\newcommand{\xBc}{\langle}
\newcommand{\xBe}{\rangle}

\newcommand{\xbD}{\Delta}
\newcommand{\xbF}{\Phi}
\newcommand{\xbG}{\Gamma}

\newcommand{\xbL}{\Lambda}
\newcommand{\xbO}{\Omega}
\newcommand{\xbP}{\Pi}

\newcommand{\xbS}{\Sigma}

\newcommand{\xba}{\alpha}
\newcommand{\xbb}{\beta}

\newcommand{\xbd}{\delta}
\newcommand{\xbe}{\in}
\newcommand{\xbf}{\phi}
\newcommand{\xbg}{\gamma}
\newcommand{\xbh}{\eta}

\newcommand{\xbk}{\kappa}
\newcommand{\xbl}{\lambda}
\newcommand{\xbm}{\mu}

\newcommand{\xbo}{\omega}
\newcommand{\xbp}{\pi}
\newcommand{\xbq}{\psi}
\newcommand{\xbr}{\rho}
\newcommand{\xbs}{\sigma}
\newcommand{\xbt}{\tau}

\newcommand{\xCI}{{\Big(}}
\newcommand{\xCJ}{{\Big)}}
\newcommand{\xCK}{\times}
\newcommand{\xCL}{\pm}

\newcommand{\xCN}{\neg}
\newcommand{\xCO}{ }
\newcommand{\xCQ}{\emptyset}

\newcommand{\xCd}{\approx}

\newcommand{\xCf}{\hspace{0.1em}}

\newcommand{\xCq}{\sim}

\newcommand{\xcA}{\forall}

\newcommand{\xcC}{\not\subseteq}
\newcommand{\xcD}{\not\supseteq}
\newcommand{\xcE}{\exists}

\newcommand{\xcH}{\not\Rightarrow}
\newcommand{\xcI}{\not\Leftarrow}
\newcommand{\xcJ}{\not\Leftrightarrow}
\newcommand{\xcK}{\not\leq}
\newcommand{\xcL}{\not\vdash}
\newcommand{\xcM}{\not\models}
\newcommand{\xcN}{\hspace{0.2em}\not\sim\hspace{-0.9em}\mid\hspace{0.8em}}
\newcommand{\xcO}{\bigvee}

\newcommand{\xcS}{\bigcap}
\newcommand{\xcT}{\bot}
\newcommand{\xcU}{\bigwedge}
\newcommand{\xcV}{\bigcup}

\newcommand{\xcX}{\Box}

\newcommand{\xcZ}{\mapsto}

\newcommand{\xca}{\infty}
\newcommand{\xcb}{\subset}
\newcommand{\xcc}{\subseteq}
\newcommand{\xcd}{\supseteq}
\newcommand{\xce}{\not\in}
\newcommand{\xcf}{\supset}
\newcommand{\xcg}{\geq}
\newcommand{\xch}{\Rightarrow}
\newcommand{\xci}{\Leftarrow}
\newcommand{\xcj}{\Leftrightarrow}
\newcommand{\xck}{\leq}
\newcommand{\xcl}{\vdash}
\newcommand{\xcm}{\models}
\newcommand{\xcn}{\hspace{0.2em}\sim\hspace{-0.9em}\mid\hspace{0.58em}}

\newcommand{\xco}{\vee}
\newcommand{\xcp}{\rightarrow}

\newcommand{\xcr}{\leftrightarrow}
\newcommand{\xcs}{\cap}
\newcommand{\xcu}{\wedge}
\newcommand{\xcv}{\cup}

\newcommand{\xcz}{\Box}

\newcommand{\xDC}{\hspace{2em}}
\newcommand{\xDH}{\item }

\newcommand{\xDN}{\ominus}
\newcommand{\xDO}{\circ}

\newcommand{\xDc}{\ll}

\newcommand{\xdA}{\mbox{\boldmath$A$}}

\newcommand{\xdC}{\mbox{\boldmath$C$}}
\newcommand{\xdD}{\mbox{\boldmath$D$}}

\newcommand{\xdL}{\mbox{\boldmath$L$}}

\newcommand{\xdO}{\mbox{\boldmath$O$}}
\newcommand{\xdP}{\mbox{\boldmath$P$}}

\newcommand{\xdR}{\Re}

\newcommand{\xda}{{\cal A}}

\newcommand{\xdc}{{\cal C}}
\newcommand{\xdd}{{\cal D}}

\newcommand{\xdf}{{\cal F}}

\newcommand{\xdi}{{\cal I}}
\newcommand{\xdj}{{\cal J}}

\newcommand{\xdl}{{\cal L}}
\newcommand{\xdm}{{\cal M}}
\newcommand{\xdn}{{\cal N}}
\newcommand{\xdo}{{\cal O}}
\newcommand{\xdp}{{\cal P}}
\newcommand{\xdr}{{\cal R}}

\newcommand{\xdu}{{\cal U}}

\newcommand{\xdx}{{\cal X}}
\newcommand{\xdy}{{\cal Y}}
\newcommand{\xdz}{{\cal Z}}

\newcommand{\xEH}{ & }
\newcommand{\xEI}{\begin{itemize}}
\newcommand{\xEJ}{\end{itemize}}
\newcommand{\xEP}{ \\ }

\newcommand{\xET}{\%}

\newcommand{\xEc}{\not<}
\newcommand{\xEd}{\neq}

\newcommand{\xEh}{\begin{enumerate}}

\newcommand{\xEj}{\end{enumerate}}

\newcommand{\xeA}{\nabla}
\newcommand{\xeB}{\not\prec}
\newcommand{\xeC}{\not\preceq}

\newcommand{\xeY}{\triangle}

\newcommand{\xeb}{\prec}
\newcommand{\xec}{\preceq}

\newcommand{\xee}{\succ}

\newcommand{\xes}{\sqsubseteq}

\newcommand{\xex}{\upharpoonright}

\newcommand{\xFO}{\parallel}

\newcommand{\xfA}{\mid}
\newcommand{\xfB}{\uparrow}

\newcommand{\xfI}{\mbox{I}}

\newcommand{\xfN}{\Uparrow}

\newcommand{\xfb}{\downarrow}

\newcommand{\Xl}{\ldots}

\newcommand{\ol}{\overline}
\newcommand{\ul}{\underline}

\newcommand{\wt}{\overbrace}

\newcommand{\xssc}{\scriptsize}

\newcommand{\xssB}{\scriptsize}

\newcommand{\bl}{\begin{lemma} \rm}
\newcommand{\el}{\end{lemma}}
\newcommand{\br}{\begin{remark} \rm}
\newcommand{\er}{\end{remark}}
\newcommand{\be}{\begin{example} \rm}
\newcommand{\ee}{\end{example}}
\newcommand{\bco}{\begin{corollary} \rm}
\newcommand{\eco}{\end{corollary}}
\newcommand{\bc}{\begin{claim} \rm}
\newcommand{\ec}{\end{claim}}
\newcommand{\bfa}{\begin{fact} \rm}
\newcommand{\efa}{\end{fact}}
\newcommand{\bp}{\begin{proposition} \rm}
\newcommand{\ep}{\end{proposition}}
\newcommand{\bd}{\begin{definition} \rm}
\newcommand{\ed}{\end{definition}}
\newcommand{\bcs}{\begin{construction} \rm}
\newcommand{\ecs}{\end{construction}}
\newcommand{\bcd}{\begin{condition} \rm}
\newcommand{\ecd}{\end{condition}}
\newcommand{\bt}{\begin{theorem} \rm}
\newcommand{\et}{\end{theorem}}
\newcommand{\bn}{\begin{notation} \rm}
\newcommand{\en}{\end{notation}}
\newcommand{\bfi}{\begin{bild} \rm}
\newcommand{\efi}{\end{bild}}
\newcommand{\bsta}{\begin{statement} \rm}
\newcommand{\esta}{\end{statement}}
\newcommand{\bcom}{\begin{comment} \rm}
\newcommand{\ecom}{\end{comment}}
\newcommand{\bdia}{\begin{diagram} \rm}
\newcommand{\edia}{\end{diagram}}

\newcommand{\bfc}{\begin{figure}[htb] \begin{center}}
\newcommand{\efc}{\end{center} \end{figure}}

\documentclass{book}

\usepackage{amssymb,latexsym,epic,eepic,rotating}

\title{
Conditionals and modularity in general logics
\thanks{
Paper 381
}
}

\author{Dov M Gabbay
\thanks{
Dov.Gabbay@kcl.ac.uk, www.dcs.kcl.ac.uk/staff/dg
} \\
King's College, London
\thanks{
Department of Computer Science, King's College London, Strand,
London WC2R 2LS, UK
} \\
and \\
Bar-Ilan University, Israel
\thanks{
Department of Computer Science,
Bar-Ilan University,
52900 Ramat-Gan, Israel
} \\
and \\
University of Luxembourg
\thanks{
Computer Science and Communications,
Faculty of Sciences,
6, rue Coudenhove-Kalergi,
L-1359 Luxembourg
} \\ \\
Karl Schlechta
\thanks{
ks@cmi.univ-mrs.fr, karl.schlechta@web.de, http://www.cmi.univ-mrs.fr/ $\sim$ ks
} \\
Laboratoire d'Informatique Fondamentale de Marseille
\thanks{
UMR 6166, CNRS and Universit\'{e} de Provence,
Address: CMI, 39, rue Joliot-Curie, F-13453 Marseille Cedex 13, France
}
}

\begin{document}

\newtheorem{lemma}{Lemma}[section]
\newtheorem{theorem}[lemma]{Theorem}
\newtheorem{proposition}[lemma]{Proposition}
\newtheorem{corollary}[lemma]{Corollary}
\newtheorem{claim}[lemma]{Claim}
\newtheorem{fact}[lemma]{Fact}
\newtheorem{remark}[lemma]{Remark}
\newtheorem{definition}{Definition}[section]
\newtheorem{construction}{Construction}[section]
\newtheorem{condition}{Condition}[section]
\newtheorem{example}{Example}[section]
\newtheorem{notation}{Notation}[section]
\newtheorem{bild}{Figure}[section]
\newtheorem{comment}{Comment}[section]
\newtheorem{statement}{Statement}[section]
\newtheorem{diagram}{Diagram}[section]

\renewcommand{\labelenumi}
  {(\arabic{enumi})}
\renewcommand{\labelenumii}
  {(\arabic{enumi}.\arabic{enumii})}
\renewcommand{\labelenumiii}
  {(\arabic{enumi}.\arabic{enumii}.\arabic{enumiii})}
\renewcommand{\labelenumiv}
  {(\arabic{enumi}.\arabic{enumii}.\arabic{enumiii}.\arabic{enumiv})}

\maketitle

\setcounter{secnumdepth}{3}
\setcounter{tocdepth}{3}

\tableofcontents

%
%
%
%

%
%
%

$ \xCO $
\clearpage
\chapter{
Introduction
}

Unless said otherwise, we work in propositional logic.
\section{
The main subjects of this book
}

This text centers around the following main subjects:

 \xEh

 \xDH
The concept of modularity and independence in
 \xEI
 \xDH
classical logic
 \xDH
non-monotonic and other non-classical logic
 \xEJ
and the consequences on
 \xEI
 \xDH
(syntactic and semantic) interpolation and
 \xDH
language change
 \xEJ
In particular, we will show the connection between interpolation for
non-monotonic logic and manipulation of an abstract notion of size.

Modularity is, for us, and essentially, the ability to put partial results
achieved independently together for a global result.

 \xDH
A uniform picture of conditionals, including

 \xEI
 \xDH
many-valued logics, and
 \xDH
structure on the language elements themselves (in contrast to structure on
the model set) and on the truth value set
 \xEJ

 \xDH
Neighbourhood semantics, their connection to independence, and their
common
points and differences for various logics, e.g.,
 \xEI
 \xDH
for defaults and deontic logic
 \xDH
for the limit version of preferential logics
 \xDH
and for general approximation.
 \xEJ
 \xEj

These subjects are not always isolated from one another, and we will
sometimes
have to go back and forth between them. For instance, a structure on the
language elements can itself be given in a modular way, and this then has
influence on the modularity of the structure on the model set.
Independence seems to be a core idea of logic, as logic is supposed to
give
the basics of reasoning, so we will not assume any
``hidden'' connections - everything which is not made explicitely otherwise,
will be assumed to be independent.

The main connections between the concepts investigate in this book are
illustrated by
Diagram \ref{Diagram Mod-Over} (page \pageref{Diagram Mod-Over}).
The left hand side concerns non-monotonic logic, the right hand side
monotonic
or classical logic. The upper part concerns mainly semantics, the lower
part
syntax.

Independence is at the core. It can be generated by Hamming distances and
relations, and can be influenced by stuctures on the language and the
truth
values. Independence is expressed by the very definition of semantics in
classical logic - a formula depends only on the value of the variables
occurring in the formula - and by suitable multiplication laws for
abstract
size in the non-monotonic case. Essentially, by these laws, a product of
two
sets is big iff the components are big. In neighbourhood semantics,
independence
is expressed by independent satisfiability of
``close'' or ``good'' along several criteria.

Semantical interpolation is the existence of
``simple'' model sets $X$ between (in the two-valued case: between by set
inclusion)
the left and the right hand model sets: $ \xbf \xcm \xbq $ results in $M(
\xbf) \xcc X \xcc M(\xbq).$
``Simple'' means here that $X$ is restricted only in the parts where both
$M(\xbf)$ and
$M(\xbq)$ are restricted, and otherwise the full product, i.e., all
truth values
may be assumed.
A - for the authors - surprising result was that monotonic logic and
antitonic
logic $ \xCf always$ have semantical interpolation. This results from the
independent definition of validity. The same is not necessarily true
for full non-monotonic logic (note, that we have $ \xbf \xcn \xbq $ iff
$M(\xbm (\xbf)) \xcc M(\xbq),$
where $M(\xbm (\xbf)) \xcc M(\xbf),$ so we have a combined downward
and upward movement:
$M(\xbf)$ - $M(\xbm (\xbf))$ - $M(\xbq)).$ The reason is that
abstract size (the sets of
``normal'' or ``important'' elements) need not be defined in an independent
way.

Semantical interpolation results also in syntactic interpolation, if the
laguage
and the operators are sufficiently rich to express these semantical
interpolants. This holds both for monotonic and non-monotonic logic.

Many aspects of independence can also be illustrated by language change.
(We did not include this in the diagram, as it might have become too
complicated to read otherwise.)

First, consider classical logic.
Let $ \xbf:=p,$ $p$ a propositional language. Consider first $L,$ with
propositional
variable just $p,$ and then $L',$ with variables $p,q.$ As long as we
know the
value of $p$ in a model $m,$ it is irrelevant whether $m$ is also defined
for
$q$ or not, and if it is, what value $q$ has. This sounds utterly trivial,
but
it is not, it is a profound conceptual property, and the whole idea of
tables of truth values is based on it. (Unfortunately, this is usually not
told to beginners, so they learn a technique (truth value tables) without
understanding that there $ \xCf is$ a problem, and how it is solved.)
This independence property is all we need to show semantical
interpolation, also
in many-valued logics,
see Proposition \ref{Proposition Gin-Pr-Int} (page \pageref{Proposition
Gin-Pr-Int}).

In non-monotonic logic, we do $ \xCf not$ have a priori such an
independence property.
Consider again $L$ and $L' $ as above, and a preferential structure. In
$L,$ the
$p-$model might be better than the $ \xCN p-$model, but this does not mean
that
in $L',$ $p-$models are better than $ \xCN p-$models. So we might have,
e.g.,
$ \xBc p \xBe  \xeb  \xBc  \xCN p \xBe $ in $M(L),$ but $ \xBc  \xCN p,q \xBe 
\xeb  \xBc  \xCN p \xCN q \xBe
\xeb  \xBc p,q \xBe  \xeb  \xBc p \xCN q \xBe,$ to give an
example. So, in $L,$ we have $TRUE \xcn p,$ in $L',$ we have $TRUE \xcn
\xCN p.$ (In terms of
abstract size, $\{ \xBc p \xBe \}$ is the smallest big subset of $M(L),$ and
$\{ \xBc  \xCN p,q \xBe \}$ the smallest big subset of $M(L'),$ see
Chapter \ref{Chapter Size-Laws} (page \pageref{Chapter Size-Laws}).) Thus,
language matters. Of course, we
can impose
additional restrictions on the various relations $ \xeb $ in $L$ and $L'
,$ e.g.,
that $p$ has precedence, like if $p \xeb \xCN p$ in $L,$ then, no matter
what the
$q-$value is in $L',$ $p-$models will be better. But this is, again,
conceptually
non-trivial, and here it might be false, contrary to the classical
case above.

We did not pursue this point in detail, yet it is very important, and
should be considered in further research - and borne in mind. For some
remarks, see Section \ref{Section Mul-Mul} (page \pageref{Section Mul-Mul})  and
Table \ref{Table Mul-Laws} (page \pageref{Table Mul-Laws}).

$ \xCO $

\vspace{10mm}

\begin{diagram}

\label{Diagram Mod-Over}
\index{Diagram Mod-Over}

\centering
\setlength{\unitlength}{1mm}
{\renewcommand{\dashlinestretch}{30}
\begin{picture}(160,190)(0,0)

\put(10,160){
{\bf Non-monotonic}
}
\put(10,156){
{\bf logic}
}

\put(110,160){
{\bf Monotonic logic}
}

\put(57,130){
Structure on language,
}
\put(63,126){
truth values
}
\path(72,122)(72,95)
\path(73,97.8)(72,95)(71,97.8)

\put(10,110){
Hamming relation,
}
\put(10,106){
Hamming distance
}
\path(42,108)(64,94)
\path(62.2,96.4)(64,94)(61.1,94.6)

\put(10,90){
Laws about size
}
\path(22,88)(22,75)
\path(22.6,76.7)(22,75)(21.4,76.7)

\put(60,90){
Independence
}
\path(58,91)(38,91)
\path(40.7,90)(38,91)(40.7,92)
\path(85,91)(108,91)
\path(105.2,92)(108,91)(105.2,90)
\path(72,88)(72,75)
\path(72.6,76.7)(72,75)(71.4,76.7)

\put(110,90){
Classical semantics
}
\path(122,88)(122,75)
\path(122.6,76.7)(122,75)(121.4,76.7)

\put(10,70){
Semantical
}
\put(10,66){
interpolation
}
\path(22,64)(22,55)
\path(22.4,56.2)(22,55)(21.6,56.2)

\put(110,70){
Semantical
}
\put(110,66){
interpolation
}
\path(122,64)(122,55)
\path(122.4,56.2)(122,55)(121.6,56.2)

\put(10,50){
Syntactical
}
\put(10,46){
interpolation
}

\put(110,50){
Syntactical
}
\put(110,46){
interpolation
}

\put(60,70){
Modularity in
}
\put(60,66){
neighbourhoods
}

\put(60,10){
Expressivity
}
\put(60,6){
of language
}
\path(69,15)(25,44)
\path(26.8,41.6)(25,44)(27.9,43.3)
\path(72,15)(120,44)
\path(117.1,43.4)(120,44)(118.1,41.7)

\put(40,1){
{\bf Connections between main concepts}
}

\end{picture}
}

\end{diagram}

\vspace{4mm}

$ \xCO $

We now give a short introduction to these main subjects.
\section{
Overview of this introduction
}

In the next sections, we give an introduction to the following chapters of
the book. In
Section \ref{Section Abstract-View} (page \pageref{Section Abstract-View}), we
try to give an abstract
definition of
independence and modularity (limited to our purposes). We conclude this
chapter
with remarks on where we used previously published material (basic
definitions
etc.), and acknowledgements.
\section{
Basic definitions
}

This chapter is relatively long, as we use a number of more or less
involved concepts, which have to be made precise. In addition, we
also want to put our work a bit more in perspective, and make it
self-contained, for the convenience of the reader.
Most of the material of this chapter (unless marked as
``new'') was published previously, see
 \cite{Sch04},
 \cite{GS08b},
 \cite{GS08c},
 \cite{GS09a}, and
 \cite{GS08f}.

We begin with basic algebraic and logical definitions, including in
particular
many laws of non-monotonic logics, in their syntactic and semantic
variants,
showing also the connections between both sides, see
Definition \ref{Definition Log-Cond-N} (page \pageref{Definition Log-Cond-N}) 
and the tables
Table \ref{Table Base2-Rules-Def-Conn-1} (page \pageref{Table
Base2-Rules-Def-Conn-1})  and
Table \ref{Table Base2-Rules-Def-Conn-2} (page \pageref{Table
Base2-Rules-Def-Conn-2}).

It seems to be a little known result that even the classical operators
permit an unusual interpretation in the infinite case, but we
claim no originality, see
Example \ref{Example Non-Standard} (page \pageref{Example Non-Standard}).

We would like to emphasize the importance of the definability preservation
(dp) property. In the infinite case, not all model sets $X$ are definable,
i.e.,
there need not necessarily be some formula $ \xbf $ or theory $T$ such
that $X=M(\xbf)$ -
the models of $ \xbf $ - or $X=M(T)$ - the models of $T.$
It is by no means evident that a model choice function $ \xbm,$ applied
to a
definable model set, gives us back again a definable model set (is
definability
preserving, or in short, dp). If $ \xbm $ does
not have this property, some representation results will not hold, which
hold if $ \xbm $ is dp, and representation results become much more
complicated,
see  \cite{Sch04} for positive and for impossibility results.
In our present context, definability is again an important concept.
Even if we have semantic interpolation, if language and operators are not
strong enough, we cannot define the semantic interpolants, so we have
semantic, but not syntactic interpolation. Examples are found in
finite Goedel logics, see
Section \ref{Section Mon-Synt-Int} (page \pageref{Section Mon-Synt-Int}).
New operators guaranteeing the definability of particularly interesting,
``universal'' interpolants, see
Definition \ref{Definition Univ-Int} (page \pageref{Definition Univ-Int}), are
discussed in
Section \ref{Section Mon-Interpol-Int} (page \pageref{Section Mon-Interpol-Int})
.
They are intricately related to the existence of
conjunctive and disjunctive normal forms, as discussed in
Section \ref{Section Analoga} (page \pageref{Section Analoga}).

We conclude this part with a - to our knowledge - unpublished result
that we can define only countably many inconsistent formulas, see
Example \ref{Example Co-Ex-Inf} (page \pageref{Example Co-Ex-Inf}). (The
question is due to D.Makinson.)

We then give a detailed introduction into the basic concepts
of many-valued logics, again, as readers might not be so familiar with
the generalizations from 2-valued to many-valued logic. In particular, the
nice correspondence between 2-valued functions and sets does not hold any
more, so we have to work with arbitrary functions, which give values to
models.
We have to re-define
what a definable model ``set'' is, and what semantical interpolation
means for many-valued logic. A formula $ \xbf $ defines such a model value
function
$f_{ \xbf },$
and we call a model value function $f$ definable iff there is some formula
$ \xbf $
such that $f=f_{ \xbf }.$
Table \ref{Table Gin-Not-Def} (page \pageref{Table Gin-Not-Def})  gives an
overview.

We then give an introduction to preferential structures and the logic
they define. These structures are among the best examined semantics for
non-monotonic logics, and
Chapter \ref{Chapter Size-Laws} (page \pageref{Chapter Size-Laws})  is also
based on the investigation of
such structures.
We first introduce the minimal variant, and then the limit variant. The
first variant is the ususal one, the second is needed to deal with
cases where there are no minimal models, due to infinite descending
chains.
(The first variant was introduced by Y.Shoham in
 \cite{Sho87b}, the second variant by P.Siegel et al. in
 \cite{BS85}. It should, however, be emphasized, that
preferential models were introduced as a semantics for deontic
logic long before they were investigated as a semantics for
non-monotonic logic, see
 \cite{Han69}).
The limit variant was further investigated in
 \cite{Sch04}, and we refer the reader there for representation
and impossibility results.
An overview of representation results for the minimal variant is given
in Table \ref{Table Base1-Pref-Rep} (page \pageref{Table Base1-Pref-Rep}).

We introduce a new concept
in this section on preferential structures,
``bubble structures'', which, we think, present a useful tool for
abstraction,
and are a semantic variant of independence in preferential structures.
Here, we have a global preferential structure between subsets
(``bubbles'') of the model set, and a fine scale structure inside those
subsets. Seen from the outside, all elements of a bubble behave the
same way, so the whole set can be treated as one element, on the
inside, we see a finer structure.

Moreover, new material on many-valued preferential structures is included.

We then go into details in the section on IBRS, introduced by D.Gabbay,
see  \cite{Gab04}, and further investigated in
 \cite{GS08b} and  \cite{GS08f}, as they
are not so much common knowledge. We also discuss here if and how
the limit version of preferential structures might be applied
to reactive structures.

We then present theory revision, as introduced by Alchorron, Gardenfors,
and Makinson, see  \cite{AGM85}. Again,
we also build on previous results by the authors,
when we discuss distance based revision, introduced by
Lehmann, Magidor, and Schlechta, see
 \cite{LMS95},
 \cite{LMS01}, and elaborated in
 \cite{Sch04}.
We also include a short paragraph on new material for theory revision
based on
many-valued logic.
\section{
Towards a uniform picture of conditionals
}

In large parts, this chapter should rather be seen more as a sketch for
future
work, than a fully elaborated theory.
\subsection{
Discussion and classification
}

It seems difficult to say what is not a conditional. The word
``condition'' suggests something like
``if  \Xl, then  \Xl'', but as the condition might be hidden in an
underlying structure, and not expressed in the object language, a
conditional
might also be an unary operator, e.g., we may read the consequence
relation
$ \xcn $ of preferential structure as
``under the condition of normality''.

Moreover, as shown at the beginning of
Section \ref{Section Mod-Uni-Int} (page \pageref{Section Mod-Uni-Int}),
Example \ref{Example Libi-Cond} (page \pageref{Example Libi-Cond}), it seems
that one can define
new conditionals ad libitum, binary, ternary, etc.

Thus, the best seems to be to say that a conditional is just any operator.
Negation, conjunction, etc., are then included, but excluded from the
discussion, as we know them well.

The classical connectives have a semantics in the boolean set operators,
but
there are other operators, like the $ \xbm -$functions of preferential
logic which
do not correspond to any such operator, and might even not preserve
definability
in the infinite case (see
Definition \ref{Definition Log-Base} (page \pageref{Definition Log-Base})). It
seems more promising to order
conditionals
by the properties of their model choice functions, e.g., whether those
functions
are idempotent,
etc., see Section \ref{Section Choice-Prop} (page \pageref{Section Choice-Prop})
.

Many conditionals can be based on binary relations, e.g. modal
conditionals
on accessibility relations, preferential consequence relations on
preference relations, counterfactuals and theory revision on distance
relations, etc. Thus, it is promising to look at those relations, and
their
properties to bring more order into the vast field of conditionals.
D.Gabbay introduced reactive structures
(see, e.g.,  \cite{Gab04}), and added supplementary expressivity to
structures based on binary relations, see
 \cite{GS08b} and  \cite{GS08f}. In particular, it
was shown there that we can have cumulativity without the
basic properties of preferential structures (e.g., OR).
This is discussed in
Section \ref{Section Cond-Bin} (page \pageref{Section Cond-Bin}).
\subsection{
Additional structure on language and truth values
}

Normally, the language elements (propositional variables) are not
structured. This is somewhat surprising, as, quite often, one variable
will be more important than another. Size or weight might often be more
important than colour for physical objects, etc. It is probably the
mathematical
tradition which was followed too closely. One of the authors gave a
semantics
to theory revision using a measure on language elements in
 \cite{Sch91-1} and  \cite{Sch91-3},
but, as far as we know, the subject was not treated in a larger
context so far. The present book often works with independence of language
elements, see in particular
Chapter \ref{Chapter Mod-Mon-Interpol} (page \pageref{Chapter Mod-Mon-Interpol})
 and
Chapter \ref{Chapter Size-Laws} (page \pageref{Chapter Size-Laws}), and Hamming
type relations and
distances
between models, where it need not be the case that all variables have
the same weight. Thus, it is obvious to discuss this subject in the
present text.
It can also be fruitful to discuss sizes of subsets of the set of
variables,
so we may, e.g., neglect differences to classical logic if they concern
only
a ``small'' set of propositional variables.

On the other hand, classical truth values have a natural order,
$FALSE<TRUE,$ and we will sometimes work with more than 2 truth values,
see in particular
Chapter \ref{Chapter Mod-Mon-Interpol} (page \pageref{Chapter Mod-Mon-Interpol})
, but also
Section \ref{Section EQ} (page \pageref{Section EQ}). So there is a natural
question: do we also
have a total
order, or a boolean order, or another order on those sets of truth values?
Or: Is there a distance between truth values, so that a change from value
$ \xCf a$ to
value $b$ is smaller than a change from $ \xCf a$ to $c?$

There is a natural correspondence between semantical structures and truth
values, which is best seen by an example: Take finite (intuitionistic)
Goedel
logics,
see Section \ref{Section Finite-Goedel} (page \pageref{Section Finite-Goedel}),
say, for simplicity with two worlds. Now, $ \xbf $ may hold nowhere,
everywhere, or only in the second world (called
``there'', in contrast to ``here'',
the first world). Thus, we can express the same situation by three truth
values: 0 for nowhere, 1 for only
``there'', 2 for everywhere.

In Section \ref{Section Softening} (page \pageref{Section Softening}), we will
make some short remarks on
``softening'' concepts, like neglecting
``small'' fragments of a language, etc. This way, we can define, e.g.,
``soft'' interpolation, where we need a small set of variables which
are not in both formulas.

Inheritance systems,
(see, e.g.,  \cite{TH89},  \cite{THT86},
 \cite{THT87},  \cite{TTH91},  \cite{Tou86},
also  \cite{Sch93} and  \cite{Sch97},
 \cite{GS08e},  \cite{GS08f}),
present many aspects of independence,
(see Section \ref{Section Inher} (page \pageref{Section Inher})).
Thus, if two nodes are not connected by valid paths, they may have
very different languages, as language elements have to be inherited,
otherwise, they are undefined. In addition, $ \xCf a$ may inherit from $b$
property $c,$
but not property $d,$ as we have a contradiction to $d$ (or, even $ \xCN
d)$ via
a different node $b'.$ Theses are among the aspects which make them
natural
for common sense reasoning,
but also quite different from traditional logics.
\subsection{
Representation for general revision, update, and counterfactuals
}

Revision (see  \cite{AGM85}, and the discussion in
Section \ref{Section Mod-TR} (page \pageref{Section Mod-TR})),
update (see  \cite{KM90}),
and counterfactuals
(see  \cite{Lew73} and  \cite{Sta68})
are special forms of conditionals, which
received much interest in the artificial intelligence community.
Explicitly or implicitly
(see  \cite{LMS95},  \cite{LMS01}),
they are based on a distance based semantics, working with
``closest worlds''.
In the case of revision, we look at those worlds which are closest to the
present $ \xCf set$ of worlds, in update and counterfactual, we look
from each present world $ \xCf individually$ to the closest worlds, and
then take
the union. Obviously, the formal properties may be very different in the
two
cases.

There are two obvious generalizations possible, and sometimes necessary.
First, ``closest'' worlds need not exist, there may be infinite descending
chains of distances without minimal elements. Second, a distance or ranked
order may force too many comparisons, when two distances or elements may
just simply not be comparable. We address representation problems for
these generalizations:

 \xEh

 \xDH
We first generalize the notion of distance for revision semantics
in Section \ref{Section Semantic-TR} (page \pageref{Section Semantic-TR}).
We mostly consider symmetrical distances, so $d(a,b)=d(b,a),$ and we
work with equivalence classes $[a,b].$ Unfortunately, one of the main
tools in  \cite{LMS01}, a loop condition, does not work any more, it
is
too close to rankedness.

We will have to work more in the spirit of general and
smooth preferential structures to obtain representation. Unfortunately,
revision does not allow many observations
(see  \cite{LMS01}, and, in particular, the impossibility results
for
revision (``Hamster Wheels'') discussed in
 \cite{Sch04}), so all we have
(see Section \ref{Section TR-Main-Rep} (page \pageref{Section TR-Main-Rep}))
are results which use more
conditions than
what can be observed from revision observations.
This problem is one of principles: we showed in
 \cite{GS08a}, see also  \cite{GS08f}, that cumulativity
suffices only
to guarantee smoothness of the structure if the domain is
closed under finite unions. But the union of two products need not
be a product any more.

To solve the problem, we use a technique employed in
 \cite{Sch96-1}, using ``witnesses'' to
testify for the conditions.

 \xDH
We then discuss the limit version (when there are no minimal distances)
for theory revision.

 \xDH
In Section \ref{Section Semantic-Up} (page \pageref{Section Semantic-Up}), we
turn to
generalized update and counterfactuals.
To solve this problem, we use a technique invented in
 \cite{MS90}, and adapt it to our situation.
The basic idea is very simple: we begin (simplified) with some world $x,$
and
arrange the other worlds around $x,$ as $x$ sees them, by their relative
distances.
Suppose we consider now one those worlds, say $y.$ Now we arrange the
worlds
around $y,$ as $y$ sees them. If we make all the new distances smaller
than the
old ones, we ``cannot look back'', etc. We continue this construction
unboundedly
(but finitely) often. If we are a little careful, everyone will only see
what he
is supposed to see. In a picture, we construct galaxies around a center,
then
planets around suns, moons around planets, etc.
The resulting construction is an $ \xda -$ranked structure, as discussed
in
 \cite{GS08d}, see also  \cite{GS08f}.

 \xDH
In Section \ref{Section TR-Up-Synt} (page \pageref{Section TR-Up-Synt}), we
discuss the corresponding
syntactic
conditions, using again ideas from
 \cite{Sch96-1}.

 \xEj
\section{
Interpolation
}
\subsection{
Introduction
}

\label{Section Sin-Introduction}

The two chapters
Chapter \ref{Chapter Mod-Mon-Interpol} (page \pageref{Chapter Mod-Mon-Interpol})
 and
Chapter \ref{Chapter Size-Laws} (page \pageref{Chapter Size-Laws})
are probably the core of the present book.

We treat very general interpolation problems for monotone and antitone,
2-valued and many-valued logics in
Chapter \ref{Chapter Mod-Mon-Interpol} (page \pageref{Chapter Mod-Mon-Interpol})
, splitting the question in two
parts,
``semantic interpolation'' and
``syntactic interpolation'', show that the first problem, existence
of semantic interpolation, has a simple and general answer, and reduce
the second question, existence of syntactic interpolation to a
definability problem. For the latter, we examine the concrete example
of finite Goedel logics. We can also show that the semantic problem
has two ``universal'' solutions, which depend only on one formula and
the shared variables.

In Chapter \ref{Chapter Size-Laws} (page \pageref{Chapter Size-Laws}), we
investigate three variants of
semantic
interpolation for non-monotonic logics, in syntactic shorthand of the
types
$ \xbf \xcn \xba \xcl \xbq,$ $ \xbf \xcl \xba \xcn \xbq,$ and $ \xbf
\xcn \xba \xcn \xbq,$ where $ \xba $ is the interpolant, and
see that two variants are closely related to multiplication laws about
abstract
size, defining (or originating from) the non-monotonic logics.
The syntactic problem is analogous to that of the monotonic case.
\subsubsection{
Background
}

Interpolation for classical logic is well-known,
see  \cite{Cra57}, and there are also
non-classical logics for which interpolation has been shown, e.g., for
Circumscription, see
 \cite{Ami02}. An extensive overview of interpolation is found
in  \cite{GM05}. Chapter 1 of this book
 \cite{GM05} gives a survey and a discussion and
the chapter puts forward that interpolation can be viewed in many
different ways
and indeed
11 points of view of interpolation are discussed. The present text
presents the
semantic interpolation, this is a new 12th point of view.
\subsection{
Problem and Method
}

In classical logic, the problem of interpolation is to find for two
formulas $ \xbf $ and $ \xbq $ such that $ \xbf \xcl \xbq $ a
``simple'' formula $ \xba $ such that $ \xbf \xcl \xba \xcl \xbq.$
``Simple'' is defined as:
``expressed in the common language of $ \xbf $ and $ \xbq $''.

Working on the semantic level has often advantages:

 \xEI

 \xDH results are robust under logically equivalent reformulations

 \xDH in many cases, the semantic level allows an easy reformulation as an
algebraic problem, whose results can be generalized to other situations

 \xDH we can split a problem in two parts: a semantical problem, and the
problem
to find a syntactic counterpart (a definability problem)

 \xDH the semantics of many non-classical logics are built on relatively
few
basic notions like size, distance, etc., and we can thus make connections
to other problems and logics

 \xDH in the case of preferential and similar logics, the very definition
of
the logic is already semantical (minimal models), so it is very natural
to proceed on this level.

 \xEJ

This strategy - translate to the semantic level, do the main work there,
and
then translate back - has proved fruitful also in the present case.

Looking back at the classical interpolation problem, and translating it
to the semantic level, it becomes: Given $M(\xbf) \xcc M(\xbq)$ (the
models sets of $ \xbf $ and
$ \xbq),$ is there a
``simple'' model set $ \xCf A$ such that $M(\xbf) \xcc A \xcc M(\xbq)?$
Or, more generally, given
model sets $X \xcc Y,$ is there ``simple'' $ \xCf A$ such that $X \xcc A
\xcc Y?$
Of course, we have to define in a natural way,
what ``simple'' means in our context.
This is discussed below in
Section \ref{Section Intro-Mon-Sema} (page \pageref{Section Intro-Mon-Sema}).

Our separation of the semantic from the syntactic question
pays immediately:

 \xEh

 \xDH
We see that monotonic (and antitonic) logics $ \xCf always$ have a
semantical
interpolant. But this interpolant may not be definable syntactically.

 \xDH
More precisely, we see that there is a whole interval of interpolants in
above situation.

 \xDH
We also see that our analysis generalizes immediately to many valued
logics, with the same result (existence of an interval of interpolants).

 \xDH
Thus, the question remains: under what conditions does a syntactic
interpolant exist?

 \xDH
In non-monotonic logics, our analysis reveals a deep connection between
semantic interpolation and questions
about (abstract) multiplication of (abstract) size.

 \xEj
\subsection{
Monotone and antitone semantic and syntactic interpolation
}

We consider here the semantic property of monotony or antitony, in the
following sense (in the two-valued case, the generalization to the
many-valued case is straightforward):

Let $ \xcl $ be some logic such that $ \xbf \xcl \xbq $ implies $M(\xbf)
\xcc M(\xbq)$
(the monotone case) or $M(\xbq) \xcc M(\xbf)$ (the antitone case).

In the many-valued case, the corresponding property is that $ \xcp $ (or $
\xcl)$
respects $ \xck,$ the order on the truth values.
\subsubsection{
Semantic interpolation
}

\label{Section Intro-Mon-Sema}

The problem (for simplicity, for the 2-valued case) reads now:

If $M(\xbf) \xcc M(\xbq)$ (or, symmetrically $M(\xbq) \xcc M(\xbf
)),$ is there a
``simple'' model set $ \xCf A,$ such that $M(\xbf) \xcc A \xcc M(\xbq),$
or $M(\xbq) \xcc A \xcc M(\xbf).$
Obviously, the problem is the same in both cases.
We will see that such $ \xCf A$ will always exist, so all such logics have
semantic interpolation (but not necessarily also syntactic interpolation).
We turn to the main conceptual problem, the definition of ``simple''.

The main conceptual problem is to define
``simple model set''. We have to look at the
syntactic problem for guidance. Suppose $ \xbf $ is defined using
propositional
variables $p$ and $q,$
$ \xbq $ using $q$ and $r.$ $ \xba $ has to be defined using only $q.$
What are the models of
$ \xba?$ By the very definition of validity in classical logic, neither
$p$ nor $r$
have any influence on whether $m$ is a model of $ \xba $ or not. Thus, if
$m$ is a model
of $ \xba,$ we can modify $m$ on $p$ and $r,$ and it will still be a
model. Classical
models are best seen as functions from the set of propositional variables
to
$\{true,false\},$ $\{t,f\},$ or so. In this terminology, any $m$ with $m
\xcm \xba $ is
``free'' to choose the value for $p$ and $r,$ and we can write the model set
A
of $ \xba $ as $\{t,f\} \xCK M_{q} \xCK \{t,f\},$ where $M_{q}$ is the set
of values for $q$ $ \xba -$models
may have $(\xCQ,$ $\{t\},$ $\{f\},$ $\{t,f\}).$

So, the semantic interpolation problem is to find sets which may be
restricted
on the common variables,
but are simply the Cartesian product of the possible values for the other
variables. To summarize: Let two model sets $X$ and $Y$ be given, where
$X$ itself
is restricted on variables $\{p_{1}, \Xl,p_{m}\}$ (i.e. the Cartesian
product for the rest),
$Y$ is restricted on $\{r_{1}, \Xl,r_{n}\},$ then we have to find a model
set $ \xCf A$ which
is restricted only on $\{p_{1}, \Xl,p_{m}\} \xcs \{r_{1}, \Xl,r_{n}\},$
and such that
$X \xcc A \xcc Y,$ of course.

Formulated this way, our approach, the problem and its solution, has two
trivial generalizations:

 \xEI

 \xDH for multi-valued logics
we take the Cartesian product of more than just $\{t,f\}.$

 \xDH $ \xbf $ may be the hypothesis, and $ \xbq $ the consequence, but
also vice versa,
there is no direction in the problem. Thus, any result for classical
logic carries over to the
core part of, e.g., preferential logics.

 \xEJ

The main result for the situation with $X \xcc Y$ is that there is always
such a
semantic interpolant $ \xCf A$ as described above
(see Proposition \ref{Proposition Sin-Interpolation} (page \pageref{Proposition
Sin-Interpolation})  for a simple case,
and Proposition \ref{Proposition Gin-Pr-Int} (page \pageref{Proposition
Gin-Pr-Int})  for the full picture).
Our proof works also for
``parallel interpolation'', a concept introduced by
Makinson et al.,  \cite{KM07}.

We explain and quote the latter result.

Suppose we have $f,g:M \xcp V,$ where, intuitively, $M$ is the set of all
models, and
$V$ the set of all truth values. Thus, $f$ and $g$ give to each model a
truth value,
and, intuitively, $f$ and $g$ each code a model set, assigning to $m$ TRUE
iff $m$ is
in the model set, and FALSE iff not. We further assume that there is an
order
on the truth value set $V.$ $ \xcA m \xbe M(f(m) \xck g(m))$ corresponds
now to $M(\xbf) \xcc M(\xbq),$
or $ \xbf \xcl \xbq $ in classical logic. Each model $m$ is itself a
function from $L,$ the
set of propositional variables to $V.$ Let now $J \xcc L.$ We say that $f$
is
insensitive to $J$ iff the values of $m$ on $J$ are irrelevant: If
$m \xex (L-J)=m' \xex (L- \xCf J),$
i.e., $m$ and $m' $ agree at least on all $p \xbe L- \xCf J,$ then
$f(m)=f(m').$ This
corresponds to the situation where the variable $p$ does not occur in the
formula $ \xbf,$ then $M(\xbf)$ is insensitive to $p,$ as the value of
any $m$ on $p$
does not matter for $m$ being a model of $ \xbf,$ or not.

We need two more definitions:

Let $J' \xcc L,$ then
$f^{+}(m_{J' }):=max\{f(m'):m' \xex J' =m \xex J' \}$ and $f^{-}(m_{J'
}):=min\{f(m'):m' \xex J' =m \xex J' \}.$

We quote now
Proposition \ref{Proposition Gin-Pr-Int} (page \pageref{Proposition Gin-Pr-Int})
:

\bp

$\hspace{0.01em}$


\label{Proposition Gin-Pr-Int-A}

Let $M$ be rich, $f,g:M \xcp V,$ $f(m) \xck g(m)$ for all $m \xbe M.$
Let $L=J \xcv J' \xcv J'',$ let $f$ be insensitive to $J,$ $g$ be
insensitive to $J''.$

Then $f^{+}(m_{J' }) \xck g^{-}(m_{J' })$ for all $m_{J' } \xbe M \xex J'
,$ and any $h:M \xex J' \xcp V$ which is
insensitive to $J \xcv J'' $ is an interpolant iff

$f^{+}(m_{J' }) \xck h(m_{J}m_{J' }m_{J'' })=h(m_{J' }) \xck g^{-}(m_{J'
})$ for all $m_{J' } \xbe M \xex J'.$

(h can be extended to the full $M$ in a unique way, as it is
insensitive to $J \xcv J''.)$

\ep

See Diagram \ref{Diagram Mon-Int} (page \pageref{Diagram Mon-Int}).
\subsubsection{
The interval of interpolants
}

Our result has an additional reading: it defines an interval of
interpolants, with lower bound $f^{+}(m_{J' })$ and upper bound
$g^{-}(m_{J' }).$
But these interpolants have a particular form. If they exist,
i.e. iff $f \xck g,$ then $f^{+}(m_{J' })$ depends only on $f$ and $J' $
(and $m),$ but $ \xCf not$ on $g,$
$g^{-}(m_{J' })$ only on $g$ and $J',$ $ \xCf not$ on $f.$ Thus, they are
universal, as we have
to look only at one function and the set of common variables.

Moreover, we will see in
Section \ref{Section Analoga} (page \pageref{Section Analoga})  that they
correspond to simple operations
on the normal
forms in classical logic. This is not surprising, as we
``simplify'' potentially complicated model sets by replacing some
coordinates with simple products. The question is, whether our logic
allows to express this simplification, classical logic does.
\subsubsection{
Syntactic interpolation
}

Recall the problem described at the beginning of
Section \ref{Section Intro-Mon-Sema} (page \pageref{Section Intro-Mon-Sema}).
We were given $M(\xbf) \xcc M(
\xbq),$ and were
looking for a ``simple'' model set $ \xCf A$ such that $M(\xbf) \xcc A
\xcc M(\xbq).$ We just saw that
such $ \xCf A$ exists, and were able to describe an interval of such $
\xCf A.$
But we have no guarantee that any such $ \xCf A$ is definable, i.e., that
there
is some $ \xba $ with $A=M(\xba).$

In classical logic, such $ \xba $ exists, see, e.g.,
Proposition \ref{Proposition Sin-Simplification-Definable} (page
\pageref{Proposition Sin-Simplification-Definable})),
but also Section \ref{Section Analoga} (page \pageref{Section Analoga}).
Basically, in classical logic, $f^{+}(m_{J' })$ and $g^{-}(m_{J' })$
correspond to
simplifications of the formulas expressed in normal form,
see Fact \ref{Fact Class-Up-Down} (page \pageref{Fact Class-Up-Down})  (in a
different notation, which we
will
explain in a moment).
This is not necessarily true in other logics, see
Example \ref{Example Sin-4-Value} (page \pageref{Example Sin-4-Value}).
(We find here again the importance of definability preservation,
a concept introduced by one of us in  \cite{Sch92}.)

If we have projections (simplifications), see
Section \ref{Section Int-Int-Mon} (page \pageref{Section Int-Int-Mon}),
we also have syntactic interpolation. At present,
we do not know whether this is a necessary condition for all natural
operators.

We can also turn the problem around, and just define suitable
operators. This is done in
Section \ref{Section Analoga} (page \pageref{Section Analoga}),
Definition \ref{Definition 2-Many-Up} (page \pageref{Definition 2-Many-Up})  and
Definition \ref{Definition 2-Many-Down} (page \pageref{Definition 2-Many-Down})
.
There is a slight problem, as one of the operands is a $ \xCf set$ of
propositional
variables, and not a formula, as usual. One, but certainly not the only
one,
possibility is to take a formula (or the corresponding model set) and
``extract'' the ``relevant'' variables from it, i.e., those, which cannot
be replaced by a product. Assume now that $f$ is one of the
generalized model ``sets'', then:

Given $f,$ define

 \xEh
 \xDH
$(f \xfB J)(m)$ $:=$ $sup\{f(m'):$ $m' \xbe M,$ $m \xex J=m' \xex J\}$
 \xDH
$(f \xfb J)(m)$ $:=$ $inf\{f(m'):$ $m' \xbe M,$ $m \xex J=m' \xex J\}$
 \xDH
$ \xbf! \xbq $ by:

$f_{ \xbf! \xbq }$ $:=$ $f_{ \xbf } \xfB (L-R(\xbq))$
 \xDH
$ \xbf? \xbq $ by:

$f_{ \xbf? \xbq }$ $:=$ $f_{ \xbf } \xfb (L-R(\xbq))$
 \xEj

We then obtain for classical logic
(see Fact \ref{Fact Class-Up-Down} (page \pageref{Fact Class-Up-Down})):

\bfa

$\hspace{0.01em}$


\label{Fact Class-Up-Down-a}

Let $J:=\{p_{1,1}, \Xl,p_{1,m_{1}}, \Xl,p_{n,1}, \Xl,p_{n,m_{n}}\}$

(1) Let $ \xbf_{i}:= \xCL p_{i,1} \xcu  \Xl  \xcu \xCL p_{i,m_{i}}$ and $
\xbq_{i}:= \xCL q_{i,1} \xcu  \Xl  \xcu \xCL q_{i,k_{i}},$
let $ \xbf:=(\xbf_{1} \xcu \xbq_{1}) \xco  \Xl  \xco (\xbf_{n} \xcu
\xbq_{n}).$
Then $ \xbf \xfB J= \xbf_{1} \xco  \Xl  \xco \xbf_{n}.$

(2) Let $ \xbf_{i}:= \xCL p_{i,1} \xco  \Xl  \xco \xCL p_{i,m_{i}}$ and $
\xbq_{i}:= \xCL q_{i,1} \xco  \Xl  \xco \xCL q_{i,k_{i}},$
let $ \xbf:=(\xbf_{1} \xco \xbq_{1}) \xcu  \Xl  \xcu (\xbf_{n} \xco
\xbq_{n}).$
Then $ \xbf \xfb J= \xbf_{1} \xcu  \Xl  \xcu \xbf_{n}.$

In a way, these operators are natural, as they simplify definable
model sets, so they can be used as a criterion of the expressive
strength of a language and logic: If $X$ is definable, and $Y$ is in some
reasonable sense simpler than $X,$ then $Y$ should also be definable.
If the language is not sufficiently strong, then we can introduce these
operators, and have also syntactic interpolation.
\subsubsection{
Finite Goedel logics
}

\efa

The semantics of finite (intuitionistic) Goedel logics is a finite chain
of
worlds, which can also be expressed by a totally ordered set of truth
values
0 \Xl n (see Section \ref{Section Finite-Goedel} (page \pageref{Section
Finite-Goedel})).
Let FALSE and TRUE be the minimal and maximal truth values.
$ \xbf $ has value false, iff it holds nowhere, and TRUE, iff it holds
everywhere,
it has value 1 iff it holds from world 2 onward, etc.
The operators are classical $ \xcu $ and $ \xco,$ negation $ \xCN $ is
defined by
$ \xCN (FALSE)=TRUE$ and $ \xCN (x)=FALSE$ otherwise.
Implication $ \xcp $ is defined by $ \xbf \xcp \xbq $ is TRUE iff $ \xbf
\xck \xbq $ (as truth values),
and the value of $ \xbq $ otherwise.

More precisely, where $f_{ \xbf }$ is the model value function of the
formula $ \xbf:$

negation $ \xCN $ is defined by:
\begin{flushleft}
\[ f_{\xCN \xbf}(m):= \left\{ \begin{array}{lcl}
TRUE \xEH iff \xEH
f_{\xbf}(m)=FALSE \xEP
\xEH \xEH \xEP
FALSE \xEH \xEH otherwise \xEP
\end{array}
\right.
\]
\end{flushleft}

implication $ \xcp $ is defined by:
\begin{flushleft}
\[ f_{\xbf \xcp \xbq}(m):= \left\{ \begin{array}{lcl}
TRUE \xEH iff \xEH
f_{\xbf}(m) \xck f_{\xbq}(m) \xEP
\xEH \xEH \xEP
f_{\xbq}(m) \xEH \xEH otherwise \xEP
\end{array}
\right.
\]
\end{flushleft}

see Definition \ref{Definition Mod-Fin-Goed} (page \pageref{Definition
Mod-Fin-Goed})  in
Section \ref{Section Finite-Goedel} (page \pageref{Section Finite-Goedel}).
We show in Section \ref{Section Example-No-Int} (page \pageref{Section
Example-No-Int})  the well-known
result that such logics for 3 worlds (and thus 4 truth values)
have no interpolation, whereas the corresponding logic for 2 worlds
has interpolation. For the latter logic, we can still find a kind of
normal form, though $ \xcp $ cannot always be reduced. At least we can
avoid
nested implications, which is not possible in the former logic for 3
worlds.

We also discuss several ``hand made'' additional operators which
allow us to define sufficiently many model sets to have syntactical
interpolation - of course, we $ \xCf know$ that we have semantical
interpolation.
A more systematic approach was discussed above, the operators $ \xbf!
\xbq $ and $ \xbf? \xbq.$
\subsection{
Laws about size and interpolation in non-monotonic logics
}

\label{Section Laws-Size-Intro}
\subsubsection{
Various concepts of size and non-monotonic logics
}

A natural interpretation of the non-monotonic rule $ \xbf \xcn \xbq $ is
that
the set of exceptional cases, i.e., those where $ \xbf $ holds, but not $
\xbq,$ is
a small subset of all the cases where $ \xbf $ holds, and the complement,
i.e.,
the set of cases where $ \xbf $ and $ \xbq $ hold, is a big subset of all
$ \xbf -$cases.

This interpretation gives an abstract semantics to non-monotonic logic,
in the sense that definitions and rules are translated to rules about
model sets, without any structural justification of those rules, as they
are given, e.g., by preferential structures, which provide structural
semantics. Yet, they are extremely useful, as they allow us to concentrate
on the essentials, forgetting about syntactical reformulations of
semantically
equivalent formulas, the laws derived from the standard proof theoretical
rules
incite to generalizations and modifications, and reveal deep connections
but
also differences. One of those insights is the connection between laws
about size and (semantical) interpolation for non-monotonic logics,
discussed in Chapter \ref{Chapter Size-Laws} (page \pageref{Chapter Size-Laws})
.

To put this abstract view a little more into perspective, we now present
three alternative systems, also working with abstract size as a semantics
for non-monotonic logics. (They will be repeated in the introduction of
Chapter \ref{Chapter Size-Laws} (page \pageref{Chapter Size-Laws}).)

 \xEI
 \xDH
the system of one of the authors for a first order setting,
published in  \cite{Sch90} and elaborated in
 \cite{Sch95-1},
 \xDH
the system of S.Ben-David and R.Ben-Eliyahu,
published in  \cite{BB94},
 \xDH
the system of N.Friedman and J.Halpern,
published in  \cite{FH96}.
 \xEJ

 \xEh

 \xDH
Defaults as generalized quantifiers:

We first recall the definition of a
``weak filter'', made official in
Definition \ref{Definition Weak-Filter} (page \pageref{Definition Weak-Filter})
:

Fix a base set $X.$
A weak filter on or over $X$ is a set $ \xdf \xcc \xdp (X),$ s.t. the
following
conditions hold:

(F1) $X \xbe \xdf $

(F2) $A \xcc B \xcc X,$ $A \xbe \xdf $ imply $B \xbe \xdf $

$(F3')$ $A,B \xbe \xdf $ imply $A \xcs B \xEd \xCQ.$

We use
weak filters on the semantical side, and add the following axioms on the
syntactical side to a FOL axiomatisation:

\index{$\xeA$}
1. $ \xeA x \xbf (x)$ $ \xcu $ $ \xcA x(\xbf (x) \xcp \xbq (x))$ $ \xcp $
$ \xeA x \xbq (x),$

2. $ \xeA x \xbf (x)$ $ \xcp $ $ \xCN \xeA x \xCN \xbf (x),$

3. $ \xcA x \xbf (x)$ $ \xcp $ $ \xeA x \xbf (x)$ and $ \xeA x \xbf (x)$ $
\xcp $ $ \xcE x \xbf (x).$

A model is now a pair, consisting of a classical FOL model $M,$ and a weak
filter over its universe. Both sides are connected by the following
definition,
where $ \xdn (M)$ is the weak filter on the universe of the classical
model $M:$

\index{$ \xBc M,\xdn (M) \xBe $}
$ \xBc M, \xdn (M) \xBe $ $ \xcm $ $ \xeA x \xbf (x)$ iff there is $A \xbe \xdn
(M)$
s.t. $ \xcA a \xbe A$ $(\xBc M, \xdn (M) \xBe $ $ \xcm $ $ \xbf [a]).$

Soundness and completeness is shown in
 \cite{Sch95-1}, see also  \cite{Sch04}.

The extension to defaults with prerequisites by restricted quantifiers is
straightforward.

 \xDH
The system of $S.$ Ben-David and $R.$ Ben-Eliyahu:

Let $ \xdn':=\{ \xdn' (A):$ $A \xcc U\}$ be a system of filters for $
\xdp (U),$ i.e. each $ \xdn' (A)$ is a
filter over A. The conditions are (in slight modification):

UC': $B \xbe \xdn' (A)$ $ \xcp $ $ \xdn' (B) \xcc \xdn' (A),$

DC': $B \xbe \xdn' (A)$ $ \xcp $ $ \xdn' (A) \xcs \xdp (B) \xcc \xdn'
(B),$

RBC': $X \xbe \xdn' (A),$ $Y \xbe \xdn' (B)$ $ \xcp $ $X \xcv Y \xbe
\xdn' (A \xcv B),$

SRM': $X \xbe \xdn' (A),$ $Y \xcc A$ $ \xcp $ $A-Y \xbe \xdn' (A)$ $
\xco $ $X \xcs Y \xbe \xdn' (Y),$

\index{$UC'$}
\index{$DC'$}
\index{$RBC'$}
\index{$SRM'$}
\index{$GTS'$}
GTS': $C \xbe \xdn' (A),$ $B \xcc A$ $ \xcp $ $C \xcs B \xbe \xdn' (B).$

 \xDH

The system of $N.$ Friedman and $J.$ Halpern:

Let $U$ be a set, $<$ a strict partial order on $ \xdp (U),$
(i.e. $<$ is transitive, and contains no cycles).
Consider the following conditions for $<:$

(B1) $A' \xcc A<B \xcc B' $ $ \xcp $ $A' <B',$

(B2) if $A,B,C$ are pairwise disjoint, then $C<A \xcv B,$ $B<A \xcv C$ $
\xcp $ $B \xcv C<A,$

(B3) $ \xCQ <X$ for all $X \xEd \xCQ,$

(B4) $A<B$ $ \xcp $ $A<B-$A,

\index{$(B1)$}
\index{$(B2)$}
\index{$(B3)$}
\index{$(B4)$}
\index{$(B5)$}
(B5) Let $X,Y \xcc A.$ If $A-X<X,$ then $Y<A-Y$ or $Y-X<X \xcs Y.$

 \xEj

The equivalence of the systems of  \cite{BB94}
and  \cite{FH96} was shown in
 \cite{Sch97-4}, see also  \cite{Sch04}.

Historical remarks:
Our own view as abstract size was inspired by the classical filter
approach,
as used e.g. in mathematical measure theory.
The first time that abstract size was related to
nonmonotonic logics was, to our knowledge,
in the second author's  \cite{Sch90} and  \cite{Sch95-1}, and,
independently, in  \cite{BB94}.
The approach to size by partial orders is
first discussed - to our knowledge - by
N.Friedman and J.Halpern, see  \cite{FH96}.
More detailed remarks can also be found in  \cite{GS08c},
 \cite{GS09a},  \cite{GS08f}.
A somewhat different approach is taken in  \cite{HM07}.

Before we introduce the connection between interpolation and
multiplicative
laws about size, we give now some comments on the laws about size
themselves.
\subsubsection{
Additive and multiplicative laws about size
}

We give here a short introduction to and some examples for additive
and multiplicative laws about size. A detailed overview is presented
in Table \ref{Table Base2-Size-Rules-1} (page \pageref{Table
Base2-Size-Rules-1}),
Table \ref{Table Base2-Size-Rules-2} (page \pageref{Table Base2-Size-Rules-2}),
and
Table \ref{Table Mul-Laws} (page \pageref{Table Mul-Laws}).
(The first two tables have to be read together, they are too
big to fit on one page.)

They show connections and how to develop a multitude of logical rules
known from nonmonotonic logics by combining a small number of principles
about size. We can use them as building blocks to construct the rules
from.
More precisely, ``size'' is to be read as
``relative size'', since it is essential to change the base sets.

In the first two tables, these principles are some basic and very natural
postulates,
$ \xCf (Opt),$ $ \xCf (iM),$ $(eM \xdi),$ $(eM \xdf),$ and a continuum
of power of the notion of
``small'', or, dually, ``big'', from $(1*s)$ to $(< \xbo *s).$
From these, we can develop the rest except, essentially, Rational
Monotony,
and thus an infinity of different rules.

The probably easiest way to see a connection between non-monotonic logics
and abstract size is by considering preferential structures.
Preferential structures
define principal filters, generated by the set of minimal elements, as
follows:
if $ \xbf \xcn \xbq $ holds in such a structure, then $ \xbm (\xbf) \xcc
M(\xbq),$ where $ \xbm (\xbf)$ is the
set of minimal elements of $M(\xbf).$ According to our ideas, we define
a principal filter $ \xdf $ over $M(\xbf)$ by $X \xbe \xdf $ iff $ \xbm
(\xbf) \xcc X \xcc M(\xbf).$ Thus,
$M(\xbf) \xcs M(\xCN \xbq)$ will be a ``small'' subset of $M(\xbf).$
(Recall that
filters contain the ``big'' sets, and ideals the ``small'' sets.)

We can now go back and forth between rules on size and logical rules,
e.g.:

(For details, see
Table \ref{Table Base2-Size-Rules-1} (page \pageref{Table Base2-Size-Rules-1}),
Table \ref{Table Base2-Size-Rules-2} (page \pageref{Table Base2-Size-Rules-2}),
and
Table \ref{Table Mul-Laws} (page \pageref{Table Mul-Laws}).)

 \xEh

 \xDH
The ``AND'' rule corresponds to the filter property (finite intersections of
big subsets are still big).

 \xDH
``Right weakening'' corresponds to the rule that supersets of big sets
are still big.

 \xDH
It is natural, but beyond filter properties themselves, to postulate that,
if $X$ is a small subset of $Y,$ and $Y \xcc Y',$ then $X$ is also a
small subset of $Y'.$
We call such properties
``coherence properties'' between filters.
This property corresponds to the logical rule $ \xCf (wOR).$

 \xDH
In the rule $(CM_{ \xbo }),$ usually called Cautious Monotony, we change
the
base set a little when going from $M(\xba)$ to $M(\xba \xcu \xbb)$
(the change is small
by the prerequisite $ \xba \xcn \xbb),$ and still have $ \xba \xcu \xbb
\xcn \xbb',$ if we had
$ \xba \xcn \xbb'.$ We see here a conceptually very different use of
``small'', as we now change the base set, over which the filter is
defined, by a small amount.

 \xDH
The rule of Rational Monotony is the last one in the first table, and
somewhat
isolated there. It
is better to be seen as a multiplicative law, as described in the third
table.
It corresponds to the rule that the product of medium (i.e, neither big
nor small) sets, still has medium size.

 \xEj
\subsubsection{
Interpolation and size
}

The connection between non-monotonic logic and
the abstract concept of size was investigated in
 \cite{GS09a}, see also  \cite{GS08f}.
There, we looked among other things at abstract addition
of size. Here, we will show a connection to abstract multiplication of
size.
Our semantic approach used decomposition of set theoretical products.
An important step was to write a set of models $ \xbS $ as a product of
some
set $ \xbS' $ (which was a restriction of $ \xbS),$ and some full
Cartesian product.
So, when we speak about size, we will have (slightly simplified) some big
subset
$ \xbS_{1}$ of one product $ \xbP_{1},$ and some big subset $ \xbS_{2}$ of
another product $ \xbP_{2},$
and will now check whether $ \xbS_{1} \xCK \xbS_{2}$ is a big subset of $
\xbP_{1} \xCK \xbP_{2}.$
In shorthand,
whether ``$big*big=big$''.
(See Definition \ref{Definition Sin-Size-Rules} (page \pageref{Definition
Sin-Size-Rules})  for precise definitions.)
Such conditions are called
coherence conditions, as they do not concern the notion of size itself,
but the way the sizes defined for different base sets are connected.
Our main results here are
Proposition \ref{Proposition Sin-Interpolation-1} (page \pageref{Proposition
Sin-Interpolation-1})  and
Proposition \ref{Proposition Mul-Mu*1-Int} (page \pageref{Proposition
Mul-Mu*1-Int}). They say that if the logic under
investigation is defined from a notion of size which satisfies
sufficiently
many multiplicative
conditions, then this logic will have interpolation of type three or even
two, see Paragraph \ref{Paragraph Three-Variants} (page \pageref{Paragraph
Three-Variants}).

Consider now some set product $X \xCK X'.$ (Intuitively, $X$ and $X' $
are model sets
on sublanguages $J$ and $J' $ of the whole language $L.)$ When we have now
a rule like: If $Y$ is a big subset of $X,$ and $Y' $ a big subset of $X'
,$ then
$Y \xCK Y' $ is a big subset of $X \xCK X',$ and conversely, we can
calculate consequences
separately in the sublanguages, and put them together to have the overall
consequences. But this is the principle behind interpolation: we can work
with independent parts.

This is made precise in
Definition \ref{Definition Sin-Size-Rules} (page \pageref{Definition
Sin-Size-Rules}),
in particular by the rule

$(\xbm *1):$ $ \xbm (X \xCK X')= \xbm (X) \xCK \xbm (X').$

(Note that the conditions $(\xbm *i)$ and $(\xbS *i)$ are equivalent, as
shown
in Proposition \ref{Proposition Sin-Product-Small} (page \pageref{Proposition
Sin-Product-Small})  (for principal filters).)

The main result is
that the multiplicative size rule $(\xbm *1)$ entails
non-monotonic interpolation of the form $ \xbf \xcn \xba \xcn \xbq,$ see
Proposition \ref{Proposition Mul-Mu*1-Int} (page \pageref{Proposition
Mul-Mu*1-Int}).

We take now a closer look at interpolation for non-monotonic
logic.
\paragraph{
The three variants of interpolation \\[2mm]
}

\label{Paragraph Three-Variants}

Consider preferential logic, a rule like $ \xbf \xcn \xbq.$ This means
that
$ \xbm (\xbf) \xcc M(\xbq).$ So we go from $M(\xbf)$ to $ \xbm (
\xbf),$ the minimal models of $ \xbf,$
and then to $M(\xbq),$ and, abstractly, we have $M(\xbf) \xcd \xbm (
\xbf) \xcc M(\xbq),$ so we have
neither necessarily $M(\xbf) \xcc M(\xbq),$ nor $M(\xbf) \xcd M(
\xbq),$ the relation between
$M(\xbf)$ and $M(\xbq)$ may be more complicated. Thus, we have neither
the monotone,
nor the antitone case. For this reason, our general results for monotone
or antitone logics do not hold any more.

But we also see here that classical logic is used, too. Suppose that
there is $ \xbf' $ which describes exactly $ \xbm (\xbf),$ then we can
write
$ \xbf \xcn \xbf' \xcl \xbq.$

So we can split preferential logic into a
core part - going from $ \xbf $ to its minimal models - and a second part,
which is
just classical logic. (Similar decompositions are
also natural for other non-monotonic logics.)
Thus, preferential logic can be seen as a combination of two logics, the
non-monotonic core, and classical logic.
It is thus natural to consider variants of the interpolation problem,
where $ \xcn $ denotes again preferential logic, and $ \xcl $ as usual
classical
logic:

Given $ \xbf \xcn \xbq,$ is there ``simple'' $ \xba $ such that
 \xEh

 \xDH
$ \xbf \xcn \xba \xcl \xbq,$ or

 \xDH
$ \xbf \xcl \xba \xcn \xbq,$ or

 \xDH
$ \xbf \xcn \xba \xcn \xbq $?

 \xEj

In most cases, we will only consider the semantical version,
as the problems of the syntactical version are very similar to those
for monotonic logics.
We turn to the variants.

 \xEh

 \xDH

The first variant, $ \xbf \xcn \xba \xcl \xbq,$ has a complete
characterization
in Proposition \ref{Proposition Sin-Non-Mon-Int-Karl} (page \pageref{Proposition
Sin-Non-Mon-Int-Karl}),
provided we have a suitable normal form (conjunctions of disjunctions).
The condition says that the relevant variables of $ \xbm (\xbf)$ have to
be
relevant for $M(\xbf).$

 \xDH

The second variant, $ \xbf \xcl \xba \xcn \xbq,$ is related to very (and
in many cases,
too) strong conditions about size.
We do not have a complete characterization, only sufficient conditions
about
size.
The size conditions we need are
(see Definition \ref{Definition Sin-Size-Rules} (page \pageref{Definition
Sin-Size-Rules})):

the abovementioned $(\xbm *1),$ and,

$(\xbm *2):$ $ \xbm (X) \xcc Y$ $ \xch $ $ \xbm (X \xex A) \xcc Y \xex A$

where $X$ need not be a product any more.

The result is given in
Proposition \ref{Proposition Sin-Interpolation-1} (page \pageref{Proposition
Sin-Interpolation-1}).

Example \ref{Example Sin-Prod-Size} (page \pageref{Example Sin-Prod-Size}) 
shows that $(\xbm *2)$ seems too
strong
when compared to probability defined size. We repeat this
example here, for the reader's convenience.

\be

$\hspace{0.01em}$


\label{Example Sin-Prod-Size-a}

Take a language of 5 propositional variables, with $X':=\{a,b,c\},$ $X''
:=\{d,e\}.$
Consider the model set $ \xbS:=\{ \xCL a \xCL b \xCL cde,$ $-a-b-c-d \xCL
e\},$ i.e. of 8 models of
$ \xCf de$ and 2 models of $- \xCf d.$ The models of $ \xCf de$ are 8/10
of all elements of $ \xbS,$ so
it is reasonable to call them a big subset of $ \xbS.$ But its projection
on
$X'' $ is only 1/3 of $ \xbS''.$

So we have a potential $ \xCf decrease$ when going to the coordinates.

This shows that weakening the prerequisite about $X$ as done in $(\xbm
*2)$ is not
innocent.

\ee

We should, however, note that sufficiently modular preferential relations
guarantee these very strong properties of the big sets, see
Section \ref{Section Ham-Rel-Dist} (page \pageref{Section Ham-Rel-Dist}).

 \xDH

We turn to the third variant, $ \xbf \xcn \xba \xcn \xbq.$ This is
probably the
most interesting one, (a) as it is more general, it loosens the connection
with classical logic, (b) it seems more natural as a rule,
and (c) it is also connected to more natural laws about size.
Again, we do not have a complete characterization, only sufficient
conditions
about size.
Here, $(\xbm *1)$ suffices, and we have our main result about
non-monotonic semanti interpolation,
Proposition \ref{Proposition Mul-Mu*1-Int} (page \pageref{Proposition
Mul-Mu*1-Int}),
that $(\xbm *1)$ entails interpolation of the type $ \xbf \xcn \xba \xcn
\xbq.$

Proposition \ref{Proposition Mod-Hamming} (page \pageref{Proposition
Mod-Hamming})  shows that $(\xbm *1)$ is
(roughly) equivalent to
the relation properties

$ \xCf (GH1)$ $ \xbs \xec \xbt $ $ \xcu $ $ \xbs' \xec \xbt' $ $ \xcu $
$(\xbs \xeb \xbt $ $ \xco $ $ \xbs' \xeb \xbt')$ $ \xch $ $ \xbs \xbs
' \xeb \xbt \xbt' $

(where $ \xbs \xec \xbt $ iff $ \xbs \xeb \xbt $ or $ \xbs = \xbt)$

$ \xCf (GH2)$ $ \xbs \xbs' \xeb \xbt \xbt' $ $ \xch $ $ \xbs \xeb \xbt $
$ \xco $ $ \xbs' \xeb \xbt' $

of a preferential relation.

$(\xCf (GH2)$ means that some compensation is possible, e.g., $ \xbt \xeb
\xbs $ might be the
case, but $ \xbs' \xeb \xbt' $ wins in the end, so $ \xbs \xbs' \xeb
\xbt \xbt'.)$

There need not always be a semantical interpolation for the third variant,
this
is shown in
Example \ref{Example Sin-Non-Mon-Int} (page \pageref{Example Sin-Non-Mon-Int}).

 \xEj

So we see that, roughly,
semantic interpolation for nonmonotonic logics works when abstract size is
defined in a modular way - and we find independence again.
In a way, this is not surprising, as we use independent definition of
validity for interpolation in classical logic, and we use independent
definition of additional structure (relations or size) for
interpolation in non-monotonic logic.
\subsubsection{
Hamming relations and size
}

\label{Section Ham-Rel-Size}

As preferential relations are determined by a relation, and give rise
to abstract notions of size and their manipulation, it is natural to
take a close look at the corresponding properties of the relation.
We already gave a few examples in the preceding sections, so we can be
concise here. Our main definitions and results on this subject are
to be found in
Section \ref{Section Ham-Rel-Dist} (page \pageref{Section Ham-Rel-Dist}), where
we also discuss distances with
similar properties.

It is not surprising that we find various types of Hamming relations and
distances in this context, as they are, by definition, modular. Neither is
it surprising that we see them again in
Chapter \ref{Chapter Neighbourhood} (page \pageref{Chapter Neighbourhood}), as
we are interested there in
independent ways
to define neighbourhoods.

Basically, these relations and distances come in two flavours, the
set and the counting variant. This is perhaps best illustrated by the
Hamming
distance of two sequence of finite, equal length. We can define the
distance
by the $ \xCf set$ of arguments where they differ, or by the $ \xCf
cardinality$
of this set. The first results in possibly incomparable distances, the
second allows ``compensation'', difference in one argument can be
compensated by equality in another argument.

For definitions and results, also those connecting them to notions of
size, see Section \ref{Section Ham-Rel-Dist} (page \pageref{Section
Ham-Rel-Dist})  in particular
Definition \ref{Definition Sin-Set-HR} (page \pageref{Definition Sin-Set-HR}).
We then show in
Proposition \ref{Proposition Mod-Hamming} (page \pageref{Proposition
Mod-Hamming})
that (smooth) Hamming relations
generate our size conditions when size is defined as above from
a relation (the set of preferred elements generates the principal
filter). Thus, Hamming relations determine logics which have
interpolation, see
Corollary \ref{Corollary Sin-Interpolation-2} (page \pageref{Corollary
Sin-Interpolation-2}).

We define Hamming relations twice, in
Section \ref{Section Ham-Rel-Dist} (page \pageref{Section Ham-Rel-Dist}), and
in
Section \ref{Section Ham-Neigh} (page \pageref{Section Ham-Neigh}),
their uses and definitions differ slightly.
\subsubsection{
Equilibrium logic
}

Equilibrium logic, due to D.Pearce, A.Valverde, see
 \cite{PV09} for motivation and further discussion, is based
on the 3-valued finite Goedel logic, also called HT logic,
HT for ``here and there''. Our results are presented in
Section \ref{Section EQ} (page \pageref{Section EQ}).

Equilibrium logic (EQ) is defined by a choice function on the model set.
First models have to be ``total'', no variable of the language
may have 1 as value. Second, if $m \xeb m',$ then $m$ is considered
better,
and $m' $ discarded, where $m \xeb m' $ iff $m$ and $m' $ give value 0 to
the same
variables, and $m$ gives value 2 to strictly less (as subset) variables
than $m' $ does.

We can define equilibrium logic by a preferential relation (taking care
also
of the first condition), but it is not smooth. Thus, our general results
from the beginning of this section will not hold, and we have to work
with ``hand knitted'' solutions. We first show that equilibrium logic
has no interpolation of the form $ \xbf \xcl \xba \xcn \xbq $ or $ \xbf
\xcn \xba \xcl \xbq,$ then that is has
interpolation of the form $ \xbf \xcn \xba \xcn \xbq,$ and that the
interpolant is also
definable, i.e., equilibrium logic has semantic and syntactic
interpolation
of this form. Essentially, semantic interpolation is due to the fact that
the preference relation is defined in a modular way, using individual
variables - as always, when we have interpolation.
\subsubsection{
Interpolation for revision and argumentation
}

We have a short and simple result
(Lemma \ref{Lemma Mul-TR} (page \pageref{Lemma Mul-TR}))
for interpolation in AGM revision.
Unfortunately, we need the variables from both sides of the
revision operator as can easily be seen by revising with TRUE.
The reader is referred to
Section \ref{Section Int-Dist-Rev} (page \pageref{Section Int-Dist-Rev})  for
details.

Somewhat surprisingly, we also have an interpolation result for one
form of argumentation, where we consider the set of arguments for a
statement
as the truth value of that statement. As we have maximum (set union), we
have the lower bound used in
Proposition \ref{Proposition Gin-Pr-Int} (page \pageref{Proposition Gin-Pr-Int})
 for the monotonic case, and can
show
Fact \ref{Fact Arg-Int} (page \pageref{Fact Arg-Int}).
See Section \ref{Section Inter-Arg} (page \pageref{Section Inter-Arg})  for
details.
\subsubsection{
Language change to obtain products
}

To achieve interpolation and other results of independence, we
often need to write a set of models as a non-trivial product.
Sometimes, this is impossible, but an equivalent reformulation
of the language can selve the problem.

As this might be interesting also for the non-specialists, we
repeat Example \ref{Example Mod-Lang-Fact} (page \pageref{Example
Mod-Lang-Fact})  here:

\be

$\hspace{0.01em}$


\label{Example Mod-Lang-Fact-a}

Consider $p=3,$ and let

$ \xCf abc,$ $a \xCN bc,$ $a \xCN b \xCN c,$ $ \xCN abc,$ $ \xCN a \xCN b
\xCN c,$ $ \xCN ab \xCN c$ be the $6=2*3$ positive cases,

$ab \xCN c,$ $ \xCN a \xCN bc$ the negative ones. (It is coincidence that
we can factorize
positive and negative cases - probably iff one of the factors is the full
product, here 2, it could also be 4 etc.)

We divide the cases by 3 new variables, grouping them together in
positive and negative cases. $a' $ is indifferent, we want this to be the
independent factor, the negative ones will be put into $ \xCN b' \xCN c'
.$
The procedure has to be made precise still. (n): negative

Let $ \xCf a' $ code the set $ \xCf abc,$ $a \xCN bc,$ $a \xCN b \xCN c,$
$ab \xCN c$ (n),

Let $ \xCN a' $ code $ \xCN a \xCN bc$ (n), $ \xCN abc,$ $ \xCN a \xCN b
\xCN c,$ $ \xCN ab \xCN c.$

Let $ \xCf b' $ code $ \xCf abc,$ $a \xCN bc,$ $ \xCN a \xCN b \xCN c,$ $
\xCN ab \xCN c$

Let $ \xCN b' $ code $a \xCN b \xCN c,$ $ab \xCN c$ (n), $ \xCN a \xCN bc$
(n), $ \xCN abc$

Let $ \xCf c' $ code $ \xCf abc,$ $a \xCN b \xCN c,$ $ \xCN abc,$ $ \xCN a
\xCN b \xCN c$

Let $ \xCN c' $ code $a \xCN bc,$ $ab \xCN c$ (n), $ \xCN a \xCN bc$ (n),
$ \xCN ab \xCN c$

Then the 6 positive instances are

$\{a', \xCN a' \} \xCK \{b' c',$ $b' \xCN c',$ $ \xCN b' c' \},$ the
negative ones

$\{a', \xCN a' \} \xCK \{ \xCN b' \xCN c' \}$

As we have 3 new variables, we code again all possible cases, so
expressivity is the same.

\ee

Crucial here is that $6=3*2,$ so we can just re-arrange the 6 models in a
different way, see
Fact \ref{Fact Mod-Lang-Fact} (page \pageref{Fact Mod-Lang-Fact}).

A similar result holds for the non-monotonic case, where the
structure must be possible, we can then redefine the language.

All details are to be found in
Section \ref{Section Lang-Manip} (page \pageref{Section Lang-Manip}).
\subsection{
Summary
}

We use our proven strategy of ``divide et impera'',
transform the problem first in a semantical question, and then
in a purely algebraic one:

 \xEI

 \xDH
Classical and basic non-monotonic logic (looking for the sharpest
consequence)
have a surprising same answer, problems show up with definability when
going back to the syntactical question.

 \xDH
Thus, we separate algebraic from logical
questions, and we see that there are logics with algebraic
interpolation, but without logical interpolation, as the
necessary sets of models are not definable in the language.
This opens the way to making the language richer to
obtain interpolation, when so desired.

 \xDH
Full non-monotonic logic is more complicated, and finds a partial answer
using the concept of size and a novel manipulation of it, justified by
certain modular relations.

 \xDH
Finally, our approach also has the advantage of short and elementary
proofs.

 \xEJ
\section{
Neighbourhood semantics
}

Neighbourhood semantics, see
Chapter \ref{Chapter Neighbourhood} (page \pageref{Chapter Neighbourhood}),
probably first introduced by
D.Scott and R.Montague in
 \cite{Sco70} and  \cite{Mon70}, and already used for deontic
logic by O.Pacheco
in  \cite{Pac07}, seem to be useful for many logics:

 \xEh

 \xDH
in preferential logics, they describe the limit variant, where we consider
neighbourhoods of an ideal, usually inexistent, situation,

 \xDH
in approximative reasoning, they describe the approximations to the final
result,

 \xDH
in deontic and default logic, they describe the
``good'' situations, i.e., deontically acceptable, or where defaults have
fired.

 \xEj

Neighbourhood semantics are used, when the ``ideal'' situation does not
exist (e.g., preferential systems without minimal elements), or
are too difficult to obtain (e.g., ``perfect'' deontic states).
\subsection{
Defining neighbourhoods
}

Neighbourhoods can be defined in various ways:

 \xEI
 \xDH
by algebraic systems, like unions of intersections of certain sets
(but not complements),
 \xDH
quality relations, which say that some points are better than others,
carrying
over to sets of points,
 \xDH
distance relations, which measure the distance to the perhaps inexistant
ideal points.
 \xEJ

The relations and distances may be given already by the underlying
structure,
e.g., in preferential structures, or they can be defined in a natural
way, e.g., from a systems of sets, as in deontic logic or default logic.
In these cases, we can define a distance between two points by the number
or set of deontic requirements or default rules which one satisfies, but
not the other. A quality relation is defined in a similar way: a point is
better, if it satisfies more requirements or rules.
\subsection{
Additional requirements
}

With these tools, we can define properties neighbourhoods should
have. E.g., we may require them to be downward closed, i.e., if $x \xbe
N,$
where $N$ is a neighbourhood, $y \xeb x,$ $y$ is better than $x,$ then $y$
should also be
in $N.$ This is a property we will certainly require in neighbourhood
semantics
for preferential structures (in the limit version). For these structures,
we will also require that for every $x \xce N,$ there should be some $y
\xbe N$ with
$y \xeb x.$ We may also require that, if $x \xbe N,$ $y \xce N,$ and $y$
is in some aspect
better than $x,$ then there must be $z \xbe N,$ which is better than both,
so
we have some kind of ``ceteris paribus'' improvement.
\subsection{
Connections between the various properties
}

There is a multitude of possible definitions (via distances, relations,
set systems), and properties, so it is not surprising that one can
investigate a multitude of connections between the different possible
definitions of neighbourhoods. We cannot cover all possible connections,
so
we compare only a few cases, and the reader is invited to complete the
picture
for the cases which interest him. The connections we examined are
presented in
Section \ref{Section Exam} (page \pageref{Section Exam}).
\subsection{
Various uses of neighbourhood semantics
}

We also distinguish the different uses of the systems of sets thus
characterized as neighbourhoods: we can look at all formulas which hold in
(all
or some) such sets (as in neighbourhood semantics for preferential
logics), or
at the formulas which exactly describe them. The latter reading avoids the
infamous Ross paradox of deontic logic.
This distinction is simple, but basic, and did probably not receive
the attention it deserves, in the literature.
\section{
An abstract view on modularity and independence
}

\label{Section Abstract-View}
\subsection{
Introduction
}

We see independence and modularity in many situations.
Roughly, it means that we can use logic as a child who plays with building
blocks and puts them together. The big picture is not more than the
elements.

One example, which seems so natural, that it is hardly ever mentioned, is
validity in classical propositonal logic. The validity of $p$ in a model
does
not depend on the values a model assigns to other propositional variables.
By induction, this property carries over to more complicated formulas.
Consequently, the validity of a formula $ \xbf $ does not depend on the
laguage:
it suffices to know the values for the fragment in which $ \xbf $ is
formulated to
decide if $ \xbf $ holds or not. This is evident, but very important, it
justifies
what we call
``semantic interpolation'': Semantic interpolation will always hold for
monotone or antitone logics. It does $ \xCf not$ follow that the language
is
sufficiently rich to describe such an interpolant, the latter will then
be ``syntactic interpolation''.
Syntactic interpolation can be guaranteed by the existence of suitable
normal forms, which allow to treat model subsets independently.

For preferential non-monotonic logic, we see conditions for
the resulting abstract notion of size and its multiplication which
guarantee
semantic interpolation also for those logics. Natural conditions for
the preference relation result in such properties of abstract size.

Independence is also at the basis of an approach to theory revision due to
Parikh and his co-authors, see
 \cite{CP00}.
Again, natural conditions on a distance relation
result in such independent ways of revision.

The rule of Rational Monotony
(see Table \ref{Table Base2-Rules-Def-Conn-2} (page \pageref{Table
Base2-Rules-Def-Conn-2}))
can also be seen as independence: we can
``cut up'' the domain, and the same rules will still hold in the fragments.
\subsection{
Abstract definition of independence
}

\label{Section Mul-Def}

$ \xCO $

\vspace{10mm}

\begin{diagram}

\label{Diagram Mul-Commut}
\index{Diagram Mul-Commut}

\centering
\setlength{\unitlength}{1mm}
{\renewcommand{\dashlinestretch}{30}
\begin{picture}(150,100)(0,0)

\path(35,90)(78,90)
\path(75.2,91)(78,90)(75.2,89)
\path(66,60)(106,60)
\path(103.2,61)(106,60)(103.2,59)
\path(35,30)(78,30)
\path(75.2,31)(78,30)(75.2,29)

\path(35,35)(55,55)
\path(52.3,53.7)(55,55)(53.7,52.3)
\path(88,35)(108,55)
\path(105.3,53.7)(108,55)(106.7,52.3)
\path(90.7,36.3)(88,35)(89.3,37.7)

\path(35,85)(55,65)
\path(53.7,67.7)(55,65)(52.3,66.3)
\path(88,85)(108,65)
\path(106.7,67.7)(108,65)(105.3,66.3)
\path(89.3,82.3)(88,85)(90.7,83.7)

\put(29,89){{\xssc $\xbS_1$}}
\put(29,29){{\xssc $\xbS_2$}}
\put(55,59){{\xssc $\xbS_1 \xDO \xbS_2$}}

\put(81,89){{\xssc $f(\xbS_1$)}}
\put(81,29){{\xssc $f(\xbS_2$)}}
\put(108,59){{\xssc $f(\xbS_1 \xDO \xbS_2)=f(\xbS_1) \xDO' f(\xbS_2$)}}

\put(35,10){{\xssc Note that $\xDO$ and $\xDO'$ might be different}}
\put(50,1){{\xssc Independence}}

\end{picture}
}

\end{diagram}

\vspace{4mm}

$ \xCO $

The right notion of independence in our context seems to be:

We have compositions $ \xDO $ and $ \xDO',$ and an operation $f.$ We can
calculate
$f(\xbS_{1} \xDO \xbS_{2})$ from $f(\xbS_{1})$ and $f(\xbS_{2}),$ but
also conversely, given
$f(\xbS_{1} \xDO \xbS_{2})$ we can calculate $f(\xbS_{1})$ and $f(
\xbS_{2}).$
Of course, in other contexts, other notions of independence
might be adequate.
More precisely:

\bd

$\hspace{0.01em}$


\label{Definition Mul-Ind}

Let $f: \xdd \xcp \xdc $ be any function from domain $ \xdd $ to co-domain
$ \xdc.$
Let $ \xDO $ be a ``composition function'' $ \xDO: \xdd \xCK \xdd \xcp \xdd
,$ likewise
for $ \xDO': \xdc \xCK \xdc \xcp \xdc.$

We say that $ \xBc f, \xDO, \xDO'  \xBe $ are independent iff for any $
\xbS_{i}
\xbe \xdd $

(1) $f(\xbS_{1} \xDO \xbS_{2})=f(\xbS_{1}) \xDO' f(\xbS_{2}),$

(2) we can recover $f(\xbS_{i})$ from $f(\xbS_{1} \xDO \xbS_{2}),$
provided we know how
$ \xbS_{1} \xDO \xbS_{2}$ splits into the $ \xbS_{i},$ without using $f$
again.
\subsubsection{
Discussion
}

\ed

 \xEh

 \xDH Ranked structures satisfy it:

Let $ \xDO = \xDO' = \xcv.$ Let $f$ be the minimal model operator $ \xbm
$ of preferential logic.
Let $X,Y \xcc X \xcv Y$ have (at least) medium size, i.e.
$X \xcs \xbm (X \xcv Y) \xEd \xCQ,$ $Y \xcs \xbm (X \xcv Y) \xEd \xCQ,$
(see below,
Section \ref{Section Mul-Nota} (page \pageref{Section Mul-Nota})).
Then $ \xbm (X \xcv Y)= \xbm (X) \xcv \xbm (Y),$ and $ \xbm (X)= \xbm (X
\xcv Y) \xcs X,$ $ \xbm (Y)= \xbm (X \xcv Y) \xcs Y.$

 \xDH Consistent classical formulas and their interpretation satisfy it:

Let $ \xDO $ be conjunction in the composed language,
$ \xDO' $ be model set intersection, $f(\xbf)=M(\xbf).$
Let $ \xbf,$ $ \xbq $ be classical formulas, defined on disjoint language
fragments
$ \xdl,$ $ \xdl' $ of some language $ \xdl''.$ Then $f(\xbf \xcu \xbq
)=M(\xbf) \xcs M(\xbq),$ and
$M(\xbf)$ is the projection of $M(\xbf) \xcs M(\xbq)$ onto the
(models of) language $ \xdl,$
likewise for $M(\xbq).$
This is due to the way validity is defined, using only variables which
occur in the formula.

As a consequence, monotonic logic has semantical interpolation -
see  \cite{GS09c}, and below,
Section \ref{Section Mul-Mon-Int} (page \pageref{Section Mul-Mon-Int}). The
definition of being insensitive
is justified by
this modularity.

 \xDH It does not hold for inconsistent classical formulas:
We cannot recover $M(a \xcu \xCN a)$ and $M(b)$ from $M(a \xcu \xCN a \xcu
b),$ as we do not
know where the inconsistency came from.
The basic reason is trivial: One empty factor suffices to make the whole
product empty, and we do not know which factor was the culprit.
See Section \ref{Section Mul-Relev} (page \pageref{Section Mul-Relev})  for the
discussion of a remedy.

 \xDH Preferential logic satisfies it under certain conditions:

If $ \xbm (X \xCK Y)= \xbm (X) \xCK \xbm (Y)$ holds for model products and
$ \xcn,$ then it
holds by definition. An important consequence is that such a logic
has interpolation of the form $ \xcn \xDO \xcn,$
see Section \ref{Section Mul-Int} (page \pageref{Section Mul-Int}).

 \xDH Modular revision a la Parikh,
see  \cite{CP00}, is based on a similar idea.

 \xEj
\subsection{
Other aspects of independence
}
\subsubsection{
Existence of normal forms
}

We may see the existence of conjunctive and disjunctive normal forms
as a form of independence: A formula may be split into elementary parts,
which are then put together by the standard operations of inf $(\xcu)$
and
sup $(\xco),$ resulting immediately in the existence of syntactic
interpolation, as both the upper and lower limits of interpolation
are definable. Note that higher finite Goedel logics do not allow
these operations, basically as we cannot always decompose nested
intuitionistic implication.
\subsubsection{
Language change
}

Independence of language fragments gives us the following perspectives:

 \xEh

 \xDH it makes independent and parallel treatment of fragments possible,
and offers thus efficient treatment in applications
(descriptive logics etc.):

Consider $X=X' \xcv X'',$ where $X',X'' $ are disjoint. Suppose size is
calculated
independently, in the following sense: Let $Y \xcc X,$ then $Z \xcc Y$ is
big iff
$Z \xcs X' \xcc Y \xcs X' $ and $Z \xcs X'' \xcc Y \xcs X'' $ both are
big. We can then calculate size
independently.

 \xDH it results in new rules similar to the classical ones like AND, OR,
Cumulativity, etc. We can thus obtain postulates about reasonable
behaviour,
but also classification by those rules, see
Table \ref{Table Mul-Laws} (page \pageref{Table Mul-Laws}), Scenario 2, Logical
property.

 \xDH it sheds light on notions like
``ceteris paribus'', which we saw in the context of obligations,
see  \cite{GS08g}, and
Definition \ref{Definition Quality-Extension} (page \pageref{Definition
Quality-Extension}),

 \xDH it clarifies notions like
``normal with respect to $ \xbf,$ but not $ \xbq $'', see
 \cite{GS08e} and  \cite{GS08f},

 \xDH it helps to understand e.g. inheritance diagrams where
arrows make other information accessible, and we need an underlying
mechanism to combine bits of information, given in different
languages, see again
 \cite{GS08e} and  \cite{GS08f}.

 \xEj
\subsubsection{
A relevance problem
}

\label{Section Mul-Relev}

Consider the formula $ \xbf:=a \xcu \xCN a \xcu b.$ Then $M(\xbf)= \xCQ
.$ But we cannot recover
where the problem came from (it might also come from $b \xcu \xCN b),$ and
this results
in the EFQ rule. We now discuss
one, purely algebraic, approach to remedy.

Consider 3 valued models, with a new value $b$ for both, in addition to
$t$ and $f.$
Above formula would then have the model $m(a)=b,$ $m(b)=t.$ So there is a
model,
EFQ fails, and we can recover the culprit.

To have the usual behaviour of $ \xcu $ as intersection, it might be good
to change
the definition so that $m(x)=b$ is always a model. Then $M(b)=\{m(b)=t,m'
(b)=b\},$
$M(\xCN b)=\{m(b)=f,m' (b)=b\},$ and $M(b \xcu \xCN b)=\{m' (b)=b\}.$

It is not yet clear which version to choose, and we have no syntactic
characterization.
\subsubsection{
Small subspaces
}

When considering small subsets in nonmonotonic logic, we neglect small
subsets of models.
What is the analogue when considering small subspaces, i.e. when
$J=J' \xcv J'',$ with $J'' $ small in $J$ in the sense of nonmonotonic
logic?

It is perhaps easiest to consider the relation based approach first.
So we have an order on $ \xbP J' $ and one on $ \xbP J'',$ $J'' $ is
small, and we want
to know how to construct a corresponding order on $ \xbP J.$ Two solutions
come
to mind:

 \xEI

 \xDH a less radical one: we make a lexicographic ordering, where the one
on $ \xbP J' $
has precedence over the one on $ \xbP J'',$

 \xDH a more radical one: we totally forget about the ordering of $ \xbP
J'',$ i.e. we
do as if the ordering on $ \xbP J'' $ were the empty set, i.e.
$ \xbs' \xbs'' \xeb \xbt' \xbt'' $ iff $ \xbs' \xeb \xbt' $ and $
\xbs'' = \xbt''.$

We call this condition $forget(J'').$

 \xEJ

The less radical one is already covered by our relation conditions
$ \xCf (GH)$ see Definition \ref{Definition Sin-Set-HR} (page
\pageref{Definition Sin-Set-HR}).
The more radical one is probably more interesting. Suppose $ \xbf' $ is
written in
language $J',$ $ \xbf'' $ in language $J'',$ we then have

$ \xbf' \xcu \xbf'' \xcn \xbq' \xcu \xbq'' $ iff $ \xbf' \xcn \xbq'
$ and $ \xbf'' \xcl \xbq''.$

This approach is of course the same as considering on the small
coordinate only ALL as a big subset, (see the lines $x*1/1*x$
in Table \ref{Table Mul-Laws} (page \pageref{Table Mul-Laws})).
\section{
Conclusion and outlook
}

In Section \ref{Section Outlook} (page \pageref{Section Outlook}), we argue
that logics which diverge
from classical
logic in the sense that they allow to conclude more or less than classical
logic concludes need an additional fundamental concept, a $ \xCf
justification.$
Classical logic has language and truth values, proof theory, and
semantics.
Here, we need more, justification, why we are allowed to conclude more or
less.
We have to show that the price we pay (divergence from truth) is
justified,
e.g., by more efficient reasoning, conjectures which ``pay'', etc.

We think that we need a new fundamental concept, which is on the same
level as proof theory and semantics.

This is an open research problem, but it seems that our tools like
abstract manipulation of abstract size are sufficient to attack it.
\section{
Previously published material, acknowledgements
}

This text builds upon previous research by the authors.
To make the text self-contained, it is therefore necessary to repeat some
previously published material.
We give now the parts concerned and their sources.

All parts of Chapter \ref{Chapter Mod-Base-Def} (page \pageref{Chapter
Mod-Base-Def})  which are not marked as
new material were published in some or all of
 \cite{Sch04},
 \cite{GS08b},
 \cite{GS08c},
 \cite{GS09a},
 \cite{GS08f}.

The additive laws on abstract size
(see Section \ref{Section Add-Size} (page \pageref{Section Add-Size})) were
published in
 \cite{GS09a} and  \cite{GS08f}.

The formal material of
Chapter \ref{Chapter Neighbourhood} (page \pageref{Chapter Neighbourhood})  was
already published in
 \cite{GS08f}, it is put here in a wider perspective.

Finally, we would like to thank D.Makinson and D.Pearce for useful
comments and
very interesting questions.
\clearpage
\chapter{
Basic definitions
}

\label{Chapter Mod-Base-Def}
\section{
Introduction
}
\subsection{
Overview of this chapter
}

This chapter contains basic definitions and results, sometimes
slightly beyond the immediate need of this book, as we
also want to put our work a bit more in perspective, and make it
self-contained, for the convenience of the reader.
Most of the material of this chapter (unless marked as
``new'') was published previously, see
 \cite{Sch04},
 \cite{GS08b},
 \cite{GS08c},
 \cite{GS09a}, and
 \cite{GS08f}.

We begin with basic algebraic and logical definitions, including in
particular
many laws of non-monotonic logics, in their syntactic and semantic
variants,
showing also the connections between both sides, see
Definition \ref{Definition Log-Cond-N} (page \pageref{Definition Log-Cond-N}) 
and the tables
Table \ref{Table Base2-Rules-Def-Conn-1} (page \pageref{Table
Base2-Rules-Def-Conn-1})  and
Table \ref{Table Base2-Rules-Def-Conn-2} (page \pageref{Table
Base2-Rules-Def-Conn-2}).

It seems to be a little known result that even the classical operators
permit an unusual interpretation in the infinite case, but we
claim no originality, see
Example \ref{Example Non-Standard} (page \pageref{Example Non-Standard}).

We would like to emphasize the importance of the definability preservation
(dp) property. In the infinite case, not all model sets $X$ are definable,
i.e.,
there is some formula $ \xbf $ or theory $T$ such that $X=M(\xbf)$ - the
models of $ \xbf $ - or
$X=M(T)$ - the models of $T.$
It is by no means evident that a model choice function $ \xbm,$ applied
to a
definable model set, gives us back again a definable model set (is
definability
preserving, or dp). If $ \xbm $ does
not have this property, some representation results will not hold, which
hold if $ \xbm $ is dp, and representation results become much more
complicated,
see  \cite{Sch04} for positive and for impossibility results.
In our present context, definability is again an important concept.
Even if we have semantic interpolation, if language and operators are not
strong enough, we cannot define the semantic interpolants, so we have
semantic, but not syntactic interpolation. Examples are found in
finite Goedel logics, see
Section \ref{Section Mon-Synt-Int} (page \pageref{Section Mon-Synt-Int}).
New operators guaranteeing the definability of particularly interesting,
``universal'' interpolants, see
Definition \ref{Definition Univ-Int} (page \pageref{Definition Univ-Int}), are
discussed in
Section \ref{Section Mon-Interpol-Int} (page \pageref{Section Mon-Interpol-Int})
.
They are intricately related to the existence of
conjunctive and disjunctive normal forms, as discussed in
Section \ref{Section Analoga} (page \pageref{Section Analoga}).

We conclude this part with a - to our knowledge - unpublished result
that we can define only countably many inconsistent formulas, see
Example \ref{Example Co-Ex-Inf} (page \pageref{Example Co-Ex-Inf}). (The
question is due to D.Makinson.)

We then give a detailed introduction into the basic concepts
of many-valued logics, again, as readers might not be so familiar with
the generalizations from 2-valued to many-valued logic. In particular, the
nice correspondence between 2-valued functions and sets does not hold any
more, so we have to work with arbitrary functions. We have to re-define
what a definable model ``set'' is, and what semantical interpolation
means for many-valued logic.
Table \ref{Table Gin-Not-Def} (page \pageref{Table Gin-Not-Def})  gives an
overview.

We then give an introduction to preferential structures and the logic
they define. These structures are among the best examined semantics for
non-monotonic logics, and
Chapter \ref{Chapter Size-Laws} (page \pageref{Chapter Size-Laws})  is also
based on the investigation of
such structures.
We first introduce the minimal variant, and then the limit variant. The
first variant is the ususal one, the second is needed to deal with
cases where there are no minimal models, due to infinite descending
chains.
(The first variant was introduced by Y.Shoham in
 \cite{Sho87b}, the second variant by P.Siegel et al. in
 \cite{BS85}. It should, however, be emphasized, that
preferential models were introduced as a semantics for deontic
logic long before they were investigated as a semantics for
non-monotonic logic, see
 \cite{Han69}).
The limit variant was further investigated in
 \cite{Sch04}, and we refer the reader there for representation
and impossibility results.
An overview of representation results for the minimal variant is given
in Table \ref{Table Base1-Pref-Rep} (page \pageref{Table Base1-Pref-Rep}).

We introduce a new concept
in this section on preferential structures,
``bubble structures'', which, we think, present a useful tool for
abstraction,
and are a semantic variant of independence in preferential structures.
Here, we have a global preferential structure between subsets
(``bubbles'') of the model set, and a fine scale structure inside those
subsets. Seen from the outside, all elements of a bubble behave the
same way, so the whole set can be treated as one element, on the
inside, we see a finer structure.

Moreover, new material on many-valued preferential structures is included.

We then go into details in the section on IBRS, introduced by D.Gabbay,
see  \cite{Gab04}, and further investigated in
 \cite{GS08b} and  \cite{GS08f}, as they
are not so much common knowledge. We also discuss here if and how
the limit version of preferential structures might be applied
to reactive structures.

We then present theory revision, as introduced by Alchorron, Gardenfors,
and Makinson, see  \cite{AGM85}. Again,
we also build on previous results by (here, one of) the authors,
when we discuss distance based revision, introduced by
Lehmann, Magidor, and Schlechta, see
 \cite{LMS95},
 \cite{LMS01}, and elaborated in
 \cite{Sch04}.
We also include a short paragraph on new material for theory revision
based on
many-valued logic.
\section{
Basic algebraic and logical definitions
}

$ \xCO $
\index{Notation FOL-Tilde}

\bn

$\hspace{0.01em}$


\label{Notation D-FOL-Tilde}

\index{FOL}
\index{NML}
We use sometimes FOL as abbreviation for first order logic, and NML for
nonmonotonic logic.
To avoid Latex complications in bigger expressions, we replace
$\widetilde{xxxxx}$ by $\wt{xxxxx}$.

\en

$ \xCO $

$ \xCO $
\index{Definition Algebraic basics}

\bd

$\hspace{0.01em}$


\label{Definition Alg-Base}

 \xEh

 \xDH

We use $ \xdp $ \index{$ \xdp $}  to denote the power set operator.
 \index{$ \xbP $}

$ \xbP \{X_{i}:i \xbe I\}$ $:=$ $\{g:$ $g:I \xcp \xcV \{X_{i}:i \xbe I\},$
$ \xcA i \xbe I.g(i) \xbe X_{i}\}$ is the general Cartesian
product, $X \xCK X' $ is the binary Cartesian product.

$card(X)$ shall denote the cardinality of $X,$
 \index{card}
and $V$ \index{$V$}  the
set-theoretic
universe we work in - the class of all sets.

Given a set of pairs $ \xdx,$ and a
set $X,$ we denote by $ \xdx \xex X:=\{ \xBc x,i \xBe  \xbe \xdx:x \xbe X\}.$
 \index{$ \xex $}
(When the context is clear, we
will sometime simply write $X$ for $ \xdx \xex X.)$

We will use the same notation $ \xex $ to denote the restriction of
functions and
in particular of sequences to a subset of the domain.

If $ \xbS $ is a set of sequences over an index set $X,$ and $X' \xcc X,$
we will abuse
notation and also write $ \xbS \xex X' $ for $\{ \xbs \xex X': \xbs \xbe
\xbS \}.$

Concatenation of sequences, e.g., of $ \xbs $ and $ \xbs',$ will be
denoted by
juxtaposition: $ \xbs \xbs'.$

 \xDH

$A \xcc B$ will denote that $ \xCf A$ is a subset of $B$ or equal to $B,$
and $A \xcb B$ that $ \xCf A$ is
a proper subset of $B,$ likewise for $A \xcd B$ and $A \xcf B.$
 \index{$ \xcc $}
 \index{$ \xcb $}
 \index{$ \xcd $}
 \index{$ \xcf $}

Given some fixed set $U$ we work in, and $X \xcc U,$ then $ \xdC (X):=U-X$
.
 \index{$ \xdC $}

 \xDH

If $ \xdy \xcc \xdp (X)$ for some
$X,$ we say that $ \xdy $ satisfies

$(\xcs)$ \index{$(\xcs)$}  iff it is closed under finite intersections,

$(\xcS)$ \index{$(\xcS)$}  iff it is closed under arbitrary
intersections,

$(\xcv)$ \index{$(\xcv)$}  iff it is closed under finite unions,

$(\xcV)$ \index{$(\xcV)$}  iff it is closed under arbitrary unions,

$(\xdC)$ \index{$(\xdC)$}  iff it is closed under complementation,

$ \xCf (-)$ \index{$ \xCf (-)$}  iff it is closed under set difference.

 \xDH

We will sometimes write $A=B \xFO C$ for: $A=B,$ or $A=C,$ or $A=B \xcv
C.$
 \index{$ \xFO $}

 \xEj

We make ample and tacit use of the Axiom of Choice.

\ed

$ \xCO $

$ \xCO $
\index{Definition Relations}

\bd

$\hspace{0.01em}$


\label{Definition Rel-Base}

$ \xeb^{*}$ \index{$ \xeb^{*}$}  will denote the transitive closure of the
relation $ \xeb.$ If a relation
$<,$
$ \xeb,$ or similar is given, $a \xcT b$ \index{$a \xcT b$}  will express
that $ \xCf a$ and $b$ are $<-$ (or $ \xeb -)$
incomparable - context will tell. Given any relation $<,$ $ \xck $ \index{$ \xck
$}  will stand for
$<$ or $=,$ conversely, given $ \xck,$ $<$ \index{$<$}  will stand
for $ \xck,$ but not $=,$ similarly
for $ \xeb $ etc.

\ed

$ \xCO $

$ \xCO $

\bd

$\hspace{0.01em}$


\label{Definition Weak-Filter}

Fix a base set $X.$

A (weak) filter on or over $X$ is a set $ \xdf \xcc \xdp (X),$ s.t.
$(F1)-(F3)$ $((F1),$ (F2),
$(F3')$ respectively) hold:

(F1) $X \xbe \xdf $

(F2) $A \xcc B \xcc X,$ $A \xbe \xdf $ imply $B \xbe \xdf $

(F3) $A,B \xbe \xdf $ imply $A \xcs B \xbe \xdf $

$(F3')$ $A,B \xbe \xdf $ imply $A \xcs B \xEd \xCQ.$

So a weak filter satisfies $(F3')$ instead of (F3).

\index{principal filter}
A filter is called a principal filter iff there is $X' \xcc X$ s.t. $ \xdf
=\{A:$ $X' \xcc A \xcc X\}.$

An (weak) ideal on or over $X$ is a set $ \xdi \xcc \xdp (X),$ s.t.
$(I1)-(I3)$ $((I1),$ (I2),
$(I3')$ respectively) hold:
(I1) $ \xCQ \xbe \xdi $

(I2) $A \xcc B \xcc X,$ $B \xbe \xdi $ imply $A \xbe \xdi $

(I3) $A,B \xbe \xdi $ imply $A \xcv B \xbe \xdi $

$(I3')$ $A,B \xbe \xdi $ imply $A \xcv B \xEd X.$

So a weak ideal satisfies $(I3')$ instead of (I3).

A filter is an abstract notion of size,
\index{big set}
elements of a filter on $X$ are called big subsets of $X,$ their
complements are
\index{small set}
\index{medium size set}
called small, and the rest have medium size. The dual applies to ideals,
this is
justified by the following trivial fact:

\ed

\bfa

$\hspace{0.01em}$


\label{Fact Weak-Filter}

If $ \xdf $ is a (weak) filter on $X,$ then $ \xdi:=\{X-A:A \xbe \xdf \}$
is a (weak) ideal on $X,$
if $ \xdi $ is a (weak) ideal on $X,$ then $ \xdf:=\{X-A:A \xbe \xdf \}$
is a (weak) filter on $X.$

\efa

$ \xCO $

$ \xCO $
\index{Definition Logic, basics}

\bd

$\hspace{0.01em}$


\label{Definition Log-Base}

 \xEh

 \xDH

$V$ will be the set of truth values when there are more than the classical
ones, TRUE and FALSE.

We work here in a classical propositional language $ \xdl,$ a theory
\index{theory}  $T$ will be
an
arbitrary set of formulas. Formulas will often be named $ \xbf,$ $ \xbq
,$ etc., theories
$T,$ $S,$ etc.

$v(\xdl)$ \index{$v(\xdl)$}  or simply
$L$ \index{$L$}  will be the set of propositional variables of $ \xdl
.$

$F(\xdl)$ \index{$F(\xdl)$}  will be the set of formulas of $ \xdl.$

A propositional model $m$ will be a function from the set of propositional
variables to the set of truth values - $V$ when we have more than two
truth values.

$M_{ \xdl }$ \index{$M_{ \xdl }$}  or simply
$M$ \index{$M$}  when the context is clear, will be the set of
(classical) models for $ \xdl,$
$M(T)$ \index{$M(T)$}  or
$M_{T}$ \index{$M_{T}$}
is the set of models of $T,$ likewise $M(\xbf)$ \index{$M(\xbf)$}  for a
formula $ \xbf.$

 \xDH

$ \xdD_{ \xdl }$ \index{$ \xdD_{ \xdl }$}  $:=\{M(T):$ $T$ a theory in $ \xdl
\},$
the set of
$ \xCf definable$ \index{$ \xCf definable$}  model sets.

Note that, in classical propositional logic, $ \xCQ,M_{ \xdl } \xbe
\xdD_{ \xdl },$ $ \xdD_{ \xdl }$ contains
singletons, is closed under arbitrary intersections and finite unions.

An operation $f: \xdy \xcp \xdp (M_{ \xdl })$ for $ \xdy \xcc \xdp (M_{
\xdl })$ is called
$ \xCf definability$ $ \xCf preserving$ \index{$ \xCf definability$ $ \xCf
preserving$},
$ \xCf (dp)$ \index{$ \xCf (dp)$}  or
$(\xbm dp)$ \index{$(\xbm dp)$}  in short, iff for all $X \xbe \xdD_{ \xdl }
\xcs \xdy $ $f(X) \xbe \xdD_{ \xdl }.$

We will also use $(\xbm dp)$ for binary functions $f: \xdy \xCK \xdy \xcp
\xdp (M_{ \xdl })$ - as needed
for theory revision - with the obvious meaning.

 \xDH

$ \xcl $ \index{$ \xcl $}  will be classical derivability, and

$ \ol{T}:=\{ \xbf:T \xcl \xbf \},$ the closure of $T$ under $ \xcl.$

 \index{$ \ol{T}$}

 \xDH

$Con(.)$ \index{$Con(.)$}  will stand for classical consistency, so $Con(
\xbf)$ will mean that
$ \xbf $ is classical consistent, likewise for $Con(T).$ $Con(T,T')$ will
stand for
$Con(T \xcv T'),$ etc.

 \xDH

Given a consequence relation $ \xcn,$ we define

$ \ol{ \ol{T} }:=\{ \xbf:T \xcn \xbf \}.$

 \index{$ \ol{ \ol{T} }$}

(There is no fear of confusion with $ \ol{T},$ as it just is not useful to
close
twice under classical logic.)

 \xDH

$T \xco T' $ \index{$T \xco T' $}  $:=\{ \xbf \xco \xbf': \xbf \xbe T, \xbf
' \xbe T' \}.$

 \xDH

If $X \xcc M_{ \xdl },$ then $Th(X)$ \index{$Th(X)$}  $:=\{ \xbf:X \xcm
\xbf \},$ likewise for
$Th(m)$ \index{$Th(m)$}, $m \xbe M_{ \xdl }.$
($ \xcm $ \index{$ \xcm $}  will usually be classical validity.)

 \xEj

\ed

$ \xCO $

In the following,
the $X_{i}$ are arbitrary (non-empty) sets, standing for the set of truth
values, and $J$ is intuitively the set of propositional variables.
So any element of $ \xbP \{X_{i}:i \xbe J\}$ stands for a propositional
model.
Inessential variables in a set $ \xbS $ are those which do not have any
influence
on the truth of the formula whose model set is $ \xbS.$

\bd

$\hspace{0.01em}$


\label{Definition Sin-Ir-Relevant}

Let $ \xbS \xcc \xbP:= \xbP \{X_{i}:i \xbe J\}.$ Define:

(1) For $ \xbs, \xbs' \xbe \xbS,$ $J' \xcc J,$ define:

$ \xbs \xCq_{J' } \xbs': \xcj \xcA x \xbe J' \xbs (x)= \xbs' (x).$

(2) $I(\xbS):=\{i \xbe J:$ $ \xbS = \xbS \xex (J-\{i\}) \xCK X_{i}\}$
(up to re-ordering), (the irrelevant
or inessential $i)$ and

$R(\xbS):=J-I(\xbS)$ (the relevant or essential $i).$

\ed

\bfa

$\hspace{0.01em}$


\label{Fact Sin-Ir-Relevant}

(1) $ \xbS = \xbS \xex R(\xbS) \xCK \xbP \xex I(\xbS)$ (up to
re-ordering)

(2) $ \xbs \xex R(\xbS)= \xbs' \xex R(\xbS)$ $ \xcu $ $ \xbs \xbe
\xbS $ $ \xch $ $ \xbs' \xbe \xbS.$

\efa

\subparagraph{
Proof
}

$\hspace{0.01em}$


(1) Enumerate $I(\xbS),$ $I(\xbS)=\{i:i< \xbk \}.$
Define $ \xbS_{j}:= \xbS \xex (R(\xbS) \xcv \{i:i<j\}).$ We show by
induction that
$ \xbS_{j}= \xbS \xex R(\xbS) \xCK \xbP \xex (I(\xbS) \xcs j)$ for $j
\xck \xbk $ (up to re-ordering).

$j=0$ is trivial - there is nothing to show.

$j \xcp j+1:$ This follows from the induction hypothesis and the
definition of $I(\xbS).$

$j$ is a limit ordinal: Any sequence of length $j$
can be written as the coherent union of shorter sequences, and these are
in both sets, as the result holds for $j' <j$ by induction hypothesis.

(2) Trivial.

$ \xcz $
\\[3ex]

We give the following example which shows that even the operators
of classical propositional logic may be interpreted in a
non-standard way (in the infinite case). We give this example,
as the basic argument is about definability, and thus fully in our
context.

\be

$\hspace{0.01em}$


\label{Example Non-Standard}

Consider a (countable for simplicity) infinite propositional language.
Consider the usual model set $M,$ and take any countable set of models
away,
resulting in $M'.$
Interpret the variables as usual, $p$ by $M(p),$ and let $M' (p):=M(p)
\xcs M'.$
Interpret $ \xCN $ and $ \xcu $ as usual, i.e., complement and
intersection, but now
in $M',$ define $ \xco,$ $ \xcp $ from $ \xCN $ and $ \xcu.$
Then $M' (\xCN \xbf):=M' -M' (\xbf)=M' -(M' \xcs M(\xbf))=M' \xcs M(
\xCN \xbf),$
$M' (\xbf \xcu \xbq)=M' (\xbf) \xcs M' (\xbq)=M' \xcs M(\xbf) \xcs
M(\xbq)$ by induction, so we have
$M' (\xbf)=M' \xcs M(\xbf).$
Consequently, the usual axioms hold, and Modus Ponens is a valid rule.
For completeness, we have to show that every consistent formula has a
model.
Let $ \xbf $ be consistent, then $M' (\xbf)=M' \xcs M(\xbf),$ but as $
\xbf $ is consistent, it
has uncountably many models, so it also has a model in $M'.$
$ \xcz $
\\[3ex]

\ee

To put our work more into perspective, we repeat now
material from  \cite{GS08c}. This gives the main definitions and
rules
for non-monotonic logics, see
Table \ref{Table Base2-Rules-Def-Conn-1} (page \pageref{Table
Base2-Rules-Def-Conn-1})  and
Table \ref{Table Base2-Rules-Def-Conn-2} (page \pageref{Table
Base2-Rules-Def-Conn-2}),
``Logical rules, definitions and connections''.

$ \xCO $

\index{Definition Logical conditions}

\bd

$\hspace{0.01em}$


\label{Definition Log-Cond-N}

The definitions are given in
Table \ref{Table Base2-Rules-Def-Conn-1} (page \pageref{Table
Base2-Rules-Def-Conn-1}),
``Logical rules, definitions and connections Part $ \xfI $'',
and
Table \ref{Table Base2-Rules-Def-Conn-2} (page \pageref{Table
Base2-Rules-Def-Conn-2}),
``Logical rules, definitions and connections Part II'',
which also show connections between different versions
of rules, the semantics, and rules about size.
(The tables are split in two, as they would not fit
onto one page otherwise.)

Explanation of the tables:

 \xEh
 \xDH
The first table gives the basic properties, the second table
those for Cumulativity and Rational Monotony.
 \xDH
The difference between the first two columns is that the first column
treats the formula version of the rule, the second the more general
theory (i.e., set of formulas) version.
 \xDH ``Corr.'' stands for ``Correspondence''.
 \xDH
The third column, ``Corr.'', is to be understood as follows:

Let a logic $ \xcn $ satisfy $ \xCf (LLE)$ and $ \xCf (CCL),$ and define a
function $f: \xdD_{ \xdl } \xcp \xdD_{ \xdl }$
by $f(M(T)):=M(\ol{ \ol{T} }).$ Then $f$ is well defined, satisfies $(
\xbm dp),$ and $ \ol{ \ol{T} }=Th(f(M(T))).$

If $ \xcn $ satisfies a rule in the left hand side,
then - provided the additional properties noted in the middle for $ \xch $
hold, too -
$f$ will satisfy the property in the right hand side.

Conversely:

If $f: \xdy \xcp \xdp (M_{ \xdl })$ is a function, with $ \xdD_{ \xdl }
\xcc \xdy,$ and we define a logic
$ \xcn $ by $ \ol{ \ol{T} }:=Th(f(M(T))),$ then $ \xcn $ satisfies $ \xCf
(LLE)$ and $ \xCf (CCL).$
If $f$ satisfies $(\xbm dp),$ then $f(M(T))=M(\ol{ \ol{T} }).$

If $f$ satisfies a property in the right hand side,
then - provided the additional properties noted in the middle for $ \xci $
hold, too -
$ \xcn $ will satisfy the property in the left hand side.

 \xDH
We use the following abbreviations for those supplementary conditions in
the
``Correspondence'' columns:

``$T= \xbf $'' means that, if one of the theories
(the one named the same way in Definition \ref{Definition Log-Cond-N} (page
\pageref{Definition Log-Cond-N}))
is equivalent to a formula, we do not need $(\xbm dp).$

$-(\xbm dp)$ stands for ``without $(\xbm dp)$''.

 \xDH
$A=B \xFO C$ will abbreviate $A=B,$ or $A=C,$ or $A=B \xcv C.$

 \xEj
Further comments:

\begin{enumerate}

\item

$(PR)$ is also called $infinite$ $conditionalization$ We choose this name for
its central role for preferential structures $(PR)$ or $(\xbm PR).$

\item

The system of rules $(AND)$ $(OR)$ $(LLE)$ $(RW)$ $(SC)$ $(CP)$ $(CM)$ $(CUM)$
is also called system $P$ (for preferential). Adding $(RatM)$ gives the system
$R$ (for rationality or rankedness).

Roughly: Smooth preferential structures generate logics satisfying system
$P$, while ranked structures generate logics satisfying system $R$.

\item

A logic satisfying $(REF)$, $(ResM)$, and $(CUT)$ is called a $consequence$
$relation$.

\item

$(LLE)$ and$(CCL)$ will hold automatically, whenever we work with model sets.

\item

$(AND)$ is obviously closely related to filters, and corresponds to closure
under finite intersections. $(RW)$ corresponds to upward closure of filters.

More precisely, validity of both depend on the definition, and the
direction we consider.

Given $f$ and $(\xbm \xcc)$, $f(X)\xcc X$ generates a principal filter:
$\{X'\xcc X:f(X)\xcc X'\}$, with
the definition: If $X=M(T)$, then $T\xcn \xbf$  iff $f(X)\xcc M(\xbf)$.
Validity of $(AND)$ and
$(RW)$ are then trivial.

Conversely, we can define for $X=M(T)$

$\xdx:=\{X'\xcc X: \xcE \xbf (X'=X\xcs M(\xbf)$ and $T\xcn \xbf)\}$.

$(AND)$ then makes $\xdx$  closed under
finite intersections, and $(RW)$ makes $\xdx$  upward
closed. This is in the infinite case usually not yet a filter, as not all
subsets of $X$ need to be definable this way.
In this case, we complete $\xdx$  by
adding all $X''$ such that there is $X'\xcc X''\xcc X$, $X'\xbe\xdx$.

Alternatively, we can define

$\xdx:=\{X'\xcc X: \xcS\{X \xcs M(\xbf): T\xcn \xbf \} \xcc X' \}$.

\item

$(SC)$ corresponds to the choice of a subset.

\item

$(CP)$ is somewhat delicate, as it presupposes that the chosen model set is
non-empty. This might fail in the presence of ever better choices, without
ideal ones; the problem is addressed by the limit versions.

\item

$(PR)$ is an infinitary version of one half of the deduction theorem: Let $T$
stand for $\xbf$, $T'$ for $\xbq$, and $\xbf \xcu \xbq \xcn \xbs$,
so $\xbf \xcn \xbq \xcp \xbs$, but $(\xbq \xcp \xbs)\xcu \xbq \xcl \xbs$.

\item

$(CUM)$ (whose more interesting half in our context is $(CM)$) may best be seen
as normal use of lemmas: We have worked hard and found some lemmas. Now
we can take a rest, and come back again with our new lemmas. Adding them to the
axioms will neither add new theorems, nor prevent old ones to hold.
(This is, of course, a meta-level argument concerning an object level rule.
But also object level rules should - at least generally - have an intuitive
justification, which will then come from a meta-level argument.)

\end{enumerate}

\begin{table}[h]

\index{$(Opt)$}
\index{$(\xbm \xcc)$}
\index{$(REF)$}
\index{Reflexivity}
\index{$(LLE)$}
\index{Left Logical Equivalence}
\index{$(RW)$}
\index{Right Weakening}
\index{$(iM)$}
\index{$(wOR)$}
\index{$(\xbm wOR)$}
\index{$(eM\xdi)$}
\index{$(disjOR)$}
\index{$(\xbm disjOR)$}
\index{$(I\xcv disj)$}
\index{$(CP)$}
\index{$(\xbm \xCQ)$}
\index{$(I_1)$}
\index{Consistency Preservation}
\index{$(\xbm \xCQ fin)$}
\index{$(AND_1)$}
\index{$(I_2)$}
\index{$(AND_n)$}
\index{$(I_n)$}
\index{$(AND)$}
\index{$(I_\xbo)$}
\index{$(CCL)$}
\index{Classical Closure}
\index{$(iM)$}
\index{$(OR)$}
\index{$(\xbm OR)$}
\index{$(PR)$}
\index{$(\xbm PR)$}
\index{$(\xbm PR')$}
\index{$(CUT)$}
\index{$ (\xbm CUT) $}

\tabcolsep=0.5pt
\caption{Logical rules, definitions and connections Part I}

\label{Table Base2-Rules-Def-Conn-1}
\begin{center}
\begin{turn}{90}
{\tiny
\begin{tabular}{|c|c|c|c|c|c|}

\hline

\multicolumn{6}{|c|}{\bf Logical rules, definitions and connections Part I}\\

\hline

\multicolumn{2}{|c|}{Logical rule}
\xEH
Corr.
\xEH
Model set
\xEH
Corr.
\xEH
Size Rules
\xEP

\hline
\hline

\multicolumn{6}{|c|}{Basics}

\xEP

\hline

$(SC)$ Supraclassicality
\index{$(SC)$}
\index{Supraclassicality}
\xEH
$(SC)$
\xEH
$\xch$
\xEH
$(\xbm \xcc)$
\xEH
trivial
\xEH
$(Opt)$
\xEP

\cline{3-3}

$ \xba \xcl \xbb $ $ \xch $ $ \xba \xcn \xbb $
\xEH
$ \ol{T} \xcc \ol{ \ol{T} }$
\xEH
$\xci$
\xEH
$f(X) \xcc X$
\xEH
\xEH
\xEP

\cline{1-1}

$(REF)$ Reflexivity
\xEH
\xEH
\xEH
\xEH
\xEH
\xEP

$ T \xcv \{\xba\} \xcn \xba $
\xEH
\xEH
\xEH
\xEH
\xEH
\xEP

\hline

$(LLE)$
\xEH
$(LLE)$
\xEH
\xEH
(trivally true)
\xEH
\xEH
\xEP

Left Logical Equivalence
\xEH
\xEH
\xEH
\xEH
\xEH
\xEP

$ \xcl \xba \xcr \xba',  \xba \xcn \xbb   \xch $
\xEH
$ \ol{T}= \ol{T' }  \xch   \ol{\ol{T}} = \ol{\ol{T'}}$
\xEH
\xEH
\xEH
\xEH
\xEP

$ \xba' \xcn \xbb $
\xEH
\xEH
\xEH
\xEH
\xEH
\xEP

\hline

$(RW)$ Right Weakening
\xEH
$(RW)$
\xEH
\xEH
(upward closure)
\xEH
trivial
\xEH
$(iM)$
\xEP

$ \xba \xcn \xbb,  \xcl \xbb \xcp \xbb'   \xch $
\xEH
$ T \xcn \xbb,  \xcl \xbb \xcp \xbb'   \xch $
\xEH
\xEH
\xEH
\xEH
\xEP

$ \xba \xcn \xbb' $
\xEH
$T \xcn \xbb' $
\xEH
\xEH
\xEH
\xEH
\xEP

\hline

$(wOR)$
\xEH
$(wOR)$
\xEH
$\xch$
\xEH
$(\xbm wOR)$
\xEH
$\xcj$
\xEH
$(eM\xdi)$
\xEP

\cline{3-3}

$ \xba \xcn \xbb,$ $ \xba' \xcl \xbb $ $ \xch $
\xEH
$ \ol{ \ol{T} } \xcs \ol{T' }$ $ \xcc $ $ \ol{ \ol{T \xco T' } }$
\xEH
$\xci$
\xEH
$f(X \xcv Y) \xcc f(X) \xcv Y$
\xEH
\xEH
\xEP

$ \xba \xco \xba' \xcn \xbb $
\xEH
\xEH
\xEH
\xEH
\xEH
\xEP

\hline

$(disjOR)$
\xEH
$(disjOR)$
\xEH
$\xch$
\xEH
$(\xbm disjOR)$
\xEH
$\xcj$
\xEH
$(I\xcv disj)$
\xEP

\cline{3-3}

$ \xba \xcl \xCN \xba',$ $ \xba \xcn \xbb,$
\xEH
$\xCN Con(T \xcv T') \xch$
\xEH
$\xci$
\xEH
$X \xcs Y= \xCQ $ $ \xch $
\xEH
\xEH
\xEP

$ \xba' \xcn \xbb $ $ \xch $ $ \xba \xco \xba' \xcn \xbb $
\xEH
$ \ol{ \ol{T} } \xcs \ol{ \ol{T' } } \xcc \ol{ \ol{T \xco T' } }$
\xEH
\xEH
$f(X \xcv Y) \xcc f(X) \xcv f(Y)$
\xEH
\xEH
\xEP

\hline

$(CP)$
\xEH
$(CP)$
\xEH
$\xch$
\xEH
$(\xbm \xCQ)$
\xEH
trivial
\xEH
$(I_1)$
\xEP

\cline{3-3}

Consistency Preservation
\xEH
\xEH
$\xci$
\xEH
\xEH
\xEH
\xEP

$ \xba \xcn \xcT $ $ \xch $ $ \xba \xcl \xcT $
\xEH
$T \xcn \xcT $ $ \xch $ $T \xcl \xcT $
\xEH
\xEH
$f(X)= \xCQ $ $ \xch $ $X= \xCQ $
\xEH
\xEH
\xEP

\hline

\xEH
\xEH
\xEH
$(\xbm \xCQ fin)$
\xEH
\xEH
$(I_1)$
\xEP

\xEH
\xEH
\xEH
$X \xEd \xCQ $ $ \xch $ $f(X) \xEd \xCQ $
\xEH
\xEH
\xEP

\xEH
\xEH
\xEH
for finite $X$
\xEH
\xEH
\xEP

\hline

\xEH
$(AND_1)$
\xEH
\xEH
\xEH
\xEH
$(I_2)$
\xEP

\xEH
$\xba\xcn\xbb \xch \xba\xcN\xCN\xbb$
\xEH
\xEH
\xEH
\xEH
\xEP

\hline

\xEH
$(AND_n)$
\xEH
\xEH
\xEH
\xEH
$(I_n)$
\xEP

\xEH
$\xba\xcn\xbb_1, \ldots, \xba\xcn\xbb_{n-1} \xch $
\xEH
\xEH
\xEH
\xEH
\xEP

\xEH
$\xba\xcN(\xCN\xbb_1 \xco \ldots \xco \xCN\xbb_{n-1})$
\xEH
\xEH
\xEH
\xEH
\xEP

\hline

$(AND)$
\xEH
$(AND)$
\xEH
\xEH
(closure under finite
\xEH
trivial
\xEH
$(I_\xbo)$
\xEP

$ \xba \xcn \xbb,  \xba \xcn \xbb'   \xch $
\xEH
$ T \xcn \xbb, T \xcn \xbb'   \xch $
\xEH
\xEH
intersection)
\xEH
\xEH
\xEP

$ \xba \xcn \xbb \xcu \xbb' $
\xEH
$ T \xcn \xbb \xcu \xbb' $
\xEH
\xEH
\xEH
\xEH
\xEP

\hline

$(CCL)$ Classical Closure
\xEH
$(CCL)$
\xEH
\xEH
(trivally true)
\xEH
trivial
\xEH
$(iM)+(I_\xbo)$
\xEP

\xEH
$ \ol{ \ol{T} }$ classically closed
\xEH
\xEH
\xEH
\xEH
\xEP

\hline

$(OR)$
\xEH
$(OR)$
\xEH
$\xch$
\xEH
$(\xbm OR)$
\xEH
$\xcj$
\xEH
$(eM\xdi)+(I_\xbo)$
\xEP

\cline{3-3}

$ \xba \xcn \xbb,  \xba' \xcn \xbb   \xch $
\xEH
$ \ol{\ol{T}} \xcs \ol{\ol{T'}} \xcc \ol{\ol{T \xco T'}} $
\xEH
$\xci$
\xEH
$f(X \xcv Y) \xcc f(X) \xcv f(Y)$
\xEH
\xEH
\xEP

$ \xba \xco \xba' \xcn \xbb $
\xEH
\xEH
\xEH
\xEH
\xEH
\xEP

\hline

\xEH
$(PR)$
\xEH
$\xch$
\xEH
$(\xbm PR)$
\xEH
$\xcj$
\xEH
$(eM\xdi)+(I_\xbo)$
\xEP

\cline{3-3}

$ \ol{ \ol{ \xba \xcu \xba' } }$ $ \xcc $ $ \ol{ \ol{ \ol{ \xba } } \xcv
\{ \xba' \}}$
\xEH
$ \ol{ \ol{T \xcv T' } }$ $ \xcc $ $ \ol{ \ol{ \ol{T} } \xcv T' }$
\xEH
$\xci (\xbm dp)+(\xbm\xcc)$
\xEH
$X \xcc Y$ $ \xch $
\xEH
\xEH
\xEP

\cline{3-3}

\xEH
\xEH
$\xcI$ $-(\xbm dp)$
\xEH
$f(Y) \xcs X \xcc f(X)$
\xEH
\xEH
\xEP

\cline{3-3}

\xEH
\xEH
$\xci (\xbm\xcc)$
\xEH
\xEH
\xEH
\xEP

\xEH
\xEH
$T'=\xbf$
\xEH
\xEH
\xEH
\xEP

\cline{3-4}

\xEH
\xEH
$\xci$
\xEH
$(\xbm PR')$
\xEH
\xEH
\xEP

\xEH
\xEH
$T'=\xbf$
\xEH
$f(X) \xcs Y \xcc f(X \xcs Y)$
\xEH
\xEH
\xEP

\hline

$(CUT)$
\xEH
$(CUT)$
\xEH
$\xch$
\xEH
$ (\xbm CUT) $
\xEH
$\xci$
\xEH
$(eM\xdi)+(I_\xbo)$
\xEP

\cline{3-3}
\cline{5-5}

$ T  \xcn \xba; T \xcv \{ \xba\} \xcn \xbb \xch $
\xEH
$T \xcc \ol{T' } \xcc \ol{ \ol{T} }  \xch $
\xEH
$\xci$
\xEH
$f(X) \xcc Y \xcc X  \xch $
\xEH
$\xcH$
\xEH
\xEP

$ T  \xcn \xbb $
\xEH
$ \ol{ \ol{T'} } \xcc \ol{ \ol{T} }$
\xEH
\xEH
$f(X) \xcc f(Y)$
\xEH
\xEH
\xEP

\hline

\end{tabular}
}
\end{turn}
\end{center}
\end{table}

\begin{table}[h]

\index{$(wCM)$}
\index{$(eM\xdf)$}
\index{$(CM_2)$}
\index{$(I_2)$}
\index{$(CM_n)$}
\index{$(I_n)$}
\index{$(CM)$}
\index{Cautious Monotony}
\index{$ (\xbm CM) $}
\index{$(\xdm^+_\xbo)$}
\index{$(ResM)$}
\index{Restricted Monotony}
\index{$(\xbm ResM)$}
\index{$(CUM)$}
\index{Cumulativity}
\index{$(\xbm CUM)$}
\index{$ (\xcc \xcd) $}
\index{$ (\xbm \xcc \xcd) $}
\index{$(eM\xdi)$}
\index{$(I_\xbo)$}
\index{($eM\xdf)$}
\index{$(RatM)$}
\index{Rational Monotony}
\index{$(\xbm RatM)$}
\index{$(\xdm^{++})$}
\index{$(RatM=)$}
\index{$(\xbm =)$}
\index{$(Log=')$}
\index{$(\xbm =')$}
\index{$(DR)$}
\index{$(Log \xFO)$}
\index{$(\xbm \xFO)$}
\index{$(Log \xcv)$}
\index{$(\xbm \xcv)$}
\index{$(Log \xcv')$}
\index{$(\xbm \xcv')$}
\index{$(\xbm \xbe)$}

\tabcolsep=0.5pt
\caption{Logical rules, definitions and connections Part II}

\label{Table Base2-Rules-Def-Conn-2}
\begin{center}
\begin{turn}{90}
{\tiny
\begin{tabular}{|c|c|c|c|c|c|}

\hline

\multicolumn{6}{|c|}{\bf Logical rules, definitions and connections Part II}\\

\hline

\multicolumn{2}{|c|}{Logical rule}
\xEH
Corr.
\xEH
Model set
\xEH
Corr.
\xEH
Size-Rule
\xEP

\hline
\hline

\multicolumn{6}{|c|}{Cumulativity}

\xEP

\hline

$(wCM)$
\xEH
\xEH
\xEH
\xEH
trivial
\xEH
$(eM\xdf)$
\xEP

$\xba\xcn\xbb, \xba'\xcl\xba, \xba\xcu\xbb\xcl\xba' \xch $
\xEH
\xEH
\xEH
\xEH
\xEH
\xEP

$\xba'\xcn\xbb$
\xEH
\xEH
\xEH
\xEH
\xEH
\xEP

\hline

$(CM_2)$
\xEH
\xEH
\xEH
\xEH
\xEH
$(I_2)$
\xEP

$\xba\xcn\xbb, \xba\xcn\xbb' \xch \xba\xcu\xbb\xcL\xCN\xbb'$
\xEH
\xEH
\xEH
\xEH
\xEH
\xEP

\hline

$(CM_n)$
\xEH
\xEH
\xEH
\xEH
\xEH
$(I_n)$
\xEP

$\xba\xcn\xbb_1, \ldots, \xba\xcn\xbb_n \xch $
\xEH
\xEH
\xEH
\xEH
\xEH
\xEP

$\xba \xcu \xbb_1 \xcu \ldots \xcu \xbb_{n-1} \xcL\xCN\xbb_n$
\xEH
\xEH
\xEH
\xEH
\xEH
\xEP

\hline

$(CM)$ Cautious Monotony
\xEH
$(CM)$
\xEH
$\xch$
\xEH
$ (\xbm CM) $
\xEH
$\xcj$
\xEH
$(\xdm^+_\xbo)(4)$
\xEP

\cline{3-3}

$ \xba \xcn \xbb,  \xba \xcn \xbb'   \xch $
\xEH
$T \xcc \ol{T' } \xcc \ol{ \ol{T} }  \xch $
\xEH
$\xci$
\xEH
$f(X) \xcc Y \xcc X  \xch $
\xEH
\xEH
\xEP

$ \xba \xcu \xbb \xcn \xbb' $
\xEH
$ \ol{ \ol{T} } \xcc \ol{ \ol{T' } }$
\xEH
\xEH
$f(Y) \xcc f(X)$
\xEH
\xEH
\xEP

\cline{1-1}

\cline{3-4}

or $(ResM)$ Restricted Monotony
\xEH
\xEH
$\xch$
\xEH
$(\xbm ResM)$
\xEH
\xEH
\xEP

\cline{3-3}

$ T  \xcn \xba, \xbb \xch T \xcv \{\xba\} \xcn \xbb $
\xEH
\xEH
$\xci$
\xEH
$ f(X) \xcc A \xcs B \xch $
\xEH
\xEH
\xEP

\xEH
\xEH
\xEH
$f(X \xcs A) \xcc B $
\xEH
\xEH
\xEP

\hline

$(CUM)$ Cumulativity
\xEH
$(CUM)$
\xEH
$\xch$
\xEH
$(\xbm CUM)$
\xEH
$\xci$
\xEH
$(eM\xdi)+(I_\xbo)+(\xdm^{+}_{\xbo})(4)$
\xEP

\cline{3-3}
\cline{5-5}

$ \xba \xcn \xbb   \xch $
\xEH
$T \xcc \ol{T' } \xcc \ol{ \ol{T} }  \xch $
\xEH
$\xci$
\xEH
$f(X) \xcc Y \xcc X  \xch $
\xEH
$\xcH$
\xEH
\xEP

$(\xba \xcn \xbb'   \xcj   \xba \xcu \xbb \xcn \xbb')$
\xEH
$ \ol{ \ol{T} }= \ol{ \ol{T' } }$
\xEH
\xEH
$f(Y)=f(X)$
\xEH
\xEH
\xEP

\hline

\xEH
$ (\xcc \xcd) $
\xEH
$\xch$
\xEH
$ (\xbm \xcc \xcd) $
\xEH
$\xci$
\xEH
$(eM\xdi)+(I_\xbo)+(eM\xdf)$
\xEP

\cline{3-3}
\cline{5-5}

\xEH
$T \xcc \ol{\ol{T'}}, T' \xcc \ol{\ol{T}} \xch $
\xEH
$\xci$
\xEH
$ f(X) \xcc Y, f(Y) \xcc X \xch $
\xEH
$\xcH$
\xEH
\xEP

\xEH
$ \ol{\ol{T'}} = \ol{\ol{T}}$
\xEH
\xEH
$ f(X)=f(Y) $
\xEH
\xEH
\xEP

\hline
\hline

\multicolumn{6}{|c|}{Rationality}

\xEP

\hline

$(RatM)$ Rational Monotony
\xEH
$(RatM)$
\xEH
$\xch$
\xEH
$(\xbm RatM)$
\xEH
$\xcj$
\xEH
$(\xdm^{++})$
\xEP

\cline{3-3}

$ \xba \xcn \xbb,  \xba \xcN \xCN \xbb'   \xch $
\xEH
$Con(T \xcv \ol{\ol{T'}})$, $T \xcl T'$ $ \xch $
\xEH
$\xci$ $(\xbm dp)$
\xEH
$X \xcc Y, X \xcs f(Y) \xEd \xCQ   \xch $
\xEH
\xEH
\xEP

\cline{3-3}

$ \xba \xcu \xbb' \xcn \xbb $
\xEH
$ \ol{\ol{T}} \xcd \ol{\ol{\ol{T'}} \xcv T} $
\xEH
$\xcI$ $-(\xbm dp)$
\xEH
$f(X) \xcc f(Y) \xcs X$
\xEH
\xEH
\xEP

\cline{3-3}

\xEH
\xEH
$\xci$ $T=\xbf$
\xEH
\xEH
\xEH
\xEP

\hline

\xEH
$(RatM=)$
\xEH
$\xch$
\xEH
$(\xbm =)$
\xEH
\xEH
\xEP

\cline{3-3}

\xEH
$Con(T \xcv \ol{\ol{T'}})$, $T \xcl T'$ $ \xch $
\xEH
$\xci$ $(\xbm dp)$
\xEH
$X \xcc Y, X \xcs f(Y) \xEd \xCQ   \xch $
\xEH
\xEH
\xEP

\cline{3-3}

\xEH
$ \ol{\ol{T}} = \ol{\ol{\ol{T'}} \xcv T} $
\xEH
$\xcI$ $-(\xbm dp)$
\xEH
$f(X) = f(Y) \xcs X$
\xEH
\xEH
\xEP

\cline{3-3}

\xEH
\xEH
$\xci$ $T=\xbf$
\xEH
\xEH
\xEH
\xEP

\hline

\xEH
$(Log=')$
\xEH
$\xch$
\xEH
$(\xbm =')$
\xEH
\xEH
\xEP

\cline{3-3}

\xEH
$Con(\ol{ \ol{T' } } \xcv T)$ $ \xch $
\xEH
$\xci$ $(\xbm dp)$
\xEH
$f(Y) \xcs X \xEd \xCQ $ $ \xch $
\xEH
\xEH
\xEP

\cline{3-3}

\xEH
$ \ol{ \ol{T \xcv T' } }= \ol{ \ol{ \ol{T' } } \xcv T}$
\xEH
$\xcI$ $-(\xbm dp)$
\xEH
$f(Y \xcs X)=f(Y) \xcs X$
\xEH
\xEH
\xEP

\cline{3-3}

\xEH
\xEH
$\xci$ $T=\xbf$
\xEH
\xEH
\xEH
\xEP

\hline

$(DR)$
\xEH
$(Log \xFO)$
\xEH
$\xch$
\xEH
$(\xbm \xFO)$
\xEH
\xEH
\xEP

\cline{3-3}

$\xba \xco \xbb \xcn \xbg \xch$
\xEH
$ \ol{ \ol{T \xco T' } }$ is one of
\xEH
$\xci$
\xEH
$f(X \xcv Y)$ is one of
\xEH
\xEH
\xEP

$\xba \xcn \xbg$ or $\xbb \xcn \xbg$
\xEH
$\ol{\ol{T}},$ or $\ol{\ol{T'}},$ or $\ol{\ol{T}} \xcs \ol{\ol{T'}}$ (by (CCL))
\xEH
\xEH
$f(X),$ $f(Y)$ or $f(X) \xcv f(Y)$
\xEH
\xEH
\xEP

\hline

\xEH
$(Log \xcv)$
\xEH
$\xch$ $(\xbm\xcc)+(\xbm=)$
\xEH
$(\xbm \xcv)$
\xEH
\xEH
\xEP

\cline{3-3}

\xEH
$Con(\ol{ \ol{T' } } \xcv T),$ $ \xCN Con(\ol{ \ol{T' } }
\xcv \ol{ \ol{T} })$ $ \xch $
\xEH
$\xci$ $(\xbm dp)$
\xEH
$f(Y) \xcs (X-f(X)) \xEd \xCQ $ $ \xch $
\xEH
\xEH
\xEP

\cline{3-3}

\xEH
$ \xCN Con(\ol{ \ol{T \xco T' } } \xcv T')$
\xEH
$\xcI$ $-(\xbm dp)$
\xEH
$f(X \xcv Y) \xcs Y= \xCQ$
\xEH
\xEH
\xEP

\hline

\xEH
$(Log \xcv')$
\xEH
$\xch$ $(\xbm\xcc)+(\xbm=)$
\xEH
$(\xbm \xcv')$
\xEH
\xEH
\xEP

\cline{3-3}

\xEH
$Con(\ol{ \ol{T' } } \xcv T),$ $ \xCN Con(\ol{ \ol{T' }
} \xcv \ol{ \ol{T} })$ $ \xch $
\xEH
$\xci$ $(\xbm dp)$
\xEH
$f(Y) \xcs (X-f(X)) \xEd \xCQ $ $ \xch $
\xEH
\xEH
\xEP

\cline{3-3}

\xEH
$ \ol{ \ol{T \xco T' } }= \ol{ \ol{T} }$
\xEH
$\xcI$ $-(\xbm dp)$
\xEH
$f(X \xcv Y)=f(X)$
\xEH
\xEH
\xEP

\hline

\xEH
\xEH
\xEH
$(\xbm \xbe)$
\xEH
\xEH
\xEP

\xEH
\xEH
\xEH
$a \xbe X-f(X)$ $ \xch $
\xEH
\xEH
\xEP

\xEH
\xEH
\xEH
$ \xcE b \xbe X.a \xce f(\{a,b\})$
\xEH
\xEH
\xEP

\hline

\end{tabular}
}
\end{turn}
\end{center}
\end{table}

\ed

$ \xCO $
\subsection{
Countably many disjoint sets
}

This might be the right place to add the following short remark
on formula sets.

$ \xCO $

We show here that - independent of the cardinality of the language -
one can define only countably many inconsistent formulas.

The question is due to D.Makinson (personal communication).

\be

$\hspace{0.01em}$


\label{Example Co-Ex-Inf}

There is a countably infinite set of formulas s.t. the defined model sets
are pairwise disjoint.

Let $p_{i}:i \xbe \xbo $ be propositional variables.

Consider $ \xbf_{i}:= \xcU \{ \xCN p_{j}:j<i\} \xcu p_{i}$ for $i \xbe
\xbo.$

Obviously, $M(\xbf_{i}) \xEd \xCQ $ for all $i.$

Let $i<i',$ we show $M(\xbf_{i}) \xcs M(\xbf_{i' })= \xCQ.$ $M(
\xbf_{i' }) \xcm \xCN p_{i},$ $M(\xbf_{i}) \xcm p_{i}.$

$ \xcz $
\\[3ex]

\ee

\bfa

$\hspace{0.01em}$


\label{Fact Co-Ex-Inf}

Any set $X$ of consistent formulas with pairwise disjoint model sets is at
most
countable.

\efa

\subparagraph{
Proof
}

$\hspace{0.01em}$


Let such $X$ be given.

(1) We may assume that $X$ consists of conjunctions of propositional
variables
or their negations.

Proof: Re-write all $ \xbf \xbe X$ as disjunctions of conjunctions $
\xbf_{j}.$ At least one of
the conjunctions $ \xbf_{j}$ is consistent. Replace $ \xbf $ by one such $
\xbf_{j}.$ Consistency
is preserved, as is pairwise disjointness.

(2) Let $X$ be such a set of formulas. Let $X_{i} \xcc X$ be the set of
formulas in $X$ with
length $i,$ i.e. a consistent conjunction of $i$ many propositional
variables or
their negations, $i>0.$

As the model sets for $X$ are pairwise disjoint, the model sets for all $
\xbf \xbe X_{i}$
have to be disjoint.

(3) It suffices now to show that each $X_{i}$ is at most countable, we
even show
that each $X_{i}$ is finite.

Proof by induction:

Consider $i=1.$ Let $ \xbf, \xbf' \xbe X_{1}.$ Let $ \xbf $ be $p$ or $
\xCN p.$ If $ \xbf' $ is not $ \xCN \xbf,$ then
$ \xbf $ and $ \xbf' $ have a common model. So one must be $p,$ the other
$ \xCN p.$ But these
are all possibilities, so $card(X_{1})$ is finite.

Let the result be shown for $k<i.$

Consider now $X_{i}.$ Take arbitrary $ \xbf \xbe X_{i}.$ Wlog, $ \xbf
=p_{1} \xcu  \Xl  \xcu p_{i}.$ Take arbitrary
$ \xbf' \xEd \xbf.$ As $M(\xbf) \xcs M(\xbf')= \xCQ,$ $ \xbf' $
must be a conjunction containing one of
$ \xCN p_{k},$ $1 \xck k \xck i.$ Consider now $X_{i,k}:=\{ \xbf' \xbe
X_{i}: \xbf' $ contains $ \xCN p_{k}\}.$
Thus $X_{i}=\{ \xbf \} \xcv \xcV \{X_{i,k}:1 \xck k \xck i\}.$ Note that
all $ \xbq, \xbq' \xbe X_{i,k}$ agree on $ \xCN p_{k},$
so the situation in $X_{i,k}$ is isomorphic to $X_{i-1}.$ So,
by induction hypothesis, $card(X_{i,k})$ is finite,
as all $ \xbf' \xbe X_{i,k}$ have to be mutually inconsistent. Thus,
$card(X_{i})$ is finite.
(Note that we did not use the fact that elements from different $X_{i,k},$
$X_{i,k' }$
also have to be mutually inconsistent, our rough proof suffices.)

$ \xcz $
\\[3ex]

Note that the proof depends very little on logic. We needed normal forms,
and used 2 truth values. Obviously, we can easily generalize to finitely
many truth values.

$ \xCO $
\subsection{
Introduction to many-valued logics
}

\label{Section Many-Val-Intro}

In 2-valued logic, we have a correspondence between sets and logic, e.g.,
we can speak abut the set of models of a formula. This has now to be
replaced by a many-valued function (or, alternatively, but generally
not pursued here, by many-valued sets).
\paragraph{
Motivation \\[2mm]
}

Preferential logics offer, indirectly, 4 truth values: classically true
and
false and defeasibly true and false. Inheritance systems offer, through
specificity, arbitrarily many truth values, with a partial order. Finite
Goedel logics offer arbitrarily many truth values, with a total order.
This motivates in our context to consider the following:

\bd

$\hspace{0.01em}$


\label{Definition Mod-Many-Log}

 \xEh

 \xDH

We assume here a finite set of truth values, $V,$ to be given, with a
partial
order $ \xck.$ We assume that a minimal and a maximal element exist,
which will
be denoted 0 and 1, TRUE and FALSE, min and max, depending on context.
In many cases we will assume that sup and inf (equivalently, max and min,
as
$V$ is supposed to be finite)
exist for any subset of $V.$
This will not always be necessary, but often it will be convenient.

Inf will correspond to classical $ \xcu,$ sup to classical $ \xco.$
Given $x \xbe V,$
$ \xdC x$ will be $inf\{y \xbe V:sup\{x,y\}=1\},$ it corresponds to
classical $ \xCN.$
We will not assume that classical $ \xCN $ is always part of the language.

 \xDH

A model is a function $m:L \xcp V.$

In classical logic, a formula $ \xbf $ defines a model set $M(\xbf) \xcc
M,$ equivalently
a function $f_{ \xbf }:M \xcp \{0,1\}$ with $f_{ \xbf }(m)=1: \xcj m \xcm
\xbf.$ A straightforward generalization
is to define in the many-valued case $f_{ \xbf }:M \xcp V.$

This definition should respect the following postulates:

 \xEh

 \xDH $f_{p}(m)=m(p)$ for $p \xbe L.$

This postulate is the basis for a seemingly trivial property, which
has far-reaching consequences: If $ \xbf $ contains only variables in $L'
\xcc L,$
then its truth value is the same in $L-$models and $L' -$models, whenever
they agree on $L'.$ Of course, validity of the operators has to be
truth functional, and again not to depend on other variables.

 \xDH $f_{ \xbf \xcu \xbq }(m)=inf\{f_{ \xbf }(m),f_{ \xbq }(m)\}$ and
$f_{ \xbf \xco \xbq }(m)=sup\{f_{ \xbf }(m),f_{ \xbq }(m)\}.$

 \xEj

For a set $T$ of formulas, we define $f_{T}(m):=inf\{f_{ \xbf }(m): \xbf
\xbe T\}.$

 \xDH

In general, a model set corresponds now to an arbitrary function $f:M \xcp
V.$
Such $f$ is called (formula) definable iff there is $ \xbf $ such that
$f=f_{ \xbf }.$
The definition of theory definable is analogous.

$ \xdp (M)$ is replaced by $V^{M},$ the set of all functions from $M$ to
$V,$ $ \xdd $ will denote
the set of all definable functions from $M$ to $V.$

 \xDH

Semantical consequence should respect $ \xck:$

 \xEh

 \xDH For $f,g:M \xcp V,$ we write $f \xcm g$ or $f \xck g$ iff $ \xcA m
\xbe M.f(m) \xck g(m),$

 \xDH we write $ \xbf \xcm \xbq $ iff $f_{ \xbf } \xcm f_{ \xbq },$

 \xDH and we assume that $ \xcp $ is compatible with $ \xcm:$

$f_{ \xbf \xcp \xbq }(m)=TRUE$ iff $f_{ \xbf }(m) \xck f_{ \xbq }(m),$ so

$f_{ \xbf \xcp \xbq }=$ (constant) TRUE iff $ \xbf \xcm \xbq.$

 \xEj

 \xEj
\subsubsection{
Definable model functions
}

\ed

In classical logic, the set of definable model sets satisfies certain
closure
conditions. We will examine them now, and generalize them.

 \xEh

 \xDH
We will assume that there are formulas TRUE and FALSE (by abuse of
language)
such that $f_{TRUE}(m)=TRUE$ and $f_{FALSE}(m)=FALSE$ for all $m,$ thus we
have at least
the two definable constant functions TRUE and FALSE, again by abuse of
language.

Not necessarily all constant functions are definable.

 \xDH
For each $x \xbe L$ there is $f_{x}$ defined by $f_{x}(m):=m(x).$

 \xDH
$ \xdd $ is closed under finite sup and finite inf.

 \xDH
$ \xdd $ will not necessarily be closed under complementation:

Given $f \xbe \xdd,$ the complement $ \xdc (f)$ is defined as above by:

$ \xdc (f)(m):=inf\{v \xbe V:sup\{v,f(m)\}=TRUE\}$ $=$ $ \xdc (f(m)).$

 \xDH
In classical logic, $ \xdd $ is closed under simplification:
If $X \xcc M$ is a definable model set, $L' \xcc L,$ then
$X' $ $:=$ $\{m \xbe M:$ $ \xcE m' \xbe X.m' \xex L' =m \xex L' \}$ is
definable.
This is a consequence of the existence of the standard normal forms,
consider e.g. the formula $p \xcu q,$ with $L=\{p,q\},$ $L' =\{p\},$ then
$X' =\{m,m' \},$
where $m(p)=m(q)=1,$ $m' (p)=1,$ $m' (q)=0,$ and the new formula is $p.$
We ``neglect'' or ``forget'' $q,$ take the projection.
It is also a sufficient condition for syntactic interpolation,
see Chapter \ref{Chapter Mod-Mon-Interpol} (page \pageref{Chapter
Mod-Mon-Interpol}).
(The following Example \ref{Example Gin-No-Def} (page \pageref{Example
Gin-No-Def})  shows that two different
formulas
might have the same model function, but should have different
projections.)

We have to define the analogon to $X' $ in many-valued logic.

Note that
 \xEI
 \xDH if $m \xex L' =m' \xex L',$ then $m \xbe X' \xcj m' \xbe X' $
 \xDH $m \xbe X' $ iff there is $m' \xbe X.m \xex L' =m' \xex L',$ thus
$f_{X' }(m)=sup\{f_{X}(m'):m \xex L' =m' \xex L' \}$
 \xEJ

So we impose the same conditions:
Let $f$ and $L' \xcc L$ be given, we look for suitable $f'.$

 \xEh
 \xDH $f' $ has to be indifferent to $L-L':$ if $m \xex L' =m' \xex L',$
then
we should have $f' (m)=f' (m').$
 \xDH $f' (m)=sup\{f(m'):m \xex L' =m' \xex L' \}.$
 \xEj

 \xDH When the model set has additional structure, we can ask whether the
resulting model choice functions preserve definability:

 \xEh
 \xDH the case of preferential structures is treated below in
Definition \ref{Definition Mod-Pref} (page \pageref{Definition Mod-Pref}),
 \xDH we can ask the same question
e.g. for modal structures: is the set of all models reachable from
some model definable, etc.
 \xEj

 \xEj

\be

$\hspace{0.01em}$


\label{Example Gin-No-Def}

This example shows that 2 different formulas $ \xbf $ and $ \xbf' $ may
define the
same $f_{ \xbf }=f_{ \xbf' },$ but neglecting a certain variable should
give
different results. We give two variants.

(1)
Set $ \xbf:=p \xco (q \xco \xCN q),$ $ \xbq:= \xCN p \xco (q \xco \xCN
q).$ So $f_{ \xbf }=f_{ \xbq },$ but negleting $q$
should result in $p$ in the first case, in $ \xCN p$ in the second case.

(2)
We work with 3 truth values, 0 for FALSE, 2 for TRUE, $ \xcu $ is as usual
interpreted by inf.
Define two new unary operators $K(x):=1$ (constant), $M(x):=min\{1,x\}.$

\vspace{5mm}
\begin{tabular}{|c|c|c|c|}

\hline

$a$ \xEH $b$ \xEH $\xbf = K(a) \xcu b$ \xEH $\xbf' = K(a) \xcu M(b)$ \xEP

\hline

0 \xEH 0 \xEH 0 \xEH 0 \xEP \hline

0 \xEH 1 \xEH 1 \xEH 1 \xEP \hline

0 \xEH 2 \xEH 1 \xEH 1 \xEP \hline

1 \xEH 0 \xEH 0 \xEH 0 \xEP \hline

1 \xEH 1 \xEH 1 \xEH 1 \xEP \hline

1 \xEH 2 \xEH 1 \xEH 1 \xEP \hline

2 \xEH 0 \xEH 0 \xEH 0 \xEP \hline

2 \xEH 1 \xEH 1 \xEH 1 \xEP \hline

2 \xEH 2 \xEH 1 \xEH 1 \xEP

\hline

\end{tabular}

\ee

So they define the same model function $f:M \xcp V.$
But when we forget about $ \xCf a,$ the first should just be $b,$ but the
second should
be $M(b).$

$ \xcz $
\\[3ex]

We may consider more systematically other operators, under which
the definable model sets should be closed:

 \xEh

 \xDH
constants for each truth value, like $FALSE=p \xcu \xCN p,$ $TRUE=p \xco
\xCN p$ in classical
logic,

 \xDH
complementation, if the complement is defined on the truth value set.
(In classical logic, this is, of course, negation.)

This might for instance be interesting for argumentation, where arguments,
or their sources, are the truth values.

 \xDH
functions similar to the basic operations SHL and SHR of computer science:
suppose a linear order be given on the truth values $0, \Xl,n,$ then
$SHR(p):=p+1$ if $p<n,$ and e.g. 0 if $p=n,$ etc.

The function $J$ of finite Goedel logics,
see Definition \ref{Definition Mod-Fin-Goed-Add} (page \pageref{Definition
Mod-Fin-Goed-Add}),
has some similarity to such shift
operations.

 \xEj
\subsubsection{
Generalization of model sets and (in)essential variables, overview
}

The Table \ref{Table Gin-Not-Def} (page \pageref{Table Gin-Not-Def})
also contains further material which will become clearer only
later.

$ \xCO $

\label{Gin-Nota-Defin-3}

The Table \ref{Table Gin-Not-Def} (page \pageref{Table Gin-Not-Def})
summarizes the situation for the 2-valued and the many-valued
case.
\begin{table}
\caption{Notation and Definitions}

\label{Table Gin-Not-Def}
\begin{center}
\tabcolsep=0.5pt.
\begin{tabular}{|c|c|c|}
\hline
\multicolumn{3}{|c|}{\bf Notation and definitions}\\
\hline
\multicolumn{3}{|c|}{

propositional language $L,$ $L' \xcc L,$ propositional variables $s, \Xl $
} \xEP
\hline

 \xEH 2-valued $\{0,1\}$ \xEH many-valued $(V, \xck)$ \xEP
\hline

definability of $f$ \xEH
\multicolumn{2}{|c|}{

$ \xcE \xbf:f_{ \xbf }=f$
} \xEP
\hline

model $m$ \xEH $m:L \xcp \{0,1\}$ \xEH $m:L \xcp V$ \xEP
$M$ set of all $L-$models \xEH \xEH \xEP
\hline

(for $ \xbG \xcc M)$ $ \xbG \xex L' $ \xEH
\multicolumn{2}{|c|}{

$ \xbG \xex L':=\{m \xex L':m \xbe \xbG \}$
} \xEP
\hline

$m \xex L' $ \xEH
\multicolumn{2}{|c|}{

like $m,$ but restricted to $L' $
} \xEP
\hline

$m \xCq_{L' }m' $ \xEH
\multicolumn{2}{|c|}{

$m \xCq_{L' }m' $ iff $ \xcA s \xbe L'.m(s)=m' (s)$
} \xEP
\hline

model $``set'' $ of formula $ \xbf $ \xEH $M(\xbf) \xcc M,$
$f_{ \xbf }:M \xcp \{0,1\}$ \xEH $f_{ \xbf }:M \xcp V$ \xEP
\hline

semantic equivalence of $ \xbf,$ $ \xbq $ \xEH
\multicolumn{2}{|c|}{

$f_{ \xbf }=f_{ \xbq }$
} \xEP
\hline

general model set \xEH $M' \xcc M,$
$f:M \xcp \{0,1\}$ \xEH $f:M \xcp V$ \xEP
\hline

$f$ insensitive to $L' $ \xEH
\multicolumn{2}{|c|}{

$ \xcA m,m' \xbe M.(m \xCq_{L-L' }m' \xch f(m)=f(m'))$
} \xEP
\hline

(ir)relevant \xEH
\multicolumn{2}{|c|}{

$s \xbe L$ is irrelevant for $f$ iff $f$ is insensitive to $s,$
} \xEP
\xEH \multicolumn{2}{|c|}{

$I(f):=\{s \xbe L:$ $s$ is irrelevant for $f\},$ $R(f):=L-I(f)$
} \xEP
\xEH \multicolumn{2}{|c|}{

$I(\xbf):=I(f_{ \xbf }),$ $R(\xbf):=R(f_{ \xbf })$
} \xEP
\hline

$f^{+}(m \xex L'),$ $f^{-}(m \xex L')$ \xEH
\multicolumn{2}{|c|}{

$f^{+}(m \xex L'):=max\{f(m'):m' \xbe M,$ $m \xCq_{L' }m' \}$
} \xEP

 \xEH
\multicolumn{2}{|c|}{

$f^{-}(m \xex L'):=min\{f(m'):m' \xbe M,$ $m \xCq_{L' }m' \}$
} \xEP
\hline

$f \xck g$ \xEH
\multicolumn{2}{|c|}{

$ \xcA m \xbe M.f(m) \xck g(m)$
} \xEP
\hline
\end{tabular}
\end{center}
\end{table}

$ \xCO $
\subsubsection{
Interpolation of many valued logics
}

\label{Section Eq-Sem}

\bd

$\hspace{0.01em}$


\label{Definition Eq-Interpol}

(3) Given $f,g,h:M \xcp V,$ we say that $h$ is a semantic interpolant for
$f$ and $g$ iff

(3.1) $ \xcA m \xbe M(f(m) \xck h(m) \xck g(m)),$

(3.2) $I(f) \xcv I(g) \xcc I(h)$

(4) Given $ \xbf,$ $ \xbq,$ we say that $ \xba $ is a syntactic
interpolant for $ \xbf $ and $ \xbq $
iff

(4.1) $ \xcA m \xbe M(f_{ \xbf }(m) \xck f_{ \xba }(m) \xck f_{ \xbq
}(m)),$

(4.2) all variables occuring in $ \xba $ occur also in $ \xbf $ and $ \xbq
.$

(5) The following will be central for constructing a semantical
interpolant:

Let $L' \xcc L,$ $m \xbe M,$ $m \xex L':L' \xcp V$ be the restriction of
a model $m$ to $L',$ $f:M \xcp V,$ then

$f^{+}(m \xex L'):=max\{f(m'):m \xCq_{L' }m' \},$ the maximal value for
any $m' $ which agrees with $m$
on $L',$ and

$f^{-}(m \xex L'):=min\{f(m'):m \xCq_{L' }m' \},$ the minimal value for
any $m' $ which agrees with $m$
on $L'.$
\clearpage
\section{
Preferential structures
}

\ed

An important part of these notes is motivated or concerns directly
preferential structures, which we now define.
\subsection{
The minimal variant
}

$ \xCO $
\index{Definition Preferential structure}

\bd

$\hspace{0.01em}$


\label{Definition Pref-Str}

Fix $U \xEd \xCQ,$ and consider arbitrary $X.$
Note that this $X$ has not necessarily anything to do with $U,$ or $ \xdu
$ below.
Thus, the functions $ \xbm_{ \xdm }$ below are in principle functions from
$V$ to $V$ - where $V$
is the set theoretical universe we work in.

Note that we work here often with copies of elements (or models).
In other areas of logic, most authors work with valuation functions. Both
definitions - copies or valuation functions - are equivalent, a copy
$ \xBc x,i \xBe $ can be seen as a state $ \xBc x,i \xBe $ with valuation $x.$
In the
beginning
of research on preferential structures, the notion of copies was widely
used, whereas e.g., [KLM90] used that of valuation functions. There is
perhaps
a weak justification of the former terminology. In modal logic, even if
two states have the same valid classical formulas, they might still be
distinguishable by their valid modal formulas. But this depends on the
fact
that modality is in the object language. In most work on preferential
stuctures, the consequence relation is outside the object language, so
different states with same valuation are in a stronger sense copies of
each other.

\ed

 \xEh

 \xDH $ \xCf Preferential$ $ \xCf models$ or $ \xCf structures.$
 \index{preferential model}
 \index{preferential structure}

 \xEh

 \xDH The version without copies:

A pair $ \xdm:= \xBc U, \xeb  \xBe $ with $U$ an arbitrary set, and $ \xeb $ an
arbitrary binary relation
on $U$ is called a $ \xCf preferential$ $ \xCf model$ or $ \xCf
structure.$

 \xDH The version with copies \index{copies}:

A pair $ \xdm:= \xBc  \xdu, \xeb  \xBe $ with $ \xdu $ an arbitrary set of
pairs,
and $ \xeb $ an arbitrary binary
relation on $ \xdu $ is called a $ \xCf preferential$ $ \xCf model$ or $
\xCf structure.$

If $ \xBc x,i \xBe  \xbe \xdu,$ then $x$ is intended to be an element of $U,$
and
$i$ the index of the
copy.

We sometimes also need copies of the relation $ \xeb.$ We will then
replace $ \xeb $
by one or several arrows $ \xba $ attacking non-minimal elements, e.g., $x
\xeb y$ will
be written $ \xba:x \xcp y$ \index{$ \xba:x \xcp y$}, $ \xBc x,i \xBe  \xeb 
\xBc y,i \xBe $ will
be written
$ \xba: \xBc x,i \xBe  \xcp  \xBc y,i \xBe $ \index{$ \xba: \xBc x,i \xBe  \xcp 
\xBc y,i \xBe $}, and
finally we might have $ \xBc  \xba,k \xBe:x \xcp y$ \index{$ \xBc  \xba,k \xBe
:x \xcp y$}  and
$ \xBc  \xba,k \xBe: \xBc x,i \xBe  \xcp  \xBc y,i \xBe $ \index{$ \xBc  \xba,k
\xBe: \xBc x,i \xBe  \xcp  \xBc y,i \xBe $}, etc.

 \xEj

 \xDH $ \xCf Minimal$ $ \xCf elements,$ the functions $ \xbm_{ \xdm }$ \index{$
\xbm_{ \xdm }$}
 \index{minimal element}

 \xEh

 \xDH The version without copies:

Let $ \xdm:= \xBc U, \xeb  \xBe,$ and define

$ \xbm_{ \xdm }(X)$ $:=$ $\{x \xbe X:$ $x \xbe U$ $ \xcu $ $ \xCN \xcE x'
\xbe X \xcs U.x' \xeb x\}.$

$ \xbm_{ \xdm }(X)$ is called the set of $ \xCf minimal$ $ \xCf elements$
of $X$ (in $ \xdm).$

Thus, $ \xbm_{ \xdm }(X)$ is the set of elements such that there is no
smaller one
in $X.$

 \xDH The version with copies:

Let $ \xdm:= \xBc  \xdu, \xeb  \xBe $ be as above. Define

$ \xbm_{ \xdm }(X)$ $:=$ $\{x \xbe X:$ $ \xcE  \xBc x,i \xBe  \xbe \xdu. \xCN
\xcE
\xBc x',i'  \xBe  \xbe \xdu (x' \xbe X$ $ \xcu $ $ \xBc x',i'  \xBe' \xeb  \xBc
x,i \xBe)\}.$

Thus, $ \xbm_{ \xdm }(X)$ is the projection on the first coordinate of the
set of elements
such that there is no smaller one in $X.$

Again, by abuse of language, we say that $ \xbm_{ \xdm }(X)$ is the set of
$ \xCf minimal$ $ \xCf elements$
of $X$ in the structure. If the context is clear, we will also write just
$ \xbm.$

We sometimes say that $ \xBc x,i \xBe $
``$ \xCf kills$'' or ``$ \xCf minimizes$'' $ \xBc y,j \xBe $ if
 \index{kill}
 \index{minimize}
$ \xBc x,i \xBe  \xeb  \xBc y,j \xBe.$ By abuse of language we also say a set
$X$ $ \xCf
kills$ or $ \xCf minimizes$ a set
$Y$ if for all $ \xBc y,j \xBe  \xbe \xdu,$ $y \xbe Y$ there is $ \xBc x,i \xBe 
\xbe \xdu,$
$x \xbe X$ s.t. $ \xBc x,i \xBe  \xeb  \xBc y,j \xBe.$

$ \xdm $ is also called $ \xCf injective$ or 1-copy \index{1-copy},
iff there is always at most one
copy
 \index{injective}
$ \xBc x,i \xBe $ for each $x.$ Note that the existence of copies corresponds to
a
non-injective labelling function - as is often used in nonclassical
logic, e.g., modal logic.

 \xEj

 \xEj

We say that $ \xdm $ is $ \xCf transitive,$ $ \xCf irreflexive,$ etc., iff
$ \xeb $ is.
 \index{transitive}
 \index{irreflexive}

Note that $ \xbm (X)$ might well be empty, even if $X$ is not.

$ \xCO $

$ \xCO $
\index{Definition Preferential logics}

\bd

$\hspace{0.01em}$


\label{Definition Pref-Log}

We define the consequence relation \index{consequence relation}  of a
preferential
structure for a
given propositional language $ \xdl.$

 \xEh

 \xDH

 \xEh

 \xDH If $m$ is a classical model of a language $ \xdl,$ we say by abuse
of language

$ \xBc m,i \xBe  \xcm \xbf $ iff $m \xcm \xbf,$

and if $X$ is any set of such pairs, that

$X \xcm \xbf $ iff for all $ \xBc m,i \xBe  \xbe X$ $m \xcm \xbf.$

 \xDH If $ \xdm $ is a preferential structure, and $X$ is a set of $ \xdl
-$models for a
classical propositional language $ \xdl,$ or a set of pairs $ \xBc m,i \xBe,$
where the $m$ are
such models, we call $ \xdm $ a $ \xCf classical$ $ \xCf preferential$ $
\xCf structure$ or $ \xCf model.$

 \xEj

 \xDH

$ \xCf Validity$ in a preferential structure, or the $ \xCf semantical$ $
\xCf consequence$ $ \xCf relation$
defined by such a structure:

Let $ \xdm $ be as above.

We define:

$T \xcm_{ \xdm } \xbf $ \index{$T \xcm_{ \xdm } \xbf $}  iff $ \xbm_{ \xdm
}(M(T)) \xcm
\xbf,$ i.e., $ \xbm_{ \xdm }(M(T)) \xcc M(\xbf).$

 \xDH

$ \xdm $ will be called $ \xCf definability$ $ \xCf preserving$ iff for
all $X \xbe \xdD_{ \xdl }$ $ \xbm_{ \xdm }(X) \xbe \xdD_{ \xdl }.$

 \xEj

As $ \xbm_{ \xdm }$ is defined on $ \xdD_{ \xdl },$ but need by no means
always result in some new
definable set, this is (and reveals itself as a quite strong) additional
property.

\ed

$ \xCO $

$ \xCO $
\index{Definition Smoothness}

\bd

$\hspace{0.01em}$


\label{Definition Smooth}

Let $ \xdy \xcc \xdp (U).$ (In applications to logic, $ \xdy $ will be $
\xdD_{ \xdl }.)$

A preferential structure $ \xdm $ is called $ \xdy -$smooth \index{$ \xdy
-$smooth}  iff for every $X \xbe \xdy $ every
element
$x \xbe X$ is either minimal in $X$ or above an element, which is minimal
in $X.$ More
precisely:

 \xEh

 \xDH The version without copies:

If $x \xbe X \xbe \xdy,$ then either $x \xbe \xbm (X)$ or there is $x'
\xbe \xbm (X).x' \xeb x.$

 \xDH The version with copies:

If $x \xbe X \xbe \xdy,$ and $ \xBc x,i \xBe  \xbe \xdu,$ then either there is
no
$ \xBc x',i'  \xBe  \xbe \xdu,$ $x' \xbe X,$
$ \xBc x',i'  \xBe  \xeb  \xBc x,i \xBe $ or there is $ \xBc x',i'  \xBe  \xbe
\xdu,$
$ \xBc x',i'  \xBe  \xeb  \xBc x,i \xBe,$ $x' \xbe X,$ s.t. there is
no $ \xBc x'',i''  \xBe  \xbe \xdu,$ $x'' \xbe X,$
with $ \xBc x'',i''  \xBe  \xeb  \xBc x',i'  \xBe.$

(Writing down all details here again might make it easier to read
applications
of the definition later on.)

 \xEj

When considering the models of a language $ \xdl,$ $ \xdm $ will be
called $ \xCf smooth$ iff
 \index{smooth}
it is $ \xdD_{ \xdl }-$smooth \index{$ \xdD_{ \xdl }-$smooth} ; $ \xdD_{ \xdl }$
is the
default.

Obviously, the richer the set $ \xdy $ is, the stronger the condition $
\xdy -$smoothness
will be.

A remark for the intuition:
Smoothness is perhaps best motivated through Gabbay's concept of
reactive diagrams, see, e.g.,
 \cite{Gab04} and  \cite{Gab08}, and also
 \cite{GS08c},  \cite{GS08f}.
In this concept, smaller, or ``better'' elements attack bigger, or
``less good'' elements. But when $ \xCf a$ attacks $b,$ and $b$ attacks $c,$
then one
might consider the attack of $b$ against $c$ weakened by the attack of
$ \xCf a$ against $b.$ In a smooth structure, for every attack against
some element $x,$
there is also an uncontested attack against $x,$ as it originates in an
element
$y,$ which is not attacked itself.

\ed

$ \xCO $

$ \xCO $

\bfa

$\hspace{0.01em}$


\label{Fact Rank-Base}

Let $ \xeb $ be an irreflexive, binary relation on $X,$ then the following
two conditions
are equivalent:

(1) There is $ \xbO $ and an irreflexive, total, binary relation $ \xeb'
$ on $ \xbO $ and a
function $f:X \xcp \xbO $ s.t. $x \xeb y$ $ \xcj $ $f(x) \xeb' f(y)$ for
all $x,y \xbe X.$

(2) Let $x,y,z \xbe X$ and $x \xcT y$ wrt. $ \xeb $ (i.e., neither $x \xeb
y$ nor $y \xeb x),$ then $z \xeb x$ $ \xch $ $z \xeb y$
and $x \xeb z$ $ \xch $ $y \xeb z.$

\efa

$ \xCO $

$ \xCO $

\subparagraph{
Proof
}

$\hspace{0.01em}$


$(1) \xch (2)$: Let $x \xcT y,$ thus neither $f(x) \xeb' f(y)$ nor $f(y)
\xeb' f(x),$ but
then $f(x)=f(y).$ Let now
$z \xeb x,$ so $f(z) \xeb' f(x)=f(y),$ so $z \xeb y.$ $x \xeb z$ $ \xch $
$y \xeb z$ is similar.

$(2) \xch (1)$: For $x \xbe X$ let $[x]:=\{x' \xbe X:x \xcT x' \},$ and $
\xbO:=\{[x]:x \xbe X\}.$ For $[x],[y] \xbe \xbO $
let $[x] \xeb' [y]: \xcj x \xeb y.$ This is well-defined: Let $x \xcT x'
,$ $y \xcT y' $ and $x \xeb y,$ then
$x \xeb y' $ and $x' \xeb y'.$ Obviously, $ \xeb' $ is an irreflexive,
total binary relation.
Define $f:X \xcp \xbO $ by $f(x):=[x],$ then $x \xeb y \xcj [x] \xeb' [y]
\xcj f(x) \xeb' f(y).$ $ \xcz $
\\[3ex]

$ \xCO $

$ \xCO $
\index{Definition Ranked relation}

\bd

$\hspace{0.01em}$


\label{Definition Rank-Rel}

We call an irreflexive, binary relation $ \xeb $ on $X,$ which satisfies
(1)
(equivalently (2)) of Fact \ref{Fact Rank-Base} (page \pageref{Fact Rank-Base})
, ranked \index{ranked}.
By abuse of language, we also call a preferential structure $ \xBc X, \xeb  \xBe
$
ranked, iff
$ \xeb $ is.

\ed

$ \xCO $

We quote from  \cite{Sch04} the following summary for
preferential structures:

$ \xCO $
\index{Proposition Pref-Representation}

Table \ref{Table Base1-Pref-Rep} (page \pageref{Table Base1-Pref-Rep}),
``Preferential representation'',
summarizes the more difficult half of a full
representation result for preferential structures. It shows equivalence
between certain abstract conditions for model choice functions and certain
preferential structures. They are shown in the respective
representation theorems.

``singletons'' means that the domain must contain all singletons,
``1 copy''
or ``$ \xcg 1$ copy'' means that the structure may contain only 1 copy for
each point,
or several, ``$(\xbm \xCQ)$'' etc. for the preferential structure mean
that the
$ \xbm -$function of the structure has to satisfy this property.

We call a characterization ``normal'' iff it is a
 \index{normal characterization}
universally quantified boolean combination (of any fixed, but perhaps
infinite, length) of rules of the usual form. We do not go into
details here.

In the second column from the left
``$ \xch $'' means, for instance for the
smooth case, that for any $ \xdy $ closed under finite unions, and any
choice function
$f$ which satisfies the conditions
in the left hand column, there is a (here $ \xdy -$smooth) preferential
structure $ \xdx $ which represents it, i.e., for all $Y \xbe \xdy $
$f(Y)= \xbm_{ \xdx }(Y),$ where
$ \xbm_{ \xdx }$ is the model choice function of the structure $ \xdx.$
The inverse arrow $ \xci $ means that the model choice function for any
smooth $ \xdx $
defined on such $ \xdy $ will satisfy the conditions on the left.

\label{Proposition Pref-Representation}

\begin{table}[h]

\index{$(\xbm \xcc)$}
\index{reactive}
\index{$(LLE)$}
\index{$(CCL)$}
\index{$(SC)$}
\index{$(\xbm CUM)$}
\index{$(\xbm \xcc \xcd)$}
\index{essentially smooth}
\index{$(SC)$}
\index{$(\xcc \xcd)$}
\index{$(\xbm PR)$}
\index{$(RW)+$}
\index{$(PR)$}
\index{$(\xbm dp)$}
\index{normal characterization}
\index{exception set}
\index{smooth}
\index{$(CUM)$}
\index{$(\xbm=)$}
\index{ranked}
\index{$\xcg 1$ copy}
\index{$(\xbm=')$}
\index{$(\xbm\xFO)$}
\index{$(\xbm\xcv)$}
\index{$(\xbm\xcv')$}
\index{$(\xbm\xbe)$}
\index{$(\xbm RatM)$}
\index{$(\xbm \xCQ)$}
\index{1 copy}
\index{$(\xbm \xCQ fin)$}
\index{$(\xcv)$}
\index{singletons}
\index{$(RatM)$}
\index{$(RatM=)$}
\index{$(Log\xcv)$}
\index{$(Log\xcv')$}

\caption{Preferential representation}

\label{Table Base1-Pref-Rep}
\begin{center}

{\tiny

\begin{tabular}{|c|c|c|c|c|}

\hline

\multicolumn{5}{|c|}{\bf Preferential representation}\\
\hline

$\xbm-$ function
\xEH
\xEH
Pref.Structure
\xEH
\xEH
Logic
\xEP

\hline

$(\xbm \xcc)$
\xEH
$\xcj$
\xEH
reactive
\xEH
$\xcj$
\xEH
$(LLE)+(CCL)+$
\xEP

\xEH
\xEH
\xEH
\xEH
$(SC)$
\xEP

\hline

$(\xbm \xcc)+(\xbm CUM)$
\xEH
$\xch$ $(\xcs)$
\xEH
reactive +
\xEH
\xEH
\xEP

\xEH
\xEH
essentially smooth
\xEH
\xEH
\xEP

\hline

$(\xbm \xcc)+(\xbm \xcc \xcd)$
\xEH
$\xch$
\xEH
reactive +
\xEH
$\xcj$
\xEH
$(LLE)+(CCL)+$
\xEP

\xEH
\xEH
essentially smooth
\xEH
\xEH
$(SC)+(\xcc \xcd)$
\xEP

\hline

$(\xbm \xcc)+(\xbm CUM)+(\xbm \xcc \xcd)$
\xEH
$\xci$
\xEH
reactive +
\xEH
\xEH
\xEP

\xEH
\xEH
essentially smooth
\xEH
\xEH
\xEP

\hline

$(\xbm \xcc)+(\xbm PR)$
\xEH
$\xci$
\xEH
general
\xEH
$\xch$ $(\xbm dp)$
\xEH
$(LLE)+(RW)+$
\xEP

\xEH
\xEH
\xEH
\xEH
$(SC)+(PR)$
\xEP

\cline{2-2}
\cline{4-4}

\xEH
$\xch$
\xEH
\xEH
$\xci$
\xEH
\xEP

\cline{2-2}
\cline{4-4}

\xEH
\xEH
\xEH
$\xcH$ without $(\xbm dp)$
\xEH
\xEP

\cline{4-5}

\xEH
\xEH
\xEH
$\xcJ$ without $(\xbm dp)$
\xEH
any ``normal''
\xEP

\xEH
\xEH
\xEH
\xEH
characterization
\xEP

\xEH
\xEH
\xEH
\xEH
of any size
\xEP

\hline

$(\xbm \xcc)+(\xbm PR)$
\xEH
$\xci$
\xEH
transitive
\xEH
$\xch$ $(\xbm dp)$
\xEH
$(LLE)+(RW)+$
\xEP

\xEH
\xEH
\xEH
\xEH
$(SC)+(PR)$
\xEP

\cline{2-2}
\cline{4-4}

\xEH
$\xch$
\xEH
\xEH
$\xci$
\xEH
\xEP

\cline{2-2}
\cline{4-4}

\xEH
\xEH
\xEH
$\xcH$ without $(\xbm dp)$
\xEH
\xEP

\cline{4-5}

\xEH
\xEH
\xEH
$\xcj$ without $(\xbm dp)$
\xEH
using ``small''
\xEP

\xEH
\xEH
\xEH
\xEH
exception sets
\xEP

\hline

$(\xbm \xcc)+(\xbm PR)+(\xbm CUM)$
\xEH
$\xci$
\xEH
smooth
\xEH
$\xch$ $(\xbm dp)$
\xEH
$(LLE)+(RW)+$
\xEP

\xEH
\xEH
\xEH
\xEH
$(SC)+(PR)+$
\xEP

\xEH
\xEH
\xEH
\xEH
$(CUM)$
\xEP

\cline{2-2}
\cline{4-4}

\xEH
$\xch$ $(\xcv)$
\xEH
\xEH
$\xci$ $(\xcv)$
\xEH
\xEP

\cline{2-2}
\cline{4-4}

\xEH
$\xcH$ without $(\xcv)$
\xEH
\xEH
$\xcH$ without $(\xbm dp)$
\xEH
\xEP

\hline

$(\xbm \xcc)+(\xbm PR)+(\xbm CUM)$
\xEH
$\xci$
\xEH
smooth+transitive
\xEH
$\xch$ $(\xbm dp)$
\xEH
$(LLE)+(RW)+$
\xEP

\xEH
\xEH
\xEH
\xEH
$(SC)+(PR)+$
\xEP

\xEH
\xEH
\xEH
\xEH
$(CUM)$
\xEP

\cline{2-2}
\cline{4-4}

\xEH
$\xch$ $(\xcv)$
\xEH
\xEH
$\xci$ $(\xcv)$
\xEH
\xEP

\cline{2-2}
\cline{4-4}

\xEH
\xEH
\xEH
$\xcH$ without $(\xbm dp)$
\xEH
\xEP

\cline{4-5}

\xEH
\xEH
\xEH
$\xcj$ without $(\xbm dp)$
\xEH
using ``small''
\xEP

\xEH
\xEH
\xEH
\xEH
exception sets
\xEP

\hline

$(\xbm\xcc)+(\xbm=)+(\xbm PR)+$
\xEH
$\xci$
\xEH
ranked, $\xcg 1$ copy
\xEH
\xEH
\xEP

$(\xbm=')+(\xbm\xFO)+(\xbm\xcv)+$
\xEH
\xEH
\xEH
\xEH
\xEP

$(\xbm\xcv')+(\xbm\xbe)+(\xbm RatM)$
\xEH
\xEH
\xEH
\xEH
\xEP

\hline

$(\xbm\xcc)+(\xbm=)+(\xbm PR)+$
\xEH
$\xcH$
\xEH
ranked
\xEH
\xEH
\xEP

$(\xbm\xcv)+(\xbm\xbe)$
\xEH
\xEH
\xEH
\xEH
\xEP

\hline

$(\xbm\xcc)+(\xbm=)+(\xbm \xCQ)$
\xEH
$\xcj$, $(\xcv)$
\xEH
ranked,
\xEH
\xEH
\xEP

\xEH
\xEH
1 copy + $(\xbm \xCQ)$
\xEH
\xEH
\xEP

\hline

$(\xbm\xcc)+(\xbm=)+(\xbm \xCQ)$
\xEH
$\xcj$, $(\xcv)$
\xEH
ranked, smooth,
\xEH
\xEH
\xEP

\xEH
\xEH
1 copy + $(\xbm \xCQ)$
\xEH
\xEH
\xEP

\hline

$(\xbm\xcc)+(\xbm=)+(\xbm \xCQ fin)+$
\xEH
$\xcj$, $(\xcv)$, singletons
\xEH
ranked, smooth,
\xEH
\xEH
\xEP

$(\xbm\xbe)$
\xEH
\xEH
$\xcg$ 1 copy + $(\xbm \xCQ fin)$
\xEH
\xEH
\xEP

\hline

$(\xbm\xcc)+(\xbm PR)+(\xbm \xFO)+$
\xEH
$\xcj$, $(\xcv)$, singletons
\xEH
ranked
\xEH
$\xcH$ without $(\xbm dp)$
\xEH
$(RatM), (RatM=)$,
\xEP

$(\xbm \xcv)+(\xbm\xbe)$
\xEH
\xEH
$\xcg$ 1 copy
\xEH
\xEH
$(Log\xcv), (Log\xcv')$
\xEP

\cline{4-5}

\xEH
\xEH
\xEH
$\xcJ$ without $(\xbm dp)$
\xEH
any ``normal''
\xEP

\xEH
\xEH
\xEH
\xEH
characterization
\xEP

\xEH
\xEH
\xEH
\xEH
of any size
\xEP

\hline

\end{tabular}

}

\end{center}
\end{table}

$ \xCO $

For more detail on preferential logics and size, the reader is referred
to,
e.g.,  \cite{GS08c} and  \cite{GS09a}.

In Section \ref{Section Mod-Many} (page \pageref{Section Mod-Many}), we will
generalize the concept of a
preferential
structure to logics with more than two truth values, see
Definition \ref{Definition Mod-Pref} (page \pageref{Definition Mod-Pref}) 
there.
\subsubsection{
New material on the minimal variant of preferential structures (Bubble
structures)
}

\label{Section New-Min}

We discuss now a variant of preferential structures which will prove to be
useful. We call them bubble structures, see
Diagram \ref{Diagram Bubble} (page \pageref{Diagram Bubble})  for illustration.

The basic idea is as follows: We have one global structure, and parts of
it, ``bubbles'', behave in a uniform way. Thus, with respect to the
``outside world'', what is inside an individual bubble is indistinguishable.

Formally, if $x$ is outside a given bubble, and $b,b' $ are inside the
bubble,
then $x \xeb b$ iff $x \xeb b',$ and $b \xeb x$ iff $b' \xeb x.$ Thus,
seen from the outside, a bubble
behaves like a layer in a ranked structure. But we do not require the
inside
of the bubble to consist of only incomparable elements, like a layer in
a ranked structure. Thus, ranked structures are special cases of bubble
structures, with the layers being the bubbles.

We may identify the bubbles with single points, and, as long as we do not
look into the bubbles, we have just a usual preferential structure.

We can see the whole structure as an abstraction, where details are
encapsulated in the bubbles. We can also work with different languages,
a ``global'' language for the big structure, and a (or several) sublanguages
which are used only inside the bubbles. When all bubbles have the same
internal structure, we can also push this structure into the truth values,
see Section \ref{Section Mod-Lang-Struc} (page \pageref{Section Mod-Lang-Struc})
. The basic idea is, of course,
another
expression of modularity, where we try to isolate different aspects of
reasoning as much as possible, to simplify the task - divide et impera!

We will not give a full formal representation result here, but point out
the
important ingredients. The cases to consider (different sublanguages,
cumulativity or not inside/outside the bubbles, etc.) may
be quite varied, and we leave it to future research to elaborate the
details.
We turn now to a list of things to consider.

 \xEh

 \xDH
Copies of the same element should probably be in the same bubble.

 \xDH
We will probably work with one copy of each bubble, so we add a 1-copy
condition to the global structure, see
 \cite{Sch04} for a discussion.

 \xDH
Singletons without copies can always be considered to be bubbles.

 \xDH
Given a global structure, the decomposition into bubbles is usually not
unique: Take a linear order, then all elements of any interval $[a,b]$
behave to the outside world in the same way, so any such interval can be
considered a bubble. But, of course, we want bubbles to be disjoint, or,
at least, forming systems of bubbles, superbubbles, etc.

 \xDH
Similar to ranked structures, any element of a bubble replaces any other
for the order relation:

Suppose $A \xcs B= \xCQ,$ $B$ a bubble, then $ \xbm (A \xcv B)= \xbm (A)
\xcv \xbm (B),$ or $= \xbm (A),$ or
$= \xbm (B),$ provided the relation is transitive.

For a counterexample to the non-transitive case, consider $b \xeb a \xeb
b',$ but not
$b \xeb b',$ with $B=\{b,b' \},$ $A=\{a\},$ then $ \xbm (B)=\{b,b' \},$ $
\xbm (A \xcv B)=\{b\}.$

 \xDH
The full rankedness conditions like $A \xcc B,$ $ \xbm (B) \xcs A \xEd
\xCQ $ $ \xch $ $ \xbm (A)= \xbm (B) \xcs A$
are generally too strong, as the inside of the bubbles need not consist of
incomparable elements.

 \xEj

$ \xCO $

\vspace{10mm}

\begin{diagram}

\label{Diagram Bubble}
\index{Diagram Bubble}

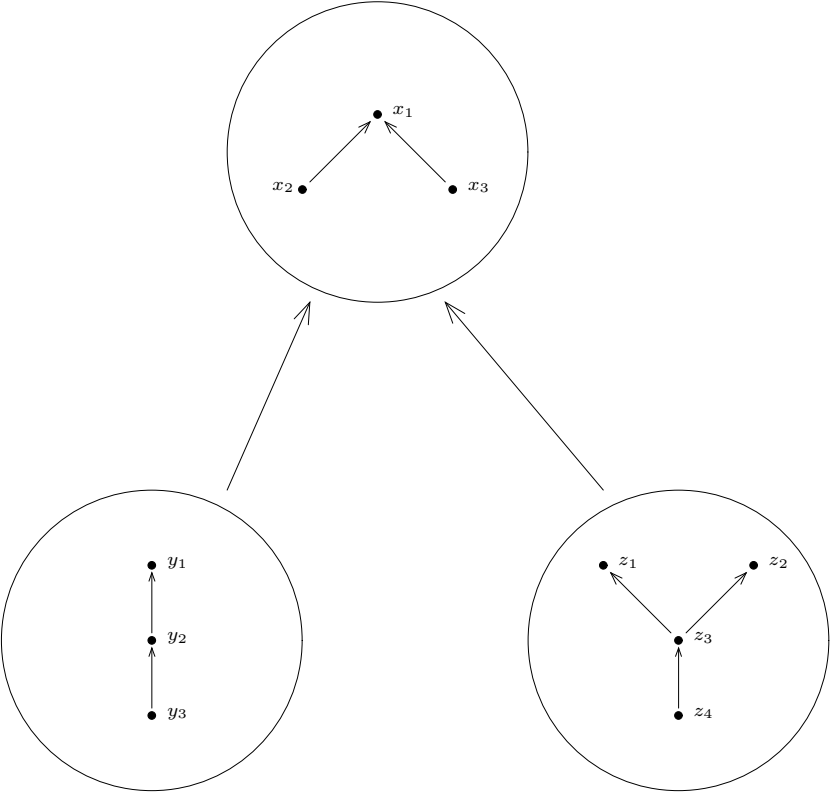
\centering
\setlength{\unitlength}{1mm}
{\renewcommand{\dashlinestretch}{30}
\begin{picture}(150,150)(0,0)

\put(50,100){\circle*{1}}
\put(40,90){\circle*{1}}
\put(60,90){\circle*{1}}

\put(52,100){{\xssc $x_1$}}
\put(36,90){{\xssc $x_2$}}
\put(62,90){{\xssc $x_3$}}

\path(41,91)(49,99)
\path(59,91)(51,99)
\path(47.5,98.3)(49,99)(48.3,97.5)
\path(51.7,97.5)(51,99)(52.5,98.3)

\put(50,95){\circle{40}}

\put(20,40){\circle*{1}}
\put(20,30){\circle*{1}}
\put(20,20){\circle*{1}}

\put(22,40){{\xssc $y_1$}}
\put(22,30){{\xssc $y_2$}}
\put(22,20){{\xssc $y_3$}}

\path(20,21)(20,29)
\path(20,31)(20,39)
\path(19.6,27.9)(20,29)(20.4,27.9)
\path(19.6,37.9)(20,39)(20.4,37.9)

\put(20,30){\circle{40}}

\put(80,40){\circle*{1}}
\put(100,40){\circle*{1}}
\put(90,30){\circle*{1}}
\put(90,20){\circle*{1}}

\put(82,40){{\xssc $z_1$}}
\put(102,40){{\xssc $z_2$}}
\put(92,30){{\xssc $z_3$}}
\put(92,20){{\xssc $z_4$}}

\path(90,21)(90,29)
\path(89,31)(81,39)
\path(91,31)(99,39)

\path(89.6,27.9)(90,29)(90.4,27.9)
\path(81.7,37.5)(81,39)(82.5,38.3)
\path(97.5,38.3)(99,39)(98.3,37.5)

\put(90,30){\circle{40}}

\path(30,50)(41,75)
\path(80,50)(59,75)
\path(38.9,72.8)(41,75)(40.8,72)
\path(60,72.2)(59,75)(61.6,73.5)

\end{picture}
}

\end{diagram}

\vspace{4mm}

$ \xCO $
\subsection{
The limit variant
}

\label{Section Pref-Lim}
\paragraph{
Motivation for the limit variant and for our approach \\[2mm]
}

$ \xCO $

\label{Motivation-MISE}

Distance based semantics give perhaps the clearest motivation for the
limit
variant. For instance,
the Stalnaker/Lewis semantics for counterfactual conditionals defines
$ \xbf > \xbq $ to hold in a (classical) model $m$ iff in those models of
$ \xbf,$ which are
closest to $m,$ $ \xbq $ holds. For this to make sense, we need, of
course, a distance
$d$ on the model set. We call this approach the minimal variant.
Usually, one makes a limit assumption: The set of $ \xbf -$models
closest to $m$ is not empty if $ \xbf $ is consistent - i.e., the $ \xbf
-$models are not
arranged around $m$ in a way that they come closer and closer, without a
minimal
distance. This is, of course, a very strong assumption, and which is
probably
difficult to justify philosophically. It seems to have its only
justification
in the fact that it avoids degenerate cases, where, in above example, for
consistent $ \xbf $ $m \xcm \xbf >FALSE$ holds. As such, this assumption
is unsatisfactory.

The limit version avoids such assumptions. It will still work in above
situation, i.e., when there are not always optimal (closest) elements, it
defines what happens when we get
``better and better'', i.e. approach the limit (the
``best'' case).

We will have to define what a suitable
``neighbourhood'' of the best cases is, in our context, this will roughly be
a set
of elements which minimizes all other elements, and is downward closed,
i.e.,
contains all elements better than some $x$ already in the set.
We call such sets MISE, for minimizing initial segment.
We will see
(Example \ref{Example D-8.2.1} (page \pageref{Example D-8.2.1})) that this
definition will not always do
what we
want it to do, and we will have to impose additional properties.

Essentially, we want MISE sets to reflect the properties of the sets of
minimal elements, if they exist. Thus, the set of minimal elements should
be
a special case of a MISE. But we also want MISE sets to be closed under
finite intersection, to have the logical (AND) property, see again
Example \ref{Example D-8.2.1} (page \pageref{Example D-8.2.1}). If our
definition is such that its
properties are
sufficiently close to those of the ideal (the minimal elements), then
we will also have the desired algebraic and logical properties, but avoid
pathologies originating from the empty set (when there are no best
elements) -
and this is what we wanted. Of course, our definition still has to
correspond
to the intuition what an approximation to the ideal case should be.

$ \xCO $

We give now the basic definitions for the limit version of preferential
and
ranked preferential structures.

\bd

$\hspace{0.01em}$


\label{Definition D-8.1.1}

(1) General preferential structures

(1.1) The version without copies:

Let $ \xdm:= \xBc U, \xeb  \xBe.$ Define

$Y \xcc X \xcc U$ is a minimizing initial segment \index{minimizing initial
segment},
or MISE \index{MISE}, of $X$ iff:

(a) $ \xcA x \xbe X \xcE x \xbe Y.y \xec x$ - where $y \xec x$ stands for
$x \xeb y$ or $x=y$
(i.e., $Y$ is minimizing)
and

(b) $ \xcA y \xbe Y, \xcA x \xbe X(x \xeb y$ $ \xch $ $x \xbe Y)$
(i.e., $Y$ is downward closed or an initial part).

(1.2) The version with copies:

Let $ \xdm:= \xBc  \xdu, \xeb  \xBe $ be as above. Define for $Y \xcc X \xcc
\xdu
$

$Y$ is a minimizing initial segment, or MISE of $X$ iff:

(a) $ \xcA  \xBc x,i \xBe  \xbe X \xcE  \xBc y,j \xBe  \xbe Y. \xBc y,j \xBe 
\xec  \xBc x,i \xBe $

and

(b) $ \xcA  \xBc y,j \xBe  \xbe Y, \xcA  \xBc x,i \xBe  \xbe X$
$(\xBc x,i \xBe  \xeb  \xBc y,j \xBe $ $ \xch $ $ \xBc x,i \xBe  \xbe Y).$

(1.3) For $X \xcc \xdu,$ let $ \xbL (X)$ \index{$ \xbL (X)$}  be the set of
MISE of $X.$

(1.4) We say that a set $ \xdx $ of MISE
is cofinal \index{cofinal}  in another set of MISE
$ \xdx' $ (for the same base set $X)$ iff for all $Y' \xbe \xdx',$
there is $Y \xbe \xdx,$ $Y \xcc Y'.$

(1.5) A MISE $X$ is called definable \index{definable}  iff
$\{x: \xcE i. \xBc x,i \xBe  \xbe X\} \xbe \xdD_{ \xdl }.$

(1.6) $T \xcm_{ \xdm } \xbf $ iff there is $Y \xbe \xbL (\xdu \xex M(T))$
such that $Y \xcm \xbf.$

$(\xdu \xex M(T):=\{ \xBc x,i \xBe  \xbe \xdu:x \xbe M(T)\}$ - if there are no
copies, we simplify in
the obvious way.)

(2) Ranked preferential structures

In the case of ranked structures, we may assume without loss of generality
that
the MISE sets have a particularly simple form:

For $X \xcc U$ $A \xcc X$ is MISE iff $X \xEd \xCQ $ and $ \xcA a \xbe A
\xcA x \xbe X(x \xeb a \xco x \xcT a$ $ \xch $ $x \xbe A).$
(A is downward and horizontally closed.)

(3) Theory Revision
 \index{theory revision}

Recall that we have a distance $d$ on the model set, and are interested
in $y \xbe Y$ which are close to $X.$

Thus, given $ \xCf X,Y,$ we define analogously:

$B \xcc Y$ is MISE iff

(1) $B \xEd \xCQ $

(2) there is $d' $ such that $B:=\{y \xbe Y: \xcE x \xbe X.d(x,y) \xck d'
\}$
(we could also have chosen $d(x,y)<d',$ this is not important).

And we define $ \xbf \xbe T*T' $ iff there is $B \xbe \xbL (M(T),M(T'))$
$B \xcm \xbf.$

\ed

There are basic problems with the limit in general preferential
structures, as we shall see now:

\be

$\hspace{0.01em}$


\label{Example D-8.2.1}

Let $a \xeb b,$ $a \xeb c,$ $b \xeb d,$ $c \xeb d$ (but $ \xeb $ not
transitive!), then $\{a,b\}$ and $\{a,c\}$ are
such $S$ and $S',$ but there is no $S'' \xcc S \xcs S' $ which is an
initial segment. If, for
instance, in $ \xCf a$ and $b$ $ \xbq $ holds, in $ \xCf a$ and $c$ $ \xbq
',$ then ``in the limit'' $ \xbq $ and $ \xbq' $
will hold, but not $ \xbq \xcu \xbq'.$ This does not seem right. We
should not be obliged
to give up $ \xbq $ to obtain $ \xbq'.$ $ \xcz $
\\[3ex]

\ee

We will therefore require it to be closed under
finite intersections, or at least, that if $S,$ $S' $ are such segments,
then there
must be $S'' \xcc S \xcs S' $ which is also such a segment.

We make this official. Let $ \xbL (X)$ be the set of initial segments of
$X,$ then
we require:

$(\xbL \xcs)$ \index{$(\xbL \xcs)$}  If $A,B \xbe \xbL (X)$ then there is $C
\xcc A \xcs B,$ $C \xbe \xbL (X).$

To familiarize the reader with the limit version, we show two easy
but important results.

\bfa

$\hspace{0.01em}$


\label{Fact D-8.2.1}

(Taken from  \cite{Sch04}, Fact 3.4.3, Proposition 3.10.16 there,
(2a) is new, but only a summary of other properties.)

Let the relation $ \xeb $ be transitive. The following hold in the limit
variant of
general preferential structures:

(1) If $A \xbe \xbL (Y),$ and $A \xcc X \xcc Y,$ then $A \xbe \xbL (X).$

(2) If $A \xbe \xbL (Y),$ and $A \xcc X \xcc Y,$ and $B \xbe \xbL (X),$
then $A \xcs B \xbe \xbL (Y).$

(2a) We summarize to make finitary semantic cumulativity evident:
Let $A \xbe \xbL (Y),$ $A \xcc X \xcc Y.$ Then, if $B \xbe \xbL (Y),$ $A
\xcs B \xbe \xbL (X).$ Conversely,
if $B \xbe \xbL (X),$ then $A \xcs B \xbe \xbL (Y).$

(3) If $A \xbe \xbL (Y),$ $B \xbe \xbL (X),$ then there is $Z \xcc A \xcv
B$ $Z \xbe \xbL (Y \xcv X).$

The following hold in the limit variant of ranked structures without
copies,
where the domain is closed under finite unions and contains all finite
sets.

(4) $A,B \xbe \xbL (X)$ $ \xch $ $A \xcc B$ or $B \xcc A,$

(5) $A \xbe \xbL (X),$ $Y \xcc X,$ $Y \xcs A \xEd \xCQ $ $ \xch $ $Y \xcs
A \xbe \xbL (Y),$

(6) $ \xbL' \xcc \xbL (X),$ $ \xcS \xbL' \xEd \xCQ $ $ \xch $ $ \xcS
\xbL' \xbe \xbL (X).$

(7) $X \xcc Y,$ $A \xbe \xbL (X)$ $ \xch $ $ \xcE B \xbe \xbL (Y).B \xcs
X=A$

\efa

\subparagraph{
Proof
}

$\hspace{0.01em}$


(1) trivial.

(2)

(2.1) $A \xcs B$ is closed in $Y:$ Let $ \xBc x,i \xBe  \xbe A \xcs B,$
$ \xBc y,j \xBe  \xeb  \xBc x,i \xBe,$ then $ \xBc y,j \xBe  \xbe A.$ If
$ \xBc y,j \xBe  \xce X,$ then $ \xBc y,j \xBe  \xce A,$ $contradiction.$
So $ \xBc y,j \xBe  \xbe X,$ but then $ \xBc y,j \xBe  \xbe B.$

(2.2) $A \xcs B$ minimizes $Y:$ Let $ \xBc a,i \xBe  \xbe Y.$

(a) If $ \xBc a,i \xBe  \xbe A-B \xcc X,$ then there is
$ \xBc y,j \xBe  \xeb  \xBc a,i \xBe,$ $ \xBc y,j \xBe  \xbe B.$ Xy closure of
A,
$ \xBc y,j \xBe  \xbe A.$

(b) If $ \xBc a,i \xBe  \xce A,$ then there is
$ \xBc a',i'  \xBe  \xbe A \xcc X,$ $ \xBc a',i'  \xBe  \xeb  \xBc a,i \xBe,$
continue by (a).

(2a) For the first part,
by (2), $A \xcs B \xbe \xbL (Y),$ so by (1), $A \xcs B \xbe \xbL (X).$ The
second part is just (2).

(3)

Let $Z$ $:=$ $\{ \xBc x,i \xBe  \xbe A$:
$ \xCN \xcE  \xBc b,j \xBe  \xec  \xBc x,i \xBe. \xBc b,j \xBe  \xbe X-B\}$ $
\xcv $
$\{ \xBc y,j \xBe  \xbe B$:
$ \xCN \xcE  \xBc a,i \xBe  \xec  \xBc y,j \xBe. \xBc a,i \xBe  \xbe Y-A\},$
where $ \xec $ stands for $ \xeb $ or $=.$

(3.1) $Z$ minimizes $Y \xcv X:$ We consider $Y,$ $X$ is symmetrical.

(a) We first show: If $ \xBc a,k \xBe  \xbe A-Z$, then there is
$ \xBc y,i \xBe  \xbe Z. \xBc a,k \xBe  \xee  \xBc y,i \xBe.$
Proof: If $ \xBc a,k \xBe  \xbe A-$Z,
then there is $ \xBc b,j \xBe  \xec  \xBc a,k \xBe,$ $ \xBc b,j \xBe  \xbe
X-$B. Then there
is $ \xBc y,i \xBe  \xeb  \xBc b,j \xBe,$ $ \xBc y,i \xBe  \xbe B.$
Xut $ \xBc y,i \xBe  \xbe Z,$ too: If not, there would be
$ \xBc a',k'  \xBe  \xec  \xBc y,i \xBe,$ $ \xBc a',k'  \xBe  \xbe Y- \xCf A,$
but
$ \xBc a',k'  \xBe  \xeb  \xBc a,k \xBe,$ contradicting closure of $ \xCf A.$

(b) If $ \xBc a'',k''  \xBe  \xbe Y- \xCf A,$ there is $ \xBc a,k \xBe  \xbe A,$
$ \xBc a,k \xBe  \xeb  \xBc a'',k''  \xBe.$ If $ \xBc a,k \xBe  \xce Z,$
continue
with (a).

(3.2) $Z$ is closed in $Y \xcv X:$ Let then $ \xBc z,i \xBe  \xbe Z,$
$ \xBc u,k \xBe  \xeb  \xBc z,i \xBe,$ $ \xBc u,k \xBe  \xbe Y \xcv X.$
Suppose $ \xBc z,i \xBe  \xbe A$ - the case
$ \xBc z,i \xBe  \xbe B$ is symmetrical.

(a) $ \xBc u,k \xBe  \xbe Y-A$ cannot be, by closure of $ \xCf A.$

(b) $ \xBc u,k \xBe  \xbe X-B$ cannot be, as $ \xBc z,i \xBe  \xbe Z,$ and by
definition of
$Z.$

(c) If $ \xBc u,k \xBe  \xbe A- \xCf Z,$ then there is
$ \xBc v,l \xBe  \xec  \xBc u,k \xBe,$ $ \xBc v,l \xBe  \xbe X-$B, so $ \xBc
v,l \xBe  \xeb  \xBc z,i \xBe,$
contradicting (b).

(d) If $ \xBc u,k \xBe  \xbe B-Z$, then there is
$ \xBc v,l \xBe  \xec  \xBc u,k \xBe,$ $ \xBc v,l \xBe  \xbe Y- \xCf A,$
contradicting (a).

(4) Suppose not, so there are $a \xbe A- \xCf B,$ $b \xbe B- \xCf A.$ But
if $a \xcT b,$ $a \xbe B$ and $b \xbe A,$
similarly if $a \xeb b$ or $b \xeb a.$

(5) As $A \xbe \xbL (X)$ and $Y \xcc X,$ $Y \xcs A$ is downward and
horizontally closed. As $Y \xcs A \xEd \xCQ,$
$Y \xcs A$ minimizes $Y.$

(6) $ \xcS \xbL' $ is downward and horizontally closed, as all $A \xbe
\xbL' $ are. As $ \xcS \xbL' \xEd \xCQ,$
$ \xcS \xbL' $ minimizes $X.$

(7) Set $B:=\{b \xbe Y: \xcE a \xbe A.a \xcT b$ or $b \xck a\}$

$ \xcz $
\\[3ex]

We have as immediate logical consequence:

\bfa

$\hspace{0.01em}$


\label{Fact D-8.2.2}

(Fact 3.4.4 of  \cite{Sch04}.)

If $ \xeb $ is transitive, then in the limit variant hold:

(1) $ \xCf (AND)$ \index{$ \xCf (AND)$},

(2) $ \xCf (OR)$ \index{$ \xCf (OR)$}.

\efa

\subparagraph{
Proof
}

$\hspace{0.01em}$


Let $ \xdz $ be the structure.

(1) Immediate by Fact \ref{Fact D-8.2.1} (page \pageref{Fact D-8.2.1}), (2) -
set $A=B.$

(2) Immediate by Fact \ref{Fact D-8.2.1} (page \pageref{Fact D-8.2.1}), (3).
$ \xcz $
\\[3ex]

We also have

\bfa

$\hspace{0.01em}$


\label{Fact Lim-Cum}

(Fact 3.4.5 in  \cite{Sch04})

Finite cumulativity holds in transitive limit structures:

If $ \xbf \xcn \xbq,$ then $ \ol{ \ol{ \xbf } }= \ol{ \ol{ \xbf \xcu \xbq
} }.$

See  \cite{Sch04} for a direct proof, or above
Fact \ref{Fact D-8.2.1} (page \pageref{Fact D-8.2.1}), (2a). $ \xcz $
\\[3ex]

\efa

We repeat now (without proof) our main logical trivialization results on
the
limit variant of general preferential structures,
Proposition 3.4.7 and Proposition 3.10.19 from  \cite{Sch04}:

\bp

$\hspace{0.01em}$


\label{Proposition D-8.3.5a}

(1)
Let the relation be transitive.
If we consider only formulas on the left of $ \xcn,$ the resulting logic
of the
limit version can also be generated by the minimal version of a (perhaps
different) preferential structure. Moreover, this structure can be chosen
smooth.

(2)
Let a logic $ \xbf \xcn \xbq $ be given by the limit variant of a ranked
structure without
copies. Then there is a ranked structure, which gives
exactly the same logic, but interpreted in the minimal variant.

$ \xcz $
\\[3ex]

\ep

(The negative results for the general not definability preserving minimal
case
apply also to the general limit case - see Section 5.2.3 in
 \cite{Sch04} for details.)
\subsubsection{
New material on the limit variant of preferential structures
}

\label{Section New-Lim}

This short section contains new material on the limit variant of
preferential
structures - a discussion of the limit variant of higher preferential
structures
will be presented below, see
Section \ref{Section Lim-High} (page \pageref{Section Lim-High}).

Consider the following analogon to $(\xbm PR)$ $(A \xcc B$ $ \xch $ $
\xbm (B) \xcs A \xcc \xbm (A)):$

\bfa

$\hspace{0.01em}$


\label{Fact Lim-MuPR}

Let $ \xeb $ be transitive, $ \xbL (X)$ the MISE systems over $X.$

Let $A \xcc B,$ $A' \xbe \xbL (A)$ $ \xch $ $ \xcE B' \xbe \xbL (B).B'
\xcs A \xcc A'.$

\efa

\subparagraph{
Proof
}

$\hspace{0.01em}$


Consider $B' $ $:=$ $\{b \xbe B:$ $b \xce A-A' $ and $ \xCN \xcE b' \xbe
A-A'.b' \xeb b\}.$
Thus, $B-B' =\{b \xbe B:$ $b \xbe A-A' $ or $ \xcE b' \xbe A-A'.b' \xeb
b\}.$

(1) $B' \xcs A \xcc A':$ Trivial.

(2) $B' $ is closed in $B:$ Let $b \xbe B',$ suppose there is $b' \xbe
B-B',$ $b' \xeb b.$ $b' \xbe A-A' $
is excluded by definition, $b' $ such that $ \xcE b'' \xbe A-A'.b'' \xeb
b' $ by transitivity.

(3) $B' $ is minimizing: Let $b \xbe B-B'.$ If $b \xbe A-A',$ then there
is $a \xbe A'.a \xeb b$ by
minimization of $ \xCf A$ by $A'.$ We have to show that $a \xbe B'.$ If
not, there must be
$b' \xbe A-A'.b' \xeb a,$ contradicting closure of $A' $ in $ \xCf A.$ If
$b$ is such that there
is $b' \xbe A-A'.b' \xeb b.$ Then there has to be $a \xbe A' $ such that
$a \xeb b' \xeb b,$ so $a \xeb b$
by transitivity, and we continue as above.

$ \xcz $
\\[3ex]

We have immediately:

\bco

$\hspace{0.01em}$


\label{Corollary Lim-PR}

$ \ol{ \ol{ \xbf \xcu \xbf' } } \xcc \ol{ \ol{ \ol{ \xbf } } \xcv \{ \xbf
' \}}.$

\eco

\subparagraph{
Proof
}

$\hspace{0.01em}$


Let $ \xbq \xbe \ol{ \ol{ \xbf \xcu \xbf' } },$ $A:=M(\xbf \xcu \xbf'
),$ $B:=M(\xbf).$ So there is $A' \xbe \xbL (A).A' \xcm \xbq,$ so
there is $B' \xbe \xbL (B).B' \xcs M(\xbf') \xcm \xbq $ by
Fact \ref{Fact Lim-MuPR} (page \pageref{Fact Lim-MuPR}), so $B' \xcm \xbf'
\xcp \xbq,$ so $ \xbf'
\xcp \xbq \xbe \ol{ \ol{ \xbf } },$ so $ \ol{ \ol{ \xbf } } \xcv \xbf'
\xcl \xbq.$
$ \xcz $
\\[3ex]
\subsection{
Preferential structures for many-valued logics
}

We can, of course, consider for given $ \xbf $ the set of models where $
\xbf $ has
maximal truth value TRUE, and then take the minimal ones as usual.
The resulting logic $ \xcn $ then makes $ \xbf \xcn \xbq $ true, iff the
minimal models
with value TRUE assign TRUE also to $ \xbq.$
See Section \ref{Section EQ} (page \pageref{Section EQ}).

But this does not seem to be the adequate way. So we adapt the definition
of
preferential structures to the many-valued situation.

\bd

$\hspace{0.01em}$


\label{Definition Mod-Pref}

Let $ \xdl $ be given with model set $M.$

Let a binary relation $ \xeb $ be given on $ \xdx,$ where $ \xdx $ is a
set of pairs
$ \xBc m,i \xBe,$ $m \xbe M,$ $i$ some index as usual.
(We use here the assumption that the truth value is independent of
indices.)

Let $f:M \xcp V$ be given, we define $ \xbm (f),$ the minimal models of
$f:$
\[ \xbm (f)(m):= \left\{ \begin{array}{lcl}
FALSE \xEH iff \xEH
\xcA \xBc m,i \xBe \xbe \xdx \xcE \xBc m',i' \xBe \xeb \xBc m,i \xBe
.f(m') \xcg f(m) \xEP
\xEH \xEH \xEP
f(m) \xEH \xEH otherwise \xEP
\end{array}
\right.
\]

This generalizes the idea that only models of $ \xbf $ can destroy models
of $ \xbf.$

Obviously, for all $v \xbe V,$ $v \xEd FALSE,$ $\{m: \xbm (f)(m)=v\} \xcc
\{m:f(m)=v\}.$

A structure is called smooth iff for all $f_{ \xbf }$
and for all $ \xBc m,i \xBe $ such that there is
$ \xBc m',i'  \xBe  \xeb  \xBc m,i \xBe $ with $f_{ \xbf }(m') \xcg f_{ \xbf
}(m),$
there is
$ \xBc m'',i''  \xBe  \xeb  \xBc m,i \xBe $ with $f_{ \xbf }(m'') \xcg f_{ \xbf
}(m),$ and
no
$ \xBc n,j \xBe  \xeb  \xBc m'',i''  \xBe $ with $f_{ \xbf }(n) \xcg f_{ \xbf
}(m'').$

A structure will be called definablity preserving iff for all $f_{ \xbf }$
$ \xbm (f_{ \xbf })$ is again the $f_{ \xbq }$ for some $ \xbq.$

\ed

\bd

$\hspace{0.01em}$


\label{Definition Many-MISE}

With these ideas, we can also define minimizing initial segments for
many-valued
structures in a straightforward way:

$F$ is a MISE with respect to $G$ iff

(0) $F \xck G.$

(1) if $F(x) \xEd 0,$ $y \xeb x,$ $G(x) \xck G(y),$ then $F(y)=G(y)$
(downward closure),

and

(2) if $G(x) \xEd 0,$ $F(x)=0,$ then there is $y$ with $y \xeb x,$ $G(x)
\xck G(y),$ $F(y)=G(y).$

\ed

We turn to representation questions.

\be

$\hspace{0.01em}$


\label{Example No-Coher}

This example shows that a suitable choice of truth values can
destroy coherence, as it is present in 2-valued preferential structures.

We want essentially $y \xeb x$ in $ \xCf A,$ $A \xcc B,$ but $y \xeB x$ in
$B.$

The solution will be to make

$F_{B}(y)<F_{B}(x),$ but $F_{A}(y) \xcg F_{A}(x),$ and $F_{A} \xck F_{B},$
e.g., we set:
$F_{B}(x)=3,$ $F_{B}(y)=2,$ $F_{A}(x)=F_{A}(y)=2.$

\ee

This example leads to the following small representation result:

\bfa

$\hspace{0.01em}$


\label{Fact Choice-Rep-Many}

Let $U$ be the universe we work in, let $ \xbm: \xdp (U) \xcp \xdp (U)$
be a function
such that

(1) $ \xbm (X) \xcc X,$

(2) there is no singleton $X=\{x\}$ with $ \xbm (X)= \xCQ.$

Then there is a many-valued preferential structure $ \xdx $ which
represents $ \xbm.$
Note that no coherence conditions are necessary.

\efa

\subparagraph{
Proof
}

$\hspace{0.01em}$


Let 0 be the smallest truth value, and $ \xcA x \xbe U$ $x$ also be a
truth value,
where for all $x \xEd y,$ $x,y \xbe U,$ $x \xcT y$ $(x,y$ as truth values
are incomparable).
Take as preference relation $x \xeb y$ for all $x,y \xbe U,$ $x \xEd y.$

Choose $X \xcc U.$ Define $F_{X}(x):=0$ iff $x \xce \xbm (X),$ and
$F_{X}(x):=x$ iff $x \xbe \xbm (X).$
Then all relations $x \xeb y$ are effective for $y \xce \xbm (X),$ as then
$F_{X}(y) \xck F_{X}(x),$ so
$y$ will not be minimal. If $y \xbe \xbm (X),$ then there is no $x \xEd
y,$ $F_{X}(y) \xck F_{X}(x).$ $ \xcz $
\\[3ex]

Above Fact \ref{Fact Choice-Rep-Many} (page \pageref{Fact Choice-Rep-Many}) 
largely solves the problem of
finding
a preferential representation for arbitrary choice functions by many
valued
structures.

But one might ask different questions in this context, e.g.:
Suppose we have a family of pairs $ \xBc F, \xbm F \xBe $ of functions giving
truth
values
to all $x \xbe U.$ Suppose $ \xcA x \xbe U. \xbm F(x) \xck F(x),$ in short
$ \xbm F \xck F.$
Suppose for simplicity that we have a minimal element 0 of truth values,
with the meaning $F(x)=0$
iff ``$x \xce F$'' (read as a set), so we will not consider $x$ with
$F(x)=0.$ Suppose
further that $ \xbm F(x)=0$ or $ \xbm F(x)=F(x).$ Then, what are the
conditions on the
family of $ \xBc F, \xbm F \xBe $ such that we can represent them by a
many-valued
preferential structure?
The answer is not as trivial as the one to the choice function
representation
problem above.

Consider the following

\be

$\hspace{0.01em}$


\label{Example Coher-F}

Consider $F,G$ with $F \xck G,$ $F(x) \xEd 0,$ $ \xbm F(x)=0,$ $ \xbm F(y)
\xEd 0,$ $F(x) \xck F(y).$
In this case, a relation $y \xeb x$ is effective for $F.$ Suppose now that
also $G(x) \xck G(y),$ then $y \xeb x$ is also effective for $G.$
We may say roughly: If not only $F \xck G,$ but for $x,y$ such that
$F(x),F(y) \xEd 0,$
also for the ``derivatives'' $F' $ and $G' $ $F' \xck G' $ holds in the
sense that
$F(x) \xck F(y)$ $ \xch $ $G(x) \xck G(y),$ then $F$ and $G$ must have the
same coherence
properties as the 2-valued choice functions in order to be
preferentially representable - as any relation effective for $F$ will also
be
effective for $G.$

\ee

This may lead us to consider the following brutal solution:

We have a global truth value relation in the sense that for all $F,G$
$F(x) \xck F(y)$ iff $G(x) \xck G(y)$ - apart from cases where, e.g.,
$F(x)=0,$ as
``$x \xce F$''. In this case, they behave just like normal 2-valued
structures,
but we could now just as well simply omit any relations $x \xeb y$ when
$F(y) \xcK F(x)$ (equivalently, for any other $G).$ So this leads us
nowhere
interesting.

We will leave the problem for further research, and only add a few
rudimentary
remarks:

 \xEh

 \xDH
We may introduce new operators in order to be able to speak about the
situation:

(a) $m \xcm O_{F} \xbf $ iff for all $m' $ such that $m' \xcm \xbf,$
$F(m) \xck F(m'),$

(b) $ \xcm O_{F}(\xbf, \xbq)$ iff for all $m$ such that $m \xcm \xbf $
and all $m' $ such that $m' \xcm \xbq,$
$F(m) \xck F(m').$

These expressions are, of course, still semi-classical, and we can
replace $ \xcm $ by a certain threshold, or consider only $m' $ such that
$F_{ \xbf }(m') \xcg F_{ \xbf }(m).$

Then, given sufficient definability power, we can express that
all models ``in'' $ \xbm F$ have truth values at least as good
as those in $F- \xbm F:=F \xcs \xdC \xbm F,$ $ \xcm O_{F}(F- \xbm F, \xbm
F),$ and use this to formulate
a coherence condition, like:
$O_{F}(F- \xbm F, \xbm F)$ $ \xch $ $O_{G}(F- \xbm F, \xbm F).$

 \xDH

To do some set theory, we will assume that the set of truth values is a
complete Boolean algebra, with symbols $ \xcu $ (or $ \xcU $ for many
arguments) for
infimum, likewise $ \xco $ and $ \xcO $ for supremum, unary - for
complement, binary
$a-b$ for $a \xcu -$b, 0 and 1. For functions $F,G,$ etc with values in
the set
of truth values, we define $ \xcu, \xco $ etc. argumentwise, e.g., $ \xcU
F_{i}$ is defined
by $(\xcU F_{i})(x):= \xcU (F_{i}(x))$ etc.

 \xDH

As an illustration, and for no other purposes, we look at some cases of
the
crucial Fact 3.3.1
in  \cite{Sch04} for representation by smooth structures, which we
repeat now here
for easier reference, together with its proof:

\bfa

$\hspace{0.01em}$


\label{Fact CS-3.3.1}

Let A, $U,$ $U',$ $Y$ and all $A_{i}$ be in $ \xdy.$

$(\xbm \xcc)$ and $(\xbm PR)$ entail:

(1) $A= \xcV \{A_{i}:i \xbe I\}$ $ \xcp $ $ \xbm (A) \xcc \xcV \{ \xbm
(A_{i}):i \xbe I\},$

(2) $U \xcc H(U),$ and $U \xcc U' \xcp H(U) \xcc H(U'),$

(3) $ \xbm (U \xcv Y)-H(U) \xcc \xbm (Y).$

$(\xbm \xcc),$ $(\xbm PR),$ $(\xbm CUM)$ entail:

(4) $U \xcc A,$ $ \xbm (A) \xcc H(U)$ $ \xcp $ $ \xbm (A) \xcc U,$

(5) $ \xbm (Y) \xcc H(U)$ $ \xcp $ $Y \xcc H(U)$ and $ \xbm (U \xcv Y)=
\xbm (U),$

(6) $x \xbe \xbm (U),$ $x \xbe Y- \xbm (Y)$ $ \xcp $ $Y \xcC H(U),$

(7) $Y \xcC H(U)$ $ \xcp $ $ \xbm (U \xcv Y) \xcC H(U).$

\efa

\subparagraph{
Proof:
}

$\hspace{0.01em}$


\label{Section Proof:}

(1) $ \xbm (A) \xcs A_{j} \xcc \xbm (A_{j}) \xcc \xcV \xbm (A_{i}),$ so by
$ \xbm (A) \xcc A= \xcV A_{i}$ $ \xbm (A) \xcc \xcV \xbm (A_{i}).$

(2) trivial.

(3) $ \xbm (U \xcv Y)-H(U)$ $ \xcc_{(2)}$ $ \xbm (U \xcv Y)-U$ $ \xcc_{(
\xbm \xcc)}$ $ \xbm (U \xcv Y) \xcs Y$ $ \xcc_{(\xbm PR)}$ $ \xbm (Y).$

(4) $ \xbm (A)$ $=$ $ \xcV \{ \xbm (A) \xcs X: \xbm (X) \xcc U\}$ $
\xcc_{(\xbm PR')}$ $ \xcV \{ \xbm (A \xcs X): \xbm (X) \xcc U\}.$ But if
$ \xbm (X) \xcc U \xcc A,$ then by
$ \xbm (X) \xcc X,$ $ \xbm (X) \xcc A \xcs X \xcc X$ $ \xcp_{(\xbm CUM)}$
$ \xbm (A \xcs X)= \xbm (X) \xcc U,$ so $ \xbm (A) \xcc U.$

(5) Let $ \xbm (Y) \xcc H(U),$ then by $ \xbm (U) \xcc H(U)$ and (1) $
\xbm (U \xcv Y) \xcc \xbm (U) \xcv \xbm (Y) \xcc H(U),$
so by (4) $ \xbm (U \xcv Y) \xcc U$ and $U \xcv Y \xcc H(U).$ Moreover, $
\xbm (U \xcv Y) \xcc U \xcc U \xcv Y$ $ \xcp_{(\xbm CUM)}$ $ \xbm (U \xcv
Y)= \xbm (U).$

(6) If not, $Y \xcc H(U),$ so $ \xbm (Y) \xcc H(U),$ so $ \xbm (U \xcv Y)=
\xbm (U)$ by (5), but $x \xbe Y- \xbm (Y)$ $ \xcp_{(\xbm PR)}$
$x \xce \xbm (U \xcv Y)= \xbm (U),$ $contradiction.$

(7) $ \xbm (U \xcv Y) \xcc H(U)$ $ \xcp_{(5)}$ $U \xcv Y \xcc H(U).$
$ \xcz $
\\[3ex]

We translate some properties and arguments:

 \xEI

 \xDH
(1)

(a) $F_{i} \xck F$ $ \xch $ $ \xbm F \xcu F_{i} \xck \xbm F_{i}$ for all
$i$ by $(\xbm PR)$
(but recall that $(\xbm PR)$ will not always hold, see
Example \ref{Example No-Coher} (page \pageref{Example No-Coher}))

(b) $ \xbm F \xck F \xck \xcO_{i}F_{i}.$

Thus $ \xbm F=$ (by $b)$ $ \xbm F \xcu \xcO_{i}F_{i}$ $=$ (distributivity)
$ \xcO_{i}(\xbm F \xcu F_{i})$ $ \xck $ (by a) $ \xcO_{i} \xbm F_{i}.$

 \xDH
(3)

We first need an analogue to $X \xcc Y \xcv Z$ $ \xch $ $X-Y \xcc Z:$

(a)
$F_{X} \xck F_{Y} \xco F_{Z}$ $ \xch $ $F_{X}-F_{Y} \xck F_{Z}.$ Proof:
$F_{X}-F_{Y}=F_{X} \xcu \xdC F_{Y}$ $ \xck $ (prerequisite)
$(F_{Y} \xcu \xdC F_{Y}) \xco (F_{Z} \xcu \xdC F_{Y})$ $=$ $0 \xco (F_{Z}
\xcu \xdC F_{Y})$ $ \xck $ $F_{Z}.$

We then need

(b)
$F_{X} \xck F_{X' }$ $ \xch $ $F_{Y}-F_{X' }$ $ \xck $ $F_{Y}-F_{X}.$
Proof: $F_{X} \xck F_{X' }$ $ \xch $ $ \xdC F_{X' } \xck \xdC F_{X},$ so
$F_{Y} \xcs \xdC F_{X' } \xck F_{Y} \xcs \xdC F_{X}.$

Thus, $ \xbm f_{U \xcv Y}-f_{H(U)}$ $ \xck $ (by (2) and (b)) $ \xbm f_{U
\xcv Y}-f_{U}$ $ \xck $ (by $(\xbm \xcc),$ (a))
$ \xbm f_{U \xcv Y} \xcu f_{Y}$ $ \xck $ $ \xbm f_{Y}$ by $(\xbm PR).$

 \xDH
(6)

$F_{Y} \xck F_{H(U)}$ $ \xch $ $ \xbm F_{Y} \xck F_{H(U)}$ $ \xch $ $ \xbm
F_{U \xcv Y}$ $=$ $ \xbm (F_{U} \xco F_{Y})$ $=$ $ \xbm F_{U}$ by (5).
$ \xbm F_{U \xcv Y} \xcu F_{Y}$ $=$ $ \xbm (F_{U} \xco F_{V}) \xcu F_{Y}$
$ \xck $ $ \xbm F_{Y}$ by $(\xbm PR),$ so
$ \xbm F_{U} \xcu F_{Y}= \xbm F_{U \xcv Y} \xcu F_{Y} \xck \xbm F_{Y}.$
Thus $ \xbm F_{U} \xcu F_{Y} \xcu \xdC \xbm F_{Y}= \xCQ,$
contradicting the prerequisite.

 \xEJ

 \xEj
\section{
IBRS and higher preferential structures
}
\subsection{
General IBRS
}

We first define IBRS:

$ \xCO $
\index{Definition IBRS}

\bd

$\hspace{0.01em}$


\label{Definition IBRS}

 \xEh
 \xDH
An $ \xCf information$ $ \xCf bearing$ $ \xCf binary$ $ \xCf relation$ $
\xCf frame$ $ \xCf IBR,$ has the form
 \index{IBR}
 \index{points}
 \index{arrows}
 \index{nodes}
$(S, \xdR),$ where $S$ is a non-empty set and $ \xdR $ is a subset of
$S_{ \xbo },$ where $S_{ \xbo }$ is defined by induction as follows:

 \xEh

 \xDH $S_{0}:=S$

 \xDH $S_{n+1}:=$ $S_{n} \xcv (S_{n} \xCK S_{n}).$

 \xDH $S_{ \xbo }=$ $ \xcV \{S_{n}:n \xbe \xbo \}$

 \xEj

We call elements from $S$ $ \xCf points$ or $ \xCf nodes,$ and elements
from $ \xdR $ $ \xCf arrows.$
Given $(S, \xdR),$ we also set $ \xdP ((S, \xdR)):=S,$ and $ \xdA ((S,
\xdR)):= \xdR.$

If $ \xba $ is an arrow, the origin and destination of $ \xba $ are
defined
as usual, and we write $ \xba:x \xcp y$ when $x$ is the origin, and $y$
the destination of the arrow $ \xba.$ We also write $o(\xba)$ and $d(
\xba)$ for
the origin and destination of $ \xba.$

 \xDH Let $Q$ be a set of atoms, and $ \xdL $ be a set of labels (usually
$\{0,1\}$ or $[0,1]).$ An $ \xCf information$ $ \xCf assignment$ $h$ on
$(S, \xdR)$
is a function $h:Q \xCK \xdR \xcp \xdL.$

 \xDH An $ \xCf information$ $ \xCf bearing$ $ \xCf system$ $ \xCf IBRS,$
has the form
$(S, \xdR,h,Q, \xdL),$ where $S,$ $ \xdR,$ $h,$ $Q,$ $ \xdL $ are as
above.

 \xEj

See Diagram \ref{Diagram IBRS} (page \pageref{Diagram IBRS})  for an
illustration.

\ed

$ \xCO $

\vspace{10mm}

\begin{diagram}

\centering
\setlength{\unitlength}{0.0006in}
{\renewcommand{\dashlinestretch}{30}
\begin{picture}(4961,5004)(0,0)
\path(1511,1583)(611,3683)
\blacken\path(685.845,3584.520)(611.000,3683.000)(630.696,3560.885)(672.451,3539.613)(685.845,3584.520)
\path(1511,1583)(2411,3683)
\blacken\path(2391.304,3560.885)(2411.000,3683.000)(2336.155,3584.520)(2349.549,3539.613)(2391.304,3560.885)
\path(3311,1583)(4361,4133)
\blacken\path(4343.050,4010.616)(4361.000,4133.000)(4287.570,4033.461)(4301.603,3988.750)(4343.050,4010.616)
\path(3316,1574)(2416,3674)
\blacken\path(2490.845,3575.520)(2416.000,3674.000)(2435.696,3551.885)(2477.451,3530.613)(2490.845,3575.520)
\path(986,2783)(2621,2783)
\blacken\path(2501.000,2753.000)(2621.000,2783.000)(2501.000,2813.000)(2465.000,2783.000)(2501.000,2753.000)
\path(2486,2783)(2786,2783)
\blacken\path(2666.000,2753.000)(2786.000,2783.000)(2666.000,2813.000)(2630.000,2783.000)(2666.000,2753.000)
\path(3311,1583)(2051,2368)
\blacken\path(2168.714,2330.008)(2051.000,2368.000)(2136.987,2279.083)(2183.406,2285.509)(2168.714,2330.008)
\path(2166,2288)(1906,2458)
\blacken\path(2022.854,2417.439)(1906.000,2458.000)(1990.019,2367.221)(2036.567,2372.629)(2022.854,2417.439)

\put(1511,1358) {{\xssc $a$}}
\put(3311,1358) {{\xssc $d$}}
\put(3311,1058)  {{\xssc $(p,q)=(1,0)$}}
\put(1511,1058)  {{\xssc $(p,q)=(0,0)$}}
\put(2411,3758){{\xssc $c$}}
\put(4361,4433){{\xssc $(p,q)=(1,1)$}}
\put(4361,4208){{\xssc $e$}}
\put(2411,3983){{\xssc $(p,q)=(0,1)$}}
\put(611,3983) {{\xssc $(p,q)=(0,1)$}}
\put(611,3758) {{\xssc $b$}}
\put(1211,2883){{\xssc $(p,q)=(1,1)$}}
\put(260,2333) {{\xssc $(p,q)=(1,1)$}}
\put(2261,1583) {{\xssc $(p,q)=(1,1)$}}
\put(1286,3233){{\xssc $(p,q)=(1,1)$}}
\put(2711,3083){{\xssc $(p,q)=(1,1)$}}
\put(3836,2633){{\xssc $(p,q)=(1,1)$}}

\put(300,700)
{{\rm\bf
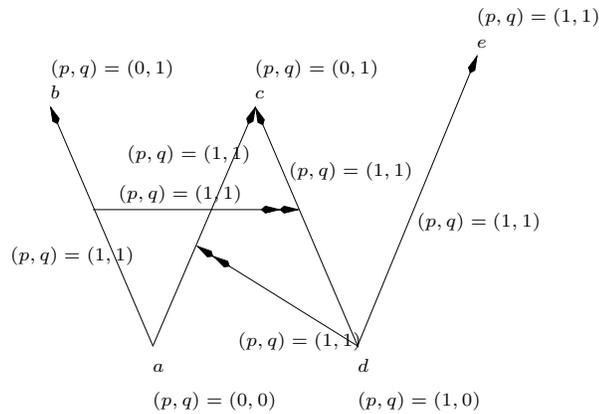
A simple example of an information bearing system.
}}

\end{picture}
}

\label{Diagram IBRS}
\index{Diagram IBRS}

\end{diagram}

We have here:
\[\begin{array}{l}
S =\{a,b,c,d,e\}.\\
\xdR = S \cup \{(a,b), (a,c), (d,c), (d,e)\} \cup \{((a,b), (d,c)),
(d,(a,c))\}.\\
Q = \{p,q\}
\end{array}
\]
The values of $h$ for $p$ and $q$ are as indicated in the figure. For
example $h(p,(d,(a,c))) =1$.

\vspace{4mm}

$ \xCO $

$ \xCO $

$ \xCO $

\bcom

$\hspace{0.01em}$


\label{Comment}

\label{Comment IBRS}

The elements in Figure Diagram \ref{Diagram IBRS} (page \pageref{Diagram IBRS}) 
can be interpreted in
many ways,
depending on the area of application.

 \xEh

 \xDH The points in $S$ can be interpreted as possible worlds, or
as nodes in an argumentation network or nodes in a neural net or
states, etc.

 \xDH The direct arrows from nodes to nodes can be interpreted as
accessibility relation, attack or support arrows in an argumentation
networks, connection in a neural nets, a preferential ordering in
a nonmonotonic model, etc.

 \xDH The labels on the nodes and arrows can be interpreted as fuzzy
values in the accessibility relation or weights in the neural net
or strength of arguments and their attack in argumentation nets, or
distances in a counterfactual model, etc.

 \xDH The double arrows can be interpreted as feedback loops to nodes
or to connections, or as reactive links changing the system which are
activated
as we pass between the nodes.

 \xEj

\ecom

$ \xCO $
\subsection{
Higher preferential structures
}

\label{Section High-Pref}
\subsubsection{
Introduction
}

We turn to the special case of higher preferential structures,
give the definitions and some results.
\subsubsection{
Definitions and facts for basic structures
}

$ \xCO $
\index{Definition Generalized preferential structure}

\bd

$\hspace{0.01em}$


\label{Definition Generalized preferential structure}

An IBR is called a $ \xCf generalized$ $ \xCf preferential$ $ \xCf
structure$ iff the origins of all
 \index{generalized preferential structure}
arrows are points. We will usually write $x,y$ etc. for points, $ \xba,$
$ \xbb $ etc. for
arrows.

\ed

$ \xCO $

$ \xCO $
\index{Definition Level-n-Arrow}

\bd

$\hspace{0.01em}$


\label{Definition Level-n-Arrow}

Consider a generalized preferential structure $ \xdx.$

(1) $ \xCf Level$ $n$ $ \xCf arrow:$

 \index{level $n$ arrow}
Definition by upward induction.

If $ \xba:x \xcp y,$ $x,y$ are points, then $ \xba $ is a level 1 arrow.

If $ \xba:x \xcp \xbb,$ $x$ is a point, $ \xbb $ a level $n$ arrow, then
$ \xba $ is a level $n+1$ arrow.
($o(\xba)$ \index{$o(\xba)$}  is the origin,
$d(\xba)$ \index{$d(\xba)$}  is the destination of $ \xba.)$

$ \xbl (\xba)$ \index{$ \xbl (\xba)$}  will denote the level of $ \xba.$

(2) $ \xCf Level$ $n$ $ \xCf structure:$

 \index{level $n$ structure}
$ \xdx $ is a level $n$ structure iff all arrows in $ \xdx $ are at most
level $n$ arrows.

We consider here only structures of some arbitrary but finite level $n.$

(3) We define for an arrow $ \xba $ by induction $O(\xba)$ and $D(\xba
).$

If $ \xbl (\xba)=1,$ then $O(\xba)$ \index{$O(\xba)$}  $:=\{o(\xba
)\},$
$D(\xba)$ \index{$D(\xba)$}  $:=\{d(\xba)\}.$

If $ \xba:x \xcp \xbb,$ then $D(\xba):=D(\xbb),$ and $O(\xba
):=\{x\} \xcv O(\xbb).$

Thus, for example, if $ \xba:x \xcp y,$ $ \xbb:z \xcp \xba,$ then $O(
\xbb):=\{x,z\},$ $D(\xbb)=\{y\}.$
Consider also the arrow $ \xbb:= \xBc  \xbb',l'  \xBe $ in
Diagram \ref{Diagram Essential-Smooth-3-1-2} (page \pageref{Diagram
Essential-Smooth-3-1-2}). There,
$D(\xbb)=\{ \xBc x,i \xBe \},$
$O(\xbb)=\{ \xBc z',m'  \xBe, \xBc y,j \xBe \}.$

\ed

$ \xCO $

$ \xCO $
\index{Example O-D}

\be

$\hspace{0.01em}$


\label{Example O-D}

Let $ \xbl (\xba)$ be finite. $ \xba $ may be of the form $ \xba:x \xcp
\xbb,$ $ \xbb:y \xcp \xbg,$
$ \xbg:z \xcp \xbh,$ $ \xbh:u \xcp v.$ $D(\xba)=\{v\},$ the last
destination in the construction,
so always a point. $O(\xba)$ is $\{x,y,z,u\},$ the set of all origins of
these
arrows, a set of points. We do $ \xCf not$ consider any other arrows
pointing
to or going from elements of this construction.

See Diagram \ref{Diagram O-D} (page \pageref{Diagram O-D}).

\ee

$ \xCO $

\vspace{10mm}

\begin{diagram}

\label{Diagram O-D}
\index{Diagram O-D}

\centering
\setlength{\unitlength}{1mm}
{\renewcommand{\dashlinestretch}{30}
\begin{picture}(150,150)(0,0)

\path(30,20)(30,52)
\path(29,49.2)(30,52)(31,49.2)

\path(14,53)(46,53)
\path(43.2,54)(46,53)(43.2,52)

\path(47,37)(47,69)
\path(46,66.2)(47,69)(48,66.2)

\path(31,70)(63,70)
\path(60.2,71)(63,70)(60.2,69)

\path(64,54)(64,86)
\path(63,83.2)(64,86)(65,83.2)

\put(29,18){{\xssc $w$}}
\put(31,35){{\xssc $\xbr$}}

\put(12,52){{\xssc $x$}}
\put(29,54){{\xssc $\xba$}}

\put(46,35){{\xssc $y$}}
\put(48,52){{\xssc $\xbb$}}

\put(29,69){{\xssc $z$}}
\put(46,72){{\xssc $\xbg$}}

\put(63,52){{\xssc $u$}}
\put(65,69){{\xssc $\xbh$}}

\put(63,87){{\xssc $v$}}

\put(30,7){{\xssc $D(x)=\{v\}$, $O(x)=\{x,y,z,u\}$}}

\end{picture}
}

\end{diagram}

\vspace{4mm}

$ \xCO $

$ \xCO $

$ \xCO $

\bcom

$\hspace{0.01em}$


\label{Comment Gen-Pref}

A counterargument to $ \xba $ is NOT an argument for $ \xCN \xba $ (this
is asking for too
much), but just showing one case where $ \xCN \xba $ holds. In
preferential structures,
an argument for $ \xba $ is a set of level 1 arrows, eliminating $ \xCN
\xba -$models. A
counterargument is one level 2 arrow, attacking one such level 1 arrow.

Of course, when we have copies, we may need many successful attacks, on
all
copies, to achieve the goal. As we may have copies of level 1 arrows, we
may need many level 2 arrows to destroy them all.

\ecom

$ \xCO $

$ \xCO $

We will not consider here diagrams with arbitrarily high levels. One
reason is
that diagrams like the following will have an unclear meaning:

\be

$\hspace{0.01em}$


\label{Example Inf-Level}

$ \xBc  \xba,1 \xBe:x \xcp y,$

$ \xBc  \xba,n+1 \xBe:x \xcp  \xBc  \xba,n \xBe $ $(n \xbe \xbo).$

Is $y \xbe \xbm (X)?$

\ee

$ \xCO $

$ \xCO $
\index{Definition Valid-Arrow}

\bd

$\hspace{0.01em}$


\label{Definition Valid-Arrow}

Let $ \xdx $ be a generalized preferential structure of (finite) level
$n.$

We define (by downward induction):

(1) $ \xCf Valid$ $X-to-Y$ $ \xCf arrow:$

 \index{valid $X-to-Y$ arrow}

Let $X,Y \xcc \xdP (\xdx).$

$ \xba \xbe \xdA (\xdx)$ is a $ \xCf valid$ $X-to-Y$ $ \xCf arrow$ iff

(1.1) $O(\xba) \xcc X,$ $D(\xba) \xcc Y,$

(1.2) $ \xcA \xbb:x' \xcp \xba.(x' \xbe X$ $ \xch $ $ \xcE \xbg:x''
\xcp \xbb.(\xbg $ is a valid $X-to-Y$ arrow)).

We will also say that $ \xba $ is a $ \xCf valid$ $ \xCf arrow$ in $X,$ or
just $ \xCf valid$ in $X,$ iff $ \xba $ is a
valid $X-to-X$ arrow.

(2) $ \xCf Valid$ $X \xch Y$ $ \xCf arrow:$
 \index{valid $X \xch Y$ arrow}

Let $X \xcc Y \xcc \xdP (\xdx).$

$ \xba \xbe \xdA (\xdx)$ is a $ \xCf valid$ $X \xch Y$ $ \xCf arrow$ iff

(2.1) $o(\xba) \xbe X,$ $O(\xba) \xcc Y,$ $D(\xba) \xcc Y,$

(2.2) $ \xcA \xbb:x' \xcp \xba.(x' \xbe Y$ $ \xch $ $ \xcE \xbg:x''
\xcp \xbb.(\xbg $ is a valid $X \xch Y$ arrow)).

Thus, any attack $ \xbb $ from $Y$ against $ \xba $ has to be countered by
a valid attack
on $ \xbb.$

(Note that in particular $o(\xbg) \xbe X,$ and that $o(\xbb)$ need not
be in $X,$ but
can be in the bigger $Y.)$

\ed

\br

$\hspace{0.01em}$


\label{Remark Valid-Valid}

Note that, in the definition of valid $X-to-Y$ arrow, $X$ and $Y$ need not
be related, but in the definition of valid $X \xch Y$ arrow, $X \xcc Y.$

Let us assume now that $X \xcc Y,$ and look at the remaining differences.

In both cases, $D(\xba) \xcc Y.$

In the $X-to-Y$ case, $O(\xba) \xcc X,$ and attacks from $X$ are
countered.

In the $X \xch Y$ case, $o(\xba) \xbe X,$ $O(\xba) \xcc Y,$ and
attacks from $Y$ are countered.

So the first condition is stronger in the $X-to-Y$ case, the second in the
$X \xch Y$ case.

\er

\be

$\hspace{0.01em}$


\label{Example Valid-Arrow}

(1)
Consider the arrow $ \xbb:= \xBc  \xbb',l'  \xBe $ in
Diagram \ref{Diagram Essential-Smooth-3-1-2} (page \pageref{Diagram
Essential-Smooth-3-1-2}).
$D(\xbb)=\{ \xBc x,i \xBe \},$
$O(\xbb)=\{ \xBc z',m'  \xBe, \xBc y,j \xBe \},$
and the only arrow attacking $ \xbb $ originates outside $X,$
so $ \xbb $ is a valid $X-to- \xbm (X)$ arrow.

(2)
Consider the arrows $ \xBc  \xba',k'  \xBe $ and $ \xBc  \xbg',n'  \xBe $ in
Diagram \ref{Diagram Essential-Smooth-3-2} (page \pageref{Diagram
Essential-Smooth-3-2}).
Both are valid $ \xbm (X) \xch X$ arrows.

\ee

\be

$\hspace{0.01em}$


\label{Example Valid-Valid}

See Diagram \ref{Diagram Valid-Valid} (page \pageref{Diagram Valid-Valid}).

Let $X \xcc Y.$

Consider the left hand part of the diagram.

The fact that $ \xbd $ originates
in $Y,$ but not in $X$ makes that $ \xba $ is not a valid $X-to-Y$ arrow,
as the
condition $O(\xba) \xcc X$ is violated. To be a valid $X \xch Y$ arrow,
we have to show
that all attacks on $ \xba $ originating from $Y$ (not only from $X)$ are
be countered by valid $X \xch Y$ arrows. This holds, as $ \xbb_{1}$ is
countered by $ \xbg_{1},$
$ \xbb_{2}$ by $ \xbg_{2}.$ All possible attacks on $ \xba,$ $ \xbg_{1},$
$ \xbg_{2}$ from outside $Y,$ like $ \xbr,$
need not be considered.

Consider the right hand part of the diagram.

The fact that there is no valid counterargument to $ \xbb'_{2}$ makes
that $ \xba' $
is not a valid $X \xch Y$ arrow. It is a valid $X-to-Y$ arrow, as
counterarguments
to $ \xba' $ like $ \xbb'_{2},$ which do not originate in $X,$ are not
considered. The
counterargument $ \xbb'_{1}$ is considered, but it is destroyed by the
valid
$X-to-Y$ arrow $ \xbg'_{1}.$

\ee

$ \xCO $

\vspace{10mm}

\begin{diagram}

\label{Diagram Valid-Valid}
\index{Diagram Valid-Valid}

\centering
\setlength{\unitlength}{1mm}
{\renewcommand{\dashlinestretch}{30}
\begin{picture}(150,150)(0,0)

\path(10,145)(130,145)(130,5)(10,5)(10,145)
\path(10,75)(130,75)
\put(135,140){$Y$}
\put(135,70){$X$}

\put(20,135){$X \xch Y$}
\put(95,135){$X-to-Y$}


\path(25,110)(55,110)
\put(39,111){{\xssc $\xbd$}}
\path(52.2,111)(55,110)(52.2,109)

\path(40,20)(40,109)
\put(41,90){{\xssc $\xba$}}
\path(39,106.2)(40,109)(41,106.2)

\path(20,85)(39,85)
\put(29,86){{\xssc $\xbb_2$}}
\path(36.4,85.9)(39,85)(36.4,84.1)

\path(30,65)(30,84)
\put(26,79){{\xssc $\xbg_2$}}
\path(29.1,81.4)(30,84)(30.9,81.4)

\path(60,55)(41,55)
\put(49,56){{\xssc $\xbb_1$}}
\path(43.6,54.1)(41,55)(43.6,55.9)

\path(50,34)(50,54)
\put(51,44){{\xssc $\xbg_1$}}
\path(49,51.3)(50,54)(51,51.3)

\path(5,55)(39,55)
\put(20,53){{\xssc $\xbr$}}
\path(36.2,56)(39,55)(36.2,54)


\path(110,20)(110,109)
\put(111,90){{\xssc $\xbd'$}}
\path(109,106.2)(110,109)(111,106.2)

\path(90,55)(109,55)
\put(91,52){{\xssc $\xba'$}}
\path(106.4,55.9)(109,55)(106.4,54.1)

\path(100,85)(100,56)
\put(101,80){{\xssc $\xbb'_2$}}
\path(101,58.8)(100,56)(99,58.8)

\path(100,35)(100,54)
\put(101,45){{\xssc $\xbb'_1$}}
\path(99.1,51.4)(100,54)(100.9,51.4)

\path(80,45)(99,45)
\put(89,42){{\xssc $\xbg'_1$}}
\path(96.4,45.9)(99,45)(96.4,44.1)

\end{picture}
}

\end{diagram}

\vspace{4mm}

$ \xCO $

$ \xCO $

$ \xCO $

\bfa

$\hspace{0.01em}$


\label{Fact Higher-Validity}

(1) If $ \xba $ is a valid $X \xch Y$ arrow, then $ \xba $ is a valid
$Y-to-Y$ arrow.

(2) If $X \xcc X' \xcc Y' \xcc Y \xcc \xdP (\xdx)$ and $ \xba \xbe \xdA
(\xdx)$ is a valid $X \xch Y$ arrow, and
$O(\xba) \xcc Y',$ $D(\xba) \xcc Y',$ then $ \xba $ is a valid $X'
\xch Y' $ arrow.

\efa

$ \xCO $

$ \xCO $

\subparagraph{
Proof
}

$\hspace{0.01em}$


Let $ \xba $ be a valid $X \xch Y$ arrow. We show (1) and (2) together by
downward induction
(both are trivial).

By prerequisite $o(\xba) \xbe X \xcc X',$ $O(\xba) \xcc Y' \xcc Y,$
$D(\xba) \xcc Y' \xcc Y.$

Case 1: $ \xbl (\xba)=n.$ So $ \xba $ is a valid $X' \xch Y' $ arrow,
and a valid $Y-to-Y$ arrow.

Case 2: $ \xbl (\xba)=n-1.$ So there is no $ \xbb:x' \xcp \xba,$ $y
\xbe Y,$ so $ \xba $ is a valid
$Y-to-Y$ arrow. By $Y' \xcc Y$ $ \xba $ is a valid $X' \xch Y' $ arrow.

Case 3: Let the result be shown down to $m,$ $n>m>1,$ let $ \xbl (\xba
)=m-1.$
So $ \xcA \xbb:x' \xcp \xba (x' \xbe Y$ $ \xch $ $ \xcE \xbg:x'' \xcp
\xbb (x'' \xbe X$ and $ \xbg $ is a valid $X \xch Y$ arrow)).
By induction hypothesis $ \xbg $ is a valid $Y-to-Y$ arrow, and a valid
$X' \xch Y' $ arrow. So $ \xba $ is a valid $Y-to-Y$ arrow, and by $Y'
\xcc Y,$ $ \xba $ is a valid
$X' \xch Y' $ arrow.

$ \xcz $
\\[3ex]

$ \xCO $

$ \xCO $

\bd

$\hspace{0.01em}$


\label{Definition Higher-Mu}

Let $ \xdx $ be a generalized preferential structure of level $n,$ $X \xcc
\xdP (\xdx).$

$ \xbm (X)$ \index{$ \xbm (X)$}  $:=\{x \xbe X:$ $ \xcE  \xBc x,i \xBe. \xCN
\xcE $
valid $X-to-X$ arrow $ \xba:x' \xcp  \xBc x,i \xBe \}.$

\ed

$ \xCO $
\subsubsection{
Definitions and facts for totally and essentially smooth structures
}

$ \xCO $

\bcom

$\hspace{0.01em}$


\label{Comment Smooth-Gen}

The purpose of smoothness \index{smoothness}  is to guarantee cumulativity.
Smoothness
achieves Cumulativity by mirroring all information present in $X$ also in
$ \xbm (X).$
Closer inspection shows that smoothness does more than necessary. This is
visible when there are copies (or, equivalently, non-injective labelling
functions). Suppose we have two copies of $x \xbe X,$
$ \xBc x,i \xBe $ and $ \xBc x,i'  \xBe,$ and there
is $y \xbe X,$ $ \xba: \xBc y,j \xBe  \xcp  \xBc x,i \xBe,$ but there is no
$ \xba': \xBc y',j'  \xBe  \xcp  \xBc x,i'  \xBe,$ $y' \xbe X.$ Then
$ \xba: \xBc y,j \xBe  \xcp  \xBc x,i \xBe $
is irrelevant, as $x \xbe \xbm (X)$ anyhow. So mirroring
$ \xba: \xBc y,j \xBe  \xcp  \xBc x,i \xBe $ in $ \xbm (X)$ is not
necessary, i.e., it is not necessary to have some
$ \xba': \xBc y',j'  \xBe  \xcp  \xBc x,i \xBe,$ $y' \xbe \xbm (X).$

On the other hand, Example \ref{Example Need-Smooth} (page \pageref{Example
Need-Smooth})  shows that,
if we want smooth structures to correspond to the property $(\xbm CUM),$
we
need at least some valid arrows from $ \xbm (X)$ also for higher level
arrows.
This ``some'' is made precise (essentially) in Definition \ref{Definition
X-Sub-X'} (page \pageref{Definition X-Sub-X'}).

From a more philosophical point of view,
when we see the (inverted) arrows of preferential structures as attacks on
non-minimal elements, then we should see smooth structures as always
having
attacks also from valid (minimal) elements. So, in general structures,
also
attacks from non-valid elements are valid; in smooth structures we always
also have attacks from valid elements.

The analogue to usual smooth structures, on level 2, is then that any
successfully attacked level 1 arrow is also attacked from a minimal point.

\ecom

$ \xCO $

$ \xCO $

\bd

$\hspace{0.01em}$


\label{Definition X-Sub-X'}

Let $ \xdx $ be a generalized preferential structure.

$X \xes X' $ \index{$X \xes X' $}  iff

(1) $X \xcc X' \xcc \xdP (\xdx),$

(2) $ \xcA x \xbe X' -X$ $ \xcA  \xBc x,i \xBe $
$ \xcE \xba:x' \xcp  \xBc x,i \xBe (\xba $ is a valid $X \xch X' $ arrow),

(3) $ \xcA x \xbe X$ $ \xcE  \xBc x,i \xBe $

$ \xDC $ $(\xcA \xba:x' \xcp  \xBc x,i \xBe (x' \xbe X' $ $ \xch $ $ \xcE \xbb
:x'' \xcp \xba.(\xbb $ is a valid $X \xch X' $ arrow))).

Note that (3) is not simply the negation of (2):

Consider a level 1 structure. Thus all level 1 arrows are valid, but the
source
of the arrows must not be neglected.

(2) reads now: $ \xcA x \xbe X' -X$
$ \xcA  \xBc x,i \xBe $ $ \xcE \xba:x' \xcp  \xBc x,i \xBe.x' \xbe X$

(3) reads: $ \xcA x \xbe X$ $ \xcE  \xBc x,i \xBe $
$ \xCN \xcE \xba:x' \xcp  \xBc x,i \xBe.x' \xbe X' $

This is intended: intuitively, read $X= \xbm (X'),$ and minimal elements
must not be
attacked at all, but non-minimals must be attacked from $X$ - which is a
modified
version of smoothness. More precisely: non-minimal elements (i.e., from
$X' -$X)
have to be validly attacked from $X,$ minimal elements must not be validly
attacked at all from $X' $ (only perhaps from the outside).

\ed

$ \xCO $

$ \xCO $

\br

$\hspace{0.01em}$


\label{Remark X-Sub-X'}

We note the special case of Definition \ref{Definition X-Sub-X'} (page
\pageref{Definition X-Sub-X'})  for level
3 structures.
We also write it immediately for the intended case
$ \xbm (X) \xes X,$ and explicitly with copies.

$x \xbe \xbm (X)$ iff

(1) $ \xcE  \xBc x,i \xBe  \xcA  \xBc  \xba,k \xBe: \xBc y,j \xBe  \xcp  \xBc
x,i \xBe $

$ \xDC (y \xbe X$ $ \xch $
$ \xcE  \xBc  \xbb',l'  \xBe: \xBc z',m'  \xBe  \xcp  \xBc  \xba,k \xBe.$

$ \xDC \xDC (z' \xbe \xbm (X)$ $ \xcu $
$ \xCN \xcE  \xBc  \xbg',n'  \xBe: \xBc u',p'  \xBe  \xcp  \xBc  \xbb',l' 
\xBe.u' \xbe X))$

See Diagram \ref{Diagram Essential-Smooth-3-1-2} (page \pageref{Diagram
Essential-Smooth-3-1-2}).

$x \xbe X- \xbm (X)$ iff

(2) $ \xcA  \xBc x,i \xBe  \xcE  \xBc  \xba',k'  \xBe: \xBc y',j'  \xBe  \xcp 
\xBc x,i \xBe $

$ \xDC (y' \xbe \xbm (X)$ $ \xcu $

$ \xDC \xDC (a)$ $ \xCN \xcE  \xBc  \xbb',l'  \xBe: \xBc z',m'  \xBe  \xcp
\xBc \xba',k'
\xBe.z' \xbe X$

$ \xDC \xDC or$

$ \xDC \xDC (b)$ $ \xcA  \xBc  \xbb',l'  \xBe: \xBc z',m'  \xBe  \xcp \xBc
\xba',k' \xBe
$

$ \xDC \xDC \xDC (z' \xbe X$ $ \xch $
$ \xcE  \xBc  \xbg',n'  \xBe: \xBc u',p'  \xBe  \xcp  \xBc  \xbb',l'  \xBe
.u' \xbe \xbm (X))$)

See Diagram \ref{Diagram Essential-Smooth-3-2} (page \pageref{Diagram
Essential-Smooth-3-2}).

\er

$ \xCO $

$ \xCO $

\vspace{15mm}

\begin{diagram}

\label{Diagram Essential-Smooth-3-1-2}
\index{Diagram Essential smoothness I}

\centering
\setlength{\unitlength}{.6mm}
{\renewcommand{\dashlinestretch}{30}
\begin{picture}(100,130)(0,0)
\put(50,80){\circle{80}}
\path(10,80)(90,80)

\put(100,100){{\xssc $X$}}
\put(100,60){{\xssc $\xbm (X)$}}

\path(50,100)(50,60)
\path(48.5,63)(50,60)(51.5,63)
\put(50,101){\circle*{0.3}}
\put(50,59){\circle*{0.3}}
\put(50,57){{\xssc $ \xBc x,i \xBe $}}
\put(50,102){{\xssc $ \xBc y,j \xBe $}}
\put(52,90){{\xssc $ \xBc \xba,k \xBe $}}

\path(20,70)(48,70)
\path(46,71)(48,70)(46,69)
\put(19,70){\circle*{0.3}}
\put(13,67){{\xssc $ \xBc z',m' \xBe $}}
\put(26,71){{\xssc $ \xBc \xbb',l' \xBe $}}

\path(30,30)(30,68)
\path(29,66)(30,68)(31,66)
\put(30,29){\circle*{0.3}}
\put(31,53){{\xssc $ \xBc \xbg',n' \xBe $}}
\put(30,26.5){{\xssc $ \xBc u',p' \xBe $}}

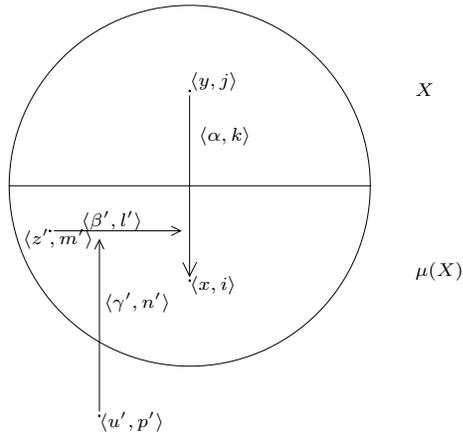
\put(30,10) {{\rm\bf Case 3-1-2}}

\end{picture}
}
\end{diagram}

\vspace{4mm}

$ \xCO $

$ \xCO $

\vspace{10mm}

\begin{diagram}

\label{Diagram Essential-Smooth-3-2}
\index{Diagram Essential smoothness II}

\centering
\setlength{\unitlength}{.6mm}
{\renewcommand{\dashlinestretch}{30}
\begin{picture}(100,110)(0,0)
\put(50,60){\circle{80}}
\path(10,60)(90,60)

\put(100,80){{\xssc $X$}}
\put(100,40){{\xssc $\xbm (X)$}}

\path(50,80)(50,40)
\path(48.5,77)(50,80)(51.5,77)
\put(50,81){\circle*{0.3}}
\put(50,39){\circle*{0.3}}
\put(50,37){{\xssc $ \xBc y',j' \xBe $}}
\put(50,82){{\xssc $ \xBc x,i \xBe $}}
\put(52,70){{\xssc $ \xBc \xba',k' \xBe $}}

\path(20,70)(48,70)
\path(46,71)(48,70)(46,69)
\put(19,70){\circle*{0.3}}
\put(13,67){{\xssc $ \xBc z',m' \xBe $}}
\put(26,71){{\xssc $ \xBc \xbb',l' \xBe $}}

\path(30,30)(30,68)
\path(29,66)(30,68)(31,66)
\put(30,29){\circle*{0.3}}
\put(31,53){{\xssc $ \xBc \xbg',n' \xBe $}}
\put(31,26.5){{\xssc $ \xBc u',p' \xBe $}}

\put(30,0) {{\rm\bf Case 3-2}}

\end{picture}
}
\end{diagram}

\vspace{4mm}

$ \xCO $

$ \xCO $

\bfa

$\hspace{0.01em}$


\label{Fact X-Sub-X'}

(1) If $X \xes X',$ then $X= \xbm (X'),$

(2) $X \xes X',$ $X \xcc X'' \xcc X' $ $ \xch $ $X \xes X''.$ (This
corresponds
to $(\xbm CUM)$ \index{$(\xbm CUM)$}.)

(3) $X \xes X',$ $X \xcc Y',$ $Y \xes Y',$ $Y \xcc X' $ $ \xch $ $X=Y.$
(This corresponds
to $(\xbm \xcc \xcd)$ \index{$(\xbm \xcc \xcd)$}.)

\efa

$ \xCO $

$ \xCO $

\subparagraph{
Proof
}

$\hspace{0.01em}$


(1) Trivial by Fact \ref{Fact Higher-Validity} (page \pageref{Fact
Higher-Validity})  (1).

(2)

We have to show

(a) $ \xcA x \xbe X'' -X$ $ \xcA  \xBc x,i \xBe $
$ \xcE \xba:x' \xcp  \xBc x,i \xBe (\xba $ is a valid $X \xch X'' $ arrow), and

(b) $ \xcA x \xbe X$ $ \xcE  \xBc x,i \xBe $
$(\xcA \xba:x' \xcp  \xBc x,i \xBe (x' \xbe X'' $ $ \xch $ $ \xcE \xbb:x'' \xcp
\xba.(\xbb $ is a valid $X \xch X'' $ arrow))).

Both follow from the corresponding condition for $X \xch X',$ the
restriction of the
universal quantifier, and Fact \ref{Fact Higher-Validity} (page \pageref{Fact
Higher-Validity})  (2).

(3)

Let $x \xbe X-$Y.

(a) By $x \xbe X \xes X',$ $ \xcE  \xBc x,i \xBe $ s.t.
$(\xcA \xba:x' \xcp  \xBc x,i \xBe (x' \xbe X' $ $ \xch $ $ \xcE \xbb:x'' \xcp
\xba.(\xbb $ is a valid $X \xch X' $ arrow))).

(b) By $x \xce Y \xes \xcE \xba_{1}:x' \xcp  \xBc x,i \xBe $ $ \xba_{1}$ is a
valid
$Y \xch Y' $ arrow, in
particular $x' \xbe Y \xcc X'.$ Moreover, $ \xbl (\xba_{1})=1.$

So by (a) $ \xcE \xbb_{2}:x'' \xcp \xba_{1}.(\xbb_{2}$ is a valid $X \xch
X' $ arrow), in particular $x'' \xbe X \xcc Y',$
moreover $ \xbl (\xbb_{2})=2.$

It follows by induction from the definition of valid $A \xch B$ arrows
that

$ \xcA n \xcE \xba_{2m+1},$ $ \xbl (\xba_{2m+1})=2m+1,$ $ \xba_{2m+1}$ a
valid $Y \xch Y' $ arrow and

$ \xcA n \xcE \xbb_{2m+2},$ $ \xbl (\xbb_{2m+2})=2m+2,$ $ \xbb_{2m+2}$ a
valid $X \xch X' $ arrow,

which is impossible, as $ \xdx $ is a structure of finite level.

$ \xcz $
\\[3ex]

$ \xCO $

$ \xCO $
\index{Definition Totally smooth}

\bd

$\hspace{0.01em}$


\label{Definition Totally-Smooth}

Let $ \xdx $ be a generalized preferential structure, $X \xcc \xdP (\xdx
).$

$ \xdx $ is called $ \xCf totally$ $ \xCf smooth$ for $X$ iff
 \index{totally smooth}

(1) $ \xcA \xba:x \xcp y \xbe \xdA (\xdx)(O(\xba) \xcv D(\xba) \xcc
X$ $ \xch $ $ \xcE \xba':x' \xcp y.x' \xbe \xbm (X))$

(2) if $ \xba $ is valid, then there must also exist such $ \xba' $ which
is valid.

(y a point or an arrow).

If $ \xdy \xcc \xdP (\xdx),$ then $ \xdx $ is called
$\xdy-totally$ $smooth$
\index{$\xdy-$totally smooth}

iff for all $X \xbe \xdy $ $ \xdx $ is totally smooth for $X.$

\ed

$ \xCO $

$ \xCO $

\be

$\hspace{0.01em}$


\label{Example Totally-Smooth}

$X:=\{ \xba:a \xcp b,$ $ \xba':b \xcp c,$ $ \xba'':a \xcp c,$ $ \xbb
:b \xcp \xba' \}$ is not totally smooth,

$X:=\{ \xba:a \xcp b,$ $ \xba':b \xcp c,$ $ \xba'':a \xcp c,$ $ \xbb
:b \xcp \xba',$ $ \xbb':a \xcp \xba' \}$ is totally smooth.

\ee

$ \xCO $

$ \xCO $

\be

$\hspace{0.01em}$


\label{Example Need-Smooth}

Consider $ \xba':a \xcp b,$ $ \xba'':b \xcp c,$ $ \xba:a \xcp c,$ $
\xbb:a \xcp \xba.$

Then $ \xbm (\{a,b,c\})=\{a\},$ $ \xbm (\{a,c\})=\{a,c\}.$
Thus, $(\xbm CUM)$ does not hold in this structure.
Note that there is no valid arrow from $ \xbm (\{a,b,c\})$ to $c.$

\ee

$ \xCO $

$ \xCO $
\index{Definition Essentially smooth}

\bd

$\hspace{0.01em}$


\label{Definition Essentially-Smooth}

Let $ \xdx $ be a generalized preferential structure, $X \xcc \xdP (\xdx
).$

$ \xdx $ is called $ \xCf essentially$ $ \xCf smooth$ for $X$ iff $ \xbm
(X) \xes X.$
 \index{essentially smooth}
If $ \xdy \xcc \xdP (\xdx),$ then $ \xdx $ is called
$\xdy-essentially$ $smooth$

 \index{$ \xdy -$essentially smooth}
iff for all $X \xbe \xdy $
$ \xbm (X) \xes X.$

\ed

$ \xCO $
\subsubsection{
Semantic representation results for higher preferential structures
}
\paragraph{
Result on not necessarily smooth structures \\[2mm]
}

We give a representation theorem, but will make it more general than
for preferential structures only. For this purpose, we will introduce some
definitions first.

$ \xCO $

\bd

$\hspace{0.01em}$


\label{Definition Eta-Rho-Structure}

Let $ \xbh, \xbr: \xdy \xcp \xdp (U).$

(1) If $ \xdx $ is a simple structure:

$ \xdx $ is called an $ \xCf attacking$ $ \xCf structure$ relative to $
\xbh $ representing $ \xbr $ iff

 \index{attacking structure relative to $ \xbh $ representing $ \xbr $}
$ \xbr (X)$ $=$ $\{x \xbe \xbh (X):$ there is no valid $X-to- \xbh (X)$
arrow $ \xba:x' \xcp x\}$

for all $X \xbe \xdy.$

(2) If $ \xdx $ is a structure with copies:

$ \xdx $ is called an $ \xCf attacking$ $ \xCf structure$ relative to $
\xbh $ representing $ \xbr $ iff

$ \xbr (X)$ $=$ $\{x \xbe \xbh (X):$ there is $ \xBc x,i \xBe $ and no valid
$X-to-
\xbh (X)$ arrow
$ \xba: \xBc x',i'  \xBe  \xcp  \xBc x,i \xBe \}$

for all $X \xbe \xdy.$

Obviously, in those cases $ \xbr (X) \xcc \xbh (X)$ for all $X \xbe \xdy
.$

Thus, $ \xdx $ is a preferential structure iff $ \xbh $ is the identity.

See Diagram \ref{Diagram Eta-Rho-1} (page \pageref{Diagram Eta-Rho-1})

\ed

$ \xCO $

\vspace{10mm}

\begin{diagram}

\label{Diagram Eta-Rho-1}
\index{Diagram Attacking structure}

\centering
\setlength{\unitlength}{.5mm}
{\renewcommand{\dashlinestretch}{30}
\begin{picture}(150,100)(0,0)

\put(30,50){\circle{50}}
\put(100,50){\circle{40}}
\put(100,50){\circle{70}}

\put(30,50){\circle*{1}}
\put(70,50){\circle*{1}}

\path(31,50)(69,50)
\path(66.6,50.9)(69,50)(66.6,49.1)

\put(30,60){\xssc{$X$}}
\put(100,80){\xssc{$\xbh(X)$}}
\put(100,60){\xssc{$\xbr(X)$}}

\put(30,10) {{\rm\bf Attacking structure}}

\end{picture}
}

\end{diagram}

\vspace{4mm}

$ \xCO $

$ \xCO $

The following result
shows that we can obtain (almost) anything with level 2 structures.

$ \xCO $

\bp

$\hspace{0.01em}$


\label{Proposition Eta-Rho-Repres}

Let $ \xbh, \xbr: \xdy \xcp \xdp (U).$ Then there is an attacking level
2 structure relative
to $ \xbh $ representing $ \xbr $ iff

(1) $ \xbr (X) \xcc \xbh (X)$ for all $X \xbe \xdy,$

(2) $ \xbr (\xCQ)= \xbh (\xCQ)$ if $ \xCQ \xbe \xdy.$

\ep

(2) is, of course, void for preferential structures.

$ \xCO $
\paragraph{
Results on essential smoothness \\[2mm]
}

$ \xCO $

\bd

$\hspace{0.01em}$


\label{Definition ODPi}

Let $ \xbm: \xdy \xcp \xdp (U)$ and $ \xdx $ be given, let
$ \xba: \xBc y,j \xBe  \xcp  \xBc x,i \xBe  \xbe \xdx.$

Define

$ \xdO (\xba)$ \index{$ \xdO (\xba)$}  $:=$ $\{Y \xbe \xdy:x \xbe Y- \xbm
(Y),y \xbe \xbm (Y)\},$

$ \xdD (\xba)$ \index{$ \xdD (\xba)$}  $:=$ $\{X \xbe \xdy:x \xbe \xbm (X),y
\xbe X\},$

$ \xbP (\xdO, \xba)$ \index{$ \xbP (\xdO, \xba)$}  $:=$ $ \xbP \{ \xbm (Y):Y
\xbe
\xdO (\xba)\},$

$ \xbP (\xdD, \xba)$ \index{$ \xbP (\xdD, \xba)$}  $:=$ $ \xbP \{ \xbm (X):X
\xbe
\xdD (\xba)\}.$

\ed

$ \xCO $

$ \xCO $

\bl

$\hspace{0.01em}$


\label{Lemma Level-3-Constr}

Let $U$ be the universe, $ \xbm: \xdy \xcp \xdp (U).$
Let $ \xbm $ satisfy $(\xbm \xcc)+(\xbm \xcc \xcd).$

Let $ \xdx $ be a level 1 preferential structure,
$ \xba: \xBc y,j \xBe  \xcp  \xBc x,i \xBe,$ $ \xdO (\xba) \xEd \xCQ,$
$ \xdD (\xba) \xEd \xCQ.$

We can modify $ \xdx $ to a level 3 structure $ \xdx' $ by introducing
level 2 and level 3
arrows s.t. no copy of $ \xba $ is valid in any $X \xbe \xdD (\xba),$
and in every $Y \xbe \xdO (\xba)$ at
least one copy of $ \xba $ is valid. (More precisely, we should write $
\xdx' \xex X$ etc.)

Thus, in $ \xdx',$

(1) $ \xBc x,i \xBe $ will not be minimal in any $Y \xbe \xdO (\xba),$

(2) if $ \xba $ is the only arrow minimizing
$ \xBc x,i \xBe $ in $X \xbe \xdD (\xba),$ $ \xBc x,i \xBe $ will now be
minimal in $X.$

The construction is made independently for all such arrows $ \xba \xbe
\xdx.$

\el

$ \xCO $

$ \xCO $

\bp

$\hspace{0.01em}$


\label{Proposition Level-3-Repr}

Let $U$ be the universe, $ \xbm: \xdy \xcp \xdp (U).$

Then any $ \xbm $ satisfying $(\xbm \xcc)$ \index{$(\xbm \xcc)$},
$(\xcs)$ \index{$(\xcs)$},
$(\xbm CUM)$ \index{$(\xbm CUM)$}  (or, alternatively,
$(\xbm \xcc)$ and
$(\xbm \xcc \xcd)$ \index{$(\xbm \xcc \xcd)$}) can be represented by a
level 3
essentially smooth structure.

\ep

$ \xCO $
\subsubsection{
Translation to logics
}

We turn to the translation to logics.

$ \xCO $

\bp

$\hspace{0.01em}$


\label{Proposition Higher-Repr}

Let $ \xcn $ be a logic for $ \xdl.$ Set $T^{ \xdm }:=Th(\xbm_{ \xdm
}(M(T))),$ where $ \xdm $ is a
generalized preferential structure, and $ \xbm_{ \xdm }$ its choice
function. Then

(1) there is a level 2 preferential structure $ \xdm $ s.t. $ \ol{ \ol{T}
}=T^{ \xdm }$ iff
$ \xCf (LLE)$ \index{$ \xCf (LLE)$},
$ \xCf (CCL)$ \index{$ \xCf (CCL)$},
$ \xCf (SC)$ \index{$ \xCf (SC)$}  hold for all $T,T' \xcc \xdl.$

(2) there is a level 3 essentially smooth preferential structure
$ \xdm $ s.t. $ \ol{ \ol{T} }=T^{ \xdm }$ iff
$ \xCf (LLE),$ $ \xCf (CCL),$ $ \xCf (SC),$
$(\xcc \xcd)$ \index{$(\xcc \xcd)$}  hold for all $T,T' \xcc \xdl.$

\ep

$ \xCO $
\subsubsection{
Discussion of the limit version of higher preferential structures
}

\label{Section Lim-High}

In usual preferential structures, the definition of a MISE in $X,$ a
minimizing
initial segment of $X,$ was clear: It should minimize all other elements
of $X,$
and it should be downward closed. This is an intuitively clear
generalization
of the set of minimal elements. The minimal elements are the best
elements, and
the least minimal ones the worst. The set of minimal elements is the
ideal, this
set may be empty, but leaving aside some of the worst is an approximation.
This is the aim, to define what holds in the limit, in a sufficiently good
approximation of the ideal. A good approximation should not exclude better
elements, so it has to be downward closed. (In addition, this property
assures that the intersection of two MISE will not be empty, which is very
important, as we do not want to conclude FALSE.) As a MISE will minimize
all other elements, we assure that finally all bad elements
may be left out. This is all so clear and simple, as the relation
expresses
a clear relation of quality. If $x \xeb y \xeb z,$ then $x$ is better than
$y,$ $y$ is better
than $z,$ and adding transitivity is natural: $x$ is better than $z.$

This intuition is not clear for higher preferential structures. If we
read $x \xeb y$ as ``x attacks $y$'', then $x \xeb y \xeb z$ may mean that
$x$ also attacks the attack $y \xeb z.$ Of course, we can argue, to
express this idea,
it suffices to explicitly add an attack from $x$ to the attack $ \xba:y
\xeb z$ itself (not
the starting point $y),$ an arrow $ \xbb:x \xcp \xba.$ Thus, we can
assume transitivity,
if we do not want it, we just write this down using $ \xbb:x \xcp \xba.$
But then we also
may have to add some destruction of the resulting arrow $ \xba':x \xcp
z,$ this can be
done, take $ \xbb':x \xcp \xba'.$ This is already somewhat cumbersome.
But sometimes
we are interested mainly in $ \xCf valid$ (in some sense) arrows, and
concatenating
them again to a valid arrow will destroy the possibility of adding such an
arrow $ \xbb'.$ So we may still destroy $ \xba,$ but not $ \xba' $ any
more, which
seems an inconvenient.

The basic problem is that we have no clear intuition about the
``quality'' of points, in the way we had it for usual preferential
structures - as described above. Asking for strong formal properties like
transitivity of ``valid'' arrows is perhaps counterintuitive. Unfortunately,
assuming transitivity was at the heart of a reasonable definition of a
MISE (assuring among other things that finite AND holds: the intersection
of two MISE is a MISE, see
Fact \ref{Fact D-8.2.1} (page \pageref{Fact D-8.2.1})  (2)). Of course, we can
put our qualms
about intuition aside and assume sufficient transitivity. But we see
immediately
a new problem. One of the nice properties of MISE was to assure a unitary
form of cumulativity, see
Fact \ref{Fact Lim-Cum} (page \pageref{Fact Lim-Cum}), again
related to transitivity: If $A \xcc B$ is a MISE in $B,$ $B \xcc C$ a MISE
in $C,$ then
$A \xcc C$ is a MISE in $C,$
see Fact \ref{Fact D-8.2.1} (page \pageref{Fact D-8.2.1}), (2a).
In higher preferential structures, we would want to have (among other
properties), that if $c \xbe C- \xCf B,$ then there is a
``valid'' arrow $ \xba $ from some $b \xbe B$ to $c,$ $ \xba:b \xcp c.$
See Diagram \ref{Diagram High-Cum} (page \pageref{Diagram High-Cum}).
Moreover, if $b \xbe B- \xCf A,$ then there should be a
``valid'' arrow from some $ \xCf a$ to $b,$ $ \xba':a \xcp b.$
Finally, we will want some ``valid'' arrow $ \xba'':a \xcp c.$
What does ``valid'' mean here. If it means that there is no (valid) attack
against the arrow in the big set, i.e., no valid attack against $ \xba $
from $C,$
no valid attack against $ \xba' $ from $B,$ then it is unclear why there
should not be
any valid attack against $ \xba'' $ from $C,$ as
``$ \xCf A$ does not know about $C,$ only about $B$''.
Thus, we have to restrict ourselves to consider attacks against $ \xba $
from $B,$
against $ \xba' $ from $ \xCf A,$ and again against $ \xba'' $ from $
\xCf A.$

In the Diagram \ref{Diagram High-Cum} (page \pageref{Diagram High-Cum}),
$ \xba $ is attacked by $ \xbb_{0}$ from $C,$ by $ \xbb $ from $B,$ but
only the latter attack
has to be countered by $ \xbg.$ The situation is similar for $ \xba'.$
By
transitivity, there is now $ \xba''.$ The attack $ \xbb''_{0}$ need not
be countered,
but the attack $ \xbb'',$ it is countered by $ \xbg''.$

But why should an attack $ \xbb''_{0}$ be admitted without counterattack?
The only
reason seems to be that it originates from elements which are attacked
themselves. So an attack against an element weakens attacks originating
from this element, something we had doubts about.

It seems difficult to find a plausible solution without a clear intuition
about strength. Thus, we limit ourselves to a very general and informal
description, and leave the question open for further research.
A MISE $A \xcc B$ for a higher preferential structure should have the
following
vague properties:

 \xEh
 \xDH
For all $b \xbe B-A$ there should be $a \xbe A$ and a valid arrow $ \xba
:a \xcp b.$

 \xDH
$ \xCf A$ should be downward closed in $B:$ If there is a valid arrow $
\xba:b \xcp a$
to some $a \xbe A,$ from some $b \xbe B,$ then $b$ should be in $ \xCf A.$
Or, in other words,
if $b \xbe B-$A, then every arrow $ \xba:b \xcp a$ into $ \xCf A$ should
be validly attacked.
 \xEj

The next section is an attempt to clarify the intuition - of one of
us (Karl Schlechta). It seems that the authors have here somewhat
diverging intuitions.
\subsubsection{
Ideas for the intuitive meaning of higher preferential structures
}

It seems we have to differentiate:
 \xEh
 \xDH relations like in preferential structures (seen as arrows)
 \xDH arrows of support
 \xDH arrows of attack
 \xEj

Examples:

 \xEh

 \xDH

Consider $a \xcp b,$ $c \xcp b.$

 \xEh
 \xDH In preferential structures, this means that
$ \xCf a$ and $c$ are better than $b.$ If we just have $a \xcp b,$ this
does not mean that
$b$ is now better than in the first situation.

 \xDH If we see arrows as support, in the first situation, $b$ is better
than in the
second situation.

 \xDH But if we see arrows as attacks,
in the second situation, there are less attacks against $b,$ so $b$ is
better.

 \xEj

 \xDH

Consider $a \xcp b \xcp c.$

 \xEh
 \xDH In preferential structures, transitivity is reasonable,
considering arrows as expressing quality.

 \xDH So it is for arrows as support, but it
seems reasonable to say that the support from $ \xCf a$ is less strong
than the
support from $b,$ being indirect support.

 \xDH In the third case, transitivity is not wanted,
on the contrary. The attack of $ \xCf a$ on $b$ attenuates the attack of
$b$ on $c.$ Still,
it seems reasonable to say that the attack does not disappear totally.

 \xDH
Suppose we add a supporting arrow $d \xcp c.$

Then, if all are supports,
the whole situation is better for $c$ than the situation with only $d \xcp
c.$

If
the arrows $a \xcp b \xcp c$ are attacks, then the whole situation is
worse for $c$
than the situation with only $d \xcp c.$

 \xEj

 \xDH

Suppose we go upward from a given set $X$ of points. So all $x \xbe X$
will
have full credibility.

 \xEh
 \xDH For preference relations, we get worse going upward.

 \xDH For support, the longer a support chain, the weaker the support it
gives. Different chains add their support (we have to define different!).

 \xDH For attack, we alternate, but attacks get weaker with longer chains.
Again,
different attacks add up.

 \xEj

 \xDH

Suppose we go downward from some point $x$ in an infinite descending
chain.

 \xEh
 \xDH In the preferential case, we just get better, without any optimum,
this makes
sense.

 \xDH In the support case, we have
``support out of thin air'', such chains should not be considered.

E.g., we have arguments that the moon is made from cheese:
The arguments go like this:
$A_{0}:$
We believe that there is a herd of at least $10^{10}$ cows hidden behind
the moon.

$A_{1}:$
We believe that there is a herd of at least $10^{10}+1$ cows hidden behind
the moon.

$A_{2}:$
We believe that there is a herd of at least $10^{10}+2$ cows hidden behind
the moon.

Etc. Of course, $A_{n+1}$ implies $A_{n},$ and $A_{0}$ makes the cheese
hypothesis seem
reasonable.

 \xDH In the attack case, we have
``attack out of thin air'', again, such chains should not be considered.

 \xEj

 \xDH
Cycles:

 \xEh
 \xDH Preferences: the intuitively best interpretation is probably to see
them as $ \xec,$ and not $ \xeb.$

 \xDH For the other cases, it seems we have to distinguish whether we
go into a cycle, or come out of it:

If we go into a cycle, we should stop after going around once, and treat
the
result as usual.

If we come out of a cycle (with no linear path going into it) we should
consider
it just as an infinite descending chain, i.e., neglect it.

 \xEj

 \xDH
Attacks on attacks etc., some attenuation, like length of path, should be
considered.
This is important for the limit version of higher preferential structures.

 \xDH
Do attacks and support only go against other arrows, or also against
points?

 \xEj

$ \xCO $

\vspace{10mm}

\begin{diagram}

\label{Diagram High-Cum}
\index{Diagram High-Cum}

\centering
\setlength{\unitlength}{1mm}
{\renewcommand{\dashlinestretch}{30}
\begin{picture}(150,150)(0,0)

\path(5,145)(120,145)(120,5)(5,5)(5,145)
\path(5,110)(120,110)
\path(5,50)(120,50)
\put(125,140){$C$}
\put(125,105){$B$}
\put(125,45){$A$}

\put(30,20){\circle*{1}}
\put(90,80){\circle*{1}}
\put(30,140){\circle*{1}}

\path(30,21)(30,137)
\path(29,134.2)(30,137)(31,134.2)
\put(25,80){{\xssc $\xba''$}}
\path(31,21)(89,79)
\path(86.3,77.7)(89,79)(87.7,76.3)
\put(60,53){{\xssc $\xba'$}}
\path(89,81)(32,137)
\path(33.3,134.3)(32,137)(34.7,135.8)
\put(60,105){{\xssc $\xba$}}

\path(9,35)(29,35)
\path(26.3,36)(29,35)(26.3,34)
\put(17,36){{\xssc $\xbb''$}}
\path(9,125)(29,125)
\path(26.3,126)(29,125)(26.3,124)
\put(17,126){{\xssc $\xbb''_0$}}

\path(46,35)(66,35)
\path(48.7,34)(46,35)(48.7,36)
\put(54,36){{\xssc $\xbb'$}}
\path(46,125)(66,125)
\path(48.7,124)(46,125)(48.7,126)
\put(54,126){{\xssc $\xbb_0$}}

\path(76,95)(96,95)
\path(78.7,94)(76,95)(78.7,96)
\put(84,96){{\xssc $\xbb$}}
\path(76,65)(96,65)
\path(78.7,64)(76,65)(78.7,66)
\put(84,66){{\xssc $\xbb'_0$}}

\path(19,19)(19,34)
\path(18.3,32)(19,34)(19.7,32)
\put(15,24){{\xssc $\xbg''$}}
\path(56,19)(56,34)
\path(55.3,32)(56,34)(56.7,32)
\put(57,24){{\xssc $\xbg'$}}

\path(90,84)(90,94)
\path(89.5,92.7)(90,94)(90.5,92.7)
\put(91,88){{\xssc $\xbg$}}

\put(29,18){{\xssc $a$}}
\put(29,141){{\xssc $c$}}
\put(91,79){{\xssc $b$}}

\end{picture}
}

\end{diagram}

\vspace{4mm}

$ \xCO $
\clearpage
\section{
Theory revision
}

\label{Section Mod-TR}

We give here the basic concepts and ideas of theory revision,
first the AGM approach, and then distance based revision.

$ \xCO $

\label{Section Toolbase1-TR-AGM}

\label{Section AGM-revision}

All material in this Section \ref{Section AGM-revision} (page \pageref{Section
AGM-revision})  is due verbatim
or in essence
to AGM - AGM for Alchourron, Gardenfors, Makinson, see e.g.,  \cite{AGM85}.


\index{Definition AGM revision}

\bd

$\hspace{0.01em}$


\label{Definition AGM}

We present in parallel the logical
and the semantic (or purely algebraic) side. For the latter, we work in
some
fixed universe $U,$ and the intuition is $U=M_{ \xdl },$ $X=M(K),$ etc.,
so, e.g., $A \xbe K$
becomes $X \xcc B,$ etc.

(For reasons of readability, we omit most caveats about definability.)

$K_{ \xcT }$ \index{$K_{ \xcT }$}  will denote the inconsistent theory.

We consider two functions, - and $*,$ taking a deductively closed theory
and a
formula as arguments, and returning a (deductively closed) theory on the
logics
side. The algebraic counterparts work on definable model sets. It is
obvious
that $ \xCf (K-1),$ $(K*1),$ $ \xCf (K-6),$ $(K*6)$ have vacuously true
counterparts on the
semantical side. Note that $K$ $ \xCf (X)$ will never change, everything
is relative
to fixed $K$ $ \xCf (X).$ $K* \xbf $ is the result of revising $K$ with $
\xbf.$ $K- \xbf $ is the result of
subtracting enough from $K$ to be able to add $ \xCN \xbf $ in a
reasonable way, called
contraction.

Moreover,
let $ \xck_{K}$ be a relation on the formulas relative to a deductively
closed theory $K$
on the formulas of $ \xdl,$ and $ \xck_{X}$ a relation on $ \xdp (U)$ or
a suitable subset of $ \xdp (U)$
relative to fixed $X.$ When the context is clear, we simply write $ \xck
.$
$ \xck_{K}$ $(\xck_{X})$ is called a relation of epistemic entrenchment
for $K$ $ \xCf (X).$

Table \ref{Table Base1-AGM-TR} (page \pageref{Table Base1-AGM-TR}),
``AGM theory revision'',
presents
``rationality postulates'' for contraction (-),
rationality postulates \index{rationality postulates}
revision $(*)$ and epistemic entrenchment. In AGM tradition, $K$ will be a
deductively closed theory, $ \xbf, \xbq $ formulas. Accordingly, $X$ will
be the set of
models of a theory, $A,B$ the model sets of formulas.

In the further development, formulas $ \xbf $ etc. may sometimes also be
full
theories. As the transcription to this case is evident, we will not go
into
details.

\renewcommand{\arraystretch}{1.2}

\begin{table}[h]

\index{Contraction}
\index{Revision}
\index{Epistemic entrenchment}
\index{$K-\xbf $}
\index{$K*\xbf $}
\index{$(K-1)$}
\index{$(K-2)$}
\index{$(K-3)$}
\index{$(K-4)$}
\index{$(K-5)$}
\index{$(K-6)$}
\index{$(K-7)$}
\index{$(K-8)$}
\index{$(K*1)$}
\index{$(K*2)$}
\index{$(K*3)$}
\index{$(K*4)$}
\index{$(K*5)$}
\index{$(K*6)$}
\index{$(K*7)$}
\index{$(K*8)$}
\index{$(EE1)$}
\index{$(EE2)$}
\index{$(EE3)$}
\index{$(EE4)$}
\index{$(EE5)$}

\caption{AGM theory revision}

\label{Table Base1-AGM-TR}
\begin{center}

{\small

\begin{tabular}{|c|c|c|c|}

\hline

\multicolumn{4}{|c|}{\bf AGM theory revision}\\
\hline

\multicolumn{4}{|c|} {Contraction, $K-\xbf $} \xEP

\hline

$(K-1)$ \xEH $K-\xbf $ is deductively closed \xEH \xEH \xEP

\hline

$(K-2)$ \xEH $K-\xbf $ $ \xcc $ $K$ \xEH $(X \xDN 2)$ \xEH $X \xcc X \xDN A$
\xEP

\hline

$(K-3)$ \xEH $\xbf  \xce K$ $ \xch $ $K-\xbf =K$ \xEH $(X \xDN 3)$ \xEH $X \xcC
A$ $
\xch $ $X \xDN A=X$ \xEP

\hline

$(K-4)$ \xEH $ \xcL \xbf $ $ \xch $ $\xbf  \xce K-\xbf $ \xEH $(X \xDN 4)$ \xEH
$A \xEd
U$ $ \xch $ $X \xDN A \xcC A$ \xEP

\hline

$(K-5)$ \xEH $K \xcc \ol{(K-\xbf) \xcv \{\xbf \}}$ \xEH $(X \xDN 5)$ \xEH $(X
\xDN
A) \xcs A$ $ \xcc $ $X$ \xEP

\hline

$(K-6)$ \xEH $ \xcl \xbf  \xcr \xbq $ $ \xch $ $K-\xbf =K-\xbq $ \xEH \xEH \xEP

\hline

$(K-7)$ \xEH $(K-\xbf) \xcs (K-\xbq)  \xcc  $ \xEH
$(X \xDN 7)$ \xEH $X \xDN (A \xcs B)  \xcc  $ \xEP

\xEH $K-(\xbf  \xcu \xbq) $ \xEH
\xEH $(X \xDN A) \xcv (X \xDN B)$ \xEP

\hline

$(K-8)$ \xEH $\xbf  \xce K-(\xbf  \xcu \xbq)  \xch  $ \xEH
$(X \xDN 8)$ \xEH $X \xDN (A \xcs B) \xcC A  \xch  $ \xEP

\xEH $K-(\xbf  \xcu \xbq) \xcc K-\xbf $ \xEH
\xEH $X \xDN A \xcc X \xDN (A \xcs B)$ \xEP

\hline
\hline

\multicolumn{4}{|c|} {Revision, $K*\xbf $} \xEP

\hline

$(K*1)$ \xEH $K*\xbf $ is deductively closed \xEH - \xEH \xEP

\hline

$(K*2)$ \xEH $\xbf  \xbe K*\xbf $ \xEH $(X \xfA 2)$ \xEH $X \xfA A \xcc A$
\xEP

\hline

$(K*3)$ \xEH $K*\xbf $ $ \xcc $ $ \ol{K \xcv \{\xbf \}}$ \xEH $(X \xfA 3)$ \xEH
$X
\xcs A \xcc X \xfA A$ \xEP

\hline

$(K*4)$ \xEH $ \xCN \xbf  \xce K  \xch $ \xEH
$(X \xfA 4)$ \xEH $X \xcs A \xEd \xCQ   \xch $ \xEP

\xEH $\ol{K \xcv \{\xbf \}}  \xcc  K*\xbf $ \xEH
\xEH $X \xfA A \xcc X \xcs A$ \xEP

\hline

$(K*5)$ \xEH $K*\xbf =K_{ \xcT }$ $ \xch $ $ \xcl \xCN \xbf $ \xEH $(X \xfA 5)$
\xEH $X \xfA A= \xCQ $ $ \xch $ $A= \xCQ $ \xEP

\hline

$(K*6)$ \xEH $ \xcl \xbf  \xcr \xbq $ $ \xch $ $K*\xbf =K*\xbq $ \xEH - \xEH
\xEP

\hline

$(K*7)$ \xEH $K*(\xbf  \xcu \xbq)  \xcc $ \xEH
$(X \xfA 7)$ \xEH $(X \xfA A) \xcs B  \xcc  $ \xEP

\xEH $\ol{(K*\xbf) \xcv \{\xbq \}}$ \xEH
\xEH $X \xfA (A \xcs B)$ \xEP

\hline

$(K*8)$ \xEH $ \xCN \xbq  \xce K*\xbf  \xch $ \xEH
$(X \xfA 8)$ \xEH $(X \xfA A) \xcs B \xEd \xCQ \xch $ \xEP

\xEH $\ol{(K*\xbf) \xcv \{\xbq \}} \xcc K*(\xbf  \xcu \xbq)$ \xEH
\xEH $ X \xfA (A \xcs B) \xcc (X \xfA A) \xcs B$ \xEP

\hline
\hline

\multicolumn{4}{|c|} {Epistemic entrenchment} \xEP

\hline

$(EE1)$ \xEH $ \xck_{K}$ is transitive \xEH
$(EE1)$ \xEH $ \xck_{X}$ is transitive \xEP

\hline

$(EE2)$ \xEH $\xbf  \xcl \xbq   \xch  \xbf  \xck_{K}\xbq $ \xEH
$(EE2)$ \xEH $A \xcc B  \xch  A \xck_{X}B$ \xEP

\hline

$(EE3)$ \xEH $ \xcA  \xbf,\xbq  $ \xEH
$(EE3)$ \xEH $ \xcA A,B $ \xEP
\xEH $ (\xbf  \xck_{K}\xbf  \xcu \xbq $ or $\xbq  \xck_{K}\xbf  \xcu \xbq)$
\xEH
\xEH $ (A \xck_{X}A \xcs B$ or $B \xck_{X}A \xcs B)$ \xEP

\hline

$(EE4)$ \xEH $K \xEd K_{ \xcT }  \xch  $ \xEH
$(EE4)$ \xEH $X \xEd \xCQ   \xch  $ \xEP
\xEH $(\xbf  \xce K$ iff $ \xcA  \xbq.\xbf  \xck_{K}\xbq)$ \xEH
\xEH $(X \xcC A$ iff $ \xcA B.A \xck_{X}B)$ \xEP

\hline

$(EE5)$ \xEH $ \xcA \xbq.\xbq  \xck_{K}\xbf   \xch   \xcl \xbf $ \xEH
$(EE5)$ \xEH $ \xcA B.B \xck_{X}A  \xch  A=U$ \xEP

\hline

\end{tabular}
}

\end{center}
\end{table}

\vspace{3mm}


\vspace{3mm}



\ed

\br

$\hspace{0.01em}$


\label{Remark TR-Rank}

(1) Note that $(X \xfA 7)$ and $(X \xfA 8)$ express a central condition
for ranked
structures: If we note $X \xfA.$ by $f_{X}(.),$ we then have:
$f_{X}(A) \xcs B \xEd \xCQ $ $ \xch $ $f_{X}(A \xcs B)=f_{X}(A) \xcs B.$

(2) It is trivial to see that AGM revision cannot be defined by an
individual
distance (see
Definition \ref{Definition Dist-Indiv-Coll} (page \pageref{Definition
Dist-Indiv-Coll})):
Suppose $X \xfA Y$ $:=$ $\{y \xbe Y:$ $ \xcE x_{y} \xbe X(\xcA y' \xbe
Y.d(x_{y},y) \xck d(x_{y},y'))\}.$
Consider $a,b,c.$ $\{a,b\} \xfA \{b,c\}=\{b\}$ by $(X \xfA 3)$ and $(X
\xfA 4),$ so $d(a,b)<d(a,c).$
But on the other hand $\{a,c\} \xfA \{b,c\}=\{c\},$ so $d(a,b)>d(a,c),$
$contradiction.$

\vspace{3mm}


\vspace{3mm}



\er

\bp

$\hspace{0.01em}$


\label{Proposition AGM-Equiv}

We refer here to
Table \ref{Table Base1-AGM-Inter} (page \pageref{Table Base1-AGM-Inter}),
``AGM interdefinability''.
Contraction, revision, and epistemic entrenchment are interdefinable by
the
following equations, i.e., if the defining side has the respective
properties,
so will the defined side. (See  \cite{AGM85}.)

\renewcommand{\arraystretch}{1.5}

\begin{table}
\caption{AGM interdefinability}

\label{Table Base1-AGM-Inter}
\begin{center}

{\scriptsize
\begin{tabular}{|c|c|}

\hline

\multicolumn{2}{|c|}{\bf AGM interdefinability}\\
\hline

$K*\xbf:= \ol{(K- \xCN \xbf)} \xcv {\xbf }$  \xEH  $X \xfA A:= (X \xDN 
\xdC(A)) \xcs A$ \xEP

\hline

$K-\xbf:= K \xcs (K* \xCN \xbf)$  \xEH  $X \xDN A:= X \xcv (X \xfA  \xdC(A))$
\xEP

\hline

$K-\xbf:=\{\xbq  \xbe K:$ $(\xbf <_{K}\xbf  \xco \xbq $ or $ \xcl \xbf)\}$ \xEH

$
X \xDN A:=
\left\{
\begin{array}{rcl}
X & iff & A=U, \\
 \xcS \{B: X \xcc B \xcc U, A<_{X}A \xcv B\} & & otherwise \\
\end{array}
\right. $
\xEP

\hline

$
\xbf  \xck_{K}\xbq: \xcr
\left\{
\begin{array}{l}
\xcl \xbf  \xcu \xbq  \\
or \\
\xbf  \xce K-(\xbf  \xcu \xbq) \\
\end{array}
\right. $
\xEH

$
A \xck_{X}B: \xcr
\left\{
\begin{array}{l}
A,B=U  \\
or \\
X \xDN (A \xcs B) \xcC A \\
\end{array}
\right. $
\xEP

\hline

\end{tabular}
}

\end{center}
\end{table}

\ep

Speaking in terms of distance defined revision, $X \xfA A$ is the set of
those
$a \xbe A,$ which are closest to $X,$ and $X \xDN A$ is the set of $y$
which are either
in $X,$ or in $ \xdC (A)$ and closest to $X$ among those in $ \xdC (A).$

\vspace{3mm}


\vspace{3mm}



\subsubsection{A remark on intuition}

The idea of epistemic entrenchment is that $ \xbf $ is more entrenched
than $ \xbq $
(relative to $K)$ iff $M(\xCN \xbq)$ is closer to $M(K)$ than $M(\xCN
\xbf)$ is to $M(K).$ In
shorthand, the more we can twiggle $K$ without reaching $ \xCN \xbf,$ the
more $ \xbf $ is
entrenched. Truth is maximally entrenched - no twiggling whatever will
reach
falsity. The more $ \xbf $ is entrenched,
the more we are certain about it. Seen this way, the properties of
epistemic
entrenchment relations are very natural (and trivial): As only the closest
points of $M(\xCN \xbf)$ count (seen from $M(K)),$ $ \xbf $ or $ \xbq $
will be as entrenched as
$ \xbf \xcu \xbq,$ and there is a logically strongest $ \xbf' $ which is
as entrenched as $ \xbf $ -
this is just the sphere around $M(K)$ with radius $d(M(K),M(\xCN \xbf
)).$

\vspace{3mm}


\vspace{3mm}

$ \xCO $

$ \xCO $

\label{Section Toolbase1-TR-DistBase}


\index{Definition Distance}

\bd

$\hspace{0.01em}$


\label{Definition Distance}

$d:U \xCK U \xcp Z$ is called a pseudo-distance \index{pseudo-distance}  on $U$
iff (d1) holds:

(d1) $Z$ is totally ordered by a relation $<.$

If, in addition, $Z$ has a $<-$smallest element 0, and (d2) holds, we say
that $d$
respects identity:

(d2) $d(a,b)=0$ iff $a=b.$

If, in addition, (d3) holds, then $d$ is called symmetric:

(d3) $d(a,b)=d(b,a).$

(For any $a,b \xbe U.)$

Note that we can force the triangle inequality to hold trivially (if we
can
choose the values in the real numbers): It suffices to choose the values
in
the set $\{0\} \xcv [0.5,1],$ i.e., in the interval from 0.5 to 1, or as
0.

\vspace{3mm}


\vspace{3mm}


\index{Definition Individual/collective distance}

\ed

\bd

$\hspace{0.01em}$


\label{Definition Dist-Indiv-Coll}

We define the collective and the individual variant of choosing the
closest
elements in the second operand by two operators,
$ \xfA, \xfB: \xdp (U) \xCK \xdp (U) \xcp \xdp (U):$

Let $d$ be a distance or pseudo-distance.

$X \xfA Y$ \index{$X \xfA Y$}  $:=$ $\{y \xbe Y:$ $ \xcE x_{y} \xbe X. \xcA
x' \xbe X, \xcA y' \xbe Y(d(x_{y},y) \xck d(x',y')\}$

(the collective variant \index{collective variant}, used in theory revision)

and

$X \xfB Y$ \index{$X \xfB Y$}  $:=$ $\{y \xbe Y:$ $ \xcE x_{y} \xbe X. \xcA
y' \xbe Y(d(x_{y},y) \xck d(x_{y},y')\}$

(the individual variant \index{individual variant}, used for counterfactual
conditionals and theory
update).

Thus, $A \xfA_{d}B$ is the subset of $B$ consisting of all $b \xbe B$ that
are closest to A.
Note that, if $ \xCf A$ or $B$ is infinite, $A \xfA_{d}B$ may be empty,
even if $ \xCf A$ and $B$ are not
empty. A condition assuring nonemptiness will be imposed when necessary.

\vspace{3mm}


\vspace{3mm}


\index{Definition Distance representation, semantics}

\ed

\bd

$\hspace{0.01em}$


\label{Definition Dist-Repr}

An operation $ \xfA: \xdp (U) \xCK \xdp (U) \xcp \xdp (U)$ is
representable iff there is a
pseudo-distance $d:U \xCK U \xcp Z$ such that

$A \xfA B$ $=$ $A \xfA_{d}B$ \index{$A \xfA_{d}B$}  $:=$ $\{b \xbe B:$ $ \xcE
a_{b} \xbe A \xcA a' \xbe A \xcA b' \xbe B(d(a_{b},b) \xck d(a',b'))\}.$

\vspace{3mm}


\vspace{3mm}


\index{Definition Distance representation, logics}

\ed

The following is the central definition, it describes the way a revision
$*_{d}$ is
attached to a pseudo-distance $d$ on the set of models.

\bd

$\hspace{0.01em}$


\label{Definition TR*d}

$T*_{d}T' $ \index{$T*_{d}T' $}  $:=$ $Th(M(T) \xfA_{d}M(T')).$

$*$ is called representable \index{representable}  iff there is a
pseudo-distance $d$ on the set of
models
s.t. $T*T' =Th(M(T) \xfA_{d}M(T')).$

\vspace{3mm}


\vspace{3mm}

\ed

$ \xCO $

$ \xCO $

\bfa

$\hspace{0.01em}$


\label{Fact AGM-In-Dist-N}

A distance based revision satisfies the AGM postulates provided:

(1) it respects identity, i.e., $d(a,a)<d(a,b)$ for all $a \xEd b,$

(2) it satisfies a limit condition: minima exist,

(3) it is definability preserving.

\efa

(All conditions are necessary.)

$ \xCO $

$ \xCO $
\index{Definition TR-Dist-N}

\bd

$\hspace{0.01em}$


\label{Definition TR-Dist}

We refer here to
Table \ref{Table Base1-TR-Dist} (page \pageref{Table Base1-TR-Dist}),
``Distance representation and revision''.
It defines the Loop Condition, and shows the correspondence
between the semantic and the syntactic side.

The prerequisites are:

Let $U \xEd \xCQ,$ $ \xdy \xcc \xdp (U)$ satisfy $(\xcs),$ $(\xcv),$
$ \xCQ \xce \xdy.$

Let $A,B,X_{i} \xbe \xdy,$ $ \xfA: \xdy \xCK \xdy \xcp \xdp (U).$

Let $*$ be a revision function defined for
arbitrary consistent theories on both sides. (This is thus a slight
extension of
the AGM framework, as AGM work with formulas only on the right of $*.)$

\begin{table}
\caption{Distance representation and revision}

\label{Table Base1-TR-Dist}
\begin{center}

\begin{turn}{90}

{\scriptsize

\begin{tabular}{|c|c|c|}

\hline

\multicolumn{3}{|c|}{\bf Distance representation and revision}\\
\hline

\xEH
\xEH
$(*Equiv)$
\xEP
\xEH
\xEH
$ \xcm T \xcr S,$ $ \xcm T' \xcr S',$ $\xch$ $T*T' =S*S',$
\xEP

\hline

\xEH
\xEH
$(*CCL)$
\xEP
\xEH
\xEH
$T*T' $ is a consistent, deductively closed theory,
\xEP

\hline

\xEH
$(\xfA Succ)$
\xEH
$(*Succ)$
\xEP
\xEH
$A \xfA B \xcc B$
\xEH
$T' \xcc T*T',$
\xEP

\hline

\xEH
$(\xfA Con)$
\xEH
$(*Con)$
\xEP
\xEH
$A \xcs B \xEd \xCQ $ $ \xch $ $A \xfA B=A \xcs B$
\xEH
$Con(T \xcv T') $ $\xch$ $T*T' = \ol{T \xcv T' },$
\xEP

\hline

Intuitively,
\xEH
$(\xfA Loop)$
\xEH
$(*Loop)$
\xEP
Using symmetry
\xEH
\xEH
\xEP
$d(X_{0},X_{1}) \xck d(X_{1},X_{2}),$
\xEH
$(X_{1} \xfA (X_{0} \xcv X_{2})) \xcs X_{0} \xEd \xCQ,$
\xEH
$Con(T_{0},T_{1}*(T_{0} \xco T_{2})),$
\xEP
$d(X_{1},X_{2}) \xck d(X_{2},X_{3}),$
\xEH
$(X_{2} \xfA (X_{1} \xcv X_{3})) \xcs X_{1} \xEd \xCQ,$
\xEH
$Con(T_{1},T_{2}*(T_{1} \xco T_{3})),$
\xEP
$d(X_{2},X_{3}) \xck d(X_{3},X_{4})$
\xEH
$(X_{3} \xfA (X_{2} \xcv X_{4})) \xcs X_{2} \xEd \xCQ,$
\xEH
$Con(T_{2},T_{3}*(T_{2} \xco T_{4}))$
\xEP
\Xl
\xEH
\Xl
\xEH
\Xl
\xEP
$d(X_{k-1},X_{k}) \xck d(X_{0},X_{k})$
\xEH
$(X_{k} \xfA (X_{k-1} \xcv X_{0})) \xcs X_{k-1} \xEd \xCQ $
\xEH
$Con(T_{k-1},T_{k}*(T_{k-1} \xco T_{0}))$
\xEP
$\xch$
\xEH
$\xch$
\xEH
$\xch$
\xEP
$d(X_{0},X_{1}) \xck d(X_{0},X_{k}),$
\xEH
$(X_{0} \xfA (X_{k} \xcv X_{1})) \xcs X_{1} \xEd \xCQ$
\xEH
$Con(T_{1},T_{0}*(T_{k} \xco T_{1}))$
\xEP

i.e., transitivity, or absence of
\xEH
\xEH
\xEP

loops involving $<$
\xEH
\xEH
\xEP

\hline

\end{tabular}
}
\end{turn}

\end{center}
\end{table}

\ed

$ \xCO $

$ \xCO $

\bp

$\hspace{0.01em}$


\label{Proposition TR-Alg-Repr}

Let $U \xEd \xCQ,$ $ \xdy \xcc \xdp (U)$ be closed under finite $ \xcs $
and finite $ \xcv,$ $ \xCQ \xce \xdy.$

(a) $ \xfA $ is representable by a symmetric pseudo-distance $d:U \xCK U
\xcp Z$ iff $ \xfA $
satisfies $(\xfA Succ)$ \index{$(\xfA Succ)$}  and
$(\xfA Loop)$ \index{$(\xfA Loop)$}  in
Definition \ref{Definition TR-Dist} (page \pageref{Definition TR-Dist}).

(b) $ \xfA $ is representable by an identity respecting symmetric
pseudo-distance
$d:U \xCK U \xcp Z$ iff $ \xfA $ satisfies
$(\xfA Succ)$ \index{$(\xfA Succ)$},
$(\xfA Con)$ \index{$(\xfA Con)$}, and
$(\xfA Loop)$ \index{$(\xfA Loop)$}
in Definition \ref{Definition TR-Dist} (page \pageref{Definition TR-Dist}).

See  \cite{LMS01} or  \cite{Sch04}.

\ep

$ \xCO $

$ \xCO $

\bp

$\hspace{0.01em}$


\label{Proposition TR-Log-Repr}

Let $ \xdl $ be a propositional language.

(a) A revision operation $*$ is representable by a symmetric consistency
and
definability preserving pseudo-distance iff $*$ satisfies
$(*Equiv)$ \index{$(*Equiv)$},
$(*CCL)$ \index{$(*CCL)$},
$(*Succ)$ \index{$(*Succ)$},
$(*Loop)$ \index{$(*Loop)$}.

(b) A revision operation $*$ is representable by a symmetric consistency
and
definability preserving, identity respecting pseudo-distance iff $*$
satisfies $(*Equiv),$ $(*CCL),$ $(*Succ),$ $(*Con)$ \index{$(*Con)$},
$(*Loop).$

See  \cite{LMS01} or  \cite{Sch04}.

\ep

$ \xCO $
\subsection{
Theory revision for many-valued logics
}

To see what we have to do for distance based revision in the case of
many-valued logics, we orient ourselves by the 2-valued case.

We considered sets $X$ and $Y,$ and were looking at pairs
$ \xBc x,y \xBe,$ $x \xbe X,$ $y \xbe Y$ and choose those with minimal
distance.
So we compared $ \xBc x,y \xBe $ and $ \xBc x',y'  \xBe,$ $x,x' \xbe X,$ $y,y'
\xbe Y,$ and
$d(x,y)$ and $d(x',y').$
Thus, we now look pairs $ \xBc x,y \xBe,$ $ \xBc x',y'  \xBe $ again, and
discard $ \xBc x,y \xBe $ if
$d(x',y')<d(x,y),$ $ \xCf and$
``the value of $ \xBc x',y'  \xBe $ is at least as good as the value of $ \xBc
x,y \xBe $''.
To make the latter precise, we postulate, setting $F$ for $X,$ $G$ for
$Y:$
$F(x) \xck F(x')$ and $G(y) \xck G(y').$
Thus, we consider

$\{ \xBc x,y \xBe $: $ \xCN \xcE  \xBc x',y'  \xBe.(d(x',y')<d(x,y)$ and (
$(F(x) \xck
F(x')$ and $G(y) \xck G(y'))$ or
$(G(x) \xck G(x')$ and $F(y) \xck F(y'))$)) $\}.$
\clearpage
\chapter{
Towards a uniform picture of conditionals
}

\label{Chapter Mod-Uni-Cond}
\section{
Introduction
}

\label{Section Mod-Uni-Int}

The word ``conditional'' contains
``condition''. So, a conditional is a structure of the
form:
``if condition $c$ holds, then so will property $p$''.
Thus, conditionals seem to be at least always binary.
But the condition may be hidden in additional structure
(like a Kripke model relation, a preferential relation, etc.).
Thus, we will also include unary structures.

But even if we restrict ourselves to binary (and why not ternary
etc.?) conditionals, it seems
that we can invent ad libitum new conditionals.
We give two (arbitrary) examples, just to illustrate the possibilities,
we do not pretend that they are very intuitive:

\be

$\hspace{0.01em}$


\label{Example Libi-Cond}

(1) A new binary conditional:
Suppose we have a distance $d$ between models, and, in addition, a real
valued
function (e.g., of utility) $f$ defined on the model set. Define now
$m \xcm \xbf > \xbq $ iff in all $ \xbf -$worlds $n,$ which are closest to
$m,$ and where $f$ is locally
constant, i.e., $f(n)$ does not change in a small $d-$neighbourhood around
$n,$
$ \xbq $ holds.

(2) A new ternary conditional:
Let again a distance $d$ be defined. $(\xbf, \xbq, \xbf')$ holds iff
in all worlds
which are equidistant to the $ \xbf $ and $ \xbf' $ worlds (defined as in
theory
revision), $ \xbq $ holds.

\ee

So it seems impossible to give an exhaustive enumeration of all possible
conditionals.

We look at different possible classifications:

We may classify conditionals as to

 \xEh

 \xDH
their arity,

 \xDH
whether they are in the object language (like counterfactual
conditionals),
or in the meta-language (like usual preferential consequences), and if so,
can we nest them, or do we run into trivialization results as for
theory revision,

 \xDH
whether they are based on classical logic, like usual modal logic, or
perhaps
some other logic,

 \xDH
according to the properties of their semantic choice functions, like the
properties
of the $ \xbm -$functions of preferential structures,

 \xDH
whether they work with $ \xCf one$ set of chosen models (as in the minimal
variant of
preferential structures), or with a $ \xCf family$ of chosen sets (as in
the
limit version of preferential structures),

 \xDH
whether they can be defined by binary relations, or perhaps some abstract
size
property, on the underlying model structure, like preferential or modal
logic,
or some higher relation, like distances between pairs of models, like
for theory revision, update, or counterfactual conditionals,

 \xDH
and if they are based on binary relations, according to the properties of
those
relations, like transitivity, smoothness, rankedness, etc.,

 \xDH
whether we have a value function on the models, like utility, whether we
have
addition, and similar operations on these values,

 \xDH
how the semantical structures are evaluated, e.g., both preferential and
modal structures work with binary relations, but the relations are
evaluated
in totally different ways,

 \xDH
whether we work with an unstructured language, or not, is there additional
structure on the truth values, or not?,

 \xEj

Thus, also an exhaustive classification and listing of ordering principles
seems
quite hopeless to obtain.

In addition, for some conditionals, especially for modal logic, there
is abundant literature, and we neither have exhaustive knowledge,
nor do we think it important to summarize this literature here.

So, what will we do?
Our purpose here is to begin to put some order into this multitude, or,
perhaps
even better, lay down some lines along which ordering is possible.
\subsection{
Overview of this chapter
}
\subsubsection{
Definition and classification
}

As argued in the introduction to this chapter,
the best seems to be to say that a conditional is just any operator.
Negation, conjunction, etc., are then included, but excluded from the
discussion, as we know them well.

The classical connectives have a semantics in the boolean set operators,
but
there are other operators, like the $ \xbm -$functions of preferential
logic which
do not correspond to any such operator, and might even not preserve
definability
in the infinite case (see
Definition \ref{Definition Log-Base} (page \pageref{Definition Log-Base})). It
seems more promising to order
conditionals
by the properties of their model choice functions, e.g., whether those
functions
are idempotent,
etc., see Section \ref{Section Choice-Prop} (page \pageref{Section Choice-Prop})
.

Many conditionals can be based on binary relations, e.g. modal
conditionals
on accessibility relations, preferential consequence relations on
preference relations, counterfactuals and theory revision on distance
relations, etc. Thus, it is promising to look at those relations, and
their
properties to bring more order into the vast field of conditionals.
D.Gabbay introduced reactive structures
(see, e.g.,  \cite{Gab04}), and added supplementary expressivity to
structures based on binary relations, see
 \cite{GS08b} and  \cite{GS08f}. In particular, it
was shown there that we can have cumulativity without the
basic properties of preferential structures (e.g., OR).
This is discussed in
Section \ref{Section Cond-Bin} (page \pageref{Section Cond-Bin}).
\subsubsection{
Additional structure on language and truth values
}

Normally, the language elements (propositional variables) are not
structured. This is somewhat surprising, as, quite often, one variable
will be more important than another. Size or weight might often be more
important than colour for physical objects, etc. It is probably the
mathematical
tradition which was followed too closely. One of the authors gave a
semantics
to theory revision using a measure on language elements in
 \cite{Sch91-1} and  \cite{Sch91-3},
but, as far as we know, the subject was not treated in a larger
context so far. The present book often works with independence of language
elements, see in particular
Chapter \ref{Chapter Mod-Mon-Interpol} (page \pageref{Chapter Mod-Mon-Interpol})
 and
Chapter \ref{Chapter Size-Laws} (page \pageref{Chapter Size-Laws}), and Hamming
type relations and
distances
between models, where it need not be the case that all variables have
the same weight. Thus, it is obvious to discuss this subject in the
present text.
It can also be fruitful to discuss sizes of subsets of the set of
variables,
so we may, e.g., neglect differences to classical logic if they concern
only
a ``small'' set of propositional variables.

On the other hand, classical truth values have a natural order,
$FALSE<TRUE,$ and we will sometimes work with more than 2 truth values,
see in particular
Chapter \ref{Chapter Mod-Mon-Interpol} (page \pageref{Chapter Mod-Mon-Interpol})
, but also
Section \ref{Section EQ} (page \pageref{Section EQ}). So there is a natural
question: do we also
have a total
order, or a boolean order, or another order on those sets of truth values?
Or: Is there a distance between truth values, so that a change from value
$ \xCf a$ to
value $b$ is smaller than a change from $ \xCf a$ to $c?$

There is a natural correspondence between semantical structures and truth
values, which is best seen by an example: Take finite (intuitionistic)
Goedel
logics,
see Section \ref{Section Finite-Goedel} (page \pageref{Section Finite-Goedel}),
say, for simplicity with two worlds. Now, $ \xbf $ may hold nowhere,
everywhere, or only in the second world (called
``there'', in contrast to ``here'',
the first world). Thus, we can express the same situation by three truth
values: 0 for nowhere, 1 for only
``there'', 2 for everywhere.

In Section \ref{Section Softening} (page \pageref{Section Softening}), we will
make some short remarks on
``softening'' concepts, like neglecting
``small'' fragments of a language, etc. This way, we can define, e.g.,
``soft'' interpolation, where we need a small set of variables which
are not in both formulas.

Inheritance systems,
(see, e.g.,  \cite{TH89},  \cite{THT86},
 \cite{THT87},  \cite{TTH91},  \cite{Tou86},
also  \cite{Sch93} and  \cite{Sch97}),
present many aspects of independence,
(see Section \ref{Section Inher} (page \pageref{Section Inher})).
Thus, if two nodes are not connected by valid paths, they may have
very different languages, as language elements have to be inherited,
otherwise, they are undefined. In addition, a may inherit from $b$
property $c,$
but not property $d,$ as we have a contradiction to $d$ (or, even $ \xCN
d)$ via
a different node $b'.$ Theses are among the aspects which make them
natural,
but also quite different from traditional logics.
\subsubsection{
Representation for general revision, update, and counterfactuals
}

Revision (see  \cite{AGM85}, and the discussion in
Section \ref{Section Mod-TR} (page \pageref{Section Mod-TR})),
update (see  \cite{KM90}),
and counterfactuals
(see  \cite{Lew73} and  \cite{Sta68})
are special forms of conditionals, which
received much interest in the artificial intelligence community.
Explicitly or implicitly
(see  \cite{LMS95},  \cite{LMS01}),
they are based on a distance based semantics, working with
``closest worlds''.
In the case of revision, we look at those worlds which are closest to the
present $ \xCf set$ of worlds, in update and counterfactual, we look
from each present world $ \xCf individually$ to the closest worlds, and
then take
the union. Obviously, the formal properties may be very different in the
two
cases.

There are two obvious generalizations possible, and sometimes necessary.
First, ``closest'' worlds need not exist, there may be infinite descending
chains of distances without minimal elements. Second, a distance or ranked
order may force too many comparisons, when two distances or elements may
just simply not be comparable. We address representation problems for
these generalizations:

 \xEh

 \xDH
We first generalize the notion of distance for revision semantics
in Section \ref{Section Semantic-TR} (page \pageref{Section Semantic-TR}).
We mostly consider symmetrical distances, so $d(a,b)=d(b,a),$ and we
work with equivalence classes $[a,b].$ Unfortunately, one of the main
tools in  \cite{LMS01}, a loop condition, does not work any more, it
is
too close to rankedness.

We will have to work more in the spirit of general and
smooth preferential structures to obtain representation. Unfortunately,
revision does not allow many observations
(see  \cite{LMS01}, and, in particular, the impossibility results
for
revision (``Hamster Wheels'') discussed in
 \cite{Sch04}), so all we have
(see Section \ref{Section TR-Main-Rep} (page \pageref{Section TR-Main-Rep}))
are results which use more
conditions than
what can be observed from revision observations.
This problem is one of principles: we showed in
 \cite{GS08a}, see also  \cite{GS08f}, that cumulativity
suffices only
to guarantee smoothness of the structure if the domain is
closed under finite unions. But the union of two products need not
be a product any more.

To solve the problem, we use a technique employed in
 \cite{Sch96-1}, using ``witnesses'' to
testify for the conditions.

 \xDH
We then discuss the limit version (when there are no minimal distances)
for theory revision.

 \xDH
In Section \ref{Section Semantic-Up} (page \pageref{Section Semantic-Up}), we
turn to
generalized update and counterfactuals.
To solve this problem, we use a technique invented in
 \cite{MS90}, and adapt it to our situation.
The basic idea is very simple: we begin (simplified) with some world $x,$
and
arrange the other worlds around $x,$ as $x$ sees them, by their relative
distances.
Suppose we consider now one those worlds, say $y.$ Now we arrange the
worlds
around $y,$ as $y$ sees them. If we make all the new distances smaller
than the
old ones, we ``cannot look back'', etc. We continue this construction
unboundedly
(but finitely) often. If we are a little careful, everyone will only see
what he
is supposed to see. In a picture, we construct galaxies around a center,
then
planets around suns, moon around planets, etc.

The resulting construction is an $ \xda -$ranked structure, as discussed
in
 \cite{GS08d}, see also  \cite{GS08f}.

 \xDH
In Section \ref{Section TR-Up-Synt} (page \pageref{Section TR-Up-Synt}), we
discuss the corresponding
syntactic
conditions, using again ideas from
 \cite{Sch96-1}.

 \xEj
\section{
An abstract view on conditionals
}
\subsection{
A general definition as arbitrary operator
}

\bd

$\hspace{0.01em}$


\label{Definition Sem-Cond}

(1) The 2-valued case:

Let $M$ be the set of models for a given language $ \xdl.$

(1.1) Single sets:

A $n-$ary semantical conditional $C$ is a $n-$ary function

$C: \xdp (M) \xCK  \Xl  \xCK \xdp (M) \xcp \xdp (M).$

(1.2) Systems of sets:

A $n-$ary semantical sys-conditional $C$ is a $n-$ary function

$C: \xdp (M) \xCK  \Xl  \xCK \xdp (M) \xcp \xdp (\xdp (M)).$

(2) The many-valued case:

Let $V$ be the set of truth values, $L$ the set of propositional variables
of $ \xdl,$
let $M$ be the set of functions $m:L \xcp V$ (such $m$ are the
generalizations of a
2-valued model). Let $ \xdm $ be the set of functions $F:M \xcp V$ (such
$F$ are the
generalization of a model set).

(2.1) Single ``sets'':

A $n-$ary many-valued semantical conditional $C$ is a $n-$ary function

$C: \xdm \xCK  \Xl  \xCK \xdm \xcp \xdm.$

(2.2) Systems of ``sets'':

A $n-$ary many-valued semantical sys-conditional $C$ is a $n-$ary function

$C: \xdm \xCK  \Xl  \xCK \xdm \xcp \xdp (\xdm).$

Note that the definition is not yet fully general, as we might have
sets of pairs of elements as a result, e.g., the pairs with minimal
distance from $X \xCK Y$ - but this might cause more confusion than
clarity.

When the context is clear, we will just speak of conditionals, without
any further precision.

\ed

\be

$\hspace{0.01em}$


\label{Example Sem-Cond}

(1) Negation (complement) $C:M(\xbf) \xcZ M(\xCN \xbf)$ is an unary
2-valued
conditional.

(2) AND (intersection) $C: \xBc M(\xbf),M(\xbq) \xBe  \xcZ M(\xbf \xcu \xbq
)$
is a binary 2-valued
conditional.

(3) Given a preferential relation on $M,$ defining the sets of minimal
elements $ \xbm (X)$ of $X,$ $ \xbm $ is an unary 2-valued conditional.

(4) Given a preferential relation on $M,$ defining the systems of
$MISE$'s, $MISE:X \xcZ \{Y \xcc X:Y$ is MISE in $X\}$ is an unary
2-valued sys-conditional.

(5)
Often, the same underlying structure can define a simple or a system
conditional: A preferential relation generates the set of minimal
elements,
and the system of MISE.

(6)
The same underlying structure may also be used totally differently, e.g.,
a binary relation may be used to find the minimal or the accessible
elements.
\subsection{
Properties of choice functions
}

\label{Section Choice-Prop}

\ee

We may classify conditionals by the abstract properties of their choice
functions:

 \xEh

 \xDH
The $ \xbm -$functions of general preferential structures obey additive
size laws,
as described in
Table \ref{Table Base2-Size-Rules-1} (page \pageref{Table Base2-Size-Rules-1}) 
and
Table \ref{Table Base2-Size-Rules-2} (page \pageref{Table Base2-Size-Rules-2}).

 \xDH
Cumulative or smooth structures obey additional laws.

 \xDH
Rational structures obey the rule
that $ \xbm (A \xcv B)= \xbm (A) \xcv \xbm (B)$ or $ \xbm (A \xcv B)= \xbm
(A)$ or $ \xbm (A \xcv B)= \xbm (B).$

 \xDH
For the usual preferential structures, we have $ \xbm (\xbm (X))= \xbm
(X).$

 \xDH
Modular structures show multiplicative rules, as described in
Table \ref{Table Mul-Laws} (page \pageref{Table Mul-Laws}).

 \xDH
Update and counterfactual conditionals are additive on the
left: $ \xbm (A \xcv B,X)= \xbm (A,X) \xcv \xbm (B,X).$

 \xEj
\subsection{
Evaluation of systems of sets
}

Note that, even if the set or system of sets (in the 2-valued case) is
fully described, it is not yet fully defined what we do with the
information, not even in the finite (and thus definability preserving)
case:

Given the system of sets, we can
 \xEh
 \xDH
determine what holds in all elements of the set, or finally in the
``good'' sets, see
Section \ref{Section Pref-Lim} (page \pageref{Section Pref-Lim}),
 \xDH
describe exactly those sets, e.g. in:

 \xEI

 \xDH
Boutelier's modal logic approach to preferential structures, see
 \cite{Bou90a},

 \xDH
the ``good'' sets of deontic logic,
see  \cite{GS08f}
and Chapter \ref{Chapter Neighbourhood} (page \pageref{Chapter Neighbourhood}).

 \xEJ

 \xEj

Such systems are not necessarily well described as neighbourhood systems,
and are usually not yet well examined.
\subsection{
Conditionals based on binary relations
}

\label{Section Cond-Bin}

Often, $C$ is determined by a binary relation, and we can try to classify
the conditionals by the properties of the relation and the way the
relation is used.

\be

$\hspace{0.01em}$


\label{Example Bin-C}

(1) Kripke structures for modal, intuitionistic, and other logics,

(2) preferential structures,

(3) distance based theory revision (the relation is on the Cartesian
product),

(4) the Stalnaker-Lewis semantics for counterfactual conditionals, and
distance
based update (again, the relation is on the Cartesian product, but used
differently).

But we can also see the following structures as given by a binary
relation:

(5) Defaults: A default $A \xcp A \xcu B$ (we take only simple examples)
is pure
defense, $A \xcu B-$models are ``better'' than $A \xcu \xCN B-$models, we do
$ \xCf not$ say that $A \xcu B-$models attack $A \xcu \xCN B-$models.
(This view can have
repercussions when chaining defaults.) The ``best'' situations are those
which satisfy a maximum of default rules.
See Chapter \ref{Chapter Neighbourhood} (page \pageref{Chapter Neighbourhood}).

(6) Obligations: An obligation can be seen the same way as a default.
See again Chapter \ref{Chapter Neighbourhood} (page \pageref{Chapter
Neighbourhood}).

(7) Contrary to duty obligations can also be seen as such rules. But here,
the best situations are those which satisfy the primary obligation, etc. -
this is $ \xCf not$ an order by specificity.

\ee

We can then classify conditionals as to the properties of these relations,
like transitivity, smoothness, etc., and their use.
\subsubsection{
Short discussion of above examples
}

Obviously, the relation of accessibility in Kripke models is one of pure
defense: if $ \xCf xRy,$
then $x$ ``supports'' or ``defends'' $y.$
This results in the trivial property

$X \xcc Y$ $ \xch $ $C(X) \xcc C(Y)$

The situation in preferential structures is more complicated:
If $x \xbe X,$ then $x$ ``defends'' itself. There is no defense of another
element,
but there is attack, if $x \xeb y,$ then $x$ attacks $y.$ The situation
is, as a matter of
fact, still more subtle. Preferential structures often work with copies.
Each copy $ \xBc x,i \xBe $ is a defense of $x.$ This becomes obvious, when we
think
that,
to destroy minimality of $x,$ we have to destroy $ \xCf all$ copies of
$x.$

As there is no attack on attacks, attacks have a property of monotony,
what is attacked in $X,$ is also attacked in any $Y \xcd X.$ This results
in:

$X \xcc Y$ $ \xch $ $C(Y) \xcs X \xcc C(X).$

Smoothness means, roughly, that, if there is an attack, there is an attack
from a valid element (or copy).

Higher preferential structures may attack attacks, which is then a
defense,
see Section \ref{Section High-Pref} (page \pageref{Section High-Pref}).

In Example \ref{Example Bin-C} (page \pageref{Example Bin-C}),
(3) and (4), we have a ranking of distances. Thus, any
$ \xBc x,y \xBe $ is a defense of $y$ (when we consider the second coordinate),
and
an
attack on bigger pairs $ \xBc x',y'  \xBe.$

Consequently, we may also classify conditionals on the role of support
and attack:
 \xEI
 \xDH
there are no conflicts, e.g.

we have pure defense, as in Kripke models,

we have pure attack, as in classical preferential structures
 \xDH
conflicts are resolved or not
e.g. by sup of the absolute values of the strength of support and attack,
when
this exists, by addition, if there is $+$ on $Y,$ etc.:
 \xEI
 \xDH
in reactivite diagrams, we may consider the biggest absolute values for
$f(x,y')$
for fixed $y',$ and
$x \xbe A:$ basic arrows are attack of force -1, attack on attack has
force $+2,$ etc.
 \xDH
in inheritance diagrams, we use the valid paths for strength comparison
(see  \cite{GS08e},  \cite{GS08f}),
 \xDH
so far, the use of support and attack does not seem to use ressources,
as is the case in Girard's Linear Logic, it remains an open research
problem
to investigate utility and properties of such approaches.
 \xEJ

 \xEJ
\section{
Conditionals and additional structure on language and truth values
}

\label{Section Mod-Many}

\label{Section Mod-Lang-Struc}
\subsection{
Introduction
}

This section gives only a very rough outline of the possible
approaches.

We can treat the set of language elements (propositional variables in our
context), or of truth values, like any other set. Thus, we can introduce
relations, operations, etc. on those sets. In addition, we can
compose sublanguages into bigger languages, small truth value sets
in bigger one, or, conversely, go from bigger to smaller sets, etc.

This section is mostly intended to open the discussion. What is
reasonable to do, will be shown when necessity arises from applications.
\subsection{
Operations on language and truth values
}

We can consider structures and operations on
 \xEI
 \xDH
language elements (propositional variables) and truth values
within one language
 \xDH
languages and truth value sets for several, different languages
 \xDH
definable model sets.
 \xEJ

We will take now a short look at the different cases.
\subsection{
Operations on language elements and truth values within one language
}
\paragraph{
Operations on language elements
}

 \xEh

 \xDH
A relation of importance between language elements:

 \xEh

 \xDH
If $p$ is more important than $q,$ then we might give
``biased'' weight e.g. in theory revision,
see  \cite{Sch95-2},
 \cite{Sch91-1}, and  \cite{Sch91-3}.

 \xDH
This might lead to constructing a modular order of models, for
non-monotonic logic and theory revision, e.g., by
a lexicographic order, where the more important language elements have
precedence.
Suppose for instance that $ \xCf a$ is considered more important than $b,$
but that we
prefer $ \xCf a$ and $b$ over their negation. This might lead to the
following order: $ \xCN a \xCN b< \xCN ab<a \xCN b<ab.$

 \xEj

 \xDH
A distance relation between language elements

 \xEh

 \xDH
See Section \ref{Section Inher} (page \pageref{Section Inher}), for the
structuring of the language in
inheritance diagrams via a distance between language elements.

 \xEj

 \xDH
Size of sets of language elements:

 \xEh

 \xDH
A notion of ``small'' sets of language elements allows to introduce
``soft'' concepts (see also below,
Section \ref{Section Softening} (page \pageref{Section Softening})):
 \xEh
 \xDH
Soft cumulativity: for each model, there is a model which differs only
on a small set of variables, and which is minimized by a minimal
model.
 \xDH
Soft modularity and independence: being modular, independent
with some exceptions.
 \xEj

 \xDH
We can define model distance by the size of elements by which they differ,
they are ``almost'' identical if this set is small.

 \xDH
If $J$ is a small subsets of $L$ (the set of propositional variables),
then we have (in classical logic) $2^{J}=2^{small}$ models, so we can work
with rules
about abstract exponentiation.

 \xEj

 \xDH
Accessibility of language elements:

 \xEh
 \xDH Consider inheritance systems: in the Nixon diamond, it is undecided
if Nixon
is a
pacifist, this can now be distinguished from the case where there no
path whatsoever from $ \xCf A$ to $B:$ language element $B$ (pazifist) is
not reachable
from
language element $ \xCf A$ (Nixon).

 \xDH Consider natural language and argumentation: constructing an
argument or
a scenario, adding language elements when necessary, and thus keeping
representation simple and efficient, introducing only what is needed.
 \xEj

 \xDH
Combination of language elements:

We may combine several language elements to one
``super-element'' and thus achieve abstraction, like all properties
for colour are grouped under ``colour''.
Conversely, we may differentiate one element into sub-elements, see also
``bubble structures'',
Section \ref{Section New-Min} (page \pageref{Section New-Min}).

 \xDH
Context dependency:

All above considerations may be context dependent in more complicated
situations, as in inheritance or argumentation.

 \xEj
\paragraph{
Operations on truth values
}

 \xEh

 \xDH
The classical truth value set is a Boolean algebra. Finite Goedel logics
(see Section \ref{Section Finite-Goedel} (page \pageref{Section Finite-Goedel})
) can be seen as having a
linearly ordered truth
value set, see below, (5).

 \xDH
We can consider more general Boolean algebras, or just partial
orders. In the case of argumentation, we may just consider the powerset of
the arguments as the set of truth values - together with the usual
operations.

 \xDH
We may combine several truth values to one
``super-value'' and thus achieve abstraction, e.g., if
the truth values are arguments. For instance, all arguments from one
source may be grouped together.
Conversely, we may differentiate one truth value into sub-values.

 \xDH
We may define a notion of size on the set of truth values.

 \xDH
Note that we have an equivalence between structure on the model set and
structure on the truth value set:

Consider a Kripke structure. We have an equivalent re-formulation when we
look at the same structure on the truth value set. $ \xbf $ is true in
world $m,$
iff the truth value $m$ for $ \xbf $ is true. Intuitionistic logic is then
a
condition about truth values, not about Kripke structures.
(See also Section \ref{Section Finite-Goedel} (page \pageref{Section
Finite-Goedel}).)

These connections could be formally expressed by isomorphy propositions.

It might be reasonable to work sometimes with mixed structures. If a
substructure is repeated everywhere, we can put it into the truth values,
so every point of the simplified structure is differentiated by the
truth value structure. This is, abstractly, similar to the
``bubble structures'' introduced in
Section \ref{Section New-Min} (page \pageref{Section New-Min}), with identical
bubbles.

 \xEj
\subsection{
Operations on several languages
}
\paragraph{
Operations on language elements \\[2mm]
}

 \xEh

 \xDH
We may compose a bigger language from sub-languages, or, conversely,
split bigger languages into sub-languages.

 \xDH
The projection (respectively the inf and sup operators of
Section \ref{Section Mon-Interpol-Int} (page \pageref{Section Mon-Interpol-Int})
) are new operators on model sets,
but also
on the language, going from a formula in a richer language to one in
a poorer language.

 \xDH
In argumentation and inheritance, we can introduce new language elements:
``just think of  \Xl.'', or by establishing a valid path upward from a
given
node in inheritance networks.

 \xDH
Note that the structures on a language $L$ need not be coherent with
structures
on $L' \xcc L,$ see Chapter \ref{Chapter Size-Laws} (page \pageref{Chapter
Size-Laws}).

 \xEj
\paragraph{
Operations on truth values \\[2mm]
}

In defeasible inheritance, we can consider accessible nodes as truth
values
(with valid paths for comparison). Thus, we $ \xCf construct$ the truth
value
structure while evaluating the net above a certain point. Likewise, we may
consider arguments for a formula as truth values, which are then
constructed
during the process.
\subsection{
Operations on definable model sets
}

We have to distinguish:

 \xEh
 \xDH Internal operators of the language, like
the classical operators AND, they have their natural interpretation in
boolean set operators.

 \xDH
External operators like the choice operator of preferential
structures $ \xbm.$

 \xDH
Algebraic operators, like
projection (respectively the inf and sup operators of
Section \ref{Section Mon-Interpol-Int} (page \pageref{Section Mon-Interpol-Int})
),
guaranteeing syntactic interpolation.

 \xEj
\subsection{
Softening concepts
}

\label{Section Softening}

There are several ways we can
``soften'' a concept.

 \xEh

 \xDH
If we have a many-valued logic, we can go down from maximal truth.

 \xDH
We can sometimes replace rules by weaker versions, e.g., $ \xbf \xcn \xbq
$ by $ \xbf \xcN \xCN \xbq.$

 \xDH
We can neglect small fragments of the language, e.g.,
 \xEI
 \xDH
instead of having
full interpolation like $ \xbf \xcn \xba \xcn \xbq $ with $ \xba $ having
only symbols common to
both $ \xbf $ and $ \xbq,$ $ \xba $ is now allowed to have
``some, but not too many'' symbols which are not common to both $ \xbf $ and
$ \xbq,$
 \xDH
strictly speaking, $ \xba \xcN \xbb,$ but if we consider a big subset of
the language,
and $ \xba' $ is the fragment of $ \xba $ in that language, likewise for
$ \xbb',$ then
$ \xba' \xcn \xbb'.$
 \xEJ

 \xDH
``Soft'' independence, where $ \xbS $ is
``almost'' $ \xbS' \xCK \xbS'',$ in several interpretations:
 \xEI
 \xDH
neglecting a small set of variables, we can write $ \xbS $ as a product

 \xDH
neglecting a small set of sequences, we can write $ \xbS $ as a product,
e.g., when $ \xbS =(\xbS' \xCK \xbS'') \xcv \{ \xbs \}.$

 \xDH
If we have a distance between language elements, then we can express that,
the more language elements are distant from each other, the more they
are independent.

 \xDH
We may have independence for small changes, but not for big ones.

A small change may be to go from
``big'' to ``medium'' size, a big change to
go from ``big'' to ``small'', etc.

 \xEJ

 \xDH
Generally, we can ``soften'' an object by considering another object, which
differs only slightly from the original one (how ever this is measured).
A concept can then be
softened by softening the properties, or the objects involved.

 \xEj

It seems premature to go into details here, which should be motivated by
concrete problems.
\subsection{
Aspects of modularity and independence in defeasible inheritance
}

\label{Section Inher}

Inheritance structures have many aspects, they are discussed in detail
in  \cite{GS08e}, see also  \cite{GS08f}.

We discuss here the following, illustrated in
Diagram \ref{Diagram Independ} (page \pageref{Diagram Independ}), and
Diagram \ref{Diagram Information-Transfer} (page \pageref{Diagram
Information-Transfer}):
the left hand side shows the diagram, the right hand side
the strength of information. E.g., in part (2), the strongest information
available at $ \xCf A$ (except $ \xCf A$ itself), is $B \xcu C.$ But there
might be
exceptions, and the most exceptional situation is $ \xCN B \xcu \xCN C.$

 \xEh

 \xDH
We can see the language as being constructed dynamically.
In (1) $A \xcp C,$ $B \xcp C,$
$ \xCf A$ knows nothing about
$B$ and vice versa, $C$ knows nothing about $ \xCf A$ or $B.$
Giving a truth value ``undecided'' would not be correct.
We are not decided about Nixon's pacifism in the Nixon Diagram,
but we have (contradictory) information about it. Here we have
$ \xCf no$ information.

Consider also Diagram \ref{Diagram Information-Transfer} (page \pageref{Diagram
Information-Transfer}).
As there is no monotonous path whatever between $e$ and $d,$
the question whether $ \xCf e$'s are $ \xCf d$'s or not, or vice versa,
does not even arise.
For the same reason, there is no question whether $ \xCf b$'s are $ \xCf
c$'s, or not.
In upward chaining formalisms, as there is no valid positive path from $
\xCf a$ to $d,$
there is no question either whether $ \xCf a$'s are $ \xCf f$'s or not.

 \xDH
We can see inheritance diagrams
as Kripke structures, where the accessibility relation
is constructed non-monotonically, it is not static.

 \xDH
We can see inheritance diagrams
as constructing a structure of truth values dynamically,
along with the relation of comparison between truth values. This truth
value
structure is not absolute, but depends on the node from which we start. (A
strongest element is always given, the information given directly at each
node. The relation between the others is decided by specificity via valid
paths.)

 \xDH
Information is given independently, $ \xCf A$ might have many normal
properties of $B,$
but not all of them.

In (2) of
Diagram \ref{Diagram Independ} (page \pageref{Diagram Independ})  $A \xcp B,$ $A
\xcp C,$ and the
information is independent. The resulting order is, seen from $ \xCf A:$
$B \xcu C<B \xcu \xCN C< \xCN B \xcu \xCN C,$ $B \xcu C< \xCN B \xcu C<
\xCN B \xcu \xCN C.$

Consider now (3): $A \xcp B \xcp C.$
The resulting order is $B \xcu C<B \xcu \xCN C< \xCN B.$ $ \xCN B$ is not
differentiated any
more, we have no information about it.

Consider again Diagram \ref{Diagram Information-Transfer} (page \pageref{Diagram
Information-Transfer}).
In our diagram, $ \xCf a$'s are $ \xCf b$'s,
but not ideal $ \xCf b$'s,
as they are not
$ \xCf d$'s, the
more specific information from $c$ wins. But they are $ \xCf e$'s, as
ideal $ \xCf b$'s are.
So they are not perfectly ideal $ \xCf b$'s, but as ideal $ \xCf b$'s as
possible. Thus, we
have graded ideality, which does not exist in preferential and similar
structures. In those structures, if an element is an ideal element, it has
all
properties of such,
if one such property is lacking, it is not ideal, and we can't say
anything any
more beyond classical logic.
Here, however, we sacrifice as little normality as possible, it is
thus a minimal change formalism.

 \xDH
Inheritance diagrams also give a distance between language elements,
and thus a structure on the language: In the simplest approach,
if the (valid) path from $ \xCf A$ to $B$ is long, then the language
elements $ \xCf A$ and $B$ are distant.
A better approach is via specificity, as usual in inheritance diagrams.

 \xEj

$ \xCO $

\vspace{10mm}

\begin{diagram}

\label{Diagram Independ}
\index{Diagram Independ}

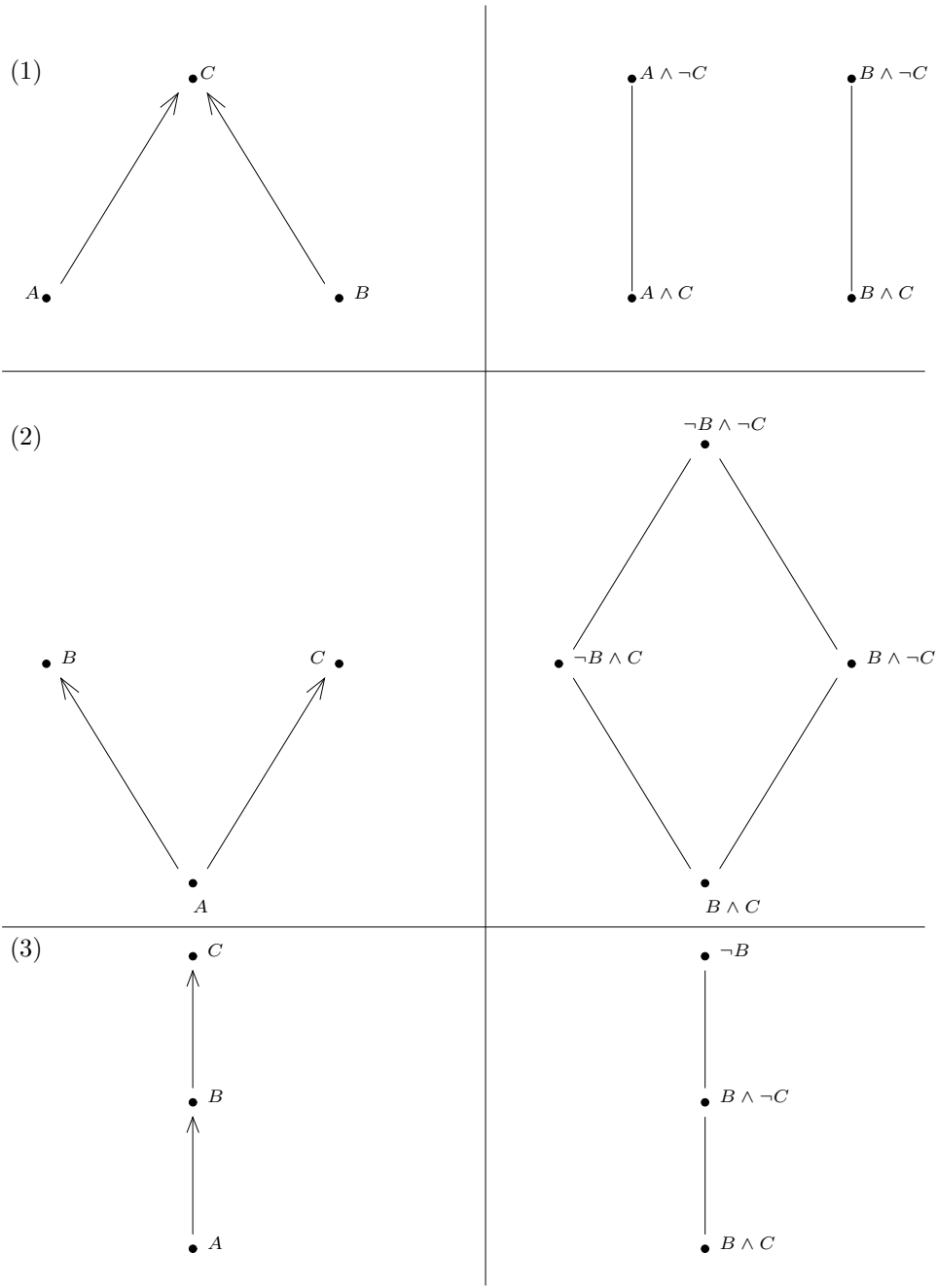
\centering
\setlength{\unitlength}{1mm}
{\renewcommand{\dashlinestretch}{30}
\begin{picture}(150,200)(0,0)

\path(70,5)(70,180)

\put(5,170){$(1)$}

\put(30,170){\circle*{1}}
\put(10,140){\circle*{1}}
\put(50,140){\circle*{1}}

\path(12,142)(28,168)
\path(25.6,166.1)(28,168)(27.4,165.1)
\path(48,142)(32,168)
\path(32.6,165.1)(32,168)(34.4,166.1)

\put(31,170){{\xssc $C$}}
\put(7,140){{\xssc $A$}}
\put(52,140){{\xssc $B$}}

\put(90,140){\circle*{1}}
\put(90,170){\circle*{1}}
\put(120,140){\circle*{1}}
\put(120,170){\circle*{1}}

\path(90,141)(90,169)
\path(120,141)(120,169)

\put(91,140){{\xssc $A \xcu C$}}
\put(91,170){{\xssc $A \xcu \xCN C$}}
\put(121,140){{\xssc $B \xcu C$}}
\put(121,170){{\xssc $B \xcu \xCN C$}}

\path(4,130)(130,130)

\put(5,120){$(2)$}

\put(30,60){\circle*{1}}
\put(10,90){\circle*{1}}
\put(50,90){\circle*{1}}

\path(28,62)(12,88)
\path(12.6,85.1)(12,88)(14.4,86.1)
\path(32,62)(48,88)
\path(45.6,86.1)(48,88)(47.4,85.1)

\put(30,56){{\xssc $A$}}
\put(12,90){{\xssc $B$}}
\put(46,90){{\xssc $C$}}

\put(100,60){\circle*{1}}
\put(80,90){\circle*{1}}
\put(120,90){\circle*{1}}
\put(100,120){\circle*{1}}

\path(98,62)(82,88)
\path(102,62)(118,88)
\path(82,92)(98,118)
\path(118,92)(102,118)

\put(100,56){{\xssc $B \xcu C$}}
\put(82,90){{\xssc $\xCN B \xcu C$}}
\put(122,90){{\xssc $B \xcu \xCN C$}}
\put(97,122){{\xssc $\xCN B \xcu \xCN C$}}

\path(4,54)(130,54)

\put(5,50){$(3)$}

\put(30,10){\circle*{1}}
\put(30,30){\circle*{1}}
\put(30,50){\circle*{1}}

\path(30,12)(30,28)
\path(29.2,25.9)(30,28)(30.8,25.9)
\path(30,32)(30,48)
\path(29.2,45.9)(30,48)(30.8,45.9)

\put(32,10){{\xssc $A$}}
\put(32,30){{\xssc $B$}}
\put(32,50){{\xssc $C$}}

\put(100,10){\circle*{1}}
\put(100,30){\circle*{1}}
\put(100,50){\circle*{1}}

\path(100,12)(100,28)
\path(100,32)(100,48)

\put(102,10){{\xssc $B \xcu C$}}
\put(102,30){{\xssc $B \xcu \xCN C$}}
\put(102,50){{\xssc $\xCN B$}}

\end{picture}
}

\end{diagram}

\vspace{4mm}

$ \xCO $

$ \xCO $

\vspace{10mm}

\begin{diagram}

\label{Diagram Information-Transfer}
\index{Diagram information transfer}

\index{information transfer}

\unitlength1.0mm
\begin{picture}(130,110)

\newsavebox{\Tweety}
\savebox{\Tweety}(130,110)[bl]
{

\put(57,18){\vector(1,1){24}}
\put(51,18){\vector(-1,1){24}}
\put(27,51){\vector(1,1){24}}
\put(81,51){\vector(-1,1){24}}

\put(67,61){\line(1,1){4}}

\put(81,47){\vector(-1,0){54}}

\put(24,51){\vector(0,1){22}}
\put(54,81){\vector(0,1){22}}

\put(53,16){$a$}
\put(23,46){$b$}
\put(83,46){$c$}
\put(53,76){$d$}
\put(23,76){$e$}
\put(53,106){$f$}

\put(10,0) {{\rm\bf Information transfer}}

}

\put(0,0){\usebox{\Tweety}}
\end{picture}

\end{diagram}

\vspace{4mm}

$ \xCO $
\section{
Representation for general revision, update, and counterfactuals
}

\label{Section Repr-TR-Up}
\subsection{
Importance of theory revision for general structures, reactivity, and its
solution
}

\label{Section TR-Importance}

Theory revision is about minimal change.
When we revise a theory with a formula, we obtain a new theory (in the
AGM approach). When we want to revise other objects, like one logic with
another logic, we may face the following problems:

 \xEh
 \xDH
We have to choose an adequate level of representation:
For instance, if we revise one logic with another logic, do we
represent the logics
 \xEh
 \xDH syntactically, as a relation between formulas, or
 \xDH semantically, e.g. by preferential structures, or
 \xDH by their abstract semantics, e.g., the sets of big subsets in
preferential logics?
 \xEj
When we chose a syntactic representation, we have to make the
approach robust under syntactic reformulations (as in AGM revision).
The same holds when we choose representation by preferential
structures, as the same logic can be represented by many different
structures,
(for instance, by playing with the number of copies). We can, of course,
choose a canonical representation, but this risks to be artificial.
Probably the best approach here is to take the abstract semantics with
big subsets.

 \xDH
Traditional revision has a limit condition, by which (semantically
speaking) closest models will always exist. It is not at all clear that
such
a condition will always be satisfied. For instance, there will often not
be
a closest real-valued continuous function with certain properties, but
only
closer and closer ones. Thus, we will have to take a limit approach,
instead of
the traditional minimal approach.

 \xDH
It will often be too strong to require a ranked order on the elements, as
AGM revision presupposes. One object may be closer in one aspect, another
object in another object, but this closeness is not comparable. So, in
general,
we will only have a partial order.

 \xDH
The closest object may not have the desired properties. E.g., we take
transitive relation, modify it minimally by adding one pair to the
relation, then the relation will not always be transitive any more.
Thus, structural properties might not be preserved.

 \xEj

Due to the importance of theory revision, we will now give new
representation
results, which address some of these problems.
\subsection{
Introduction
}

We work here on representation for not necessarily ranked orders for
distance based revision, and for update/counterfactual conditionals.

We do first the general, and then the smooth case for revision, then
the same for update/counterfactual conditionals.
We also examine the limit versions of revision and
update/counterfactual conditionals.

Recall from Table \ref{Table Base1-Pref-Rep} (page \pageref{Table
Base1-Pref-Rep})  the semantic
representation results
for general preferential structures. We will use them now. The
complication
of our situation lies in the fact that we cannot directly observe the
``best'' elements (which are now pairs), but only their projections. In the
case
of ranked structures, representation as shown in
 \cite{LMS01} was easier, as, in this case, any minimal element was
comparable with
any non-minimal element, but this is not the case any more.
\subsection{
Semantic representation for generalized distance based theory revision
}

\label{Section Semantic-TR}
\subsubsection{
Description of the problem
}

We characterized distance based revision in  \cite{LMS01}.
Our aim is to generalize the result to more general notions of distance,
where the abstract distance is just a partial order, and also to cases
where the limit condition, i.e., there are always closest elements,
may be violated.

The central characterizing condition in  \cite{LMS01} was an elegant
loop
condition, due to M.Magidor. Unfortunately, this condition seems to
be closely related to the rankedness condition, which we do not have
any more. (More precisely, any minimal element was comparable to any
non-minimal element, now a minimal and a non-minimal element need not be
comparable. Even in the smooth case, a non-minimal element has to be
comparable
only to $ \xCf one$ minimal element, but not to all minimal elements.) We
do not see
how to find a similar nice condition in our more general situation.
(The loop condition would have to be modified anyway, but this is not the
main
problem:
$(X \xfA (Y \xcv Y')) \xcs Y \xEd \xCQ $ has to be changed to:
$(X \xfA (Y \xcv Y')) \xcs Y' = \xCQ $ and $(X \xfA (Y \xcv Y')) \xcc
Y.)$

We treated in the past a problem similar to the present one in  \cite{Sch96-1}.
There, the existence of (at least) one minimal comparable element was
captured
by the idea of ``witnesses'', described by a formula
$ \xbf (\xBc a,b \xBe,A \xCK B),$ which expresses
that $ \xBc a,b \xBe $ is a valid candidate not only for $A \xCK B$ but also for
all
$A' \xCK B',$ $A' \xcc A,$ $B' \xcc B.$ This gives the main condition for
preferential representation. Another paper treating somehow
similar problems was  \cite{BLS99}, where we worked with
``patches''.

The idea pursued here is the same as in  \cite{Sch96-1}.
We only isolate the different layers of
the problem better, so our approach is more flexible for various
conditions
imposed on the abstract distance. This allows us to treat many variants
in a very general way, mostly just citing the general representation
results,
and putting things together only at the end.

 \xEh
 \xDH
Thus, we first give a general
machinery for treating symmetrical distances by a translation to
equivalence classes.
See Section \ref{Section TR-Equiv} (page \pageref{Section TR-Equiv}).
 \xDH
We then consider representation for the equivalence
classes, with various additional conditions like
smoothness/semantical cumulativity, translating them back to conditions
about the original pairs $ \xBc a,b \xBe.$
See Section \ref{Section TR-Main-Rep} (page \pageref{Section TR-Main-Rep}).
 \xDH
In the next step, we see how much we can
do for the original problem, where we do not have arbitrary $A \xcc U \xCK
U$ $(\xCf U$ is the
universe we work in), but only sets of the type $A \xCK B,$ $A,B \xcc U.$
In particular,
the union of $A \xCK B$ and $A' \xCK B' $ need not be again such a
product.
Moreoever, we can
observe only the projection of $ \xbm (A \xCK B)$ onto the second
coordinate:
$\{b \xbe B: \xcE  \xBc a,b \xBe  \xbe \xbm (A \xCK B)\},$
see Section \ref{Section TR-Proj} (page \pageref{Section TR-Proj}).
More precisely, we look for $ \xbm (Z)_{X}:=\{x \xbe X: \xcE y \xbe
Y. \xBc x,y \xBe  \xbe \xbm (Z)\}$ and
$ \xbm (Z)_{Y}:=\{y \xbe Y: \xcE x \xbe X. \xBc x,y \xBe  \xbe \xbm (Z)\}$.
Of course, $x \xbe \xbm (Z)_{X}$, $y \xbe \xbm (Z)_{Y}$ does not imply
$ \xBc x,y \xBe  \xbe \xbm (Z).$
We only know that for all $x \xbe \xbm (Z)_{X}$ there is at least one $y
\xbe \xbm (Z)_{Y}$ such that
$ \xBc x,y \xBe  \xbe \xbm (Z),$ and conversely.
Thus, in the original problem, we have less input, and also less
observation
of the output. The implications of scarce information for representation
was,
for example, illustrated in
 \cite{Sch04}, by
``Hamster Wheels'', where we showed that this might make finite
characterizations impossible.
 \xDH
In a final step, we will go from the semantic
side to the descriptive, logical side.
See Section \ref{Section TR-Up-Synt} (page \pageref{Section TR-Up-Synt}).
 \xEj
\subsubsection{
Symmetrical distances and translation to equivalence classes
}

\label{Section TR-Equiv}

We work in some universe $U.$ We want to express that the distance from $
\xCf a$ to $b$
is the same as the distance form $b$ to $ \xCf a.$ So, when pairs stand
for distance,
$ \xBc a,b \xBe $ and $ \xBc b,a \xBe $ will be equivalent, and we consider
equivalence
classes
$[a,b]$ instead of pairs $ \xBc a,b \xBe.$
The $[a,b]$'s will be the abstract distances, not the $ \xBc a,b \xBe $'s.

We define

\bd

$\hspace{0.01em}$


\label{Definition A-equiv}

Let $A,A' \xcc U \xCK U.$

(1) $A \xCq A' $ iff $ \xcA  \xBc a,b \xBe  \xbe A(\xBc a,b \xBe  \xbe A' $ or
$ \xBc b,a \xBe  \xbe A'
)$ and
$ \xcA  \xBc a,b \xBe  \xbe A' (\xBc a,b \xBe  \xbe A$ or $ \xBc b,a \xBe  \xbe
A).$

(2) $[A]:=\{[a,b]:$ $ \xBc a,b \xBe  \xbe A$ or $ \xBc b,a \xBe  \xbe A\}$

\ed

Note that we might lose elements when going from $ \xCf A$ to $[ \xCf A],$
precisely
when both $ \xBc a,b \xBe  \xbe A$ and $ \xBc b,a \xBe  \xbe A,$ then both are
in the class
$[a,b].$

We have

\bfa

$\hspace{0.01em}$


\label{Fact A-Equiv}

$A \xCq A' $ iff $[A]=[A' ].$

\efa

\subparagraph{
Proof
}

$\hspace{0.01em}$


Trivial. $ \xcz $
\\[3ex]

We suppose a choice function $ \xbm $ to be given, chosing pairs of
minimal
distance, more precisely
$ \xbm (A):=\{ \xBc a,b \xBe  \xbe A:$ $ \xcA  \xBc a',b'  \xBe  \xbe A.d(a,b)
\xck d(a',b'
)\}.$
We try to represent $ \xbm $ by a suitable ordering relation $ \xeb $ on
pairs.
As the distance relation is supposed to be symmetric, we will work with
equivalence classes, considering $ \xbm' ([A]),$ instead of $ \xbm (A),$
and
representing $ \xbm' $ instead of $ \xbm.$ We have to define $ \xbm' $
from $ \xbm,$ and make sure
that the definition is independent of the choice of the particular element
we work with. For this to work, we
consider the following axiom $(\xbm S1)$ about symmetry of distance. We
also
add immediately an axiom $(\xbm S2)$ which says that, if both
$ \xBc a,b \xBe $ and $ \xBc b,a \xBe $ are present, then both or none are
minimal.

\bd

$\hspace{0.01em}$


\label{Definition Mu-Symm}

We define two conditions about symmetrical $ \xbm:$

$(\xbm S1)$ Let $A,A' \xcc U \xCK U.$

If $A \xCq A',$ then $\{[a,b]:$ $ \xBc a,b \xBe  \xbe \xbm (A)$ or $ \xBc b,a
\xBe  \xbe \xbm
(A)\}$ $=$
$\{[a,b]:$ $ \xBc a,b \xBe  \xbe \xbm (A')$ or $ \xBc b,a \xBe  \xbe \xbm (A'
)\}.$

$(\xbm S2)$ Let $A \xcc U \xCK U,$ and both $ \xBc a,b \xBe, \xBc b,a \xBe 
\xbe A.$

Then $ \xBc a,b \xBe  \xbe \xbm (A)$ iff $ \xBc b,a \xBe  \xbe \xbm (A).$

Axiom $(\xbm S1)$ says that it is unimportant for $ \xbm $ in which form
$[a,b]$ is
present, as $ \xBc a,b \xBe,$ $ \xBc b,a \xBe,$ or as both.

\ed

The main definition is now:

\bd

$\hspace{0.01em}$


\label{Definition Mu'}

$ \xbm' ([A]):=\{[a,b] \xbe [A]:$ $ \xBc a,b \xBe  \xbe \xbm (A)$ or $ \xBc b,a
\xBe  \xbe
\xbm (A)\}.$

\ed

We have to show that $ \xbm' $ is well-defined, i.e. the following Fact
holds.

\bfa

$\hspace{0.01em}$


\label{Fact Mu'-Well-Def}

If $(\xbm S1)$ holds, and $[A]=[B],$ then $ \xbm' ([A])= \xbm' ([B]).$

\efa

\subparagraph{
Proof
}

$\hspace{0.01em}$


If $[A]=[B],$ then by Fact \ref{Fact A-Equiv} (page \pageref{Fact A-Equiv})  $A
\xCq B,$ so by
$(\xbm S1)$ $\{[a,b]:$ $ \xBc a,b \xBe  \xbe \xbm (A)$ or $ \xBc b,a \xBe  \xbe
\xbm (A)\}$
$=$
$\{[a,b]:$ $ \xBc a,b \xBe  \xbe \xbm (B)$ or $ \xBc b,a \xBe  \xbe \xbm (B)\},$
and by
Definition \ref{Definition Mu'} (page \pageref{Definition Mu'})  $ \xbm' ([A])=
\xbm' ([B]).$ $ \xcz $
\\[3ex]

Suppose now we have a representation result for $ \xbm',$ i.e. some
relation $ \xeb $ such that (simplified, without copies):

$ \xbm' ([A])$ $=$ $\{[a,b] \xbe [A]:$ $ \xCN \xcE [a',b' ] \xbe [A].[a'
,b' ] \xeb [a,b]\}.$

Recall that $[a,b]$ is the (abstract) symmetrical distance between $ \xCf
a$ and $b.$
Then the relation $ \xeb $ describes those pairs $ \xBc a,b \xBe  \xbe A$ which
have
minimal distance, more precisely, all $ \xBc a,b \xBe $ or $ \xBc b,a \xBe $
such that
$[a,b] \xbe [A].$ So our representation problem is solved.

We may want to impose additional conditions.
They may come from $ \xbm',$ where, e.g., we may want semantic
cumulativity,
they have then to be translated back to conditions about $ \xbm.$
They may come from requirements about the distance, e.g.,
that all distances from $x$ to $x$ are minimal, they have to be translated
to conditions for $ \xbm'.$
We may also have conditions about the domain, e.g., we do not consider
all $A \xcc U \xCK U,$ but only some, so we work on some $ \xdy \xcc \xdp
(U \xCK U),$ which may have
some conditions like closure under finite unions - or, conditions on the
sets of equivalence classes to be considered.
Such additional conditions will be considered below.
\subsubsection{
General representation results and their translation
}

\label{Section TR-Main-Rep}
\paragraph{
The general and smooth case \\[2mm]
}

We continue to work in $U,$ but with the equivalence classes. Let
$[U]:=\{[a,b]:a,b \xbe U\},$ and $ \xdy \xcc [U].$
For the smooth case, let $ \xdy $ be closed under finite unions.

Consider $ \xbm': \xdy \xcp \xdp ([U]).$

We know from previous work, see
Table \ref{Table Base1-Pref-Rep} (page \pageref{Table Base1-Pref-Rep})
that

 \xEh

 \xDH
such $ \xbm' $ can be represented by a preferential structure iff

$(\xbm \xcc)$ $ \xbm' (X) \xcc X$

and

$(\xbm PR)$ $X \xcc Y$ $ \xch $ $ \xbm' (Y) \xcs X \xcc \xbm' (X)$

hold (the structure can be chosen transitive),

 \xDH
if $ \xdy $ is closed under finite unions,
such $ \xbm' $ can be represented by a preferential structure iff,
in addition

$(\xbm CUM)$ $ \xbm' (X) \xcc Y \xcc X$ $ \xch $ $ \xbm' (X)= \xbm'
(Y)$

holds. (Again, the structure can be chosen transitive.)

 \xEj

We have to consider the translation of the conditions back from $ \xbm' $
to
$ \xbm,$ where $ \xbm' $ is defined from $ \xbm $ as in
Definition \ref{Definition Mu'} (page \pageref{Definition Mu'}).

$(\xbm \xcc):$ $ \xbm' ([X]) \xcc [X]$ does not necessarily entail $
\xbm (X) \xcc X,$ as, e.g.,
the possibility
$ \xBc a,b \xBe  \xbe X,$ $ \xBc b,a \xBe  \xce X,$ but $ \xBc b,a \xBe  \xbe
\xbm (X)$ is left open.
But the converse is true, and
natural, so we impose $ \xbm (X) \xcc X.$ Note, however, that this
observation destroys
full equivalence of the conditions for $ \xbm,$ we have to bear in mind
that they
only hold for $ \xbm'.$

$(\xbm PR)$ for $ \xbm $ entails $(\xbm PR)$ for $ \xbm':$

Let $ \xbm' [X]:= \xbm' ([X]).$

We have to show $[X] \xcc [Y]$ $ \xch $ $ \xbm' [Y] \xcs [X] \xcc \xbm'
[X].$

Let $[X] \xcc [Y].$
Consider $X$ and $Y$ such that $X \xcc Y.$ This is possible, e.g., add
inverse pairs
to $Y$ as necessary, so by prerequisite $ \xbm (Y) \xcs X \xcc \xbm (X).$
Note that we used
here some prerequisite about domain closure of the original set $ \xdy.$
Let $[a,b] \xbe \xbm' [Y] \xcs [X],$ so e.g.
$ \xBc a,b \xBe  \xbe \xbm (Y) \xcc Y,$ and
$ \xBc a,b \xBe  \xbe X$ or $ \xBc b,a \xBe  \xbe X.$ In the latter case,
$ \xBc b,a \xBe  \xbe Y,$ and by $(\xbm S2)$ also $ \xBc b,a \xBe  \xbe \xbm
(Y).$ So let
without loss of generality
$ \xBc a,b \xBe  \xbe \xbm (Y) \xcs X \xcc \xbm (X),$ so $[a,b] \xbe \xbm'
[X].$

$(\xbm CUM)$ for $ \xbm $ entails $(\xbm CUM)$ for $ \xbm':$

We have to show $ \xbm' [X] \xcc [Y] \xcc [X]$ $ \xch $ $ \xbm' [X]=
\xbm' [Y].$

Assume as above $ \xbm (X) \xcc Y \xcc X$ by $(\xbm S2)$ and sufficient
closure conditions
of the original domain, so $ \xbm (X)= \xbm (Y),$ and $ \xbm' [X]= \xbm'
[Y].$
\paragraph{
Additional properties \\[2mm]
}

\label{Section TR-Add}

The following additional condition is very important, it says that $x$ has
minimal distance to itself, and all $[x,x]$ are minimal. It corresponds to
a universal distance 0.

\bd

$\hspace{0.01em}$


\label{Definition Mu-Id}

$(\xbm Id)$ If there is for some $x$ $ \xBc x,x \xBe  \xbe A,$ then $ \xbm
(A)=\{ \xBc x,x \xBe  \xbe A\}$

\ed

This condition translates directly to a condition about $ \xbm':$

If $[x,x] \xbe [A],$ then $ \xbm' [A]=\{[x,x] \xbe [A]\}.$

The following conditions are taken from
 \cite{Sch04} (Section 3.2.3 there).
The first two impose that 1 copy for each point suffice, in
the finite and infinite case, the third imposes transitivity:

\bd

$\hspace{0.01em}$


\label{Definition Additional-Cond}

(1-fin) Let $X=A \xcv B_{1} \xcv B_{2}$ and $A \xcs \xbm (X)= \xCQ.$ Then
$A \xcc (A \xcv B_{1}- \xbm (A \xcv B_{1})) \xcv (A \xcv B_{2}- \xbm (A
\xcv B_{2})).$

(1-infin) Let $X=A \xcv \xcV \{B_{i}:i \xbe I\},$ and $A \xcs \xbm (X)=
\xCQ.$ Then
$A \xcc \xcV \{A \xcv B_{i}- \xbm (A \xcv B_{i})\}.$

(T) $ \xbm (A \xcv B) \xcc A,$ $ \xbm (B \xcv C) \xcc B$ $ \xch $ $ \xbm
(A \xcv C) \xcc A.$

\ed

It is also discussed there that in the 1-copy case, not all structures can
be made transitive, contrary to the situation when we allow arbitrary many
copies. The reader is referred there for further reference.
\paragraph{
The limit version \\[2mm]
}

\label{Section TR-Limit}

We know from
Section \ref{Section Pref-Lim} (page \pageref{Section Pref-Lim}),
Fact \ref{Fact D-8.2.1} (page \pageref{Fact D-8.2.1})  that for $ \xCf
transitive$ relations $ \xeb
$ on $ \xdy,$
the MISE system $ \xbL ([X])$ has the following properties, corresponding
to the
system of sets of minimal elements:

(1) If $A \xbe \xbL (Y),$ and $A \xcc X \xcc Y,$ then $A \xbe \xbL (X).$

(2) If $A \xbe \xbL (Y),$ and $A \xcc X \xcc Y,$ and $B \xbe \xbL (X),$
then $A \xcs B \xbe \xbL (Y).$

(2a) Let $A \xbe \xbL (Y),$ $A \xcc X \xcc Y.$ Then, if $B \xbe \xbL (Y),$
$A \xcs B \xbe \xbL (X).$ Conversely,
if $B \xbe \xbL (X),$ then $A \xcs B \xbe \xbL (Y).$

(3) If $A \xbe \xbL (Y),$ $B \xbe \xbL (X),$ then there is $Z \xcc A \xcv
B$ $Z \xbe \xbL (Y \xcv X).$

Using $(\xbm S2)$ and closure of the domain as above, we see that these
conditions
carry over to conditions about the MISE systems $ \xbL (X)$ of the
original domain.
\paragraph{
Higher order structures, and their limit version \\[2mm]
}

We know from previous work, see
Section \ref{Section High-Pref} (page \pageref{Section High-Pref}),
Proposition \ref{Proposition Eta-Rho-Repres} (page \pageref{Proposition
Eta-Rho-Repres})  and
Proposition \ref{Proposition Level-3-Repr} (page \pageref{Proposition
Level-3-Repr})
that any function $ \xbm' $ satisfying $(\xbm \xcc)$ can be represented
by a higher
preferential structure, and if $ \xbm' $ satisfies also $(\xbm CUM)$
(and $(\xbm \xcS)),$ it can
be represented by an essentially smooth higher preferential structure.

As this, too, was a very abstract result, independent of logic, it
also carries over to our situation.

We discussed in
Section \ref{Section Lim-High} (page \pageref{Section Lim-High})  that more
research has to go into
higher preferential structures before we can present
an intuitively valid limit version. So this stays an
open problem.
\subsubsection{
The original problem: products and projections
}

\label{Section TR-Proj}

We know from previous work, see, e.g.,
 \cite{GS08f}, Section 4.2.2.3, Example 4.2.4
(see also  \cite{Sch96-3},  \cite{Sch04},
 \cite{GS08a}),
that we need closure of the domain under finite unions for
representation of cumulative functions by smooth structures.
Obviously, the domain of products $X \xCK Y$ is usually not closed
under finite unions. It is therefore hopeless to obtain nice
representation results for many cases when we work directly
with products, we have to work with more general sets of
pairs, at least for testing the conditions. This was one of the reasons
why we
split the representation problem immediately into several layers.
Recall also that we do not see the resulting sets of minimal elements,
but only their projections - so we are quite
``blind''. We have shown in
 \cite{Sch04} that such kind of blindness can have very serious
consequences on representation results. See the results on
absence of definability preservation there, in particular
Section 5.2.3, Proposition 5.2.15 there, but also the
``Hamster Wheels'' described in
 \cite{Sch04}.

So let us see how we go from products to more arbitrary sets of pairs,
and treat our problem with projections. We adopt the same strategy as
we did in  \cite{Sch96-1}.
(Our solution is more general, as we also manage more general input,
to solve the smooth case.)
In this article, we ``hid'' the problem in a formula
$ \xbf (\xBc a,b \xBe,A \xCK B),$ we make it now explicit. Obviously, our
solution is not very nice, but we see no good alternative.
To put this into perspective, the reader may want to look back
to Section \ref{Section Mod-TR} (page \pageref{Section Mod-TR}).

A semantical revision function $ \xfA,$ which
assigns to every pair of (intuitively: model) sets
$ \xBc A,B \xBe $ a subset $A \xfA B \xcc B$ is representable by a
preferential, smooth, higher preferential, etc., structure iff the
following
condition $(\xbm \xcE)$ of existence of a more general function $f$
holds:

\bd

$\hspace{0.01em}$


\label{Definition Mu-Ex}

$(\xbm \xcE)$ There is a representing function $f$ for more arbitrary
sets of pairs of elements
$ \xBc a,b \xBe $ with the required properties, accepting all original
product sets $A \xCK B$ as arguments, and such that the
projection of the set of minimal elements on the second coordinate
is $A \xfA B,$ i.e., $ \xbp_{2}(f(A \xCK B))=A \xfA B.$

\ed

Again, we pay generality with lack of elegance.
\subsection{
Semantic representation for generalized update and counterfactuals
}

\label{Section Semantic-Up}
\subsubsection{
Introduction
}

We work here in the spirit of the Stalnaker/Lewis semantics for
counterfactual
conditionals,
see  \cite{Sta68},  \cite{Lew73}.
We treat update and counterfactual conditionals together, as it has
been tradition since the seminal paper by Katsuno and Mendelson,
 \cite{KM90}.

At first sight, everything seems simple, as we
``simply add up what we see from the points on the left'',
i.e., $X \xfB Y$ $=$ $ \xcV \{\{x\} \xfB Y:x \xbe X\}$ $=$ $ \xcV \{\{y
\xbe Y: \xCN \xcE y' \xbe Y(d(x,y') \xeb d(x,y)\}:x \xbe X\}.$
So it is tempting to do as if all distances from fixed $x$
can be chosen independently (up to symmetry) for all other distances.

But things are a bit more complicated, as we
``look into the pairs''
$ \xBc x,y \xBe,$ they all have to have $x$ on the left. So situations about
transitivity
like $[ab]<[bc]<[cd]<[ad],$ thus $[ab]<[ad],$ an example discussed in
 \cite{MS90} have to be considered. In this article, we constructed
a common real
valued distance
from a set of independent distances (one for each point, showing the
world as ``it sees it''), by multiplying copies. The same strategy
works again here.

In more detail:

 \xEI

 \xDH
We first go from distances to equivalence classes respecting symmetry
as we did in
Section \ref{Section Semantic-TR} (page \pageref{Section Semantic-TR}).

 \xDH
We then work with finite sequences of points like
$ \xBc a,b,x,y \xBe $ and let essentially
$d(\xBc a,b,x \xBe, \xBc a,b,x,y \xBe):=d(x,y),$
the latter the original distance, and
make all other distances infinite (except when both sequences are
identical, of course).
Thus, we compare in a non-trivial way only sequences $ \xbs, \xbt,$
where $ \xbt $ is
$ \xbs,$ with one element added at the end.

 \xDH
In addition, we make all distances between shorter sequences bigger
than all distances between longer sequences (unless they are already
infinite or 0). This gives us a layered structure, repeating
$ \xbo $ many times (in descending order) the distance structure between
elements.

 \xDH
Finally, we identify sequences with their end points, and have again that
every point sees the world according to its own distances, but we work
with
a common abstract distance.

 \xEJ
\subsubsection{
The general and smooth case
}

\label{Section Up-Gen}

We describe now formally the idea.

We consider sequences of points from $U$ without direct repetitions,
i.e., $ \xbs_{i} \xEd \xbs_{i+1}.$ We then define
\begin{flushleft}
\[ d(\xbs,\xbt):= \left\{ \begin{array}{lcll}
0 \xEH iff \xEH \xbs=\xbt \xEH \xEP
\xEH \xEH \xEH \xEP
d(\xbs_n,\xbs_{n+1}) \xEH iff \xEH \xbs=\xbs_0,..,\xbs_{n+1},
\xbt=\xbs_0,..,\xbs_{n} \xEH or \xEP
\xEH \xEH \xbt=\xbt_0,..,\xbt_{n+1},
\xbs=\xbt_0,..,\xbt_{n} \xEH \xEP
\xEH \xEH \xEH \xEP
\xca \xEH \xEH otherwise \xEH \xEP
\end{array}
\right.
\]
\end{flushleft}

This gives us distances between sequences of length $n+1$ and $n+2.$
The distance between sequences is symmetrical iff the base distance is
symmetrical.

Now, we have to arrange all these layers in a way to make distances
between
longer sequences smaller, resulting in the following picture,
see Diagram \ref{Diagram Update-S} (page \pageref{Diagram Update-S}). In this
diagram, the distance
situation
between points is shown on the right, between corresponding sequences on
the left. 0 is the smallest distance, $ \xca $ the biggest, and we have an
infinite descending chain of layers. The distance between sequences is
symmetrical iff the original distance was, it is transitive, iff the
original
distance was.
(Formally, $d(\xbs, \xbs +x)>d(\xbt, \xbt +x)$ if $ \xbs $ and $ \xbt
$ have the same end point, and $ \xbs $ is
longer than $ \xbt $ is. Here, $+$ is appending one more element.)

$ \xCO $

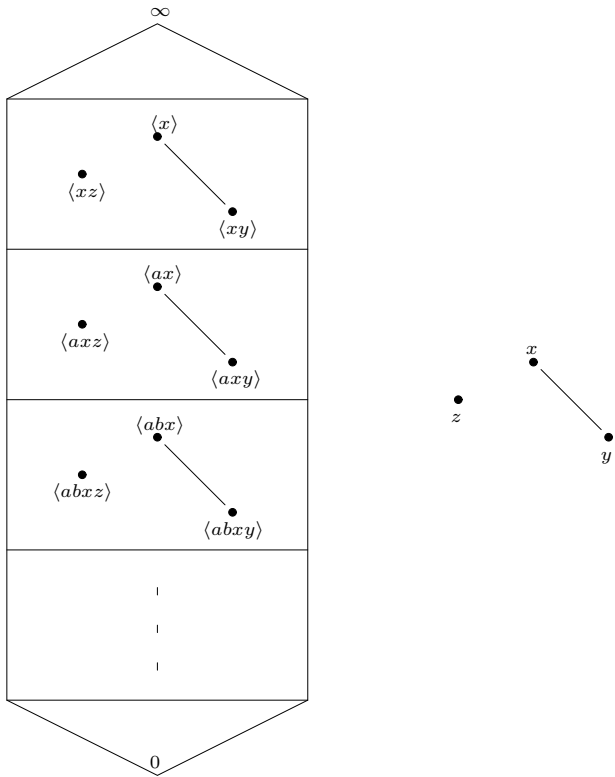
\vspace{10mm}

\begin{diagram}

\label{Diagram Update-S}
\index{Diagram Update-S}

\centering
\setlength{\unitlength}{1mm}
{\renewcommand{\dashlinestretch}{30}
\begin{picture}(150,150)(0,0)

\put(90,60){\circle*{1}}
\put(89,61){{\xssc $x$}}
\put(80,55){\circle*{1}}
\put(79,52){{\xssc $z$}}
\put(100,50){\circle*{1}}
\put(99,47){{\xssc $y$}}
\path(91,59)(99,51)

\put(40,90){\circle*{1}}
\put(39,91){{\xssc $\xBc x \xBe$}}
\put(30,85){\circle*{1}}
\put(28,82){{\xssc $\xBc xz \xBe$}}
\put(50,80){\circle*{1}}
\put(48,77){{\xssc $\xBc xy \xBe$}}
\path(41,89)(49,81)

\put(40,70){\circle*{1}}
\put(38,71){{\xssc $\xBc ax \xBe$}}
\put(30,65){\circle*{1}}
\put(27,62){{\xssc $\xBc axz \xBe$}}
\put(50,60){\circle*{1}}
\put(47,57){{\xssc $\xBc axy \xBe$}}
\path(41,69)(49,61)

\put(40,50){\circle*{1}}
\put(37,51){{\xssc $\xBc abx \xBe$}}
\put(30,45){\circle*{1}}
\put(26,42){{\xssc $\xBc abxz \xBe$}}
\put(50,40){\circle*{1}}
\put(46,37){{\xssc $\xBc abxy \xBe$}}
\path(41,49)(49,41)

\path(40,5)(20,15)(20,95)(40,105)(60,95)(60,15)(40,5)
\path(20,95)(60,95)
\path(20,75)(60,75)
\path(20,55)(60,55)
\path(20,35)(60,35)
\path(20,15)(60,15)
\put(39,106){{\xssc $\xca$}}
\put(39,6){{\xssc $0$}}

\path(40,30)(40,29)
\path(40,25)(40,24)
\path(40,20)(40,19)

\end{picture}
}

\end{diagram}

\vspace{4mm}

$ \xCO $

This is an $ \xda -$ranked structure, as discussed in
 \cite{GS08d}, see also  \cite{GS08f}.

When we evaluate the structure, we are interested in the set of (original)
points $Y$ individually closest to some set of original points $X.$ We
look now at
all sequences whose end points are in $X,$ and from those to all sequences
whose
end points are in $Y.$ The closest ones have the same end points as the
(indiviually) closest to $X$ elements of $Y.$ So we are done, and can work
with
the new structure, representing it by some relation $ \xeb.$

Not yet, quite: We have destroyed cumulativity, if this was a property.
But only
superficially, as we argue now. $ \xCf Inside$ each layer, smoothness was
preserved, but each layer gives the same answer, so, in the end,
cumulativity
$ \xCf is$ preserved.

More precisely,
the property of cumulativity in our context is the following:

\bd

$\hspace{0.01em}$


\label{Definition Up-Cum}

$ \xCf (Up- \xCf Cum)$
Fix $X,$ consider $Y$ and $Y'.$ If $Y' \xcc Y,$ and $X \xfB Y \xcc Y',$
then $X \xfB Y=X \xfB Y'.$

\ed

As the new structure was equivalent to the old structure in each layer,
in each layer cumulativity is preserved, and the structure can thus be
represented by a smooth relation in each layer, and we just put
the layers together in the final construction, we obtain local smoothness.
Details about cumulativity are left to the reader.

We now give a formal representation result for the basic case.

Consider the following property:

\bd

$\hspace{0.01em}$


\label{Definition Mu-Cover}

$(\xbm \xcV)$ If $ \xdx $ is a cover of $X,$ then $ \xbm (X)= \xcV \{
\xbm (X'):X' \xbe \xdx \}.$

\ed

We postulate: $ \xfB $ satisfies $(\xbm \xcV)$ in the following sense on
the left:

All $ \xbm_{Y}(X):=X \xfB Y$ satisfy $(\xbm \xcV)$ for any $Y.$

We can now define:

\bd

$\hspace{0.01em}$


\label{Definition Indiv-Result}

(1) $x \xfN Y:= \xcS \{X \xfB Y:x \xbe X\},$

(2) $A \xfN B:= \xcV \{a \xfN B:a \xbe A\}.$

\ed

Note that $\{x\} \xfB Y$ need not be defined, this way we avoid
postulating that
singletons are definable. We have:

\bfa

$\hspace{0.01em}$


\label{Fact Indiv-Result}

Let the domain be closed under finite intersections, let $(\xbm \xcV)$
hold for
$ \xfB $ on the left, then $A \xfN B=A \xfB B.$

\efa

\subparagraph{
Proof
}

$\hspace{0.01em}$


``$ \xcc $'': $b \xbe A \xfN B$ $ \xch $ $ \xcE a \xbe A.b \xbe a \xfN B$
$ \xch $ $b \xbe A \xfB B,$ both implications by
Definition \ref{Definition Indiv-Result} (page \pageref{Definition
Indiv-Result}).

``$ \xcd $'': $b \xce A \xfN B$ $ \xch $ (by definition) $ \xcA a \xbe
A.b \xce a \xfN B$ $ \xch $ (by definition)
$ \xcA a \xbe A \xcE X(a \xbe X \xcu b \xce X \xfB B)$ $ \xch $ (by $(
\xbm \xcV))$ $ \xcA a \xbe A \xcE X(a \xbe X \xcs A \xcu b \xce (X \xcs
A) \xfB B).$
As there is such $X \xcs A$ for all $a \xbe A,$ such $X \xcs A$ form a
cover of A, and by $(\xbm \xcV)$
again, $a \xce A \xfB B.$

$ \xcz $
\\[3ex]

Note that closure under finite intersections is here a very weak
prerequisite,
usually satisfied, and, in addition, it is $ \xCf not$ required on
products or so,
only on one coordinate.

We try to show now that relevant properties carry over from $ \xfB $ to $
\xfN.$
It seems, however, that we need an additional property:

\bd

$\hspace{0.01em}$


\label{Definition Supp-Cond}

(Approx) If $x \xfN Y \xcc A,$ then there is $X$ such that $x \xbe X,$ $X
\xfB Y \xcc A.$

\ed

This condition is not as strong as requiring that singletons are in the
domain,
it says that we can approximate the results from singletons sufficiently
well.

We have now:

\bfa

$\hspace{0.01em}$


\label{Fact Cond-Transfer}

$(\xbm \xcc)$ and $(\xbm PR)$ carry over from $ \xfB $ to $ \xfN,$ if
(Approx) holds, the same
is true for $(\xbm CUM).$

More precisely:

(1) If $X \xfB Y \xcc Y,$ then also $x \xfN Y \xcc Y$ (for $x \xbe X).$

(2) If $Y \xcc Y' \xch (X \xfB Y') \xcs Y \xcc X \xfB Y,$ then also $Y
\xcc Y' \xch (x \xfN Y') \xcs Y \xcc x \xfN Y.$

(3) Let (Approx) and $(\xcs)$ hold. Then:
If $(X \xfB Y') \xcc Y \xcc Y' \xch (X \xfB Y)=(X \xfB Y'),$ then also
$(x \xfN Y') \xcc Y \xcc Y' \xch (x \xfN Y)=(x \xfN Y').$

(For all $X,Y,Y',x,$ etc.)

\efa

\subparagraph{
Proof
}

$\hspace{0.01em}$


(1) Trivial.

(2) Let $Y \xcc Y'.$ $y \xbe (x \xfN Y') \xcs Y$ $ \xcj $ $y \xbe Y$ $
\xcu $ $ \xcA X(x \xbe X \xch y \xbe X \xfB Y')$ $ \xch $
(by $(\xbm PR)$ for $ \xfB)$ $ \xcA X(x \xbe X \xch y \xbe X \xfB Y)$ $
\xch $ $y \xbe x \xfN Y.$

(3) Let $(x \xfN Y') \xcc Y \xcc Y'.$
By (2) $(x \xfN Y') \xcs Y \xcc x \xfN Y,$ by prerequisite $x \xfN Y'
\xcc Y,$ so
$x \xfN Y' \xcc x \xfN Y.$ It remains to show $x \xfN Y \xcc x \xfN Y'.$
Let $y \xbe x \xfN Y,$ so for all $X$ such that
$x \xbe X,$ $y \xbe X \xfB Y.$ Assume $y \xce x \xfN Y',$ so there is $X'
,$ $x \xbe X',$ $y \xce X' \xfB Y'.$
We use now (Approx) to chose suitable $X'' $ with $x \xbe X'',$ such that
$X'' \xfB Y' \xcc Y.$ Now $x \xbe X' \xcs X'',$ so $y \xbe (X' \xcs X'')
\xfB Y.$ By $(\xbm \xcV)$
$(X' \xcs X'') \xfB Y' \xcc X' \xfB Y',$ so $y \xce (X' \xcs X'') \xfB
Y'.$ By $X'' \xfB Y' \xcc Y$ and $(\xbm \xcV)$
$(X' \xcs X'') \xfB Y' \xcc Y,$ so by $(\xbm CUM)$ for $ \xfB,$ $(X'
\xcs X'') \xfB Y=(X' \xcs X'') \xfB Y',$ contradiction.

$ \xcz $
\\[3ex]

We put our ideas together:

Let $ \xfB $ satisfy:

(1) $(\xbm \xcV)$ on the left,

(2) $(\xbm \xcc)$ and $(\xbm PR)$ on the right,

(3) let the domain be closed under $(\xcs).$

By Fact \ref{Fact Cond-Transfer} (page \pageref{Fact Cond-Transfer}), $x \xfN
Y$ has suitable conditions,
so we can represent
$x \xfN Y$ individually for each $x$ by a suitable relation. The
construction of
Section \ref{Section Up-Gen} (page \pageref{Section Up-Gen})  shows how to
combine the individual
relations.
Fact \ref{Fact Indiv-Result} (page \pageref{Fact Indiv-Result})  shows that $X
\xfN Y=X \xfB Y.$
So we can define $X \xfN Y$ from the individual $x \xfN Y,$ having $X \xfB
Y=X \xfN Y.$
This results in a suitable choice function for each element $x,$ which
can,
individually, be represented by a relation $ \xeb_{x}$ (for each, fixed,
$x).$
The construction with $ \xbo $ many layers from the beginning of
Section \ref{Section Up-Gen} (page \pageref{Section Up-Gen})  shows how to
construct a global relation $
\xeb $ on a suitable
space of sequences, whose evaluation gives the same results as the
individual
relations $ \xeb_{x}.$ Conversely, any such order gives a choice function
satisfying $(\xbm \xcc)$ and $(\xbm PR)$ on the right, and, in suitable
interpretation,
$(\xbm \xcV)$ on the left.

If $ \xfB $ satisfies $(\xbm CUM)$ on the right, and if we can
approximate the results
of singletons sufficiently well, i.e., if (Approx) holds,
the choice functions $x \xfN (.)$ will also
have the property $(\xbm CUM),$ they can represented by a smooth
structure for
each $x,$ and the individual structures can be put together as above. The
resulting global structure is not globally smooth, but locally smooth,
which
is sufficient. Conversely, if the structure is locally smooth, $ \xfB $
and $ \xfN $
will satisfy $(\xbm CUM).$
\subsubsection{
Further conditions and representation questions
}

 \xEI

 \xDH
Note that we built the 0-property already into the basic construction, so
there is nothing to do here.

 \xDH
The 1-copy case was already discussed in
Section \ref{Section TR-Add} (page \pageref{Section TR-Add}), we refer the
reader there.

 \xDH
For the limit version, all the important algebraic material stays valid,
so we can use it to show that the same (unitary versions of) laws hold
there as in the minimal variant, analogous to the TR case,
see Section \ref{Section TR-Limit} (page \pageref{Section TR-Limit}).

 \xDH
Higher preferential structures need not satisfy $(\xbm PR),$ so even
distances not satisfying this condition can be represented by a
suitable combination of techniques we are now familiar with.

 \xEJ
\subsection{
Syntactic representation for generalized revision, update, counterfactuals
}

\label{Section TR-Up-Synt}

Our aim in this
Section \ref{Section Repr-TR-Up} (page \pageref{Section Repr-TR-Up})
is not only to generalize the notions of distance, but also to generalize
to
what revision and update can be applied, and in which way (minimal or
limit
version). Consequently, the syntactic side is a less important part than
in classical revision and update (and counterfactuals). So we will only
indicate how to proceed, and will leave the rest to the reader - or future
research.

First, some introductory remarks.

Translation from semantics to logics and back is discussed extensively
in  \cite{Sch04}, see, e.g., Section 3.4, 4.2.3 there. Particular
attention
is given to definability preservation and the problems arising from the
lack
of it, see there Chapter 5. Further discussion can be found in
 \cite{GS08f}, see, e.g., Section 5.4 there.
\paragraph{
Discussion of the Theory Revision situation \\[2mm]
}

We have to look back at
Definition \ref{Definition Mu-Ex} (page \pageref{Definition Mu-Ex}).

The input, $A \xCK B,$ is simple, and will be
$M(\xbf) \xCK M(\xbq)$ for formulas, $M(S) \xCK M(T)$ for full
theories.

The output will be (exactly - to obtain definability preservation)
some $M(\xbq')$ or $M(T').$

Problems are in-between. We have to find a language (which will not just
be
classical propositional language) rich enough to

 \xEI

 \xDH
Express finite unions of products of the type $M(\xbf) \xCK M(\xbq).$

 \xDH
Rich enough to describe the result of the choice function in
the desired situation. Usually, this will not just be a product,
or even a finite union of products.

 \xDH
Going to the projection will probably not be very difficult.

 \xEJ

It seems that one needs here in many cases specially tailored
languages.
\paragraph{
Discussion of the update/counterfactual situation \\[2mm]
}

The case for update and counterfactuals is easier, essentially because
we do not need any cross-comparison between
$ \xBc a,b \xBe $ and $ \xBc a',b'  \xBe.$

The conditions $(\xbm \xcV),$ $(\xbm \xcc),$ and $(\xbm PR)$
translate directly
to the syntactical side - as usual. Care has to be taken to make the
results
of minimization definability preserving.
\paragraph{
Example of syntactic conditions \\[2mm]
}

To give a flavour of a full set of conditions, we quote from
 \cite{Sch96-1} the logical counterpart of the semantical
representation
result, Theorem 3.1 there.

We first introduce a definition (Definition 3.3 there).

We use the abbreviation
``cct'' for ``consistent complete theory'',
i.e., corresponding to a single model.

Fix a propositional language $ \xdl.$

\bd

$\hspace{0.01em}$


\label{Definition PDS-3}

We consider two logics, $.^{i}: \xdp (\xdl) \xCK \xdp (\xdl) \xcp \xdp
(\xdl),$ $ \xBc S,T \xBe  \xcZ  \xBc S,T \xBe ^{i} \xcc \xdl.$

We say that both logics $ \xBc S,T \xBe ^{i}$ are given by a function $f: \xdD
\xCK
\xdD \xcp \xdp (M_{ \xdl } \xCK M_{ \xdl })$
iff for all theories $S,T$

$ \xBc S,T \xBe ^{i}=\{ \xbf: \xcA m \xbe \xbp_{i}(f(M_{S} \xCK M_{T})).m \xcm
\xbf
\},$ where $ \xbp_{i}$ is the projection on the
$i-$th coordinate. Analogously, we say that they are given by a
preferential
structure on $M_{ \xdl } \xCK M_{ \xdl }$ iff $f$ is given by such a
structure.

We call a function $f$ on pairs of models definability preserving (dp)
iff for all theories $S,$ $T$ $ \xbp_{1}(f(M_{S} \xCK M_{T}))=M_{U}$ and $
\xbp_{2}(f(M_{S} \xCK M_{T}))=M_{V}$
for some theories $U,$ $V,$ where $M_{S}$ is the set of $S-$models etc.

Note that then $ \xbp_{i}(f(M_{S} \xCK M_{T}))=M_{ \xBc S,T \xBe ^{i}}$.

For $U,V$ complete, consistent theories (cct), $S,T$ any theories, we
abbreviate

$ \xbq (U,V,S,T):= \xcA S',T' (U \xcl S' \xcl S$ $ \xcu $ $V \xcl T' \xcl
T$ $ \xcp $ $U \xcl  \xBc S',T'  \xBe ^{1}$ $ \xcu $ $V \xcl  \xBc S',T'  \xBe
^{2}).$

\ed

We quote now the syntactic result, Theorem 3.4 of
 \cite{Sch96-1}:

\bp

$\hspace{0.01em}$


\label{Proposition PDS-4}

Let $ \xBc S,T \xBe ^{1}$, $ \xBc S,T \xBe ^{2}$ be two logics on pairs of
theories. Then
$ \xBc S,T \xBe ^{i}$ are given by
a dp preferential structure iff

(1) $ \ol{S}= \ol{S' }$ $ \xcu $ $ \ol{T}= \ol{T' }$ $ \xcp $ $ \xBc S,T \xBe
^{1}$
$=$ $ \xBc S',T'  \xBe ^{1}$ $ \xcu $ $ \xBc S,T \xBe ^{2}$ $=$ $ \xBc S',T' 
\xBe ^{2}$

(2) $ \xBc S,T \xBe ^{i}$ is classically closed

(3) $ \xBc S,T \xBe ^{1} \xcl S$, $ \xBc S,T \xBe ^{2} \xcl T$

(4S) If $U$ is a cct with $U \xcl S,$ then
$U \xcl  \xBc S,T \xBe ^{1}$ iff there is a cct $V$ such that $V \xcl T$ and $
\xbq
(U,V,S,T)$

(4T) If $V$ is a cct with $V \xcl T,$ then
$V \xcl  \xBc S,T \xBe ^{2}$ iff there is a cct $U$ such that $U \xcl S$ and $
\xbq
(U,V,S,T).$

\ep

The proof is straightforward, and the reader is referred to
 \cite{Sch96-1}.

$ \xCO $

$ \xCO $
\chapter{
Monotone and antitone semantic and syntactic interpolation
}

\label{Chapter Mod-Mon-Interpol}
\section{
Introduction
}

\label{Section Inter-Intro}
\subsection{
Overview
}

The two chapters
Chapter \ref{Chapter Mod-Mon-Interpol} (page \pageref{Chapter Mod-Mon-Interpol})
 and
Chapter \ref{Chapter Size-Laws} (page \pageref{Chapter Size-Laws})
are probably the core of the present book.

We treat very general interpolation problems for monotone and antitone,
2-valued and many-valued logics in the present chapter,
splitting the question in two parts,
``semantic interpolation'' and
``syntactic interpolation'', show that the first problem, existence
of semantic interpolation, has a simple and general answer, and reduce
the second question, existence of syntactic interpolation to a
definability problem. For the latter, we examine the concrete example
of finite Goedel logics. We can also show that the semantic problem
has two ``universal'' solutions, which depend only on one formula and
the shared variables.

In Chapter \ref{Chapter Size-Laws} (page \pageref{Chapter Size-Laws}), we
investigate three variants of
semantic
interpolation for non-monotonic logics, in syntactic shorthand of the
types
$ \xbf \xcn \xba \xcl \xbq,$ $ \xbf \xcl \xba \xcn \xbq,$ and $ \xbf
\xcn \xba \xcn \xbq,$ where $ \xba $ is the interpolant, and
see that two variants are closely related to multiplication laws about
abstract
size, defining (or originating from) the non-monotonic logics.
The syntactic problem is analogous to that of the monotonic case.
\subsubsection{
Background
}

Interpolation for classical logic is well-known,
see  \cite{Cra57}, and there are also
non-classical logics for which interpolation has been shown, e.g., for
Circumscription, see
 \cite{Ami02}. An extensive overview of interpolation is found
in  \cite{GM05}. Chapter 1 of this book gives a survey and a
discussion and
the chapter puts forward that interpolation can be viewed in many
different ways
and indeed
11 points of view of interpolation are discussed. The present text
presents the
semantic interpolation, this is a new 12th point of view.
\subsection{
Problem and Method
}

In classical logic, the problem of interpolation is to find for two
formulas $ \xbf $ and $ \xbq $ such that $ \xbf \xcl \xbq $ a
``simple'' formula $ \xba $ such that $ \xbf \xcl \xba \xcl \xbq.$
``Simple'' is defined as:
``expressed in the common language of $ \xbf $ and $ \xbq $''.

Working on the semantic level has often advantages:

 \xEI

 \xDH results are robust under logically equivalent reformulations

 \xDH in many cases, the semantic level allows an easy reformulation as an
algebraic problem, whose results can be generalized to other situations

 \xDH we can split a problem in two parts: a semantical problem, and the
problem
to find a syntactic counterpart (a definability problem)

 \xDH the semantics of many non-classical logics are built on relatively
few
basic notions like size, distance, etc., and we can thus make connections
to other problems and logics

 \xDH in the case of preferential and similar logics, the very definition
of
the logic is already semantical (minimal models), so it is very natural
to proceed on this level

 \xEJ

This strategy - translate to the semantic level, do the main work there,
and
then translate back - has proved fruitful also in the present case.

Looking back at the classical interpolation problem, and translating it
to the semantic level, it becomes: Given $M(\xbf) \xcc M(\xbq)$ (the
models sets of $ \xbf $ and
$ \xbq),$ is there a
``simple'' model set $ \xCf A$ such that $M(\xbf) \xcc A \xcc M(\xbq)?$
Or, more generally, given
model sets $X \xcc Y,$ is there ``simple'' $ \xCf A$ such that $X \xcc A
\xcc Y?$

Of course, we have to define in a natural way,
what ``simple'' means in our context.
This is discussed below in
Section \ref{Section Intro-Mon-Sema-b} (page \pageref{Section Intro-Mon-Sema-b})
.

Our separation of the semantic from the syntactic question
pays immediately:

 \xEh

 \xDH
We see that monotonic (and antitonic) logics $ \xCf always$ have a
semantical
interpolant. But this interpolant may not be definable syntactically.

 \xDH
More precisely, we see that there is a whole interval of interpolants in
above situation.

 \xDH
We also see that our analysis generalizes immediately to many valued
logics, with the same result (existence of an interval of interpolants).

 \xDH
Thus, the question remains: under what conditions does a syntactic
interpolant exist?

 \xDH
In non-monotonic logics, our analysis reveals a deep connection between
semantic interpolation and questions
about (abstract) multiplication of (abstract) size.

 \xEj
\subsection{
Monotone and antitone semantic and syntactic interpolation
}

We consider here the semantic property of monotony or antitony, in the
following sense (in the two-valued case, the generalization to the
many-valued case is straightforward):

Let $ \xcl $ be some logic such that $ \xbf \xcl \xbq $ implies $M(\xbf)
\xcc M(\xbq)$
(the monotone case) or $M(\xbq) \xcc M(\xbf)$ (the antitone case).

In the many-valued case, the corresponding property is that $ \xcp $ (or $
\xcl)$
respects $ \xck,$ the order on the truth values.
\subsubsection{
Semantic interpolation
}

\label{Section Intro-Mon-Sema-b}

The problem (for simplicity, for the 2-valued case) reads now:

If $M(\xbf) \xcc M(\xbq)$ (or, symmetrically $M(\xbq) \xcc M(\xbf
)),$ is there a
``simple'' model set $ \xCf A,$ such that $M(\xbf) \xcc A \xcc M(\xbq),$
or $M(\xbq) \xcc A \xcc M(\xbf).$
Obviously, the problem is the same in both cases.
We will see that such $ \xCf A$ will always exist, so all such logics have
semantic interpolation (but not necessarily also syntactic interpolation).

The main conceptual problem is to define
``simple model set''. We have to look at the
syntactic problem for guidance. Suppose $ \xbf $ is defined using
propositional
variables $p$ and $q,$
$ \xbq $ using $q$ and $r.$ $ \xba $ has to be defined using only $q.$
What are the models of
$ \xba?$ By the very definition of validity in classical logic, neither
$p$ nor $r$
have any influence on whether $m$ is a model of $ \xba $ or not. Thus, if
$m$ is a model
of $ \xba,$ we can modify $m$ on $p$ and $r,$ and it will still be a
model. Classical
models are best seen as functions from the set of propositional variables
to
$\{true,false\},$ $\{t,f\},$ or so. In this terminology, any $m$ with $m
\xcm \xba $ is
``free'' to choose the value for $p$ and $r,$ and we can write the model set
A
of $ \xba $ as $\{t,f\} \xCK M_{q} \xCK \{t,f\},$ where $M_{q}$ is the set
of values for $q$ $ \xba -$models
may have $(\xCQ,$ $\{t\},$ $\{f\},$ $\{t,f\}).$

So, the semantic interpolation problem is to find sets which may be
restricted
on the common variables,
but are simply the Cartesian product of the possible values for the other
variables. To summarize: Let two model sets $X$ and $Y$ be given, where
$X$ itself
is restricted on variables $\{p_{1}, \Xl,p_{m}\}$ (i.e. the Cartesian
product for the rest),
$Y$ is restricted on $\{r_{1}, \Xl,r_{n}\},$ then we have to find a model
set $ \xCf A$ which
is restricted only on $\{p_{1}, \Xl,p_{m}\} \xcs \{r_{1}, \Xl,r_{n}\},$
and such that
$X \xcc A \xcc Y,$ of course.

Formulated this way, our approach, the problem and its solution, has two
trivial generalizations:

 \xEI

 \xDH for multi-valued logics
we take the Cartesian product of more than just $\{t,f\}.$

 \xDH $ \xbf $ may be the hypothesis, and $ \xbq $ the consequence, but
also vice versa,
there is no direction in the problem. Thus, any result for classical
logic carries over to the
core part of, e.g., preferential logics.

 \xEJ

The main result for the situation with $X \xcc Y$ is that there is always
such a
semantic interpolant $ \xCf A$ as described above
(see Proposition \ref{Proposition Sin-Interpolation} (page \pageref{Proposition
Sin-Interpolation})  for a simple case,
and Proposition \ref{Proposition Gin-Pr-Int} (page \pageref{Proposition
Gin-Pr-Int})  for the full picture).
Our proof works also for
``parallel interpolation'', a concept introduced by
Makinson et al.,  \cite{KM07}.

We explain and quote the latter result.

Suppose we have $f,g:M \xcp V,$ where, intuitively, $M$ is the set of all
models, and
$V$ the set of all truth values. Thus, $f$ and $g$ give to each model a
truth value,
and, intuitively, $f$ and $g$ each code a model set, assigning to $m$ TRUE
iff $m$ is
in the model set, and FALSE iff not. We further assume that there is an
order
on the truth value set $V.$ $ \xcA m \xbe M(f(m) \xck g(m))$ corresponds
now to $M(\xbf) \xcc M(\xbq),$
or $ \xbf \xcl \xbq $ in classical logic. Each model $m$ is itself a
function from $L,$ the
set of propositional variables to $V.$ Let now $J \xcc L.$ We say that $f$
is
insensitive to $J$ iff the values of $m$ on $J$ are irrelevant: If $m \xex
(L-J)=m' \xex (L-$J),
i.e., $m$ and $m' $ agree at least on all $p \xbe L-$J, then $f(m)=f(m'
).$ This
corresponds to the situation where the variable $p$ does not occur in the
formula $ \xbf,$ then $M(\xbf)$ is insensitive to $p,$ as the value of
any $m$ on $p$
does not matter for $m$ being a model of $ \xbf,$ or not.

We need two more definitions:

Let $J' \xcc L,$ then
$f^{+}(m_{J' }):=max\{f(m'):m' \xex J' =m \xex J' \}$ and $f^{-}(m_{J'
}):=min\{f(m'):m' \xex J' =m \xex J' \}.$

We quote now
Proposition \ref{Proposition Gin-Pr-Int} (page \pageref{Proposition Gin-Pr-Int})
, slightly simplified:

\bp

$\hspace{0.01em}$


\label{Proposition Gin-Pr-Int-B}

Let $f,g:M \xcp V,$ $f(m) \xck g(m)$ for all $m \xbe M.$
Let $L=J \xcv J' \xcv J'',$ let $f$ be insensitive to $J,$ $g$ be
insensitive to $J''.$

Then $f^{+}(m_{J' }) \xck g^{-}(m_{J' })$ for all $m_{J' } \xbe M \xex J'
,$ and any $h:M \xex J' \xcp V$ which is
insensitive to $J \xcv J'' $ is an interpolant iff

$f^{+}(m_{J' }) \xck h(m_{J}m_{J' }m_{J'' })=h(m_{J' }) \xck g^{-}(m_{J'
})$ for all $m_{J' } \xbe M \xex J'.$

(h can be extended to the full $M$ in a unique way, as it is
insensitive to $J \xcv J''.)$

\ep

See Diagram \ref{Diagram Mon-Int} (page \pageref{Diagram Mon-Int}).
\subsubsection{
The interval of interpolants
}

Our result has an additional reading: it defines an interval of
interpolants, with lower bound $f^{+}(m_{J' })$ and upper bound
$g^{-}(m_{J' }).$
But these interpolants have a particular form. If they exist,
i.e. iff $f \xck g,$ then $f^{+}(m_{J' })$ depends only on $f$ and $J' $
(and $m),$ but $ \xCf not$ on $g,$
$g^{-}(m_{J' })$ only on $g$ and $J',$ $ \xCf not$ on $f.$ Thus, they are
universal, as we have
to look only at one function and the set of common variables.

Moreover, we will see in
Section \ref{Section Analoga} (page \pageref{Section Analoga})  that they
correspond to simple operations
on the normal
forms in classical logic. This is not surprising, as we
``simplify'' potentially complicated model sets by replacing some
coordinates with simple products. The question is, whether our logic
allows to express this simplification, classical logic does.
\subsubsection{
Syntactic interpolation
}

Recall the problem described at the beginning of
Section \ref{Section Intro-Mon-Sema-b} (page \pageref{Section Intro-Mon-Sema-b})
. We were given $M(\xbf) \xcc M(
\xbq),$ and were
looking for a ``simple'' model set $ \xCf A$ such that $M(\xbf) \xcc A
\xcc M(\xbq).$ We just saw that
such $ \xCf A$'s exists, and were able to describe an interval of such $
\xCf A$'s.
But we have no guarantee that any such $ \xCf A$ is definable, i.e., that
there
is some $ \xba $ with $A=M(\xba).$

In classical logic, such $ \xba $ exists, see, e.g.,
Proposition \ref{Proposition Sin-Simplification-Definable} (page
\pageref{Proposition Sin-Simplification-Definable})),
but also Section \ref{Section Analoga} (page \pageref{Section Analoga}).
Basically, in classical logic, $f^{+}(m_{J' })$ and $g^{-}(m_{J' })$
correspond to
simplifications of the formulas expressed in normal form,
see Fact \ref{Fact Class-Up-Down} (page \pageref{Fact Class-Up-Down})  (in a
different notation, which we
will
explain in a moment).
This is not necessarily true in other logics, see
Example \ref{Example Sin-4-Value} (page \pageref{Example Sin-4-Value}).
We find here again the importance of definability preservation,
a concept introduced by one of us in  \cite{Sch92}.

If we have projections (simplifications), see
Section \ref{Section Int-Int-Mon} (page \pageref{Section Int-Int-Mon}),
we also have syntactic interpolation. At present,
we do not know whether this is a necessary condition for all natural
operators.

We can also turn the problem around, and just define suitable
operators. This is done in
Section \ref{Section Analoga} (page \pageref{Section Analoga}),
Definition \ref{Definition 2-Many-Up} (page \pageref{Definition 2-Many-Up})  and
Definition \ref{Definition 2-Many-Down} (page \pageref{Definition 2-Many-Down})
.
There is a slight problem, as one of the operands is a $ \xCf set$ of
propositional
variables, and not a formula, as usual. One, but certainly not the only
one,
possibility is to take a formula (or the corresponding model set) and
``extract'' the ``relevant'' variables from it, i.e., those, which cannot
be replaced by a product. Assume now that $f$ is one of the
generalized model ``sets'', then:

Given $f,$ define

 \xEh
 \xDH
$(f \xfB J)(m)$ $:=$ $sup\{f(m'):$ $m' \xbe M,$ $m \xex J=m' \xex J\}$
 \xDH
$(f \xfb J)(m)$ $:=$ $inf\{f(m'):$ $m' \xbe M,$ $m \xex J=m' \xex J\}$
 \xDH
$ \xbf! \xbq $ by:

$f_{ \xbf! \xbq }$ $:=$ $f_{ \xbf } \xfB (L-R(\xbq))$
 \xDH
$ \xbf? \xbq $ by:

$f_{ \xbf? \xbq }$ $:=$ $f_{ \xbf } \xfb (L-R(\xbq))$
 \xEj

We then obtain for classical logic
(see Fact \ref{Fact Class-Up-Down} (page \pageref{Fact Class-Up-Down})):

\bfa

$\hspace{0.01em}$


\label{Fact Class-Up-Down-b}

Let $J:=\{p_{1,1}, \Xl,p_{1,m_{1}}, \Xl,p_{n,1}, \Xl,p_{n,m_{n}}\}$

(1) Let $ \xbf_{i}:= \xCL p_{i,1} \xcu  \Xl  \xcu \xCL p_{i,m_{i}}$ and $
\xbq_{i}:= \xCL q_{i,1} \xcu  \Xl  \xcu \xCL q_{i,k_{i}},$
let $ \xbf:=(\xbf_{1} \xcu \xbq_{1}) \xco  \Xl  \xco (\xbf_{n} \xcu
\xbq_{n}).$
Then $ \xbf \xfB J= \xbf_{1} \xco  \Xl  \xco \xbf_{n}.$

(2) Let $ \xbf_{i}:= \xCL p_{i,1} \xco  \Xl  \xco \xCL p_{i,m_{i}}$ and $
\xbq_{i}:= \xCL q_{i,1} \xco  \Xl  \xco \xCL q_{i,k_{i}},$
let $ \xbf:=(\xbf_{1} \xco \xbq_{1}) \xcu  \Xl  \xcu (\xbf_{n} \xco
\xbq_{n}).$
Then $ \xbf \xfb J= \xbf_{1} \xcu  \Xl  \xcu \xbf_{n}.$

In a way, these operators are natural, as they simplify definable
model sets, so they can be used as a criterion of the expressive
strength of a language and logic: If $X$ is definable, and $Y$ is in some
reasonable sense simpler than $X,$ then $Y$ should also be definable.
If the language is not sufficiently strong, then we can introduce these
operators, and have also syntactic interpolation.
\subsubsection{
Finite Goedel logics
}

\efa

The semantics of finite (intuitionistic) Goedel logics is a finite chain
of
worlds, which can also be expressed by a totally ordered set of truth
values
0 \Xl n (see Section \ref{Section Finite-Goedel} (page \pageref{Section
Finite-Goedel})).
Let FALSE and TRUE be the minimal and maximal truth values.
$ \xbf $ has value false, iff it holds nowhere, and TRUE, iff it holds
everywhere,
it has value 1 iff it holds from world 2 onward, etc.
The operators are classical $ \xcu $ and $ \xco,$ negation $ \xCN $ is
defined by
$ \xCN (FALSE)=TRUE$ and $ \xCN (x)=FALSE$ otherwise.
Implication $ \xcp $ is defined by $ \xbf \xcp \xbq $ is TRUE iff $ \xbf
\xck \xbq $ (as truth values),
and the value of $ \xbq $ otherwise.

More precisely, where $f_{ \xbf }$ is the model value function of the
formula $ \xbf:$

negation $ \xCN $ is defined by:
\begin{flushleft}
\[ f_{\xCN \xbf}(m):= \left\{ \begin{array}{lcl}
TRUE \xEH iff \xEH
f_{\xbf}(m)=FALSE \xEP
\xEH \xEH \xEP
FALSE \xEH \xEH otherwise \xEP
\end{array}
\right.
\]
\end{flushleft}

implication $ \xcp $ is defined by:
\begin{flushleft}
\[ f_{\xbf \xcp \xbq}(m):= \left\{ \begin{array}{lcl}
TRUE \xEH iff \xEH
f_{\xbf}(m) \xck f_{\xbq}(m) \xEP
\xEH \xEH \xEP
f_{\xbq}(m) \xEH \xEH otherwise \xEP
\end{array}
\right.
\]
\end{flushleft}

see Definition \ref{Definition Mod-Fin-Goed} (page \pageref{Definition
Mod-Fin-Goed})  in
Section \ref{Section Finite-Goedel} (page \pageref{Section Finite-Goedel}).
We show in Section \ref{Section Example-No-Int} (page \pageref{Section
Example-No-Int})  the well-known
result that such logics for 3 worlds (and thus 4 truth values)
have no interpolation, whereas the corresponding logic for 2 worlds
has interpolation. For the latter logic, we can still find a kind of
normal form, though $ \xcp $ cannot always be reduced. At least we can
avoid
nested implications, which is not possible in the logic for 3
worlds.

We also discuss several ``hand made'' additional operators which
allow us to define sufficiently many model sets to have syntactical
interpolation - of course, we $ \xCf know$ that we have semantical
interpolation.
A more systematic approach was discussed above, the operators $ \xbf!
\xbq $ and $ \xbf? \xbq.$
\section{
Monotone and antitone semantic interpolation
}

\label{Section Mon-Sem-Int}

We explain what is happening here. Assume $A(p,q)$ proves $B(q,r)$ where
the
set of models $M(A(p,q))$ is a subset of $M(B(q,r))$ $ \xCf or$ $ \xCf
vice$ $ \xCf versa,$ i.e.,
$B$ proves $ \xCf A$ and the subset relation is also the inverse. We
discuss here the
first variant. - The common language is $p.$ This means for all $p,q,r$
$(A(p,q)$ proves $B(q,r))$ or equivalently
for all q(for some $p$ $ \xCf A)$ proves (for all $r$ $B)).$ Semantically
this means, in the
monotonic case, the (union on $p$ models of $ \xCf A)$ is a subset of
(intersection on $r$ models of $B).$ We want to extract from this a set of
models of $q$ interpolating in between. This is what the set theoretical
manipulation below does. The result is formulated in
Proposition \ref{Proposition Sin-Interpolation} (page \pageref{Proposition
Sin-Interpolation}).

Once we find the semantic interpolant we ask under what conditions
can we find a syntactic $C$ to do the job. This we investigate in the
rest of the section.
In classical logic, the semantic result carries over immediately
to the syntactic level, as is shown in
Proposition \ref{Proposition Sin-Simplification-Definable} (page
\pageref{Proposition Sin-Simplification-Definable}).
$C$ can be found in many cases by enriching the language.
There are papers in the literature with a title
``repairing interpolation for logic $X$'',
e.g., by Areces, Blackburn, and Marx,
see  \cite{ABM03}, this is
what they do for some particular logic $X.$

Thus, the following interpolation results can be read upward
(monotonic logic) or downward (the core of non-monotonic logic, in the
following
sense:
$ \xbg $ is the theory of the minimal models of $ \xba,$ and not just any
formula
which holds in the set of minimal models - which would be downward,
and then upward again in the sense of model set inclusion), in the latter
case
we have to be careful: we usually cannot go upward again, so we have the
sharpest possible case in mind. The case of mixed movement - down and then
up -
as in full non-monotonic logic is treated
in Section \ref{Section Sem-Int-NML} (page \pageref{Section Sem-Int-NML}).

As a warming up exercise, we do first a simplified version of the
two-valued
case, giving only the lower bound.
Parts (2) and (3) of the following Proposition concern
``parallel interpolation'', see  \cite{KM07}.
\subsection{
The two-valued case
}

Recall that we can work here with sets of models, which are named $ \xbS $
etc.,
to remind us that model sets are sets of sequences. Part (2) and (3)
concern
``parallel interpolation'', a terminology used by Makinson et al. in
 \cite{KM07}. The proofs of these parts are straightforward
generalizations
of the simple case, they are mentioned for completeness' sake.

\bp

$\hspace{0.01em}$


\label{Proposition Sin-Interpolation}

Let $ \xbS' \xcc \xbS \xcc \xbP,$ where $ \xbP = \xbP \{X_{i}:i \xbe
L\}.$

Recall Definition \ref{Definition Sin-Ir-Relevant} (page \pageref{Definition
Sin-Ir-Relevant})  for
the definitions of $I$ and $R.$

(1) Let $ \xbS'':= \xbS' \xex (R(\xbS) \xcs R(\xbS')) \xCK \xbP
\xex (I(\xbS) \xcv I(\xbS')).$

Then $ \xbS' \xcc \xbS'' \xcc \xbS.$

The following two results
concern ``parallel interpolation'', terminology introduced by
D.Makinson in  \cite{KM07}. Thus in the first case, $ \xbS' $ is a
product, in the second
case, $ \xbS $ is a product. We do interpolation for a whole family of
partial lower
or upper bounds in parallel, thus its name.

(2) Let $ \xdj $ be a disjoint cover of $L.$

Let $ \xbS' = \xbP \{ \xbS'_{K}:K \xbe \xdj \}$ with $ \xbS'_{K} \xcc
\xbP \{X_{i}:i \xbe K\}.$

Let $ \xbS''_{K}:= \xbS'_{K} \xex (R(\xbS) \xcs R(\xbS'_{K})) \xCK
\xbP \{X_{i}:i \xbe K,i \xbe I(\xbS) \xcv I(\xbS'_{K})\}.$

Let $ \xbS'':= \xbP \{ \xbS''_{K}:K \xbe \xdj \}$ (re-ordered).

Then $ \xbS' \xcc \xbS'' \xcc \xbS.$

(3) Let $ \xdj $ be a disjoint cover of $L.$

Let $ \xbS = \xbP \{ \xbS_{K}:K \xbe \xdj \}$ with $ \xbS_{K} \xcc \xbP
\{X_{i}:i \xbe K\}.$

Let $ \xbS''_{K}:= \xbS' \xex (R(\xbS') \xcs R(\xbS_{K})) \xCK \xbP
\{X_{i}:i \xbe K,i \xbe I(\xbS_{K}) \xcv I(\xbS')\}.$

Let $ \xbS'':= \xbP \{ \xbS''_{K}:K \xbe \xdj \}$ (re-ordered).

Then $ \xbS' \xcc \xbS'' \xcc \xbS.$

\ep

\subparagraph{
Proof
}

$\hspace{0.01em}$


(1)

(1.1) $ \xbS' \xcc \xbS'' $ is trivial.

(1.2) $ \xbS'' \xcc \xbS:$

We can use Fact \ref{Fact Sin-Ir-Relevant} (page \pageref{Fact Sin-Ir-Relevant})
 (2), or argue directly. We
will do the
latter.

Let $m \xbe \xbS'',$ so, by definition, there is $m' \xbe \xbS' $ such
that
$m \xex (R(\xbS) \xcs R(\xbS'))=m' \xex (R(\xbS) \xcs R(\xbS'
)).$
Define $m'' $ as follows: On $(R(\xbS) \xcs R(\xbS')) \xcv I(\xbS'
),$ $m'' $ is like $m,$ on
the other $i,$ $m'' $ is like $m'.$ $m' $ differs from $m'' $ at most on
$I(\xbS'),$ so
by definition of $I(\xbS'),$ $m'' \xbe \xbS' \xcc \xbS.$ $m'' $ is
like $m$ at least on
$(R(\xbS) \xcs R(\xbS')) \xcv I(\xbS') \xcd R(\xbS),$ so by
definition of $R(\xbS),$ $m \xbe \xbS.$

(2)

(2.1) $ \xbS' \xcc \xbS''.$

$ \xbS'_{K} \xcc \xbS''_{K},$ so by $ \xbS'' = \xbP \xbS''_{K}$ the
result follows.

(2.2) $ \xbS'' \xcc \xbS.$

Let $m'' \xbe \xbS'',$ $m'' = \xDO \{m''_{K}:K \xbe \xdj \}$ for
suitable $m''_{K} \xbe \xbS''_{K},$ where $ \xDO $ stands for
the composition (or concatenation) of partial sequences or models.
Consider $m''_{K}.$ By definition of $ \xbS''_{K},$ there is $m'_{K} \xbe
\xbS'_{K}$ s.t.
$m''_{K} \xex (R(\xbS'_{K}) \xcs R(\xbS))=m'_{K} \xex (R(\xbS'_{K})
\xcs R(\xbS)),$ so there is $n'_{K} \xbe \xbS'_{K}$ s.t.
$m''_{K} \xex R(\xbS)=n'_{K} \xex R(\xbS).$
Let $n':= \xDO \{n'_{K}:K \xbe \xdj \},$ so by $ \xbS' = \xbP \{ \xbS
'_{K}:K \xbe \xdj \}$ $n' \xbe \xbS' \xcc \xbS.$
But $m'' \xex R(\xbS)=n' \xex R(\xbS),$ so $m'' \xbe \xbS.$

(3)

We first show $R(\xbS_{K})=R(\xbS) \xcs K.$

Let $i \xbe I(\xbS_{K}),$ then $ \xbS_{K}= \xbS_{K} \xex (K-\{i\}) \xCK
X_{i},$ but $ \xbS = \xbP \{ \xbS_{K}:K \xbe \xdj \},$ so
$ \xbS = \xbS \xex (L-\{i\}) \xCK X_{i},$ and $i \xbe I(\xbS).$
Conversely, let $i \xbe I(\xbS) \xcs K,$ then
$ \xbS = \xbS \xex (L-\{i\}) \xCK X_{i},$ so $ \xbS \xex K= \xbS \xex
(K-\{i\}) \xCK X_{i},$ so $i \xbe I(\xbS_{K}).$

(3.1) $ \xbS' \xcc \xbS''.$

$ \xbS' \xex K \xcc \xbS''_{K},$ so by $ \xbS'' = \xbP \{ \xbS''_{K}:K
\xbe \xdj \},$ $ \xbS' \xcc \xbS''.$

(3.2) $ \xbS'' \xcc \xbS.$

By $ \xbS = \xbP \{ \xbS_{K}:K \xbe \xdj \},$ it suffices to show $ \xbS
''_{K} \xcc \xbS_{K}.$

Let $m''_{K} \xbe \xbS''_{K}.$ So there is $m' \xbe \xbS' $ s.t.
$m' \xex (R(\xbS_{K}) \xcs R(\xbS'))=m''_{K} \xex (R(\xbS_{K}) \xcs
R(\xbS')),$ so there is $n' \xbe \xbS' \xcc \xbS $ s.t.
$n' \xex R(\xbS_{K})=m''_{K} \xex R(\xbS_{K}),$ so by $R(\xbS_{K})=R(
\xbS) \xcs K,$
there is $m \xbe \xbS $ s.t. $m \xex K=m''_{K},$ so $m''_{K} \xbe
\xbS_{K}.$

$ \xcz $
\\[3ex]
\subsection{
The many-valued case
}

For the basic definitions and facts, see
Section \ref{Section Many-Val-Intro} (page \pageref{Section Many-Val-Intro}).

We assume here that the max and min of arbitrary sets of truth values
will always exist. In
Section \ref{Section Inter-Arg} (page \pageref{Section Inter-Arg}), we consider
argumentation, and see the
set of
arguments for some formula $ \xbf $ as its truth value. There, only the
max
(i.e., union) will always exist, we will see there that this is also
sufficient
for some form of interpolation.

Recall that we do not work with sets of models or sequences any more, but
with arbitrary functions $f,g,h:M \xcp V,$ where each $m \xbe M$ is a
function $m:L \xcp V,$
where, intuitively, $L$ stands
for the set of propositional variables, $V$ for the set of truth values,
$M$ is the set of many-valued models,
and $f$ etc. are functions (intuitively, $f=f_{ \xbf }$ etc.) assigning
each $m$ a value,
intuitively, the value of $ \xbf $ in $m.$
Again, we will consider $f \xck g$ and look for some $h$ with $f \xck h
\xck g,$ where
$I(f) \xcv I(g) \xcc I(h).$

\bd

$\hspace{0.01em}$


\label{Definition Gin-Base}

(1) Let $J \xcc L,$ $f:M \xcp V.$ Define

$f^{+}(m \xex J):=max\{f(m'):m \xex J=m' \xex J\}$ and

$f^{-}(m \xex J):=min\{f(m'):m \xex J=m' \xex J\}.$

(Similarly, if $m$ is defined only
on $J,$ the condition is $m' \xex J=m,$ instead of $m \xex J=m' \xex J.)$

(2) Call $M$ rich iff for all $m,m' \xbe M,$ $J \xcc L,$ $(m \xex J) \xcv
(m' \xex (L-J)) \xbe M.$
(I.e., we may cut and paste models.)

This assumption is usually given, it is mainly here to remind the
reader that it is not trivial, and we have to make sure it really
holds. A problem might, e.g., arise when we consider only subsets of
all models, i.e., some $M' \xcc M,$ and not the full $M.$

Note that the possibility of arbitrary combinations of models is also an
aspect
of independence.

(3) A reminder: Call $f:M \xcp V$ insensitive to $J \xcc L$ iff for all
$m,n:$
$m \xex (L-J)=n \xex (L- \xCf J)$
implies $f(m)=f(n)$ - i.e., the values of $m$ on $J$ have no importance
for $f.$
See Section \ref{Section Many-Val-Intro} (page \pageref{Section Many-Val-Intro})
,
Table \ref{Table Gin-Not-Def} (page \pageref{Table Gin-Not-Def}).

(4)
We will sometimes write $m_{J}$ for $m \xex J.$

\ed

\br

$\hspace{0.01em}$


\label{Remark Gin-Base}

Let $J \xcc L,$ $m \xbe M,$ where $m:L \xcp V,$ and $f:M \xcp V.$

 \xEh

 \xDH
Obviously, if $J=L,$ then $f^{+}(m_{J})=f^{-}(m_{J})=f(m).$

 \xDH
We define the following ternary operators $+$ and -:

$+(f,m,J):=f^{+}(m_{J}) \xbe V,$ $-(f,m,J):=f^{-}(m_{J}) \xbe V.$

 \xDH
Usually, we fix $f$ and $J,$ and are interested in the value
$f^{+}(m_{J})$ and $f^{-}(m_{J})$
for various $m.$ Seen this way, we have binary operators

$+(f,J):M \xcp V,$ defined by $+(f,J)(m):=f^{+}(m_{J})$ and

$-(f,J):M \xcp V,$ defined by $-(f,J)(m):=f^{-}(m_{J}),$

i.e., the results $+(f,J)$ and $-(f,J)$ are new model functions like $f$
again.

 \xDH
When we fix now $J,$ we have unary functions $+(f)$ and $-(f)$ which
assign
to the old functions $f$ new functions $f^{+}$ and $f^{-}.$ This is
probably the
best way to consider our new operators.

By definition,
$f^{-}(m) \xck f(m) \xck f^{+}(m)$ for all $m.$

 \xDH
In the two-valued case, $f$ is a model set $X,$ and for any model $m,$ $m
\xbe X$ or $m \xce X.$
Fix now $J.$

In above sense, $m' \xbe X^{+}$ iff there is $m$ such that $m \xex J=m'
\xex J$ and $m \xbe X,$ and
$m' \xbe X^{-}$ iff for all $m$ such that $m \xex J=m' \xex J,$ also $m
\xbe X$ holds.

Thus, $X^{-} \xcc X \xcc X^{+}.$

 \xDH
When we thus fix $J,$
$+$ and - are new operators on model set functions, or on
model sets in the two-valued case. As such, they are similar to other
new operators, like the $ \xbm -$operator of preferential structures. But
they
are simpler, in the following sense: they do not require an additional
structure
like a relation $ \xeb,$ but are
``built into'' the model structure itself, as they only need the algebraic
structure already present in products.
Thus, in a certain way, they are more elementary.

But, just as $ \xbm $ need not preserve definability, neither do $+$ or -;
if $f=f_{ \xbf }$ for some $ \xbf,$ then there need not be $ \xbq $ such
that $f^{+}=f_{ \xbq }$ or $f^{-}=f_{ \xbq }.$

 \xEj

\er

Let $L=J \xcv J' \xcv J'' $ be a disjoint union. If $f:M \xcp V$ is
insensitive to $J \xcv J'',$ we
can define for $m_{J' }:J' \xcp V$ $f(m_{J' })$ as any $f(m')$ such that
$m' \xex J' =m_{J' }.$

\bp

$\hspace{0.01em}$


\label{Proposition Gin-Pr-Int}

Let $M$ be rich, $f,g:M \xcp V,$ $f(m) \xck g(m)$ for all $m \xbe M.$
Let $L=J \xcv J' \xcv J'',$ let $f$ be insensitive to $J,$ $g$ be
insensitive to $J''.$

Then $f^{+}(m_{J' }) \xck g^{-}(m_{J' })$ for all $m_{J' } \xbe M \xex J'
,$ and any $h:M \xex J' \xcp V$ which is
insensitive to $J \xcv J'' $ is an interpolant iff

$f^{+}(m_{J' }) \xck h(m_{J}m_{J' }m_{J'' })=h(m_{J' }) \xck g^{-}(m_{J'
})$ for all $m_{J' } \xbe M \xex J'.$

(h can be extended to the full $M$ in a unique way, as it is
insensitive to $J \xcv J'',$ so it does not really matter whether we
define
$h$ on $L$ or on $J'.)$

\ep

See Diagram \ref{Diagram Mon-Int} (page \pageref{Diagram Mon-Int}).

$ \xCO $

\vspace{10mm}

\begin{diagram}

\label{Diagram Mon-Int}
\index{Diagram Mon-Int}

\centering
\setlength{\unitlength}{1mm}
{\renewcommand{\dashlinestretch}{30}
\begin{picture}(150,100)(0,0)

\path(20,90)(90,90)(90,20)(20,20)(20,90)

\path(20,20)(20,18)
\path(90,20)(90,18)

\path(40,90)(40,18)
\path(70,90)(70,18)

\path(20,75)(90,75)
\path(20,45)(90,45)

\multiput(40,45)(-2,1){10}{\line(-2,1){0.7}}
\multiput(70,45)(2,1){10}{\line(2,1){0.7}}

\multiput(40,75)(-2,-1){10}{\line(-2,-1){0.7}}
\multiput(70,75)(2,-1){10}{\line(2,-1){0.7}}

\put(13,75){{\xssc $g(m)$}}
\put(9,65){{\xssc $g^-(m_{J'})$}}
\put(92,45){{\xssc $f(m)$}}
\put(92,55){{\xssc $f^+(m_{J'})$}}

\put(30,12){{\xssc $J$}}
\put(54,12){{\xssc $J'$}}
\put(80,12){{\xssc $J"$}}

\put(36,5){{\xssc Semantic interpolation, $f \xck g$}}

\end{picture}
}

\end{diagram}

\vspace{4mm}

$ \xCO $

\subparagraph{
Proof
}

$\hspace{0.01em}$


Let $L=J \xcv J' \xcv J'' $ be a pairwise disjoint union. Let $f$ be
insensitive to $J,$ $g$ be
insensitive to $J''.$

$h:M \xcp V$ will have to be insensitive to $J \xcv J'',$ so we will have
to define
$h$ on $M \xex J',$ the extension to $M$ is then trivial.

Fix arbitrary $m_{J' }:J' \xcp V,$ $m_{J' }=m \xex J' $ for some $m \xbe
M.$
We first show $f^{+}(m_{J' }) \xck g^{-}(m_{J' }).$

Proof:
Choose $m_{J'' }$ such that $f^{+}(m_{J' })=f(m_{J}m_{J' }m_{J'' })$ for
any $m_{J}.$
(Recall that $f$ is insensitive to $J.)$
Let $n_{J'' }$ be one such $m_{J'' }.$
Likewise,
choose $m_{J}$ such that $g^{-}(m_{J' })=g(m_{J}m_{J' }m_{J'' })$ for any
$m_{J'' }.$
Let $n_{J}$ be one such $m_{J}.$
Consider $n_{J}m_{J' }n_{J'' } \xbe M$ (recall that $M$ is rich).
By definition, $f^{+}(m_{J' })=f(n_{J}m_{J' }n_{J'' })$ and $g^{-}(m_{J'
})=g(n_{J}m_{J' }n_{J'' }),$ but
by prerequisite $f(n_{J}m_{J' }n_{J'' }) \xck g(n_{J}m_{J' }n_{J'' }),$ so
$f^{+}(m_{J' }) \xck g^{-}(m_{J' }).$

Thus, any $h$ such that $h$ is insensitive to $J \xcv J'' $ and

(Int) $f^{+}(m_{J' }) \xck h(m):=h(m_{J' }) \xck g^{-}(m_{J' })$

is an interpolant for $f$ and $g.$ The definition $h(m):=h(m_{J' })$ is
possible,
as $h$ is insensitive to $J \xcv J''.$

$f(m_{J}m_{J' }m_{J'' })$ $ \xck $ $h(m_{J}m_{J' }m_{J'' })$ $ \xck $
$g(m_{J}m_{J' }m_{J'' })$ follows
trivially, using above notation: $f(m_{J}m_{J' }m_{J'' })$ $ \xck $
$f(m_{J}m_{J' }n_{J'' })$ $=$
$f^{+}(m_{J' })$ $ \xck $ $h(m_{J' })$ $=$ $h(m_{J}m_{J' }m_{J'' })$ $
\xck $ $g^{-}(m_{J' })$ $=$ $g(n_{J}m_{J' }m_{J'' })$ $ \xck $
$g(m_{J}m_{J' }m_{J'' }).$

But (Int) is also a necessary condition.

Proof:

Suppose $h$ is insensitive to $J \xcv J'' $ and $h(m_{J}m_{J' }m_{J''
})=h(m_{J' })<f^{+}(m_{J' }).$
Let $n_{J'' }$ be as above, i.e.,
$f(m_{J}m_{J' }n_{J'' })=f^{+}(m_{J' })$ for any $m_{J}.$ Then
$h(m_{J}m_{J' }n_{J'' })=h(m_{J' })<f^{+}(m_{J' })=f(m_{J}m_{J' }n_{J''
}),$ so $h$ is not an interpolant.

The proof that $h(m_{J}m_{J' }m_{J'' })=h(m_{J' })$ has to be $ \xck
g^{-}(m_{J' })$ is analogous.

We summarize:

$f$ and $g$ have an interpolant $h,$ and
$h$ is an interpolant for $f$ and $g$ iff $h$ is insensitive to $J \xcv
J'' $ and
for any $m_{J' } \xbe M \xex J' $
$f^{+}(m_{J' }) \xck h(m_{J}m_{J' }m_{J'' })=h(m_{J' }) \xck g^{-}(m_{J'
}).$

$ \xcz $
\\[3ex]

\bd

$\hspace{0.01em}$


\label{Definition Standard-Int}

It is thus justified to call in above situation

$f^{+}(m_{J' })$ and $g^{-}(m_{J' })$ the standard interpolants,

more precisely, we call
$h$ such that $h(m)=f^{+}(m_{J' })$ or $h(m)=g^{-}(m_{J' })$ with $R(h)
\xcc J' $ a standard
interpolant.

\ed

It seems that the same technique can be used to show many-valued
semantic interpolation for modal and first-order logic, but we have
not checked in detail.
\section{
The interval of interpolants in monotonic or antitonic logics
}

\label{Section Int-Int-Mon}

\label{Section Mon-Interpol-Int}
\subsection{
Introduction
}

We take now a closer look at the interval of interpolants, with already
some remarks on the syntactic side.

By Proposition \ref{Proposition Gin-Pr-Int} (page \pageref{Proposition
Gin-Pr-Int}),
we have an interval of interpolants, and the extremes, the
standard interpolants,
see Definition \ref{Definition Standard-Int} (page \pageref{Definition
Standard-Int}),
are particularly
interesting, for the following reasons:

 \xEh

 \xDH
precisely because they are the extremes

 \xDH
they are universal in the sense that they depend only one of the two
functions
and the variable set, more precisely:

\bd

$\hspace{0.01em}$


\label{Definition Univ-Int}

Let $L=J \xcv J' \xcv J'',$ let $f:M \xcp V$ be insensitive to $J.$

(1) $h$ is an upper universal interpolant for $f$ and $J' $ iff for all
$g$ such that
$g$ is insensitive to $J'',$ and $f \xck g,$ $h$ is an interpolant for
$f$ and $g,$ i.e.,
$f \xck h \xck g$ and $h$ is insensitive to $J \xcv J''.$

(2) $h$ is a lower universal interpolant for $f$ and $J' $ iff for all $g$
such that
$g$ is insensitive to $J'',$ and $g \xck f,$ $h$ is an interpolant for
$f$ and $g,$ i.e.,
$g \xck h \xck f$ and $h$ is insensitive to $J \xcv J''.$

\ed

So $f^{+}(m_{J' })$ is (the only - provided we have sufficiently many $g$
to consider)
upper for $f,$ $g^{-}(m_{J' })$ the (again, only) lower interpolant for
$g.$

 \xDH
They have a particularly simple structure. E.g., we could also choose for
one $m_{J' }$ the lower boudary, for another $m'_{J' }$ the upper
boundary, etc., and
compose an interpolant in this somewhat arbitrary way.

 \xDH
The operators $+$ and - are themselves interesting and elementary.

 \xDH
It might be difficult to find necessary and sufficient criteria for the
existence of some syntactic interpolant in the interval, whereas it might
be possible to find such criteria for the existence of one or both of the
universal interpolants.

 \xDH
In particular, we can ask the following questions:

If $f$ is given by some $ \xbf,$ i.e., $f=f_{ \xbf },$ and $J' \xcc L$ as
above, defining
$f^{+}_{ \xbf }(m):=f^{+}_{ \xbf }(m_{J' })$ and $f^{-}_{ \xbf
}(m):=f^{-}_{ \xbf }(m_{J' }),$ can we find formulas
$ \xbq $ and $ \xbq' $ containing only variables from $J' $
such that $f_{ \xbq }=f^{+}_{ \xbf }$ and $f_{ \xbq' }=f^{-}_{ \xbf }$?
In more detail:

 \xEI

 \xDH
Is this possible in classical logic?

 \xDH
Is this possible in other logics?

 \xDH
Do we find criteria for the language (operators, truth values, etc.) which
guarantee that such $ \xbq $ and $ \xbq' $ exist?

 \xEJ

 \xEj

We turn to some examples and a simple fact.
\subsection{
Examples and a simple fact
}

\be

$\hspace{0.01em}$


\label{Example Min=Max}

This trivial, many-valued example shows that
infimum and supremum might coincide. Take 3 truth values, $0,1,2,$ in this
order, and 3 propositional variables, $p,q,r.$ We want to fix $q,$ and let
$p$ and $r$ float. $ \xbf:=p \xcu q \xcu r,$ $ \xbf':=p \xco q \xco r,$
so $ \xbf \xck \xbf'.$ We fix $q$ at 1,
so $sup\{F_{ \xbf }(m):m(q)=1\}=inf\{F_{ \xbf' }(m):m(q)=1\}=1.$

\ee

We present now two examples for the interval of interpolants.

\be

$\hspace{0.01em}$


\label{Example Int-Int-Class}

Consider classical logic with 4 propositional variables, $p,q,r,s.$

Let $ \xbf:=p \xcu q \xcu r,$ $ \xbq:=q \xco r \xco s.$ Obviously, $
\xbf \xcm \xbq.$

Let $f:=f_{ \xbf }$ assign to a model $m$
the truth value $ \xbf $ has in $m,$ likewise for $g:=f_{ \xbq }$ for $
\xbq.$

Let $f' (m):=sup\{f(m'):m \xex \{q,r\}=m' \xex \{q,r\}\},$
$g' (m):=inf\{g(m'):m \xex \{q,r\}=m' \xex \{q,r\}\},$
then $f' (m)=f_{q \xcu r}(m),$ $g' (m)=f_{q \xco r}(m).$

$f' $ and $g' $ are the bounds of the
interval, and, e.g., $q$ and $r$ are both inside the interval, we have
$f' \xck f_{q},f_{r} \xck g'.$

\ee

\be

$\hspace{0.01em}$


\label{Example Int-Int-Many}

(See Section \ref{Section Finite-Goedel} (page \pageref{Section Finite-Goedel}) 
for motivation.)

Consider a finite intuitionistic Goedel language over $p,q,r,$ with two
worlds,
and an additional operator $E \xbf,$ which says that $ \xbf $ holds
everywhere.
We have 3 truth values, $0<1<2,$ $ \xcu $ and $ \xco $ with the usual
interpretation of inf
and sup, intuitionistic negation with $ \xCN \xbf $ has value 2 iff $ \xbf
$ has value 0,
and 0 otherwise,
so $ \xCN \xCN \xbf $ has value 0 iff $ \xbf $ has value 0, and 2
otherwise,
and $E \xbf $ has value 2 iff $ \xbf $ has value 2, and 0 otherwise.

Consider $ \xbf:=p \xcu Eq$ and $ \xbq:= \xCN \xCN q \xco r.$ Obviously,
$ \xbf \xcm \xbq.$

$f:=f_{ \xbf }$ assigns to a model $m$
the truth value $ \xbf $ has in $m,$ likewise for $g:=f_{ \xbq }$ for $
\xbq.$

Let $f' (m):=sup\{f(m'):m \xex \{q\}=m' \xex \{q\}\},$
$g' (m):=inf\{g(m'):m \xex \{q\}=m' \xex \{q\}\}.$
It suffices to consider the truth values for $q:$
If $m(q)=0,$ then $f' (m)=g' (m)=0,$
if $m(q)=1,$ then $f' (m)=0,$ $g' (m)=2,$
If $m(q)=2,$ then $f' (m)=g' (m)=2.$

So $f' (m)=f_{Eq}(m),$ $g' (m)=f_{ \xCN \xCN q}(m).$

$f' $ and $g' $ are the bounds of the
interval, and, $q$ is inside the interval, we have
$f' \xck f_{q} \xck g'.$
\subsection{
The analoga of + and - as new semantic and syntactic operators
}

\label{Section Analoga}
\subsubsection{
Motivation
}

\ee

We introduced and discussed the operators $+$ and - in
Definition \ref{Definition Gin-Base} (page \pageref{Definition Gin-Base})  and
Remark \ref{Remark Gin-Base} (page \pageref{Remark Gin-Base}).
In the binary version,
see Remark \ref{Remark Gin-Base} (page \pageref{Remark Gin-Base}),
they were $+(f,J)$ and $-(f,J),$ where $f$ may be the
model function of some formula $ \xbf,$ $f=f_{ \xbf }.$ So, intuitively,
we can
read $+(\xbf,J)$ and $-(\xbf,J)$ as new syntactic operators.

A usual binary operator between formulas has
two formulas as arguments.
$J,$ however, is just a set of propositional variables. Of course, we may
take the conjunction of $J$ in the finite case. As formulas are finite,
most variables are ``free'' anyway, so finiteness is not a real restriction.
Still, so far, one of the formulas has to be a conjunction of variables,
which is bizarr for a normal operator. We have to generalize, but are
relatively
free how we do this, as we really just need the set of variables. The
perhaps
simplest way to do this is to take $R(\xbq),$ the set of relevant
variables
of $ \xbq.$ Then $+(\xbf, \xbq):=f^{+}_{ \xbf }(m \xex R(\xbq))$ and
$-(\xbf, \xbq):=f^{-}_{ \xbf }(m \xex R(\xbq))$ look more
respectable. Note, however, that, e.g., $R(p \xcu q)=R(p \xco q),$ so the
precise form
of $ \xbq $ enters only slightly into the picture.

We make this formal now.
\subsubsection{
Formal definition and results
}

Let $J \xcc L,$ $X \xcc M.$

\bd

$\hspace{0.01em}$


\label{Definition 2-Many-Up}

The supremum (the $ \xCf lower$ bound in interpolation):

 \xEh
 \xDH
The 2-valued case:

 \xEh
 \xDH
$X \xfB J$ $:=$ $\{m \xbe M$: $ \xcE m' (m \xex J=m' \xex J$ $ \xcu $ $m'
\xbe X)\}$
 \xDH
$M(\xbf! \xbq)$ $:=$ $M(\xbf) \xfB (L-R(\xbq))$

(It is another question, whether $M(\xbf! \xbq)$ is definable - unless,
of course,
we give the new operator full right as an operator in the language, then
$M(\xbf! \xbq)$ is,
by definition, definable. This will also hold in the other cases.)

Remark: $ \xbf $ is a finite formula, so $R(\xbf)$ is finite, so whether
we fix
$M(\xbf)$ on a finite or infinite set $J$ does not matter. Thus, we
could just
as well have defined $M(\xbf! \xbq)=M(\xbf) \xfB R(\xbq).$ This is
a matter of taste.

 \xEj

 \xDH
The many-valued case:

Given $f,$ define

 \xEh
 \xDH
$(f \xfB J)(m)$ $:=$ $sup\{f(m'):$ $m' \xbe M,$ $m \xex J=m' \xex J\}$
 \xDH
$f_{ \xbf! \xbq }$ $:=$ $f_{ \xbf } \xfB (L-R(\xbq))$
 \xEj
 \xEj

\ed

\bd

$\hspace{0.01em}$


\label{Definition 2-Many-Down}

The infimum (the $ \xCf upper$ bound in interpolation):

 \xEh
 \xDH
The 2-valued case:

 \xEh
 \xDH
$X \xfb J$ $:=$ $\{m \xbe M$: $ \xcA m' (m \xex J=m' \xex J$ $ \xch $ $m'
\xbe X)\}$
 \xDH
$M(\xbf? \xbq)$ $:=$ $M(\xbf) \xfb (L-R(\xbq))$
 \xEj

 \xDH
The many-valued case:

Given $f,$ define

 \xEh
 \xDH
$(f \xfb J)(m)$ $:=$ $inf\{f(m'):$ $m' \xbe M,$ $m \xex J=m' \xex J\}$
 \xDH
$f_{ \xbf? \xbq }$ $:=$ $f_{ \xbf } \xfb (L-R(\xbq))$
 \xEj
 \xEj

\ed

\br

$\hspace{0.01em}$


\label{Remark Up-Down}

(1) $ \xbf! \xbq $ and $ \xbf? \xbq $ will give immediately syntactic
interpolation,
when they are defined (or equivalent to another formula).

(2) The existence (or equivalence) of $ \xbf! \xbq $ and $ \xbf? \xbq $
can be seen as
as well-behaviour of a language and logic, as simplified model sets,
obtained by simple algebraic operations, are again definable.

$ \xcz $
\\[3ex]
\subsubsection{
The special case of classical logic.
}

\label{Section Class-Proj}

\er

\bfa

$\hspace{0.01em}$


\label{Fact Class-Up-Down-Pre}

(1) $(X \xcv X') \xfB J$ $=$ $(X \xfB J)$ $ \xcv $ $(X' \xfB J)$

(2) $(X \xcs X') \xfb J$ $=$ $(X \xfb J)$ $ \xcs $ $(X' \xfb J)$

(3) Let $ \xbf $ and $ \xbq $ be consistent, but no tautologies, and let $
\xbf $ and $ \xbq $
have no common variables. Let $J$ contain all
variables in $ \xbf,$ but no variable in $ \xbq.$

Then $M(\xbf \xcu \xbq) \xfB J$ $=$ $M(\xbf),$ and $M(\xbf \xco \xbq
) \xfb J$ $=$ $M(\xbf).$

\efa

\subparagraph{
Proof
}

$\hspace{0.01em}$


(1) is shown in Fact \ref{Fact 2-Up} (page \pageref{Fact 2-Up}), (8),
(2) is shown in Fact \ref{Fact 2-Down} (page \pageref{Fact 2-Down}), (9).
As the arguments are simple, we give them here already, for the reader's
convenience:
(1): Let $m \xbe (X \xcv X') \xfB J,$
so there is $m' \xbe X \xcv X' $ with $m \xex J=m' \xex J,$ so $m \xbe (X
\xfB J) \xcv (X' \xfB J).$ The converse is
even easier.
(2): Let $m \xbe (X \xcs X') \xfb J,$ so for all $m' $ such that $m' \xex
J=m \xex m$ $m' \xbe X \xcs X',$ so
$m \xbe (X \xfb J) \xcs (X' \xfb J).$ Again, the converse is even easier.

We turn to (3).

We first show $M(\xbf \xcu \xbq) \xfB J$ $=$ $M(\xbf).$
$m \xbe M(\xbf \xcu \xbq) \xfB J$ iff there is $m',$ $m \xex J=m' \xex
J,$ $m' \xcm \xbf \xcu \xbq.$ Suppose
$m \xbe M(\xbf \xcu \xbq) \xfB J,$ but $m \xcM \xbf.$ Then any $m' $
such that $m \xex J=m' \xex J$ will also
make $ \xbf $ false, contradiction. Conversely, let $m \xcm \xbf.$ Then
we can modify $m$
outside $J$ to make $ \xbq $ true, so $m \xbe M(\xbf \xcu \xbq) \xfB J.$

We now show $M(\xbf \xco \xbq) \xfb J$ $=$ $M(\xbf).$
$m \xbe M(\xbf \xco \xbq) \xfb J$ iff for all $m' $ such that $m' \xex
J=m \xex J,$ $m' \xcm \xbf \xco \xbq.$
Suppose $m \xbe M(\xbf \xco \xbq) \xfb J,$ but $m \xcM \xbf.$ Change
now $m$ outside $J$ to make $ \xbq $ false,
this is possible, contradiction. Conversely, let $m \xcm \xbf.$ Then, no
matter how we
change $m$ outside $J,$ $m' $ will still be a model of $ \xbf,$ and thus
of $ \xbf \xco \xbq.$

$ \xcz $
\\[3ex]

We use the disjunctive and conjunctive normal forms in classical logic
to obtain the following result:

\bfa

$\hspace{0.01em}$


\label{Fact Class-Up-Down}

Let the logic be classical propositional logic.

Let $J:=\{p_{1,1}, \Xl,p_{1,m_{1}}, \Xl,p_{n,1}, \Xl,p_{n,m_{n}}\}$

(1) Let $ \xbf_{i}:= \xCL p_{i,1} \xcu  \Xl  \xcu \xCL p_{i,m_{i}}$ and $
\xbq_{i}:= \xCL q_{i,1} \xcu  \Xl  \xcu \xCL q_{i,k_{i}},$
let $ \xbf:=(\xbf_{1} \xcu \xbq_{1}) \xco  \Xl  \xco (\xbf_{n} \xcu
\xbq_{n}).$
Then $ \xbf \xfB J= \xbf_{1} \xco  \Xl  \xco \xbf_{n}.$

(2) Let $ \xbf_{i}:= \xCL p_{i,1} \xco  \Xl  \xco \xCL p_{i,m_{i}}$ and $
\xbq_{i}:= \xCL q_{i,1} \xco  \Xl  \xco \xCL q_{i,k_{i}},$
let $ \xbf:=(\xbf_{1} \xco \xbq_{1}) \xcu  \Xl  \xcu (\xbf_{n} \xco
\xbq_{n}).$
Then $ \xbf \xfb J= \xbf_{1} \xcu  \Xl  \xcu \xbf_{n}.$

\efa

\subparagraph{
Proof
}

$\hspace{0.01em}$


(1)
By Fact \ref{Fact Class-Up-Down-Pre} (page \pageref{Fact Class-Up-Down-Pre}) 
(1)
$M(\xbf) \xfB J$ $=$ $M(\xbf_{1} \xcu \xbq_{1}) \xfB J$ $ \xcv $  \Xl
$ \xcv $ $M(\xbf_{n} \xcu \xbq_{n}) \xfB J.$
By Fact \ref{Fact Class-Up-Down-Pre} (page \pageref{Fact Class-Up-Down-Pre}) 
(3), $M(\xbf_{i} \xcu \xbq_{i})
\xfB J$ $=$ $M(\xbf_{i}).$

(2)
By Fact \ref{Fact Class-Up-Down-Pre} (page \pageref{Fact Class-Up-Down-Pre}) 
(2)
$M(\xbf) \xfb J$ $=$ $M(\xbf_{1} \xco \xbq_{1}) \xfb J$ $ \xcs $  \Xl
$ \xcs $ $M(\xbf_{n} \xco \xbq_{n}) \xfb J.$
By Fact \ref{Fact Class-Up-Down-Pre} (page \pageref{Fact Class-Up-Down-Pre}) 
(3), $M(\xbf_{i} \xco \xbq_{i})
\xfb J$ $=$ $M(\xbf_{i}).$

$ \xcz $
\\[3ex]
\subsubsection{
General results on the new operators
}

\label{Section General-New-Oper}

We turn to some general results on the new operators.
Note that we work with some redundancy, as the positive cases for the
many-valued situation imply the 2-valued result. Still, as the operators
are somewhat unusual, we prefer to proceed slowly, and first treat the
two-valued case.

\bfa

$\hspace{0.01em}$


\label{Fact 2-Up}

Consider the 2-valued case.

(1) $X$ $ \xcc $ $X \xfB J$ $ \xcc $ $M.$

(2) $X \xfB L$ $=$ $X.$

(3) $X \xfB \xCQ $ $=$ $M$ iff $X \xEd \xCQ.$

(4) $ \xCQ \xfB J$ $=$ $ \xCQ.$

(5) $M \xfB J$ $=$ $M.$

(6) $X \xcc X' $ $ \xch $ $X \xfB J \xcc X' \xfB J.$

(7) $J \xcc J' $ $ \xch $ $X \xfB J' \xcc X \xfB J.$

(8) $(X \xcv X') \xfB J$ $=$ $(X \xfB J) \xcv (X' \xfB J).$

(9) $(X \xcs X') \xfB J$ $ \xcc $ $(X \xfB J) \xcs (X' \xfB J).$
The converse (i.e., $ \xcd)$ is not always true.

(10) $ \xdC (X \xfB J)$ $ \xcc $ $(\xdC X) \xfB J.$
The converse is not always true.

(11) $X \xfB (J \xcv J')$ $ \xcc $ $(X \xfB J) \xcs (X \xfB J').$
The converse is not always true.

(12) $(X \xfB J) \xcv (X \xfB J')$ $ \xcc $ $X \xfB (J \xcs J').$
The converse is not always true.

(13) In general, $ \xdC (X \xfB J) \xcC X \xfB (\xdC J).$
In general, $X \xfB (\xdC J) \xcC \xdC (X \xfB J).$

(14) $(X \xfB J) \xfB J' $ $=$ $X \xfB (J \xcs J').$

(15) $(\xbf!p)!q= \xbf!p \xcu q.$

(16) $M(\xbf!TRUE)$ $=$ $M(\xbf).$

(17)
Let $ \xbg_{X}$ be the characteristic function for $X.$
Then $ \xbg_{X \xfB J}(m)$ $=$ $sup\{ \xbg_{X}(m'):$ $m' \xbe M,$ $m \xex
J=m' \xex J\}.$

\efa

\subparagraph{
Proof
}

$\hspace{0.01em}$


$(1)-(4)$ and $(6)-(8)$ are trivial.

(5) follows from (1).

(9) follows from (6).

For the converse:
Let $L:=\{p,q\},$ $J:=\{p\}.$
Let $m:=p \xcu q,$ $m':=p \xcu \xCN q,$ $X:=\{m\},$ $X':=\{m' \},$ so $X
\xcs X' = \xCQ,$
$m \xex J=m' \xex J,$ but $m,m' \xbe (X \xfB J) \xcs (X' \xfB J).$

(10)
Let $m \xbe M,$ $m \xce X \xfB J$ $ \xch $ $ \xcA m' \xbe X.m \xex J \xEd
m' \xex J$ $ \xch $ $m \xce X$ $ \xch $ $m \xbe \xdC X \xcc ((\xdC X)
\xfB J).$

For the converse:
Let $L:=\{p,q\},$ $J:=\{p\}.$
Let $X=\{p \xcu q\}.$
So $X \xfB J=\{p \xcu q,p \xcu \xCN q\},$
$ \xdC (X \xfB J)$ $=$ $\{ \xCN p \xcu q, \xCN p \xcu \xCN q\},$ $ \xdC
(X)$ $=$ $\{p \xcu \xCN q, \xCN p \xcu q, \xCN p \xcu \xCN q\},$
$(\xdC X) \xfB J=M.$

(11) follows from (7).

For the converse:
Let $L:=\{p,q\},$ $J:=\{p\},$ $J':=\{q\}.$
Let $X:=\{p \xcu q, \xCN p \xcu \xCN q\},$ so $p \xcu \xCN q \xbe (X \xfB
J) \xcs (X \xfB J'),$ but $p \xcu \xCN q \xce X=X \xfB (J \xcv J').$

(12) follows from (7).

For the converse:
Let $L:=\{p,q\},$ $J:=\{p\},$ $J':=\{q\}.$
Let $X:=\{p \xcu q\}.$
$X \xfB \xCQ =M$ by (3), but $ \xCN p \xcu \xCN q \xce (X \xfB J) \xcv (X
\xfB J').$

(13)
Consider $L:=\{p,q\},$ $J:=\{p\},$ $X:=\{p \xcu q\}.$ Then
$X \xfB J=\{p \xcu q,p \xcu \xCN q\},$ so $ \xdc (X \xfB J)=\{ \xCN p \xcu
q, \xCN p \xcu \xCN q\}.$
$X \xfB \xdc J=X \xfB \{q\}=\{p \xcu q, \xCN p \xcu q\}.$

(14)

``$ \xcd $'': Let $m \xbe X \xfB (J \xcs J'),$ so there is $m' \xbe X,$
$m \xex J \xcs J' $ $=$ $m' \xex J \xcs J'.$
Define $m'' $ on $J$ like $m',$ on $L-J$ like $m.$
(Note that we combine here arbitrarily!)
So, as $m' \xbe X,$ $m'' \xbe X \xfB J.$ We want to show $m \xbe (X \xfB
J) \xfB J',$ it suffices to
show $m'' \xex J' =m \xex J'.$ But $m'' $ is on $L-J$ like $m,$ on $J
\xcs J' $ like $m',$
where $m' $ is like $m,$ so $m'' \xex J' =m \xex J'.$

``$ \xcc $'': Let $m \xbe (X \xfB J) \xfB J',$ we have to find $m' \xbe
X.m \xex (J \xcs J')=m' \xex (J \xcs J').$ As
$m \xbe (X \xfB J) \xfB J',$ there is $m_{0} \xbe X \xfB J,$ $m_{0} \xex
J' =m \xex J'.$ As $m_{0} \xbe X \xfB J,$ there is $m_{1} \xbe X,$
$m_{1} \xex J=m_{0} \xex J.$ So $m_{1} \xex (J \xcs J')=m_{0} \xex (J
\xcs J')=m \xex (J \xcs J').$

(15)
Let $X:=M(\xbf).$
$(X! \xCf p)! \xCf q$ $=$ $(X \xfB (L-\{p\})) \xfB (L-\{q\})$ $=$ (by
(14)) $X \xfB (L-\{p,q\})$ $=$ $ \xbf!p \xcu q.$

(16) $M(\xbf!TRUE)$ $=$ $M(\xbf) \xfB L$ $=$ $M(\xbf).$

(17)
$ \xbg_{X \xfB J}(m)=1$ iff $m \xbe X \xfB J$ iff $ \xcE m' \xbe X.m \xex
J=m' \xex J$ iff
$ \xcE m' (\xbg_{X}(m')=1$ and $m \xex J=m' \xex J)$ iff $sup\{
\xbg_{X}(m'):$ $m' \xbe M,$ $m \xex J=m' \xex J\}$ $=$ 1.

$ \xcz $
\\[3ex]

\bfa

$\hspace{0.01em}$


\label{Fact Many-Up}

Consider the many-valued case.

Let $0:M \xcp V$ be the constant 0 $(=$ minimal value) function,
$1:M \xcp V$ be the constant 1 $(=$ maximal value) function.

The same counterexamples as for the 2-valued case work, we just
take the characteristic functions.

(1) $f$ $ \xck $ $f \xfB J$ $ \xck $ 1.

(2) $f \xfB L$ $=$ $f.$

(4) $0 \xfB J$ $=$ 0.

(5) $1 \xfB J$ $=$ 1.

(6) $f \xck g$ $ \xch $ $f \xfB J \xck g \xfB J.$

(7) $J \xcc J' $ $ \xch $ $f \xfB J' \xck f \xfB J.$

(8) $(sup(f,g)) \xfB J$ $=$ $sup(f \xfB J,g \xfB J).$

(9) $(inf(f,g)) \xfB J$ $ \xck $ $inf(f \xfB J,g \xfB J).$
The converse is not always true.

(11) $f \xfB (J \xcv J')$ $ \xck $ $inf(f \xfB J,f \xfB J').$
The converse is not always true.

(12) $sup(f \xfB J,f \xfB J')$ $ \xck $ $f \xfB (J \xcs J').$
The converse is not always true.

(14) $(f \xfB J) \xfB J' $ $=$ $f \xfB (J \xcs J').$

(15) $(\xbf!p)!q= \xbf!p \xcu q.$

(16) $f_{ \xbf!TRUE}$ $=$ $f_{ \xbf }.$

\efa

\subparagraph{
Proof
}

$\hspace{0.01em}$


$(1)-(2),$ (4), $(6)-(8)$ are trivial.

(5) follows from (1).

(9) follows from (6).

(11) by (7).

(12) by (7).

(14)

For simplicity, we write just $f \xfB J \xcs J' (m)$ for $(f \xfB J \xcs
J')(m),$
$(f \xfB J) \xfB J' (m)$ for $((f \xfB J) \xfB J')(m),$ etc.

``$ \xcg $'': Consider $f \xfB J \xcs J' (m).$ We have to show $(f \xfB
J) \xfB J' (m) \xcg f \xfB (J \xcs J')(m).$
Let $m' $ be such that $m' \xex J \xcs J' =m \xex J \xcs J',$ $f \xfB J
\xcs J' (m)=f(m')$
(i.e., $m' $ has maximal value among such $m').$
(Use again arbitrary combinations.)
Define $m'' $ on $J$ like $m',$ on $L-J$ like $m,$ so, as $m'' $ is like
$m' $ on $J,$
$f \xfB J(m'')$ $:=$ $max\{f(m_{0}):m_{0} \xex J=m'' \xex J\} \xcg f(m'
).$ Note that $m'' $ is on $L-J$ like $m,$
on $J \xcs J' $ like $m',$ where $m' $ is like $m,$ so $m'' \xex J' =m
\xex J',$ so
$(f \xfB J) \xfB J' (m)$ $:=$ $max\{f \xfB J(m_{0}):m_{0} \xex J' =m \xex
J' \}$ $ \xcg $ $f \xfB J(m'').$
Thus, $(f \xfB J) \xfB J' (m)$ $ \xcg $ $f \xfB J(m'')$ $ \xcg $ $f(m')$
$=$ $f \xfB J \xcs J' (m).$

``$ \xck $'': Consider $(f \xfB J) \xfB J' (m).$ Let $m' $ be such that
$m \xex J' =m' \xex J',$
$(f \xfB J) \xfB J' (m)=(f \xfB J)(m'),$ i.e., $m' $ has maximal value.
Let $m'' $ be such that $m'' \xex J=m' \xex J,$ and $(f \xfB J)(m')=f(m''
),$ i.e., $m'' $ has maximal
value. Then $m'' \xex J \xcs J' =m' \xex J \xcs J' =m \xex J \xcs J',$ so
$f \xfB J \xcs J' (m)$ $ \xcg $ $f(m'')$ $=$ $f \xfB J(m')$ $=$ $(f \xfB
J) \xfB J' (m).$

(15)
Let $f:=f_{ \xbf }.$
$(f! \xCf p)! \xCf q$ $=$ $(f \xfB (L-\{p\})) \xfB (L-\{q\})$ $=$ (by
(14)) $f \xfB (L-\{p,q\})$ $=$ $f!p \xcu q.$

(16) $f_{ \xbf!TRUE}$ $=$ $f_{ \xbf } \xfB L$ $=$ $f_{ \xbf }.$

$ \xcz $
\\[3ex]

\bfa

$\hspace{0.01em}$


\label{Fact 2-Down}

Consider the 2-valued case.

(1) $X \xfb J \xcc X.$

(2) $X \xfb L=X.$

(3) $X \xfb \xCQ $ $=$ $ \xCQ $ iff $X \xEd M,$ and $M$ iff $X=M.$

(4) $ \xCQ \xfb J= \xCQ.$

(5) $M \xfb J=M.$

(6) $X \xcc X' $ $ \xch $ $X \xfb J \xcc X' \xfb J.$

(7) $J \xcc J' $ $ \xch $ $X \xfb J \xcc X \xfb J'.$

(8) $(X \xfb J) \xcv (X' \xfb J)$ $ \xcc $ $(X \xcv X') \xfb J.$
The converse is not always true.

(9) $(X \xcs X') \xfb J$ $=$ $(X \xfb J) \xcs (X' \xfb J).$

(10) $(\xdC X) \xfb J$ $ \xcc $ $ \xdC (X \xfb J).$
The converse is not always true.

(11) $(X \xfb J) \xcv (X \xfb J')$ $ \xcc $ $X \xfb (J \xcv J').$
The converse is not always true.

(12) $X \xfb (J \xcs J)$ $ \xcc $ $(X \xfb J) \xcs (X \xfb J').$
The converse is not always true.

(13) In general, $ \xdC (X \xfb J) \xcC X \xfb (\xdC J)$ and $X \xfb (
\xdC J) \xcC \xdC (X \xfb J).$

(14) $(X \xfb J) \xfb J' $ $=$ $X \xfb (J \xcs J')$

(15) $(\xbf?p)?q= \xbf?p \xcu q.$

(16) $M(\xbf?TRUE)$ $=$ $M(\xbf).$

(17)
Let $ \xbg_{X}$ be the characteristic function for $X.$
Then $ \xbg_{X \xfb J}(m)$ $=$ $inf\{ \xbg_{X}(m'):$ $m' \xbe M,$ $m \xex
J=m' \xex J\}.$

\efa

\subparagraph{
Proof
}

$\hspace{0.01em}$


$(1)-(6)$ are trivial.

(7)
$m \xbe X \xfb J$ $ \xch $ $ \xcA m' (m \xex J=m' \xex J$ $ \xch $ $m'
\xbe X).$ But $m \xex J' =m' \xex J' $ $ \xch $ $m \xex J=m' \xex J$ $
\xch $ $m' \xbe X,$
so $m \xbe X \xfb J'.$

(8) by (6).

For the converse:
Let $L:=\{p,q\},$ $J:=\{p\}.$
Let $X=\{p \xcu q\}.$
Let $X':=\{p \xcu \xCN q\}.$ So $X \xfb J=X' \xfb J= \xCQ,$ $(X \xcv X'
) \xfb J$ $=$ $\{p \xcu q,p \xcu \xCN q\}.$

(9)

``$ \xcc $'' by (6).

``$ \xcd $'' Let $m \xbe (X \xfb J) \xcs (X' \xfb J),$ so all $m' $ such
that $m' \xex J=m \xex J$ are in $X,$ and in $X',$
so they are in $X \xcs X',$ so $m \xbe (X \xcs X') \xfb J.$

(10)
$m \xbe (\xdC X) \xfb J$ $ \xch $ $m \xce X$ $ \xch $ $ \xcE m' (m \xex
J=m' \xex J$ $ \xcu $ $m' \xce X),$ so $m \xce X \xfb J.$

For the converse:
Let $L:=\{p,q\},$ $J:=\{p\}.$
Let $X=\{p \xcu q\}.$
$X \xfb J= \xCQ,$ so $ \xdC (X \xfb J)=M.$ $ \xdC (X)$ $=$ $\{p \xcu \xCN
q, \xCN p \xcu q, \xCN p \xcu \xCN q\},$
$(\xdC X) \xfb J$ $=$ $\{ \xCN p \xcu q, \xCN p \xcu \xCN q\}.$

(11) by (7).

For the converse:
Let $L:=\{p,q\},$ $J:=\{p\},$ $J':=\{q\}.$
Let $X=\{p \xcu q\}.$
$X \xfb (J \xcv J')=X,$ $X \xfb J=X \xfb J' = \xCQ.$

(12) by (7).

For the converse:
Let $L:=\{p,q\},$ $J:=\{p\},$ $J':=\{q\}.$
Let $X=\{p \xcu q,p \xcu \xCN q, \xCN p \xcu q\}.$
$X \xfb \xCQ = \xCQ $ by (3). $X \xfb J=\{p \xcu q,p \xcu \xCN q\},$ $X
\xfb J' =\{p \xcu q, \xCN p \xcu q\},$
$(X \xfb J) \xcs (X \xfb J')=\{p \xcu q\}.$

(13)

``$ \xcC $'': Let $J:= \xCQ,$ $X \xEd M,$ so $ \xdC J=L,$ $X \xfb \xdC
J=X,$ as $X \xfb \xCQ = \xCQ,$ $ \xdC (X \xfb J)=M.$

``$ \xcD $'': Let $J:= \xCQ.$ $M \xfb (\xdC J)=M,$ $M \xfb J=M,$ so $
\xdC (M \xfb J)= \xCQ.$

(14)

``$ \xcc $'': Let $m \xbe (X \xfb J) \xfb J',$ we have to show $m \xbe X
\xfb (J \xcs J'),$ i.e.,
for all $m'' $ such that $m'' \xex J \xcs J' =m \xex J \xcs J' $ $m'' \xbe
X.$ Take $m'' $ such that
$m \xex (J \xcs J')=m'' \xex (J \xcs J').$
Define $m' $ such that $m' $ on $J' $ is like $m,$ $m' $ on $L-J' $ is
like $m''.$
Then by $m \xbe (X \xfb J) \xfb J' $ and $m \xex J' =m' \xex J',$ $m'
\xbe X \xfb J.$ Note that
$m'' $ on $J$ is like $m':$ On $L-J' $ $m' $ is like $m''.$
On $J \xcs J',$ $m' $ is like $m,$ and $m$ is like $m'',$ so $m'' $ is
like $m'.$
As $m' \xbe X \xfb J,$ $m'' \xex J=m' \xex J,$ $m'' \xbe X$ by definition.

``$ \xcd $'': Let $m \xbe X \xfb (J \xcs J').$ So, if $m' \xex (J \xcs
J')=m \xex (J \xcs J'),$ $m' \xbe X.$
We have to show $m'' \xex J' =m \xex J' $ $ \xch $ $m'' \xbe X \xfb J,$
i.e. if $m'' \xex J' =m \xex J' $ and
$m_{0} \xex J=m'' \xex J,$ then $m_{0} \xbe X.$
But $m_{0} \xex (J \xcs J')=m'' \xex (J \xcs J')=m \xex (J \xcs J'),$
so $m_{0} \xbe X$ by prerequisite.

(15)
Let $X:=M(\xbf).$
$(\xbf?p)?q$ $=$ $(X \xfb (L-\{p\})) \xfb (L-\{q\})$ $=$ (by (14)) $X
\xfb (L-\{p,q\})$ $=$ $ \xbf \xfb p \xcu q.$

(16) $M(\xbf?TRUE)$ $=$ $M(\xbf) \xfb L$ $=$ $M(\xbf).$

(17)
$ \xbg_{X \xfb J}(m)=1$ iff $m \xbe X \xfb J$ iff $ \xcA m' \xbe M(m \xex
J=m' \xex J$ $ \xch $ $m' \xbe X)$ iff
$inf\{ \xbg_{X}(m'):$ $m' \xbe M,$ $m \xex J=m' \xex J\}$ $=$ 1.

$ \xcz $
\\[3ex]

\bfa

$\hspace{0.01em}$


\label{Fact Many-Down}

Consider the many-valued case.

Let $0:M \xcp V$ be the constant 0 $(=$ minimal value) function,
$1:M \xcp V$ be the constant 1 $(=$ maximal value) function.

The same counterexamples as for the 2-valued case work, we just
take the characteristic functions.

(1) $f \xfb J \xck f.$

(2) $f \xfb L=f.$

(4) $0 \xfb J=0.$

(5) $1 \xfb J=1.$

(6) $f \xck g$ $ \xch $ $f \xfb J \xck g \xfb J.$

(7) $J \xcc J' $ $ \xch $ $f \xfb J \xck f \xfb J'.$

(8) $sup(f \xfb J,g \xfb J)$ $ \xck $ $sup(f,g) \xfb J,$
the converse is not always true.

(9) $inf(f,g) \xfb J$ $=$ $inf(f \xfb J,g \xfb J).$

(11) $sup(f \xfb J,f \xfb J')$ $ \xck $ $f \xfb (J \xcv J'),$
the converse is not always true.

(12) $f \xfb (J \xcs J)$ $ \xck $ $inf(f \xfb J,f \xfb J'),$
the converse is not always true.

(14) $(f \xfb J) \xfb J' $ $=$ $f \xfb (J \xcs J')$

(15) $(\xbf?p)?q= \xbf?p \xcu q$

(16) $f_{ \xbf?TRUE}$ $=$ $f_{ \xbf }.$

\efa

\subparagraph{
Proof
}

$\hspace{0.01em}$


$(1)-(6)$ are trivial.

(7)
$m \xbe X \xfb J$ $ \xch $ $ \xcA m' (m \xex J=m' \xex J$ $ \xch $ $m'
\xbe X).$ But $m \xex J' =m' \xex J' $ $ \xch $ $m \xex J=m' \xex J$ $
\xch $ $m' \xbe X,$
so $m \xbe X \xfb J'.$

(8) by (6).

(9)

``$ \xck $'' by (6).

``$ \xcg $'':
Suppose $(inf(f,g) \xfb J)(m)<inf(f \xfb J(m),g \xfb J(m)).$ Choose $m' $
such that
$m \xex J=m' \xex J$ and $inf(f,g)(m')=inf(f,g) \xfb J(m).$ Note that $f
\xfb J(m) \xck f(m'),$
$g \xfb J(m) \xck g(m').$
Suppose without loss of generality $f(m') \xck g(m').$ Then
$f(m')=inf(f,g)(m')=(inf(f,g) \xfb J)(m)<inf(f \xfb J(m),g \xfb J(m))
\xck f \xfb J(m),$ a
contradiction.

(11) by (7).

(12) by (7).

(14)

``$ \xck $'': We have to show $(f \xfb J) \xfb J' (m)$ $ \xck $ $f \xfb
(J \xcs J')(m).$
Let $m'' $ be such that $m \xex J \xcs J' =m'' \xex J \xcs J',$ we have
to show
$(f \xfb J) \xfb J' (m)$ $ \xck $ $f(m''),$ then $(f \xfb J) \xfb J' (m)$
$ \xck $
$f \xfb J \xcs J' (m)$ $:=$ $min\{f(m''):m \xex J \xcs J' =m'' \xex J
\xcs J' \}.$
Take such $m''.$
Define $m' $ such that $m' $ on $J' $ is like $m,$ $m' $ on $L-J' $ is
like $m''.$
As $m' \xex J' =m \xex J',$ $(f \xfb J) \xfb J' (m) \xck f \xfb J(m').$
Note that $m'' $ is like $m' $ on $J:$ On $L-J',$ $m' $ is like $m''.$
On $J \xcs J',$ $m' $
is like $m,$ and $m$ is like $m'',$ so $m'' $ is like $m'.$
Thus $f \xfb J(m') \xck f(m''),$ so $(f \xfb J) \xfb J' (m)$ $ \xck $ $f
\xfb J(m')$ $ \xck $ $f(m'').$

``$ \xcg $'': We show $f \xfb J \xcs J' (m)$ $ \xck $ $(f \xfb J) \xfb J'
(m).$ We have to show that
$min\{f \xfb J(m'):m' \xex J' =m \xex J' \}$ $ \xcg $ $f \xfb J \xcs J'
(m).$
So take any $m' $ such that $m' \xex J' =m \xex J'.$ We have to show that
$f \xfb J(m') \xcg f \xfb J \xcs J' (m),$ i.e., for all $m'' $ such that
$m'' \xex J=m' \xex J,$ $f(m'')$ $ \xcg $ $f \xfb J \xcs J' (m)$ $:=$
$min\{f(m_{0}):m_{0} \xex J \xcs J' =m \xex J \xcs J' \}.$
If $m'' \xex J=m' \xex J,$ then $m'' \xex J \xcs J' =m' \xex J \xcs J' =m
\xex J \xcs J',$ so
$m'' $ is one of those $m_{0},$ and we are done.

(15)
Let $f:=f_{ \xbf }.$
(f?p)?q $=$ $(f \xfb (L-\{p\})) \xfb (L-\{q\})$ $=$ (by (14)) $f \xfb
(L-\{p,q\})$ $=$ $f \xfb p \xcu q.$

(16) $f_{ \xbf?TRUE}$ $=$ $f_{ \xbf } \xfb L$ $=$ $f_{ \xbf }.$

$ \xcz $
\\[3ex]
\section{
Monotone and antitone syntactic interpolation
}

\label{Section Mon-Synt-Int}
\subsection{
Introduction
}

We have shown $ \xCf semantic$ interpolation, this is not yet
$ \xCf syntactic$ interpolation. We still need that the set of
sequences is definable.
(The importance of definability in the context of non-monotonic
logics was examined by one of the authors
in  \cite{Sch92}.) Note that we
``simplified'' the set of sequences when going from $ \xbS' $ to $ \xbS''
,$
in Proposition \ref{Proposition Sin-Interpolation} (page \pageref{Proposition
Sin-Interpolation}), but perhaps the logic
at hand does not share this perspective. Consider, e.g.,
intuitionistic logic with three worlds, linearly ordered.
This is a monotonic logic, so by our results, it has
semantic interpolation. But it has no syntactic
interpolation, so the created set of models must not be
definable, see
Example \ref{Example Sin-4-Value} (page \pageref{Example Sin-4-Value}).
In classical propositional logic, the created set
$ \xCf is$ definable, as we will see in
Proposition \ref{Proposition Sin-Simplification-Definable} (page
\pageref{Proposition Sin-Simplification-Definable}).

But first a side remark.

\bcom

$\hspace{0.01em}$


\label{Comment:}

Usual approaches to repair interpolation construct a chain of ever richer
languages for interpolation:

We can go further with a logic in language L0 for which there is no
interpolation. For every pair of formulas which give a counterexample to
interpolation we introduce a new connective which corresponds to the
semantic
interpolant. Now we have a language L1 which allows for interpolation
for formulas in the original language L0. L1 itself may or may not
have interpolation. So we might have to continue to L2, L3, etc.

\ecom

This leads to the following definition for the cases where we have
semantic, but not necessarily syntactic interpolation.
We have, however, not examined this notion.

\bd

$\hspace{0.01em}$


\label{Definition Sin-Level}

A logic has level 0 semantic interpolation, iff it has interpolation.

A logic has level $n+1$ semantic interpolation iff it has no level $n$
semantic interpolation, but introducing new elements into the language
(of level $n)$ results in interpolation also for the new
language.

\ed

The case of full non-monotonic logic is, of course, different, as
the logics might not even have semantic interpolation,
so above repairing might not be possible.
\subsection{
The classical propositional case
}

We saw already in
Section \ref{Section Analoga} (page \pageref{Section Analoga})  the general
result for classical logic, we
give here a more special argument.
We show that semantic interpolation carries over to syntactic
interpolation in classical propositional logic. What we see as
simplification of sets of sequences, is seen the same way by classical
logic. If the first is definable, so is the second.

Nicht schon irgendwo???

\bp

$\hspace{0.01em}$


\label{Proposition Sin-Simplification-Definable}

Simplification (i.e. sup) preserves definability in classical
propositional
logic:

Let $ \xbG = \xbS \xex X' \xCK \xbP X''.$
Then, if $ \xbS $ is formula definable, so is $ \xbG.$

\ep

\subparagraph{
Proof
}

$\hspace{0.01em}$


As $ \xbS $ is formula definable, it is defined by $ \xbf_{1} \xco  \Xl
\xco \xbf_{n},$ where
$ \xbf_{i}= \xbq_{i,1} \xcu  \Xl  \xcu \xbq_{i,n_{i}}.$ Let $ \xbF_{i}:=\{
\xbq_{i,1}, \Xl, \xbq_{i,n_{i}}\},$
$ \xbF'_{i}:=\{ \xbq \xbe \xbF_{i}: \xbq \xbe X' \}$ (more precisely, $
\xbq $ or $ \xCN \xbq \xbe X'),$ $ \xbF''_{i}:= \xbF_{i}- \xbF'_{i}.$
Let $ \xbf'_{i}:= \xcU \xbF'_{i}.$ Thus $ \xbf_{i} \xcl \xbf'_{i}.$ We
show that $ \xbf'_{1} \xco  \Xl  \xco \xbf'_{n}$ defines $ \xbG.$
(Alternatively, we may replace all $ \xbq \xbe \xbF''_{i}$ by TRUE.)

(1) $ \xbG \xcm \xbf'_{1} \xco  \Xl  \xco \xbf'_{n}$

Lt $ \xbs \xbe \xbG,$ then there is $ \xbt \xbe \xbS $ s.t. $ \xbs \xex
X' = \xbf \xex X'.$ By prerequisite,
$ \xbt \xcm \xbf_{1} \xco  \Xl  \xco \xbf_{n},$ so $ \xbt \xcm \xbf'_{1}
\xco  \Xl  \xco \xbf'_{n},$ so $ \xbs \xcm \xbf'_{1} \xco  \Xl  \xco
\xbf'_{n}.$

(2) Suppose $ \xbs \xce \xbG,$ we have to show $ \xbs \xcM \xbf'_{1}
\xco  \Xl  \xco \xbf'_{n}.$

Suppose then $ \xbs \xce \xbG,$ but $ \xbs \xcm \xbf'_{1} \xco  \Xl
\xco \xbf'_{n},$ without loss of
generality
$ \xbs \xcm \xbf'_{1}= \xcU \xbF'_{1}.$
As $ \xbs \xce \xbG,$ there is no $ \xbt \xbe \xbS $ $ \xbt \xex X' =
\xbs \xex X'.$
Choose $ \xbt $ s.t. $ \xbs \xex X' = \xbt \xex X' $ and $ \xbt \xcm \xbq
$ for all $ \xbq \xbe \xbF''_{1}.$
By $ \xbs \xcm \xbq $ for $ \xbq \xbe \xbF'_{1},$ and $ \xbs \xex X' =
\xbt \xex X' $ $ \xbt \xcm \xbq $ for $ \xbq \xbe \xbF'_{1}.$
By prerequisite, $ \xbt \xcm \xbq $ for $ \xbq \xbe \xbF''_{1},$ so $
\xbt \xcm \xbq $ for all $ \xbq \xbe \xbF_{1},$ so
$ \xbt \xcm \xbf_{1} \xco  \Xl  \xco \xbf_{n},$ and $ \xbt \xbe \xbS,$ as
$ \xbf \xco  \Xl  \xco \xbf $ defines $ \xbS,$ contradiction.

$ \xcz $
\\[3ex]

\bco

$\hspace{0.01em}$


\label{Corollary Sin-Simplification-Definable}

The same result holds if $ \xbS $ is theory definable.

\eco

\subparagraph{
Proof
}

$\hspace{0.01em}$


(Outline). Define $ \xbS $ by a (possibly infinite) conjunction of
(finite)
disjunctions. Transform this into a possibly infinite disjunction of
possibly infinite conjunctions. Replace all $ \xbf \xbe \xbF''_{i}$ by
TRUE. The same
proof as above shows that this defines $ \xbG $ (finiteness was nowhere
needed). Transform backward into a conjunction of finite disjunctions,
where
the $ \xbf \xbe \xbF''_{i}$ are replaced by TRUE.

$ \xcz $
\\[3ex]

We could have used the following trivial fact for the proof of
Proposition \ref{Proposition Sin-Simplification-Definable} (page
\pageref{Proposition Sin-Simplification-Definable}):

\bfa

$\hspace{0.01em}$


\label{Fact Sin-Union}

Let $ \xbS = \xbS_{1} \xcv \xbS_{2}.$ Then $ \xbS \xex X' \xCK \xbP X'' =(
\xbS_{1} \xex X' \xCK \xbP X'') \xcv (\xbS_{2} \xex X' \xCK \xbP X'').$

\efa

\subparagraph{
Proof
}

$\hspace{0.01em}$


``$ \xcc $'': Let $ \xbs \xbe \xbS \xex X' \xCK \xbP X'',$ then there is
$ \xbs' \xbe \xbS $ s.t. $ \xbs \xex X' = \xbs' \xex X',$
so $ \xbs' \xbe \xbS_{1}$ or $ \xbs' \xbe \xbS_{2}.$ If $ \xbs' \xbe
\xbS_{1},$ then $ \xbs \xbe \xbS_{1} \xex X' \xCK \xbP X'',$ likewise
for $ \xbs' \xbe \xbS_{2}.$

The converse is even more trivial.

$ \xcz $
\\[3ex]

\br

$\hspace{0.01em}$


\label{Remark Sin-Intersection}

An analogous result about interesection does not hold, of course:
$(\xbS_{1} \xex X' \xCK X'') \xcs (\xbS_{2} \xex X' \xCK X'')$ might
well be $ \xEd \xCQ,$ but $ \xbS_{1} \xcs \xbS_{2}= \xCQ.$
\subsection{
General finite (intuionistic) Goedel logics
}

\label{Section Finite-Goedel}

\er

The semantics is a linearly ordered finite set of worlds,
with increasing truth, as usual in intuitionistic logics.
Let $n$ be the number of worlds, then $ \xbf $ has truth value
$k$ in the structure iff $ \xbf $ holds from world $n+1-k$ onward
(and 0 iff it never holds). Thus, if it holds everywhere, it
has truth value $n,$ if it holds from world 2 onward, it has value
$n-1,$ etc.
It is well known that such logics have (syntactic) interpolation
for $n=2,$ but not for $n>2.$ This is the reason we treat them here.
We will connect the interpolation problem to the existence of
normal forms. The connection is incomplete, as we will show that
suitable normal forms entail interpolation, but we do not know if
this condition is necessary.

\bd

$\hspace{0.01em}$


\label{Definition Mod-Fin-Goed}

Finite intuitionistic Goedel logics with $n+1$ truth values
$FALSE=0 \xck 1 \xck  \Xl  \xck n=TRUE$ are defined as follows:

 \xEh
 \xDH
$f_{ \xbf \xcu \xbq }(m):=inf\{f_{ \xbf }(m),f_{ \xbq }(m)\},$
 \xDH
$f_{ \xbf \xco \xbq }(m):=sup\{f_{ \xbf }(m),f_{ \xbq }(m)\},$
 \xDH negation $ \xCN $ is defined by:
\begin{flushleft}
\[ f_{\xCN \xbf}(m):= \left\{ \begin{array}{lcl}
TRUE \xEH iff \xEH
f_{\xbf}(m)=FALSE \xEP
\xEH \xEH \xEP
FALSE \xEH \xEH otherwise \xEP
\end{array}
\right.
\]
\end{flushleft}

 \xDH implication $ \xcp $ is defined by:
\begin{flushleft}
\[ f_{\xbf \xcp \xbq}(m):= \left\{ \begin{array}{lcl}
TRUE \xEH iff \xEH
f_{\xbf}(m) \xck f_{\xbq}(m) \xEP
\xEH \xEH \xEP
f_{\xbq}(m) \xEH \xEH otherwise \xEP
\end{array}
\right.
\]
\end{flushleft}

 \xEj

\ed

Thus, for $n+1=2,$ this is classical logic. So we assume now $n \xcg 2.$

\bd

$\hspace{0.01em}$


\label{Definition Mod-Fin-Goed-Add}

We will also consider the following additional operators:

 \xEh
 \xDH $J$ is defined by:
\begin{flushleft}
\[ f_{J \xbf}(m):= \left\{ \begin{array}{lclcl}
f_{\xbf}(m) \xEH iff \xEH
f_{\xbf}(m)=FALSE \xEH or \xEH f_{\xbf}(m)=TRUE \xEP
\xEH \xEH \xEH \xEH \xEP
f_{\xbf}(m)+1 \xEH \xEH otherwise \xEH \xEH \xEP
\end{array}
\right.
\]
\end{flushleft}

The intuitive meaning is:
``it holds in the next moment''

 \xDH $ \xCf A$ is defined by:
\begin{flushleft}
\[ f_{A(\xbf)}(m):= \left\{ \begin{array}{lcl}
TRUE \xEH iff \xEH
f_{\xbf}(m)=TRUE \xEP
\xEH \xEH \xEP
FALSE \xEH \xEH otherwise \xEP
\end{array}
\right.
\]
\end{flushleft}

Thus, $ \xCf A$ is the dual of negation, we might call it affirmation.

 \xDH $F$ is defined by:
\begin{flushleft}
\[ f_{F \xbf}(m):= \left\{ \begin{array}{lclcl}
FALSE \xEH iff \xEH
f_{\xbf}(m)=FALSE \xEH or \xEH f_{\xbf}(m)=TRUE \xEP
\xEH \xEH \xEH \xEH \xEP
f_{\xbf}(m)+1 \xEH \xEH otherwise \xEH \xEH \xEP
\end{array}
\right.
\]
\end{flushleft}

The intuitive meaning is:
``it begins to hold in the next moment''

 \xDH $Z$ (cyclic addition of 1) is defined by:
\begin{flushleft}
\[ f_{Z \xbf}(m):= \left\{ \begin{array}{lcl}
FALSE \xEH iff \xEH f_{\xbf}(m)=TRUE \xEP
\xEH \xEH \xEP
f_{\xbf}(m)+1 \xEH \xEH otherwise \xEP
\end{array}
\right.
\]
\end{flushleft}

Note that $Z$ is slightly different from $J.$
We do not know if there is an intuitive meaning.
 \xEj

To help the intuition, we give the truth tables of the basic operators
for $n=3,$ and of $ \xCN, \xcp,J,A,F,Z$ for $n=4$ and $n=6.$

\vspace{5mm}
\begin{tabular}{|c||c||c||c|c|c|c||c|c|c|c||c|c|c|c||c|c|c|c|}

\hline

\xEH \xEH $b$ \xEH \xEH 0 \xEH 1 \xEH 2 \xEH
\xEH 0 \xEH 1 \xEH 2 \xEH
\xEH 0 \xEH 1 \xEH 2 \xEH
\xEH 0 \xEH 1 \xEH 2 \xEP

\hline
\hline

$a$ \xEH $\xCN a$ \xEH \xEH $a \xcp b$ \xEH \xEH \xEH \xEH
$a \xcu b$ \xEH \xEH \xEH \xEH
$a \xco b$ \xEH \xEH \xEH \xEH
$a \xcr b$ \xEH \xEH \xEH \xEP

\hline

0 \xEH 2 \xEH \xEH \xEH 2 \xEH 2 \xEH 2 \xEH
\xEH 0 \xEH 0 \xEH 0 \xEH
\xEH 0 \xEH 1 \xEH 2 \xEH
\xEH 2 \xEH 0 \xEH 0 \xEP

\hline

1 \xEH 0 \xEH \xEH \xEH 0 \xEH 2 \xEH 2 \xEH
\xEH 0 \xEH 1 \xEH 1 \xEH
\xEH 1 \xEH 1 \xEH 2 \xEH
\xEH 0 \xEH 2 \xEH 1 \xEP

\hline

2 \xEH 0 \xEH \xEH \xEH 0 \xEH 1 \xEH 2 \xEH
\xEH 0 \xEH 1 \xEH 2 \xEH
\xEH 2 \xEH 2 \xEH 2 \xEH
\xEH 0 \xEH 1 \xEH 2 \xEP

\hline

\end{tabular}

\vspace{5mm}
\begin{tabular}{|c||c|c|c|c|c||c||c|c|c|c|c|}

\hline

\xEH \xEH \xEH \xEH \xEH \xEH $b$ \xEH \xEH 0 \xEH 1 \xEH 2 \xEH 3 \xEP

\hline
\hline

$a$ \xEH $\xCN a$ \xEH $Ja$ \xEH $Aa$ \xEH $Fa$ \xEH $Za$ \xEH \xEH $a \xcp b$
\xEH \xEH \xEH \xEH \xEP

\hline

0 \xEH 3 \xEH 0 \xEH 0 \xEH 0 \xEH 1 \xEH \xEH \xEH 3 \xEH 3 \xEH 3 \xEH 3 \xEP

\hline

1 \xEH 0 \xEH 2 \xEH 0 \xEH 2 \xEH 2 \xEH \xEH \xEH 0 \xEH 3 \xEH 3 \xEH 3 \xEP

\hline

2 \xEH 0 \xEH 3 \xEH 0 \xEH 3 \xEH 3 \xEH \xEH \xEH 0 \xEH 1 \xEH 3 \xEH 3 \xEP

\hline

3 \xEH 0 \xEH 3 \xEH 3 \xEH 0 \xEH 0 \xEH \xEH \xEH 0 \xEH 1 \xEH 2 \xEH 3 \xEP

\hline

\end{tabular}

\vspace{5mm}
\begin{tabular}{|c||c|c|c|c|c||c||c|c|c|c|c|c|c|}

\hline

\xEH \xEH \xEH \xEH \xEH \xEH $b$ \xEH \xEH 0 \xEH 1 \xEH 2 \xEH 3 \xEH 4 \xEH 5
\xEP

\hline
\hline

$a$ \xEH $\xCN a$ \xEH $Ja$ \xEH $Aa$ \xEH $Fa$ \xEH $Za$ \xEH \xEH $a \xcp b$
\xEH \xEH \xEH \xEH \xEH \xEH \xEP

\hline

0 \xEH 5 \xEH 0 \xEH 0 \xEH 0 \xEH 1 \xEH \xEH \xEH 5 \xEH 5 \xEH 5 \xEH 5 \xEH
5 \xEH 5 \xEP

\hline

1 \xEH 0 \xEH 2 \xEH 0 \xEH 2 \xEH 2 \xEH \xEH \xEH 0 \xEH 5 \xEH 5 \xEH 5 \xEH
5 \xEH 5 \xEP

\hline

2 \xEH 0 \xEH 3 \xEH 0 \xEH 3 \xEH 3 \xEH \xEH \xEH 0 \xEH 1 \xEH 5 \xEH 5 \xEH
5 \xEH 5 \xEP

\hline

3 \xEH 0 \xEH 4 \xEH 0 \xEH 4 \xEH 4 \xEH \xEH \xEH 0 \xEH 1 \xEH 2 \xEH 5 \xEH
5 \xEH 5 \xEP

\hline

4 \xEH 0 \xEH 5 \xEH 0 \xEH 5 \xEH 5 \xEH \xEH \xEH 0 \xEH 1 \xEH 2 \xEH 3 \xEH
5 \xEH 5 \xEP

\hline

5 \xEH 0 \xEH 5 \xEH 5 \xEH 0 \xEH 0 \xEH \xEH \xEH 0 \xEH 1 \xEH 2 \xEH 3 \xEH
4 \xEH 5 \xEP

\hline

\end{tabular}
\vspace{5mm}

\subsubsection{
The basic operators $\xcu, \xco, \xcp, \xCN$
}

\ed

We work now towards a suitable normal form, even though we cannot
obtain it for $n>3.$ This will also indicate a way to repair those
logics by introducing suitable additional operators, which allow
to obtain such normal forms.

We have the following fact:

\bfa

$\hspace{0.01em}$


\label{Fact Mod-Fin-Goed}

(0) With one variable $ \xCf a$ we can define up to semantical equivalence
exactly
the following 6 different formulas:

$ \xCf a,$ $ \xCN a,$ $ \xCN \xCN a,$ $TRUE=a \xcp a,$ $FALSE= \xCN (a
\xcp a),$ $ \xCN \xCN a \xcp a.$

The following semantic equivalences hold:

(Note: all except (14) hold also for 4 and 6 truth values, so
probably for arbitrarily many truth values, but this is not checked so
far.)

Triple negation can be simplified:

(1) $ \xCN \xCN \xCN a$ $ \xcr $ $ \xCN a$

Disjunction and conjunction combine classically:

(2) $ \xCN (a \xco b)$ $ \xcr $ $ \xCN a \xcu \xCN b$

(3) $ \xCN (a \xcu b)$ $ \xcr $ $ \xCN a \xco \xCN b$

(4) $a \xcu (b \xco c)$ $ \xcr $ $(a \xcu b) \xco (a \xcu c)$

(5) $a \xco (b \xcu c)$ $ \xcr $ $(a \xco b) \xcu (a \xco c)$

Implication can be eliminated from combined negation and implication:

(6) $ \xCN (a \xcp b)$ $ \xcr $ $ \xCN \xCN a \xcu \xCN b$

(7) $(a \xcp \xCN b)$ $ \xcr $ $(\xCN a \xco \xCN b)$

(8) $(\xCN a \xcp b)$ $ \xcr $ $(\xCN \xCN a \xco b)$

Implication can be put inside when combined with $ \xcu $ and $ \xco:$

(9) $(a \xco b \xcp c)$ $ \xcr $ $((a \xcp c) \xcu (b \xcp c))$

(10) $(a \xcu b \xcp c)$ $ \xcr $ $((a \xcp c) \xco (b \xcp c))$

(11) $(a \xcp b \xcu c)$ $ \xcr $ $((a \xcp b) \xcu (a \xcp c))$

(12) $(a \xcp b \xco c)$ $ \xcr $ $((a \xcp b) \xco (a \xcp c))$

Nested implication can be flattened for nesting on the right:

(13) $(a \xcp (b \xcp c))$ $ \xcr $ $((a \xcu b \xcp c) \xcu (a \xcu \xCN
c \xcp \xCN b))$

\efa

\subparagraph{
Proof
}

$\hspace{0.01em}$


We use $T$ for TRUE, $F$ for FALSE.

(0)
The truth table for the 6 formulas is given by the following table:

\begin{center}
\begin{tabular}{|c||c|c|c|c|c|c|c|}
\hline

$ \xCf a$ \xEH $ \xCf a$ \xEH $ \xCN a$ \xEH $ \xCN \xCN a$ \xEH $T$ $(=a
\xcp a)$ \xEH $F$ $(= \xCN (a \xcp a))$ \xEH $ \xCN \xCN a \xcp a$ \xEP
\hline

0 \xEH 0 \xEH $n$ \xEH 0 \xEH $n$ \xEH 0 \xEH $n$ \xEP
\hline

1 \xEH 1 \xEH 0 \xEH $n$ \xEH $n$ \xEH 0 \xEH 1 \xEP
\hline

2 \xEH 2 \xEH 0 \xEH $n$ \xEH $n$ \xEH 0 \xEH 2 \xEP
\hline

 \Xl  \xEH  \Xl  \xEH  \Xl  \xEH  \Xl  \xEH  \Xl  \xEH  \Xl  \xEH.. \xEP
\hline

$n$ \xEH $n$ \xEH 0 \xEH $n$ \xEH $n$ \xEH 0 \xEH $n$ \xEP
\hline
\end{tabular}
\end{center}

We see that the first line takes the values 0 and $n,$ and the $n-1$ other
lines take the vectors of values $(1, \Xl n),$ $(0, \Xl,0),$ $(n, \Xl
,n),$ and that all
combinations of first line values and those vectors occur. Thus, we can
check closure separately for the first line and the other lines, which is
now
trivial.

(1) Trivial.

$(2)+(3)$ Both sides can only be $T$ or $F.$
(2): Suppose $a \xck b,$ then $ \xCN (a \xco b)=T$ iff $b=F,$ and $ \xCN a
\xcu \xCN b=T$ iff $b=F.$
The case $b \xck a$ is symmetrical.
(3): similar: $ \xCN (a \xcu b)=T$ iff $a=F$ and $ \xCN a \xco \xCN b=T$
iff $a=F.$

(4) Suppose $b \xck c.$ Thus $(a \xcu b) \xco (a \xcu c)=a \xcu c.$ If $a
\xck c,$ then $a \xcu (b \xco c)=a,$ else
$a \xcu (b \xco c)=c.$ The case $c \xck b$ is symmetrical.

(5) Suppose $b \xck c.$ Then $(a \xco b) \xcu (a \xco c)=a \xco b.$ If $a
\xck b,$ then $a \xco (b \xcu c)=b,$ else
$a \xco (b \xcu c)=a.$

(6) Both sides are $T$ or $F.$ $ \xCN (a \xcp b)=T$ iff $a>b$ and $b=F.$
$ \xCN \xCN a \xcu \xCN b=T$ iff $b=F$ and $ \xCN \xCN a=F.$ $ \xCN \xCN
a=F$ iff $a>F.$

(7) Again, both sides are $T$ or $F.$ $a \xcp \xCN b=T$ iff $a \xck \xCN
b$ iff $b=F$ or $a=F.$

(8) $ \xCN a \xcp b=T$ iff $a>F$ or $b=T.$ If $a>F,$ then $ \xCN a \xcp
b=T$ and $ \xCN \xCN a \xco b=T.$
If $a=F,$ then $ \xCN a \xcp b=b,$ and $ \xCN \xCN a \xco b=b.$

(9) $a \xco b \xcp c$ is $a \xcp c$ or $b \xcp c.$
If $a \xck b,$ then it is $b \xcp c,$ and $a \xcp c \xcg b \xcp c.$ The
case $b \xck a$ is symmetrical.

$(10)-(12)$ are similar to (9), e.g., (11): If $b \xck c,$ then $(a \xcp
b) \xcu c=a \xcp b,$ and
$(a \xcp b) \xcu (a \xcp c)=a \xcp b.$

(13)
Case 1. $b \xck c:$ Then $ \xCN c \xck \xCN b,$ and $a \xcp (b \xcp c)=T,$
$a \xcu b \xcp c=T,$ $a \xcu \xCN c \xcp \xCN b=T.$

Case 2. $b>c:$ So $(a \xcp (b \xcp c))=(a \xcp c),$ and $ \xCN b=F.$

Case 2.1, $a \xck b:$ So $a \xcu b \xcp c=a \xcp c.$

Case 2.1.1. $ \xCN c=F:$ so $a \xcu \xCN c \xcp \xCN b=T,$ and we are
done.

Case 2.1.2. $ \xCN c=T:$ So $c=F,$ $a \xcp c=a \xcp F,$ and $a \xcu \xCN c
\xcp \xCN b=a \xcp \xCN b=a \xcp F,$ and
we are done again.

Case 2.2. $a>b:$ So $a>b>c,$ $a \xcp c=c,$ $a \xcu b \xcp c=b \xcp c=c,$
and $ \xCN b=F.$

Case 2.2.1. $ \xCN c=F:$ So $a \xcu \xCN c \xcp \xCN b=T,$ and we are
done.

Case 2.2.2. $ \xCN c=T:$ So $c=F.$ Thus $a \xcu \xCN c \xcp \xCN b=a \xcp
\xCN b=a \xcp F.$ But also $a \xcp c=a \xcp F,$
and we are done again.

$ \xcz $
\\[3ex]

We assume now that

(Assumption) Any formula of the type $(\xbf \xcp \xbf') \xcp \xbq $ is
equivalent to
a formula $ \xbF $ containing only flat $ \xcp $'s.

We will later show that this is true for $n=2.$

\bfa

$\hspace{0.01em}$


\label{Fact Eq-Normal}

Let above assumption be true. Then:

Every formula $ \xbf $ can be transformed into a semantically equivalent
formula $ \xbq $
of the following form:

(1) $ \xbq $ has the form $ \xbf_{1} \xco  \Xl  \xco \xbf_{n}$

(2) every $ \xbf_{i}$ has the form $ \xbf_{i,1} \xcu  \Xl  \xcu
\xbf_{i,m}$

(3) every $ \xbf_{i,m}$ has one of the following forms:

$p,$ or $ \xCN p,$ or $ \xCN \xCN p,$ or $p \xcp q$ - where $p$ and $q$
are propositional variables.

Note that also $ \xbf \xcp \xbf =TRUE$ can be replaced by $ \xCN a \xco
\xCN \xCN a.$

\efa

\subparagraph{
Proof
}

$\hspace{0.01em}$


The numbers refer to Fact \ref{Fact Mod-Fin-Goed} (page \pageref{Fact
Mod-Fin-Goed}).

We first push $ \xCN $ downward, towards the interior:

 \xEI
 \xDH $ \xCN (\xbf \xcu \xbq)$ is transformed to $ \xCN \xbf \xco \xCN
\xbq $ by (3).
 \xDH $ \xCN (\xbf \xco \xbq)$ is transformed to $ \xCN \xbf \xcu \xCN
\xbq $ by (2).
 \xDH $ \xCN (\xbf \xcp \xbq)$ is transformed to $ \xCN \xCN \xbf \xcu
\xCN \xbq $ by (6).
 \xEJ

We next eliminate any $ \xbf \xcp \xbq $ where $ \xbf $ and $ \xbq $ are
not propositional variables:

 \xEI
 \xDH $ \xCN \xbf \xcp \xbq $ is transformed to $ \xCN \xCN \xbf \xco \xbq
$ by (8).
 \xDH $ \xbf \xcu \xbf' \xcp \xbq $ is transformed to $(\xbf \xcp \xbq)
\xco (\xbf' \xcp \xbq)$ by (10).
 \xDH $ \xbf \xco \xbf' \xcp \xbq $ is transformed to $(\xbf \xcp \xbq)
\xcu (\xbf' \xcp \xbq)$ by (9).
 \xDH $(\xbf \xcp \xbf') \xcp \xbq $ is transformed to flat $ \xbF $ by
the assumption.
 \xEJ

 \xEI
 \xDH $ \xbf \xcp \xCN \xbq $ is transformed to $ \xCN \xbf \xco \xCN \xbq
$ by (7).
 \xDH $ \xbf \xcp \xbq \xcu \xbq' $ is transformed to $(\xbf \xcp \xbq)
\xcu (\xbf \xcp \xbq')$ by (11).
 \xDH $ \xbf \xcp \xbq \xco \xbq' $ is transformed to $(\xbf \xcp \xbq)
\xco (\xbf \xcp \xbq')$ by (12).
 \xDH $ \xbf \xcp (\xbq \xcp \xbq')$ is transformed to $(\xbf \xcu
\xbq \xcp \xbq') \xcu (\xbf \xcu \xCN \xbq' \xcp \xCN \xbq)$ by (13).
 \xEJ

Finally, we push $ \xcu $ inside:

$ \xbf \xcu (\xbq \xco \xbq')$ is transformed to $(\xbf \xcu \xbq)
\xco (\xbf \xcu \xbq')$ by (4).

The exact proof is, of course, by induction.

$ \xcz $
\\[3ex]

This normal form allows us to use the following facts:

\bfa

$\hspace{0.01em}$


\label{Fact Eq-Equiv-1}

We will now work for syntactic interpolation. For this purpose, we show
that, if $f$ is definable in
Proposition \ref{Proposition Gin-Pr-Int} (page \pageref{Proposition Gin-Pr-Int})
, i.e. there is $ \xbf $ with $f=f_{
\xbf },$ then
$f^{+}$ in the same Proposition is also definable. Recall that $f^{+}(m)$
was defined as
the maximal $f(m')$ for $m' \xex J' =m \xex J'.$ We use the normal form
just shown,
to show that conjuncts and disjuncts can be treated separately.

Our aim is to find a formula which characterizes the maximum. More
precisely,
if $f=f_{ \xbf }$ for some $ \xbf,$ we look for $ \xbf' $ such that
$f_{ \xbf' }(m)=max\{f_{ \xbf }(m'):$ $m' \xbe M,m \xex J=m' \xex J\}.$

First, a trivial fact, which shows that we can treat the elements of $J$
(or
$L-$J) one after
the other: $max\{g(x,y):x \xbe X,y \xbe Y\}=max\{max\{g(x,y):x \xbe X\}:y
\xbe Y\}.$ (Proof: The
interior max on the right hand side range over subsets of $X \xCK Y,$ so
they are all
$ \xck $ than the left hand side. Conversely, the left hand max is assumed
for some
$ \xBc x,y \xBe,$ which also figures on the right hand side. A full proof would
be
an
induction.)

Next, we show that we can treat disjunctions separately for one $x \xbe
L,$ and also
conjunctions, as long as $x$ occurs only in one of the conjuncts. Again, a
full
proof would be by induction, we only show the crucial arguments.
First, some notation:

\efa

\bn

$\hspace{0.01em}$


\label{Notation Eq-x}

(1) We write $m=_{(x)}m' $ as shorthand for $m \xex (L-\{x\})=m' \xex
(L-\{x\})$

(2) Let $f:M \xcp V,$ $x \xbe L,$ then $f_{(x)}(m):=max\{f(m'):$ $m' \xbe
M,m=_{(x)}m' \}.$

(3) Let $f_{ \xbf }:M \xcp V,$ and $(f_{ \xbf })_{(x)}=f_{ \xbf' }$ for
some $ \xbf',$ then we write $ \xbf_{(x)}$ for
(some such) $ \xbf'.$

\en

\bfa

$\hspace{0.01em}$


\label{Fact Eq-Lib}

(1) If $ \xbf = \xbf' \xco \xbf'',$ and $ \xbf'_{(x)},$ $ \xbf
''_{(x)}$ both exist, then so does $ \xbf_{(x)},$
and $ \xbf_{(x)}= \xbf'_{(x)} \xco \xbf''_{(x)}.$

(2) If $ \xbf = \xbf' \xcu \xbf'',$ $ \xbf'_{(x)}$ exists, and $ \xbf
'' $ does not contain $x,$ then
$ \xbf_{(x)}$ exists, and $ \xbf_{(x)}= \xbf'_{(x)} \xcu \xbf''.$

\efa

\subparagraph{
Proof
}

$\hspace{0.01em}$


(1) We have to show $f_{ \xbf_{(x)}}$ $=f_{(\xbf'_{(x)} \xco \xbf
''_{(x)})}.$

By definition of validity of $ \xco,$ we have
$f_{(\xbf'_{(x)} \xco \xbf''_{(x)})}(m)=max\{f_{ \xbf'_{(x)}}(m),f_{
\xbf''_{(x)}}(m)\}.$
$f_{ \xbf_{(x)}}(m):=max\{f_{ \xbf }(m'):m' =_{(x)}m\},$ so
$f_{(\xbf'_{(x)} \xco \xbf''_{(x)})}(m)$ $=$ $max\{max\{f_{ \xbf' }(m'
):m' =_{(x)}m\},max\{f_{ \xbf'' }(m'):m' =_{(x)}m\}\}$ $=$
$max\{max\{f_{ \xbf' }(m'),f_{ \xbf'' }(m')\}:m' =_{(x)}m\}$ $=$
(again by definition of validity of $ \xco)$
$max\{f_{ \xbf' \xco \xbf'' }(m'):m' =_{(x)}m\}$ $=$ $max\{f_{ \xbf
}(m'):m' =_{(x)}m\}$ $=$ $f_{(\xbf_{(x)})}(m).$

(2) We have to show $f_{ \xbf_{(x)}}=f_{(\xbf'_{(x)} \xcu \xbf
''_{(x)})}.$
By definition of validity of $ \xcu,$ we have
$f_{(\xbf'_{(x)} \xcu \xbf''_{(x)})}(m)=inf\{f_{ \xbf'_{(x)}}(m),f_{
\xbf''_{(x)}}(m)\}.$ So
$f_{(\xbf'_{(x)} \xcu \xbf''_{(x)})}(m)$ $=$ $inf\{max\{f_{ \xbf' }(m'
):m' =_{(x)}m\},max\{f_{ \xbf'' }(m'):m' =_{(x)}m\}\}$ $=$
(as $ \xbf'' $ does not contain $x)$
$inf\{max\{f_{ \xbf' }(m'):m' =_{(x)}m\},f_{ \xbf'' }(m)\}$ $=$
$max\{inf\{f_{ \xbf' }(m'),f_{ \xbf'' }(m)\}:m' =_{(x)}m\}$ $=$ (again
by definition of validity of $ \xcu,$
and by the fact that $ \xbf'' $ does not contain $x)$
$max\{f_{ \xbf' \xcu \xbf'' }(m'):m' =_{(x)}m\}$ $=$ $max\{f_{ \xbf
}(m'):m' =_{(x)}m\}$ $=$ $f_{(\xbf }$ $)^{(m).}$

$ \xcz $
\\[3ex]

Thus, we can calculate disjunctions separately, and also conjunctions,
as long as the latter have no variables in common. In classical logic,
we are finished, as we can break down conjunctions into parts which have
no variables in common.
The problem here are formulas of the type $a \xcp b,$ as they may have
variables in
common with other conjuncts, and, as we will see
in Fact \ref{Fact Eq-Def} (page \pageref{Fact Eq-Def})  (2) and (3), they cannot
be eliminated.

Thus, we have to consider situations like $(a \xcp b) \xcu (b \xcp c),$ $a
\xcu (a \xcp b),$ etc.,
where without loss of generality none is of the form $a \xcp a,$ as this
can
be replaced by TRUE.

To do as many cases together as possible, it is useful to use
Fact \ref{Fact Mod-Fin-Goed} (page \pageref{Fact Mod-Fin-Goed})  (9) and (11)
backwards, to obtain
general formulas.
We then see that the cases to examine are of the form:

$ \xbf $ $=$ $((b_{1} \xco  \Xl  \xco b_{n}) \xcp a) \xcu (a \xcp (c_{1}
\xcu  \Xl  \xcu c_{m})) \xcu \xbs a \xcu \xbt a \xcu \xbr a,$
where none of the $b_{i}$ or $c_{i}$ are $ \xCf a,$ and
where $n,m$ may be 0, and $ \xbs, \xbt, \xbr $ are absence $(\xCQ,$ no
$ \xCf a),$ $ \xCf a,$ $ \xCN a,$ or
$ \xCN \xCN a.$

We have the following equalities (for $F=FALSE,$ $T=TRUE):$

$a \xcu \xCN a=F,$ $a \xcu \xCN \xCN a=a,$ $ \xCN a \xcu \xCN \xCN a=F,$
$a \xcu \xCN a \xcu \xCN \xCN a=F.$

Thus, it suffices to consider $ \xbs $ as empty, $ \xCf a,$ $ \xCN a,$ $
\xCN \xCN a,$ which leaves us
with 4 cases. Moreover, we see that we always treat $b_{1} \xco  \Xl  \xco
b_{n}$ and
$c_{1} \xcu  \Xl  \xcu c_{m}$ as one block, so we can without loss of
generality restrict
the consideration to the 12 cases:

$ \xbf_{1,1}$ $:=$ $(b \xcp a)$

$ \xbf_{1,2}$ $:=$ $(b \xcp a) \xcu a$

$ \xbf_{1,3}$ $:=$ $(b \xcp a) \xcu \xCN a$

$ \xbf_{1,4}$ $:=$ $(b \xcp a) \xcu \xCN \xCN a$

$ \xbf_{2,1}$ $:=$ $(a \xcp c)$

$ \xbf_{2,2}$ $:=$ $(a \xcp c) \xcu a$

$ \xbf_{2,3}$ $:=$ $(a \xcp c) \xcu \xCN a$

$ \xbf_{2,4}$ $:=$ $(a \xcp c) \xcu \xCN \xCN a$

$ \xbf_{3,1}$ $:=$ $(b \xcp a) \xcu (a \xcp c)$

$ \xbf_{3,2}$ $:=$ $(b \xcp a) \xcu (a \xcp c) \xcu a$

$ \xbf_{3,3}$ $:=$ $(b \xcp a) \xcu (a \xcp c) \xcu \xCN a$

$ \xbf_{3,4}$ $:=$ $(b \xcp a) \xcu (a \xcp c) \xcu \xCN \xCN a$

We consider now the maximum, when we let $ \xCf a$ float, i.e., consider
all
$m' $ such that $m \xex L-\{a\}=m' \xex L-\{a\}.$ Let $ \xbf'_{i,j}$ be
this maximum.
For $ \xbf_{1,1},$ $ \xbf_{1,2},$ $ \xbf_{1,4},$ $ \xbf'_{i,j}=T$ (take
$a=T).$

For $ \xbf_{2,1},$ $ \xbf_{2,3},$ $ \xbf'_{i,j}=T$ (take $a=F).$

Next, we consider the remaining simple cases $ \xbf_{1,i}$ and $
\xbf_{2,i}.$
We show $ \xbf'_{1,3}= \xCN b,$ $ \xbf'_{2,2}=c,$ $ \xbf'_{2,4}= \xCN
\xCN c,$ see
Table \ref{Table Eq-Neglecting-1} (page \pageref{Table Eq-Neglecting-1}).
(We abbreviate e.g. $m(a)<m(b)$ by $a<b.)$

.
\begin{table}

\label{Table Eq-Neglecting-1}
\begin{center}
\caption{Table Neglecting a variable - Part 1}
\begin{tabular}{|c||c|c|c|c|}
\hline

$ \xbf_{1,3}$ \xEH \xEH $b \xcp a$ \xEH $ \xCN a$ \xEH $(b \xcp a) \xcu
\xCN a$ \xEP
\hline

 \xEH $a<b$ \xEH $ \xCf a$ \xEH \xEH $F$ \xEP
\hline

 \xEH $a=b$ \xEH $T$ \xEH $ \xCN b$ \xEH $ \xCN b$ \xEP
\hline

 \xEH $a>b$ \xEH $T$ \xEH $F$ \xEH $F$ \xEP
\hline
\hline

$ \xbf_{2,2}$ \xEH \xEH $a \xcp c$ \xEH $ \xCf a$ \xEH $(a \xcp c) \xcu a$
\xEP
\hline

 \xEH $a<c$ \xEH $T$ \xEH $ \xCf a$ \xEH $ \xCf a$ \xEP
\hline

 \xEH $a=c$ \xEH $T$ \xEH $c$ \xEH $c$ \xEP
\hline

 \xEH $a>c$ \xEH $c$ \xEH $ \xCf a$ \xEH $c$ \xEP
\hline
\hline

$ \xbf_{2,4}$ \xEH \xEH $a \xcp c$ \xEH $ \xCN \xCN a$ \xEH $(a \xcp c)
\xcu \xCN \xCN a$ \xEP
\hline

 \xEH $a<c$ \xEH $T$ \xEH $ \xCN \xCN a \xck \xCN \xCN c$ \xEH $ \xCN \xCN
a$ \xEP
\hline

 \xEH $a=c$ \xEH $T$ \xEH $ \xCN \xCN c$ \xEH $ \xCN \xCN c$ \xEP
\hline

 \xEH $a>c$ \xEH $c$ \xEH $ \xCN \xCN a$ $(\xcg \xCN \xCN c,c)$ \xEH $c$
$(\xck \xCN \xCN c)$ \xEP
\hline
\end{tabular}
\end{center}
\end{table}

We show now $ \xbf'_{3,1}=b \xcp c,$ $ \xbf'_{3,2}=c,$ $ \xbf'_{3,3}=(b
\xcp c) \xcu \xCN b,$ $ \xbf'_{3,4}=(b \xcp c) \xcu \xCN \xCN c,$ see
Table \ref{Table Eq-Neglecting-2} (page \pageref{Table Eq-Neglecting-2}).
\begin{table}

\label{Table Eq-Neglecting-2}
\begin{center}
\caption{Table Neglecting a variable - Part 2}
{\xssB
\begin{tabular}{|c||l|c|c|c|c||c|c|c|}
\hline

$ \xbf_{3,1}$ \xEH \xEH $b \xcp a$ \xEH $a \xcp c$ \xEH \xEH $(b \xcp a)
\xcu (a \xcp c)$ \xEH $b \xcp c$ \xEH \xEH \xEP
\hline

 \xEH Case 1: $b \xck c$ \xEH \xEH \xEH \xEH \xEH \xEH \xEH \xEP
\hline

 \xEH 1.1: $a<b$ \xEH $ \xCf a$ \xEH $T$ \xEH \xEH $ \xCf a$ \xEH $T$ \xEH
\xEH \xEP
\hline

 \xEH 1.2: $a=b$ \xEH $T$ \xEH $T$ \xEH \xEH $T$ \xEH \xEH \xEH \xEP
\hline

 \xEH 1.3: $b \xck a \xck c$ \xEH $T$ \xEH $T$ \xEH \xEH $T$ \xEH \xEH
\xEH \xEP
\hline

 \xEH 1.4: $c<a$ \xEH $T$ \xEH $c$ \xEH \xEH $c$ \xEH \xEH \xEH \xEP
\hline

 \xEH Case 2: $c<b$ \xEH \xEH \xEH \xEH \xEH \xEH \xEH \xEP
\hline

 \xEH 2.1: $a \xck c$ \xEH $ \xCf a$ \xEH $T$ \xEH \xEH $ \xCf a$ $(\xck
c)$ \xEH $c$ \xEH \xEH \xEP
\hline

 \xEH 2.2: $c<a<b$ \xEH $ \xCf a$ \xEH $c$ \xEH \xEH $c$ \xEH \xEH \xEH
\xEP
\hline

 \xEH 2.3: $b \xck a$ \xEH $T$ \xEH $c$ \xEH \xEH $c$ \xEH \xEH \xEH \xEP
\hline
\hline

$ \xbf_{3,2}$ \xEH \xEH $b \xcp a$ \xEH $a \xcp c$ \xEH $ \xCf a$ \xEH $(b
\xcp a) \xcu (a \xcp c) \xcu a$ \xEH $c$ \xEH \xEH \xEP
\hline

 \xEH Case 1: $b \xck c$ \xEH \xEH \xEH \xEH \xEH \xEH \xEH \xEP
\hline

 \xEH 1.1: $a<b$ \xEH $ \xCf a$ \xEH $T$ \xEH \xEH $ \xCf a$ $(<c)$ \xEH
\xEH \xEH \xEP
\hline

 \xEH 1.2: $a=b$ \xEH $T$ \xEH $T$ \xEH \xEH $ \xCf a$ $(\xck c)$ \xEH
\xEH \xEH \xEP
\hline

 \xEH 1.3: $b \xck a \xck c$ \xEH $T$ \xEH $T$ \xEH \xEH $ \xCf a$ $(\xck
c)$ \xEH \xEH \xEH \xEP
\hline

 \xEH 1.4: $c<a$ \xEH $T$ \xEH $c$ \xEH \xEH $c$ \xEH \xEH \xEH \xEP
\hline

 \xEH Case 2: $c<b$ \xEH \xEH \xEH \xEH \xEH \xEH \xEH \xEP
\hline

 \xEH 2.1: $a \xck c$ \xEH $ \xCf a$ \xEH $T$ \xEH \xEH $ \xCf a$ $(\xck
c)$ \xEH \xEH \xEH \xEP
\hline

 \xEH 2.2: $c<a<b$ \xEH $ \xCf a$ \xEH $c$ \xEH \xEH $c$ \xEH \xEH \xEH
\xEP
\hline

 \xEH 2.3: $b \xck a$ \xEH $T$ \xEH $c$ \xEH \xEH $c$ \xEH \xEH \xEH \xEP
\hline
\hline

$ \xbf_{3,3}$ \xEH \xEH $b \xcp a$ \xEH $a \xcp c$ \xEH $ \xCN a$ \xEH $(b
\xcp a) \xcu (a \xcp c) \xcu \xCN a$ \xEH $b \xcp c$ \xEH $ \xCN b$ \xEH
$(b \xcp c) \xcu \xCN b$ \xEP
\hline

 \xEH Case 1: $b \xck c$ \xEH \xEH \xEH \xEH \xEH \xEH \xEH \xEP
\hline

 \xEH 1.1: $a<b$ \xEH $ \xCf a$ \xEH $T$ \xEH $ \xCN a$ \xEH $F$ \xEH $T$
\xEH $ \xCN b$ \xEH $ \xCN b$ \xEP
\hline

 \xEH 1.2: $a=b$ \xEH $T$ \xEH $T$ \xEH $ \xCN a$ \xEH $ \xCN a= \xCN b$
\xEH \xEH \xEH \xEP
\hline

 \xEH 1.3: $b \xck a \xck c$ \xEH $T$ \xEH $T$ \xEH $ \xCN a$ $(\xck \xCN
b)$ \xEH $ \xCN a$ \xEH \xEH \xEH \xEP
\hline

 \xEH 1.4: $c<a$ \xEH $T$ \xEH $c$ \xEH $ \xCN a$ $(\xck \xCN b, \xCN c)$
\xEH $ \xCN a \xcu c=F$ \xEH \xEH \xEH \xEP
\hline

 \xEH Case 2: $c<b$ \xEH \xEH \xEH \xEH \xEH \xEH \xEH \xEP
\hline

 \xEH 2.1: $a \xck c$ \xEH $ \xCf a$ \xEH $T$ \xEH $ \xCN a$ \xEH $F$ \xEH
$c$ \xEH $ \xCN b=F$ \xEH $F$ \xEP
\hline

 \xEH 2.2: $c<a<b$ \xEH $ \xCf a$ \xEH $c$ \xEH $ \xCN a$ \xEH $F$ \xEH
\xEH \xEH \xEP
\hline

 \xEH 2.3: $b \xck a$ \xEH $T$ \xEH $c$ \xEH $ \xCN a$ $(\xck \xCN c)$
\xEH $F$ \xEH \xEH \xEH \xEP
\hline
\hline

$ \xbf_{3,4}$ \xEH \xEH $b \xcp a$ \xEH $a \xcp c$ \xEH $ \xCN \xCN a$
\xEH $(b \xcp a) \xcu (a \xcp c) \xcu \xCN \xCN a$ \xEH $b \xcp c$ \xEH $
\xCN \xCN c$ \xEH $(b \xcp c) \xcu
 \xCN \xCN c$ \xEP
\hline

 \xEH Case 1: $b \xck c$ \xEH \xEH \xEH \xEH \xEH \xEH \xEH \xEP
\hline

 \xEH 1.1: $a<b$ \xEH $ \xCf a$ \xEH $T$ \xEH $ \xCN \xCN a$ \xEH $a \xck
\xCN \xCN c$ \xEH $T$ \xEH $ \xCN \xCN c$ \xEH $ \xCN \xCN c$ \xEP
\hline

 \xEH 1.2: $a=b$ \xEH $T$ \xEH $T$ \xEH $ \xCN \xCN a$ \xEH $ \xCN \xCN a
\xck \xCN \xCN c$ \xEH \xEH \xEH \xEP
\hline

 \xEH 1.3: $b \xck a \xck c$ \xEH $T$ \xEH $T$ \xEH $ \xCN \xCN a$ \xEH $
\xCN \xCN a \xck \xCN \xCN c$ \xEH \xEH \xEH \xEP
 \xEH \xEH \xEH \xEH \xEH $= \xCN \xCN c$ if $a=c$ \xEH \xEH \xEH \xEP
\hline

 \xEH 1.4: $c<a$ \xEH $T$ \xEH $c$ \xEH $T$ \xEH $c \xck \xCN \xCN c$ \xEH
\xEH \xEH \xEP
\hline

 \xEH Case 2: $c<b$ \xEH \xEH \xEH \xEH \xEH \xEH \xEH \xEP
\hline

 \xEH 2.1: $a \xck c$ \xEH $ \xCf a$ \xEH $T$ \xEH $ \xCN \xCN a$ \xEH $a
\xck c$ \xEH $c$ \xEH $ \xCN \xCN c$ \xEH $c$ \xEP
\hline

 \xEH 2.2: $c<a<b$ \xEH $ \xCf a$ \xEH $c$ \xEH $T$ \xEH $c$ \xEH \xEH
\xEH \xEP
\hline

 \xEH 2.3: $b \xck a$ \xEH $T$ \xEH $c$ \xEH $T$ \xEH $c$ \xEH \xEH \xEH
\xEP
\hline
\end{tabular}
}
\end{center}
\end{table}

\br

$\hspace{0.01em}$


\label{Remark Eq-Trans}

We cannot improve the value of $ \xbf \xcp \xbq $ by taking a detour $
\xbf \xcp \xba_{1} \xcp  \Xl  \xcp \xba_{n} \xcp \xbq $
because the destination determines the value: in any column of $ \xcp,$
there
is only max and a constant value. And if we go further down than needed,
we get only worse, going from right to left deteriorates the values in the
lines.
$ \xcz $
\\[3ex]

\er

We can achieve the same result by first closing under the following
rules, and then erasing all formulas containing $ \xCf a:$

(1) $ \xcp $ under transitivity, i.e.

$((b_{1} \xco  \Xl  \xco b_{n}) \xcp a) \xcu (a \xcp (c_{1} \xcu  \Xl
\xcu c_{m}))$ $ \xch $ $((b_{1} \xco  \Xl  \xco b_{n}) \xcp (c_{1} \xcu
\Xl  \xcu c_{m}))$

(2) $ \xbs' a$ and $ \xcp $ as follows:

$((b_{1} \xco  \Xl  \xco b_{n}) \xcp a) \xcu (a \xcp (c_{1} \xcu  \Xl
\xcu c_{m})),$ $ \xCf a$ $ \xch $
$((b_{1} \xco  \Xl  \xco b_{n}) \xcp (c_{1} \xcu  \Xl  \xcu c_{m})) \xcu
c_{1} \xcu  \Xl  \xcu c_{m}$

$((b_{1} \xco  \Xl  \xco b_{n}) \xcp a) \xcu (a \xcp (c_{1} \xcu  \Xl
\xcu c_{m})),$ $ \xCN \xCN a$ $ \xch $
$((b_{1} \xco  \Xl  \xco b_{n}) \xcp (c_{1} \xcu  \Xl  \xcu c_{m})) \xcu
\xCN \xCN c_{1} \xcu  \Xl  \xcu \xCN \xCN c_{m}$

$((b_{1} \xco  \Xl  \xco b_{n}) \xcp a) \xcu (a \xcp (c_{1} \xcu  \Xl
\xcu c_{m})),$ $ \xCN a$ $ \xch $
$((b_{1} \xco  \Xl  \xco b_{n}) \xcp (c_{1} \xcu  \Xl  \xcu c_{m})) \xcu
\xCN b_{1} \xcu  \Xl  \xcu \xCN b_{n}$

In summary: the semantical interpolant constructed in
Section \ref{Section Eq-Sem} (page \pageref{Section Eq-Sem})  is definable, if
the assumption holds, so
the HT logic
(see Section \ref{Section HT} (page \pageref{Section HT})) has
also syntactic interpolation. This result is well-known,
but we need the techniques for the next section.

In some cases, introducing
new constants analogous to TRUE, FALSE - in the cited case e.g. ONE, TWO
when truth starts at world one or two - might help, but we did not
investigate this question.
This question is also examined in  \cite{ABM03}.
\subsubsection{
An important example for non-existence of interpolation
}

\label{Section Example-No-Int}

We turn now to an important example. It shows that the logic with
3 worlds, and thus 4 truth values, has no interpolation. But first, we
show as much as possible for the general case (arbitrarily many truth
values).

\be

$\hspace{0.01em}$


\label{Example Sin-4-Value}

Let

$ \xba (p,q,r):= \xCI p \xcp (((q \xcp r) \xcp q) \xcp q) \xCJ \xcp p,$

$ \xbb (p,s):=((s \xcp p) \xcp s) \xcp s$

\ee

We will show that
$ \xba (p,q,r) \xcp \xbb (p,s)$ holds in the case of 3 worlds, but that
there is no syntactic
interpolant (which could use only $p).$

Introducing a new operator $ \xCf Jp$ meaning
``from next moment onwards $p$ holds and if now is the last moment then $p$
holds
now''
gives enough definability to have also syntactic interpolation for $ \xba
$ and
$ \xbb $ above. This will be shown in
Section \ref{Section 4-Val} (page \pageref{Section 4-Val}). First, we give some
general results for
above example.

\bfa

$\hspace{0.01em}$


\label{Fact Mod-Fin-Goed-Impl}

Let $T,$ truth, be the maximal truth value.

 \xEh
 \xDH $ \xbf:=((\Xl (((a \xcp b) \xcp b) \xcp b) \Xl) \xcp b)$ has the
following truth value $v(\xbf)$ in a
model $m:$

 \xEh
 \xDH if the number $n$ of $b$ on the right of the first $ \xcp $ is odd:

if $v(a) \xck v(b),$ then $v(\xbf)=T,$ otherwise $v(\xbf)=v(b),$

 \xDH if the number $m$ of $b$ on the right of the first $ \xcp $ is even:

if $v(a) \xck v(b),$ then $v(\xbf)=v(b),$ otherwise $v(\xbf)=T.$
 \xEj

 \xDH $ \xbf:=((\Xl (((a \xcp b) \xcp a) \xcp a) \Xl) \xcp a)$ has the
following truth value $v(\xbf)$ in a
model $m:$

 \xEh
 \xDH if the number $n$ of $ \xCf a$ on the right of the first $ \xcp $ is
odd:

if $v(b)<v(a),$ then $v(\xbf)=T,$ otherwise $v(\xbf)=v(a),$

 \xDH if the number $m$ of $ \xCf a$ on the right of the first $ \xcp $ is
even:

if $v(b)<v(a),$ then $v(\xbf)=v(a),$ otherwise $v(\xbf)=T.$
 \xEj

 \xEj

\efa

\subparagraph{
Proof
}

$\hspace{0.01em}$


(1)

We proceed by induction.

(1.1) For $n=1,$ it is the definition of $ \xcp.$

(1.2) Case $n=2:$ If $v(a) \xck v(b),$ then $v(a \xcp b)=T,$ so $v((a \xcp
b) \xcp b)=v(b).$
If $v(a)>v(b),$ then $v(a \xcp b)=v(b),$ so $v((a \xcp b) \xcp b)=T.$

The general induction works as for the step from $n=1$ to $n=2.$

(2)

$n=1:$

$ \xbf =(a \xcp b) \xcp a.$ If $v(b)<v(a),$ then $v(a \xcp b)=v(b),$ so
$v(\xbf)=T.$
If $v(b) \xcg v(a),$ then $v(a \xcp b)=T,$ so $ \xbf (\xbf)=v(a).$

$n \xcp n+1:$

$ \xbf = \xbq \xcp a.$ If $v(b)<v(a),$ then, if $n$ is odd, $v(\xbq)=T,$
so $v(\xbf)=v(a),$ if $n$ is even,
$v(\xbq)=v(a),$ so $v(\xbf)=T.$
If $v(b) \xcg v(a),$ then, if $n$ is odd, $v(\xbq)=v(a),$ so $v(\xbf
)=T,$ if $n$ is even, $v(\xbq)=T,$
so $v(\xbf)=v(a).$

$ \xcz $
\\[3ex]

\bco

$\hspace{0.01em}$


\label{Corollary Mod-Int-Fail}

Let $T$ be the maximal truth value TRUE.

Consider again the formulas of
Example \ref{Example Sin-4-Value} (page \pageref{Example Sin-4-Value}),
$ \xba (p,q,r):= \xCI p \xcp (((q \xcp r) \xcp q) \xcp q) \xCJ \xcp p,$
$ \xbb (p,s):=((s \xcp p) \xcp s) \xcp s.$

We use Fact \ref{Fact Mod-Fin-Goed-Impl} (page \pageref{Fact Mod-Fin-Goed-Impl})
, the numbers refer to this
fact.

Let $f:=f_{ \xbb }.$
Let $f' (m):=min\{f(m'):m \xex p=m' \xex p\}.$ Fix $m.$
By (2.2), if $m(p)=T,$ then $f' (m)=T,$ if $m(p)<T,$ then $f' (m)=m(p)+1.$

Let $g:=f_{ \xba }.$
Let $g' (m):=max\{g(m'):m \xex p=m' \xex p\}.$ Fix $m.$
$ \xba $ is of the form $(p \xcp \xbf) \xcp p,$ so by (2.1),
if $m' (\xbf)<m' (p),$ then $m' (\xba)=T,$ if $m' (\xbf) \xcg m'
(p),$ then $m' (\xba)=m' (p).$
By (2.2), we have:
if $m' (r)<m' (q),$ then $m' (\xbf)=m' (q),$ if $m' (r) \xcg m' (q),$
then $m' (\xbf)=T.$
Note that $m' (\xbf)>0.$

Table \ref{Table Alpha} (page \pageref{Table Alpha})  shows that for $T=3$ $
\xba \xcl \xbb,$ for
$T=4$ $ \xba \xcL \xbb.$

Thus, an interpolant $h$ must have $h(0)=0$ or 1, $h(1)=1$ or 2,
$h(2)=h(3)=3$
in the case $T=3.$ This is impossible by
Fact \ref{Fact Mod-Fin-Goed} (page \pageref{Fact Mod-Fin-Goed}).

$ \xcz $
\\[3ex]

\eco

.
\begin{table}

\label{Table Alpha}
\begin{center}
\caption{Table $\xba \xcl \xbb$}
\begin{tabular}{|c|c|c||c|c|c|}
\hline
\multicolumn{3}{|c|}{T=3} \xEH
\multicolumn{3}{||c|}{T=4} \xEP
\hline

$m(p)$ \xEH $f' (m)$ \xEH $g' (m)$ \xEH $m(p)$ \xEH $f' (m)$ \xEH $g' (m)$
\xEP
\hline

0 \xEH 1 \xEH 0 \xEH 0 \xEH 1 \xEH 0 \xEP
\hline

1 \xEH 2 \xEH 1 \xEH 1 \xEH 2 \xEH 1 \xEP
\hline

2 \xEH 3 \xEH 3 \xEH 2 \xEH 3 \xEH 4 \xEP
\hline

3 \xEH 3 \xEH 3 \xEH 3 \xEH 4 \xEH 4 \xEP
\hline

 \xEH \xEH \xEH 4 \xEH 4 \xEH 4 \xEP
\hline
\end{tabular}
\end{center}
\end{table}
\subsubsection{
The additional operators $J$ and $A$
}

The following was checked with a small computer program:

(1) $ \xCf A$ alone will generate 12 semantically different formulas with
1 variable,
but it does not suffice to obtain interpolation.

(2) $ \xCf J$ alone will generate 8 semantically different formulas with 1
variable,
and it will solve the interpolation problem for $ \xba (p,q,r)$ and $ \xbb
(p,s)$
of Example \ref{Example Sin-4-Value} (page \pageref{Example Sin-4-Value})

(3) $ \xCf A$ and $ \xCf J$ will generate 48 semantically different
formulas with 1
variable.
\vspace{5mm}
\subsubsection{
The additional operator $F$
}

For $n+1$ truth values, let for $k<n$

$ \xbf'_{k}:= \xCN (F^{k}(a) \xcp F^{k+1}(a)).$

$ \xbf_{k}:=a \xcu \xbf'_{k}.$

(For $k=n,$ we take $ \xCN a.)$

Then
\begin{flushleft}
\[ f_{\xbf_k}(m):= \left\{ \begin{array}{lcl}
n-k \xEH iff \xEH m=n-k \xEP
\xEH \xEH \xEP
FALSE \xEH \xEH otherwise \xEP
\end{array}
\right.
\]
\end{flushleft}

Applying $F$ again, we can increase the value from $n-k$ up to TRUE.

We give the table of $ \xbf'_{k}(\xba)$ for 6 truth values.
\\ \\
\vspace{5mm}
\begin{tabular}{|c|c|c|c|c|c|c|}
\hline

$ \xba $ \xEH 0 \xEH 1 \xEH 2 \xEH 3 \xEH 4 \xEH 5 \xEP
\hline

$ \xbf'_{5}= \xCN \xba $ \xEH 5 \xEH 0 \xEH 0 \xEH 0 \xEH 0 \xEH 0 \xEP
\hline

$ \xbf'_{4}= \xCN (F^{4}(\xba) \xcp F^{5}(\xba))$ \xEH 0 \xEH 5 \xEH
0 \xEH 0 \xEH 0 \xEH 0 \xEP
\hline

$ \xbf'_{3}= \xCN (F^{3}(\xba) \xcp F^{4}(\xba))$ \xEH 0 \xEH 0 \xEH
5 \xEH 0 \xEH 0 \xEH 0 \xEP
\hline

$ \xbf'_{2}= \xCN (F^{2}(\xba) \xcp F^{3}(\xba))$ \xEH 0 \xEH 0 \xEH
0 \xEH 5 \xEH 0 \xEH 0 \xEP
\hline

$ \xbf'_{1}= \xCN (F(\xba) \xcp F^{2}(\xba))$ \xEH 0 \xEH 0 \xEH 0
\xEH 0 \xEH 5 \xEH 0 \xEP
\hline

$ \xbf'_{0}= \xCN (\xba \xcp F(\xba))$ \xEH 0 \xEH 0 \xEH 0 \xEH 0
\xEH 0 \xEH 5 \xEP
\hline
\end{tabular}

This allows definition by cases.
Suppose we want $ \xba $ to have the same result as $ \xbq $ if $ \xbq $
has value $p,$ and the same
result as $ \xbq' $ otherwise, more precisely:
\begin{flushleft}
\[ F_{\xba}(m):= \left\{ \begin{array}{lcl}
F_{\xbq}(m) \xEH iff \xEH F_{\xbq}(m)=p \xEP
\xEH \xEH \xEP
F_{\xbq'}(m) \xEH \xEH otherwise \xEP
\end{array}
\right.
\]
\end{flushleft}

then we define:

$ \xba:=(\xbf'_{n-p}(\xbq) \xcu \xbq) \xco ((\xCN \xbf'_{n-p}(
\xbq)) \xcu \xbq').$

As $ \xbf_{k}$ contains $ \xcp,$ but only 1 variable, there is no problem
for projections here.

Note the following important fact:

More generally,
we can construct in the same way new functions from old ones by cases,
like:
\begin{flushleft}
\[ F_{\xbq}(m):= \left\{ \begin{array}{lcl}
F_{\xbs}(m) \xEH iff \xEH F_{\xba}(m)=s \xEP
\xEH \xEH \xEP
F_{\xbt}(m) \xEH iff \xEH F_{\xba}(m)=t \xEP
\xEH \xEH \xEP
F_{\xbr}(m) \xEH \xEH otherwise \xEP
\end{array}
\right.
\]
\end{flushleft}

But we can $ \xCf not$ attribute arbitrary values as in

if condition1 holds, then x1, if condition2 holds, then x2, etc.

This is also reflected by the
fact that by the above, for 4 truth values, we can (using $ \xco)$ obtain
$\{0,3\} \xCK \{0,1,2,3\} \xCK \{0,2,3\} \xCK \{0,3\}$ $=$ $2*4*3*2$ $=$
48 semantically different
formulas with 1 variable, and no more (checked with a computer program).
Therefore, $F$ is weaker than $Z,$ which can generate arbitrary functions.

Thus, we can also define $ \xbs \xcp \xbt $ by cases, in a uniform way for
all $n.$
Consider, e.g., the case $v(\xbf)=3,$ $v(\xbq)=2.$ $v(\xbf \xcp \xbq
)$ should be equal to $v(\xbq),$
we take:
$ \xbf_{n-3}(\xbf) \xcu \xbf_{n-2}(\xbq) \xcu \xbq.$

We conclude by the following

\bfa

$\hspace{0.01em}$


\label{Fact F-J-A}

$F$ and $J+A$ are interdefinable:

$Aa \xcr \xCN (a \xcp Fa),$

$Ja \xcr Fa \xco Aa,$

$Fa \xcr Ja \xcu \xCN Aa.$

$ \xcz $
\\[3ex]
\subsubsection{
The additional operator $Z$
}

\efa

We turn to the operator $Z.$

\bd

$\hspace{0.01em}$


\label{Definition Rep-Z}

We introduce the following, derived, auxiliary, operators:

$S_{i}(\xbf):=n$ iff $v(\xbf)=i,$ and 0 otherwise, for $i=0, \Xl,n$

$K_{i}(\xbf):=i$ for any $ \xbf,$ for $i=0, \Xl,n$

\ed

\be

$\hspace{0.01em}$


\label{Example Rep-5}

We give here the example for $n=5.$

\vspace{5mm}
\begin{tabular}{|c|c|c|c|c|c|c|c|}

\hline

\xEH a \xEH 0 \xEH 1 \xEH 2 \xEH 3 \xEH 4 \xEH 5 \xEP

\hline

$Z$ \xEH   \xEH 1 \xEH 2 \xEH 3 \xEH 4 \xEH 5 \xEH 0 \xEP

\hline
\hline

       \xEH a \xEH 0 \xEH 1 \xEH 2 \xEH 3 \xEH 4 \xEH 5 \xEP

\hline

$K_0$  \xEH   \xEH 0 \xEH 0 \xEH 0 \xEH 0 \xEH 0 \xEH 0 \xEP

\hline

$K_1$  \xEH   \xEH 1 \xEH 1 \xEH 1 \xEH 1 \xEH 1 \xEH 1 \xEP

\hline

$K_2$  \xEH   \xEH 2 \xEH 2 \xEH 2 \xEH 2 \xEH 2 \xEH 2 \xEP

\hline

$K_3$  \xEH   \xEH 3 \xEH 3 \xEH 3 \xEH 3 \xEH 3 \xEH 3 \xEP

\hline

$K_4$  \xEH   \xEH 4 \xEH 4 \xEH 4 \xEH 4 \xEH 4 \xEH 4 \xEP

\hline

$K_5$  \xEH   \xEH 5 \xEH 5 \xEH 5 \xEH 5 \xEH 5 \xEH 5 \xEP

\hline
\hline

       \xEH a \xEH 0 \xEH 1 \xEH 2 \xEH 3 \xEH 4 \xEH 5 \xEP

\hline

$S_0$  \xEH   \xEH 5 \xEH 0 \xEH 0 \xEH 0 \xEH 0 \xEH 0 \xEP

\hline

$S_1$  \xEH   \xEH 0 \xEH 5 \xEH 0 \xEH 0 \xEH 0 \xEH 0 \xEP

\hline

$S_2$  \xEH   \xEH 0 \xEH 0 \xEH 5 \xEH 0 \xEH 0 \xEH 0 \xEP

\hline

$S_3$  \xEH   \xEH 0 \xEH 0 \xEH 0 \xEH 5 \xEH 0 \xEH 0 \xEP

\hline

$S_4$  \xEH   \xEH 0 \xEH 0 \xEH 0 \xEH 0 \xEH 5 \xEH 0 \xEP

\hline

$S_5$  \xEH   \xEH 0 \xEH 0 \xEH 0 \xEH 0 \xEH 0 \xEH 5 \xEP

\hline

\end{tabular}

\ee

\bfa

$\hspace{0.01em}$


\label{Fact Rep-Z}

(1) We can define $S_{i}(\xbf)$ and $K_{i}(\xbf)$ for $0 \xck i \xck
n$ from $ \xCN, \xcu, \xco,Z.$

(2) We can define any $m-$ary truth function from $ \xCN, \xcu, \xco
,Z.$

\efa

\subparagraph{
Proof
}

$\hspace{0.01em}$


(1)

$S_{i}(\xbf)= \xCN Z^{n-i}(\xbf),$ $K_{i}(\xbf)=Z^{i+1}(\xCN (\xbf
\xcu \xCN \xbf))$

(2)

Suppose $ \xBc i_{1}, \Xl,i_{m} \xBe $ should have value $i,$
$ \xBc i_{1}, \Xl,i_{m} \xBe  \xcZ i,$ we can express this by

$S_{i_{1}}(x_{1}) \xcu, \Xl, \xcu S_{i_{m}}(x_{m}) \xcu K_{i}.$

We then take the disjunction of all such expressions:

$ \xcO \{S_{i_{1}}(x_{1}) \xcu, \Xl, \xcu S_{i_{m}}(x_{m}) \xcu K_{i}:$
$ \xBc i_{1}, \Xl,i_{m} \xBe  \xcZ i\}$

\bco

$\hspace{0.01em}$


\label{Corollary Rep-Log-Int}

Any model function is definable from $Z,$ so any semantical interpolant is
also
a syntactic one.
$ \xcz $
\\[3ex]
\subsection{
The three valued intuitionistic logic Here/There HT
}

\label{Section HT}

\eco

We now give a short introduction to the well-known 3-valued
intuitionistic logic HT (Here/There), with some results also for similar
logics with more than 3 values. Many of these properties were found and
checked
with a small computer program. In particular, we show the existence of a
normal
form, similar to classical propositional logic, but $ \xcp $ cannot be
eliminated.
Consequently, we cannot always separate propositional variables easily.

Our main result here (which is probably well known, we claim no priority)
is
that ``forgetting'' a variable preserves definability in the following
sense:
Let, e.g., $ \xbf =a \xcu b,$ and $M(\xbf)$ be the set of models where $
\xbf $ has maximal
truth value (2 here), then there is $ \xbf' $ such that the set of models
where
$ \xbf' $ has value 2 is the set of all models which agree with a model
of $ \xbf $ on,
e.g., $b.$ We ``forget'' about $ \xCf a.$
Our $ \xbf' $ is here, of course, $b.$ Here, the problem is trivial, it
is a
bit less so when $ \xcp $ is involved, as we cannot always separate the
two parts.
For example, the result of ``forgetting'' about $ \xCf a$ in the formula $a
\xcu (a \xcp b)$ is
$b,$ in the formula $a \xcp b$ it is TRUE. Thus, forgetting about a
variable
preserves definability, and the abovementioned semantical interpolation
property
carries over to the syntactic side, similarly to the result on classical
logic, see Section \ref{Section Class-Proj} (page \pageref{Section Class-Proj})
.

\bfa

$\hspace{0.01em}$


\label{Fact Eq-Def}

These results were checked with a small computer program:

(1) With 2 variables $ \xCf a,b$ are definable, using the operators $ \xCN
,$ $ \xcp,$ $ \xcu,$ $ \xco,$
174 semantically different formulas.
$ \xco $ is not needed, i.e. with or without $ \xco $ we have the
same set of definable formulas. We have, e.g.,
$a \xco b$ $ \xcr $ $ \xCI \xCI b \xcp (\xCN \xCN a \xcp a) \xCJ \xcp
\xCI (\xCN a \xcp b) \xcu (\xCN \xCN a \xcp a) \xCJ \xCJ.$

(2) With the operators $ \xCN,$ $ \xcu,$ $ \xco $ only 120
semantically different formulas are definable. Thus, $ \xcp $ cannot be
expressed by the other operators.

\efa

\bfa

$\hspace{0.01em}$


\label{Fact Eq-Equiv-2}

$((a \xcp b) \xcp c)$ $ \xcr $ $((\xCN a \xcp c) \xcu (b \xcp c) \xcu (a
\xco \xCN b \xco c))$

holds in the 3 valued case.

(Thanks to D.Pearce for telling us.)

\efa

This Fact is well known, we have verified it by computer, but not by hand.

\bco

$\hspace{0.01em}$


\label{Corollary Eq-Equiv-2}

The 3 valued case has interpolation.

\eco

\subparagraph{
Proof
}

$\hspace{0.01em}$


By Fact \ref{Fact Eq-Equiv-2} (page \pageref{Fact Eq-Equiv-2}), we can flatten
nested $ \xcp $'s, so
we have projection.
$ \xcz $
\\[3ex]
\subsection{
Finite Goedel logics with 4 truth values
}

\label{Section 4-Val}

Interpolation fails for:

$(d \xcp (((a \xcp b) \xcp a) \xcp a) \xcp d$ $ \xcl $ $(((c \xcp d) \xcp
c) \xcp c.$

Consider the table in
Corollary \ref{Corollary Mod-Int-Fail} (page \pageref{Corollary Mod-Int-Fail}),
and the comment about possible
interpolants
after the table.
By the proof of Fact \ref{Fact Mod-Fin-Goed} (page \pageref{Fact Mod-Fin-Goed})
, (0), we see that above
formulas
have no interpolant. But it is trivial to see that $ \xCf Jp$ will be an
interpolant,
see Definition \ref{Definition Mod-Fin-Goed-Add} (page \pageref{Definition
Mod-Fin-Goed-Add}).

Note that the implication is not true for more than 4 truth values, as we
saw
in Corollary \ref{Corollary Mod-Int-Fail} (page \pageref{Corollary
Mod-Int-Fail}),
so this example will not be a counterexample to interpolation any
more.

We have checked with a computer program, but not by hand:

Introducing a new constant, 1, which has always truth value 1, and is thus
simpler than above operator $J,$ gives exactly 2 different interpolants
for
the formulas of Example \ref{Example Sin-4-Value} (page \pageref{Example
Sin-4-Value}):
$(p \xcp 1) \xcp p,$ and $(p \xcp 1) \xcp 1.$
Introducing an additional constant 2 will give still other interpolants.
But if one is permitted in classical logic to use
the constants TRUE and FALSE for interpolation, why not 1 and 2 here?

$ \xCO $

$ \xCO $
\clearpage
\chapter{
Laws about size and interpolation in non-monotonic logics
}

\label{Chapter Size-Laws}
\section{
Introduction
}
\subsection{
Various concepts of size and non-monotonic logics
}

A natural interpretation of the non-monotonic rule $ \xbf \xcn \xbq $ is
that
the set of exceptional cases, i.e., those where $ \xbf $ holds, but not $
\xbq,$ is
a small subset of all the cases where $ \xbf $ holds, and the complement,
i.e.,
the set of cases where $ \xbf $ and $ \xbq $ hold, is a big subset of all
$ \xbf -$cases.

This interpretation gives an abstract semantics to non-monotonic logic,
in the sense that definitions and rules are translated to rules about
model sets, without any structural justification of those rules, as they
are given, e.g., by preferential structures, which provide structural
semantics. Yet, they are extremely useful, as they allow us to concentrate
on the essentials, forgetting about syntactical reformulations of
semantically
equivalent formulas, the laws derived from the standard proof theoretical
rules
incite to generalizations and modifications, and reveal deep connections
but
also differences. One of those insights is the connection between laws
about size and (semantical) interpolation for non-monotonic logics.

To put this abstract view a little more into perspective, we present
three alternative systems, also working with abstract size as a semantics
for non-monotonic logics. (They were already mentioned in
Section \ref{Section Laws-Size-Intro} (page \pageref{Section Laws-Size-Intro})
.)

 \xEI
 \xDH
the system of one of the authors for a first order setting,
published in  \cite{Sch90} and elaborated in
 \cite{Sch95-1},
 \xDH
the system of S.Ben-David and R.Ben-Eliyahu,
published in  \cite{BB94},
 \xDH
the system of N.Friedman and J.Halpern,
published in  \cite{FH96}.
 \xEJ

 \xEh

 \xDH
Defaults as generalized quantifiers:

We first recall the definition of a
``weak filter'', made official in
Definition \ref{Definition Weak-Filter} (page \pageref{Definition Weak-Filter})
:

Fix a base set $X.$
A weak filter on or over $X$ is a set $ \xdf \xcc \xdp (X),$ s.t. the
following
conditions hold:

(F1) $X \xbe \xdf $

(F2) $A \xcc B \xcc X,$ $A \xbe \xdf $ imply $B \xbe \xdf $

$(F3')$ $A,B \xbe \xdf $ imply $A \xcs B \xEd \xCQ.$

We use
weak filters on the semantical side, and add the following axioms on the
syntactical side to a FOL axiomatisation:

\index{$\xeA$}
1. $ \xeA x \xbf (x)$ $ \xcu $ $ \xcA x(\xbf (x) \xcp \xbq (x))$ $ \xcp $
$ \xeA x \xbq (x),$

2. $ \xeA x \xbf (x)$ $ \xcp $ $ \xCN \xeA x \xCN \xbf (x),$

3. $ \xcA x \xbf (x)$ $ \xcp $ $ \xeA x \xbf (x)$ and $ \xeA x \xbf (x)$ $
\xcp $ $ \xcE x \xbf (x).$

A model is now a pair, consisting of a classical FOL model $M,$ and a weak
filter over its universe. Both sides are connected by the following
definition,
where $ \xdn (M)$ is the weak filter on the universe of the classical
model $M:$

\index{$ \xBc M,\xdn (M) \xBe $}
$ \xBc M, \xdn (M) \xBe $ $ \xcm $ $ \xeA x \xbf (x)$ iff there is $A \xbe \xdn
(M)$
s.t. $ \xcA a \xbe A$ $(\xBc M, \xdn (M) \xBe $ $ \xcm $ $ \xbf [a]).$

Soundness and completeness is shown in
 \cite{Sch95-1}, see also  \cite{Sch04}.

The extension to defaults with prerequisites by restricted quantifiers is
straightforward.

 \xDH
The system of $S.$ Ben-David and $R.$ Ben-Eliyahu:

Let $ \xdn':=\{ \xdn' (A):$ $A \xcc U\}$ be a system of filters for $
\xdp (U),$ i.e. each $ \xdn' (A)$ is a
filter over A. The conditions are (in slight modification):

UC': $B \xbe \xdn' (A)$ $ \xcp $ $ \xdn' (B) \xcc \xdn' (A),$

DC': $B \xbe \xdn' (A)$ $ \xcp $ $ \xdn' (A) \xcs \xdp (B) \xcc \xdn'
(B),$

RBC': $X \xbe \xdn' (A),$ $Y \xbe \xdn' (B)$ $ \xcp $ $X \xcv Y \xbe
\xdn' (A \xcv B),$

SRM': $X \xbe \xdn' (A),$ $Y \xcc A$ $ \xcp $ $A-Y \xbe \xdn' (A)$ $
\xco $ $X \xcs Y \xbe \xdn' (Y),$

\index{$UC'$}
\index{$DC'$}
\index{$RBC'$}
\index{$SRM'$}
\index{$GTS'$}
GTS': $C \xbe \xdn' (A),$ $B \xcc A$ $ \xcp $ $C \xcs B \xbe \xdn' (B).$

 \xDH

The system of $N.$ Friedman and $J.$ Halpern:

Let $U$ be a set, $<$ a strict partial order on $ \xdp (U),$
(i.e. $<$ is transitive, and contains no cycles).
Consider the following conditions for $<:$

(B1) $A' \xcc A<B \xcc B' $ $ \xcp $ $A' <B',$

(B2) if $A,B,C$ are pairwise disjoint, then $C<A \xcv B,$ $B<A \xcv C$ $
\xcp $ $B \xcv C<A,$

(B3) $ \xCQ <X$ for all $X \xEd \xCQ,$

(B4) $A<B$ $ \xcp $ $A<B-$A,

\index{$(B1)$}
\index{$(B2)$}
\index{$(B3)$}
\index{$(B4)$}
\index{$(B5)$}
(B5) Let $X,Y \xcc A.$ If $A-X<X,$ then $Y<A-Y$ or $Y-X<X \xcs Y.$

 \xEj

The equivalence of the systems of  \cite{BB94}
and  \cite{FH96} was shown in
 \cite{Sch97-4}, see also  \cite{Sch04}.

Historical remarks:
Our own view as abstract size was inspired by the classical filter
approach,
as used e.g. in mathematical measure theory.
The first time that abstract size was related to
nonmonotonic logics was, to our knowledge,
in the second author's  \cite{Sch90} and  \cite{Sch95-1}, and,
independently, in  \cite{BB94}.
The approach to size by partial orders is
first discussed - to our knowledge - by
N.Friedman and J.Halpern, see  \cite{FH96}.
More detailed remarks can also be found in  \cite{GS08c},
 \cite{GS09a},  \cite{GS08f}.
A somewhat different approach is taken in  \cite{HM07}.

Before we introduce the connection between interpolation and
multiplicative
laws about size, we give now some comments on the laws about size
themselves.
\subsection{
Additive and multiplicative laws about size
}

We give here a short introduction to and some examples for additive
and multiplicative laws about size. A detailed overview is presented
in Table \ref{Table Base2-Size-Rules-1} (page \pageref{Table
Base2-Size-Rules-1}),
Table \ref{Table Base2-Size-Rules-2} (page \pageref{Table Base2-Size-Rules-2}),
and
Table \ref{Table Mul-Laws} (page \pageref{Table Mul-Laws}).
(The first two tables have to be read together, they are too
big to fit on one page.)

They show connections and how to develop a multitude of logical rules
known from nonmonotonic logics by combining a small number of principles
about size. We can use them as building blocks to construct the rules
from.
More precisely, ``size'' is to be read as
``relative size'', since it is essential to change the base sets.

In the first two tables, these principles are some basic and very natural
postulates,
$ \xCf (Opt),$ $ \xCf (iM),$ $(eM \xdi),$ $(eM \xdf),$ and a continuum
of power of the notion of
``small'', or, dually, ``big'', from $(1*s)$ to $(< \xbo *s).$
From these, we can develop the rest except, essentially, Rational
Monotony,
and thus an infinity of different rules.

The probably easiest way to see a connection between non-monotonic logics
and abstract size is by considering preferential structures.
Preferential structures
define principal filters, generated by the set of minimal elements, as
follows:
if $ \xbf \xcn \xbq $ holds in such a structure, then $ \xbm (\xbf) \xcc
M(\xbq),$ where $ \xbm (\xbf)$ is the
set of minimal elements of $M(\xbf).$ According to our ideas, we define
a principal filter $ \xdf $ over $M(\xbf)$ by $X \xbe \xdf $ iff $ \xbm
(\xbf) \xcc X \xcc M(\xbf).$ Thus,
$M(\xbf) \xcs M(\xCN \xbq)$ will be a ``small'' subset of $M(\xbf).$
(Recall that
filters contain the ``big'' sets, and ideals the ``small'' sets.)

We can now go back and forth between rules on size and logical rules,
e.g.:

(For details, see
Table \ref{Table Base2-Size-Rules-1} (page \pageref{Table Base2-Size-Rules-1}),
Table \ref{Table Base2-Size-Rules-2} (page \pageref{Table Base2-Size-Rules-2}),
and
Table \ref{Table Mul-Laws} (page \pageref{Table Mul-Laws}).)

 \xEh

 \xDH
The ``AND'' rule corresponds to the filter property (finite intersections of
big subsets are still big).

 \xDH
``Right weakening'' corresponds to the rule that supersets of big sets
are still big.

 \xDH
It is natural, but beyond filter properties themselves, to postulate that,
if $X$ is a small subset of $Y,$ and $Y \xcc Y',$ then $X$ is also a
small subset of $Y'.$
We call such properties
``coherence properties'' between filters.
This property corresponds to the logical rule $ \xCf (wOR).$

 \xDH
In the rule $(CM_{ \xbo }),$ usually called Cautious Monotony, we change
the
base set a little when going from $M(\xba)$ to $M(\xba \xcu \xbb)$
(the change is small
by the prerequisite $ \xba \xcn \xbb),$ and still have $ \xba \xcu \xbb
\xcn \xbb',$ if we had
$ \xba \xcn \xbb'.$ We see here a conceptually very different use of
``small'', as we now change the base set, over which the filter is
defined, by a small amount.

 \xDH
The rule of Rational Monotony is the last one in the first table, and
somewhat
isolated there. It
is better to be seen as a multiplicative law, as described in the third
table.
It corresponds to the rule that the product of medium (i.e, neither big
nor small) sets, has still medium size.

 \xEj
\subsection{
Interpolation and size
}

The connection between non-monotonic logic and
the abstract concept of size was investigated in
 \cite{GS09a}, see also  \cite{GS08f}.
There, we looked among other things at abstract addition
of size. Here, we will show a connection to abstract multiplication of
size.
Our semantic approach used decomposition of set theoretical products.
An important step was to write a set of models $ \xbS $ as a product of
some
set $ \xbS' $ (which was a restriction of $ \xbS),$ and some full
Cartesian product.
So, when we speak about size, we will have (slightly simplified) some big
subset
$ \xbS_{1}$ of one product $ \xbP_{1},$ and some big subset $ \xbS_{2}$ of
another product $ \xbP_{2},$
and will now check whether $ \xbS_{1} \xCK \xbS_{2}$ is a big subset of $
\xbP_{1} \xCK \xbP_{2}.$
In shorthand,
whether ``$big*big=big$''.
(See Definition \ref{Definition Sin-Size-Rules} (page \pageref{Definition
Sin-Size-Rules})  for precise definitions.)
Such conditions are called
coherence conditions, as they do not concern the notion of size itself,
but the way the sizes defined for different base sets are connected.
Our main results here are
Proposition \ref{Proposition Sin-Interpolation-1} (page \pageref{Proposition
Sin-Interpolation-1})  and
Proposition \ref{Proposition Mul-Mu*1-Int} (page \pageref{Proposition
Mul-Mu*1-Int}). They say that if the logic under
investigation is defined from a notion of size which satisfies
sufficiently
many conditions, then this logic will have interpolation of type one or
even
two.

Consider now some set product $X \xCK X'.$ (Intuitively, $X$ and $X' $
are model sets
on sublanguages $J$ and $J' $ of the whole language $L.)$ When we have now
a rule like: If $Y$ is a big subset of $X,$ and $Y' $ a big subset of $X'
,$ then
$Y \xCK Y' $ is a big subset of $X \xCK X',$ and conversely, we can
calculate consequences
separately in the sublanguages, and put them together to have the overall
consequences. But this is the principle behind interpolation: we can work
with independent parts.

This is made precise in
Definition \ref{Definition Sin-Size-Rules} (page \pageref{Definition
Sin-Size-Rules}),
in particular by the rule

$(\xbm *1):$ $ \xbm (X \xCK X')= \xbm (X) \xCK \xbm (X').$

(Note that the conditions $(\xbm *i)$ and $(\xbS *i)$ are equivalent, as
shown
in Proposition \ref{Proposition Sin-Product-Small} (page \pageref{Proposition
Sin-Product-Small})  (for principal filters).)

The main result is
that the multiplicative size rule $(\xbm *1)$ entails
non-monotonic interpolation of the form $ \xbf \xcn \xba \xcn \xbq,$ see
Proposition \ref{Proposition Mul-Mu*1-Int} (page \pageref{Proposition
Mul-Mu*1-Int}).

We take now a closer look at interpolation for non-monotonic
logic.
\paragraph{
The three variants of interpolation \\[2mm]
}

Consider preferential logic, a rule like $ \xbf \xcn \xbq.$ This means
that
$ \xbm (\xbf) \xcc M(\xbq).$ So we go from $M(\xbf)$ to $ \xbm (
\xbf),$ the minimal models of $ \xbf,$
and then to $M(\xbq),$ and, abstractly, we have $M(\xbf) \xcd \xbm (
\xbf) \xcc M(\xbq),$ so we have
neither necessarily $M(\xbf) \xcc M(\xbq),$ nor $M(\xbf) \xcd M(
\xbq),$ the relation between
$M(\xbf)$ and $M(\xbq)$ may be more complicated. Thus, we have neither
the monotone,
nor the antitone case. For this reason, our general results for monotone
or antitone logics do not hold any more.

But we also see here that classical logic is used, too. Suppose that
there is $ \xbf' $ which describes exactly $ \xbm (\xbf),$ then we can
write
$ \xbf \xcn \xbf' \xcl \xbq.$

So we can split preferential logic into a
core part - going from $ \xbf $ to its minimal models - and a second part,
which is
just classical logic. (Similar decompositions are
also natural for other non-monotonic logics.)
Thus, preferential logic can be seen as a combination of two logics, the
non-monotonic core, and classical logic.
It is thus natural to consider variants of the interpolation problem,
where $ \xcn $ denotes again preferential logic, and $ \xcl $ as usual
classical
logic:

Given $ \xbf \xcn \xbq,$ is there ``simple'' $ \xba $ such that
 \xEh

 \xDH
$ \xbf \xcn \xba \xcl \xbq,$ or

 \xDH
$ \xbf \xcl \xba \xcn \xbq,$ or

 \xDH
$ \xbf \xcn \xba \xcn \xbq?$

 \xEj

In most cases, we will only consider the semantical version,
as the problems of the syntactical version are very similar to those
for monotonic logics.
We turn to the variants.

 \xEh

 \xDH

The first variant, $ \xbf \xcn \xba \xcl \xbq,$ has a complete
characterization
in Proposition \ref{Proposition Sin-Non-Mon-Int-Karl} (page \pageref{Proposition
Sin-Non-Mon-Int-Karl}),
provided we have a suitable normal form (conjunctions of disjunctions).
The condition says that the relevant variables of $ \xbm (\xbf)$ have to
be
relevant for $M(\xbf).$

 \xDH

The second variant, $ \xbf \xcl \xba \xcn \xbq,$ is related to very (and
in many cases,
too) strong conditions about size.
We do not have a complete characterization, only sufficient conditions
about
size.
The size conditions we need are
(see Definition \ref{Definition Sin-Size-Rules} (page \pageref{Definition
Sin-Size-Rules})):

the abovementioned $(\xbm *1),$ and,

$(\xbm *2):$ $ \xbm (X) \xcc Y$ $ \xch $ $ \xbm (X \xex A) \xcc Y \xex A$

where $X$ need not be a product any more.

The result is given in
Proposition \ref{Proposition Sin-Interpolation-1} (page \pageref{Proposition
Sin-Interpolation-1}).

Example \ref{Example Sin-Prod-Size} (page \pageref{Example Sin-Prod-Size}) 
shows that $(\xbm *2)$ seems too
strong
when compared to probability defined size.

We should, however, note that sufficiently modular preferential relations
guarantee these very strong properties of the big sets, see
Section \ref{Section Ham-Rel-Dist} (page \pageref{Section Ham-Rel-Dist}).

 \xDH

We turn to the third variant, $ \xbf \xcn \xba \xcn \xbq.$ This is
probably the
most interesting one, as it is more general, loosens the connection
with classical logic, seems more natural as a rule,
and is also connected to more natural laws about size.
Again, we do not have a complete characterization, only sufficient
conditions
about size.
Here, $(\xbm *1)$ suffices, and we have our main result about
non-monotonic semanti interpolation,
Proposition \ref{Proposition Mul-Mu*1-Int} (page \pageref{Proposition
Mul-Mu*1-Int}),
that $(\xbm *1)$ entails interpolation of the type $ \xbf \xcn \xba \xcn
\xbq.$

Proposition \ref{Proposition Mod-Hamming} (page \pageref{Proposition
Mod-Hamming})  shows that $(\xbm *1)$ is
(roughly) equivalent to
the properties

$ \xCf (GH1)$ $ \xbs \xec \xbt $ $ \xcu $ $ \xbs' \xec \xbt' $ $ \xcu $
$(\xbs \xeb \xbt $ $ \xco $ $ \xbs' \xeb \xbt')$ $ \xch $ $ \xbs \xbs
' \xeb \xbt \xbt' $

(where $ \xbs \xec \xbt $ iff $ \xbs \xeb \xbt $ or $ \xbs = \xbt)$

$ \xCf (GH2)$ $ \xbs \xbs' \xeb \xbt \xbt' $ $ \xch $ $ \xbs \xeb \xbt $
$ \xco $ $ \xbs' \xeb \xbt' $

of a preferential relation.

$(\xCf (GH2)$ means that some compensation is possible, e.g., $ \xbt \xeb
\xbs $ might be the
case, but $ \xbs' \xeb \xbt' $ wins in the end, so $ \xbs \xbs' \xeb
\xbt \xbt'.)$

There need not always be a semantical interpolation for the third variant,
this
is shown in
Example \ref{Example Sin-Non-Mon-Int} (page \pageref{Example Sin-Non-Mon-Int}).

 \xEj

So we see that, roughly,
semantic interpolation for nonmonotonic logics works when abstract size is
defined in a modular way - and we see independence again.
In a way, this is not surprising, as we use independent definition of
validity for interpolation in classical logic, and we use independent
definition of additional structure (relations or size) for
interpolation in non-monotonic logic.
\subsection{
Hamming relations and size
}

As preferential relations are determined by a relation, and give rise
to abstract notions of size and their manipulation, it is natural to
take a close look at the corresponding properties of the relation.
We already gave a few examples in the preceding sections, so we can be
concise here. Our main definitions and results on this subject are
to be found in
Section \ref{Section Ham-Rel-Dist} (page \pageref{Section Ham-Rel-Dist}), where
we also discuss distances with
similar properties.

It is not surprising that we find various types of Hamming relations and
distances in this context, as they are, by definition, modular. Neither is
it surprising that we see them again in
Chapter \ref{Chapter Neighbourhood} (page \pageref{Chapter Neighbourhood}), as
we are interested there in
independent ways
to define neighbourhoods.

Basically, these relations and distances come in two flavours, the
set and the counting variant. This is perhaps best illustrated by the
Hamming
distance of two sequence of finite, equal length. We can define the
distance
by the $ \xCf set$ of arguments where they differ, or by the $ \xCf
cardinality$
of this set. The first results in possibly incomparable distances, the
second allows ``compensation'', difference in one argument can be
compensated by equality in another argument.

For definitions and results, also those connecting them to notions of
size, see Section \ref{Section Ham-Rel-Dist} (page \pageref{Section
Ham-Rel-Dist})  in particular
Definition \ref{Definition Sin-Set-HR} (page \pageref{Definition Sin-Set-HR}).
We then show in
Proposition \ref{Proposition Mod-Hamming} (page \pageref{Proposition
Mod-Hamming})
that (smooth) Hamming relations
generate our size conditions when size is defined as above from
a relation (the set of preferred elements generates the principal
filter). Thus, Hamming relations determine logics which have
interpolation, see
Corollary \ref{Corollary Sin-Interpolation-2} (page \pageref{Corollary
Sin-Interpolation-2}).
\subsection{
Equilibrium logic
}

Equilibrium logic, due to D.Pearce, A.Valverde, see
 \cite{PV09} for motivation and further discussion, is based
on the 3-valued finite Goedel logic, also called HT logic,
HT for ``here and there''. Our results are presented in
Section \ref{Section EQ} (page \pageref{Section EQ}).

Equilibrium logic (EQ) is defined by a choice function on the model set.
First models have to be ``total'', no variable of the language
may have 1 as value. Second, if $m \xeb m',$ then $m$ is considered
better,
and $m' $ discarded, where $m \xeb m' $ iff $m$ and $m' $ give value 0 to
the same
variables, and $m$ gives value 2 to strictly less (as subset) variables
than $m' $ does.

We can define equilibrium logic by a preferential relation (taking care
also
of the first condition), but it is not smooth. Thus, our general results
from the beginning of this section will not hold, and we have to work
with ``hand knitted'' solutions. We first show that equilibrium logic
has no interpolation of the form $ \xbf \xcl \xba \xcn \xbq $ or $ \xbf
\xcn \xba \xcl \xbq,$ then that is has
interpolation of the form $ \xbf \xcn \xba \xcn \xbq,$ and that the
interpolant is also
definable, i.e., equilibrium logic has semantic and syntactic
interpolation
of this form. Essentially, semantic interpolation is due to the fact that
the preference relation is defined in a modular way, using individual
variables - as always, when we have interpolation.
\subsection{
Interpolation for revision and argumentation
}

We have a short and simple result
(Lemma \ref{Lemma Mul-TR} (page \pageref{Lemma Mul-TR}))
for interpolation in AGM revision.
Unfortunately, we need the variables from both sides of the
revision operator as can easily be seen by revising with TRUE.
The reader is referred to
Section \ref{Section Int-Dist-Rev} (page \pageref{Section Int-Dist-Rev})  for
details.

Somewhat surprisingly, we also have an interpolation result for one
form of argumentation, where we consider the set of arguments for a
statement
as the truth value of that statement. As we have maximum (set union), we
have the lower bound used in
Proposition \ref{Proposition Gin-Pr-Int} (page \pageref{Proposition Gin-Pr-Int})
 for the monotonic case, and can
show
Fact \ref{Fact Arg-Int} (page \pageref{Fact Arg-Int}).
See Section \ref{Section Inter-Arg} (page \pageref{Section Inter-Arg})  for
details.
\subsection{
Language change to obtain products
}

To achieve interpolation and other results of independence, we
often need to write a set of models as a non-trivial product.
Sometimes, this is impossible, but an equivalent reformulation
of the language can solve the problem, see
Example \ref{Example Mod-Lang-Fact} (page \pageref{Example Mod-Lang-Fact}).

Crucial there is that $6=3*2,$ so we can just re-arrange the 6 models in a
different way, see
Fact \ref{Fact Mod-Lang-Fact} (page \pageref{Fact Mod-Lang-Fact}).

A similar result holds for the non-monotonic case, where the
structure must be possible, we can then redefine the language.

All details are to be found in
Section \ref{Section Lang-Manip} (page \pageref{Section Lang-Manip}).
\section{
Laws about size
}
\subsection{
Additive laws about size
}

\label{Section Add-Size}

We now give the main additive rules for
manipulation of abstract size from  \cite{GS09a}, see
Table \ref{Table Base2-Size-Rules-1} (page \pageref{Table Base2-Size-Rules-1}) 
and
Table \ref{Table Base2-Size-Rules-2} (page \pageref{Table Base2-Size-Rules-2}),
``Rules on size''.

The notation is explained with some redundancy, so the reader will not
have to leaf back and forth to
Chapter \ref{Chapter Mod-Base-Def} (page \pageref{Chapter Mod-Base-Def}).

$ \xCO $
\subsubsection{
Notation
}

\label{Section Mul-Nota}

 \xEh

 \xDH

$ \xdp (X)$ is the power set of $X,$ $ \xcc $ is the subset relation, $
\xcb $ the strict part of
$ \xcc,$ i.e. $A \xcb B$ iff $A \xcc B$ and $A \xEd B.$
The operators $ \xcu,$ $ \xCN,$ $ \xco,$ $ \xcp $ and $ \xcl $ have
their usual, classical interpretation.

 \xDH

$ \xdi (X) \xcc \xdp (X)$ and $ \xdf (X) \xcc \xdp (X)$ are dual abstract
notions of size, $ \xdi (X)$ is the
set of ``small'' subsets of $X,$ $ \xdf (X)$ the set of ``big'' subsets of
$X.$ They are
dual in the sense that $A \xbe \xdi (X) \xcj X-A \xbe \xdf (X).$ ``$ \xdi
$'' evokes ``ideal'',
``$ \xdf $'' evokes ``filter'' though the full strength of both is reached
only
in $(< \xbo *s).$ ``s'' evokes ``small'', and ``$(x*s)$'' stands for
``$x$ small sets together are still not everything''.

 \xDH

If $A \xcc X$ is neither in $ \xdi (X),$ nor in $ \xdf (X),$ we say it has
medium size, and
we define $ \xdm (X):= \xdp (X)-(\xdi (X) \xcv \xdf (X)).$ $
\xdm^{+}(X):= \xdp (X)- \xdi (X)$ is the set of subsets
which are not small.

 \xDH

$ \xeA x \xbf $ is a generalized first order quantifier, it is read
``almost all $x$ have property $ \xbf $''. $ \xeA x(\xbf: \xbq)$ is the
relativized version, read:
``almost all $x$ with property $ \xbf $ have also property $ \xbq $''. To
keep the table
``Rules on size''
simple, we write mostly only the non-relativized versions.
Formally, we have $ \xeA x \xbf: \xcj \{x: \xbf (x)\} \xbe \xdf (U)$
where $U$ is the universe, and
$ \xeA x(\xbf: \xbq): \xcj \{x:(\xbf \xcu \xbq)(x)\} \xbe \xdf (\{x:
\xbf (x)\}).$
Soundness and completeness results on $ \xeA $ can be found in
 \cite{Sch95-1}.

 \xDH

Analogously, for propositional logic, we define:

$ \xba \xcn \xbb $ $: \xcj $ $M(\xba \xcu \xbb) \xbe \xdf (M(\xba)),$

where $M(\xbf)$ is the set of models of $ \xbf.$

 \xDH

In preferential structures, $ \xbm (X) \xcc X$ is the set of minimal
elements of $X.$
This generates a principal filter by $ \xdf (X):=\{A \xcc X: \xbm (X) \xcc
A\}.$ Corresponding
properties about $ \xbm $ are not listed systematically.

 \xDH

The usual rules $ \xCf (AND)$ etc. are named here $(AND_{ \xbo }),$ as
they are in a
natural ascending line of similar rules, based on strengthening of the
filter/ideal properties.

 \xDH

For any set of formulas $T,$ and any consequence relation $ \xcn,$ we
will use
$ \ol{T}:=\{ \xbf:T \xcl \xbf \},$ the set of classical consequences of
$T,$ and
$ \ol{ \ol{T} }:=\{ \xbf:T \xcn \xbf \},$ the set of consequences of $T$
under the relation $ \xcn.$

 \xDH

We say that a set $X$ of models is definable by a formula (or a theory)
iff
there is a formula $ \xbf $ (a theory $T)$ such that $X=M(\xbf),$ or
$X=M(T),$ the set of
models of $ \xbf $ or $T,$ respectively.

 \xDH

Most rules are explained in the
table ``Logical rules'',
and ``RW'' stands for Right Weakening.

 \xEj


\subsubsection{
The groupes of rules
}

The rules concern properties of $ \xdi (X)$ or $ \xdf (X),$ or
dependencies
between such properties for different $X$ and $Y.$ All $X,Y,$ etc. will
be subsets of some universe, say $V.$ Intuitively, $V$ is the set of
all models of some fixed propositional language. It is not
necessary to consider all subsets of $V,$ the intention is to consider
subsets of $V,$ which are definable by a formula or a theory.
So we assume all $X,Y$ etc. taken from some $ \xdy \xcc \xdp (V),$ which
we call the domain.
In the former case, $ \xdy $ is closed under set difference, in
the latter case not necessarily so. (We will mention it when we need some
particular closure property.)

The rules are divided into 5 groups:

 \xEh

 \xDH $ \xCf (Opt),$ which says that ``All'' is optimal - i.e. when there
are no
exceptions, then a soft rule $ \xcn $ holds.

 \xDH 3 monotony rules:

 \xEh
 \xDH $ \xCf (iM)$ is inner monotony, a subset of a small set is small,
 \xDH $(eM \xdi)$ external monotony for ideals: enlarging the base set
keeps small
sets small,
 \xDH $(eM \xdf)$ external monotony for filters: a big subset stays big
when the base
set shrinks.
 \xEj

These three rules are very natural if ``size'' is anything coherent over
change
of base sets. In particular, they can be seen as weakening.

 \xDH $(\xCd)$ keeps proportions, it is here mainly to point the
possibility out.

 \xDH a group of rules $x*s,$ which say how many small sets will not yet
add to
the base set. The notation
``$(< \xbo *s)$'' is an allusion to the full filter property, that
filters are closed under $ \xCf finite$ intersections.
 \xDH Rational monotony, which can best be understood as robustness of $
\xdm^{+},$
see $(\xdm^{++})(3).$

 \xEj

We will assume all base sets to be non-empty in order to avoid pathologies
and
in particular clashes between $ \xCf (Opt)$ and $(1*s).$

Note that the full strength of the usual definitions of
a filter \index{filter}  and an
ideal \index{ideal}  are reached only in line $(< \xbo *s).$
\paragraph{
Regularities
}

 \xEh

 \xDH

The group of rules $(x*s)$ use ascending strength of $ \xdi / \xdf.$

 \xDH

The column $(\xdm^{+})$ contains interesting algebraic properties. In
particular,
they show a strengthening from $(3*s)$ up to Rationality. They are not
necessarily
equivalent to the corresponding $(I_{x})$ rules, not even in the presence
of the basic rules. The examples show that care has to be taken when
considering
the different variants.

 \xDH

Adding the somewhat superflous $(CM_{2}),$ we have increasing cautious
monotony from $ \xCf (wCM)$ to full $(CM_{ \xbo }).$

 \xDH

We have increasing ``or'' from $ \xCf (wOR)$ to full $(OR_{ \xbo }).$

 \xDH

The line $(2*s)$ is only there because there seems to be no $(
\xdm^{+}_{2}),$ otherwise
we could begin $(n*s)$ at $n=2.$

 \xEj

$ \xCO $
\paragraph{
Summary \\[2mm]
}

We can obtain all rules except $ \xCf (RatM)$ and $(\xCd)$ from $ \xCf
(Opt),$ the monotony
rules - $ \xCf (iM),$ $(eM \xdi),$ $(eM \xdf)$ -, and $(x*s)$ with
increasing $x.$
\subsubsection{
Table
}

$ \xCO $

The following table is split in two, as it is too big for printing in
one page.

$ \xCO $

\label{Definition Gen-Filter-2-Teile}
\index{Definition Size rules}

(See
Table \ref{Table Base2-Size-Rules-1} (page \pageref{Table Base2-Size-Rules-1}),
``Rules on size - Part $ \xfI $''
and
Table \ref{Table Base2-Size-Rules-2} (page \pageref{Table Base2-Size-Rules-2}),
``Rules on size - Part II''.

\begin{table}[h]

\index{Ideal}
\index{Filter}
\index{$ \xdm^+ $}
\index{$ \xeA $}
\index{Optimal proportion}
\index{$(Opt)$}
\index{Monotony}
\index{Improving proportions}
\index{$(iM)$}
\index{internal monotony}
\index{$(eM \xdi)$}
\index{external monotony for ideals}
\index{$(eM \xdf)$}
\index{external monotony for filters}
\index{$(iM)$}
\index{$(eM \xdi)$}
\index{$(eM \xdf)$}
\index{Keeping proportions}
\index{$(\xCd)$}
\index{$(\xdi \xcv disj)$}
\index{$(\xdf \xcv disj)$}
\index{$(\xdm^+ \xcv disj)$}
\index{Robustness of proportions}
\index{$(1*s)$}
\index{$(\xdi_1)$}
\index{$(\xdf_1)$}
\index{$(\xeA_1)$}
\index{$(2*s)$}
\index{$(\xdi_2)$}
\index{$(\xdf_2)$}
\index{$(\xeA_2)$}
\index{$(n*s)$}
\index{$(\xdi_n)$}
\index{$(\xdf_n)$}
\index{$(\xdm^{+}_{n})$}
\index{$(\xeA_n)$}
\index{$(< \xbo*s)$}
\index{$(\xdi_\xbo)$}
\index{$(\xdf_\xbo)$}
\index{$(\xdm^{+}_{ \xbo })$}
\index{$(\xeA_{\xbo})$}
\index{Robustness of $\xdm^+$}
\index{$(\xdm^{++})$}

\caption{Rules on size - Part I}

\label{Table Base2-Size-Rules-1}

\tabcolsep=0.5pt

\begin{center}

{\tiny

\begin{tabular}{|c|c@{.}c|c|c|}

\hline

\multicolumn{5}{|c|}{\bf Rules on size - Part I}\\
\hline

\xEH
``Ideal''
\xEH
``Filter''
\xEH
$ \xdm^+ $
\xEH
$ \xeA $
\xEP

\hline
\hline

\multicolumn{5}{|c|}{Optimal proportion} \xEP

\hline

$(Opt)$
\xEH
$ \xCQ \xbe \xdi (X)$
\xEH
$X \xbe \xdf (X)$
\xEH
\xEH
$ \xcA x \xba \xcp \xeA x \xba$
\xEP

\xEH
\xEH
\xEH
\xEH
\xEP

\hline
\hline

\multicolumn{5}{|c|}
{Monotony (Improving proportions). $(iM)$: internal monotony,}
\xEP
\multicolumn{5}{|c|}
{$(eM \xdi)$: external monotony for ideals,
$(eM \xdf)$: external monotony for filters}
\xEP

\hline

$(iM)$
\xEH
$A \xcc B \xbe \xdi (X)$
\xEH
$A \xbe \xdf (X)$,
\xEH
\xEH
$\xeA x \xba \xcu \xcA x (\xba \xcp \xba')$
\xEP

\xEH
$ \xch $
\xEH
$A \xcc B \xcc X$
\xEH
\xEH
$ \xcp $ $ \xeA x \xba'$
\xEP

\xEH
$A \xbe \xdi (X)$
\xEH
$ \xch $ $B \xbe \xdf (X)$
\xEH
\xEH
\xEP

\hline

$(eM \xdi)$
\xEH
$X \xcc Y \xch$
\xEH
\xEH
\xEH
$\xeA x (\xba: \xbb) \xcu$
\xEP

\xEH
$\xdi (X) \xcc \xdi (Y)$
\xEH
\xEH
\xEH
$\xcA x (\xba' \xcp \xbb) \xcp$
\xEP

\xEH
\xEH
\xEH
\xEH
$\xeA x (\xba \xco \xba': \xbb)$
\xEP

\xEH
\xEH
\xEH
\xEH
\xEP

\xEH
\xEH
\xEH
\xEH
\xEP

\xEH
\xEH
\xEH
\xEH
\xEP

\xEH
\xEH
\xEH
\xEH
\xEP

\hline

$(eM \xdf)$
\xEH
\xEH
$X \xcc Y \xch$
\xEH
\xEH
$\xeA x (\xba: \xbb) \xcu$
\xEP

\xEH
\xEH
$\xdf (Y) \xcs \xdp (X) \xcc $
\xEH
\xEH
$\xcA x (\xbb \xcu \xba \xcp \xba') \xcp$
\xEP

\xEH
\xEH
$ \xdf (X)$
\xEH
\xEH
$\xeA x (\xba \xcu \xba': \xbb)$
\xEP

\xEH
\xEH
\xEH
\xEH
\xEP

\hline
\hline

\multicolumn{5}{|c|}{Keeping proportions} \xEP

\hline

$(\xCd)$
\xEH
$(\xdi \xcv disj)$
\xEH
$(\xdf \xcv disj)$
\xEH
$(\xdm^+ \xcv disj)$
\xEH
$\xeA x(\xba: \xbb) \xcu$
\xEP

\xEH
$A \xbe \xdi (X),$
\xEH
$A \xbe \xdf (X),$
\xEH
$A \xbe \xdm^+ (X),$
\xEH
$\xeA x(\xba': \xbb) \xcu$
\xEP

\xEH
$B \xbe \xdi (Y),$
\xEH
$B \xbe \xdf (Y),$
\xEH
$B \xbe \xdm^+ (Y),$
\xEH
$\xCN \xcE x(\xba \xcu \xba') \xcp$
\xEP

\xEH
$X \xcs Y= \xCQ $ $ \xch $
\xEH
$X \xcs Y= \xCQ $ $ \xch $
\xEH
$X \xcs Y= \xCQ $ $ \xch $
\xEH
$\xeA x(\xba \xco \xba': \xbb)$
\xEP

\xEH
$A \xcv B \xbe \xdi (X \xcv Y)$
\xEH
$A \xcv B \xbe \xdf (X \xcv Y)$
\xEH
$A \xcv B \xbe \xdm^+ (X \xcv Y)$
\xEH
\xEP

\xEH
\xEH
\xEH
\xEH
\xEP

\xEH
\xEH
\xEH
\xEH
\xEP

\hline
\hline

\multicolumn{5}{|c|}{Robustness of proportions: $n*small \xEd All$} \xEP

\hline

$(1*s)$
\xEH
$(\xdi_1)$
\xEH
$(\xdf_1)$
\xEH
\xEH
$(\xeA_1)$
\xEP

\xEH
$X \xce \xdi (X)$
\xEH
$ \xCQ \xce \xdf (X)$
\xEH
\xEH
$ \xeA x \xba \xcp \xcE x \xba $
\xEP

\hline

$(2*s)$
\xEH
$(\xdi_2)$
\xEH
$(\xdf_2)$
\xEH
\xEH
$(\xeA_2)$
\xEP

\xEH
$A,B \xbe \xdi (X) \xch $
\xEH
$A,B \xbe \xdf (X) \xch $
\xEH
\xEH
$ \xeA x \xba \xcu \xeA x \xbb $
\xEP

\xEH
$ A \xcv B \xEd X$
\xEH
$A \xcs B \xEd \xCQ $
\xEH
\xEH
$ \xcp $ $ \xcE x(\xba \xcu \xbb)$
\xEP

\hline

$(n*s)$
\xEH
$(\xdi_n)$
\xEH
$(\xdf_n)$
\xEH
$(\xdm^{+}_{n})$
\xEH
$(\xeA_n)$
\xEP

$(n \xcg 3)$
\xEH
$A_{1},.,A_{n} \xbe \xdi (X) $
\xEH
$A_{1},.,A_{n} \xbe \xdi (X) $
\xEH
$X_{1} \xbe \xdf (X_{2}),., $
\xEH
$ \xeA x \xba_{1} \xcu.  \xcu \xeA x \xba_{n} $
\xEP

\xEH
$ \xch $
\xEH
$ \xch $
\xEH
$ X_{n-1} \xbe \xdf (X_{n})$ $ \xch $
\xEH
$ \xcp $
\xEP

\xEH
$ A_{1} \xcv.  \xcv A_{n} \xEd X $
\xEH
$A_{1} \xcs.  \xcs A_{n} \xEd \xCQ$
\xEH
$X_{1} \xbe \xdm^{+}(X_{n})$
\xEH
$ \xcE x (\xba_{1} \xcu.  \xcu \xba_{n}) $
\xEP

\xEH
\xEH
\xEH
\xEH
\xEP

\hline

$(< \xbo*s)$
\xEH
$(\xdi_\xbo)$
\xEH
$(\xdf_\xbo)$
\xEH
$(\xdm^{+}_{ \xbo })$
\xEH
$(\xeA_{\xbo})$
\xEP

\xEH
$A,B \xbe \xdi (X) \xch $
\xEH
$A,B \xbe \xdf (X) \xch $
\xEH
(1)
\xEH
$ \xeA x \xba \xcu \xeA x \xbb \xcp $
\xEP

\xEH
$ A \xcv B \xbe \xdi (X)$
\xEH
$ A \xcs B \xbe \xdf (X)$
\xEH
$A \xbe \xdf (X),$ $X \xbe \xdm^{+}(Y)$
\xEH
$ \xeA x(\xba \xcu \xbb)$
\xEP

\xEH
\xEH
\xEH
$ \xch $ $A \xbe \xdm^{+}(Y)$
\xEH
\xEP

\xEH
\xEH
\xEH
(2)
\xEH
\xEP

\xEH
\xEH
\xEH
$A \xbe \xdm^{+}(X),$ $X \xbe \xdf (Y)$
\xEH
\xEP

\xEH
\xEH
\xEH
$ \xch $ $A \xbe \xdm^{+}(Y)$
\xEH
\xEP

\xEH
\xEH
\xEH
(3)
\xEH
\xEP

\xEH
\xEH
\xEH
$A \xbe \xdf (X),$ $X \xbe \xdf (Y)$
\xEH
\xEP

\xEH
\xEH
\xEH
$ \xch $ $A \xbe \xdf (Y)$
\xEH
\xEP

\xEH
\xEH
\xEH
(4)
\xEH
\xEP

\xEH
\xEH
\xEH
$A,B \xbe \xdi (X)$ $ \xch $
\xEH
\xEP

\xEH
\xEH
\xEH
$A-B \xbe \xdi (X-$B)
\xEH
\xEP

\hline
\hline

\multicolumn{5}{|c|}{Robustness of $\xdm^+$} \xEP

\hline

$(\xdm^{++})$
\xEH
\xEH
\xEH
$(\xdm^{++})$
\xEH
\xEP

\xEH
\xEH
\xEH
(1)
\xEH
\xEP

\xEH
\xEH
\xEH
$A \xbe \xdi (X),$ $B \xce \xdf (X)$
\xEH
\xEP

\xEH
\xEH
\xEH
$ \xch $ $A-B \xbe \xdi (X-B)$
\xEH
\xEP

\xEH
\xEH
\xEH
(2)
\xEH
\xEP

\xEH
\xEH
\xEH
$A \xbe \xdf (X), B \xce \xdf (X)$
\xEH
\xEP

\xEH
\xEH
\xEH
$ \xch $ $A-B \xbe \xdf (X-B)$
\xEH
\xEP

\xEH
\xEH
\xEH
(3)
\xEH
\xEP

\xEH
\xEH
\xEH
$A \xbe \xdm^+ (X),$
\xEH
\xEP

\xEH
\xEH
\xEH
$X \xbe \xdm^+ (Y)$
\xEH
\xEP

\xEH
\xEH
\xEH
$ \xch $ $A \xbe \xdm^+ (Y)$
\xEH
\xEP

\hline

\end{tabular}
}
\end{center}
\end{table}

\newpage

\begin{table}[h]

\index{AND}
\index{OR}
\index{Cautious Monotony}
\index{Rational Monotony}
\index{Optimal proportion}
\index{$(Opt)$}
\index{$(SC)$}
\index{Monotony}
\index{Improving proportions}
\index{$(iM)$}
\index{$(RW)$}
\index{$(eM \xdi)$}
\index{$(PR')$}
\index{$(wOR)$}
\index{$(\xbm wOR)$}
\index{$(\xbm PR)$}
\index{$(eM \xdf)$}
\index{$(wCM)$}
\index{Keeping proportions}
\index{$(\xCd)$}
\index{$(NR)$}
\index{$(disjOR)$}
\index{$(\xbm disjOR)$}
\index{Robustness of proportions}
\index{$(1*s)$}
\index{$(CP)$}
\index{$(AND_{1})$}
\index{$(2*s)$}
\index{$(AND_{2})$}
\index{$(OR_{2})$}
\index{$(CM_{2})$}
\index{$(n*s)$}
\index{$(AND_{n})$}
\index{$(OR_{n})$}
\index{$(CM_{n})$}
\index{$(< \xbo*s)$}
\index{$(AND_{ \xbo })$}
\index{$(OR_{ \xbo })$}
\index{$(CM_{ \xbo })$}
\index{$(\xbm OR)$}
\index{$(\xbm CM)$}
\index{Robustness of $\xdm^+$}
\index{$(\xdm^{++})$}
\index{$(RatM)$}
\index{$(\xbm RatM)$}

\caption{Rules on size - Part II}

\label{Table Base2-Size-Rules-2}

\begin{center}

{\tiny

\begin{tabular}{|c|c|c|c|c|}

\hline

\multicolumn{5}{|c|}{\bf Rules on size - Part II}\\
\hline

\xEH
various rules
\xEH
AND
\xEH
OR
\xEH
Caut./Rat.Mon.
\xEP

\hline
\hline

\multicolumn{5}{|c|}{Optimal proportion} \xEP

\hline

$(Opt)$
\xEH
$(SC)$
\xEH
\xEH
\xEH
\xEP

\xEH
$ \xba \xcl \xbb \xch \xba \xcn \xbb $
\xEH
\xEH
\xEH
\xEP

\hline
\hline

\multicolumn{5}{|c|}
{Monotony (Improving proportions)}
\xEP
\multicolumn{5}{|c|}
{}
\xEP

\hline

$(iM)$
\xEH
$(RW)$
\xEH
\xEH
\xEH
\xEP

\xEH
$ \xba \xcn \xbb, \xbb \xcl \xbb' \xch $
\xEH
\xEH
\xEH
\xEP

\xEH
$ \xba \xcn \xbb' $
\xEH
\xEH
\xEH
\xEP

\hline

$(eM \xdi)$
\xEH
$(PR')$
\xEH
\xEH
$(wOR)$
\xEH
\xEP

\xEH
$\xba \xcn \xbb, \xba \xcl \xba',$
\xEH
\xEH
$ \xba \xcn \xbb,$ $ \xba' \xcl \xbb $ $ \xch $
\xEH
\xEP

\xEH
$\xba' \xcu \xCN \xba \xcl \xbb \xch$
\xEH
\xEH
$ \xba \xco \xba' \xcn \xbb $
\xEH
\xEP

\xEH
$\xba' \xcn \xbb$
\xEH
\xEH
$(\xbm wOR)$
\xEH
\xEP

\xEH
$(\xbm PR)$
\xEH
\xEH
$\xbm(X \xcv Y) \xcc \xbm(X) \xcv Y$
\xEH
\xEP

\xEH
$X \xcc Y \xch$
\xEH
\xEH
\xEH
\xEP

\xEH
$\xbm(Y) \xcs X \xcc \xbm(X)$
\xEH
\xEH
\xEH
\xEP

\hline

$(eM \xdf)$
\xEH
\xEH
\xEH
\xEH
$(wCM)$
\xEP

\xEH
\xEH
\xEH
\xEH
$\xba \xcn \xbb, \xba' \xcl \xba,$
\xEP

\xEH
\xEH
\xEH
\xEH
$\xba \xcu \xbb \xcl \xba' \xch$
\xEP

\xEH
\xEH
\xEH
\xEH
$\xba' \xcn \xbb$
\xEP

\hline
\hline

\multicolumn{5}{|c|}{Keeping proportions} \xEP

\hline

$(\xCd)$
\xEH
$(NR)$
\xEH
\xEH
$(disjOR)$
\xEH
\xEP

\xEH
$\xba \xcn \xbb \xch$
\xEH
\xEH
$ \xba \xcn \xbb,$ $ \xba' \xcn \xbb' $
\xEH
\xEP

\xEH
$\xba \xcu \xbg \xcn \xbb$
\xEH
\xEH
$ \xba \xcl \xCN \xba',$ $ \xch $
\xEH
\xEP

\xEH
or
\xEH
\xEH
$ \xba \xco \xba' \xcn \xbb \xco \xbb' $
\xEH
\xEP

\xEH
$\xba \xcu \xCN \xbg \xcn \xbb$
\xEH
\xEH
$(\xbm disjOR)$
\xEH
\xEP

\xEH
\xEH
\xEH
$X \xcs Y = \xCQ \xch$
\xEH
\xEP

\xEH
\xEH
\xEH
$\xbm(X \xcv Y) \xcc \xbm(X) \xcv \xbm(Y)$
\xEH
\xEP

\hline
\hline

\multicolumn{5}{|c|}{Robustness of proportions: $n*small \xEd All$} \xEP

\hline

$(1*s)$
\xEH
$(CP)$
\xEH
$(AND_{1})$
\xEH
\xEH
\xEP

\xEH
$\xba\xcn\xcT \xch \xba\xcl\xcT$
\xEH
$ \xba \xcn \xbb $ $ \xch $ $ \xba \xcL \xCN \xbb $
\xEH
\xEH
\xEP

\hline

$(2*s)$
\xEH
\xEH
$(AND_{2})$
\xEH
$(OR_{2})$
\xEH
$(CM_{2})$
\xEP

\xEH
\xEH
$ \xba \xcn \xbb,$ $ \xba \xcn \xbb' $ $ \xch $
\xEH
$ \xba \xcn \xbb \xch \xba \xcN \xCN \xbb $
\xEH
$ \xba \xcn \xbb \xch \xba \xcN \xCN \xbb $
\xEP

\xEH
\xEH
$ \xba \xcL \xCN \xbb \xco \xCN \xbb' $
\xEH
\xEH
\xEP

\hline

$(n*s)$
\xEH
\xEH
$(AND_{n})$
\xEH
$(OR_{n})$
\xEH
$(CM_{n})$
\xEP

$(n \xcg 3)$
\xEH
\xEH
$ \xba \xcn \xbb_{1},., \xba \xcn \xbb_{n}$
\xEH
$ \xba_{1} \xcn \xbb,., \xba_{n-1} \xcn \xbb $
\xEH
$ \xba \xcn \xbb_{1},., \xba \xcn \xbb_{n-1}$
\xEP

\xEH
\xEH
$ \xch $
\xEH
$ \xch $
\xEH
$ \xch $
\xEP

\xEH
\xEH
$ \xba \xcL \xCN \xbb_{1} \xco.  \xco \xCN \xbb_{n}$
\xEH
$ \xba_{1} \xco.  \xco \xba_{n-1} \xcN \xCN \xbb $
\xEH
$ \xba \xcu \xbb_1 \xcu.  \xcu \xbb_{n-2} \xcN $
\xEP

\xEH
\xEH
\xEH
\xEH
$ \xCN \xbb_{n-1}$
\xEP

\hline

$(< \xbo*s)$
\xEH
\xEH
$(AND_{ \xbo })$
\xEH
$(OR_{ \xbo })$
\xEH
$(CM_{ \xbo })$
\xEP

\xEH
\xEH
$ \xba \xcn \xbb,$ $ \xba \xcn \xbb' $ $ \xch $
\xEH
$ \xba \xcn \xbb,$ $ \xba' \xcn \xbb $ $ \xch $
\xEH
$ \xba \xcn \xbb,$ $ \xba \xcn \xbb' $ $ \xch $
\xEP

\xEH
\xEH
$ \xba \xcn \xbb \xcu \xbb' $
\xEH
$ \xba \xco \xba' \xcn \xbb $
\xEH
$ \xba \xcu \xbb \xcn \xbb' $
\xEP

\xEH
\xEH
\xEH
$(\xbm OR)$
\xEH
$(\xbm CM)$
\xEP

\xEH
\xEH
\xEH
$\xbm(X \xcv Y) \xcc \xbm(X) \xcv \xbm(Y)$
\xEH
$\xbm(X) \xcc Y \xcc X \xch$
\xEP

\xEH
\xEH
\xEH
\xEH
$\xbm(Y) \xcc \xbm(X)$
\xEP

\xEH
\xEH
\xEH
\xEH
\xEP

\xEH
\xEH
\xEH
\xEH
\xEP

\xEH
\xEH
\xEH
\xEH
\xEP

\xEH
\xEH
\xEH
\xEH
\xEP

\xEH
\xEH
\xEH
\xEH
\xEP

\xEH
\xEH
\xEH
\xEH
\xEP

\xEH
\xEH
\xEH
\xEH
\xEP

\hline
\hline

\multicolumn{5}{|c|}{Robustness of $\xdm^+$} \xEP

\hline

$(\xdm^{++})$
\xEH
\xEH
\xEH
\xEH
$(RatM)$
\xEP

\xEH
\xEH
\xEH
\xEH
$ \xba \xcn \xbb,  \xba \xcN \xCN \xbb'   \xch $
\xEP

\xEH
\xEH
\xEH
\xEH
$ \xba \xcu \xbb' \xcn \xbb $
\xEP

\xEH
\xEH
\xEH
\xEH
$(\xbm RatM)$
\xEP

\xEH
\xEH
\xEH
\xEH
$X \xcc Y,$
\xEP

\xEH
\xEH
\xEH
\xEH
$X \xcs \xbm(Y) \xEd \xCQ \xch$
\xEP

\xEH
\xEH
\xEH
\xEH
$\xbm(X) \xcc \xbm(Y) \xcs X$
\xEP

\xEH
\xEH
\xEH
\xEH
\xEP

\xEH
\xEH
\xEH
\xEH
\xEP

\xEH
\xEH
\xEH
\xEH
\xEP

\xEH
\xEH
\xEH
\xEH
\xEP

\hline

\end{tabular}
}
\end{center}
\end{table}

$ \xCO $

$ \xCO $
\subsection{
Multiplicative laws about size
}

We are mainly interested in nonmonotonic logic. In this domain,
independence is strongly connected to multiplication of
abstract size, and an important part of the present text treats this
connection
and its repercussions.

We have at least two scenarios for multiplication, one is decribed in
Diagram \ref{Diagram Mul-Add} (page \pageref{Diagram Mul-Add}), the second in
Diagram \ref{Diagram Mul-Prod} (page \pageref{Diagram Mul-Prod}).
In the first scenario, we have nested sets, in the second, we have set
products.
In the first scenario, we consider subsets which behave as the big set
does,
in the second scenario we consider subspaces, and decompose the behaviour
of the big space into behaviour of the subspaces.
In both cases, this results naturally in multiplication of abstract sizes.
When we look at the corresponding relation properties, they are quite
different (rankedness vs. some kind of modularity). But this is perhaps
to be expected, as the two scenarios are quite different.

Other scenarios which might be interesting to consider in our framework
are:

 \xEI

 \xDH

When we have more than two truth values, say 3, and 2 is considered a big
subset, and we have $n$ propositional variables, and $m$ of them are
considered
many, then $2^{m}$ might give a
``big'' subset of the total of $3^{n}situations.$

 \xDH

Similarly, when we fix 1 variable, consider 2 cases of the possible 3, and
multiply this with a ``big'' set of models.

 \xDH

We may also consider the utility or cost of a situation, and work with
a ``big'' utility, and ``many'' situations, etc.

 \xDH

Note that, in the case of distances, subspaces add distances, and do not
multiply them: $d(xy,x' y')=d(x,x')+d(y,y').$

 \xEJ

These questions are left for further research,
see also Section \ref{Section Mod-Lang-Struc} (page \pageref{Section
Mod-Lang-Struc}).
\subsubsection{
Multiplication of size for subsets
}

Here we have nested sets, $A \xcc X \xcc Y,$ $ \xCf A$ is a certain
proportion of $X,$ and $X$ of $Y,$
resulting in a multiplication of relative size or proportions. This is a
classical subject of nonmonotonic logic,
see the last section, taken from  \cite{GS09a}, it is partly
repeated here
to stress the common points with the other scenario.

$ \xCO $

\vspace{10mm}

\begin{diagram}

\label{Diagram Mul-Add}
\index{Diagram Mul-Add}

\centering
\setlength{\unitlength}{1mm}
{\renewcommand{\dashlinestretch}{30}
\begin{picture}(150,100)(0,0)

\path(20,90)(90,90)(90,30)(20,30)(20,90)

\path(55,90)(55,30)
\path(55,60)(90,60)

\path(55,20)(90,20)
\path(55,22)(55,18)
\path(90,22)(90,18)

\path(20,10)(90,10)
\path(20,12)(20,8)
\path(90,12)(90,8)

\path(100,60)(100,30)
\path(98,60)(102,60)
\path(98,30)(102,30)

\put(54,6){{\xssc $Y$}}
\put(72,16){{\xssc $X$}}
\put(102,44){{\xssc $A$}}

\put(50,1){{\xssc Scenario 1}}

\end{picture}
}

\end{diagram}

\vspace{4mm}

$ \xCO $
\paragraph{
Properties \\[2mm]
}

Diagram \ref{Diagram Mul-Add} (page \pageref{Diagram Mul-Add})  is to be read as
follows:
The whole set $Y$ is split in $X$ and $Y- \xCf X,$ $X$ is split in $ \xCf
A$ and $X- \xCf A.$
$X$ is a small/medium/big part of $Y,$
$ \xCf A$ is a small/medium/big part of $X.$
The question is: is $ \xCf A$ a small/medium/big part of $Y?$

Note that the relation of $ \xCf A$ to $X$ is conceptually different from
that of
$X$ to $ \xCf Y,$ as we change the base set by going from $X$ to $ \xCf
Y,$ but not when going
from $ \xCf A$ to $X.$ Thus, in particular, when we read the diagram as
expressing
multiplication, commutativity is not necessarily true.

We looked at this scenario already in  \cite{GS09a}, but there from
an
additive point of view, using various basic properties like $ \xCf (iM),$
$(eM \xdi),$
$(eM \xdf),$ see Section \ref{Section Add-Size} (page \pageref{Section
Add-Size}).
Here, we use just multiplication - except sometimes for motivation.

We examine different rules:

If $Y=X$ or $X=A,$ there is nothing to show, so 1 is the neutral element
of
multiplication.

If $X \xbe \xdi (Y)$ or $A \xbe \xdi (X),$ then we should have $A \xbe
\xdi (Y).$ (Use for motivation
$ \xCf (iM)$ or $(eM \xdi)$ respectively.)

So it remains to look at the following cases, with the ``natural'' answers
given already:

(1) $X \xbe \xdf (Y),$ $A \xbe \xdf (X)$ $ \xch $ $A \xbe \xdf (Y),$

(2) $X \xbe \xdm^{+}(Y),$ $A \xbe \xdf (X)$ $ \xch $ $A \xbe \xdm^{+}(Y),$

(3) $X \xbe \xdf (Y),$ $A \xbe \xdm^{+}(X)$ $ \xch $ $A \xbe \xdm^{+}(Y),$

(4) $X \xbe \xdm^{+}(Y),$ $A \xbe \xdm^{+}(X)$ $ \xch $ $A \xbe
\xdm^{+}(Y).$

But (1) is case (3) of $(\xdm^{+}_{ \xbo })$ in  \cite{GS09a}, see
Table ``Rules on size'' in
Section \ref{Section Add-Size} (page \pageref{Section Add-Size}).

(2) is case (1) of $(\xdm^{+}_{ \xbo })$ there,

(3) is case (2) of $(\xdm^{+}_{ \xbo })$ there, finally,

(4) is $(\xdm^{++})$ there.

So the first three correspond to various expressions of $(AND_{ \xbo }),$
$(OR_{ \xbo }),$
$(CM_{ \xbo }),$ the last one to $ \xCf (RatM).$

But we can read them also the other way round, e.g.:

(1) corresponds to: $ \xba \xcn \xbb,$ $ \xba \xcu \xbb \xcn \xbg $ $
\xch $ $ \xba \xcn \xbg,$

(2) corresponds to: $ \xba \xcN \xCN \xbb,$ $ \xba \xcu \xbb \xcn \xbg $
$ \xch $ $ \xba \xcN \xCN (\xbb \xcu \xbg),$

(3) corresponds to: $ \xba \xcn \xbb,$ $ \xba \xcu \xbb \xcN \xCN \xbg $
$ \xch $ $ \xba \xcN \xCN (\xbb \xcu \xbg).$

All these rules might be seen as too idealistic, so just as we did in
 \cite{GS09a}, we can consider milder versions:
We might for instance consider a rule which says that $big* \Xl *big,$ $n$
times,
is not small. Consider for instance the case $n=2.$
So we would conclude that $ \xCf A$ is not small in $Y.$ In terms of
logic, we
then have: $ \xba \xcn \xbb,$ $ \xba \xcu \xbb \xcn \xbg $ $ \xch $ $
\xba \xcN (\xCN \xbb \xco \xCN \xbg).$ We can obtain the same
logical property from $3*small \xEd all.$
\subsubsection{
Multiplication of size for subspaces
}

Our main interest here is multiplication for subspaces, which we discuss
now.
\paragraph{
Properties \\[2mm]
}

\label{Section Mul-Mul}

$ \xCO $

\vspace{10mm}

\begin{diagram}

\label{Diagram Mul-Prod}
\index{Diagram Mul-Prod}

\centering
\setlength{\unitlength}{1mm}
{\renewcommand{\dashlinestretch}{30}
\begin{picture}(150,100)(0,0)

\path(20,90)(90,90)(90,30)(20,30)(20,90)

\path(55,90)(55,30)

\path(55,70)(90,70)
\path(20,75)(55,75)

\path(20,20)(90,20)

\path(20,22)(20,18)
\path(55,22)(55,18)
\path(90,22)(90,18)

\path(10,75)(10,30)
\path(8,75)(12,75)
\path(8,30)(12,30)

\path(100,70)(100,30)
\path(98,70)(102,70)
\path(98,30)(102,30)

\put(37,16){{\xssc $\xbS_1$}}
\put(72,16){{\xssc $\xbS_2$}}

\put(5,50){{\xssc $\xbG_1$}}
\put(101,50){{\xssc $\xbG_2$}}

\put(50,1){{\xssc Scenario 2}}

\end{picture}
}

\end{diagram}

\vspace{4mm}

$ \xCO $

In this scenario, $ \xbS_{i}$ are sets of sequences,
see Diagram \ref{Diagram Mul-Prod} (page \pageref{Diagram Mul-Prod}),
corresponding, intuitively, to a set of models in language $ \xdl_{i},$
$ \xbS_{i}$ will be the set of $ \xba_{i}-$models, and the subsets $
\xbG_{i}$ are to be seen as
the ``best'' models, where $ \xbb_{i}$ will hold. The languages are supposed
to
be disjoint sublanguages of a common language $ \xdl.$
As the $ \xbS_{i}$ have symmetrical roles,
there is no intuitive reason for multiplication not to be commutative.

We can interpret the situation twofold:

First, we work separately in sublanguage $ \xdl_{1}$ and $ \xdl_{2},$ and,
say, $ \xba_{i}$ and $ \xbb_{i}$
are both defined in $ \xdl_{i},$ and we look at $ \xba_{i} \xcn \xbb_{i}$
in the sublanguage $ \xdl_{i},$
or, we consider both $ \xba_{i}$ and $ \xbb_{i}$ in the big language $
\xdl,$ and look at
$ \xba_{i} \xcn \xbb_{i}$ in $ \xdl.$ These two ways are a priori
completely different.
Speaking in preferential terms, it is not at all clear why the orderings
on the submodels should have anything to do with the orderings on the
whole models. It seems a very desirable property, but we have to postulate
it, which we do now (an overview is given
in Table \ref{Table Mul-Laws} (page \pageref{Table Mul-Laws})). We first give
informally a list of
such rules,
mainly to show the connection with the first scenario. Later,
see Definition \ref{Definition Sin-Size-Rules} (page \pageref{Definition
Sin-Size-Rules}), we will
introduce formally some rules for which we show a connection with
interpolation.
Here, e.g.,
``$(big*big \xch big)$'' stands for
``if both factors are big, so will be the product'',
this will be abbreviated by
``$b*b \xch b$'' in Table \ref{Table Mul-Laws} (page \pageref{Table Mul-Laws})
.

$(big*1 \xch big)$ Let $ \xbG_{1} \xcc \xbS_{1},$ if $ \xbG_{1} \xbe \xdf
(\xbS_{1}),$ then $ \xbG_{1} \xCK \xbS_{2} \xbe \xdf (\xbS_{1} \xCK
\xbS_{2}),$
(and the dual rule for $ \xbS_{2}$ and $ \xbG_{2}).$

This property preserves proportions, so it seems intuitively quite
uncontested, whenever we admit coherence over products. (Recall that there
was nothing to show in the first scenario.)

When we re-consider above case: suppose $ \xba \xcn \xbb $ in the
sublanguage, so
$M(\xbb) \xbe \xdf (M(\xba))$ in the sublanguage, so by $(big*1 \xch
big),$ $M(\xbb) \xbe \xdf (M(\xba))$
in the big language $ \xdl.$

We obtain the dual rule for small (and likewise, medium size) sets:

$(small*1 \xch small)$ Let $ \xbG_{1} \xcc \xbS_{1},$ if $ \xbG_{1} \xbe
\xdi (\xbS_{1}),$ then $ \xbG_{1} \xCK \xbS_{2} \xbe \xdi (\xbS_{1} \xCK
\xbS_{2}),$
(and the dual rule for $ \xbS_{2}$ and $ \xbG_{2}),$

establishing $ \xCf All=1$ as the neutral element for multiplication.

We look now at other, plausible rules:

$(small*x \xch small)$ $ \xbG_{1} \xbe \xdi (\xbS_{1}),$ $ \xbG_{2} \xcc
\xbS_{2}$ $ \xch $ $ \xbG_{1} \xCK \xbG_{2} \xbe \xdi (\xbS_{1} \xCK
\xbS_{2})$

$(big*big \xch big)$ $ \xbG_{1} \xbe \xdf (\xbS_{1}),$ $ \xbG_{2} \xbe
\xdf (\xbS_{2})$ $ \xch $ $ \xbG_{1} \xCK \xbG_{2} \xbe \xdf (\xbS_{1}
\xCK \xbS_{2})$

$(big*medium \xch medium)$ $ \xbG_{1} \xbe \xdf (\xbS_{1}),$ $ \xbG_{2}
\xbe \xdm^{+}(\xbS_{2})$ $ \xch $ $ \xbG_{1} \xCK \xbG_{2} \xbe \xdm^{+}(
\xbS_{1} \xCK \xbS_{2})$

$(medium*medium \xch medium)$ $ \xbG_{1} \xbe \xdm^{+}(\xbS_{1}),$ $
\xbG_{2} \xbe \xdm^{+}(\xbS_{2})$ $ \xch $ $ \xbG_{1} \xCK \xbG_{2} \xbe
\xdm^{+}(\xbS_{1} \xCK \xbS_{2})$

When we accept all above rules, we can invert $(big*big \xch big),$ as a
big product
must be composed of big components. Likewise, at least one component of a
small
product has to be small - see
Proposition \ref{Proposition Sin-Product-Small} (page \pageref{Proposition
Sin-Product-Small}).

We see that these properties give a lot of modularity. We can calculate
the consequences of $ \xba $ and $ \xba' $ separately - provided $ \xba
,$ $ \xba' $ use disjoint
alphabets - and put the results together afterwards. Such properties are
particularly interesting for classification purposes, where subclasses
are defined with disjoint alphabets.

Recall that
we work here with a notion of ``big'' and
``small'' subsets, which may be thought of as defined by a filter (ideal),
though we usually will not need the full strength of a filter (ideal). But
assume as usual that $A \xcc B \xcc C$ and $A \xcc C$ is big imply $B \xcc
C$ is big, that $C \xcc C$ is
big, and define
$A \xcc B$ is small iff $(B-A) \xcc B$ is big, call all subsets which are
neither big
nor small medium size. For an extensive discussion, see
 \cite{GS09a}.

Let $X' \xcv X'' =X$ be a disjoint cover, so $ \xbP X= \xbP X' \xCK \xbP
X''.$ We consider subsets $ \xbS $
etc. of $ \xbP X.$ If not said otherwise, $ \xbS $ etc. need not be a
product $ \xbS' \xCK \xbS''.$
We will sometimes write $ \xbP':= \xbP X',$ $ \xbS'':= \xbS \xex X''
.$ The roles of $X' $ and $X'' $ are
interchangeable, e.g., instead of $ \xbG \xex X' \xcc \xbS \xex X',$ we
may also write
$ \xbG \xex X'' \xcc \xbS \xex X''.$

We consider here the following two sets of three finite product rules
about size
and $ \xbm.$ Both sets will be shown to be equivalent in
Proposition \ref{Proposition Sin-Product-Small} (page \pageref{Proposition
Sin-Product-Small}).

\bd

$\hspace{0.01em}$


\label{Definition Sin-Size-Rules}

$(S*1)$ $ \xbD \xcc \xbS' \xCK \xbS'' $ is big iff there is $ \xbG =
\xbG' \xCK \xbG'' \xcc \xbD $ such that $ \xbG' \xcc \xbS' $ and
$ \xbG'' \xcc \xbS'' $ are big

$(S*2)$ $ \xbG \xcc \xbS $ is big $ \xch $ $ \xbG \xex X' \xcc \xbS \xex
X' $ is big - where
$ \xbS $ is not necessarily a product.

$(S*3)$ $A \xcc \xbS $ is big $ \xch $ there is $B \xcc \xbP' \xCK \xbS
'' $ big such that $B \xex X'' \xcc A \xex X'' $ - again,
$ \xbS $ is not necessarily a product.

$(\xbm *1)$ $ \xbm (\xbS' \xCK \xbS'')= \xbm (\xbS') \xCK \xbm (
\xbS'')$

$(\xbm *2)$ $ \xbm (\xbS) \xcc \xbG $ $ \xch $ $ \xbm (\xbS \xex X')
\xcc \xbG \xex X' $

$(\xbm *3)$ $ \xbm (\xbP X' \xCK \xbS'') \xex X'' \xcc \xbm (\xbS)
\xex X'' $

$(s*s)$ Let $ \xbG_{i} \xcc \xbS_{i},$ then $ \xbG_{1} \xCK \xbG_{2} \xcc
\xbS_{1} \xCK \xbS_{2}$ is small iff $ \xbG_{1} \xcc \xbS_{1}$ is small or
$ \xbG_{1} \xcc \xbS_{1}$ is small.

A generalization to more than two factors is obvious.

One can also consider weakenings, e.g.,

$(S*1')$ $ \xbG' \xCK \xbS'' \xcc \xbS' \xCK \xbS'' $ is big iff $
\xbG' \xcc \xbS' $ is big.

\ed

\bp

$\hspace{0.01em}$


\label{Proposition Sin-Product-Small}

 \xEh
 \xDH Let $(S*1)$ hold. Then:

$ \xbG' \xCK \xbG'' \xcc \xbS' \xCK \xbS'' $ is small iff $ \xbG'
\xcc \xbS' $ or $ \xbG'' \xcc \xbS'' $ is small.

 \xDH If the filters over $ \xCf A$ are principal filters, generated by $
\xbm (A),$
i.e. $B \xcc A$ is big iff $ \xbm (A) \xcc B \xcc A$ for some $ \xbm (A)
\xcc A,$ then:

$(S*i)$ is equivalent to $(\xbm *i),$ $i=1,2,3.$

 \xDH Let the notion of size satisfy $ \xCf (Opt),$ $ \xCf (iM),$ and $(<
\xbo *s),$ see
the tables ``Rules on size'' in
Section \ref{Section Add-Size} (page \pageref{Section Add-Size}).
Then $(\xbm *1)$ and $(s*s)$ are equivalent.
 \xEj

\ep

\subparagraph{
Proof
}

$\hspace{0.01em}$


(1)

``$ \xci $'':

Suppose $ \xbG' \xcc \xbS' $ is small. Then $ \xbS' - \xbG' \xcc \xbS
' $ is big and
$(\xbS' - \xbG') \xCK \xbS'' \xcc \xbS' \xCK \xbS'' $ is big by
$(S*1).$ But
$(\xbG' \xCK \xbG'') \xcs ((\xbS' - \xbG') \xCK \xbS'')= \xCQ,$
so $ \xbG' \xCK \xbG'' \xcc \xbS' \xCK \xbS'' $ is small.

``$ \xch $'':

For the converse, suppose that neither $ \xbG' \xcc \xbS' $ nor $ \xbG
'' \xcc \xbS'' $ are small. Let
$A \xcc \xbS' \xCK \xbS'' $ be big, we show that $A \xcs (\xbG' \xCK
\xbG'') \xEd \xCQ.$
By $(S*1)$ there are $B' \xcc \xbS' $ and $B'' \xcc \xbS'' $ big, and
$B' \xCK B'' \xcc A.$ Then $B' \xcs \xbG' \xEd \xCQ,$
$B'' \xcs \xbG'' \xEd \xCQ,$ so there is $ \xBc x',x''  \xBe  \xbe (B' \xCK B''
)
\xcs (\xbG' \xCK \xbG'') \xcc A \xcs (\xbG' \xCK \xbG'').$

(2.1)

``$ \xch $''

``$ \xcc $'': $ \xbm (\xbS') \xcc \xbS' $ and $ \xbm (\xbS'') \xcc
\xbS'' $ are big, so by $(S*1)$
$ \xbm (\xbS') \xCK \xbm (\xbS'') \xcc \xbS' \xCK \xbS'' $ is big,
so $ \xbm (\xbS' \xCK \xbS'') \xcc \xbm (\xbS') \xCK \xbm (\xbS''
).$

``$ \xcd $'': $ \xbm (\xbS' \xCK \xbS'') \xcc \xbS' \xCK \xbS'' $
is big $ \xch $ by $(S*1)$ there is $ \xbG' \xCK \xbG'' \xcc \xbm (\xbS
' \xCK \xbS'')$
and $ \xbG' \xcc \xbS',$ $ \xbG'' \xcc \xbS'' $ big $ \xch $ $ \xbm (
\xbS') \xcc \xbG',$ $ \xbm (\xbS'') \xcc \xbG'' $ $ \xch $
$ \xbm (\xbS') \xCK \xbm (\xbS'') \xcc \xbm (\xbS' \xCK \xbS''
).$

``$ \xci $''

Let $ \xbG' \xcc \xbS' $ be big, $ \xbG'' \xcc \xbS'' $ be big, $ \xbG
' \xCK \xbG'' \xcc \xbD,$ then $ \xbm (\xbS') \xcc \xbG',$
$ \xbm (\xbS'') \xcc \xbG'',$ so
by $(\xbm *1)$ $ \xbm (\xbS)= \xbm (\xbS') \xCK \xbm (\xbS'')
\xcc \xbG' \xCK \xbG'' \xcc \xbD,$ so $ \xbD $ is big.

Let $ \xbD \xcc \xbS $ be big, then by $(\xbm *1)$ $ \xbm (\xbS') \xCK
\xbm (\xbS'')= \xbm (\xbS) \xcc \xbD.$

(2.2)

``$ \xch $''

$ \xbm (\xbS) \xcc \xbG $ $ \xch $ $ \xbG \xcc \xbS $ big $ \xch $ by
$(S*2)$ $ \xbG \xex X' \xcc \xbS \xex X' $ big $ \xch $ $ \xbm (\xbS \xex
X') \xcc \xbG \xex X'.$

``$ \xci $''

$ \xbG \xcc \xbS $ big $ \xch $ $ \xbm (\xbS) \xcc \xbG $ $ \xch $ by $(
\xbm *2)$ $ \xbm (\xbS \xex X') \xcc \xbG \xex X' $ $ \xch $ $ \xbG \xex
X' \xcc \xbS \xex X' $ big.

(2.3)

``$ \xch $''

$ \xbm (\xbS) \xcc \xbS $ big $ \xch $ $ \xcE B \xcc \xbP X' \xCK \xbS
'' $ big such that $B \xex X'' \xcc \xbm (\xbS) \xex X'' $ by $(S*3),$
thus in
particular $ \xbm (\xbP X' \xCK \xbS'') \xex X'' \xcc \xbm (\xbS)
\xex X''.$

``$ \xci $''

$A \xcc \xbS $ big $ \xch $ $ \xbm (\xbS) \xcc A.$ $ \xbm (\xbP X' \xCK
\xbS'') \xcc \xbP X' \xCK \xbS'' $ is big, and by $(\xbm *3)$
$ \xbm (\xbP X' \xCK \xbS'') \xex X'' \xcc \xbm (\xbS) \xex X'' \xcc
A \xex X''.$

(3)

``$ \xch $'':

(1) Let $ \xbG' \xcc \xbS' $ be small, we show that $ \xbG' \xCK \xbG
'' \xcc \xbS' \xCK \xbS'' $ is small.
So $ \xbS' - \xbG' \xcc \xbS' $ is big, so by $ \xCf (Opt)$ and $(\xbm
*1)$
$(\xbS' - \xbG') \xCK \xbS'' \xcc \xbS' \xCK \xbS'' $ is big, so
$ \xbG' \xCK \xbS'' $ $=$ $(\xbS' \xCK \xbS'')-((\xbS' - \xbG')
\xCK \xbS'')$ $ \xcc $ $ \xbS' \xCK \xbS'' $ is small,
so by $ \xCf (iM)$ $ \xbG' \xCK \xbG'' \xcc \xbS' \xCK \xbS'' $ is
small.

(2) Suppose $ \xbG' \xcc \xbS' $ and $ \xbG'' \xcc \xbS'' $ are not
small, we show that
$ \xbG' \xCK \xbG'' \xcc \xbS' \xCK \xbS'' $
is not small. So $ \xbS' - \xbG' \xcc \xbS' $ and $ \xbS'' - \xbG''
\xcc \xbS'' $ are not big.
We show that $Z$ $:=$ $((\xbS' \xCK \xbS'')-(\xbG' \xCK \xbG''))$
$ \xcc $ $ \xbS' \xCK \xbS'' $ is not big.
$Z$ $=$ $(\xbS' \xCK (\xbS'' - \xbG'')) \xcv ((\xbS' - \xbG')
\xCK \xbS'').$

Suppose $X' \xCK X'' $ $ \xcc $ $Z,$ then $X' \xcc \xbS' - \xbG' $ or
$X'' \xcc \xbS'' - \xbG''.$ Proof:
Let $X' \xcC \xbS' - \xbG' $ and $X'' \xcC \xbS'' - \xbG'',$ but $X'
\xCK X'' \xcc Z.$ Let $ \xbs' \xbe X' -(\xbS' - \xbG'),$
$ \xbs'' \xbe X'' -(\xbS'' - \xbG''),$ consider $ \xbs' \xbs''.$
$ \xbs' \xbs'' \xce (\xbS' - \xbG') \xCK \xbS'',$ as $ \xbs' \xce
\xbS' - \xbG',$
$ \xbs' \xbs'' \xce \xbS' \xCK (\xbS'' \xCK \xbG''),$ as $ \xbs''
\xce \xbS'' - \xbG'',$ so $ \xbs' \xbs'' \xce Z.$

By prerequisite, $ \xbS' - \xbG' \xcc \xbS' $ is not big, $ \xbS'' -
\xbG'' \xcc \xbS'' $ is not big, so by
$ \xCf (iM)$ no $X' $ with $X' \xcc \xbS' - \xbG' $ is big, no $X'' $
with $X'' \xcc \xbS'' - \xbG'' $ is big, so
by $(\xbm *1)$ or $(S*1)$ $Z \xcc \xbS' \xCK \xbS'' $ is not big, so $
\xbG' \xCK \xbG'' \xcc \xbS' \xCK \xbS'' $
is not small.

``$ \xci $'':

(1) Suppose $ \xbG' \xcc \xbS' $ is big, $ \xbG'' \xcc \xbS'' $ is
big, we have to show
$ \xbG' \xCK \xbG'' \xcc \xbS' \xCK \xbS'' $ is big.
$ \xbS' - \xbG' \xcc \xbS' $ is small, $ \xbS'' - \xbG'' \xcc \xbS''
$ is small, so by $(s*s)$
$(\xbS' - \xbG') \xCK \xbS'' \xcc \xbS' \xCK \xbS'' $ is small and
$ \xbS' \xCK (\xbS'' - \xbG'') \xcc \xbS' \xCK \xbS'' $ is small,
so by $(< \xbo *s)$
$(\xbS' \xCK \xbS'')-(\xbG' \xCK \xbG'')$ $=$
$((\xbS' - \xbG') \xCK \xbS'') \xcv (\xbS' \xCK (\xbS'' - \xbG
''))$ $ \xcc $ $ \xbS' \xCK \xbS'' $ is small,
so $ \xbG' \xCK \xbG'' \xcc \xbS' \xCK \xbS'' $ is big.

(2) Suppose $ \xbG' \xCK \xbG'' \xcc \xbS' \xCK \xbS'' $ is big, we
have to show $ \xbG' \xcc \xbS' $ is big, and
$ \xbG'' \xcc \xbS'' $ is big. By prerequisite,
$(\xbS' \xCK \xbS'')-(\xbG' \xCK \xbG'')$ $=$
$((\xbS' - \xbG') \xCK \xbS'') \xcv (\xbS' \xCK (\xbS'' - \xbG
''))$ $ \xcc $ $ \xbS' \xCK \xbS'' $ is small,
so by $ \xCf (iM)$ $ \xbS' \xCK (\xbS'' - \xbG'') \xcc \xbS' \xCK
\xbS'' $ is small, so by $ \xCf (Opt)$ and $(s*s)$
$ \xbS'' - \xbG'' \xcc \xbS'' $ is small, so $ \xbG'' \xcc \xbS'' $
is big, and likewise $ \xbG' \xcc \xbS' $
is big.

$ \xcz $
\\[3ex]

\paragraph{
Discussion
}

$\hspace{0.01em}$


\label{Section Discussion}

We compare these rules to probability defined size.

Let ``big'' be defined by
``more than $50 \xET $''. If $ \xbP X' $ and $ \xbP X'' $ have 3 elements
each, then subsets of $ \xbP X' $
or $ \xbP X'' $ of $card \xcg 2$
are big. But taking the product may give $4/9<1/2.$ So the product rule
``$big*big=big$'' will not hold there. One direction will hold, of course.

Next, we discuss the prerequisite $ \xbS = \xbS' \xCK \xbS''.$ Consider
the following
example:

\be

$\hspace{0.01em}$


\label{Example Sin-Prod-Size}

Take a language of 5 propositional variables, with $X':=\{a,b,c\},$ $X''
:=\{d,e\}.$
Consider the model set $ \xbS:=\{ \xCL a \xCL b \xCL cde,$ $-a-b-c-d \xCL
e\},$ i.e. of 8 models of
$ \xCf de$ and 2 models of $- \xCf d.$ The models of $ \xCf de$ are 8/10
of all elements of $ \xbS,$ so
it is reasonable to call them a big subset of $ \xbS.$ But its projection
on
$X'' $ is only 1/3 of $ \xbS''.$

So we have a potential $ \xCf decrease$ when going to the coordinates.

This shows that weakening the prerequisite about $ \xbS $ as done in
$(S*2)$ is not
innocent.

\ee

\br

$\hspace{0.01em}$


\label{Remark Sin-Boolean}

When we set small sets to 0, big sets to 1, we have the following boolean
rules
for filters:

(1) $0+0=0$

(2) $1+x=1$

(3) $-0=1,$ $-1=0$

(4) $0*x=0$

(5) $1*1=1$

There are no such rules for medium size sets, as the union of two medium
size sets may be big, but also stay medium.

Such multiplication rules capture the behaviour of Reiter defaults and of
defeasible inheritance.
\subsection{
Hamming relations and distances
}

\label{Section Ham-Rel-Dist}
\subsubsection{
Hamming relations and multiplication of size
}

\er

We now define Hamming relations in various flavours, and then
(see Proposition \ref{Proposition Mod-Hamming} (page \pageref{Proposition
Mod-Hamming}))
show that (smooth) Hamming relations
generate a notion of size which satisfies our conditions, defined in
Definition \ref{Definition Sin-Size-Rules} (page \pageref{Definition
Sin-Size-Rules}).
Corollary \ref{Corollary Sin-Interpolation-2} (page \pageref{Corollary
Sin-Interpolation-2})  will put our results
together, and show that (smooth) Hamming relations generate
preferential logics with interpolation.

We will conclude this section by showing that our conditions
$(\xbm *1)$ and $(\xbm *2)$ essentially characterise Hamming relations.

Note that we re-define Hamming relations in
Section \ref{Section Ham-Neigh} (page \pageref{Section Ham-Neigh}), as already
announced in
Section \ref{Section Ham-Rel-Size} (page \pageref{Section Ham-Rel-Size}).

\bd

$\hspace{0.01em}$


\label{Definition Sin-Set-HR}

We abuse notation, and define a relation $ \xcc $ on $ \xbP X',$ $ \xbP
X'',$ and $ \xbP X' \xCK \xbP X''.$

(1) Define $x \xeb y: \xcj x \xec y$ and $x \xEd y,$ thus $ \xbs \xec \xbt
$ iff $ \xbs \xeb \xbt $ or $ \xbs = \xbt.$

(2) We say that a relation $ \xec $ satisfies (GH3) iff

(GH3) $ \xbs' \xbs'' \xec \xbt' \xbt'' $ $ \xcj $ $ \xbs' \xec \xbt'
$ and $ \xbs'' \xec \xbt'',$

(Thus, $ \xbs' \xbs'' \xeb \xbt' \xbt'' $ iff $ \xbs' \xbs'' \xec
\xbt' \xbt'' $ and $(\xbs' \xeb \xbt' $ or $ \xbs'' \xeb \xbt''
).)$

(3) Call a relation $ \xeb $ a $ \xCf GH$ $(=$ general Hamming) relation
iff the following two
conditions hold:

$ \xCf (GH1)$ $ \xbs' \xec \xbt' $ $ \xcu $ $ \xbs'' \xec \xbt'' $ $
\xcu $ $(\xbs' \xeb \xbt' $ $ \xco $ $ \xbs'' \xeb \xbt'')$ $ \xch $
$ \xbs' \xbs'' \xeb \xbt' \xbt'' $

(where $ \xbs' \xec \xbt' $ iff $ \xbs' \xeb \xbt' $ or $ \xbs' =
\xbt')$

$ \xCf (GH2)$ $ \xbs' \xbs'' \xeb \xbt' \xbt'' $ $ \xch $ $ \xbs'
\xeb \xbt' $ $ \xco $ $ \xbs'' \xeb \xbt'' $

$ \xCf (GH2)$ means that some compensation is possible, e.g., $ \xbt'
\xeb \xbs' $ might be the
case, but $ \xbs'' \xeb \xbt'' $ wins in the end, so $ \xbs' \xbs''
\xeb \xbt' \xbt''.$

We use $ \xCf (GH)$ for $(GH1)+(GH2).$

\ed

\be

$\hspace{0.01em}$


\label{Example Sin-Int-Circum}

The circumscription relation satisfies (GH3) with $ \xCN p \xck p$ and
$ \xcU \xCL q_{i} \xck \xcU \xCL q'_{i}$ iff $ \xcA i(\xCL q_{i} \xck
\xCL q'_{i}).$

\ee

\br

$\hspace{0.01em}$


\label{Remark Sin-Independence}

(1) The independence makes sense because the concept of models, and thus
the
usual
interpolation for classical logic relie on the independence of the
assignments.

(2) This corresponds to social choice for many independent dimensions.

(3) We can also consider such factorisation as an approximation:
we can do part of the reasoning independently.

\er

\bd

$\hspace{0.01em}$


\label{Definition Sin-Relation-to-Filter}

Given a relation $ \xec,$ define as usual a principal filter $ \xdf (X)$
generated by the
$ \xec -$minimal elements:

$ \xbm (X):=\{x \xbe X: \xCN \xcE x' \xeb x.x' \xbe X\},$

$ \xdf (X):=\{A \xcc X: \xbm (X) \xcc A\}.$

\ed

The following proposition summarizes various properties for the different
Hamming relations:

\bp

$\hspace{0.01em}$


\label{Proposition Mod-Hamming}

Let $ \xbP X= \xbP X' \xCK \xbP X'',$ $ \xbS \xcc \xbP X.$

(1)
Let $ \xec $ be a smooth relation satisfying (GH3). Then
$(\xbm *2)$ holds, and thus $(S*2)$ by Proposition \ref{Proposition
Sin-Product-Small} (page \pageref{Proposition Sin-Product-Small}), (2).

(2) Let again $ \xbS':= \xbS \xex X',$ $ \xbS'':= \xbS \xex X''.$
Let $ \xec $ be a smooth relation satisfying (GH3). Then:

$ \xbm (\xbS') \xCK \xbm (\xbS'') \xcc \xbS $ $ \xch $ $ \xbm (\xbS
)= \xbm (\xbS') \xCK \xbm (\xbS'').$

(Here $ \xbS = \xbS' \xCK \xbS'' $ will not necessarily hold.)

(3) Let again $ \xbS':= \xbS \xex X',$ $ \xbS'':= \xbS \xex X''.$
Let $ \xec $ be a relation satisfying (GH3), and $ \xbS = \xbS' \xCK \xbS
''.$ Then
$(\xbm *1)$ holds, and thus, by Proposition \ref{Proposition Sin-Product-Small}
(page \pageref{Proposition Sin-Product-Small}), (2) $(S*1).$

(4)
Let $ \xec $ be a smooth relation satisfying (GH3), then
$(\xbm *3)$ holds, and thus by Proposition \ref{Proposition Sin-Product-Small}
(page \pageref{Proposition Sin-Product-Small})  (2) $(S*3).$

(5)

$(\xbm *1)$ and $(\xbm *2)$ and the usual axioms for smooth relations
characterize smooth relations satisfying (GH3).

(6)

Let $ \xbs \xeb \xbt \xcj \xbt \xce \xbm (\{ \xbs, \xbt \})$ and $ \xeb $
be smooth.
Then $ \xbm $ satisfies $(\xbm *1)$ (or, by
Proposition \ref{Proposition Sin-Product-Small} (page \pageref{Proposition
Sin-Product-Small})
equivalently $(s*s))$
iff $ \xeb $ is a $ \xCf GH$ relation.

(7) Let $ \xbG' \xcc \xbS',$ $ \xbG'' \xcc \xbS'',$ $ \xbG' \xCK
\xbG'' \xcc \xbS' \xCK \xbS'' $ be small, let $ \xCf (GH2)$ hold, then
$ \xbG' \xcc \xbS' $ is small or $ \xbG'' \xcc \xbS'' $ is small.

(8) Let $ \xbG' \xcc \xbS' $ be small, $ \xbG'' \xcc \xbS'',$ let $
\xCf (GH1)$ hold, then
$ \xbG' \xCK \xbG'' \xcc \xbS' \xCK \xbS'' $ is small.

\ep

\subparagraph{
Proof
}

$\hspace{0.01em}$


(1)
Suppose $ \xbm (\xbS) \xcc \xbG $ and $ \xbs' \xbe \xbS \xex X' - \xbG
\xex X',$ we show $ \xbs' \xce \xbm (\xbS \xex X').$

Let $ \xbs = \xbs' \xbs'' \xbe \xbS,$ then $ \xbs \xce \xbG,$ so $
\xbs \xce \xbm (\xbS).$ So here is $ \xbr \xeb \xbs,$
$ \xbr \xbe \xbm (\xbS) \xcc \xbG $ by
smoothness. Let $ \xbr = \xbr' \xbr''.$
We have $ \xbr' \xec \xbs' $ by (GH3). $ \xbr' = \xbs' $ cannot be, as
$ \xbr' \xbe \xbG \xex X',$ and
$ \xbs' \xce \xbG \xex X'.$ So $ \xbr' \xeb \xbs',$ and $ \xbs' \xce
\xbm (\xbS \xex X').$

(2)

``$ \xcd $'': Let $ \xbs' \xbe \xbm (\xbS'),$ $ \xbs'' \xbe \xbm (
\xbS'').$ By prerequisite, $ \xbs' \xbs'' \xbe \xbS.$ Suppose
$ \xbt \xeb \xbs' \xbs'',$ then $ \xbt' \xeb \xbs' $ or $ \xbt''
\xeb \xbs'',$ contradiction.

``$ \xcc $'': Let $ \xbs \xbe \xbm (\xbS),$ suppose $ \xbs' \xce \xbm
(\xbS')$ or $ \xbs'' \xce \xbm (\xbS'').$
So there are $ \xbt' \xec \xbs',$ $ \xbt'' \xec \xbs'' $ with $ \xbt
' \xbe \xbm (\xbS'),$ $ \xbt'' \xbe \xbm (\xbS'')$ by
smoothness. Moreoever, $ \xbt' \xeb \xbs' $ or $ \xbt'' \xeb \xbs''.$
By prerequisite $ \xbt' \xbt'' \xbe \xbS,$
and $ \xbt' \xbt'' \xeb \xbs,$ so $ \xbs \xce \xbm (\xbS).$

(3)

``$ \xcd $'': As in (2), the prerequisite holds trivially.

``$ \xcc $'': As in (2), but we do not need $ \xbt' \xbe \xbm (\xbS'
),$ $ \xbt'' \xbe \xbm (\xbS''),$ as $ \xbt' \xbt'' $ will
be in $ \xbS $ trivially. So smoothness is not needed.

(4)

Let again $ \xbS'' = \xbS \xex X''.$

Let $ \xbD:= \xbP X' \xCK \xbS'',$ $ \xbs = \xbs' \xbs'' \xbe \xbm (
\xbD).$
Suppose $ \xbs'' \xce \xbm (\xbS) \xex X''.$ There cannot be any $
\xbt \xeb \xbs,$ $ \xbt \xbe \xbS,$ by $ \xbS \xcc \xbD.$
So $ \xbs \xce \xbS,$ but $ \xbs'' \xbe \xbS'',$ so there is $ \xbt
\xbe \xbS $ $ \xbt'' = \xbs''.$ As $ \xbt $ is not minimal,
there must be minimal $ \xbr = \xbr' \xbr'' \xeb \xbt,$ $ \xbr \xbe
\xbS $ by smoothness. As $ \xbr $ is minimal,
$ \xbr'' \xEd \xbs'',$ and as
$ \xbr \xeb \xbt,$ $ \xbr'' \xeb \xbs'' $ by (GH3). By prerequisite $
\xbs' \xbr'' \xbe \xbD,$ and $ \xbs' \xbr'' \xeb \xbs,$
contradiction.

Note that smoothness is essential. Otherwise, there might be an infinite
descending chain $ \xbt_{i}$ below $ \xbt,$ all with $ \xbt''_{i}= \xbs
'',$ but none below $ \xbs.$

(5)

If $ \xec $ is smooth and satisfies (GH3), then $(\xbm *1)$ and $(\xbm
*2)$ hold by
(1) and (3). For the converse:

Define as usual $ \xbs \xeb \xbt: \xcj \xbt \xce \xbm (\{ \xbs, \xbt
\}).$ Let $ \xbs = \xbs' \xbs'',$ $ \xbt = \xbt' \xbt''.$

We have to show:

$ \xbs \xeb \xbt $ iff $ \xbs' \xec \xbt' $ and $ \xbs'' \xec \xbt'' $
and $(\xbs' \xeb \xbt' $ or $ \xbs'' \xeb \xbt'').$

``$ \xci $'':

Suppose $ \xbs' \xeb \xbt' $ and $ \xbs'' \xec \xbt''.$ Then $ \xbm
(\{ \xbs', \xbt' \})=\{ \xbs' \},$ and
$ \xbm (\{ \xbs'', \xbt'' \})=\{ \xbs'' \}$ (either $ \xbs'' \xeb
\xbt'' $ or $ \xbs'' = \xbt'',$ so in both cases
$ \xbm (\{ \xbs'', \xbt'' \})=\{ \xbs'' \}).$ As $ \xbt' \xce \xbm
(\{ \xbs', \xbt' \}),$
$ \xbt $ $ \xce $ $ \xbm (\{ \xbs', \xbt' \} \xCK \{ \xbs'', \xbt''
\})$ $=$ (by $(\xbm *1))$ $ \xbm (\{ \xbs', \xbt' \}) \xCK \xbm (\{
\xbs'', \xbt'' \})$
$=$ $\{ \xbs' \} \xCK \{ \xbs'' \}=\{ \xbs \},$
so by smoothness $ \xbs \xeb \xbt.$

``$ \xch $'':

Conversely, if $ \xbs \xeb \xbt,$ so $ \xbG:=\{ \xbs \}= \xbm (\xbS)$
for $ \xbS:=\{ \xbs, \xbt \},$ so by $(\xbm *2)$
$ \xbm (\xbS \xex X')= \xbm (\{ \xbs', \xbt' \}) \xcc \xbG \xex X'
=\{ \xbs' \},$ so $ \xbs' \xec \xbt',$
analogously $ \xbm (\xbS \xex X'')= \xbm (\{ \xbs'', \xbt'' \}) \xcc
\xbG \xex X'' =\{ \xbs'' \},$ so $ \xbs'' \xec \xbt'',$ but both
cannot be equal.

(6)

(6.1) $(\xbm *1)$ entails the $ \xCf GH$ relation conditions

$ \xCf (GH1):$
Suppose $ \xbs' \xeb \xbt' $ and $ \xbs'' \xec \xbt''.$ Then $ \xbt'
\xce \xbm (\{ \xbs', \xbt' \})=\{ \xbs' \},$ and
$ \xbm (\{ \xbs'', \xbt'' \})=\{ \xbs'' \}$ (either $ \xbs'' \xeb
\xbt'' $ or $ \xbs'' = \xbt'',$ so in both cases
$ \xbm (\{ \xbs'', \xbt'' \})=\{ \xbs'' \}).$ As $ \xbt' \xce \xbm
(\{ \xbs', \xbt' \}),$
$ \xbt' \xbt'' \xce \xbm (\{ \xbs', \xbt' \} \xCK \{ \xbs'', \xbt
'' \})$ $=_{(\xbm *1)}$ $ \xbm (\{ \xbs', \xbt' \}) \xCK \xbm (\{ \xbs
'', \xbt'' \})$
$=\{ \xbs' \} \xCK \{ \xbs'' \}=\{ \xbs' \xbs'' \},$ so by smoothness
$ \xbs' \xbs'' \xeb \xbt' \xbt''.$

$ \xCf (GH2):$
Let $X:=\{ \xbs', \xbt' \},$ $Y:=\{ \xbs'', \xbt'' \},$ so $X \xCK
Y=\{ \xbs' \xbs'', \xbs' \xbt'', \xbt' \xbs'', \xbt' \xbt''
\}.$
Suppose $ \xbs' \xbs'' \xeb \xbt' \xbt'',$ so $ \xbt' \xbt'' \xce
\xbm (X \xCK Y)$ $=_{(\xbm *1)}$ $ \xbm (X) \xCK \xbm (Y).$ If $ \xbs'
\xeB \xbt',$ then
$ \xbt' \xbe \xbm (X),$ likewise if $ \xbs'' \xeB \xbt'',$ then $ \xbt
'' \xbe \xbm (Y),$ so $ \xbt' \xbt'' \xbe \xbm (X \xCK Y),$
contradiction.

(6.2) The $ \xCf GH$ relation conditions generate $(\xbm *1).$

$ \xbm (X \xCK Y) \xcc \xbm (X) \xCK \xbm (Y):$
Let $ \xbt' \xbe X,$ $ \xbt'' \xbe Y,$ $ \xbt' \xbt'' \xce \xbm (X)
\xCK \xbm (Y),$ then $ \xbt' \xce \xbm (X)$ or $ \xbt'' \xce \xbm (Y).$
Suppose
$ \xbt' \xce \xbm (X),$ let $ \xbs' \xbe X,$ $ \xbs' \xeb \xbt',$ so
by condition $ \xCf (GH1)$ $ \xbs' \xbt'' \xeb \xbt' \xbt'',$ so
$ \xbt' \xbt'' \xce \xbm (X \xCK Y).$

$ \xbm (X) \xCK \xbm (Y) \xcc \xbm (X \xCK Y):$
Let $ \xbt' \xbe X,$ $ \xbt'' \xbe Y,$ $ \xbt' \xbt'' \xce \xbm (X
\xCK Y),$ so there is $ \xbs' \xbs'' \xeb \xbt' \xbt'',$ $ \xbs'
\xbe X,$ $ \xbs'' \xbe Y,$
so by $ \xCf (GH2)$ either $ \xbs' \xeb \xbt' $ or $ \xbs'' \xeb \xbt
'',$ so $ \xbt' \xce \xbm (X)$ or $ \xbt'' \xce \xbm (Y),$ so
$ \xbt' \xbt'' \xce \xbm (X) \xCK \xbm (Y).$

(7)

Suppose $ \xbG' \xcc \xbS' $ is not small, so there is $ \xbg' \xbe
\xbG' $ and no $ \xbs' \xbe \xbS' $ with $ \xbs' \xeb \xbg'.$
Fix this $ \xbg'.$
Consider $\{ \xbg' \} \xCK \xbG''.$ As $ \xbG' \xCK \xbG'' \xcc \xbS
' \xCK \xbS'' $ is small, there is for each
$ \xbg' \xbg'',$ $ \xbg'' \xbe \xbG'' $ some $ \xbs' \xbs'' \xbe
\xbS' \xCK \xbS'',$ $ \xbs' \xbs'' \xeb \xbg' \xbg''.$ By $ \xCf
(GH2)$
$ \xbs' \xeb \xbg' $ or $ \xbs'' \xeb \xbg'',$
but $ \xbs' \xeb \xbg' $ was excluded, so for all $ \xbg'' \xbe \xbG''
$ there is $ \xbs'' \xbe \xbS'' $ with
$ \xbs'' \xeb \xbg'',$ so $ \xbG'' \xcc \xbS'' $ is small.

(8)

Let $ \xbg' \xbe \xbG',$ so there is $ \xbs' \xbe \xbS' $ and $ \xbs
' \xeb \xbg'.$ By $ \xCf (GH1),$ for any
$ \xbg'' \xbe \xbG'' $ $ \xbs' \xbg'' \xeb \xbg' \xbg'',$ so no $
\xbg' \xbg'' \xbe \xbG' \xCK \xbG'' $ is minimal.

$ \xcz $
\\[3ex]

\be

$\hspace{0.01em}$


\label{Example Sin-Component-Inverse}

Even for smooth relations satisfying (GH3), the converse of $(\xbm *2)$
is not
necessarily true:

Let $ \xbs' \xeb \xbt',$ $ \xbt'' \xeb \xbs'',$ $ \xbS:=\{ \xbs,
\xbt \},$ then $ \xbm (\xbS)= \xbS,$ but $ \xbm (\xbS')=\{ \xbs'
\},$
$ \xbm (\xbS'')=\{ \xbt'' \},$ so $ \xbm (\xbS) \xEd \xbm (\xbS')
\xCK \xbm (\xbS'').$

We need the additional assumption that $ \xbm (\xbS') \xCK \xbm (\xbS
'') \xcc \xbS,$ see
Proposition \ref{Proposition Mod-Hamming} (page \pageref{Proposition
Mod-Hamming})  (2).

\ee

\be

$\hspace{0.01em}$


\label{Example Mul-GH-Rel}

The following are examples of $ \xCf GH$ relations:

Define on all components $X_{i}$ a relation $ \xeb_{i}.$

(1) The set variant Hamming relation:

Let the relation $ \xeb $ be defined on $ \xbP \{X_{i}:i \xbe I\}$ by $
\xbs \xeb \xbt $ iff for all $j$ $ \xbs_{j} \xec_{j} \xbt_{j},$
and there is at least one $i$ such that $ \xbs_{i} \xeb_{i} \xbt_{i}.$

(2) The counting variant Hamming relation:

Let the relation $ \xeb $ be defined on $ \xbP \{X_{i}:i \xbe I\}$ by $
\xbs \xeb \xbt $ iff the number of $i$
such that $ \xbs_{i} \xeb_{i} \xbt_{i}$ is bigger than the number of $i$
such that $ \xbt_{i} \xeb_{i} \xbs_{i}.$

(3) The weighed counting Hamming relation:

Like the counting relation, but we give different (numerical) importance
to different $i.$ E.g., $ \xbs_{1} \xeb \xbt_{1}$ may count 1, $ \xbs_{2}
\xeb \xbt_{2}$ may count 2, etc.

$ \xcz $
\\[3ex]

\ee

\paragraph{
Note
}

$\hspace{0.01em}$


\label{Section Note}

Note that already $(\xbm *1)$ results in a strong independence result in
the second
scenario:
Let $ \xbs \xbr' \xeb \xbt \xbr',$ then $ \xbs \xbr'' \xeb \xbt \xbr
'' $ for all $ \xbr''.$ Thus, whether $\{ \xbr'' \}$ is small,
or medium size (i.e. $ \xbr'' \xbe \xbm (\xbS')),$ the behaviour of $
\xbS \xCK \{ \xbr'' \}$ is the same.
This we do not have in the first scenario, as small sets may behave very
differently from medium size sets. (But, still, their internal structure
is the same, only the minimal elements change.)
When $(\xbm *2)$ holds, then if $ \xbs \xbs' \xeb \xbt \xbt' $ and $
\xbs \xEd \xbt,$ then $ \xbs \xeb \xbt,$ i.e. we
need not have $ \xbs' = \xbt'.$
\subsubsection{
Hamming distances and revision
}

\label{Section Sin-Hamming-Dist}

This short Section is mainly intended to put our work in a broader
perspective, by showing a connection of Hamming distances to
modular revision as introduced by Parikh and his co-authors.
The main result here is
Corollary \ref{Corollary Sin-Parikh-Revision} (page \pageref{Corollary
Sin-Parikh-Revision}).
We will not go into details of motivation here, and refer the reader to,
e.g.,
 \cite{Par96} for further discussion.

Thus, we have modular distances and relations, i.e., Hamming
distances and relations, we have modular revision as described below,
and we have modular logic, which has the (semantic) interpolation
property.
We want to point out here in particular this cross reference from modular
revision to modular logic, i.e., logic with interpolation.

We recall:

\bd

$\hspace{0.01em}$


\label{Definition Sin-Bar}

Given a distance $d,$ define for two sets $X,Y$

$X \xfA Y:=\{y \xbe Y: \xcE x \xbe X(\xCN \xcE x' \xbe X,y' \xbe Y.d(x'
,y')<d(x,y))\}.$

We assume that $X \xfA Y \xEd \xCQ $ if $X,Y \xEd \xCQ.$ Note that this
is related to the
consistency axiom of AGM theory revision: revising by a consistent formula
gives a consistent result. The assumption may be wrong due to infinite
descending chains of distances.

\ed

\bd

$\hspace{0.01em}$


\label{Definition Sin-TR}

Given $ \xfA,$ we can define an AGM revision operator $*$ as follows:

$T* \xbf:=Th(M(T) \xfA M(\xbf))$

where $T$ is a theory, and $Th(X)$ is the set of formulas which hold in
all $x \xbe X.$

It was shown in  \cite{LMS01} that a revision operator thus defined
satisfies the AGM revision postulates.

\ed

\bd

$\hspace{0.01em}$


\label{Definition Mul-GHD}

Let $d$ be an abstract distance on some product space $X \xCK Y,$ and its
components.
(We require of distances only that they are comparable, that $d(x,y)=0$
iff
$x=y,$ and that $d(x,y) \xcg 0.)$

$d$ is called a generalized Hamming distance $ \xCf (GHD)$
iff it satisfies the following two properties:

$ \xCf (GHD1)$ $d(\xbs, \xbt) \xck d(\xba, \xbb)$ and $d(\xbs',
\xbt') \xck d(\xba', \xbb')$ and
$(d(\xbs, \xbt)<d(\xba, \xbb)$ or $d(\xbs', \xbt')<d(\xba',
\xbb'))$ $ \xch $ $d(\xbs \xbs', \xbt \xbt')<d(\xba \xba', \xbb
\xbb')$

$ \xCf (GHD2)$ $d(\xbs \xbs', \xbt \xbt')<d(\xba \xba', \xbb \xbb
')$ $ \xch $
$d(\xbs, \xbt)<d(\xba, \xbb)$ or $d(\xbs', \xbt')<d(\xba',
\xbb')$

(Compare this definition to
Definition \ref{Definition Sin-Set-HR} (page \pageref{Definition Sin-Set-HR}).)

\ed

We have a result analogous to the relation case:

\bfa

$\hspace{0.01em}$


\label{Fact Mul-GHD-Bar}

Let $ \xfA $ be defined by a generalized Hamming distance, then $ \xfA $
satisfies

(1)

$(\xfA *)$ $(\xbS_{1} \xCK \xbS_{1}') \xfA (\xbS_{2} \xCK \xbS_{2}'
)=(\xbS_{1} \xfA \xbS_{2}) \xCK (\xbS_{1}' \xfA \xbS_{2}').$

(2) $(\xbS_{1}' \xfA \xbS_{2}') \xCK (\xbS_{1}'' \xfA \xbS_{2}'') \xcc
\xbS_{2}$ and $(\xbS_{2}' \xfA \xbS_{1}') \xCK (\xbS_{2}'' \xfA
\xbS_{1}'') \xcc \xbS_{1}$ $ \xch $
$(\xbS_{1}) \xfA (\xbS_{2})=(\xbS_{1}' \xfA \xbS_{2}') \xCK (
\xbS_{1}'' \xfA \xbS_{2}''),$ if the distance is symmetric

(where $ \xbS_{i}$ here is not necessarily $ \xbS_{i}' \xCK \xbS_{i}'',$
etc.).

\efa

\subparagraph{
Proof
}

$\hspace{0.01em}$


(1) and (2).

``$ \xcc $'':

Suppose $d(\xbs \xbs', \xbt \xbt')$ is minimal. If there is
$ \xba \xbe \xbS_{1},$ $ \xbb \xbe \xbS_{2}$ such that $d(\xba, \xbb
)<d(\xbs, \xbt),$ then
$d(\xba \xbs', \xbb \xbt')<d(\xbs \xbs', \xbt \xbt')$ by $ \xCf
(GHD1),$ so $d(\xbs, \xbt)$ and $d(\xbs', \xbt')$ have to be
minimal.

``$ \xcd $'':

For the converse, suppose $d(\xbs, \xbt)$ and $d(\xbs', \xbt')$
are minimal, but
$d(\xbs \xbs', \xbt \xbt')$ is not, so $d(\xba \xba', \xbb \xbb'
)<d(\xbs \xbs', \xbt \xbt')$ for some $ \xba \xba', \xbb \xbb',$
then $d(\xba, \xbb)<d(\xbs, \xbt)$ or $d(\xba', \xbb')<d(\xbs
', \xbt')$ by $ \xCf (GHD2),$ contradiction.

$ \xcz $
\\[3ex]

These properties translate to logic as follows:

\bco

$\hspace{0.01em}$


\label{Corollary Mul-GHD-TR}

If $ \xbf $ and $ \xbq $ are defined on a separate language from that of $
\xbf' $ and $ \xbq',$
and the distance satisfies $ \xCf (GHD1)$ and $ \xCf (GHD2),$ then for
revision holds:

$(\xbf \xcu \xbf')*(\xbq \xcu \xbq')=(\xbf * \xbq) \xcu (\xbf' *
\xbq').$

\eco

\bco

$\hspace{0.01em}$


\label{Corollary Sin-Parikh-Revision}

By Corollary \ref{Corollary Mul-GHD-TR} (page \pageref{Corollary Mul-GHD-TR}),
Hamming distances
generate decomposable revision operators a la Parikh,
see  \cite{Par96},
also in the generalized form of variable $K$ and $ \xbf.$

\eco

We conclude with a small result on partial (semantical) revision:

\bfa

$\hspace{0.01em}$


\label{Fact Sin-HD-Part}

Let $ \xfA $ be defined by a Hamming distance, then:

$ \xbP X \xfA \xbS \xcc \xbG $ $ \xch $ $ \xbP X' \xfA (\xbS \xex X')
\xcc \xbG \xex X'.$

(Recall that $ \xbP' $ is the restriction of $ \xbP $ to $X'.)$

\efa

\subparagraph{
Proof
}

$\hspace{0.01em}$


Let $t \xbe \xbS \xex X' - \xbG \xex X',$ we show $t \xce \xbP X' \xfA (
\xbS \xex X').$ Let $ \xbt \xbe \xbS $ be such that
$ \xbt' =t,$ then $ \xbt \xce \xbG $ (otherwise $t \xbe \xbG \xex X'),$
so $ \xbt \xce \xbP X \xfA \xbS,$ so there is
$ \xba = \xba' \xba'' \xbe \xbP X,$
$ \xbb = \xbb' \xbb'' \xbe \xbS,$ with $d(\xba, \xbb)$ minimal, so
$d(\xba, \xbb)<d(\xbs, \xbt)$ for all $ \xbs \xbe \xbP X.$ If
$d(\xbs', \xbt')$ were minimal
for some $ \xbs,$ then we would consider $ \xbs' \xba'',$ $ \xbt'
\xbb'',$
then $d(\xba' \xba'', \xbb' \xbb'')<d(\xbs' \xba'', \xbt' \xbb
'')$ is impossible by (GHD2), so
$ \xbt' \xbb'' \xbe \xbP X \xfA \xbS,$ so $ \xbt' \xbb'' \xbe \xbG,$
and $t \xbe \xbG \xex X',$ contradiction.

$ \xcz $
\\[3ex]
\subsubsection{
Discussion of representation
}

It would be nice to have a representation result like the one for
Hamming relations, see
Proposition \ref{Proposition Mod-Hamming} (page \pageref{Proposition
Mod-Hamming}), (5).
But this is impossible, for the following reason:

In constructing the representing distance from revision results, we made
arbitrary choices (see the proofs in  \cite{LMS01} or
 \cite{Sch04}). I.e., we choose sometimes
arbitrarily $d(x,y) \xck d(x',y'),$ when we do not have enough
information
to decide. (This is an example of the fact that the problem of
``losing ignorance'' should not be underestimated, see e.g.
 \cite{GS08f}.)
As we do not follow the same procedure for all cases, there
is no guarantee that the different representations will fit together.

Of course, it might be possible to come to a uniform choice, and one
could then attempt a representation result. This is left as an open
problem.
\subsection{
Summary of properties
}

\label{Section Mul-Sum}

We summarize in this section properties related to multiplicative laws.

They are collected in
Table \ref{Table Mul-Laws} (page \pageref{Table Mul-Laws}).

$pr(b)=b$ means: the projection of a big set on one of its coordinates
is big again.

$ \xCO $

\label{Table-Mul-Laws}
\begin{table}
\caption{Multiplication laws}

\label{Table Mul-Laws}
\begin{center}
\tabcolsep=0.5pt.
\begin{turn}{90}
{\xssB
\begin{tabular}{|c|c|c|c|c|c|c|c|c|}
\hline
\multicolumn{9}{|c|}{

Multiplication laws
} \xEP
\hline
Multiplication \xEH
\multicolumn{3}{|c|}{

Scenario 1
} \xEH
\multicolumn{5}{|c|}{

Scenario 2 $(*$ symmetrical, only 1 side shown)
} \xEP

law \xEH
\multicolumn{3}{|c|}{

(see Diagram \ref{Diagram Mul-Add} (page \pageref{Diagram Mul-Add}))
} \xEH
\multicolumn{5}{|c|}{

(see Diagram \ref{Diagram Mul-Prod} (page \pageref{Diagram Mul-Prod}))
} \xEP
\cline{2-9}

 \xEH Corresponding algebraic \xEH Logical property \xEH Relation \xEH
Algebraic property \xEH
Logical property \xEH
\multicolumn{3}{|c|}{

Interpolation
} \xEP
\cline{7-9}

 \xEH addition property \xEH \xEH property \xEH $(\xbG_{i} \xcc
\xbS_{i})$ \xEH
$(\xba, \xbb $ in $L_{1},$ $ \xba', \xbb' $ in $L_{2}$ \xEH
Multiplic. \xEH Relation \xEH Inter- \xEP
 \xEH \xEH \xEH \xEH \xEH
$L=L_{1} \xcv L_{2}$ (disjoint)) \xEH law \xEH property \xEH polation \xEP
\hline
\multicolumn{9}{|c|}{

Non-monotonic logic
} \xEP
\hline

$x*1 \xch x$ \xEH
\multicolumn{3}{|c|}{trivial}

 \xEH
$ \xbG_{1} \xbe \xdf (\xbS_{1})$ $ \xch $
 \xEH $ \xba \xcn_{L_{i}} \xbb $ $ \xch $ $ \xba \xcn_{L} \xbb $ \xEH \xEH
\xEH \xEP
\cline{1-4}

$1*x \xch x$ \xEH
\multicolumn{3}{|c|}{trivial}

 \xEH
$ \xbG_{1} \xCK \xbS_{2} \xbe \xdf (\xbS_{1} \xCK \xbS_{2})$
 \xEH \xEH \xEH \xEH \xEP
\hline

$x*s \xch s$ \xEH $ \xCf (iM)$ \xEH $ \xba \xcn \xCN \xbb $ $ \xch $ \xEH
-
 \xEH
$ \xbG_{1} \xbe I(\xbS_{1}) \xch $
 \xEH $ \xba \xcn_{L_{1}} \xbb,$ $ \xbb' \xcl_{L_{2}} \xba' $ $ \xch $
\xEH \xEH \xEH \xEP
 \xEH
$A \xcc B \xbe \xdi (X)$ $ \xch $ $A \xbe \xdi (X)$
 \xEH $ \xba \xcn \xCN \xbb \xco \xbg $ \xEH \xEH
$ \xbG_{1} \xCK \xbG_{2} \xbe I(\xbS_{1} \xCK \xbS_{2})$
 \xEH
$ \xba \xcu \xba' \xcn (\xbb \xcu \xba') \xco (\xba \xcu \xbb')$
 \xEH \xEH \xEH \xEP
\cline{1-4}

$s*x \xch s$ \xEH $(eM \xdi)$ \xEH $ \xba \xcu \xbb \xcn \xCN \xbg $ $
\xch $ \xEH - \xEH \xEH \xEH \xEH \xEH \xEP
 \xEH
$X \xcc Y$ $ \xch $ $ \xdi (X) \xcc \xdi (Y),$
 \xEH $ \xba \xcn \xCN \xbb \xco \xCN \xbg $ \xEH \xEH \xEH \xEH \xEH \xEH
\xEP
 \xEH
$X \xcc Y$ $ \xch $
 \xEH \xEH \xEH \xEH \xEH \xEH \xEH \xEP
 \xEH
$ \xdf (Y) \xcs \xdp (X) \xcc \xdf (X)$
 \xEH \xEH \xEH \xEH \xEH \xEH \xEH \xEP
\hline

$b*b \xch b$ \xEH $(< \xbo *s),$ $(\xdm^{+}_{ \xbo })$ (3) \xEH $ \xba
\xcn \xbb,$ $ \xba \xcu \xbb \xcn \xbg $ \xEH - (Filter)
 \xEH
$ \xbG_{1} \xbe \xdf (\xbS_{1}), \xbG_{2} \xbe \xdf (\xbS_{2})$ $ \xch $
 \xEH $ \xba \xcn_{L_{1}} \xbb,$ $ \xba' \xcn_{L_{2}} \xbb' $ $ \xch $
\xEH $b*b \xcj b$
 \xEH $ \xCf (GH)$ \xEH $ \xcn \xDO \xcn $ \xEP
$(\xbm *1)$ \xEH
$A \xbe \xdf (X),X \xbe \xdf (Y)$ $ \xch $
 \xEH $ \xch $ $ \xba \xcn \xbg $ \xEH \xEH
$ \xbG_{1} \xCK \xbG_{2} \xbe \xdf (\xbS_{1} \xCK \xbS_{2})$
 \xEH
$ \xba \xcu \xba' \xcn_{L} \xbb \xcu \xbb' $
 \xEH $(\xbm *1)$ \xEH \xEH \xEP
 \xEH $A \xbe \xdf (Y)$ \xEH \xEH \xEH \xEH \xEH \xEH \xEH \xEP
\cline{1-6}

$b*m \xch m$ \xEH $(< \xbo *s),$ $(\xdm^{+}_{ \xbo })$ (2) \xEH $ \xba
\xcn \xbb,$ $ \xba \xcu \xbb \xcN \xCN \xbg $
 \xEH - (Filter) \xEH
$ \xbG_{1} \xbe \xdf (\xbS_{1}), \xbG_{2} \xbe \xdm^{+}(\xbS_{2})$ $
\xch $
 \xEH
$ \xba \xcN_{L_{1}} \xCN \xbb,$ $ \xba' \xcn_{L_{2}} \xbb' $ $ \xch $
 \xEH \xEH \xEH \xEP
 \xEH
$A \xbe \xdm^{+}(X),X \xbe \xdf (Y)$ $ \xch $
 \xEH $ \xch $ $ \xba \xcN \xCN \xbb \xco \xCN \xbg $
 \xEH \xEH
$ \xbG_{1} \xCK \xbG_{2} \xbe \xdm^{+}(\xbS_{1} \xCK \xbS_{2})$
 \xEH
$ \xba \xcu \xba' \xcN_{L} \xCN (\xbb \xcu \xbb')$
 \xEH \xEH \xEH \xEP
 \xEH $A \xbe \xdm^{+}(Y)$ \xEH \xEH \xEH \xEH \xEH \xEH \xEH \xEP
\cline{1-4}

$m*b \xch m$ \xEH $(< \xbo *s),$ $(\xdm^{+}_{ \xbo })$ (1) \xEH
$ \xba \xcN \xCN \xbb,$ $ \xba \xcu \xbb \xcn \xbg $
 \xEH - (Filter) \xEH
 \xEH \xEH \xEH \xEH \xEP
 \xEH
$A \xbe \xdf (X),X \xbe \xdm^{+}(Y)$ $ \xch $
 \xEH $ \xch $ $ \xba \xcN \xCN \xbb \xco \xCN \xbg $
 \xEH \xEH \xEH \xEH \xEH \xEH \xEP
 \xEH $A \xbe \xdm^{+}(Y)$ \xEH \xEH \xEH \xEH \xEH \xEH \xEH \xEP
\cline{1-6}

$m*m \xch m$ \xEH $(\xdm^{++})$ \xEH Rational Monotony \xEH ranked \xEH
$ \xbG_{1} \xbe \xdm^{+}(\xbS_{1}), \xbG_{2} \xbe \xdm^{+}(\xbS_{2})$
 \xEH
$ \xba \xcN_{L_{1}} \xCN \xbb,$ $ \xba' \xcN_{L_{2}} \xCN \xbb' $ $
\xch $
 \xEH \xEH \xEH \xEP
 \xEH
$A \xbe \xdm^{+}(X),X \xbe \xdm^{+}(Y)$
 \xEH \xEH \xEH $ \xch $
 \xEH
$ \xba \xcu \xba' \xcN_{L} \xCN (\xbb \xcu \xbb')$
 \xEH \xEH \xEH \xEP
 \xEH $ \xch $ $A \xbe \xdm^{+}(Y)$ \xEH \xEH \xEH
$ \xbG_{1} \xCK \xbG_{2} \xbe \xdm^{+}(\xbS_{1} \xCK \xbS_{2})$
 \xEH \xEH \xEH \xEH \xEP
\hline

$b*b \xcj b,$ \xEH \xEH \xEH \xEH \xEH
$ \xba \xcn \xbb $ $ \xch $ $ \xba \xex L_{1} \xcn \xbb \xex L_{1}$
 \xEH $(\xbm *1)$ $+$
 \xEH (GH3) \xEH $ \xcl \xDO \xcn $ \xEP
$pr(b)=b$ \xEH \xEH \xEH \xEH \xEH
and
 \xEH $(\xbm *2)$
 \xEH \xEH \xEP
$(\xbm *2)$ \xEH \xEH \xEH \xEH \xEH
$ \xba \xcn_{L_{1}} \xbb,$ $ \xba' \xcn_{L_{2}} \xbb' $ $ \xch $
 \xEH \xEH \xEH \xEP
 \xEH \xEH \xEH \xEH \xEH
$ \xba \xcu \xba' \xcn_{L} \xbb \xcu \xbb' $
 \xEH \xEH \xEH \xEP
\hline

$J' $ small \xEH \xEH \xEH \xEH
 \xEH $ \xba \xcu \xba' \xcn \xbb \xcu \xbb' $ $ \xcj $
 \xEH \xEH $forget(J')$ \xEH - \xEP
 \xEH \xEH \xEH \xEH \xEH
$ \xba \xcn \xbb,$ $ \xba' \xcl \xbb' $ \xEH \xEH \xEH \xEP
\hline
\multicolumn{9}{|c|}{

Theory revision
} \xEP
\hline

 \xEH \xEH \xEH \xEH $(\xfA *):$ \xEH \xEH \xEH $ \xCf (GHD)$ \xEH $(
\xbf \xcu \xbf')*(\xbq \xcu \xbq') \xcl \xbr $ \xEP
 \xEH \xEH \xEH \xEH $(\xbS_{1} \xCK \xbS_{1}') \xfA (\xbS_{2} \xCK
\xbS_{2}')=$ \xEH $(\xbf \xcu \xbf')*(\xbq \xcu \xbq')=$
 \xEH \xEH \xEH $ \xch $ $ \xbf' * \xbq' \xcl \xbr $ \xEP
 \xEH \xEH \xEH \xEH $(\xbS_{1} \xfA \xbS_{2}) \xCK (\xbS_{1}' \xfA
\xbS_{2}')$ \xEH
$(\xbf * \xbq) \xcu (\xbf' * \xbq')$
 \xEH \xEH \xEH $ \xbf, \xbq $ in $J,$ \xEP
 \xEH \xEH \xEH \xEH
 \xEH \xEH \xEH \xEH $ \xbf', \xbq', \xbr $ in $L- \xCf J$ \xEP
\hline
\end{tabular}
}
\end{turn}
\end{center}
\end{table}

$ \xCO $

Note that $A \xCK B \xcc X \xCK Y$ big $ \xch $ $A \xcc X$ big etc. is
intuitively better justified
than the other direction, as the proportion might increase in the latter,
decrease in the former. Cf.
the table ``Rules on size'',
Section \ref{Section Add-Size} (page \pageref{Section Add-Size}),
``increasing proportions''.
\clearpage
\subsection{
Language change in classical and non-monotonic logic
}

\label{Section Lang-Manip}

\bfa

$\hspace{0.01em}$


\label{Fact Mod-Lang-Fact}

We can obtain factorization by language change, provided cardinalities
permit this.

\efa

\subparagraph{
Proof
}

$\hspace{0.01em}$


Consider $k$ variables, suppose we have $p=m*n$ positive instances, and
that
we can divide $k$ into $k' $ and $k'' $ such that $2^{k' } \xcg m,$
$2^{k'' } \xcg n,$ then we can factorize:

Choose $m$ sequences of 0/1 of length $k',$ $n$ sequences of length $k''
.$ They will code
the positive instances: there are $p=m*n$ pairs of the chosen sequences.
Take any bijection between these pairs and the positive instances, and
continue the bijection arbitrarily between other pairs and negative
instances.

We can do the same also for 2 sets, corresponding to $K,$ $ \xbf $ to have
a common
factorization, they both have to admit common factors like $m,n$ above.
We then choose first the pairs for e.g. $K,$ then for $ \xbf,$ then the
rest.

$ \xcz $
\\[3ex]

\be

$\hspace{0.01em}$


\label{Example Mod-Lang-Fact}

Consider $p=3,$ and let

$ \xCf abc,$ $a \xCN bc,$ $a \xCN b \xCN c,$ $ \xCN abc,$ $ \xCN a \xCN b
\xCN c,$ $ \xCN ab \xCN c$ be the $6=2*3$ positive cases,

$ab \xCN c,$ $ \xCN a \xCN bc$ the negative ones. (It is coincidence that
we can factorize
positive and negative cases - probably iff one of the factors is the full
product, here 2, it could also be 4 etc.)

We divide the cases by 3 new variables, grouping them together in
positive and negative cases. $a' $ is indifferent, we want this to be the
independent factor, the negative ones will be put into $ \xCN b' \xCN c'
.$
The procedure has to be made precise still. (n): negative

Let $ \xCf a' $ code the set $ \xCf abc,$ $a \xCN bc,$ $a \xCN b \xCN c,$
$ab \xCN c$ (n),

Let $ \xCN a' $ code $ \xCN a \xCN bc$ (n), $ \xCN abc,$ $ \xCN a \xCN b
\xCN c,$ $ \xCN ab \xCN c.$

Let $ \xCf b' $ code $ \xCf abc,$ $a \xCN bc,$ $ \xCN a \xCN b \xCN c,$ $
\xCN ab \xCN c$

Let $ \xCN b' $ code $a \xCN b \xCN c,$ $ab \xCN c$ (n), $ \xCN a \xCN bc$
(n), $ \xCN abc$

Let $ \xCf c' $ code $ \xCf abc,$ $a \xCN b \xCN c,$ $ \xCN abc,$ $ \xCN a
\xCN b \xCN c$

Let $ \xCN c' $ code $a \xCN bc,$ $ab \xCN c$ (n), $ \xCN a \xCN bc$ (n),
$ \xCN ab \xCN c$

Then the 6 positive instances are

$\{a', \xCN a' \} \xCK \{b' c',$ $b' \xCN c',$ $ \xCN b' c' \},$ the
negative ones

$\{a', \xCN a' \} \xCK \{ \xCN b' \xCN c' \}$

As we have 3 new variables, we code again all possible cases, so
expressivity is the same.

$ \xcz $
\\[3ex]

\ee

The same holds for non-monotonic logic.

We give an example:

\be

$\hspace{0.01em}$


\label{Example Non-Mon-Lang-Change}

Suppose we have the rule that
``positive is better than negative''
(all other things equal).
Then, for two variables, $ \xCf a$ and $b,$ we have the comparisons
$ab \xeb a \xCN b \xeb \xCN a \xCN b,$ and $ab \xeb \xCN ab \xeb \xCN a
\xCN b.$

Suppose now we are given the situation $ \xCN cd \xeb c \xCN d \xeb cd$
and $ \xCN cd \xeb \xCN c \xCN d \xeb cd,$
which has the same order structure, but with negations not fitting.
We put $c \xCN d$ and $ \xCN cd$ into a new variable $a',$ cd and $ \xCN
c \xCN d$ into $ \xCN a',$
$ \xCN cd$ and $ \xCN c \xCN d$ into $b',$ and $c \xCN d$ and cd into $
\xCN b'.$
Then $a' b' $ corresponds to $ \xCN cd,$ $a' \xCN b' $ to $c \xCN d,$ $
\xCN a' b' $ to $ \xCN c \xCN d,$
$ \xCN a' \xCN b' $ to cd - and we have the desired structure. Thus, if
the
geometric structure is possible, then we can change language and obtain
the desired pattern.

But we cannot obtain by language change a pattern of the type
$ab \xeb a \xCN b$ without any other comparison, if it is supposed to be
based
on a component-wise comparison.

\ee

We summarize:

We can cut the model set as we like:

Choose half to go into $p_{0},$ half into $ \xCN p_{0,}$ again half of
$p_{0}$ into $p_{0} \xcu p_{1},$ half
into $p_{0} \xcu \xCN p_{1},$ etc.
\section{
Semantic interpolation for non-monotonic logic
}

\label{Section Sem-Int-NML}
\subsection{
Discussion
}

We discuss here the full non-monotonic case, i.e., downward and upward.
We consider here a non-monotonic logic $ \xcn.$
We look at the interpolation problem in 3 ways.

Given $ \xbf \xcn \xbq,$ there is an interpolant $ \xba $ such that

(1) $ \xbf \xcn \xba \xcl \xbq,$
see Section \ref{Section Sin-Non-Mon-Int-Karl} (page \pageref{Section
Sin-Non-Mon-Int-Karl}),

(2) $ \xbf \xcl \xba \xcn \xbq,$
see Section \ref{Section Sin-Non-Mon-Int-Dov} (page \pageref{Section
Sin-Non-Mon-Int-Dov}),

(3) $ \xbf \xcn \xba \xcn \xbq,$
see Section \ref{Section Mul-Int} (page \pageref{Section Mul-Int}).

The first variant will be fully characterized below.

The second and third variant have no full characterization at the time of
writing (to the authors' knowledge), but are connected to very interesting
properties about multiplication of size and componentwise independent
relations.

We begin with the following negative result:

\be

$\hspace{0.01em}$


\label{Example Sin-Non-Mon-Int}

Full non-monotonic logics, i.e. down and up, has not necessarily
interpolation.

Consider the model order $pq \xeb p \xCN q \xeb \xCN p \xCN q \xeb \xCN
pq.$ Then $ \xCN p \xcn \xCN q,$ there are no
common variables, and $true \xcn q$ (and, of course, $ \xCN p \xcN
false).$
(Full consequence of $ \xCN p$ is $ \xCN p \xCN q,$ so this has
trivial interpolation.)
\subsection{
Interpolation of the form $\xbf \xcn \xba \xcl \xbq $
}

\label{Section Sin-Non-Mon-Int-Karl}

\ee

\bfa

$\hspace{0.01em}$


\label{Fact Sin-Int-Aux}

Let $var(\xbf)$ be the set of relevant variables of $ \xbf.$

Let $ \xbS \xcc \xbP X,$ $var(\xba) \xcs var(\xbb)= \xCQ,$ $var(\xbb
) \xcs R(\xbS)= \xCQ,$ $ \xbb $ not a tautology, then
$ \xbS \xcc M(\xba \xco \xbb)$ $ \xch $ $ \xbS \xcc M(\xba).$

\efa

\subparagraph{
Proof
}

$\hspace{0.01em}$


Suppose not, so there is $ \xbs \xbe \xbS $ such that $ \xbs \xcm \xba
\xco \xbb,$ $ \xbs \xcM \xba.$ As $ \xbb $ is not
a tautology, there is an assignment to $var(\xbb)$ which makes $ \xbb $
wrong. Consider
$ \xbt $ such that $ \xbs = \xbt $ except on $var(\xbb),$ where $ \xbt $
makes $ \xbb $ wrong, using this assignment.
By $var(\xba) \xcs var(\xbb)= \xCQ,$ $ \xbt \xcm \xCN \xba.$ By
$var(\xbb) \xcs R(\xbS)= \xCQ,$ $ \xbt \xbe \xbS.$
So $ \xbt \xcM \xba \xco \xbb $ for some $ \xbt \xbe \xbS,$
contradiction.

$ \xcz $
\\[3ex]

$ \xCO $

\vspace{10mm}

\begin{diagram}

\label{Diagram Xcn-Xcl}
\index{Diagram Xcn-Xcl}

\centering
\setlength{\unitlength}{1mm}
{\renewcommand{\dashlinestretch}{30}
\begin{picture}(150,100)(0,0)

\path(20,90)(90,90)(90,30)(20,30)(20,90)

\path(40,80)(80,80)
\path(30,70)(60,70)
\path(45,55)(60,55)
\path(45,50)(75,50)

\path(40,79)(40,81)
\path(80,79)(80,81)
\path(30,69)(30,71)
\path(60,69)(60,71)
\path(45,54)(45,56)
\path(60,54)(60,56)
\path(45,49)(45,51)
\path(75,49)(75,51)

\put(92,80){{\xssc $\xbS = M(\xbf)$}}
\put(92,70){{\xssc $\xbS' = M(\xbq)$}}
\put(92,55){{\xssc $\xbS'' = M(\xba)$}}
\put(92,50){{\xssc $\xbm(\xbS) = \xbm(\xbf)$}}

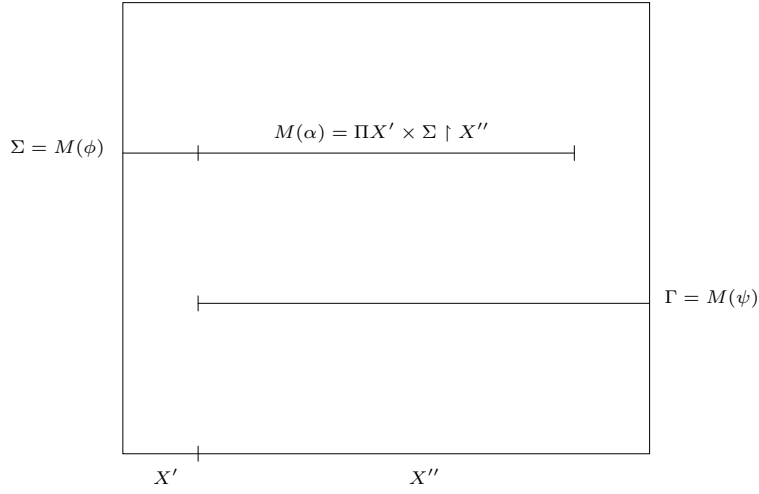
\put(24,14){{\xssc Non-monotonic interpolation, $\xbf \xcn \xba \xcl \xbq$}}

\end{picture}
}

\end{diagram}

\vspace{4mm}

$ \xCO $

\bp

$\hspace{0.01em}$


\label{Proposition Sin-Non-Mon-Int-Karl}

We use here normal forms (conjunctions of disjunctions).

Consider a finite language. Let a semantic choice function $ \xbm $ be
given,
as discussed in
Section \ref{Section Sin-Introduction} (page \pageref{Section Sin-Introduction})
, defined for sets of sequences
(read: models).

$ \xcn $ has interpolation iff for all $ \xbS $ $I(\xbS) \xcc I(\xbm (
\xbS))$ holds.

In the infinite case, we need as additional prerequisite that
$ \xbm (\xbS)$ is definable if $ \xbS $ is.

\ep

\subparagraph{
Proof
}

$\hspace{0.01em}$


Work with reformulations of $ \xbS $ etc. which use only essential $(=$
relevant)
variables.

``$ \xch $'':

Suppose the condition is wrong. Then $X:=I(\xbS)-I(\xbm (\xbS))=I(
\xbS) \xcs R(\xbm (\xbS)) \xEd \xCQ.$
Thus there is some $ \xbs' \xbe \xbm (\xbS) \xex R(\xbS)$ which
cannot be continued by some choice $ \xbr $ in $X \xcv (I(\xbS) \xcs I(
\xbm (\xbS)))$ in $ \xbm (\xbS),$ i.e.
$ \xbs' \xbr \xce \xbm (\xbS).$

We first do the finite case:
We identify models with their describing formulas.
Consider the
formula $ \xbf:= \xbs' \xcp \xCN \xbr = \xCN \xbs' \xco \xCN \xbr.$ We
have $Th(\xbS) \xcn \xbf,$ as $ \xbm (\xbS) \xcc M(\xbf).$
Suppose $ \xbS'' $ is a semantical interpolant for $ \xbS $ and $ \xbf.$
So $ \xbm (\xbS) \xcc \xbS'' \xcc M(\xbf),$
and $ \xbS'' $ does not contain any variables in $ \xbr $ as essential
variables.
By Fact \ref{Fact Sin-Int-Aux} (page \pageref{Fact Sin-Int-Aux}), $ \xbm (\xbS
) \xcc \xbS'' \xcc M(
\xCN \xbs'),$ but $ \xbs' \xbe \xbm (\xbS) \xex R(\xbS),$
contradiction.

We turn to the infinite case. Consider again $ \xbs' \xbr.$ As $ \xbs'
\xbr \xce \xbm (\xbS),$
and $ \xbm (\xbS)$ is definable,
there is some formula $ \xbf $ which holds in $ \xbm (\xbS),$ but fails
in $ \xbs' \xbr.$
Thus, $Th(\xbS) \xcn \xbf.$ Write $ \xbf $ as a disjunction of
conjunctions.
Let $ \xbS'' $ be an interpolant for $ \xbS $ and $M(\xbf).$ Thus $
\xbm (\xbS) \xcc \xbS'' \xcc M(\xbf),$ and
$ \xbs' \xbr \xce M(\xbf),$ so $ \xbm (\xbS) \xcc \xbS'' \xcc M(
\xbf) \xcc M(\xCN \xbs' \xco \xCN \xbr),$ so
$ \xbS'' \xcc M(\xCN \xbs')$ by Fact \ref{Fact Sin-Int-Aux} (page
\pageref{Fact Sin-Int-Aux}).
So $ \xbm (\xbS) \xcm \xCN \xbs',$ contradiction, as $ \xbs' \xbe
\xbm (\xbS) \xex R(\xbS).$
(More precisely, we have to argue here with not necessarily definable
model sets.)

``$ \xci $'':

Let $I(\xbS) \xcc I(\xbm (\xbS)).$ Let $ \xbS \xcn \xbS',$ i.e. $
\xbm (\xbS) \xcc \xbS'.$ Write $ \xbm (\xbS)$ as a (possibly
infinite) conjunction of
disjunctions, using only relevant variables. Form $ \xbS'' $ from $ \xbm
(\xbS)$
by omitting all variables in this description which are not in $R(\xbS'
).$
Note that all remaining variables are in $R(\xbm (\xbS)) \xcc R(\xbS
),$ so
$ \xbS'' $ is a candidate for interpolation.
See Diagram \ref{Diagram Xcn-Xcl} (page \pageref{Diagram Xcn-Xcl}).

(1) $ \xbm (\xbS) \xcc \xbS'':$ Trivial.

(2) $ \xbS'' \xcc \xbS':$
Let $ \xbs \xbe \xbS''.$ Then there is $ \xbt \xbe \xbm (\xbS) \xcc
\xbS' $ such that $ \xbs \xex R(\xbS')= \xbt \xex R(\xbS'),$ so
$ \xbs \xbe \xbS'.$

A shorter argument is as follows:
$ \xbm (\xbS) \xcm M(\xbS')$ has a semantic interpolant by
Section \ref{Section Mon-Sem-Int} (page \pageref{Section Mon-Sem-Int}), which
is by prerequisite
also an interpolant for $ \xbS $ and $ \xbS'.$

It remains to show in the infinite case that $ \xbS'' $ is definable.
This can be
shown as in Proposition \ref{Proposition Sin-Simplification-Definable} (page
\pageref{Proposition Sin-Simplification-Definable}).

$ \xcz $
\\[3ex]
\subsection{
Interpolation of the form $\xbf \xcl \xba \xcn \xbq $
}

\label{Section Sin-Non-Mon-Int-Dov}

This situation is much more interesting than the last one,
discussed in
Section \ref{Section Sin-Non-Mon-Int-Karl} (page \pageref{Section
Sin-Non-Mon-Int-Karl}). In this section, and the
next one,
Section \ref{Section Mul-Int} (page \pageref{Section Mul-Int}), we connect
abstract multiplication laws
for
size to interpolation. To our knowledge, such multiplication laws
are considered here for the first time, and so also their
connection to interpolation problems.

We introduced two sets of three conditions about abstract size
(see Definition \ref{Definition Sin-Size-Rules} (page \pageref{Definition
Sin-Size-Rules})), and then showed
in Proposition \ref{Proposition Sin-Product-Small} (page \pageref{Proposition
Sin-Product-Small})  that both sets
are equivalent.

We show now that the first two, or the last condition entail
interpolation,
see Proposition \ref{Proposition Sin-Interpolation-1} (page \pageref{Proposition
Sin-Interpolation-1}).

Recall that,
in preferential structures, size is generated by a relation. $A \xcc B$ is
a big
subset iff A contains all minimal elements of $B$ (with respect to this
relation). Hamming relations, see
Definition \ref{Definition Sin-Set-HR} (page \pageref{Definition Sin-Set-HR}),
generate a notion
of size which satisfies our multiplicative conditions (if they are
smooth).

Thus, if a preferential logic is defined by a smooth Hamming relation,
it has semantic interpolation (in our sense here). This is summarized in
Corollary \ref{Corollary Sin-Interpolation-2} (page \pageref{Corollary
Sin-Interpolation-2}).

$ \xCO $

\vspace{10mm}

\begin{diagram}

\label{Diagram Xcl-Xcn}
\index{Diagram Xcl-Xcn}

\centering
\setlength{\unitlength}{1mm}
{\renewcommand{\dashlinestretch}{30}
\begin{picture}(150,100)(0,0)

\path(20,90)(90,90)(90,30)(20,30)(20,90)

\path(20,70)(80,70)
\path(30,50)(90,50)

\path(30,69)(30,71)
\path(30,49)(30,51)
\path(30,29)(30,31)
\path(80,69)(80,71)

\put(5,70){{\xssc $\xbS = M(\xbf)$}}
\put(92,50){{\xssc $\xbG= M(\xbq)$}}
\put(40,72){{\xssc $M(\xba) = \xbP X' \xCK \xbS \xex X''$}}
\put(24,26){{\xssc $X'$}}
\put(58,26){{\xssc $X''$}}

\put(24,14){{\xssc Non-monotonic interpolation, $\xbf \xcl \xba \xcn \xbq$}}

\end{picture}
}

\end{diagram}

\vspace{4mm}

$ \xCO $

\bp

$\hspace{0.01em}$


\label{Proposition Sin-Interpolation-1}

We assume definability as shown in
Proposition \ref{Proposition Sin-Simplification-Definable} (page
\pageref{Proposition Sin-Simplification-Definable}).

Interpolation of the form $ \xbf \xcl \xba \xcn \xbq $ exists, if

(1) both $(S*1)$ and $(S*2),$

or

(2) $(S*3)$ hold,

when $ \xbb \xcn \xbg $ is defined by:

$ \xbb \xcn \xbg: \xcj \xbm (\xbb)= \xbm (M(\xbb)) \xcc M(\xbg),$
and

$ \xbm (X)$ is the generator of the principal filter over $X.$

(We saw in Example \ref{Example Sin-Prod-Size} (page \pageref{Example
Sin-Prod-Size})  that $(S*2),$ and thus
also
$(S*3),$ will often be too strong.)

\ep

\subparagraph{
Proof
}

$\hspace{0.01em}$


Let $ \xbS:=M(\xbf),$ $ \xbG:=M(\xbq),$ $X' $ the set of variables
only in $ \xbf,$ so $ \xbG = \xbP X' \xCK \xbG \xex X'',$
where $ \xbP X' = \xbP X'.$
Set $ \xba:=Th(\xbP X' \xCK \xbS''),$ where $ \xbS'' = \xbS \xex X''
.$
Note that variables only in $ \xbq $ are automatically taken care of,
as $ \xbS'' $ can be written as a product without mentioning them.
See Diagram \ref{Diagram Xcl-Xcn} (page \pageref{Diagram Xcl-Xcn}).

By prerequisite, $ \xbm (\xbS) \xcc \xbG,$ we have to show $ \xbm (
\xbP X' \xCK \xbS'') \xcc \xbG.$

(1)

$ \xbm (\xbP X' \xCK \xbS'')= \xbm (\xbP X') \xCK \xbm (\xbS'')$
by $(S*1)$ and
Proposition \ref{Proposition Sin-Product-Small} (page \pageref{Proposition
Sin-Product-Small})  (2).
By $ \xbm (\xbS) \xcc \xbG,$ $(S*2),$ and
Proposition \ref{Proposition Sin-Product-Small} (page \pageref{Proposition
Sin-Product-Small})  (2), $ \xbm (\xbS'')=
\xbm (\xbS \xex X'') \xcc \xbG \xex X'',$ so
$ \xbm (\xbP X' \xCK \xbS'')= \xbm (\xbP X') \xCK \xbm (\xbS'')
\xcc \xbm (\xbP X') \xCK \xbG \xex X'' \xcc \xbG.$

(2)

$ \xbm (\xbP X' \xCK \xbS'') \xex X'' \xcc \xbm (\xbS) \xex X'' \xcc
\xbG \xex X'' $ by $(S*3)$ and
Proposition \ref{Proposition Sin-Product-Small} (page \pageref{Proposition
Sin-Product-Small})  (2).
So $ \xbm (\xbP X' \xCK \xbS'') \xcc \xbP X' \xCK (\xbm (\xbP X' \xCK
\xbS'') \xex X'') \xcc \xbP X' \xCK (\xbG \xex X'')= \xbG.$

$ \xcz $
\\[3ex]

The following Corollary puts our results together.

\bco

$\hspace{0.01em}$


\label{Corollary Sin-Interpolation-2}

Interpolation in the form $ \xbf \xcl \xba \xcn \xbq $ exists, when $ \xcn
$ is defined by a smooth
Hamming relation, more precisely, a smooth relation satisfying (GH3).

\eco

\subparagraph{
Proof
}

$\hspace{0.01em}$


We give two proofs:

(1)

By
Proposition \ref{Proposition Mod-Hamming} (page \pageref{Proposition
Mod-Hamming})  $(S*1)$ and $(S*2)$ hold.
Thus, by Proposition \ref{Proposition Sin-Interpolation-1} (page
\pageref{Proposition Sin-Interpolation-1})  (1),
interpolation exists.

(2)

By Proposition \ref{Proposition Mod-Hamming} (page \pageref{Proposition
Mod-Hamming}), $(S*3)$ holds, so by
Proposition \ref{Proposition Sin-Interpolation-1} (page \pageref{Proposition
Sin-Interpolation-1})  (2), interpolation exists.

$ \xcz $
\\[3ex]
\subsection{
Interpolation of the form $\xbf \xcn \xba \xcn \xbq $
}

\label{Section Mul-Int}

The following result is perhaps the main result of the book. The
conditions are natural, and not too strong, and the connection between
those multiplicative properties and interpolation gives quite deep
insights into the basics of non-monotonic logics.

\bp

$\hspace{0.01em}$


\label{Proposition Mul-Mu*1-Int}

$(\xbm *1)$ entails semantical interpolation of the form $ \xbf \xcn \xba
\xcn \xbq $ in 2-valued
non-monotonic logic generated by minimal model sets.
(As the model sets might not be definable, syntactic interpolation does
not follow automatically.)

\ep

\subparagraph{
Proof
}

$\hspace{0.01em}$


Let the product be defined on $J \xcv J' \xcv J'' $ (i.e., $J \xcv J' \xcv
J'' $ is the set
of propositional variables in the intended application).
Let $ \xbf $ be defined on $J \xcv J',$ $ \xbq $ on $J' \xcv J''.$
See Diagram \ref{Diagram Mul-Base-2-b} (page \pageref{Diagram Mul-Base-2-b}).

We abuse notation and write $ \xbf \xcn \xbS $ if $ \xbm (\xbf) \xcc
\xbS.$ As usual,
$ \xbm (\xbf)$ abbreviates $ \xbm (M(\xbf)).$

For clarity, even if it clutters up notation, we will be precise about
where
$ \xbm $ is formed. Thus, we write $ \xbm_{J \xcv J' \xcv J'' }(X)$ when
we take the minimal elements
in the full product, $ \xbm_{J}(X)$ when we consider only the product on
$J,$ etc.
$X_{J}$ will be shorthand for $ \xbP \{X_{j}:j \xbe J\}.$

Let $ \xbf \xcn \xbq,$ i.e., $ \xbm_{J \xcv J' \xcv J'' }(\xbf) \xcc M(
\xbq).$
We show that $X_{J} \xCK (\xbm_{J \xcv J' \xcv J'' }(\xbf) \xex J')
\xCK X_{J'' }$ is a semantical interpolant,
i.e., that
$ \xbm_{J \xcv J' \xcv J'' }(\xbf)$ $ \xcc $ $X_{J} \xCK (\xbm_{J \xcv
J' \xcv J'' }(\xbf) \xex J') \xCK X_{J'' },$
and that
$ \xbm_{J \xcv J' \xcv J'' }(X_{J} \xCK (\xbm_{J \xcv J' \xcv J'' }(\xbf
) \xex J') \xCK X_{J'' }) \xcc M(\xbq).$

The first property is trivial, we turn to the second.

(1) As $M(\xbf)=M(\xbf) \xex (J \xcv J') \xCK X_{J'' },$
$ \xbm_{J \xcv J' \xcv J'' }(\xbf)$ $=$ $ \xbm_{J \xcv J' }(M(\xbf)
\xex (J \xcv J')) \xCK \xbm_{J'' }(X_{J'' })$ by $(\xbm *1).$

(2) By $(\xbm *1)$ again,
$ \xbm_{J \xcv J' \xcv J'' }(X_{J} \xCK (\xbm_{J \xcv J' \xcv J'' }(\xbf
) \xex J') \xCK X_{J'' })$ $=$
$ \xbm_{J}(X_{J}) \xCK \xbm_{J' }(\xbm_{J \xcv J' \xcv J'' }(\xbf) \xex
J') \xCK \xbm_{J'' }(X_{J'' }).$

So it suffices to show $ \xbm_{J}(X_{J}) \xCK \xbm_{J' }(\xbm_{J \xcv J'
\xcv J'' }(\xbf) \xex J') \xCK \xbm_{J'' }(X_{J'' }) \xcm \xbq.$

Proof:
Let $ \xbs = \xbs_{J} \xbs_{J' } \xbs_{J'' }$ $ \xbe $
$ \xbm_{J}(X_{J}) \xCK \xbm_{J' }(\xbm_{J \xcv J' \xcv J'' }(\xbf) \xex
J') \xCK \xbm_{J'' }(X_{J'' }),$
so $ \xbs_{J} \xbe \xbm_{J}(X_{J}).$

By $ \xbm_{J' }(\xbm_{J \xcv J' \xcv J'' }(\xbf) \xex J') \xcc \xbm_{J
\xcv J' \xcv J'' }(\xbf) \xex J',$
there is $ \xbs' = \xbs_{J}' \xbs_{J' }' \xbs_{J'' }' \xbe \xbm_{J \xcv
J' \xcv J'' }(\xbf)$
such that $ \xbs'_{J' }= \xbs_{J' },$ i.e. $ \xbs' = \xbs'_{J} \xbs_{J'
} \xbs'_{J'' }.$ As $ \xbs' \xbe \xbm_{J \xcv J' \xcv J'' }(\xbf),$ $
\xbs' \xcm \xbq.$

By (1) and $ \xbs_{J'' } \xbe \xbm_{J'' }(X_{J'' })$ also $ \xbs'_{J}
\xbs_{J' } \xbs_{J'' } \xbe \xbm_{J \xcv J' \xcv J'' }(\xbf),$ so also
$ \xbs'_{J} \xbs_{J' } \xbs_{J'' } \xcm \xbq.$

But $ \xbq $ does not depend on $J,$ so also
$ \xbs $ $=$ $ \xbs_{J} \xbs_{J' } \xbs_{J'' }$ $ \xcm $ $ \xbq.$

$ \xcz $
\\[3ex]

$ \xCO $

\vspace{10mm}

\begin{diagram}

\label{Diagram Mul-Base-2-b}

\index{Diagram Mul-Base-2}

\centering
\setlength{\unitlength}{1mm}
{\renewcommand{\dashlinestretch}{30}
\begin{picture}(150,100)(0,0)

\path(20,90)(90,90)(90,30)(20,30)(20,90)

\path(20,89.5)(40,89.5)
\path(70,89.5)(90,89.5)
\path(40,49.5)(70,49.5)

\path(20,30)(20,28)
\path(90,30)(90,28)

\path(40,90)(40,28)
\path(70,90)(70,28)

\path(20,75)(70,75)
\path(40,80)(70,70)(90,70)
\path(20,50)(90,50)

\put(16,75){{\xssc $\xbf$}}
\put(93,70){{\xssc $\xbq$}}
\put(12,50){{\xssc $\xbm(\xbf)$}}

\put(30,24){{\xssc $J$}}
\put(54,24){{\xssc $J'$}}
\put(80,24){{\xssc $J"$}}

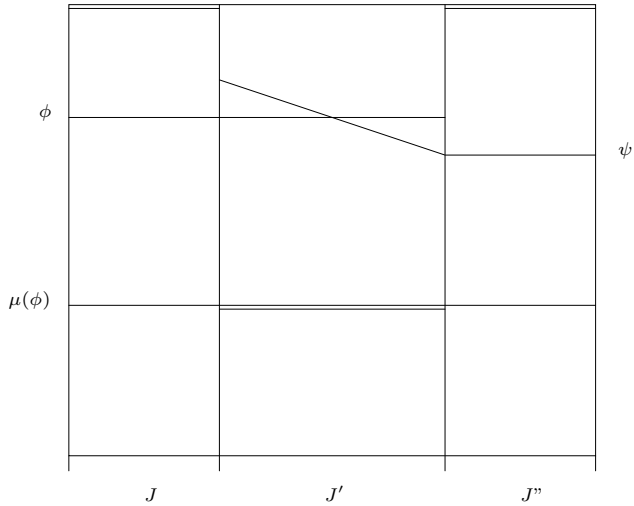
\put(20,14){{\xssc Non-monotonic interpolation, $\xbf \xcn \xba \xcn \xbq$}}
\put(18,10){{\xssc Double lines: interpolant $\xbP J \xCK (\xbm(\xbf)\xex J')
\xCK \xbP J"$}}
\put(12,6){{\xssc Alternative interpolants (in center part): $\xbf\xex J'$ or
$(\xbf\xcu\xbq)\xex J'$}}

\end{picture}
}

\end{diagram}

\vspace{4mm}

$ \xCO $
\subsubsection{
Remarks for the converse: from interpolation to $(\xbm *1)$
}

\be

$\hspace{0.01em}$


\label{Example Mul-Mu}

We show here in (1.1) and (1.2) that half of the condition $(\xbm *1)$ is
not
sufficient for interpolation, and in (2) that interpolation may hold, even
if
$(\xbm *1)$ fails. When looking closer, the latter is not surprising: $
\xbm $ of
sub-products may be defined in a funny way, which has nothing to do with
the
way $ \xbm $ on the big product is defined.

Consider the language based on $p,q,r.$
We define two orders:

(a) $ \xeb $ on sequences of length 3 by
$ \xCN p \xCN q \xCN r \xeb p \xCN q \xCN r,$ and leave all other
sequences of length 3 $ \xeb -$incomparabel.

(b) $<$ on sequences of the same length by
$ \xbs < \xbt $ iff there is $ \xCN x$ is $ \xbs,$ $x$ in $ \xbt,$ but
for no $y$ $y$ in $ \xbs,$ $ \xCN y$ in $ \xbt.$
E.g., $ \xCN p<p,$ $ \xCN pq<pq,$ $ \xCN p \xCN q<pq,$ but $ \xCN pq \xEc
p \xCN q.$

Work now with $ \xeb.$
Let $ \xbf = \xCN q \xcu \xCN r,$ $ \xbq = \xCN p \xcu \xCN q,$ so $ \xbm
(\xbf)= \xCN p \xcu \xCN q \xcu \xCN r,$ and $ \xbf \xcn \xbq.$
Suppose there is $ \xba,$ $ \xbf \xcn \xba \xcn \xbq,$ $ \xba $ written
with $q$ only, so $ \xba $ is equivalent
to FALSE, TRUE, $q,$ or $ \xCN q.$ $ \xbf \xcN FALSE,$ $ \xbf \xcN q.$
$TRUE \xcN \xbq,$ $ \xCN q \xcN \xbq.$ Thus, there
is no such $ \xba,$ and $ \xcn $ has no interpolation.

We show in (1.1) and (1.2) that
we can make both directions of $(\xbm *1)$ true separately, so they do
not
suffice to obtain interpolation.

(1.1) We make $ \xbm (X \xCK Y) \xcc \xbm (X) \xCK \xbm (Y)$ true, but not
the converse.
Take the order $ \xeb $ on sequences of length 3 as described above.
Do not order any sequences of length 2 or 1, i.e. $ \xbm $ is there always
identity. Thus, $ \xbm (X \xCK Y) \xcc X \xCK Y= \xbm (X) \xCK \xbm (Y)$
holds trivially.

(1.2) We make $ \xbm (X \xCK Y) \xcd \xbm (X) \xCK \xbm (Y)$ true, but not
the converse.
We order all sequences of length 1 or 2 by $<,$ and all sequences of
length 3
by $ \xeb.$
Suppose $ \xbs \xbe X \xCK Y- \xbm (X \xCK Y).$ Case 1: $X \xCK Y$
consists of sequences of length 2.
Then, by definition, $ \xbs \xce \xbm (X) \xCK \xbm (Y).$ Case 2: $X \xCK
Y$ consists of sequences
of length 3. Then $ \xbs =p \xCN q \xCN r,$ and $ \xCN p \xCN q \xCN r
\xbe X \xCK Y.$ So
$\{p, \xCN p\} \xcc X$ or $\{p \xCN q, \xCN p \xCN q\} \xcc X,$ but in
both cases $ \xbs \xex X \xce \xbm (X),$ so
$ \xbs \xce \xbm (X) \xCK \xbm (Y).$
Finally, note that $ \xbm (TRUE) \xcC \{ \xCN p \xCN q \xCN r\}$ $=$
$ \xbm (\{p, \xCN p\})$ $ \xCK $
$ \xbm (\{ \xBc q,r \xBe, \xBc q, \xCN r \xBe, \xBc  \xCN q,r \xBe, \xBc 
\xCN q, \xCN r \xBe \}),$
so full $(\xbm *1)$ does not hold.

(2) We make interpolation hold, but $ \xbm (X) \xCK \xbm (Y) \xcC \xbm (X
\xCK Y):$
We order all sequences of length 3 by $<.$ Shorter sequences are made
incomparabel, so for shorter sequences $ \xbm (X)=X.$
Obviously, in general $ \xbm (X) \xCK \xbm (Y) \xcC \xbm (X \xCK Y).$
But the proof of
Proposition \ref{Proposition Mul-Mu*1-Int} (page \pageref{Proposition
Mul-Mu*1-Int})  goes through as above, only
directly, without the
use of factorizing and taking $ \xbm $ of the factors.

The reason is that the order on sequences of length 3 behaves in a modular
way.

Thus, we can have the result, even if subspaces behave differently, if the
completions of the subspace in the whole space behave as they should.

$ \xcz $
\\[3ex]
\subsubsection{
The extreme interpolants in non-monotonic logics
}

\ee

Two possible ways to solve the problem of interpolation
are illustrated in
Diagram \ref{Diagram Nml-Int-Ways} (page \pageref{Diagram Nml-Int-Ways}).
The main differences between the left and the right hand part,
$ \xbm (\xbf) \xex J' $ vs. $ \xbf \xex J',$ are underlined.

$ \xCO $

\vspace{10mm}

\begin{diagram}

\label{Diagram Nml-Int-Ways}
\index{Diagram Nml-Int-Ways}

\centering
\setlength{\unitlength}{1mm}
{\renewcommand{\dashlinestretch}{30}
\begin{picture}(150,150)(0,0)

\put(45,120){
$\xbf$ = $(\xbf \xex J \xcv J') \xCK  \xbP J''$
}

\path(66,118)(100,85)
\path(64,118)(30,105)

\put(27,100){
$\xbm(\xbf)$
}
\path(30,98)(30,85)

\put(10,80){
$\xbP J \xCK {\ul{\xbm(\xbf) \xex J'}} \xCK  \xbP J''$
}
\put(80,80){
$\xbP J \xCK {\ul{\xbf \xex J'}} \xCK  \xbP J''$
}
\path(30,78)(30,65)
\path(100,78)(100,65)

\put(10,60){
$\xbm(\xbP J \xCK \xbm(\xbf) \xex J' \xCK  \xbP J'')$
}
\put(80,60){
$\xbm(\xbP J \xCK \xbf \xex J' \xCK  \xbP J'')$
}

\put(25,55){
$= \hspace{2mm} (?)$
}
\put(95,55){
$= \hspace{2mm} (?)$
}

\put(10,50){
$\xbm(\xbP J) \xCK \xbm(\xbf) \xex J' \xCK \xbm(\xbP J'')$
}
\put(80,50){
$\xbm(\xbP J) \xCK \xbm(\xbf) \xex J' \xCK \xbm(\xbP J'')$
}

\path(100,48)(66,25)
\path(30,48)(64,25)

\put(40,20){
$\xbq$ =  $\xbP J \xCK (\xbq \xex J' \xcv J'')$
}

\end{picture}
}

\end{diagram}

\vspace{4mm}

$ \xCO $

We may take the left path (as done in the proof of
Proposition \ref{Proposition Mul-Mu*1-Int} (page \pageref{Proposition
Mul-Mu*1-Int}))
from top to bottom, apply first $ \xbm,$ then
generalize by taking the product, and then apply $ \xbm $ again.
Note that we apply $ \xbm $ to a product in the last step. This is all
semantical routine (the generalization may be a definability problem when
looking at syntax). It depends on multiplication laws whether the outcome
is $ \xbm (\xbP J) \xCK \xbm (\xbf) \xex J' \xCK \xbm (\xbP J''),$
and whether this entails classically $ \xbq.$

The path on the right generalizes first, and then applies $ \xbm $ to a
product,
but it seems more unclear if this has still anything to do with the
original $ \xbm (\xbf),$ as $ \xbf $ was not a product of 3 components,
only of 2 components.
This has a similarity to the $ \xbf \xcl \xba \xcn \xbq $ problem
discussed in
Section \ref{Section Sin-Non-Mon-Int-Dov} (page \pageref{Section
Sin-Non-Mon-Int-Dov}).
We avoided this by taking $ \xbm $ of the original $ \xbf $ first on the
left hand path.
The rest is similar to the left hand path.

We can also consider a variant which is in a sense between the 2 ways
discussed here (as the underlined part $ \xbm (\xbf) \xex J' $ implies
the underlined part
$ \xbf \xex J'),$ considering $ \xbP J \xCK (\xbf \xcu \xbq \xex J')
\xCK \xbP J''.$ And if we look for a most general
interpolant, leaving aside considerations about multiplications, we may
use the (OR) rule:
If $ \xba \xex J' \xcn \xbq \xex J',$ $ \xba' \xex J' \xcn \xbq \xex J'
,$ then by the (OR) rule, also
$(\xba \xco \xba') \xex J' \xcn \xbq \xex J'.$
\subsection{
Interpolation for distance based revision
}

\label{Section Int-Dist-Rev}

The limiting condition (consistency) imposes a strong restriction:
Even for $ \xbf *TRUE,$ the result may need many variables (those in $
\xbf).$

\bl

$\hspace{0.01em}$


\label{Lemma Mul-TR}

Let $ \xfA $ satisfy $(\xfA *),$ as defined in
Fact \ref{Fact Mul-GHD-Bar} (page \pageref{Fact Mul-GHD-Bar}).

Let $J \xcc L,$ $ \xbr $ be written in sublanguage $J,$ let $ \xbf, \xbq
$ be written in $L- \xCf J,$
let $ \xbf', \xbq' $ be written in $J' \xcc J.$

Let $(\xbf \xcu \xbf')*(\xbq \xcu \xbq') \xcl \xbr,$ then $ \xbf'
* \xbq' \xcl \xbr.$

(This is suitable interpolation, but we also need to factorize the
revision construction.)

\el

\subparagraph{
Proof
}

$\hspace{0.01em}$


$(\xbf \xcu \xbf')*(\xbq \xcu \xbq')$ $=$ $(\xbf * \xbq) \xcu (
\xbf' * \xbq')$ by $(\xfA *).$
So $(\xbf \xcu \xbf')*(\xbq \xcu \xbq') \xcl \xbr $ iff $(\xbf *
\xbq) \xcu (\xbf' * \xbq') \xcl \xbr,$ but
$(\xbf * \xbq) \xcu (\xbf' * \xbq') \xcl \xbr $ iff $(\xbf' * \xbq
') \xcl \xbr,$ as $ \xbr $ contains no variables
of $ \xbf' $ or $ \xbq'.$
$ \xcz $
\\[3ex]
\subsection{
The equilibrium logic EQ
}

\label{Section EQ}
\subsubsection{
Introduction and outline
}

\label{Section Mul-Mon-Int}

Equilibrium logic, EQ, introduced by Pearce et al., see
 \cite{PV09}, is based on the three valued Goedel logic, based on
two worlds,
here and there, also called HT logic. It is defined through a model choice
function on HT models, and non-monotonic. For motivation and context, the
reader is referred to
 \cite{PV09}.

We look at the problem in two different ways.
First, we investigate the more classical case where we are interested
only in models which have maximal value (i.e., 2), and then look at
the more refined situation, as it was described in
Section \ref{Section Eq-Sem} (page \pageref{Section Eq-Sem}). Again, we will
consider, for both
approaches,
the three types of interpolation for $ \xbf \xcn \xbq:$ $ \xbf \xcl \xba
\xcn \xbq,$ $ \xbf \xcn \xba \xcl \xbq,$
and $ \xbf \xcn \xba \xcn \xbq.$ We will show that there is not always
interpolation of the
first two types, but that there is always interpolation of the third type.

In both cases, we have to make a decision: when going from all models to
minimal models, which value will we give to those models we eliminate as
non-minimal? In the classical case, it was trivial, they get value 0. Now,
we could also lower their value by some smaller amount. But we decide to
give them value 0, here, too.
\subsubsection{
Basic definition, and definability of chosen models
}

We first give the basic definition (the model choice function):

Fix a finite propositional language $ \xdl.$

\bd

$\hspace{0.01em}$


\label{Definition Eq-Eq}

(The definition is due to
 \cite{PV09}.)

(1) A model $m$ is total iff for no $a \xbe L$ $m(a)=1.$

(2) $m \xeb m' $ iff

(2.1) for all $a \xbe L$ $m(a)=0$ $ \xcj $ $m' (a)=0$ and

(2.2) $\{a \xbe L:m(a)=2\} \xcb \{a \xbe L:m' (a)=2\}$

$(\xbs \xeb \xbt $ iff $T$ is preserved, and $H$ goes down.
Thus, only changes from 2 to 1 are possible when $ \xbs \xeb \xbt.)$

(3) $m \xbe X$ is an equilibrium model of $X$ iff $m$ is total, and
there is no $m' \xeb m,$ $m' \xbe X.$ (We can add $m \xeb m$ if $m(x)=1$
for some $x \xbe L,$
so we can define equilibrium models as minimal HT models for some
relation $ \xeb.)$

(4) $ \xbm (X)$ will be the ``set'' of equilibrium models of $X.$

Recall that,
if $f_{X}:M \xcp V$ is the model function of $X,$ then
$f_{ \xbm (X)}:M \xcp V$ will be defined as follows:
\begin{flushleft}
\[ f_{\xbm(X)}(m):= \left\{ \begin{array}{lcl}
f_{X}(m) \xEH iff \xEH
m \xDC is \xDC chosen \xEP
\xEH \xEH \xEP
0 \xEH \xEH otherwise \xEP
\end{array}
\right.
\]
\end{flushleft}

\ed

\br

$\hspace{0.01em}$


\label{Remark Eq-Eq}

Consequently, we have a sort of ``anti-smoothnes'': if a model is not
minimal,
then any model below it is NOT chosen. Thus, we cannot use general
results based on smoothness.

\er

We now show:

\bfa

$\hspace{0.01em}$


\label{Fact Eq-Mu-Def}

Work in a finite language.

If $f_{X}$ is definable, then so is $f_{ \xbm (X)}.$

\efa

\subparagraph{
Proof
}

$\hspace{0.01em}$


Note: We do not need a uniform way to find the defining
formula, we just need the formula - even if it is
``handcrafted''. We will consider all models one by one. Finally,
we will AND the definition of $f_{X}$ with suitable
formulae.

If there is no model $m$ such that $f_{X}(m)>0$ and $m(p)=1$ for some $p,$
then
$f_{ \xbm (X)}=f_{X}.$ Consider now all $m$ such that $f_{X}(m)>0$ and
$m(p)=1$ for some $p,$
and define $ \xbf_{m}$ for this $m.$ As the language is finite, the set of
such
$m$ is finite, and if $f_{X}=f_{ \xbf },$ then $f_{ \xbf \xcu \xcU \{
\xbf_{m}: \xcE p.m(p)=1,f_{X}(m)>0\}}$ will
define $f_{ \xbm (X)}.$ A model is eliminated iff it contains $m(p)=1$ for
some $p,$
or there is some $m' \xeb m$ which contains some such $p.$
So, if we can show for $m$ with $m(p)=1$ for some $p$ $ \xbf_{m}(m)=0,$
and
$ \xbf_{m}(m')=0$ if $m \xeb m',$ and perhaps $ \xbm_{m}(m')$ if $m' $
contains some other $m' (p)=1,$
but $ \xbf_{m}(m')=2$ for all other $m',$ then we are done.

Consider now $m$ with $m(p)=1$ for some $p.$ Let $p_{1}, \Xl,p_{k}$ be
such that $m(p_{i})=0,$
and $q_{1}, \Xl,q_{n}$ be such that $m(q_{i})>0.$ $k$ may be 0, but not
$n.$
Define $ \xbf_{m}:= \xCN \xCN p_{1} \xco  \Xl  \xco \xCN \xCN p_{k} \xco
\xCN q_{1} \xco  \Xl  \xco \xCN q_{n}.$

First, consider $m.$ $m(\xCN \xCN p_{i})=0$ for all $i,$ and also $m(
\xCN q_{i})=0$ for all $i,$
so $ \xbf_{m}(m)=0,$ as desired.

Second, suppose $m \xeb m'.$ As $m(x)=0$ iff $m' (x)=0,$ the values for
$p_{i}$ did not
change, so still $m' (\xCN \xCN p_{i})=0,$ and the values for $q_{i}$ may
have changed,
but not to 0, so still $m' (\xCN q_{i})=0,$ and $ \xbf_{m}(m')=0.$

Suppose now $m \xeB m'.$
Case 1: $m' (p_{i}) \xEd 0$ for some $i,$ then $m' (\xCN \xCN p_{i})=2,$
and we are done.
Case 2: $m' (q_{i})=0$ for some $i,$ then $m' (\xCN q_{i})=2,$ and we are
done.
Thus, if $ \xbf_{m}(m')=0,$ then $m(x)=0$ iff $m' (x)=0.$

So $m(x)$ can only differ from $m' (x)$ on the $q_{i},$ and $m' (q_{i})
\xEd 0.$
If for all $q_{i}$ $m' (q_{i})=2,$ then $m \xeb m',$ and we are done.
If not, then $m' (q_{i})=1$ for some $i,$ and it should be eliminated
anyway.

$ \xcz $
\\[3ex]
\subsubsection{
The approach with models of value 2
}

We first define formally what we want to do:

\bd

$\hspace{0.01em}$


\label{Definition Eq-Conseq}

(1) Set $M_{2}(\xbf):=\{m \xbe M:m(\xbf)=2=TRUE\}.$

(2) Set $ \xbm_{2}(\xbf):=\{m \xbe M:m(\xbf)=2,$ and $m$ is an
equilibrium model $\}.$

(3) Set $ \xbf \xcn \xbq $ iff $ \xbm_{2}(\xbf) \xcc M_{2}(\xbq).$

(4) Set $ \xbf \xcl \xbq $ iff $ \xcA m \xbe M.m(\xbf) \xck m(\xbq).$

\ed

We will show that interpolation of type
(a) $ \xbf \xcl \xba \xcn \xbq $ and (b) $ \xbf \xcn \xba \xcl \xbq $ may
fail, but interpolation of type
(c) $ \xbf \xcn \xba \xcn \xbq $ will exist.

Definability of the interpolant will be shown using the definability
results
for HT. We can use the techniques
and results developped there
(``neglecting'' some variables), and see that
the semantical interpolant is definable, so we have also syntactical
interpolation.

\be

$\hspace{0.01em}$


\label{Example Eq-xcl-xcn}

(EQ has no interpolation of the form $ \xbf \xcl \xba \xcn \xbq.)$

Work with 3 variables, $a,b,c.$

Consider $ \xbS:=\{ \xBc 0,2,2 \xBe, \xBc 2,1,0 \xBe, \xBc 2,2,0 \xBe \}.$

By the above, and classical behaviour of ``or'' and
``and'', $ \xbS $ is definable by
$ \xbf:=(\xCN a \xcu b \xcu c) \xco (a \xcu \xCN \xCN b \xcu \xCN c),$
i.e. $ \xbS =\{m:m(\xbf)=2\}.$

Note that $ \xBc 2,2,0 \xBe $ is total, but
$ \xBc 2,1,0 \xBe  \xeb  \xBc 2,2,0 \xBe,$ thus
$ \xbm (\xbS)=\{ \xBc 0,2,2 \xBe \}.$

So $ \xbS \xcn c=2$ (or $ \xbS \xcn \xcX c).$ Let $X':=\{a,b\},$ $X''
:=\{c\}.$

All possible interpolants $ \xbG $ must not contain $ \xCf a$ or $b$ as
essential variables,
and they must contain $ \xbS.$ The smallest candidate $ \xbG $ is $ \xbP
X' \xCK \{0,2\}.$
But $ \xbs:= \xBc 0,0,0 \xBe  \xbe \xbG,$ $ \xbs $ is total, and there cannot be
any $ \xbt \xeb \xbs,$ so
$ \xbs \xbe \xbm (\xbG),$ so $ \xbG \xcN c=2.$

For completeness' sake, we write all elements of $ \xbG:$

$ \xBc 0,0,0 \xBe $ $ \xBc 0,0,2 \xBe $

$ \xBc 0,1,0 \xBe $ $ \xBc 0,1,2 \xBe $

$ \xBc 0,2,0 \xBe $ $ \xBc 0,2,2 \xBe $

$ \xBc 1,0,0 \xBe $ $ \xBc 1,0,2 \xBe $

$ \xBc 1,1,0 \xBe $ $ \xBc 1,1,2 \xBe $

$ \xBc 1,2,0 \xBe $ $ \xBc 1,2,2 \xBe $

$ \xBc 2,0,0 \xBe $ $ \xBc 2,0,2 \xBe $

$ \xBc 2,1,0 \xBe $ $ \xBc 2,1,2 \xBe $

$ \xBc 2,2,0 \xBe $ $ \xBc 2,2,2 \xBe $

Recall that no sequence containing 1 is total, and when we go from 2 to 1,
we have a smaller model.
Thus, $ \xbm (\xbG)=\{ \xBc 0,0,0 \xBe, \xBc 0,0,2 \xBe \}.$

\ee

\be

$\hspace{0.01em}$


\label{Example Eq-xcn-xcl}

(EQ has no interpolation of the form $ \xbf \xcn \xba \xcl \xbq.)$

Consider 2 variables, $a,b,$ and $ \xbS:=\{0,2\} \xCK \{0,1,2\}$

No $ \xbs $ containing 1 can be in $ \xbm (\xbS),$ as a matter of fact,
$ \xbm (\xbS)=\{ \xBc 0,0 \xBe, \xBc 2,0 \xBe \}.$
$ \xbS $ is defined by $a \xco \xCN a,$ $ \xbm (\xbS)$ is defined by $(a
\xco \xCN a) \xcu \xCN b.$

So we have $a \xco \xCN a \xcn b \xco \xCN b,$ even $a \xco \xCN a \xcn
\xCN b.$

The only possible interpolants are TRUE or FALSE.
$a \xco \xCN a \xcN FALSE,$ and $TRUE \xcL \xCN b.$

\ee

\bfa

$\hspace{0.01em}$


\label{Fact EQ-2-Int}

EQ has interpolation of the form $ \xbf \xcn \xba \xcn \xbq.$

\efa

\subparagraph{
Proof
}

$\hspace{0.01em}$


Let $ \xbf \xcn \xbq,$ i.e., $ \xbm_{2}(\xbf) \xcc M_{2}(\xbq).$
We have to find $ \xba $ such that $ \xbm_{2}(\xbf) \xcc M_{2}(\xba),$
and $ \xbm_{2}(\xba) \xcc M_{2}(\xbq).$

Let $J=I(\xbf),$ $J'' =I(\xbq).$ Consider $X:= \xbP J \xCK (\xbm_{2}(
\xbf) \xex J') \xCK \xbP J''.$
By the same arguments (``neglecting'' $J$ and $J''),$ $X$ is definable as
$M_{2}(\xba)$ for some $ \xba.$

Obviously, $ \xbm_{2}(\xbf) \xcc M_{2}(\xba).$ Consider now $
\xbm_{2}(\xba),$ we have to show
$ \xbm_{2}(\xba) \xcc M_{2}(\xbq).$ If $ \xbm_{2}(\xba)= \xCQ,$ we
are done, so suppose there is
$m \xbe \xbm_{2}(\xba).$ Suppose $m \xce M_{2}(\xbq).$ There is $m'
\xbe \xbm_{2}(\xbf),$ $m' \xex J' =m \xex J'.$

We use now $+$ for concatenation.

Consider $m'' =(m \xex J)+m' \xex (J' \xcv J'').$ As $m' \xbe \xbm_{2}(
\xbf) \xcc M_{2}(\xbf),$ and
$M_{2}(\xbf)= \xbP J \xCK M_{2}(\xbf) \xex (J' \xcv J''),$ $m'' \xbe
M_{2}(\xbf).$ $m \xex (J \xcv J')=m'' \xex (J \xcv J'),$ thus
by $J'' \xcc I(\xbq),$ $m'' \xce M_{2}(\xbq).$ Thus, $m'' \xce
\xbm_{2}(\xbf).$ So either there is $n \xbe M_{2}(\xbf)$
such that $n(y)=0$ iff $m'' (y)=0$ and $\{y:n(y)=2\} \xcb \{y:m'' (y)=2\}$
or $m'' (y)=1$ for some
$y \xbe L.$ Suppose $m'' (y)=1$ for some $y.$ $y$ cannot be in $J' \xcv
J'',$ as
$m'' \xex (J' \xcv J'')=m' \xex (J' \xcv J''),$ and $m' \xbe \xbm_{2}(
\xbf).$ $y$ cannot be in $J,$ as
$m'' \xex J=m \xex J,$ and $m \xbe \xbm_{2}(X).$

So there must be $n \xbe M_{2}(\xbf)$ as above. Case 1:
$\{y \xbe J' \xcv J'':n(y)=2\} \xcb \{y \xbe J' \xcv J'':m'' (y)=2\}.$
Then $n' =m' \xex J+n \xex (J' \xcv J'')$ would
eliminate $m' $ from $ \xbm_{2}(\xbf),$ so this cannot be. Thus,
$n \xex (J' \xcv J'')=m'' \xex (J' \xcv J'').$ So
$\{y \xbe J:n(y)=2\} \xcb \{y \xbe J:m'' (y)=2\}=\{y \xbe J:m(y)=2\}$ by
$m'' \xex J=m \xex J.$
Consider now $n' =n \xex J+m \xex (J' \xcv J'').$ $n' \xbe \xbP J \xCK
\xbm_{2}(\xbf) \xex J' \xCK \xbP J''.$
$n' (y)=0$ iff $m(y)=0$ by construction of $n' $ and $n.$ So $n' \xeb m,$
and
$m \xce \xbm_{2}(\xbP J \xCK (\xbm_{2}(\xbf) \xex J') \xCK \xbP J''
),$ contradiction.

$ \xcz $
\\[3ex]
\subsubsection{
The refined approach
}

We consider now more truth values, in the sense that
$ \xbf \xcn \xbq $ iff $f_{ \xbm (f_{ \xbf })} \xck f_{ \xbq }$ - and not
only restricted to value 2, as
in Definition \ref{Definition Eq-Conseq} (page \pageref{Definition Eq-Conseq}).
The arguments and examples will be the same, they are given for
completess' sake
only.

Again, we show that there need not be interpolation of the forms
$ \xbf \xcl \xba \xcn \xbq $ or $ \xbf \xcn \xba \xcl \xbq,$ but there
will be interpolation of the type
$ \xbf \xcn \xba \xcn \xbq.$

\be

$\hspace{0.01em}$


\label{Example Eq-xcl-xcn-ref}

(EQ has no interpolation of the form $ \xbf \xcl \xba \xcn \xbq.)$

Work with 3 variables, $a,b,c.$ Models will be written as $ \xBc 0,0,0 \xBe,$
etc.,
in the obvious meaning.

Consider $ \xbf:=(\xCN a \xcu b \xcu c) \xco (a \xcu \xCN \xCN b \xcu
\xCN c).$

For $ \xBc 0,1,1 \xBe,$ $ \xBc 0,1,2 \xBe,$ $ \xBc 0,2,1 \xBe,$ $ \xBc 1,1,0
\xBe,$ $ \xBc 1,2,0 \xBe,$ $f_{ \xbf }$
has value 1,
for $ \xBc 0,2,2 \xBe,$ $ \xBc 2,1,0 \xBe,$ $ \xBc 2,2,0 \xBe $ $f_{ \xbf }$
has value 2, all other
values are 0.
The only chosen model is $ \xBc 0,2,2 \xBe,$ all others contain 1, or are
minimized.
So $f_{ \xbm (\xbf)}$ has value 2 for $ \xBc 0,2,2 \xBe,$ all other values
are 0.
Obviously, $f_{ \xbm (\xbf)} \xck f_{c},$ so $ \xbf \xcn c.$

As shown in Fact \ref{Fact Mod-Fin-Goed} (page \pageref{Fact Mod-Fin-Goed}), we
can define with $c$ only
$c,$ $ \xCN c,$ $ \xCN \xCN c,$ $c \xcp c,$ $ \xCN (c \xcp c),$ $ \xCN
\xCN c \xcp c$ (up to semantical equivalence).
But none is an interpolant of the type $ \xbf \xcl \xba \xcn \xbq:$
The left hand condition fails for $c,$ $ \xCN c,$ $ \xCN \xCN c,$ $ \xCN
(c \xcp c),$
the right hand condition fails for $c \xcp c$ and $ \xCN \xCN c \xcp c,$
as
$f_{ \xbm (c \xcp c)}(\xBc 0,0,0 \xBe)=f_{ \xbm (\xCN \xCN c \xcp
c)}(\xBc 0,0,0 \xBe)=2.$

\ee

\be

$\hspace{0.01em}$


\label{Example Eq-xcn-xcl-ref}

(EQ has no interpolation of the form $ \xbf \xcn \xba \xcl \xbq.)$

Consider 2 variables, $a,b,$ $ \xbf:=a \xco \xCN a.$

$f_{ \xbf }(m)=1$ iff $m(a)=1,$ and 2 otherwise.
Note that (the model) $ \xBc 1,0 \xBe  \xeb  \xBc 2,0 \xBe,$ but $f_{ \xbf }(
\xBc 1,0 \xBe)=1,$ $f_{
\xbf }(\xBc 2,0 \xBe)=2,$
so this minimization does
not ``count''.
Consequently, $f_{ \xbm (\xbf)}(m)=2$ iff $m= \xBc 0,0 \xBe $ or $m= \xBc 2,0
\xBe,$ and
$f_{ \xbm (\xbf)}(m)=0$
otherwise. Thus, $ \xbf \xcn \xCN b.$ But $ \xbf \xcN FALSE,$ and $TRUE
\xcL \xCN b.$

$ \xcz $
\\[3ex]

\ee

\bfa

$\hspace{0.01em}$


\label{Fact EQ-2-Int-ref}

EQ has interpolation of the form $ \xbf \xcn \xba \xcn \xbq.$

\efa

\subparagraph{
Proof
}

$\hspace{0.01em}$


See Diagram \ref{Diagram Mul-Base-2-b} (page \pageref{Diagram Mul-Base-2-b}) 
for illustration.

Let $L=J \xcv J' \xcv J'',$ $J'' =I(\xbf),$ $J=I(\xbq).$
As $ \xbf $ does not contain any variables in $J'' $ in an essential way,
$f_{ \xbf }(m)=f_{ \xbf }(m')$ if $m \xex J \xcv J' =m' \xex J \xcv J'.$
Thus, if $a \xbe J'',$ and
$m \xex L-\{a\}=m' \xex L-\{a\}=m'' \xex L-\{a\},$ and $m(a)=0,$ $m'
(a)=1,$ $m'' (a)=2,$ then by
$f_{ \xbf }(m)=f_{ \xbf }(m')=f_{ \xbf }(m''),$ neither $m' $ nor $m'' $
survives minimization, i.e.,
$f_{ \xbm (\xbf)}(m')=f_{ \xbm (\xbf)}(m'')=0.$ Thus, if $m(a) \xEd
0$ for some $a \xbe J'',$ then $f_{ \xbm (\xbf)}(m)=0.$
On the other hand, if $m \xex J' \xcv J'' =m' \xex J' \xcv J'',$ then
$f_{ \xbq }(m)=f_{ \xbq }(m').$

Define now the semantic interpolant by
$h(m):=sup\{f_{ \xbm (\xbf)}(m'):m' \xex J' =m \xex J' \}.$ Obviously,
$f_{ \xbm (\xbf)} \xck h,$ so, if $h=f_{ \xba }$ for
some $ \xba,$ then $ \xbf \xcn \xba.$ It remains to show that $f_{ \xbm
(h)} \xck f_{ \xbq },$ then
$ \xba \xcn \xbq,$ and we are done.

For the same reasons as discussed above, $f_{ \xbm (h)}(m)=0$ if $m(a)
\xEd 0$ for
some $a \xbe J \xcv J''.$ Take now arbitrary $m,$ we have to show $f_{
\xbm (h)}(m) \xck f_{ \xbq }(m).$
If $m(a) \xEd 0$ for some $a \xbe J \xcv J'',$ there is nothing to show.
So suppose $m(a)=0$ for all $a \xbe J \xcv J''.$
By the above, $h(m)=sup\{f_{ \xbm (\xbf)}(m')$: $m' \xex J' =m \xex J'
$ $ \xcu $ $ \xcA a \xbe J''.m' (a)=0\},$ so,
as $m(a)=0$ for $a \xbe J'',$
$h(m)=sup\{f_{ \xbm (\xbf)}(m')$: $m' \xex J' \xcv J'' =m \xex J' \xcv
J'' \}.$
By prerequisite, $f_{ \xbm (\xbf)}(m') \xck f_{ \xbq }(m')$ for all
$m',$ but $ \xbq $ does not
contain essential variables in $J,$ so if $m' \xex J' \xcv J'' =m \xex J'
\xcv J'',$ then
$f_{ \xbq }(m)=f_{ \xbq }(m'),$ thus $h(m) \xck f_{ \xbq }(m),$ but $f_{
\xbm (h)} \xck h,$ so $f_{ \xbm (h)} \xck h(m) \xck f_{ \xbq }(m).$

$ \xcz $
\\[3ex]
\section{
Context and structure
}

\label{Section Sin-Context}

The discussion in this Section is intended to open the perspective
and separate support from attack, and, even more broadly, separate
logic from manipulation of model sets. But this is not pursued here, and
intended to be looked at in future research.

We take the importance of condition $(\xbm *3)$ (or $(S*3))$
as occasion for a broader remark.

 \xEh

 \xDH
This condition points to a weakening of the Hamming condition:

Adding new ``branches'' in $X' $ will not give new minimal elements in $X''
,$
but may destroy other minimal elements in $X''.$ This can be achieved
by a sort of semi-rankedness:
If $ \xbr $ and $ \xbs $ are different only in the $X' -$part, then $ \xbt
\xeb \xbr $ iff $ \xbt \xeb \xbs,$
but not necessarily $ \xbr \xeb \xbt $ iff $ \xbs \xeb \xbt.$

 \xDH
In more abstract terms:

When we separate support from attack (support: a branch $ \xbs' $ in $X'
$ supports
a continuation $ \xbs'' $ in $X'' $ iff $ \xbs \xDO \xbs'' $ is minimal,
i.e. not attacked,
attack:
a branch $ \xbt $ in $X' $ attacks a continuation $ \xbs'' $ in $X'' $
iff it prevents all $ \xbs \xDO \xbs'' $
to be minimal), we see that new branches will not support any new
continuations,
but may well attack continuations.

More radically, we can consider paths $ \xbs'' $ as positive information,
$ \xbs' $ as potentialy negative information. Thus, $ \xbP' $ gives
maximal
negative information, and thus smallest set of accepted models.

 \xDH
We can interpret this as follows:
$X'' $ determines the base set.
$X' $ is the context. This determines the choice (subset of the base set).
We compare to preferential structures:
In preferential structures, $ \xeb $ is not part of the
language either, it is context. And we have the same behaviour as
shown in the
fundamental property of preferential structures:
the bigger the set, the more attacks are possible.

 \xDH
The concept of size looks only at the result of support and attack, so it
is necessarily somewhat coarse. Future research should also investigate
both
concepts separately.

 \xEj

We broaden this.

Following a tradition begun by Kripke, one has added structure to the set
of
classical models, reachability, preference, etc. Perhaps one should
emphasize a
more abstract approach,
as done by one the authors e.g. in  \cite{Sch92}, and elaborated in
 \cite{Sch04}, see in particular the distinction between structural
and algebraic
semantics in the latter. Our suggestion is to separate structure from
logic in
the semantics, and to treat what we called context above by a separate
``machinery''. So, given a set $X$ of models, we have some abstract
function $f,$ which chooses the models where the consequences hold,
$f(X).$

Now, we can put into this ``machinery'' whatever we want.

The abstract properties of preferential or modal structures are well
known.

But we can also investigate non-static $f,$ where $f$ changes in function
of what we already did - ``reacting'' to the past.

We can look at usual properties of $f,$ complexity,
generation by some simple structure like a special
machine, etc.

So we advocate the separation of usual, classical semantics,
from the additional properties, which are treated
``outside''.
It might be interesting to forget altogether about logic,
classify those functions or more complicated devices which
correspond to some logical property, and investigate them and
their properties.
\section{
Interpolation for argumentation
}

\label{Section Inter-Arg}

Arguments (e.g., in inheritance), are sometimes ordered by a partial order
only.
We may define:
$ \xbq $ follows from $ \xbf $ in argumentation iff for every argument for
$ \xbf $ there is
a better or equal argument for $ \xbq.$ It is not sufficient to give just
one
argument, there might not be a best one. We have to consider the $ \xCf
set$ of
all arguments.

Consequently, if $V$ $=$ truth value set $=$ set of arguments with a
partial order
$ \xec,$ we have to look at functions $f:M \xcp \xdp (V),$ where to each
model $m$ (M the model
set) is assigned a set of arguments (which support that
``m belongs to $f$''.) Example: is $m$ a weevil? Yes, it has a long nose.
Yes, it
has articulate antennae \Xl. Thus, $f_{weevil}(m)=\{$ long nose,
articulate antennae $\}.$
We have to define $ \xck $ on $ \xdp (V).$ We think a good definition is:

\bd

$\hspace{0.01em}$


\label{Definition Val-Set-Order}

For $A,B \xcc V$ (V with partial order $ \xec)$ we define:

$A \xck B$ iff $ \xcA a \xbe A \xcE b \xbe B.a \xec b.$

\ed

This seems to be a decent definition of comparison of argument sets.
Why not conversely? Suppose we have a very shaky argument $b$ for $ \xbq
,$ then
to say that arguments for $ \xbq $ are better than arguments for $ \xbf,$
we would need
an even worse argument for $ \xbf.$ This does not seem right.

Thus, we define for arbitrary model functions $f$ and $g:$

\bd

$\hspace{0.01em}$


\label{Definition f-g-comparison}

Let $f,g:M \xcp \xdp (V).$ We say $f$ entails $g$ iff:

$f \xck g$ iff $ \xcA a \xbe f(m) \xcE b \xbe g(m).a \xec b.$

\ed

In total orders, sup and inf were defined.
We want for sup: $A,B \xck sup(A,B),$ and $A,B \xck C$ $ \xch $ $sup(A,B)
\xck C.$
So we define:

\bd

$\hspace{0.01em}$


\label{Definition Arg-Sup}

For a set $ \xda $ of argument sets, define $sup(\xda):= \xcV \xda.$

\ed

\bfa

$\hspace{0.01em}$


\label{Fact Arg-Sup}

We have

(1) For all $A \xbe \xda $ $A \xck sup(\xda).$

(2) If for all $A \xbe \xda $ $A \xck B,$ then $sup(\xda) \xck B.$

\efa

\subparagraph{
Proof
}

$\hspace{0.01em}$


Trivial by definition of $ \xck.$
$ \xcz.$
\\[3ex]

What is the inf? A definition should also work if the order is empty.
Then, $inf(A,B)=A \xcs B,$ which may be empty. This is probably not what
we want.
It is probably best to leave inf undefined.

But we can replace $A \xck inf(B,C)$ by $A \xck B$ and $A \xck C,$ so we
can work with inf
on the right of $ \xck $ without a definition, replacing it by the
universal
quantifier (or, equivalently, by AND).

For interpolation,
for $L=J \xcv J' \xcv J'',$ $f$ insensitive to $J,$ $g$ insensitive to
$J'',$ $f(m) \xck g(m)$ for all
$m \xbe M,$ we looked at
$f^{+}(m_{J' })$ $:=$ $sup\{f(m'):m \xex J' =m' \xex J' \}$ and
$g^{-}(m_{J' })$ $:=$ $inf\{g(m'):m \xex J' =m' \xex J' \}.$
We showed that $f^{+}(m_{J' }) \xck g^{-}(m_{J' }).$

We have to modify and show:

\bfa

$\hspace{0.01em}$


\label{Fact F-Min-G-Max}

$sup\{f(m'):$ $m \xex J' =m' \xex J' \}$ $:=$ $ \xcV \{f(m'):$ $m \xex
J' =m' \xex J' \}$ $ \xck $ $g(m'')$ for all $m'' $ such
that $m \xex J' =m'' \xex J'.$

\efa

\subparagraph{
Proof
}

$\hspace{0.01em}$


By definition of $ \xck,$ it suffices to show that
$ \xcA m' \xcA m'' (m \xex J' =m' \xex J' $ $ \xcu $ $m \xex J' =m'' \xex
J' $ $ \xch $ $f(m') \xck g(m'')).$

Take $m' $ and $m'' $ as above, so $m' \xex J' =m \xex J' =m'' \xex J'.$
Define $m_{0}$ by
$m_{0} \xex J=m'' \xex J,$ $m_{0} \xex J' \xcv J'' =m' \xex J' \xcv J''.$
As $f$ is insensitive to $J,$ $f(m')=f(m_{0}) \xck g(m_{0})$
by prerequisite. Note that $m_{0} \xex J \xcv J' =m'' \xex J \xcv J',$ as
$m_{0} \xex J' =m' \xex J' =m'' \xex J'.$
As $g$ is insensitive to $J'',$ $g(m_{0})=g(m'').$ So we have $f(m'
)=f(m_{0}) \xck g(m_{0})=g(m'').$

$ \xcz $
\\[3ex]

\bfa

$\hspace{0.01em}$


\label{Fact Arg-Int}

$f^{+}(m_{J' })$ is an interpolant for $f$ and $g$ under above
prerequisites.

\efa

\subparagraph{
Proof
}

$\hspace{0.01em}$


Define $h(m):=f^{+}(m_{J' }).$ We have to show that $h$ is an interpolant.
$f(m) \xck h(m)$ is trivial by definition. It remains to show that $h(m)
\xck g(m)$ for
all $m.$ $h(m)$ $:=$ $sup\{f(m'):m \xex J' =m' \xex J' \}$ $ \xck $
$g(m)$ iff
$ \xcA m' (m' \xex J' =m \xex J' $ $ \xch $ $f(m') \xck g(m)),$ but this
is a special case of the proof of
Fact \ref{Fact F-Min-G-Max} (page \pageref{Fact F-Min-G-Max}). $ \xcz $
\\[3ex]

Note that the same approach may also be used in other contexts, e.g.
considering worlds in Kripke structures as truth values, $w \xbe M(\xbf
)$ iff
$w \xbe f_{ \xbf }.$
All we really need is some kind of sup and inf.
\clearpage
\chapter{
Neighbourhood semantics
}

\label{Chapter Neighbourhood}
\section{
Introduction
}

Neighbourhood semantics,
probably first introduced by
D.Scott and R.Montague in
 \cite{Sco70} and  \cite{Mon70}, and already used for deontic
logic by O.Pacheco
in  \cite{Pac07} to avoid unwanted weakening of obligations, seem to
be useful for
many logics:

 \xEh

 \xDH
in preferential logics, they describe the limit variant, where we consider
neighbourhoods of an ideal, usually inexistent, situation,

 \xDH
in approximative reasoning, they describe the approximations to the final
result,

 \xDH
in deontic and default logic, they describe the
``good'' situations, i.e., deontically acceptable, or where defaults have
fired.

 \xEj

Neighbourhood semantics are used, when the ``ideal'' situation does not
exist (e.g., preferential systems without minimal elements), or
are too difficult to obtain (e.g., ``perfect'' deontic states).
\subsection{
Defining neighbourhoods
}

Neighbourhoods can be defined in various ways:

 \xEI
 \xDH
by algebraic systems, like unions of intersections of certain sets
(but not complements),
 \xDH
quality relations, which say that some points are better than others,
carrying
over to sets of points,
 \xDH
distance relations, which measure the distance to the perhaps inexistant
ideal points.
 \xEJ

The relations and distances may be given already by the underlying
structure,
e.g., in preferential structures, or they can be defined in a natural
way, e.g., from a systems of sets, as in deontic logic or default logic.
In these cases, we can define a distance between two points by the number
or set of deontic requirements or default rules which one satisfies, but
not the other. A quality relation is defined in a similar way: a point is
better, if it satisfies more requirements or rules.
\subsection{
Additional requirements
}

With these tools, we can define properties neighbourhoods should
have. E.g., we may require them to be downward closed, i.e., if $x \xbe
N,$
where $N$ is a neighbourhood, $y \xeb x,$ $y$ is better than $x,$ then $y$
should also be
in $N.$ This is a property we will certainly require in neighbourhood
semantics
for preferential structures (in the limit version). For these structures,
we will also require that for every $x \xce N,$ there should be some $y
\xbe N$ with
$y \xeb x.$ We may also require that, if $x \xbe N,$ $y \xce N,$ and $y$
is in some aspect
better than $x,$ then there must be $z \xbe N,$ which is better than both,
so
we have some kind of ``ceteris paribus'' improvement.
\subsection{
Connections between the various properties
}

There is a multitude of possible definitions (via distances, relations,
set systems), and properties, so it is not surprising that one can
investigate a multitude of connections between the different possible
definitions of neighbourhoods. We cannot cover all possible connections,
so
we compare only a few cases, and the reader is invited to complete the
picture
for the cases which interest him. The connections we examined are
presented in
Section \ref{Section Exam} (page \pageref{Section Exam}).
\subsection{
Various uses of neighbourhood semantics
}

We also distinguish the different uses of the systems of sets thus
characterized as neighbourhoods: we can look at all formulas which hold in
(all
or some) such sets (as in neighbourhood semantics for preferential
logics), or
at the formulas which exactly describe them. The latter reading avoids the
infamous Ross paradox of deontic logic.
This distinction is simple, but basic, and did probably not receive
the attention it deserves, in the literature.
\section{
Detailed overview
}

Our starting point was to give the
``derivation'' in deontic systems a precise semantical meaning.
We extend this now to encompass the following situations:

 \xEh
 \xDH
Deontic systems, including contrary-to-duty obligations
 \xDH
Default systems a la Reiter.
 \xDH
The limit version of preferential structures
 \xDH
Approximative logic
 \xEj

We borrow the word ``neighbourhood'' from analysis and topology, but
should be aware that our use will be, at least partly, different.

Common to topology and our domain is that - for reasons to be discussed -
we are not only interested in one ideal point or one ideal set, but in
sets
which are in some sense bigger, and whose elements are in some sense close
to
the ``ideal''.
\subsection{
Motivation
}

What are the reasons to consider some kind of
``approximation''?

 \xEh

 \xDH
First, the ``ideal'' might not exist, e.g.:
 \xEh
 \xDH
In preferential structures,
minimal models are the ideal, but there might be none, due to infinite
descending chains.
So the usual approach via minimal models leads
to inconsistency, we have to take the limit approach,
see Definition \ref{Definition D-8.1.1} (page \pageref{Definition D-8.1.1}).
The same
holds, e.g., for theory revision or counterfactual conditionals
without closest worlds.

 \xDH
Default rules might be contradictory, so the ideal (all defaults are
satisfied) is impossible to obtain.
 \xEj

 \xDH
Second, the ideal might exist, but be too difficult to obtain, e.g.:
 \xEh

 \xDH
In deontic logic, the requirements to lead a perfectly moral life might
just be beyond human power. The same may hold for other imperative
systems.
E.g., we might be obliged to post the letter
and to water the plants, but we have not time for both, so doing one or
the other is certainly better than nothing (so the
``or'' in the Ross paradox is $ \xCf not$ the problem).

 \xDH
It might be too costly to obtain perfect cleanliness, so we have to settle
with sufficiently clean.
 \xDH
Approximate reasoning will try to find better and better answers, perhaps
without hope to find an ideal answer,

 \xEj

 \xDH
Things might be even more complicated by a (partial or total) hierarchy
between aims. E.g., it is a ``stronger'' law not to kill than
not to steal.

 \xEj
\subsection{
Tools to define neighbourhoods
}

To define a suitable notion of neighbourhood, we may have various tools:

 \xEh
 \xDH
We may have a quality relation between points, where $a \xeb b$ says that
$b$ is
in some sense ``better'' than $ \xCf a.$

Such relations are, e.g., given in:

 \xEh

 \xDH
in preferential structures, it is the preference relation

 \xDH
in defaults, a (normal) default gives a quality relation:
the situations which satisfy it, are ``better'' than those which do
not - here, a default gives a quality relation not only to two
situations, but usually between two $ \xCf sets$ of situations -
the same holds for deontic logics and other imperative systems,

 \xDH
We may have borders separating subsets, with a direction, where
it is ``better'' inside (or outside) the border, as in $X-$logic,
see  \cite{BS85}, and  \cite{Sch04}.

 \xDH
in approximation, one situation might be closer than the other
to the ideal.

 \xEj

 \xDH
We may have several such relations, which may also partly contradict each
other, and we may have a relation of importance between different $ \xeb $
and $ \xeb',$
(as in the example of not to kill or not to steal). The better a situation
is,
the closer it should be to our ideal.
 \xDH
We may have a distance relation between points, and these distances may be
partially or totally ordered. With the distance relation, we can measure
the distance of a point to the ideal points (all of them, the closest one,
the most distant ideal point, etc.). Even if the ideal points do not
exist,
we can perhaps find a reasonable measure of distance to them.

This can be found
in distance semantics for theory revision and counterfactuals. There,
it is the ``closeness'' relation, we are interested only in the
closest models, and if they do not exist, in the sufficiently
close ones (the limit approach),

 \xEj
\subsection{
Additional requirements
}

But there might still be other requirements:

 \xEh

 \xDH
We might postulate that neighbourhoods do not only contain all
sufficiently good
points, but also do $ \xCf not$ contain any points which are too bad.

 \xDH
We may require that they are closed under certain operations, e.g.:
 \xEh

 \xDH
If $x$ is in a neighbourhood $X,$ and $y$ better than $x,$ then it should
also be in $X,$
(closure under improvement,
see Definition \ref{Definition Closed} (page \pageref{Definition Closed})).
 \xDH
For all $y$ and any neighbourhood $X,$ there should be some $x \xbe X,$
which is better
than $y.$

(This and the preceeding requirement are those of MISE, see
Definition \ref{Definition D-8.1.1} (page \pageref{Definition D-8.1.1}).)
 \xDH
If we have a notion of distance, and $x,x' $ are in a neighbourhood $X,$
then
anything ``between'' $x$ and $x' $ should be in $X$ $(\xCf X$ is convex).
Thus, when we move
in $X,$ we do not risk to leave $X.$
 \xDH
Similarly, if $x \xbe X,$ and $y$ is an ideal point, then everything
between $x$ and $y$
is in $X.$ Or, if $x \xbe X,$ and $y$ is an ideal point closest to $x,$
then everything
between $x$ and $y$ should be in $X.$ So, when we improve our situation,
we will
not leave the neighbourhood. (A star shaped set around an ideal point may
satisfy this requirement, without being convex.)
(See Definition \ref{Definition Neighbourhood} (page \pageref{Definition
Neighbourhood}).)
 \xDH
If we have to satisfy several requirements, we can ask whether this is
possible independently for those requirements, or if we have to sacrify
one requirement in order to satisfy another. If the latter is the case, is
there a hierarchy of requirements?

 \xDH
Elements in a ``good'' neighbourhood should be better than the others:

If $x \xbe X$ and the closest (to $x)$ $y \xbe \xdC (X),$ then $x \xeb y$
should hold, and,
conversely, if $y \xbe \xdC (X),$ and the closest (to $y)$ $x \xbe X,$
then $x \xeb y$ should hold,
see Definition \ref{Definition Quality-Extension} (page \pageref{Definition
Quality-Extension}),

 \xDH
It might be desirable to improve quality by moving into a good
neighbourhood,
without sacrificing anything achieved already, this is supposed to capture
the
``ceteris paribus'' idea:

If $x \xbe X$ and $y \xce X$ satisfy a set $R$ of rules, then there is $x'
\xbe X,$ which
also satisfies $R,$ and which is better than $y,$
see Definition \ref{Definition Delta-O} (page \pageref{Definition Delta-O}).

 \xDH
Given a set of ``good'' sets, we might be able to construct all
good neighbourhoods by simple algebraic operations:
Any good neighbourhood $X$ is a union of intersections of the ``good'' sets,
see Definition \ref{Definition ui} (page \pageref{Definition ui}),

 \xDH
Finally, if this exists, the set of ideal points should probably satisfy
our criteria, the set of ideal points should be
a ``good'' neighbourhood.

 \xEj

 \xEj

Of particular interest are requirements which are in some
sense independent:

 \xEh

 \xDH
We should try to satisfy an obligation, without violating another
obligation, which was not violated before.

 \xDH
The idea behind the Stalnaker/Lewis semantics for counterfactuals,
(see  \cite{Sta68},  \cite{Lew73}),
is to look at the closest, i.e., minimally changed situations.
``If it were to rain, $ \xfI $ would use an umbrella'' means something like:
``If it were to rain, and there were not a very strong wind'' (there is no
such
wind now), ``if $ \xfI $ had an umbrella'' (I have one now), etc., i.e. if
things
were mostly as they are now, with the exception that now it does not rain,
and in the situation $ \xfI $ speak about it rains, then $ \xfI $ will use
an umbrella.

 \xDH
The distance semantics for theory revision looks also (though with a
slightly different formal approach) at the closest, minimally changed,
situations.

 \xDH
This idea of ``ceteris paribus'' is the attempt to isolate a
necessary change from the rest of the situation, and is thus
intimately related to the concept of independence. Of course,
a minimal change might not be possible, but small enough changes
might do. Consider, e.g., the constant function $f:[0,1] \xcp [0,1],$
$f(x):=0,$ and we look for a minimally changed $ \xCf continous$ function
with $f(0.5):=1.$ This does not exist. So we have to do with
approximation, and look at functions ``sufficiently''
close to the first function. This is one of the reasons we have
to look at the limit variant of theory revision and counterfactuals.

 \xEj

Remark: We do not look here at $ \xCf paths$ which lead (efficiently?)
to better and better situations.
\subsection{
Interpretation of the neighbourhoods
}

Once we have identified our ``good'' neighbourhoods,
we can interpret the result in several ways:

 \xEh
 \xDH
We can ask what holds in $ \xCf all$ good neighbourhoods.
 \xDH
We can ask what holds (finally) in $ \xCf some$ good neighbourhood - this
is the approach for limit preferential structures and similar situations.
 \xDH
We may be not so much interested in what holds in all or some
good neighbourhoods, but to describe them: This is
the problem of the semantics of a system of obligations.
In short: what distinguishes a good from a bad set of
situations.

Such characterization of ``good'' situations
will give us a new semantics not only for deontic logics,
and thus a precise semantical meaning for the
``derivation'' in deontic systems,
see below for a justification, but
also for defaults, preferential structures, etc.
In particular, such descriptions will not necessarily be
closed under arbitrary classical weakening - see the
infamous Ross Paradox,
Example \ref{Example Ross-Paradox} (page \pageref{Example Ross-Paradox}).

 \xEj
\subsection{
Overview of the different lines of reasoning
}

This chapter is conceptually somewhat complicated, therefore we
give now an overview of the different aspects:

 \xEh
 \xDH
We look at different tools and ways to define neighbourhoods, using
distances, quality relations, and perhaps combining them, or purely
algebraic ways like unions, intersections, etc.
 \xDH
We look at additional requirements for neighbourhoods, using such tools,
like closure principles.
 \xDH
We investigate how to obtain such natural relations, distances, etc.
from different structures, e.g., from obligations, defaults, preferential
models, etc.
 \xDH
We look at various possibilities to interpret the neighbourhood systems we
have constructed, e.g., we can ask what holds in all or some
neighbourhoods,
what finally holds in neighbourhoods (when we have a grading of the
neighbourhoods), or what characterizes the neighbourhoods, e.g., in the
case of deontic logic.
 \xDH
We conclude (unsystematically) with connections between the different
concepts.
 \xEj
\subsection{
Extensions
}

This might be the place to make a remark on extensions.

An extension is, roughly, a maximal consistent set of information, or, a
smallest non-empty set of models. In default logic, we can follow
contradictory
default to the end of reasoning, and obtain perhaps vcontradictory
information,
likewise in inheritance nets, etc.
Usually, one then takes the intersection of extensions, what is true in
all
extensions, which is - provided the laguage is adequate - the
``OR'' of the extensions.

But we can also see preferential structures as resulting in extensions,
where every minimal model is an extension:

Consider a preferential structure with 4 models, say pq, $p \xCN q,$ $
\xCN pq,$ $ \xCN p \xCN q,$
ordered by $pq \xeb p \xCN q,$ $ \xCN pq \xeb p \xCN q.$ Then we can see
the relation roughly as two
defaults: $p \xCN q:pq,$ and $p \xCN q: \xCN pq,$ with two extensions: pq
and $ \xCN pq.$
So, we can see a preferential structure as having usually many extensions
(unless there is a single best model, of course), and we take as
result the intersection of extensions, i.e., the theory which holds
in $ \xCf all$ minimal models.

In preferential structures, the construction of the set of minimal models
is a one-step process: a model is in or out. In defaults, for instance,
the construction is more complicated, we branch in the process. This is
what
may make the construction problematic, and gives rise to different
approaches
like taking immediately the intersection of extensions in inheritance
networks, etc. But this difference to preferential structures is in the
$ \xCf process$ of construction, it is not in the outcome.

These questions are intimately related to our neighbourhood semantics,
as the constructions can be seen as an approximation to the ideal, the
final outcome.
\section{
Tools and requirements for neighbourhoods and how to obtain them
}
\subsection{
Tools to define neighbourhoods
}
\paragraph{
Background \\[2mm]
}

We often work with an additional structure, some $ \xdo \xcc \xdp (U),$
where $U$ is
the universe (intuitively, $U$ is a set of propositional models), which
allows to define distances and quality relations in a natural way.
Intuitively, $ \xdo $ is a base set of
``good'' sets, from which we will construct other
``good'' sets.

Basically, $x$ is better than $y,$ iff $x$ is in more (as a set or by
counting)
$O \xbe \xdo $ than $y$ is, and
the distance between $x$ and $y$ is the set (or cardinality) of $O \xbe
\xdo $ where
$x \xbe O,$ $y \xce O,$ or vice versa.

Sometimes, it is more appropriate to work with sequences of 0/1, where
1 stands for $O,$ 0 for $ \xdC (O)$ for $O \xbe \xdo.$

Thus, we work with sets $ \xbS $ of sequences. Note that $ \xbS $ need not
contain all
possible sequences, corresponding to the possibility that, e.g., $O \xcs
O' = \xCQ $
for $O,O' \xbe \xdo.$

Moreover, we may have a difference in quality between $O$ and $ \xdC (O):$
if
$O$ is an obligation, then $x \xbe O$ is - at least for this obligation -
better than
$x' \xce O.$ The same holds for defaults of the type $: \xbf / \xbf,$
with $O=M(\xbf).$
We will follow the tradition of preferential structures, and
``smaller'' will mean ``better''.
\subsubsection{
Algebraic tools
}

Let here again $ \xdo \xcc \xdp (U).$

\bd

$\hspace{0.01em}$


\label{Definition L-And}

Given a finite propositional laguage $ \xdl $ defined by the set $v(\xdl
)$ of
propositional
variables, let $ \xdl_{ \xcu }$ be the set of all consistent conjunctions
of
elements from $v(\xdl)$ or their negations. Thus, $p \xcu \xCN q \xbe
\xdl_{ \xcu }$ if $p,q \xbe v(\xdl),$ but
$p \xco q,$ $ \xCN (p \xcu q) \xce \xdl_{ \xcu }.$ Finally, let $ \xdl_{
\xco \xcu }$ be the set of all (finite)
disjunctions of formulas from $ \xdl_{ \xcu }.$ (As we will later not
consider all
formulas from $ \xdl_{ \xcu },$ this will be a real restriction.)

Given a set of models $M$ for a finite language $ \xdl,$ define
$ \xbf_{M}$ \index{$ \xbf_{M}$}
$:= \xcU \{p \xbe v(\xdl): \xcA m \xbe M.m(p)=v\} \xcu \xcU \{ \xCN p:p
\xbe v(\xdl), \xcA m \xbe M.m(p)=f\} \xbe \xdl_{ \xcu }.$
(If there are no such $p,$ set $ \xbf_{M}:=TRUE.)$

This is the strongest $ \xbf \xbe \xdl_{ \xcu }$ which holds in $M.$

\ed

\bd

$\hspace{0.01em}$


\label{Definition ui}

$X \xcc U' $ is $ \xCf (ui)$ (for union of intersections) iff there is a
family $ \xdo_{i} \xcc \xdo,$
 \index{(ui)}
$i \xbe I$ such that $X=(\xcV \{ \xcS \xdo_{i}:i \xbe I\}) \xcs U'.$

\ed

Unfortunately, as we will see later, this definition is not very useful
for
simple relativization.

\bd

$\hspace{0.01em}$


\label{Definition Validity}

Let $ \xdo' \xcc \xdo.$ Define for $m \xbe U$ and $ \xbd: \xdo' \xcp
2=\{0,1\}$

$m \xcm \xbd $ $: \xcj $ $ \xcA O \xbe \xdo' (m \xbe O \xcj \xbd (O)=1)$

\ed

\bd

$\hspace{0.01em}$


\label{Definition Independence}

$ \xdo $ is independent \index{independent}  iff $ \xcA \xbd: \xdo \xcp 2.
\xcE m \xbe U.m \xcm \xbd.$

\ed

Obviously, independence does not inherit downward to subsets of $U.$

\bd

$\hspace{0.01em}$


\label{Definition Delta-O}

$ \xdd (\xdo)$ \index{$ \xdd (\xdo)$}  $:=\{X \xcc U':$ $ \xcA \xdo' \xcc
\xdo $ $ \xcA \xbd: \xdo' \xcp 2$

$ \xDC $ $((\xcE m,m' \xbe U,$ $m,m' \xcm \xbd,$ $m \xbe X,m' \xce X)$ $
\xch $ $(\xcE m'' \xbe X.m'' \xcm \xbd \xcu m'' \xeb_{s}m'))\}$

This property expresses that we can satisfy obligations independently: If
we
respect $O,$ we can, in addition,
respect $O',$ and if we are hopeless kleptomaniacs, we may still not be a
murderer. If $X \xbe \xdd (\xdo),$ we can go from $U-X$
into $X$ by improving on all $O \xbe \xdo,$ which we have not fixed by $
\xbd,$ if $ \xbd $ is
not too rigid.
\subsubsection{
Relations
}

\ed

We may have an abstract relation $ \xec $ of quality on the domain, but we
may also
define it from the structure of the sequences, as we will do now.

\bd

$\hspace{0.01em}$


\label{Definition Quality}

Consider the case of sequences.

Given a relation $ \xec $ (of quality) on the codomain, we extend this to
sequences
in $ \xbS:$

$x \xCq y$ \index{$x \xCq y$}  $: \xcj $ $ \xcA i \xbe I(x(i) \xCq y(i))$

$x \xec y$ \index{$x \xec y$}  $: \xcj $ $ \xcA i \xbe I(x(i) \xec y(i))$

$x \xeb y$ \index{$x \xeb y$}  $: \xcj $ $ \xcA i \xbe I(x(i) \xec y(i))$
and $ \xcE i \xbe I(x(i) \xeb y(i))$

In the $ \xbe -$case, we will consider $x \xbe i$ better than $x \xce i.$
As we have only two
values, true/false, it is easy to count the positive and negative cases
(in more complicated situations, we might be able to multiply), so we have
an analogue of the two Hamming distances, which we might call the
Hamming quality \index{Hamming quality}  relations.

Let $ \xdo \xcc \xdp (U)$ be given now.

(Recall that we follow the preferential tradition, ``smaller'' will
mean ``better''.)

$x \xCq_{s}y$ $: \xcj $ $ \xdo (x)= \xdo (y),$

$x \xec_{s}y$ $: \xcj $ $ \xdo (y) \xcc \xdo (x),$

$x \xeb_{s}y$ $: \xcj $ $ \xdo (y) \xcb \xdo (x),$

$x \xCq_{c}y$ $: \xcj $ $card(\xdo (x))=card(\xdo (y)),$

$x \xec_{c}y$ $: \xcj $ $card(\xdo (y)) \xck card(\xdo (x)),$

$x \xeb_{c}y$ $: \xcj $ $card(\xdo (y))<card(\xdo (x)).$
\subsubsection{
Distances
}

\label{Section Ham-Neigh}

\ed

Note that we defined Hamming relations already in
Section \ref{Section Ham-Rel-Dist} (page \pageref{Section Ham-Rel-Dist}), as
announced in
Section \ref{Section Ham-Rel-Size} (page \pageref{Section Ham-Rel-Size}).

$ \xCO $

\label{Hamming}


\index{Hamming distance}

\bd

$\hspace{0.01em}$


\label{Definition Hamming-Distance}

Given $x,y \xbe \xbS,$ a set of sequences over an index set $I,$
the Hamming distance \index{Hamming distance}
comes in two flavours:

$d_{s}(x,y)$ \index{$d_{s}(x,y)$}  $:=\{i \xbe I:x(i) \xEd y(i)\},$ the set
variant,

$d_{c}(x,y)$ \index{$d_{c}(x,y)$}  $:=card(d_{s}(x,y)),$ the counting variant.

We define $d_{s}(x,y) \xck d_{s}(x',y')$ iff $d_{s}(x,y) \xcc d_{s}(x'
,y'),$

thus, $s-$distances are not always comparabel.
Consequently, readers should be aware that $d_{s}$-values are $ \xCf not$
always comparable, even though $<$ and $ \xck $ may suggest a linear
order.
We use these symbols to be in line with other distances.

There are straightforward generalizations of the counting variant:

We can also give different importance to different $i$ in the counting
variant, so e.g.,
$d_{c}(\xBc x,x'  \xBe, \xBc y,y'  \xBe)$ might be 1 if $x \xEd y$ and $x'
=y',$ but 2 if
$x=y$ and $x' \xEd y'.$

If the $x \xbe \xbS $ may have more than 2 different values, then a
varying
individual distance may also reflect to the distances in $ \xbS.$
So, (for any distance $d)$ if $d(x(i),x' (i))<d(x(i),x'' (i)),$ then (the
rest being
equal), we may have $d(x,x')<d(x,x'').$

\vspace{3mm}


\vspace{3mm}

\ed

\bfa

$\hspace{0.01em}$


\label{Fact Hamming-Distance}

(1) If the $x \xbe \xbS $ have only 2 values, say TRUE and FALSE, then
$d_{s}(x,y)=\{i \xbe I:x(i)=TRUE\} \xeY \{i \xbe I:y(i)=TRUE\},$ where $
\xeY $ is the symmetric
set difference.

(2) $d_{c}$ has the normal addition, set union takes the role of addition
for $d_{s},$
$ \xCQ $ takes the role of 0 for $d_{s},$
both are distances in the following sense:

(2.1) $d(x,y)=0$ iff $x=y,$

(2.2) $d(x,y)=d(y,x),$

(2.3) the triangle inequality holds, for the set variant in the form
$d_{s}(x,z) \xcc d_{s}(x,y) \xcv d_{s}(y,z).$

\efa

\subparagraph{
Proof
}

$\hspace{0.01em}$


(2.3) If $i \xce d_{s}(x,y) \xcv d_{s}(y,z),$ then $x(i)=y(i)=z(i),$ so
$x(i)=z(i)$ and
$i \xce d_{s}(x,z).$

The others are trivial.

$ \xcz $
\\[3ex]

$ \xCO $

Recall that the $ \xbs \xbe \xbS $ will often stand for a sequence of
possibilities $O/ \xdC (O)$ with $O \xbe \xdo.$ Thus, the distance
between two
such sequences $ \xbs $ and $ \xbs' $ is the number or set of $O,$ where
$ \xbs $ codes being in $O$ and $ \xbs' $ codes being in $ \xdC (O),$ or
vice versa.

\br

$\hspace{0.01em}$


\label{Remark Hamming-Equivalence}

If the $x(i)$ are equivalence classes, one has to be careful not to
confound the distance between the classes and the resulting distance
between
elements of the classes, as two different elements in the same class have
distance 0. So in Fact \ref{Fact Hamming-Distance} (page \pageref{Fact
Hamming-Distance})  2.1
only one direction holds.

\er

\bd

$\hspace{0.01em}$


\label{Definition Between}

(1) We can define for any distance $d$ with some minimal requirements a
notion of
``between''.
 \index{between}

If the codomain of $d$ has an ordering $ \xck,$ but no addition, we
define:

$ \xBc x,y,z \xBe _{d}$ \index{$ \xBc x,y,z \xBe _{d}$}  $: \xcj $ $d(x,y) \xck
d(x,z)$ and $d(y,z)
\xck d(x,z).$

If the codomain has a commutative addition, we define

$ \xBc x,y,z \xBe _{d}$ \index{$ \xBc x,y,z \xBe _{d}$}  $: \xcj $
$d(x,z)=d(x,y)+d(y,z)$ - in
$d_{s}$ $+$ will be replaced by $ \xcv,$ i.e.

$ \xBc x,y,z \xBe _{s}$ \index{$ \xBc x,y,z \xBe _{s}$}  $: \xcj $
$d(x,z)=d(x,y) \xcv d(y,z).$

For above two Hamming distances, we will write
$ \xBc x,y,z \xBe _{s}$ and $ \xBc x,y,z \xBe _{c}$ \index{$ \xBc x,y,z \xBe
_{c}$}.

(2) We further define:

$[x,z]_{d}$ \index{$[x,z]_{d}$}  $:=\{y \xbe X: \xBc x,y,x \xBe _{d}\}$ - where
$X$ is
the set we work in.

We will write $[x,z]_{s}$ \index{$[x,z]_{s}$}
and $[x,z]_{c}$ \index{$[x,z]_{c}$}  when appropriate.

(3) For $x \xbe U,$ $X \xcc U$ set $x \xFO_{d}X$ \index{$x \xFO_{d}X$}
$:=$ $\{x' \xbe X: \xCN \xcE x'' \xEd x' \xbe X.d(x,x') \xcg d(x,x''
)\}.$

Note that, if $X \xEd \xCQ,$ then $x \xFO X \xEd \xCQ.$

We omit the index when this does not cause confusion. Again, when
adequate,
we write $ \xFO_{s}$ \index{$ \xFO_{s}$}
and $ \xFO_{c}$ \index{$ \xFO_{c}$}.

\ed

For problems with characterizing ``between'' see
 \cite{Sch04}.

\bfa

$\hspace{0.01em}$


\label{Fact Between}

(0) $ \xBc x,y,z \xBe _{d}$ $ \xcj $ $ \xBc z,y,x \xBe _{d}.$

Consider the situation of a set of sequences $ \xbS,$ with $ \xbs:I \xcp
S$ for $ \xbs \xbe \xbS $

Let $A:=A_{ \xbs, \xbs'' }:=\{ \xbs': \xcA i \xbe I(\xbs (i)= \xbs''
(i) \xcp \xbs' (i)= \xbs (i)= \xbs'' (i))\}.$
Then

(1) If $ \xbs' \xbe A,$ then $d_{s}(\xbs, \xbs'')=d_{s}(\xbs, \xbs
') \xcv d_{s}(\xbs', \xbs''),$ so $ \xBc  \xbs, \xbs', \xbs''  \xBe
_{s}.$

(2) If $ \xbs' \xbe A$ and $S$ consists of 2 elements (as in classical
2-valued
logic), then $d_{s}(\xbs, \xbs')$ and $d_{s}(\xbs', \xbs'')$ are
disjoint.

(3) $[ \xbs, \xbs'' ]_{s}=A.$

(4) If, in addition, $S$ consists of 2 elements, then $[ \xbs, \xbs''
]_{c}=A.$

\efa

\subparagraph{
Proof
}

$\hspace{0.01em}$


(0) Trivial.

(1) ``$ \xcc $'' follows from Fact \ref{Fact Hamming-Distance} (page
\pageref{Fact Hamming-Distance}),
(2.3).

Conversely, if e.g. $i \xbe d_{s}(\xbs, \xbs'),$ then
by prerequisite $i \xbe d_{s}(\xbs, \xbs'').$

(2) Let $i \xbe d_{s}(\xbs, \xbs') \xcs d_{s}(\xbs', \xbs''),$
then $ \xbs (i) \xEd \xbs' (i)$ and $ \xbs' (i) \xEd \xbs'' (i),$
but then by $card(S)=2$ $ \xbs (i)= \xbs'' (i),$ but $ \xbs' \xbe A,$
$contradiction.$

We turn to (3) and (4):

If $ \xbs' \xce A,$ then there is $i' $ such that $ \xbs (i')= \xbs''
(i') \xEd \xbs' (i').$ On the other hand,
for all $i$ such that $ \xbs (i) \xEd \xbs'' (i)$ $i \xbe d_{s}(\xbs,
\xbs') \xcv d_{s}(\xbs', \xbs'').$ Thus:

(3) By (1) $ \xbs' \xbe A$ $ \xch $ $ \xBc  \xbs, \xbs', \xbs''  \xBe _{s}.$
Suppose $ \xbs' \xce A,$ so there is $i' $ such that
$i' \xbe d_{s}(\xbs, \xbs')-d_{s}(\xbs, \xbs''),$ so
$ \xBc  \xbs, \xbs', \xbs''  \xBe _{s}$ cannot be.

(4) By (1) and (2) $ \xbs' \xbe A$ $ \xch $
$ \xBc  \xbs, \xbs', \xbs''  \xBe _{c}.$ Conversely, if $ \xbs' \xce A,$ then
$card(d_{s}(\xbs, \xbs'))+card(d_{s}(\xbs', \xbs'')) \xcg
card(d_{s}(\xbs, \xbs''))+2.$

$ \xcz $
\\[3ex]
\subsection{
Obtaining such tools
}

We consider a set of sequences $ \xbS,$ for $x \xbe \xbS $ $x:I \xcp S,$
$I$ a finite index set, $S$
some set.
Often, $S$ will be $\{0,1\},$ $x(i)=1$ will mean that $x \xbe i,$ when $I
\xcc \xdp (U)$ and $x \xbe U.$
For abbreviation, we will call this (unsystematically, often context will
tell) the $ \xbe -$case \index{$ \xbe -$case}.
Often, $I$ will be written $ \xdo,$ intuitively, $O \xbe \xdo $ is then
an obligation, and
$x(O)=1$ means $x \xbe O,$ or $x$ ``satisfies'' the obligation $O.$

\bd

$\hspace{0.01em}$


\label{Definition O-x}

In the $ \xbe -$case, set $ \xdo (x):=\{O \xbe \xdo:x \xbe O\}.$
\subsection{
Additional requirements for neighbourhoods
}

\ed

\bd

$\hspace{0.01em}$


\label{Definition Closed}

Given any relation $ \xec $ (of quality), we say that $X \xcc U$ is
(downward) closed (with
respect to $ \xec)$ iff $ \xcA x \xbe X \xcA y \xbe U(y \xec x$ $ \xch $
$y \xbe X).$

\ed

(Warning, we follow the preferential tradition, ``smaller'' will
mean ``better''.)

\bfa

$\hspace{0.01em}$


\label{Fact Subset-Closure}

Let $ \xec $ be given.

(1) Let $D \xcc U' \xcc U'',$ $D$ closed in $U'',$ then $D$ is also
closed in $U'.$

(2) Let $D \xcc U' \xcc U'',$ $D$ closed in $U',$ $U' $ closed in $U''
,$ then $D$ is
closed in $U''.$

(3) Let $D_{i} \xcc U' $ be closed for all $i \xbe I,$ then so are $ \xcV
\{D_{i}:i \xbe I\}$ and $ \xcS \{D_{i}:i \xbe I\}.$

\efa

\subparagraph{
Proof
}

$\hspace{0.01em}$


(1) Trivial.

(2) Let $x \xbe D \xcc U',$ $x' \xec x,$ $x' \xbe U'',$ then $x' \xbe U'
$ by closure of $U'',$ so $x' \xbe D$ by
closure of $U'.$

(3) Trivial.

$ \xcz $
\\[3ex]

\bd

$\hspace{0.01em}$


\label{Definition Quality-Extension}

Given a quality relation $ \xeb $ between elements, and a distance $d,$ we
extend the
quality relation to sets and define:

(1) $x \xeb Y$ \index{$x \xeb Y$}  $: \xcj $ $ \xcA y \xbe (x \xFO Y).x \xeb
y.$ (The closest elements - i.e. there are no
closer ones - of $Y,$ seen from $x,$ are less good than $x.)$

analogously $X \xeb y$ $: \xcj $ $ \xcA x \xbe (y \xFO X).x \xeb y$

(2) $X \xeb_{l}Y$ \index{$X \xeb_{l}Y$}  $: \xcj $ $ \xcA x \xbe X.x \xeb Y$
and $ \xcA y \xbe Y.X \xeb y$
(X is locally better than $Y).$

When necessary, we will write $ \xeb_{l,s}$ \index{$ \xeb_{l,s}$}
or $ \xeb_{l,c}$ \index{$ \xeb_{l,c}$}  to distinguish the
set from the counting variant.

For the next definition, we use the notion of size: $ \xeA \xbf $ iff for
almost all $ \xbf $
holds i.e. the set of exceptions is small.

(3) $X \xDc_{l}Y$ \index{$X \xDc_{l}Y$}  $: \xcj $ $ \xeA x \xbe X.x \xeb Y$
and $ \xeA y \xbe Y.X \xeb y.$

We will likewise write $ \xDc_{l,s}$ \index{$ \xDc_{l,s}$}  etc.

This definition is supposed to capture quality difference under minimal
change,
the ``ceteris paribus \index{ceteris paribus}'' idea:
$X \xeb_{l} \xdC X$ should hold for an obligation $X.$
Minimal change is coded by $ \xFO,$ and ``ceteris paribus'' by minimal
change.

\ed

\bfa

$\hspace{0.01em}$


\label{Fact General-Obligation}

If $X \xeb_{l} \xdC X,$ and $x \xbe U$ an optimal point (there is no
better one), then $x \xbe X.$

\efa

\subparagraph{
Proof
}

$\hspace{0.01em}$


If not, then take $x' \xbe X$ closest to $x,$ this must be better than
$x,$ contradiction.
$ \xcz $
\\[3ex]

\bd

$\hspace{0.01em}$


\label{Definition Neighbourhood}

Given a distance $d,$ we define:

(1) Let $X \xcc Y \xcc U',$ then $Y$ is a neighbourhood of $X$ in $U' $
iff

$ \xcA y \xbe Y \xcA x \xbe X(x$ is closest to $y$ among all $x' $ with
$x' \xbe X$ $ \xch $ $[x,y] \xcs U' \xcc Y).$

(Closest means that there are no closer ones.)

When we also have a quality relation $ \xeb,$ we define:

(2) Let $X \xcc Y \xcc U',$ then $Y$ is an improving neighbourhood of $X$
in $U' $ iff

$ \xcA y \xbe Y \xcA x((x$ is closest to $y$ among all $x' $ with $x' \xbe
X$ and $x' \xec y)$ $ \xch $ $[x,y] \xcs U' \xcc Y).$

When necessary, we will have to say for (3) and (4) which variant, i.e.
set or counting, we mean.

\ed

\bd

$\hspace{0.01em}$


\label{Definition Ham-Umgeb}

Given a Hamming distance and a Hamming relation,
$X$ is called a Hamming neighbourhood of the best cases iff for any $x
\xbe X$
and $y$ a best case with minimal distance from $x,$ all elements between
$x$ and $y$ are in $X.$

\ed

\bfa

$\hspace{0.01em}$


\label{Fact Neighbourhood}

(1) If $X \xcc X' \xcc \xbS,$ and $d(x,y)=0$ $ \xch $ $x=y,$ then $X$ and
$X' $ are Hamming neighbourhoods
of $X$ in $X'.$

(2) If $X \xcc Y_{j} \xcc X' \xcc \xbS $ for $j \xbe J,$ and all $Y_{j}$
are Hamming Neighbourhoods
of $X$ in $X',$ then so are $ \xcV \{Y_{j}:j \xbe J\}$ and $ \xcS
\{Y_{j}:j \xbe J\}.$

\efa

\subparagraph{
Proof
}

$\hspace{0.01em}$


(1) is trivial (we need here that $d(x,y)=0$ $ \xch $ $x=y).$

(2) Trivial.

$ \xcz $
\\[3ex]
\subsection{
Connections between the various concepts
}

\label{Section Exam}

\bfa

$\hspace{0.01em}$


\label{Fact Hamming-Neighbourhood}

If $x,y$ are models, then $[x,y]=M(\xbf_{\{x,y\}}).$
(See Definition \ref{Definition Between} (page \pageref{Definition Between}) 
and
Definition \ref{Definition L-And} (page \pageref{Definition L-And}).)

\efa

\subparagraph{
Proof
}

$\hspace{0.01em}$


$m \xbe [x,y]$ $ \xcj $ $ \xcA p(x \xcm p,y \xcm p \xch m \xcm p$ and $x
\xcM p,y \xcM p \xch m \xcM p),$
$m \xcm \xbf_{\{x,y\}}$ $ \xcj $ $m \xcm \xcU \{p:x(p)=y(p)=v\} \xcu \xcU
\{ \xCN p:x(p)=y(p)=f\}.$
$ \xcz $
\\[3ex]

The requirement of closure causes a problem for
the counting approach: Given e.g. two obligations $O,$ $O',$ then any two
elements
in just one obligation have the same quality, so if one is in, the other
should
be, too. But this prevents now any of the original obligations to have the
desirable property of closure. In the counting case, we will obtain a
ranked
structure, where elements satisfy
0, 1, 2, etc. obligations, and we are unable to differentiate inside those
layers. Moreover, the set variant seems to be closer to logic, where we do
not
count the propositional variables which hold in a model, but consider them
individually. For these reasons, we will not pursue the counting approach
as
systematically as the set approach. One should, however, keep in mind that
the
counting variant gives a ranking relation of quality, as all qualities are
comparable, and the set variant does not. A ranking seems to be
appreciated
sometimes in the literature, though we are not really sure why.

Of particular interest is the combination of $d_{s}$ and $ \xec_{s}$
$(d_{c}$ and $ \xec_{c})$
respectively - where by $ \xec_{s}$ we also mean $ \xeb_{s}$ and $
\xCq_{s},$ etc. We turn to
this now.

\bfa

$\hspace{0.01em}$


\label{Fact Quality-Distance}

We work in the $ \xbe -$case.

(1) $x \xec_{s}y$ $ \xch $ $d_{s}(x,y)= \xdo (x)- \xdo (y)$

Let $a \xeb_{s}b \xeb_{s}c.$ Then

(2) $d_{s}(a,b)$ and $d_{s}(b,c)$ are not comparable,

(3) $d_{s}(a,c)=d_{s}(a,b) \xcv d_{s}(b,c),$ and thus $b \xbe [a,c]_{s}.$

This does not hold in the counting variant, as Example \ref{Example Count} (page
\pageref{Example Count})  shows.

(4) Let $x \xeb_{s}y$ and $x' \xeb_{s}y$ with $x,x' $ $
\xeb_{s}-$incomparabel. Then $d_{s}(x,y)$ and $d_{s}(x',y)$
are incomparable.

(This does not hold in the counting variant, as then all distances are
comparable.)

(5) If $x \xeb_{s}z,$ then for all $y \xbe [x,z]_{s}$ $x \xec_{s}y
\xec_{s}z.$

\efa

\subparagraph{
Proof
}

$\hspace{0.01em}$


(1) Trivial.

(2) We have $ \xdo (c) \xcb \xdo (b) \xcb \xdo (a),$ so the results
follows from (1).

(3) By definition of $d_{s}$ and (1).

(4) $x$ and $x' $ are $ \xec_{s}$-incomparable, so there are $O \xbe \xdo
(x)- \xdo (x'),$
$O' \xbe \xdo (x')- \xdo (x).$

As $x,x' \xeb_{s}y,$ $O,O' \xce \xdo (y),$ so $O \xbe d_{s}(x,y)-d_{s}(x'
,y),$
$O' \xbe d_{s}(x',y)-d_{s}(x,y).$

(5) $x \xeb_{s}z$ $ \xch $ $ \xdo (z) \xcb \xdo (x),$ $d_{s}(x,z)= \xdo
(x)- \xdo (z).$ By prerequisite
$d_{s}(x,z)=d_{s}(x,y) \xcv d_{s}(y,z).$
Suppose $x \xeC_{s}y.$ Then there is $i \xbe \xdo (y)- \xdo (x) \xcc
d_{s}(x,y),$ so
$i \xce \xdo (x)- \xdo (z)=d_{s}(x,z),$ $contradiction.$

Suppose $y \xeC_{s}z.$ Then there is $i \xbe \xdo (z)- \xdo (y) \xcc
d_{s}(y,z),$ so
$i \xce \xdo (x)- \xdo (z)=d_{s}(x,z),$ $contradiction.$

$ \xcz $
\\[3ex]

\be

$\hspace{0.01em}$


\label{Example Count}

In this and similar examples, we will use the model notation. Some
propositional variables $p,$ $q,$ etc. are given, and models are described
by
$p \xCN qr,$ etc. Moreover, the propositional variables are the
obligations, so
in this example we have the obligations $p,$ $q,$ $r.$

Consider $x:= \xCN p \xCN qr,$ $y:=pq \xCN r,$ $z:= \xCN p \xCN q \xCN r.$
Then $y \xeb_{c}x \xeb_{c}z,$ $d_{c}(x,y)=3,$
$d_{c}(x,z)=1,$ $d_{c}(z,y)=2,$ so $x \xce [y,z]_{c}.$ $ \xcz $
\\[3ex]

\ee

\bfa

$\hspace{0.01em}$


\label{Fact Local-Global}

Take the set version.

If $X \xeb_{l,s} \xdC X,$ then $X$ is downward $ \xeb_{s}$-closed.

\efa

\subparagraph{
Proof
}

$\hspace{0.01em}$


Suppose $X \xeb_{l,s} \xdC X,$ but $X$ is not downward closed.

Case 1: There are $x \xbe X,$ $y \xce X,$ $y \xCq_{s}x.$ Then $y \xbe x
\xFO_{s} \xdC X,$ but $x \xeB y,$ $contradiction.$

Case 2:
There are $x \xbe X,$ $y \xce X,$
$y \xeb_{s}x.$ By $X \xeb_{l,s} \xdC X,$ the elements in $X$ closest to
$y$ must be better than $y.$
Thus, there is $x' \xeb_{s}y,$ $x' \xbe X,$ with minimal distance from
$y.$ But then
$x' \xeb_{s}y \xeb_{s}x,$ so $d_{s}(x',y)$ and $d_{s}(y,x)$ are
incomparable by
Fact \ref{Fact Quality-Distance} (page \pageref{Fact Quality-Distance}), so $x$
is among those
with minimal distance from $y,$ so $X \xeb_{l,s} \xdC X$ does not hold. $
\xcz $
\\[3ex]

\be

$\hspace{0.01em}$


\label{Example Dependent-2}

We work with the set variant.

This example shows that $ \xec_{s}-$closed does not imply $X \xeb_{l,s}
\xdC X,$ even
if $X$ contains the best elements.

Let $ \xdo:=\{p,q,r,s\},$ $U':=\{x:=p \xCN q \xCN r \xCN s,$ $y:= \xCN
pq \xCN r \xCN s,$ $x':=pqrs\},$ $X:=\{x,x' \}.$
$x' $ is the best element of $U',$ so $X$ contains the best elements,
and $X$ is downward closed in $U',$
as $x$ and $y$ are not comparable. $d_{s}(x,y)=\{p,q\},$ $d_{s}(x'
,y)=\{p,r,s\},$ so the
distances from $y$ are not comparable, so $x$ is among the closest
elements in $X,$
seen from $y,$ but $x \xeB_{s}y.$

The lack of comparability is essential here, as the following Fact shows.

$ \xcz $
\\[3ex]

\ee

We have, however, for the counting variant:

\bfa

$\hspace{0.01em}$


\label{Fact Count-Closed}

Consider the counting variant. Then

If $X$ is downward closed, then $X \xeb_{l,c} \xdC X.$

\efa

\subparagraph{
Proof
}

$\hspace{0.01em}$


Take any $x \xbe X,$ $y \xce X.$ We have $y \xec_{c}x$ or $x \xeb_{c}y,$
as any two elements are
$ \xec_{c}-$comparabel. $y \xec_{c}x$ contradicts closure, so $x
\xeb_{c}y,$ and $X \xeb_{l,c} \xdC X$ holds
trivially. $ \xcz $
\\[3ex]

\be

$\hspace{0.01em}$


\label{Example Not-Global}

Let $U':=\{x,x',y,y' \}$ with $x':=pqrs,$ $y':=pqr \xCN s,$ $x:= \xCN
p \xCN qr \xCN s,$ $y:= \xCN p \xCN q \xCN r \xCN s.$

Consider $X:=\{x,x' \}.$

The counting version:

Then $x' $ has quality 4 (the best), $y' $ has quality 3, $x$ has 1, $y$
has 0.

$d_{c}(x',y')=1,$ $d_{c}(x,y)=1,$ $d_{c}(x,y')=2.$

\ee

Then above ``ceteris paribus'' criterion is satisfied, as $y' $ and $x$ do
not
``see'' each other, so
$X \xeb_{l,c} \xdC X.$

But $X$ is not downward closed, below $x \xbe X$ is a better element $y'
\xce X.$

This seems an argument against $X$ being an obligation.

The set version:

We still have $x' \xeb_{s}y' \xeb_{s}x \xeb_{s}y.$ As shown in
Fact \ref{Fact Quality-Distance} (page \pageref{Fact Quality-Distance}),
$d_{s}(x,y)$ (and also
$d_{s}(x',y'))$ and $d_{s}(x,y')$ are not comparable, so our argument
collapses.

As a matter of fact, we have the result that the ``ceteris paribus''
criterion entails downward closure in the set variant, see
Fact \ref{Fact Local-Global} (page \pageref{Fact Local-Global}).

$ \xcz $
\\[3ex]

In the following
Section \ref{Section Not-Indep} (page \pageref{Section Not-Indep})  and
Section \ref{Section Indep} (page \pageref{Section Indep}),
we will assume a set $ \xdo $ of obligations to be given.
We define the relation $ \xeb:= \xeb_{ \xdo }$ as described in
Definition \ref{Definition Quality} (page \pageref{Definition Quality}), and
the distance $d$ is the
Hamming distance \index{Hamming distance}
based on $ \xdo,$ see
Definition \ref{Definition Hamming-Distance} (page \pageref{Definition
Hamming-Distance}).

We work here mostly in the set version, the $ \xbe -$case, only
in the final
Section \ref{Section Counting-Case} (page \pageref{Section Counting-Case}), we
will look at the counting case.
\subsubsection{
The not necessarily independent case
}

\label{Section Not-Indep}

\be

$\hspace{0.01em}$


\label{Example Dependent-1}

Work in the set variant. We show that $X$ $ \xec_{s}$-closed does not
necessarily imply that $X$ contains all $ \xec_{s}$-best elements.

Let $ \xdo:=\{p,q\},$ $U':=\{p \xCN q, \xCN pq\},$ then all elements of
$U' $ have best quality
in $U',$ $X:=\{p \xCN q\}$ is closed, but does not contain all best
elements. $ \xcz $
\\[3ex]

\ee

\be

$\hspace{0.01em}$


\label{Example Dependent-3}

Work in the set variant. We show that $X$ $ \xec_{s}$-closed does not
necessarily imply that $X$ is a neighbourhood of the best elements, even
if $X$ contains them.

Consider $x:=pq \xCN rstu,$ $x':= \xCN pqrs \xCN t \xCN u,$ $x'':=p \xCN
qr \xCN s \xCN t \xCN u,$
$y:=p \xCN q \xCN r \xCN s \xCN t \xCN u,$ $z:=pq \xCN r \xCN s \xCN t
\xCN u.$ $U:=\{x,x',x'',y,z\},$ the $ \xeb_{s}-$best elements
are $x,x',x'',$ they are contained in $X:=\{x,x',x'',z\}.$
$d_{s}(z,x)=\{s,t,u\},$
$d_{s}(z,x')=\{p,r,s\},$ $d_{s}(z,x'')=\{q,r\},$ so $x'' $ is one of the
best
elements closest to $z.$ $d(z,y)=\{q\},$ $d(y,x'')=\{r\},$ so $[z,x''
]=\{z,y,x'' \},$ $y \xce X,$ but $X$
is downward closed. $ \xcz $
\\[3ex]

\ee

\bfa

$\hspace{0.01em}$


\label{Fact Global-Dependent-S}

Work in the set variant.

Let $X \xEd \xCQ,$ $X$ $ \xec_{s}$-closed. Then

(1) $X$ does not necessarily contain all best elements.

Assume now that $X$ contains, in addition, all best elements. Then

(2) $X \xeb_{l,s} \xdC X$ does not necessarily hold.

(3) $X$ is $ \xCf (ui).$

(4) $X \xbe \xdd (\xdo)$ does not necessarily hold.

(5) $X$ is not necessarily a neighbourhood of the best elements.

(6) $X$ is an improving neighbourhood of the best elements.

\efa

\subparagraph{
Proof
}

$\hspace{0.01em}$


(1) See Example \ref{Example Dependent-1} (page \pageref{Example Dependent-1})

(2) See Example \ref{Example Dependent-2} (page \pageref{Example Dependent-2})

(3)
If there is $m \xbe X,$ $m \xce O$ for all $O \xbe \xdo,$ then by closure
$X=U,$ take $ \xdo_{i}:= \xCQ.$

For $m \xbe X$ let $ \xdo_{m}:=\{O \xbe \xdo:m \xbe O\}.$ Let $X':= \xcV
\{ \xcS \xdo_{m}:m \xbe X\}.$

$X \xcc X':$ trivial, as $m \xbe X \xcp m \xbe \xcS \xdo_{m} \xcc X'.$

$X' \xcc X:$ Let $m' \xbe \xcS \xdo_{m}$ for some $m \xbe X.$ It suffices
to show that $m' \xec_{s}m.$
$m' \xbe \xcS \xdo_{m}= \xcS \{O \xbe \xdo:m \xbe O\},$ so for all $O
\xbe \xdo $ $(m \xbe O \xcp m' \xbe O).$

(4) Consider Example \ref{Example Dependent-2} (page \pageref{Example
Dependent-2}),
let $dom(\xbd)=\{r,s\},$ $ \xbd (r)= \xbd (s)=0.$ Then $x,y \xcm \xbd,$
but $x' \xcM \xbd $ and $x \xbe X,$ $y \xce X,$
but there is no $z \xbe X,$ $z \xcm \xbd $ and $z \xeb y,$ so $X \xce \xdd
(\xdo).$

(5) See Example \ref{Example Dependent-3} (page \pageref{Example Dependent-3}).

(6) By Fact \ref{Fact Quality-Distance} (page \pageref{Fact Quality-Distance}),
(5).

$ \xcz $
\\[3ex]

\bfa

$\hspace{0.01em}$


\label{Fact Global-Dependent-Converse-S}

Work in the set variant

(1.1) $X \xeb_{l,s} \xdC X$ implies that $X$ is $ \xec_{s}$-closed.

(1.2) $X \xeb_{l,s} \xdC X$ $ \xch $ $X$ contains all best elements

(2.1) $X$ is $ \xCf (ui)$ $ \xch $ $X$ is $ \xec_{s}$-closed.

(2.2) $X$ is $ \xCf (ui)$ does not necessarily imply that $X$ contains all
$ \xec_{s}$-best
elements.

(3.1) $X \xbe \xdd (\xdo)$ $ \xch $ $X$ is $ \xec_{s}$-closed

(3.2) $X \xbe \xdd (\xdo)$ implies that $X$ contains all $
\xec_{s}$-best elements.

(4.1) $X$ is an improving neighbourhood of the $ \xec_{s}$-best elements $
\xch $ $X$ is
$ \xec_{s}$-closed.

(4.2) $X$ is an improving neighbourhood of the best elements $ \xch $ $X$
contains all best elements.

\efa

\subparagraph{
Proof
}

$\hspace{0.01em}$


(1.1) By Fact \ref{Fact Local-Global} (page \pageref{Fact Local-Global}).

(1.2) By Fact \ref{Fact General-Obligation} (page \pageref{Fact
General-Obligation}).

(2.1)
Let $O \xbe \xdo,$ then $O$ is downward closed (no $y \xce O$ can be
better than $x \xbe O).$
The rest follows from Fact \ref{Fact Subset-Closure} (page \pageref{Fact
Subset-Closure})  (3).

(2.2) Consider Example \ref{Example Dependent-1} (page \pageref{Example
Dependent-1}), $p$ is $ \xCf (ui)$
(formed in $U!),$ but $p \xcs X$ does
not contain $ \xCN pq.$

(3.1)
Let $X \xbe \xdd (\xdo),$ but let $X$ not be closed.
Thus, there are $m \xbe X,$ $m' \xec_{s}m,$ $m' \xce X.$

Case 1: Suppose $m' \xCq m.$ Let $ \xbd_{m}: \xdo \xcp 2,$ $
\xbd_{m}(O)=1$ iff $m \xbe O.$ Then $m,m' \xcm \xbd_{m},$
and there cannot be any $m'' \xcm \xbd_{m},$ $m'' \xeb_{s}m',$ so $X \xce
\xdd (\xdo).$

Case 2: $m' \xeb_{s}m.$ Let $ \xdo':=\{O \xbe \xdo:m \xbe O \xcj m'
\xbe O\},$ $dom(\xbd)= \xdo',$ $ \xbd (O):=1$ iff $m \xbe O$ for
$O \xbe \xdo'.$ Then $m,m' \xcm \xbd.$ If there is $O \xbe \xdo $ such
that $m' \xce O,$ then by $m' \xec_{s}m$ $m \xce O,$
so $O \xbe \xdo'.$ Thus for all $O \xce dom(\xbd).m' \xbe O.$ But then
there is no $m'' \xcm \xbd,$ $m'' \xeb_{s}m',$
as $m' $ is already optimal among the $n$ with $n \xcm \xbd.$

(3.2) Suppose $X \xbe \xdd (\xdo),$ $x' \xbe U-X$ is a best element,
take $ \xbd:= \xCQ,$ $x \xbe X.$
Then there must be $x'' \xeb x',$ $x'' \xbe X,$ but this is impossible as
$x' $ was best.

(4.1) By Fact \ref{Fact Quality-Distance} (page \pageref{Fact Quality-Distance})
, (4) all minimal elements
have incomparabel
distance. But if $z \xec y,$ $y \xbe X,$ then either $z$ is
minimal or it is above a minimal element, with minimal distance from $y,$
so $z \xbe X$
by Fact \ref{Fact Quality-Distance} (page \pageref{Fact Quality-Distance})  (3).

(4.2) Trivial.

$ \xcz $
\\[3ex]
\subsubsection{
The independent case
}

\label{Section Indep}

Assume now the system to be independent, i.e. all combinations of $ \xdo $
are present.

Note that there is now only one minimal element, and the notions of
Hamming neighbourhood of the best elements
and improving Hamming neighbourhood of the best elements coincide.

\bfa

$\hspace{0.01em}$


\label{Fact Global-Independent-S}

Work in the set variant.

Let $X \xEd \xCQ,$ $X$ $ \xec_{s}$-closed. Then

(1) $X$ contains the best element.

(2) $X \xeb_{l,s} \xdC X$

(3) $X$ is $ \xCf (ui).$

(4) $X \xbe \xdd (\xdo)$

(5) $X$ is a (improving) Hamming neighbourhood of the best elements.

\efa

\subparagraph{
Proof
}

$\hspace{0.01em}$


(1) Trivial.

(2) Fix $x \xbe X,$ let $y$ be closest to $x,$ $y \xce X.$ Suppose $x \xeB
y,$ then there must be
$O \xbe \xdo $ such that $y \xbe O,$ $x \xce O.$ Choose $y' $ such that
$y' $ is like $y,$ only $y' \xce O.$ If $y' \xbe X,$ then
by closure $y \xbe X,$ so $y' \xce X.$ But $y' $ is closer to $x$ than $y$
is, $contradiction.$

Fix $y \xbe U-$X. Let $x$ be closest to $y,$ $x \xbe X.$ Suppose $x \xeB
y,$ then there is $O \xbe \xdo $
such that $y \xbe O,$ $x \xce O.$ Choose $x' $ such that $x' $ is like
$x,$ only $x' \xbe O.$ By closure of $X,$
$x' \xbe X,$ but $x' $ is closer to $y$ than $x$ is, $contradiction.$

(3) By Fact \ref{Fact Global-Dependent-S} (page \pageref{Fact
Global-Dependent-S})  (3)

(4)
Let $X$ be closed, and $ \xdo' \xcc \xdo,$ $ \xbd: \xdo' \xcp 2,$
$m,m' \xcm \xbd,$ $m \xbe X,$ $m' \xce X.$
Let $m'' $ be such that $m'' \xcm \xbd,$ and for all $O \xbe \xdo -dom(
\xbd)$ $m'' \xbe O.$ This exists by
independence.
Then $m'' \xec_{s}m',$ but also $m'' \xec_{s}m,$ so $m'' \xbe X.$ Suppose
$m'' \xCq m',$ then
$m' \xec_{s}m'',$ so $m' \xbe X,$ contradiction, so $m'' \xeb_{s}m'.$

(5) Trivial by (1), the remark preceding this Fact, and
Fact \ref{Fact Global-Dependent-S} (page \pageref{Fact Global-Dependent-S}) 
(6).

\bfa

$\hspace{0.01em}$


\label{Fact Global-Independent-Converse-S}

Work in the set variant.

(1) $X \xeb_{l,s} \xdC x$ $ \xch $ $X$ is $ \xec_{s}$-closed,

(2) $X$ is $ \xCf (ui)$ $ \xch $ $X$ is $ \xec_{s}$-closed,

(3) $X \xbe \xdd (\xdo)$ $ \xch $ $X$ is $ \xec_{s}$-closed,

(4) $X$ is a (improving) neighbourhood of the best elements $ \xch $ $X$
is $ \xec_{s}$-closed.

\efa

\subparagraph{
Proof
}

$\hspace{0.01em}$


(1) Suppose there are $x \xbe X,$ $y \xbe U-$X, $y \xeb x.$ Choose them
with minimal distance.
If $card(d_{s}(x,y))>1,$ then there is
$z,$ $y \xeb_{s}z \xeb_{s}x,$ $z \xbe X$ or $z \xbe U-$X, contradicting
minimality. So $card(d_{s}(x,y))=1.$
So $y$ is among the closest elements of $U-X$ seen from $x,$ but then by
prerequisite
$x \xeb y,$ $contradiction.$

(2) By Fact \ref{Fact Global-Dependent-Converse-S} (page \pageref{Fact
Global-Dependent-Converse-S})  (2.1).

(3) By Fact \ref{Fact Global-Dependent-Converse-S} (page \pageref{Fact
Global-Dependent-Converse-S})  (3.1).

(4) There is just one best element $z,$ so if $x \xbe X,$ then $[x,z]$
contains all
$y$ $y \xeb x$ by Fact \ref{Fact Quality-Distance} (page \pageref{Fact
Quality-Distance})  (3).

$ \xcz $
\\[3ex]

The $ \xdd (\xdo)$ condition seems to be adequate only for the
independent situation,
so we stop considering it now.

\bfa

$\hspace{0.01em}$


\label{Fact Int-Union}

Let $X_{i} \xcc U,$ $i \xbe I$ a family of sets, we note the following
about closure under unions and intersections:

(1) If the $X_{i}$ are downward closed, then so are their unions and
intersections.

(2) If the $X_{i}$ are $ \xCf (ui),$ then so are their unions and
intersections.

\efa

\subparagraph{
Proof
}

$\hspace{0.01em}$


Trivial. $ \xcz $
\\[3ex]

We do not know whether $ \xeb_{l,s}$ is preserved under unions and
intersections,
it does not seem an easy problem.

\bfa

$\hspace{0.01em}$


\label{Fact Relativization}

(1) Being downward closed is preserved while going to subsets.

(2) Containing the best elements is not preserved (and thus neither the
neighbourhood property).

(3) The $ \xdd (\xdo)$ property is not preserved.

(4) $ \xec_{l,s}$ is not preserved.

\efa

\subparagraph{
Proof
}

$\hspace{0.01em}$


(4) Consider Example \ref{Example Not-Global} (page \pageref{Example
Not-Global}), and
eliminate $y$ from $U',$ then the closest to $x$ not in $X$ is $y',$
which is better.

$ \xcz $
\\[3ex]
\subsubsection{
Remarks on the counting case
}

\label{Section Counting-Case}

\br

$\hspace{0.01em}$


\label{Remark}

In the counting variant all qualities are comparabel. So if $X$
is closed, it will contain all minimal elements.

\er

\be

$\hspace{0.01em}$


\label{Example H-N-Local}

We measure distance by counting.

Consider $a:= \xCN p \xCN q \xCN r \xCN s,$ $b:= \xCN p \xCN q \xCN rs,$
$c:= \xCN p \xCN qr \xCN s,$ $d:=pqr \xCN s,$
let $U:=\{a,b,c,d\},$ $X:=\{a,c,d\}.$ $d$ is the best element,
$[a,d]=\{a,d,c\},$
so $X$ is an improving Hamming neighbourhood, but $b \xeb a,$ so $X
\xeB_{l,c} \xdC X.$

$ \xcz $
\\[3ex]

\ee

\bfa

$\hspace{0.01em}$


\label{Fact Local-H-N}

We measure distances by counting.

$X \xeb_{l,c} \xdC X$ does not necessarily imply that $X$ is an improving
Hamming
neighbourhood of the best elements.

\efa

\subparagraph{
Proof
}

$\hspace{0.01em}$


Consider Example \ref{Example Not-Global} (page \pageref{Example Not-Global}).
There $X \xeb_{l,c} \xdC X.$
$x' $ is the best element, and
$y' \xbe [x',x],$ but $y' \xce X.$ $ \xcz $
\\[3ex]
\section{
Neighbourhoods in deontic and default logic
}
\subsection{
Introduction
}

Deontic and default logic have very much in common.
Both have a built-in quality relation, where situations which satisfy
the deontic rules are better than those which do not, or closer to the
normal case in the default situation.

They differ in the interpretation of the result. In default logic, we want
to know what holds in the ``best'' or most normal situations, in
deontic logic, we want to characterize the
``good'' situations, and avoid paradoxa like the Ross-paradox.

Note that our treatment concern only obligations and defaults without
prerequisites, but this suffices for our purposes: to construct
neighbourhood semantics for both. When we work with prerequisites,
we have to consider the possibilities of branching into different
``extensions'', which is an independent problem.

We discussed MISE extensively in
Section \ref{Section Pref-Lim} (page \pageref{Section Pref-Lim}), so it will
not be necessary to
repeat the presentation.
\subsection{
Two important examples for deontic logic
}

\be

$\hspace{0.01em}$


\label{Example Ross-Paradox}

The original version of the Ross paradox reads: If we have the obligation
to post the letter, then we have the obligation to post or burn the
letter.
Implicit here is the background knowledge that burning the letter implies
not
to post it, and is even worse than not posting it.

We prefer a modified version, which works with two independent
obligations:
We have the obligation to post the letter, and we have the obligation to
water
the plants. We conclude by unrestricted weakening that we have the
obligation to
post the letter or $ \xCf not$ to water the plants. This is obvious
nonsense.

\ee

\be

$\hspace{0.01em}$


\label{Example Considerate-Assassin}

Normally, one should not offer a cigarette to someone, out of respect for
his health. But the considerate assassin might do so nonetheless, on the
cynical reasoning that the victim's health is going to suffer anyway:

(1) One should not kill, $ \xCN k.$

(2) One should not offer cigarettes, $ \xCN o.$

(3) The assassin should offer his victim a cigarette before killing him,
if $k,$ then $o.$

Here, globally, $ \xCN k$ and $ \xCN o$ is best, but among $k-$worlds, $o$
is better than $ \xCN o.$
The model ranking is $ \xCN k \xcu \xCN o \xeb \xCN k \xcu o \xeb k \xcu o
\xeb k \xcu \xCN o.$
\subsection{
Neighbourhoods for deontic systems
}

\ee

A set $ \xdr $ of deontic or default rules defines naturally quality and
distance
relations:

 \xEh
 \xDH
A situation (model) $m$ is better than a model $m' $ iff
$m$ satisfies more rules than $m' $ does.
``More'' can be defined by counting, or by the superset relation.
In both cases, we will note the relation here by $ \xeb.$
(See Definition \ref{Definition Quality} (page \pageref{Definition Quality}).)

 \xDH
The distance between two models $m,m' $ is the number - or the set - of
rules satisfied by one, but not by the other.
In both cases, we will note the distance here by $d.$
Given a distance, we can define ``between'':
a is between $b$ and $c$ iff $d(b,c)=d(b,a)+d(a,c)$
(in the case of sets, $+$ will be $ \xcv).$
See Definition \ref{Definition Hamming-Distance} (page \pageref{Definition
Hamming-Distance})  and
Definition \ref{Definition Between} (page \pageref{Definition Between}).

 \xEj

We have here in each case two variants of Hamming relations or distances.

With these ideas, we can define ``good'' sets $X$
in a number of ways:

Then, if $ \xdr $ is a family of rules, and
if $x$ and $x' $ are in the same subset $ \xdr' \xcc \xdr $ of rules,
then a
rule derived from $ \xdr $ should not separate them. More precisely, if
$x \xbe O \xbe \xdr \xcj x' \xbe O \xbe \xdr,$ and $D$ is a derived rule,
then $x \xbe D \xcj x' \xbe D.$

We think that being closed is a desirable property for obligations: what
is
at least as good as one element in the obligation should be ``in'', too.

But it may be a matter of debate which of the different possible
notions of neighbourhood should be chosen for a given deontic system.
It seems, however, that we should use the characterization
of neighbourhoods to describe acceptable situations.
Thus, we determine the ``best'' situations, and all
neighbourhoods of the best situations are reasonable approximations
to the best situations, and can thus be considered
``derived'' from the original system of obligations.
\clearpage
\chapter{
Conclusion and outlook
}
\section{
Conclusion
}

An important part of this book concerns the concept of independence.
One of the authors described it - mostly as questions - as
``homogenousness'' in his first book,
 \cite{Sch97-2}. But it took quite some time and detours to
find a reasonable, and in hindsight obvious, answer.
\subsection{
Semantic and syntactic interpolation
}

We usually try to decompose a logical problem - often formulated in the
language
of this logic - into a semantical part, and then the translation to
syntax.
This has proven fruitful in the past, and does so here, too. The reason is
that the semantical problems are often very different from the translation
problems. The latter concern usually definabilty questions, which tend to
be similar for various logical problems.

Here, we were able to see that $ \xCf semantical$ interpolation will
always
exist for monotonic or antitonic logics, but that the language and the
operators may be too weak to define the interpolants syntactically. In
contrast, semantical interpolants for non-monotonic logics need not
always exist. We detail this now briefly.
\subsection{
Independence and interpolation for monotonic logic
}

Independence is closely related to (semantic) interpolation, as a matter
of
fact, in monotonic logic, the very definition of validity is based
on independence, and guarantees semantic interpolation, also for
many-valued
logics, provided the order on the truth values is sufficiently strong, as
we saw
in Chapter \ref{Chapter Mod-Mon-Interpol} (page \pageref{Chapter
Mod-Mon-Interpol}), wheras the expressive
strength of the language
determines whether we can $ \xCf define$ the semantic interpolant (or
interpolants),
see the same chapter. We also saw that we often have an interval of
interpolants, and the upper and lower limits are $ \xCf universal$
interpolants
in the following sense: they depend only on $ \xCf one$ formula, plus the
set of
propositional variables common to both formulas, but $ \xCf not$ on the
second
formula itself.
\subsection{
Independence and interpolation for non-monotonic logic
}

Perhaps the central chapter of the book is
Chapter \ref{Chapter Size-Laws} (page \pageref{Chapter Size-Laws}), where we
connect interpolation of
non-monotonic
logics to multiplicative laws about abstract size. The entrance to this
is to see that logics like preferential logics define an abstract notion
of size by considering the exceptional cases as forming a small subset of
all cases, and, dually, the normal cases as a big subset of all cases.
Then, laws for non-monotonic logics can be seen as laws about $ \xCf
addition$
of small and big subsets, and about general well-behaviour of the notions
of big and small. E.g., if $X$ is a small subset of $Y,$ and $Y \xcc Y',$
then $X$
should also be a small subset of $Y'.$ These connections were
investigated
systematically first (to our knowledge) in
 \cite{GS09a}, see also  \cite{GS08f}. It was pointed out
there that the rule
of rational monotony does not fit well into laws about addition, and that
it has to be seen rather as a rule about independence.
It fits now well into our laws about $ \xCf multiplication,$ which express
independence for non-monotonic logics. It is natural that the laws we need
now for semantic interpolation are laws about multiplication, as we speak
about
products of model sets. Interpolation for non-monotonic logic has
(at least) three different forms, where we may mix the non-monotonic
consequence
relation $ \xcn $ with the classical relation $ \xcl.$ We saw that two
variants are
connected to such multiplicative laws, and, especially the weakest form
has
a translation to a very natural law about size. We can go on and relate
these laws to natural laws about preferential relations, when the
logic is preferential.

The problem of syntactic interpolation is the same as for the monotonic
case.

These multiplicative laws about size have repercussions beyond
interpolation, as
they also say what should happen when we change the language, e.g., have a
rule
$ \xbf \xcn \xbq $ in language $L,$ and now chnage to a bigger language
$L',$ whether we
can still expect $ \xbf \xcn \xbq $ to hold. This seems a trivial problem,
it is not,
and somehow seems to have escaped attention so far.
\subsection{
Neighbourhood semantics
}

The concluding chapter of the book concerns neighbourhood semantics,
see Chapter \ref{Chapter Neighbourhood} (page \pageref{Chapter Neighbourhood}).
Such semantics are ubiquitous
in
non-classical logic, they can be found as systems of ever better
sets in the limit version of preferential logics, as a semantics
for deontic and default logics, for approximative logics, etc.
We looked at the common points, how to define them, and what properties
to require for them. This chapter should be seen as a toolbox, where
one finds the tools to construct the semantics one needs for the
particular
case at hand.
\section{
Outlook
}

\label{Section Outlook}

We think that further research should concern the dynamic aspects of
reasoning,
like iterated revision, revising one non-monotonic logic with another
non-monotonic logic.

Moreover, it seems to us that any non-classical logic (which is not an
extension of the former, like modal logic, but diverges in its results
from classical logic) needs a $ \xCf justification,$ so such logics do not
only consist of language, proof theory, and semantics, but of
language, proof theory, semantics, $ \xCf and$ $ \xCf justification.$
\subsection{
The dynamics of reasoning
}

So far, most work on non-monotonic and related logics concern one step
in a reasoning process only. Notable exceptions are, e.g.,
 \cite{Spo88}, and  \cite{DP94},  \cite{DP97}.

It seems quite certain that there is no universal formalism, if, e.g., in
a theory revision task,
we are given $ \xbf,$ and then $ \xCN \xbf $ as information, we can
imagine situations
where we should come back to the original state, and others, where this
should
not be the case.

So, the dynamics of reasoning need further investigation.

We took already
some steps here, when we investigated generalized revision (see
Chapter \ref{Chapter Mod-Uni-Cond} (page \pageref{Chapter Mod-Uni-Cond}),
Section \ref{Section Repr-TR-Up} (page \pageref{Section Repr-TR-Up})),
as the results can be applied to revising one preferential logic with
another one. (Usually, such logics will not have a natural ranked order,
so traditional revision will not work.) But we need more than tools,
we need satisfactory systems.
\subsection{
A revision of basic concepts of logic: justification
}

$ \xCO $

Some logics like inductive logics
(``proving'' a theory from a limited number of cases),
non-monotonic logics, revision and
update logics go beyond classical logic, they allow to derive formulas
which
cannot be derived in classical logic. Some might also be more modest,
allowing
less derivations, and some might be a mixture, e.g. approximative logics,
allowing to derive some formulas
which cannot be derived in classical logic, and not allowing to derive
other
formulas which can be derived in classical logic.

Let us call all those logics
``bold logics''.

Suppose that we agree that classical logic corresponds to
``truth''.

But then we need a $ \xCf justification$ to do other logic than classical
logic,
as we know or suspect - or someone else knows or suspects - that our
reasoning
is in some cases false. (Let us suppose for simplicity that we know this
erroneousness ourselves.)

Whatever this justification may be, we have now a fundamentally new
situation.

Classical logic has language, proof theory, and semantics. Representation
theorems say that the latter correspond. Non-monotonic logic also has
(language and) proof theory, and semantics. But something is missing:
the justification - which we do $ \xCf not$ need for classical logic, as
we do not
have any false reasoning to justify.

Thus,
 \xEI
 \xDH
classical logic consists of
 \xEh
 \xDH
language (variables and operators),
 \xDH
proof theory,
 \xDH
semantics.
 \xEj

 \xDH
bold logic consists of
 \xEh
 \xDH
language (variables and operators),
 \xDH
proof theory,
 \xDH
semantics,
 \xDH
justification.
 \xEj

 \xEJ

If a bold logic has no justification - whatever that may be - it is just
foolishness, and the bolder it is (the more it diverges from classical
logic), the more foolish it is.

So let us consider justifications - in a far from exhaustive list.

 \xEh
 \xDH
First, on the negative side, costs.
 \xEh
 \xDH
A false result has a cost. This cost depends
on the problem we try to solve. Suppose we have a case
``man, blond''. Classifying this case falsely as
``man, black hair'', has a different cost when we try to
determine the amount of hair dyes to buy, and when we are on the
lookout for a blond serial killer on the run.
 \xDH
Calculating our bold logic has a cost, too (time and space).
Usually, this will also depend on the case, the cost is not the same
for all cases.
E.g., let $T=p \xco (\xCN p \xco q),$ then the cost to determine whether
$m \xcm T$ is smaller
for $p-$models, than for $ \xCN p-$models, as we have to check now in
addition $q.$

In addition, there may be a global cost of calculation.

 \xEj

 \xDH
Second, on the positive side, benefits:

 \xEh
 \xDH
Classical logic also has its costs of calculation, similar to the above.
 \xDH
In some cases, classical logic may not be strong enough to decide the case
at hand. Hearing a vague noise in the jungle may not be enough to decide
whether it is a lion or not, but we climb up the tree nonetheless. Bold
logic allows us to avoid desaster, by precaution.
 \xDH
Parsimony, elegance, promises to future elaboration, may also be
considered
benefits.
 \xEj

 \xEj

We can then say that a bold logic is justified, iff the benefits
(summarized over all cases to consider)
are at least as big as the costs
(summarized over all cases to consider).

Diverging more from classical logic incurs a bigger cost, so the bolder
a logic is, the stronger its benefits must be, otherwise it is not
rational to choose it for reasoning.

When we look at preferential logic and its abstract semantics of
``big'' and ``small'' sets
(see, e.g.,  \cite{GS09a},  \cite{GS08f}),
we can consider this semantics as being
an $ \xCf implicit$ justification:
The cases wrongly treated are together just a
``small'' or ``unimportant'' subset of all cases.
(Note that this says nothing about the benefits.)
But we want to make concepts clear, and explicit, however they
may be treated in the end.

$ \xCO $

See also  \cite{GW08} for related reflections.
\clearpage

$ \xCO $

\begin{theindex}

\addcontentsline{toc}{chapter}{Index}

\item $(*CCL)$, 93
\item $(*Con)$, 93
\item $(*Equiv)$, 93
\item $(*Loop)$, 93
\item $(*Succ)$, 93
\item $(-)$, 41
\item $(1*s)$, 189, 190
\item $(2*s)$, 189, 190
\item $(< \xbo*s)$, 189, 190
\item $(AND)$, 52, 65
\item $(AND_1)$, 52, 190
\item $(AND_{2})$, 190
\item $(AND_n)$, 52, 190
\item $(AND_{ \xbo })$, 190
\item $(B1)$, 25, 166
\item $(B2)$, 25, 166
\item $(B3)$, 25, 166
\item $(B4)$, 25, 166
\item $(B5)$, 25, 166
\item $(CCL)$, 52, 59, 85
\item $(CM)$, 53
\item $(CM_2)$, 53, 190
\item $(CM_n)$, 53, 190
\item $(CM_{ \xbo })$, 190
\item $(CP)$, 52, 190
\item $(CUM)$, 53, 59
\item $(CUT)$, 52
\item $(disjOR)$, 52, 190
\item $(dp)$, 43
\item $(DR)$, 53
\item $(EE1)$, 94
\item $(EE2)$, 94
\item $(EE3)$, 94
\item $(EE4)$, 94
\item $(EE5)$, 94
\item $(eM\xdf)$, 53, 53, 189, 189, 190
\item $(eM\xdi)$, 52, 53, 189, 189, 190
\item $(iM)$, 52, 52, 189, 189, 190
\item $(I\xcv disj)$, 52
\item $(I_1)$, 52
\item $(I_2)$, 52, 53
\item $(I_n)$, 52, 53
\item $(I_\xbo)$, 52, 53
\item $(K*1)$, 94
\item $(K*2)$, 94
\item $(K*3)$, 94
\item $(K*4)$, 94
\item $(K*5)$, 94
\item $(K*6)$, 94
\item $(K*7)$, 94
\item $(K*8)$, 94
\item $(K-1)$, 94
\item $(K-2)$, 94
\item $(K-3)$, 94
\item $(K-4)$, 94
\item $(K-5)$, 94
\item $(K-6)$, 94
\item $(K-7)$, 94
\item $(K-8)$, 94
\item $(LLE)$, 52, 59, 85
\item $(Log=')$, 53
\item $(Log \xcv)$, 53, 59
\item $(Log \xcv')$, 53, 59
\item $(Log \xFO)$, 53
\item $(n*s)$, 189, 190
\item $(NR)$, 190
\item $(Opt)$, 52, 189, 190
\item $(OR)$, 52, 65
\item $(OR_{2})$, 190
\item $(OR_{n})$, 190
\item $(OR_{ \xbo })$, 190
\item $(PR)$, 52, 59
\item $(PR')$, 190
\item $(RatM)$, 53, 59, 190
\item $(RatM=)$, 53, 59
\item $(REF)$, 52
\item $(ResM)$, 53
\item $(RW)$, 52, 190
\item $(RW)+$, 59
\item $(SC)$, 52, 59, 59, 85, 190
\item (ui), 219
\item $(wCM)$, 53, 190
\item $(wOR)$, 52, 190
\item $(\xbL \xcs)$, 63
\item $(\xbm =)$, 53, 59
\item $(\xbm =')$, 53, 59
\item $ (\xbm CM) $, 53, 190
\item $(\xbm CUM)$, 53, 59, 81, 85
\item $ (\xbm CUT) $, 52
\item $(\xbm disjOR)$, 52, 190
\item $(\xbm dp)$, 43, 59
\item $(\xbm OR)$, 52, 190
\item $(\xbm PR)$, 52, 59, 190
\item $(\xbm PR')$, 52
\item $(\xbm RatM)$, 53, 59, 190
\item $(\xbm ResM)$, 53
\item $(\xbm wOR)$, 52, 190
\item $(\xbm \xbe)$, 53, 59
\item $(\xbm \xcc)$, 52, 59, 85
\item $ (\xbm \xcc \xcd) $, 53, 59, 81, 85
\item $(\xbm \xCQ)$, 52, 59
\item $(\xbm \xCQ fin)$, 52, 59
\item $(\xbm \xcv)$, 53, 59
\item $(\xbm \xcv')$, 53, 59
\item $(\xbm \xFO)$, 53, 59
\item $ (\xcc \xcd) $, 53, 59, 85
\item $(\xCd)$, 189, 190
\item $(\xcs)$, 41
\item $(\xcS)$, 41
\item $(\xcs)$, 85
\item $(\xcV)$, 41
\item $(\xcv)$, 41, 59
\item $(\xdC)$, 41
\item $(\xdf \xcv disj)$, 189
\item $(\xdf_1)$, 189
\item $(\xdf_2)$, 189
\item $(\xdf_n)$, 189
\item $(\xdf_\xbo)$, 189
\item $(\xdi \xcv disj)$, 189
\item $(\xdi_1)$, 189
\item $(\xdi_2)$, 189
\item $(\xdi_n)$, 189
\item $(\xdi_\xbo)$, 189
\item $(\xdm^{++})$, 53, 189, 190
\item $(\xdm^+ \xcv disj)$, 189
\item $(\xdm^{+}_{n})$, 189
\item $(\xdm^+_\xbo)$, 53, 189
\item $(\xeA_1)$, 189
\item $(\xeA_2)$, 189
\item $(\xeA_n)$, 189
\item $(\xeA_{\xbo})$, 189
\item $(\xfA Con)$, 93
\item $(\xfA Loop)$, 93, 93
\item $(\xfA Succ)$, 93, 93
\item 1-copy, 56
\item 1 copy, 59
\item $<$, 41
\item AND, 190
\item arrows, 71
\item attacking structure relative to $ \xbh $ representing $ \xbr $, 83
\item $a \xcT b$, 41
\item $A \xfA_{d}B$, 92
\item between, 222
\item big set, 42
\item card, 41
\item Cautious Monotony, 53, 190
\item ceteris paribus, 224
\item Classical Closure, 52
\item cofinal, 63
\item collective variant, 92
\item $Con(.)$, 43
\item consequence relation, 56
\item Consistency Preservation, 52
\item Contraction, 94
\item copies, 55
\item Cumulativity, 53
\item $d(\xba)$, 73
\item $D(\xba)$, 73
\item $DC'$, 25, 166
\item $definability$ $ \xCf preserving$, 43
\item $definable$, 43, 63
\item Definition AGM revision, 90
\item Definition Algebraic basics, 41
\item Definition Distance, 91
\item Definition Distance representation, logics, 92
\item Definition Distance representation, semantics, 92
\item Definition Essentially smooth, 82
\item Definition Generalized preferential structure, 72
\item Definition IBRS, 70
\item Definition Individual/collective distance, 91
\item Definition Level-n-Arrow, 73
\item Definition Logic, basics, 42
\item Definition Logical conditions, 44
\item Definition Preferential logics, 56
\item Definition Preferential structure, 55
\item Definition Ranked relation, 58
\item Definition Relations, 41
\item Definition Size rules, 173
\item Definition Smoothness, 56
\item Definition Totally smooth, 82
\item Definition TR-Dist-N, 92
\item Definition Valid-Arrow, 75
\item Diagram Attacking structure, 83
\item Diagram Bubble, 60
\item Diagram Essential smoothness I, 80
\item Diagram Essential smoothness II, 80
\item Diagram High-Cum, 88
\item Diagram IBRS, 71
\item Diagram Independ, 108
\item Diagram information transfer, 110
\item Diagram Mod-Over, 13
\item Diagram Mon-Int, 134
\item Diagram Mul-Add, 174
\item Diagram Mul-Base-2, 199
\item Diagram Mul-Commut, 33
\item Diagram Mul-Prod, 176
\item Diagram Nml-Int-Ways, 201
\item Diagram O-D, 73
\item Diagram Update-S, 118
\item Diagram Valid-Valid, 76
\item Diagram Xcl-Xcn, 196
\item Diagram Xcn-Xcl, 194
\item $d_{c}(x,y)$, 221
\item $d_{s}(x,y)$, 221
\item Epistemic entrenchment, 94
\item essentially smooth, 59, 83
\item Example O-D, 73
\item exception set, 59
\item external monotony for filters, 189
\item external monotony for ideals, 189
\item $F(\xdl)$, 42
\item filter, 173
\item Filter, 189
\item FOL, 40
\item generalized preferential structure, 72
\item $GTS'$, 25, 166
\item Hamming distance, 221, 221, 228
\item Hamming quality, 220
\item IBR, 71
\item ideal, 173
\item Ideal, 189
\item Improving proportions, 189, 190
\item independent, 220
\item individual variant, 92
\item information transfer, 110
\item injective, 56
\item internal monotony, 189
\item irreflexive, 56
\item $K*\xbf $, 94
\item $K-\xbf $, 94
\item Keeping proportions, 189, 190
\item kill, 56
\item $K_{ \xcT }$, 90
\item $L$, 42
\item Left Logical Equivalence, 52
\item level $n$ arrow, 73
\item level $n$ structure, 73
\item $M$, 42
\item $M(T)$, 42
\item $M(\xbf)$, 42
\item medium size set, 42
\item minimal element, 55
\item minimize, 56
\item minimizing initial segment, 62
\item MISE, 62
\item Monotony, 189, 190
\item $M_{T}$, 42
\item $M_{ \xdl }$, 42
\item NML, 40
\item nodes, 71
\item normal characterization, 58, 59
\item Notation FOL-Tilde, 40
\item $O(\xba)$, 73
\item $o(\xba)$, 73
\item Optimal proportion, 189, 190
\item OR, 190
\item points, 71
\item preferential model, 55
\item preferential structure, 55
\item principal filter, 42
\item Proposition Pref-Representation, 58
\item pseudo-distance, 91
\item ranked, 58, 59
\item rationality postulates, 90
\item Rational Monotony, 53, 190
\item $RBC'$, 25, 166
\item reactive, 59
\item Reflexivity, 52
\item representable, 92
\item Restricted Monotony, 53
\item Revision, 94
\item Right Weakening, 52
\item Robustness of proportions, 189, 190
\item Robustness of $\xdm^+$, 189, 190
\item singletons, 59
\item small set, 42
\item smooth, 57, 59
\item smoothness, 78
\item $SRM'$, 25, 166
\item Supraclassicality, 52
\item $T*_{d}T' $, 92
\item $Th(m)$, 43
\item $Th(X)$, 43
\item theory, 42
\item theory revision, 63
\item totally smooth, 82
\item transitive, 56
\item $T \xcm_{ \xdm } \xbf $, 56
\item $T \xco T' $, 43
\item $UC'$, 25, 166
\item $V$, 41
\item $v(\xdl)$, 42
\item valid $X-to-Y$ arrow, 75
\item valid $X \xch Y$ arrow, 75
\item $x \xCq y$, 220
\item $X \xDc_{l}Y$, 224
\item $x \xeb y$, 220
\item $x \xeb Y$, 224
\item $X \xeb_{l}Y$, 224
\item $x \xec y$, 220
\item $X \xes X' $, 79
\item $X \xfA Y$, 91
\item $X \xfB Y$, 92
\item $x \xFO_{d}X$, 222
\item $[x,z]_{c}$, 222
\item $[x,z]_{d}$, 222
\item $[x,z]_{s}$, 222
\item $ \ol{T}$, 43
\item $ \ol{ \ol{T} }$, 43
\item $ \xba:x \xcp y$, 55
\item $ \xba: \xBc x,i \xBe  \xcp \xBc y,i \xBe $, 55
\item $ \xBc M,\xdn (M) \xBe $, 25, 166
\item $ \xBc x,y,z \xBe _{c}$, 222
\item $ \xBc x,y,z \xBe _{d}$, 222, 222
\item $ \xBc x,y,z \xBe _{s}$, 222
\item $ \xBc  \xba,k \xBe :x \xcp y$, 55
\item $ \xBc  \xba,k \xBe: \xBc x,i \xBe  \xcp  \xBc y,i \xBe $, 55
\item $ \xbe -$case, 223
\item $ \xbf_{M}$, 219
\item $ \xbL (X)$, 63
\item $ \xbl (\xba)$, 73
\item $ \xbm (X)$, 78
\item $ \xbm_{ \xdm }$, 55
\item $ \xbP $, 41
\item $ \xbP (\xdD, \xba)$, 84
\item $ \xbP (\xdO, \xba)$, 84
\item $ \xcb $, 41
\item $ \xcc $, 41
\item $ \xcd $, 41
\item $ \xcf $, 41
\item $\xcg 1$ copy, 59
\item $ \xck $, 41
\item $ \xcl $, 43
\item $ \xcm $, 43
\item $ \xdC $, 41
\item $ \xDc_{l,s}$, 224
\item $ \xdD (\xba)$, 84
\item $ \xdd (\xdo)$, 220
\item $ \xdD_{ \xdl }$, 43
\item $ \xdD_{ \xdl }-$smooth, 57
\item $ \xdm^+ $, 189
\item $ \xdO (\xba)$, 84
\item $ \xdp $, 41
\item $ \xdy -$essentially smooth, 83
\item $ \xdy -$smooth, 57
\item $\xdy-$totally smooth, 82
\item $\xeA$, 24, 166, 189
\item $ \xeb^{*}$, 41
\item $ \xeb_{l,c}$, 224
\item $ \xeb_{l,s}$, 224
\item $ \xex $, 41
\item $ \xFO $, 41
\item $ \xFO_{c}$, 222
\item $ \xFO_{s}$, 222

\end{theindex}

\end{document}